\documentclass[11pt]{article}
\pdfoutput=1
\usepackage[utf8]{inputenc}
\usepackage[T1]{fontenc}
\usepackage{lmodern}
\usepackage{microtype}
\usepackage[margin=1in]{geometry}
\usepackage[shortlabels]{enumitem}

\usepackage{twisted-double-macros}
\numberwithin{equation}{section}

\usepackage[breaklinks,hidelinks,bookmarksnumbered]{hyperref}
\usepackage[capitalize,noabbrev,nameinlink]{cleveref}

\usepackage{thmtools}
\makeatletter % XXX: https://github.com/muzimuzhi/thmtools/issues/75
\@ifundefined{newcounteralias}{}{%
  \renewcommand\thmt@autorefsetup{\@xa\def\csname\thmt@envname autorefname\@xa\endcsname\@xa{\thmt@thmname}}%
}
\makeatother
\declaretheorem[within=section]{theorem}
\declaretheorem[sibling=theorem]{proposition}
\declaretheorem[sibling=theorem]{corollary}
\declaretheorem[sibling=theorem]{lemma}
\declaretheorem[sibling=theorem,style=definition,qed=\qedsymbol]{definition}
\declaretheorem[sibling=theorem,style=definition,qed=\qedsymbol]{construction}
\declaretheorem[sibling=theorem,style=remark,qed=\qedsymbol]{remark}
\declaretheorem[sibling=theorem,style=remark,qed=\qedsymbol]{example}
\declaretheorem[within=section]{conjecture}
\crefname{construction}{Construction}{Constructions}
\Crefname{construction}{Construction}{Constructions}
\crefname{conjecture}{Conjecture}{Conjectures}
\Crefname{conjecture}{Conjecture}{Conjectures}

\usepackage{nicematrix}
\NiceMatrixOptions{cell-space-limits=3pt}
\newenvironment{dblArray}[1]{\begin{NiceArray}{#1}[hvlines]}{\end{NiceArray}}

\usepackage[style=alphabetic,maxbibnames=10]{biblatex}
\renewbibmacro{in:}{}

\title{Twisted double functors and loosely discrete opfibrations}
\author{Michael Lambert \and David Jaz Myers \and Evan Patterson}
\date{}

\begin{document}
\maketitle

\begin{abstract}
  Various situations in the theory and applications of double categories,
  ranging from a loose Yoneda theory and loose compact closure to
  double-operadic systems theory, require a notion of double copresheaf in which
  the action is by loose morphisms rather than tight ones. In this paper, we
  develop and compare several models for loose copresheaves on double
  categories. First, we introduce a new notion of morphism between double
  categories, called \emph{twisted double functors}, which send tight morphisms
  to loose morphisms and vice versa, and use these to define \emph{twisted
    copresheaves}. We exhibit numerous examples of twisted double functors,
  starting with the twisted Hom functor and the twisted representables on a
  double category. Corresponding to this functorial notion of loose copresheaf
  is a fibrational one, an internal version of a discrete opfibration that we
  call a \emph{loosely discrete opfibration}. We prove that twisted copresheaves
  and cloven loosely discrete opfibrations are equivalent via an elements
  construction. Finally, we compare twisted bimodules with double categories
  over the walking loose arrow, or \emph{double barrels}, via a collage
  construction.
\end{abstract}

\tableofcontents

\section{Introduction}

The Yoneda embedding and accompanying Yoneda lemma are cornerstones of category
theory, often capable of reducing calculations in an arbitrary category to the
category of sets. In a now seminal paper \cite{pare2011}, Robert Paré introduced
a Yoneda theory for double categories. The paper is important not only for its
double-categorical Yoneda lemma, but also for the structures that it introduces
to even formulate this result. The Hom double functor, and the representables
derived from it, naturally take the form of span-valued lax double functors. The
receiving structure for a double-categorical Yoneda embedding must therefore be
a double category whose objects are lax double functors. It turns out that this
double category is, in general, only virtual. Its loose arrows and cells are
modules and modulations, concepts which Paré introduces, building on their
bicategorical analogues \cite{wood2003}. In the first and last authors' work on
double-categorical logic \cite{lambert2024}, the virtual double category of lax
functors plays an essential role as the virtual double category of models of a
double theory.

On objects, Paré's Hom double functor for a double category $\dbl{D}$ sends a
pair of objects in $\dbl{D}$ to the set of tight morphisms between them. Of
course, a double category has two kinds of morphisms, tight and loose. More
precisely, then, Paré's theory is a \emph{tight} Yoneda theory for double
categories. An impetus for this work is to develop the theory in the other
direction, a \emph{loose} Yoneda theory. Here too we expect that formulating the
results will uncover important new structures on double categories. In fact,
this paper is solely concerned with investigating these new structures, ending
roughly where a loose Yoneda theory can begin.

To see the obstacles involved, consider the putative loose Hom functor on a
double category $\dbl{D}$. On objects, it should send a pair of objects in
$\dbl{D}$ to the set (or category) of proarrows between them. It should send a
pair of proarrows in $\dbl{D}$ to the mapping that acts by pre- and
post-composition with those proarrows, and a pair of arrows in $\dbl{D}$ to the
span (or profunctor) of cells in $\dbl{D}$ with those arrows as source and
target. But notice the peculiarity of this assignment. In the receiving double
category, be it $\Span$ or $\Prof$, it is the tight morphisms that are
function-like and the loose morphisms that are span-like. Such a loose Hom
functor must therefore be \emph{twisted}, exchanging the tight and loose
directions by sending tight morphisms to loose morphisms and vice versa. And
unless $\dbl{D}$ happens to be strict, this cannot be accomplished by simply
taking an ordinary double functor out of the transpose of $\dbl{D}$, as the
transpose does not make sense.

We thus begin the paper by axiomatizing the concept of a \emph{twisted double
  functor}. Considering the loose Hom, a twisted double functor must be lax in
the tight-to-loose direction, as those composition comparisons will encode the
tight composition of cells, dual to how the composition comparisons in the
(tight) Hom double functor encode the loose composition of cells
\cite[\S{2.1}]{pare2011}. On the other hand, since composition of loose
morphisms in a (pseudo) double category obeys the category axioms only up to
natural isomorphism, a twisted functor must also be at least pseudo in the
loose-to-tight direction. So, unlike ordinary double functors, twisted double
functors allow composition comparisons in both directions.

This requirement has implications for the receiving double category. Since the
only cells in $\Span$ bounded by loose identities are themselves loose
identities, to obtain nontrivial loose-to-tight comparisons, we must take
$\Prof$ as the receiving structure for the twisted Hom functor. In contrast, the
tight Hom functor can be regarded as a span-valued lax functor, as Paré does, or
equivalently as a profunctor-valued \emph{normal} lax functor, as Fr{\"o}hlich
and Moser do in their alternative formulation of a tight Yoneda theory
\cite{frohlich2024}.

In short, then, a \emph{twisted copresheaf} will be a twisted normal lax functor
valued in profunctors. Our first major construction, the \emph{twisted Hom},
will be a twisted copresheaf on the product of a double category and its loose
dual. It sends each pair of objects to the category of loose morphisms between
them and globular cells between those.

From the twisted Hom many additional examples of twisted functors are readily
generated. In general, a twisted double functor can be precomposed with an
ordinary double functor to yield another twisted functor. Particularly important
are the \emph{twisted representables} obtained by precomposing the twisted Hom
with a constant in one slot. Examples of twisted representables include the
indexed families of slice categories, coslice categories, and subobject posets.
Categories of ordinary copresheaves are also indexed by a twisted functor,
valued in adjunctions. In fact, when the indexing double category is an
equipment, any twisted copresheaf, and so in particular any twisted
representable, automatically furnishes a bi-indexed category, i.e., a
pseudo-functorially indexed family of adjunctions. From the viewpoint of
categorical logic, we thus regard twisted copresheaves as generalizing
structures like hyperdoctrines by putting both sides on an equal footing: rather
than having a 1-category in the domain and a 2-category in the codomain, both
domain and codomain are double categories.

As further examples, the composition of generalized Moore machines via lenses
and of structured cospan systems via hypergraph category operations can be
recovered as the action of a twisted representable (suitably restricted) on
double categories of structured (co)spans. Such applications to systems theory
are discussed further in \cref{intro:systems}. In a different direction, the
twisted Hom functor is central to the definition of a \emph{twisted adjunction}
between double categories. The need for such a definition arose in our
recognition of the centrality of \emph{compactness} and \emph{$*$-autonomy} for
double categories but also that of the inadequacy of the existing theory for a
proper treatment. See \cref{intro:adjunctions} for further discussion.

In order to better understand twistedness, we treat twisted copresheaves as
representations of a fixed double category and investigate the fibration concept
corresponding to such representations. This process continues a pattern of
development whereby a fibration concept is extracted from an elements
construction appropriate to the given type of representation of a 2-, bi- or
double category. This move has formerly recovered \emph{discrete double
  fibrations} \cite{lambert2021} from the original development of ``tight double
copresheaves'' \cite{pare2011}; then \emph{double fibrations} in
\cite{cruttwell2022} as a generalization; and finally \emph{discrete
  2-fibrations} in \cite{lambert2024a}. Here, the difference is that our
representations are profunctor-valued and are in our sense \emph{twisted},
adding certain layers of complication. For this reason our elements
correspondence and the representation theorem are only at the level of objects,
that is, between fibrations and representations, without considering the full
two-dimensional structure.

In a not altogether unexpected twist, the fibration concept recovered here is
essentially the transpose of the discrete concept discovered in
\cite{lambert2021}. That is, we show that profunctor-valued twisted
representations correspond to a weak version of \emph{internal discrete
  opfibrations} in the 2-category of categories. Internal discrete opfibrations
can ordinarily be made sense of in any category with sufficient finite limits.
Such a thing is an internal functor for which the appropriate square involving
internal domain morphisms presents the pullback of the object part of the
functor along the domain morphism \emph{up to isomorphism}. This is precisely
the approach taken in \cite{johnstone2002,maclane1992} where discrete
opfibrations in an elementary topos are defined to be such internal functors and
are shown to be algebras for a pull-push monad on the slice of the topos over
the object of objects of the base category. They are seen as an elementary
axiomatization of copresheaves in the topos.

Simply instantiating this definition in the 2-category of categories would be
too strict, however. For on the one hand, the resulting lifting properties are
simply not enjoyed by our examples. And on the other hand, morally speaking,
asking that such a square presents the pullback in the sense of
\emph{isomorphism} is too strong for categorical structures or properties.
Rather, we ask that the accompanying square presents the appropriate pullback
only up to weak equivalence of categories. This gives all the appropriate
lifting properties of our examples and implies that such double functors are
\emph{discrete in the loose direction}. For this reason we dub these double
functors \emph{loosely discrete opfibrations}.

To obtain a twisted copresheaf from such a loosely discrete opfibration, certain
algebraic constructions need to be introduced in the form of
\emph{cleavages} for such opfibrations. This might be unexpected as ordinary
discrete opfibrations are always cloven for trivial reasons. However, our
categorification is only discrete in one direction and in a particular sense. A
cleavage satisfying further algebraic assumptions must be supplied in order to
return to representations.

That cleavages for genuine fibrations supply algebraic structure was certainly
known in the 1960s, at least as reported in \cite{gray1966}. The cotensor of the
base category with the two-element ordinal category provides a pull-push monad
for which any cloven fibration over that base is an algebra. This approach was
generalized in \cite{lambert2024a} to show that cloven discrete 2-fibrations are
certain algebras. Thus, we already have an instance where discreteness in one
direction nonetheless requires a cleavage to supply the algebraic data needed
for the construction of a representation.

In the case of loosely discrete opfibrations, a cleavage is one half of a
LARI split equivalence. That is, to be a loosely discrete opfibration is to say
that the canonical functor from the category of proarrows and cells of the total
double category to the pullback of the object part of the double functor along
the source functor is a weak equivalence
(\cref{def:loosely-discrete-opfibration}). A cleavage is a choice of
pseudo-inverse that supplies a left-adjoint-right-inverse (LARI)
(\cref{def:cleavage-for-disc-opfibration}). The composition of the cleavage
functor together with the target functor in the total double category is then
asked to be the structure morphism of a pseudo-algebra for the pull-push monad
associated to the base double category. The pseudo-algebra cell is exactly what
proves the comparison cells in the associated twisted profunctor-valued
representation. And the algebra associativity and unitality assumptions are
exactly what makes the two constructions appropriately round-trip, showing that
there is a correspondence between profunctor-valued representations of a double
category and cloven loosely discrete opfibrations over it.

We also compare \emph{twisted bimodules} --- certain profunctor-valued twisted functors --- with the double barrels considered in \cite{brown2025} as models of loose bimodules in \cref{sec:double-barrel}. A \emph{double barrel} is a double functor into the walking loose arrow, generalizing the well-known presentation of a bimodule (i.e., a profunctor) of categories as a functor into the walking arrow. Any double barrel may be restricted along the twisted Hom functor of its domain to obtain an associated twisted bimodule of \emph{sections} of the double barrel. We show that this assignment is 2-functorial as a corollary of an investigation into the 2-functoriality of the twisted Hom construction (\cref{prop:embedding-into-twisted-copresheaves}). While the 2-category of double barrels is simply the slice of the 2-category of double categories, lax functors, and tight transformations over the walking loose arrow, the 2-category of twisted bimodules requires us to define an interesting notion of \emph{relative modification} (\cref{defn:relative.modification}) for the 2-cells.

In \cref{subsection:collage}, we define the \emph{collage} construction for twisted bimodules, giving us a 2-functor back from twisted bimodules to double barrels. We conjecture that these two 2-functors are mutually inverse.

\subsection{Applications to loose compact closure} \label{intro:adjunctions}

It is observed in \cite{carboni1987} that any \emph{bicategory of relations} is
compact closed with a trivial duality involution on objects. This is carried
over to \emph{double categories of relations} since the loose structures of such
double categories are in particular bicategories of relations
\cite{lambert2022}. The defining adjunctions are phrased in terms of Hom
categories and package proarrow structure, like our twisted Hom functor.
However, owing to the assumption of being \emph{locally posetal}, bicategories
of relations are more like monoidal 2-categories and not properly bicategorical.
A genuine bicategorical generalization of these notions is carried out in the
important paper \cite{daystreet1997}, which defines autonomy and duality in
terms of parameterized biadjunctions of homomorphisms of bicategories. Homs in
this environment are of course categories and \emph{a priori} package the loose
structure of bicategories, since that is all there is. Important results in this
direction are then mere applications of the \emph{bicategorical Yoneda lemma},
which simply has no analogue in the theater of category-valued representations
of double categories, hence our project.

Given that there already is a well-founded bicategorical theory of such
parameterized adjunctions, why do we pursue a double-categorical version? Our
general viewpoint is genuinely double-categorical. We view as fundamental the
theoretical importance and utility of axiomatizing two different but interacting
directions of morphisms, one \emph{functional}, and the other \emph{relational}
or \emph{module-like}. Further, we view structured double categories such as
equipments in the framework of \emph{fibrational semantics}, that is, as natural
models of logics over type theories where arrows model terms and proarrows model
propositions. So, viewing structures purely bicategorically is like studying a
logic without an accompanying type theory. Within this perspective, we have the
goal of axiomatizing the double category of profunctors as a ``double topos''
whose internal logic should be a linear logic over a directed type theory.
Compactness should be a key feature of this axiomatization as it is with
\emph{2-toposes} \cite{weber2007}. However, as noted above, our known structures
of bicategories of relations and double categories of relations are both compact
but have a \emph{trivial} duality involution on objects. The candidate duality
involution for profunctors is clearly taking the \emph{opposite category}. Thus,
we intend to axiomatize compactness for profunctors viewing \emph{opposite} as
our non-trivial duality involution and to do this in a genuinely double
categorical way, hence the need for parameterized loose adjunctions. These are
phrased in terms of our twisted Hom functor. Other potential double structures
with non-trivial duality arise in directed homotopy theory, Lawvere metric
spaces, vector spaces, Hilbert spaces, and even ordinary ring theory. Each of
these could provide further instances of the developed framework.

Our interest also stems in part from our intended applications and the
recurrence of compactness in these settings. Our recent work \cite{lambert2025}
on \emph{data queries as functorial semantics} shows that data operations can be
realized through interaction between separately tracked relational propositions
on the one hand and genuinely functional assignments on the other. That is, in a
finitely presented double category of relations viewed as the schema for a
relational database, a proarrow encodes a relational proposition instanced by
data tables. In this context, compactness, with its trivial duality involution,
says that the arrangement of the columns across the two sides of the proarrow
does not matter. Our view, however, is that compactness has utility in arranging
and executing queries, optimizing such queries, and in presenting knowledge
bases with preferred semantic biases. Compactness is thus in these contexts seen
to be non-trivial, and keeping both orthogonal directions separate in the
double-categorical viewpoint affords maximal flexibility in the choice of types,
terms and propositions in setting up our relational schemas.

\subsection{Applications to compositional systems theory} \label{intro:systems}

A theme in applied category theory has been the study of open systems and their
composition. Fong and Spivak observe that the hypergraph categories used in
compositional systems theory may be equivalently described as algebras over a
(colored) operad of finite sets and cospans \cite{fong2019}. A cospan
with a coproduct source of some higher arity is interpreted as an undirected
wiring diagram (in the sense of \cite{spivak2013}) which may couple a collection
of systems by equating their respective exposed variables
according to the left leg of the cospan, and choosing to re-expose those
variables coming from the right leg. Formally, we may suppose that our systems
of interest form a finitely cocomplete category and that we have a
finite colimit preserving functor including finite sets as ``discrete'' or 
unconstrained systems. We may then equip a system with interface realized as a 
finite set by interpreting the (unconstrained) variables in this interface into 
the (constrained) variables of the system via a certain morphism in the
category of systems, or equivalently a cospan with one trivial leg given by an 
arrow from the empty system. A wiring diagram may then act on a family of
systems simply by cospan composition.

This procedure essentially gives a hypergraph category structure to
\emph{structured cospans} \cite{baez2020}, exhibiting them as algebras over the
operad of undirected wiring diagrams, except that of course this action is only
associative and unital up to coherent isomorphism. Rather, we should think of
this action not as one by an operad or a symmetric monoidal category, but by a
(symmetric monoidal) double category. Indeed, the above action is nothing more
than the restriction of the twisted Hom functor of the double category of 
cospans over the systems category along the empty system constant on the left 
and the double functor between cospan double categories induced from the
inclusion of finite sets on the right.

The idea of seeing categories of systems as indexed by a \emph{double category}
of interfaces and interaction patterns is the centerpiece of the second author's
work on compositional systems theory. In the paper \emph{Double categories of
  open dynamical systems} \cite{myers2021} and the unpublished book
\cite{myers2023}, the compositionality of dynamical systems is organized as
``indexed double categories'', or normal profunctor-valued lax double functors. 
However, the examples given in that paper take the represented double 
categories to be (supposedly strict) double categories of spans \emph{in the
  tight direction}; clearly, the second author should have used twisted double
functors on (pseudo-)double categories of spans, correctly placed in the loose
direction. With the technology developed in this paper, we can well understand
these twisted double functors of dynamical systems.

Libkind and the second author later resolved the strictness problem by using
\emph{double barrels} to model loose bimodules \cite{libkind2025}; they describe
a procedure for constructing symmetric monoidal loose right modules of dynamical
systems over symmetric monoidal double categories whose loose morphisms give
coupling laws for those systems. The procedure given there is, generally, to
construct a symmetric monoidal double category whose loose morphisms
may be understood as partially-coupled systems, and then to restrict its loose
Hom bimodule appropriately. While \cite{libkind2025} uses
\emph{double barrels} as its model of loose Hom bimodule, the approach taken
there is nicely understood in terms of twisted double functors. We formally
compare the two approaches in \cref{sec:double-barrel}.

\subsection{Conventions}

Double categories and double functors are assumed to be pseudo unless otherwise
stated; nevertheless, we typically elide associators and unitors in pasting
diagrams. The \define{tight morphisms} or \define{arrows} of a double category
$\dbl{D}$ are the morphisms of $\dbl{D}_0$ and the \define{loose morphisms} or
\define{proarrows} of $\dbl{D}$ are the objects of $\dbl{D}_1$. A tight
composite $x \xto{f} y \xto{g} z$ is written in applicative order as
$gf = g \circ f$ or in diagrammatic order as $f \cdot g$. Loose composition
$\odot: \dbl{D}_1 \times_{\dbl{D}_0} \dbl{D}_1 \to \dbl{D}_1$ is always written in
diagrammatic order, and loose identities are denoted
$\id: \dbl{D}_0 \to \dbl{D}_1$.

The \define{tight dual} or \define{opposite} of a double category $\dbl{D}$,
reversing its arrows, is denoted $\dbl{D}^\op$ and is defined by
$(\dbl{D}^\op)_0 \coloneqq (\dbl{D}_0)^\op$ and
$(\dbl{D}^\op)_1 \coloneqq (\dbl{D}_1)^\op$. The \define{loose dual} of a double
category $\dbl{D}$, reversing its proarrows, is denoted $\dbl{D}^\co$ and
exchanges the source and target functors
$\dbl{D}_1 \rightrightarrows \dbl{D}_0$.

\section{Twisted double functors}

The \define{transpose} of a strict double category $\dbl{D}$ is another strict
double category $\dbl{D}^\top$ that exchanges the categories of horizontal and
vertical morphisms and transposes the cells. Due to the asymmetry in the
definition of a (pseudo) double category, being strict in one direction and weak
in the other, transposition is possible only under the assumption of strictness.
Yet there are still occasions when it is useful to consider ``double functors'' of
the form $\dbl{D}^\top \to \dbl{E}$ or $\dbl{D} \to \dbl{E}^\top$, even when neither
transpose makes formal sense. We call these ``twisted double functors'' and denote
them $\dbl{D} \twistto \dbl{E}$. We begin by defining twisted \emph{doubly lax}
functors, though we will be primarily interested in laxness with respect to
arrows in the domain.

\begin{definition}[Twisted double functor] \label{defn:twisted-double-functor}
  Let $\dbl{D}$ and $\dbl{E}$ be double categories. A \define{twisted doubly lax
    functor} $F: \dbl{D} \twistto \dbl{E}$ consists of:
  \begin{itemize}
    \item for each object $x$ in $\dbl{D}$, an object $Fx$ in $\dbl{E}$;
    \item for each arrow $f: x \to y$ in $\dbl{D}$, a proarrow $Ff: Fx \proto Fy$ in
      $\dbl{E}$;
    \item for each proarrow $m: x \proto y$ in $\dbl{D}$, an arrow $Fm: Fx \to Fy$ in
      $\dbl{E}$; and
    \item for each cell $\stdInlineCell{\alpha}$ in $\dbl{D}$, a cell
      $\inlineCell{Fx}{Fw}{Fy}{Fz}{Ff}{Fg}{Fm}{Fn}{F\alpha}$ in $\dbl{E}$.
  \end{itemize}
  In addition, there are \define{tight-to-loose} comparison cells: for each
  pair of composable arrows $x \xto{f} y \xto{g} z$ in $\dbl{D}$, a
  \define{composition comparison} $F_{f,g}$, and for each object $x$ in
  $\dbl{D}$, an \define{identity comparison} $F_x$, which are cells in $\dbl{E}$
  of the form
  \begin{equation*}
    % https://q.uiver.app/#q=WzAsNSxbMCwwLCJGeCJdLFsxLDAsIkZ5Il0sWzIsMCwiRnoiXSxbMCwxLCJGeCJdLFsyLDEsIkZ6Il0sWzMsNCwiRihmIFxcY2RvdCBnKSIsMix7InN0eWxlIjp7ImJvZHkiOnsibmFtZSI6ImJhcnJlZCJ9fX1dLFswLDMsIiIsMCx7ImxldmVsIjoyLCJzdHlsZSI6eyJoZWFkIjp7Im5hbWUiOiJub25lIn19fV0sWzIsNCwiIiwyLHsibGV2ZWwiOjIsInN0eWxlIjp7ImhlYWQiOnsibmFtZSI6Im5vbmUifX19XSxbMCwxLCJGZiIsMCx7InN0eWxlIjp7ImJvZHkiOnsibmFtZSI6ImJhcnJlZCJ9fX1dLFsxLDIsIkZnIiwwLHsic3R5bGUiOnsiYm9keSI6eyJuYW1lIjoiYmFycmVkIn19fV0sWzEsNSwiRl97ZixnfSIsMSx7ImxhYmVsX3Bvc2l0aW9uIjo0MCwic2hvcnRlbiI6eyJ0YXJnZXQiOjIwfSwic3R5bGUiOnsiYm9keSI6eyJuYW1lIjoibm9uZSJ9LCJoZWFkIjp7Im5hbWUiOiJub25lIn19fV1d
    \begin{tikzcd}
      Fx & Fy & Fz \\
      Fx && Fz
      \arrow["Ff", "\shortmid"{marking}, from=1-1, to=1-2]
      \arrow[Rightarrow, no head, from=1-1, to=2-1]
      \arrow["Fg", "\shortmid"{marking}, from=1-2, to=1-3]
      \arrow[Rightarrow, no head, from=1-3, to=2-3]
      \arrow[""{name=0, anchor=center, inner sep=0}, "{F(f \cdot g)}"', "\shortmid"{marking}, from=2-1, to=2-3]
      \arrow["{F_{f,g}}"{description, pos=0.4}, draw=none, from=1-2, to=0]
    \end{tikzcd}
    \qquad\text{and}\qquad
    % https://q.uiver.app/#q=WzAsNCxbMCwwLCJGeCJdLFsxLDAsIkZ4Il0sWzAsMSwiRngiXSxbMSwxLCJGeCJdLFswLDEsIlxcbWF0aHJte2lkfV97Rnh9IiwwLHsic3R5bGUiOnsiYm9keSI6eyJuYW1lIjoiYmFycmVkIn19fV0sWzIsMywiRigxX3gpIiwyLHsic3R5bGUiOnsiYm9keSI6eyJuYW1lIjoiYmFycmVkIn19fV0sWzAsMiwiIiwwLHsibGV2ZWwiOjIsInN0eWxlIjp7ImhlYWQiOnsibmFtZSI6Im5vbmUifX19XSxbMSwzLCIiLDAseyJsZXZlbCI6Miwic3R5bGUiOnsiaGVhZCI6eyJuYW1lIjoibm9uZSJ9fX1dLFs0LDUsIkZfeCIsMSx7InNob3J0ZW4iOnsic291cmNlIjoyMCwidGFyZ2V0IjoyMH0sInN0eWxlIjp7ImJvZHkiOnsibmFtZSI6Im5vbmUifSwiaGVhZCI6eyJuYW1lIjoibm9uZSJ9fX1dXQ==
    \begin{tikzcd}
      Fx & Fx \\
      Fx & Fx
      \arrow[""{name=0, anchor=center, inner sep=0}, "{\mathrm{id}_{Fx}}", "\shortmid"{marking}, from=1-1, to=1-2]
      \arrow[Rightarrow, no head, from=1-1, to=2-1]
      \arrow[Rightarrow, no head, from=1-2, to=2-2]
      \arrow[""{name=1, anchor=center, inner sep=0}, "{F(1_x)}"', "\shortmid"{marking}, from=2-1, to=2-2]
      \arrow["{F_x}"{description}, draw=none, from=0, to=1]
    \end{tikzcd}.
  \end{equation*}
  Finally, there are \define{loose-to-tight} comparison cells: for each pair of
  composable proarrows $x \xproto{m} y \xproto{n} z$ in $\dbl{D}$, a
  \define{composition comparison} $F^{m,n}$, and for each object $x$ in
  $\dbl{D}$, an \define{identity comparison} $F^x$, which are cells in $\dbl{E}$
  of the form:
  \begin{equation*}
    % https://q.uiver.app/#q=WzAsNSxbMCwwLCJGeCJdLFswLDEsIkZ5Il0sWzAsMiwiRnoiXSxbMSwwLCJGeCJdLFsxLDIsIkZ6Il0sWzAsMSwiRm0iLDJdLFsxLDIsIkZuIiwyXSxbMyw0LCJGKG0gXFxvZG90IG4pIl0sWzAsMywiIiwyLHsibGV2ZWwiOjIsInN0eWxlIjp7ImJvZHkiOnsibmFtZSI6ImJhcnJlZCJ9LCJoZWFkIjp7Im5hbWUiOiJub25lIn19fV0sWzIsNCwiIiwwLHsibGV2ZWwiOjIsInN0eWxlIjp7ImJvZHkiOnsibmFtZSI6ImJhcnJlZCJ9LCJoZWFkIjp7Im5hbWUiOiJub25lIn19fV0sWzEsNywiRl57bSxufSIsMSx7ImxhYmVsX3Bvc2l0aW9uIjo0MCwic2hvcnRlbiI6eyJ0YXJnZXQiOjIwfSwic3R5bGUiOnsiYm9keSI6eyJuYW1lIjoibm9uZSJ9LCJoZWFkIjp7Im5hbWUiOiJub25lIn19fV1d
    \begin{tikzcd}
      Fx & Fx \\
      Fy \\
      Fz & Fz
      \arrow["\shortmid"{marking}, equals, from=1-1, to=1-2]
      \arrow["Fm"', from=1-1, to=2-1]
      \arrow[""{name=0, anchor=center, inner sep=0}, "{F(m \odot n)}", from=1-2, to=3-2]
      \arrow["Fn"', from=2-1, to=3-1]
      \arrow["\shortmid"{marking}, equals, from=3-1, to=3-2]
      \arrow["{F^{m,n}}"{description, pos=0.4}, draw=none, from=2-1, to=0]
    \end{tikzcd}
    \qquad\text{and}\qquad
    % https://q.uiver.app/#q=WzAsNCxbMCwwLCJGeCJdLFswLDEsIkZ4Il0sWzEsMCwiRngiXSxbMSwxLCJGeCJdLFswLDEsIjFfe0Z4fSIsMl0sWzIsMywiRihcXG1hdGhybXtpZH1feCkiXSxbMCwyLCIiLDIseyJsZXZlbCI6Miwic3R5bGUiOnsiYm9keSI6eyJuYW1lIjoiYmFycmVkIn0sImhlYWQiOnsibmFtZSI6Im5vbmUifX19XSxbMSwzLCIiLDAseyJsZXZlbCI6Miwic3R5bGUiOnsiYm9keSI6eyJuYW1lIjoiYmFycmVkIn0sImhlYWQiOnsibmFtZSI6Im5vbmUifX19XSxbNCw1LCJGXngiLDEseyJzaG9ydGVuIjp7InNvdXJjZSI6MjAsInRhcmdldCI6MjB9LCJzdHlsZSI6eyJib2R5Ijp7Im5hbWUiOiJub25lIn0sImhlYWQiOnsibmFtZSI6Im5vbmUifX19XV0=
    \begin{tikzcd}
      Fx & Fx \\
      Fx & Fx
      \arrow["\shortmid"{marking}, equals, from=1-1, to=1-2]
      \arrow[""{name=0, anchor=center, inner sep=0}, "{1_{Fx}}"', from=1-1, to=2-1]
      \arrow[""{name=1, anchor=center, inner sep=0}, "{F(\mathrm{id}_x)}", from=1-2, to=2-2]
      \arrow["\shortmid"{marking}, equals, from=2-1, to=2-2]
      \arrow["{F^x}"{description}, draw=none, from=0, to=1]
    \end{tikzcd}.
  \end{equation*}
  This data must satisfy the following axioms.
  \begin{itemize}
    \item Naturality of tight-to-loose composition comparisons: for cells
      $\inlineCell{x}{x'}{y}{y'}{m}{n}{f}{f'}{\alpha}$ and
      $\inlineCell{y}{y'}{z}{z'}{n}{p}{g}{g'}{\beta}$ in $\dbl{D}$, we have
      \begin{equation*}
        % https://q.uiver.app/#q=WzAsOCxbMCwwLCJGeCJdLFsxLDAsIkZ5Il0sWzAsMSwiRngnIl0sWzEsMSwiRnknIl0sWzIsMCwiRnoiXSxbMiwxLCJGeiciXSxbMCwyLCJGeCciXSxbMiwyLCJGeiciXSxbMCwxLCJGZiIsMCx7InN0eWxlIjp7ImJvZHkiOnsibmFtZSI6ImJhcnJlZCJ9fX1dLFswLDIsIkZtIiwyXSxbMSwzLCJGbiIsMV0sWzIsMywiRmYnIiwyLHsic3R5bGUiOnsiYm9keSI6eyJuYW1lIjoiYmFycmVkIn19fV0sWzEsNCwiRmciLDAseyJzdHlsZSI6eyJib2R5Ijp7Im5hbWUiOiJiYXJyZWQifX19XSxbMyw1LCJGZyciLDIseyJzdHlsZSI6eyJib2R5Ijp7Im5hbWUiOiJiYXJyZWQifX19XSxbNCw1LCJGcCJdLFsyLDYsIiIsMSx7ImxldmVsIjoyLCJzdHlsZSI6eyJoZWFkIjp7Im5hbWUiOiJub25lIn19fV0sWzUsNywiIiwxLHsibGV2ZWwiOjIsInN0eWxlIjp7ImhlYWQiOnsibmFtZSI6Im5vbmUifX19XSxbNiw3LCJGKGYnIFxcY2RvdCBnJykiLDIseyJzdHlsZSI6eyJib2R5Ijp7Im5hbWUiOiJiYXJyZWQifX19XSxbOSwxMCwiRlxcYWxwaGEiLDEseyJzaG9ydGVuIjp7InNvdXJjZSI6MjAsInRhcmdldCI6MjB9LCJzdHlsZSI6eyJib2R5Ijp7Im5hbWUiOiJub25lIn0sImhlYWQiOnsibmFtZSI6Im5vbmUifX19XSxbMywxNywiRl97ZicsZyd9IiwxLHsibGFiZWxfcG9zaXRpb24iOjQwLCJzaG9ydGVuIjp7InRhcmdldCI6MjB9LCJzdHlsZSI6eyJib2R5Ijp7Im5hbWUiOiJub25lIn0sImhlYWQiOnsibmFtZSI6Im5vbmUifX19XSxbMTAsMTQsIkZcXGJldGEiLDEseyJzaG9ydGVuIjp7InNvdXJjZSI6MjAsInRhcmdldCI6MjB9LCJzdHlsZSI6eyJib2R5Ijp7Im5hbWUiOiJub25lIn0sImhlYWQiOnsibmFtZSI6Im5vbmUifX19XV0=
        \begin{tikzcd}
          Fx & Fy & Fz \\
          {Fx'} & {Fy'} & {Fz'} \\
          {Fx'} && {Fz'}
          \arrow["Ff", "\shortmid"{marking}, from=1-1, to=1-2]
          \arrow[""{name=0, anchor=center, inner sep=0}, "Fm"', from=1-1, to=2-1]
          \arrow["Fg", "\shortmid"{marking}, from=1-2, to=1-3]
          \arrow[""{name=1, anchor=center, inner sep=0}, "Fn"{description}, from=1-2, to=2-2]
          \arrow[""{name=2, anchor=center, inner sep=0}, "Fp", from=1-3, to=2-3]
          \arrow["{Ff'}"', "\shortmid"{marking}, from=2-1, to=2-2]
          \arrow[Rightarrow, no head, from=2-1, to=3-1]
          \arrow["{Fg'}"', "\shortmid"{marking}, from=2-2, to=2-3]
          \arrow[Rightarrow, no head, from=2-3, to=3-3]
          \arrow[""{name=3, anchor=center, inner sep=0}, "{F(f' \cdot g')}"', "\shortmid"{marking}, from=3-1, to=3-3]
          \arrow["{F\alpha}"{description}, draw=none, from=0, to=1]
          \arrow["{F\beta}"{description}, draw=none, from=1, to=2]
          \arrow["{F_{f',g'}}"{description, pos=0.4}, draw=none, from=2-2, to=3]
        \end{tikzcd}
        \quad=\quad
        % https://q.uiver.app/#q=WzAsNyxbMCwwLCJGeCJdLFsyLDAsIkZ6Il0sWzIsMSwiRnoiXSxbMSwwLCJGeSJdLFswLDEsIkZ4Il0sWzAsMiwiRngnIl0sWzIsMiwiRnonIl0sWzAsNCwiIiwxLHsibGV2ZWwiOjIsInN0eWxlIjp7ImhlYWQiOnsibmFtZSI6Im5vbmUifX19XSxbMSwyLCIiLDEseyJsZXZlbCI6Miwic3R5bGUiOnsiaGVhZCI6eyJuYW1lIjoibm9uZSJ9fX1dLFswLDMsIkZmIiwwLHsic3R5bGUiOnsiYm9keSI6eyJuYW1lIjoiYmFycmVkIn19fV0sWzMsMSwiRmciLDAseyJzdHlsZSI6eyJib2R5Ijp7Im5hbWUiOiJiYXJyZWQifX19XSxbNCwyLCJGKGYgXFxjZG90IGcpIiwyLHsic3R5bGUiOnsiYm9keSI6eyJuYW1lIjoiYmFycmVkIn19fV0sWzQsNSwiRm0iLDJdLFsyLDYsIkZwIl0sWzUsNiwiRihmJyBcXGNkb3QgZycpIiwyLHsic3R5bGUiOnsiYm9keSI6eyJuYW1lIjoiYmFycmVkIn19fV0sWzMsMTEsIkZfe2YsZ30iLDEseyJsYWJlbF9wb3NpdGlvbiI6NDAsInNob3J0ZW4iOnsidGFyZ2V0IjoyMH0sInN0eWxlIjp7ImJvZHkiOnsibmFtZSI6Im5vbmUifSwiaGVhZCI6eyJuYW1lIjoibm9uZSJ9fX1dLFsxMSwxNCwiRihcXGFscGhhIFxcY2RvdCBcXGJldGEpIiwxLHsibGFiZWxfcG9zaXRpb24iOjYwLCJzaG9ydGVuIjp7InNvdXJjZSI6MjAsInRhcmdldCI6MjB9LCJzdHlsZSI6eyJib2R5Ijp7Im5hbWUiOiJub25lIn0sImhlYWQiOnsibmFtZSI6Im5vbmUifX19XV0=
        \begin{tikzcd}
          Fx & Fy & Fz \\
          Fx && Fz \\
          {Fx'} && {Fz'}
          \arrow["Ff", "\shortmid"{marking}, from=1-1, to=1-2]
          \arrow[Rightarrow, no head, from=1-1, to=2-1]
          \arrow["Fg", "\shortmid"{marking}, from=1-2, to=1-3]
          \arrow[Rightarrow, no head, from=1-3, to=2-3]
          \arrow[""{name=0, anchor=center, inner sep=0}, "{F(f \cdot g)}"', "\shortmid"{marking}, from=2-1, to=2-3]
          \arrow["Fm"', from=2-1, to=3-1]
          \arrow["Fp", from=2-3, to=3-3]
          \arrow[""{name=1, anchor=center, inner sep=0}, "{F(f' \cdot g')}"', "\shortmid"{marking}, from=3-1, to=3-3]
          \arrow["{F_{f,g}}"{description, pos=0.4}, draw=none, from=1-2, to=0]
          \arrow["{F(\alpha \cdot \beta)}"{description, pos=0.6}, draw=none, from=0, to=1]
        \end{tikzcd}.
      \end{equation*}
    \item Naturality of tight-to-loose identity comparisons: for each proarrow
      $m: x \proto y$ in $\dbl{D}$,
      \begin{equation*}
        % https://q.uiver.app/#q=WzAsNixbMCwxLCJGeSJdLFsxLDEsIkZ5Il0sWzAsMiwiRnkiXSxbMSwyLCJGeSJdLFswLDAsIkZ4Il0sWzEsMCwiRngiXSxbMCwxLCJcXG1hdGhybXtpZH1fe0Z5fSIsMCx7InN0eWxlIjp7ImJvZHkiOnsibmFtZSI6ImJhcnJlZCJ9fX1dLFsyLDMsIkYgMV95IiwyLHsic3R5bGUiOnsiYm9keSI6eyJuYW1lIjoiYmFycmVkIn19fV0sWzAsMiwiIiwwLHsibGV2ZWwiOjIsInN0eWxlIjp7ImhlYWQiOnsibmFtZSI6Im5vbmUifX19XSxbMSwzLCIiLDAseyJsZXZlbCI6Miwic3R5bGUiOnsiaGVhZCI6eyJuYW1lIjoibm9uZSJ9fX1dLFs0LDAsIkZtIiwyXSxbNSwxLCJGbSJdLFs0LDUsIlxcbWF0aHJte2lkfV97Rnh9IiwwLHsic3R5bGUiOnsiYm9keSI6eyJuYW1lIjoiYmFycmVkIn19fV0sWzYsNywiRl95IiwxLHsic2hvcnRlbiI6eyJzb3VyY2UiOjIwLCJ0YXJnZXQiOjIwfSwic3R5bGUiOnsiYm9keSI6eyJuYW1lIjoibm9uZSJ9LCJoZWFkIjp7Im5hbWUiOiJub25lIn19fV0sWzEyLDYsIlxcbWF0aHJte2lkfV97Rm19IiwxLHsibGFiZWxfcG9zaXRpb24iOjQwLCJzaG9ydGVuIjp7InNvdXJjZSI6MjAsInRhcmdldCI6MjB9LCJzdHlsZSI6eyJib2R5Ijp7Im5hbWUiOiJub25lIn0sImhlYWQiOnsibmFtZSI6Im5vbmUifX19XV0=
        \begin{tikzcd}
          Fx & Fx \\
          Fy & Fy \\
          Fy & Fy
          \arrow[""{name=0, anchor=center, inner sep=0}, "{\mathrm{id}_{Fx}}", "\shortmid"{marking}, from=1-1, to=1-2]
          \arrow["Fm"', from=1-1, to=2-1]
          \arrow["Fm", from=1-2, to=2-2]
          \arrow[""{name=1, anchor=center, inner sep=0}, "{\mathrm{id}_{Fy}}", "\shortmid"{marking}, from=2-1, to=2-2]
          \arrow[Rightarrow, no head, from=2-1, to=3-1]
          \arrow[Rightarrow, no head, from=2-2, to=3-2]
          \arrow[""{name=2, anchor=center, inner sep=0}, "{F 1_y}"', "\shortmid"{marking}, from=3-1, to=3-2]
          \arrow["{\mathrm{id}_{Fm}}"{description, pos=0.4}, draw=none, from=0, to=1]
          \arrow["{F_y}"{description}, draw=none, from=1, to=2]
        \end{tikzcd}
        \quad=\quad
        % https://q.uiver.app/#q=WzAsNixbMCwwLCJGeCJdLFsxLDAsIkZ4Il0sWzAsMSwiRngiXSxbMSwxLCJGeCJdLFswLDIsIkZ5Il0sWzEsMiwiRnkiXSxbMCwxLCJcXG1hdGhybXtpZH1fe0Z4fSIsMCx7InN0eWxlIjp7ImJvZHkiOnsibmFtZSI6ImJhcnJlZCJ9fX1dLFsyLDMsIkYgMV94IiwwLHsic3R5bGUiOnsiYm9keSI6eyJuYW1lIjoiYmFycmVkIn19fV0sWzAsMiwiIiwwLHsibGV2ZWwiOjIsInN0eWxlIjp7ImhlYWQiOnsibmFtZSI6Im5vbmUifX19XSxbMSwzLCIiLDAseyJsZXZlbCI6Miwic3R5bGUiOnsiaGVhZCI6eyJuYW1lIjoibm9uZSJ9fX1dLFsyLDQsIkZtIiwyXSxbMyw1LCJGbSJdLFs0LDUsIkYgMV95IiwyLHsic3R5bGUiOnsiYm9keSI6eyJuYW1lIjoiYmFycmVkIn19fV0sWzYsNywiRl94IiwxLHsibGFiZWxfcG9zaXRpb24iOjQwLCJzaG9ydGVuIjp7InNvdXJjZSI6MjAsInRhcmdldCI6MjB9LCJzdHlsZSI6eyJib2R5Ijp7Im5hbWUiOiJub25lIn0sImhlYWQiOnsibmFtZSI6Im5vbmUifX19XSxbNywxMiwiRjFfbSIsMSx7InNob3J0ZW4iOnsic291cmNlIjoyMCwidGFyZ2V0IjoyMH0sInN0eWxlIjp7ImJvZHkiOnsibmFtZSI6Im5vbmUifSwiaGVhZCI6eyJuYW1lIjoibm9uZSJ9fX1dXQ==
        \begin{tikzcd}
          Fx & Fx \\
          Fx & Fx \\
          Fy & Fy
          \arrow[""{name=0, anchor=center, inner sep=0}, "{\mathrm{id}_{Fx}}", "\shortmid"{marking}, from=1-1, to=1-2]
          \arrow[Rightarrow, no head, from=1-1, to=2-1]
          \arrow[Rightarrow, no head, from=1-2, to=2-2]
          \arrow[""{name=1, anchor=center, inner sep=0}, "{F 1_x}", "\shortmid"{marking}, from=2-1, to=2-2]
          \arrow["Fm"', from=2-1, to=3-1]
          \arrow["Fm", from=2-2, to=3-2]
          \arrow[""{name=2, anchor=center, inner sep=0}, "{F 1_y}"', "\shortmid"{marking}, from=3-1, to=3-2]
          \arrow["{F_x}"{description, pos=0.4}, draw=none, from=0, to=1]
          \arrow["{F1_m}"{description}, draw=none, from=1, to=2]
        \end{tikzcd}.
      \end{equation*}
    \item Naturality of loose-to-tight composition comparisons: for cells
      $\inlineCell{x}{y}{x'}{y'}{m}{m'}{f}{g}{\alpha}$ and
      $\inlineCell{y}{z}{y'}{z'}{n}{n'}{g}{h}{\beta}$ in $\dbl{D}$, we have
      \begin{equation*}
        % https://q.uiver.app/#q=WzAsOCxbMCwwLCJGeCJdLFswLDEsIkZ5Il0sWzEsMCwiRngnIl0sWzAsMiwiRnoiXSxbMSwxLCJGeSciXSxbMSwyLCJGeiciXSxbMiwwLCJGeCciXSxbMiwyLCJGeiciXSxbMCwxLCJGbSIsMl0sWzAsMiwiRmYiLDAseyJzdHlsZSI6eyJib2R5Ijp7Im5hbWUiOiJiYXJyZWQifX19XSxbMSwzLCJGbiIsMl0sWzEsNCwiRmciLDAseyJzdHlsZSI6eyJib2R5Ijp7Im5hbWUiOiJiYXJyZWQifX19XSxbMyw1LCJGaCIsMix7InN0eWxlIjp7ImJvZHkiOnsibmFtZSI6ImJhcnJlZCJ9fX1dLFsyLDQsIkZtJyJdLFs0LDUsIkZuJyJdLFsyLDYsIiIsMSx7ImxldmVsIjoyLCJzdHlsZSI6eyJib2R5Ijp7Im5hbWUiOiJiYXJyZWQifSwiaGVhZCI6eyJuYW1lIjoibm9uZSJ9fX1dLFs1LDcsIiIsMSx7ImxldmVsIjoyLCJzdHlsZSI6eyJib2R5Ijp7Im5hbWUiOiJiYXJyZWQifSwiaGVhZCI6eyJuYW1lIjoibm9uZSJ9fX1dLFs2LDcsIkYobScgXFxvZG90IG4nKSJdLFs0LDE3LCJGXnttJyxuJ30iLDEseyJzaG9ydGVuIjp7InRhcmdldCI6MjB9LCJzdHlsZSI6eyJib2R5Ijp7Im5hbWUiOiJub25lIn0sImhlYWQiOnsibmFtZSI6Im5vbmUifX19XSxbMTEsMTIsIkZcXGJldGEiLDEseyJvZmZzZXQiOi0xLCJzaG9ydGVuIjp7InNvdXJjZSI6MjAsInRhcmdldCI6MjB9LCJzdHlsZSI6eyJib2R5Ijp7Im5hbWUiOiJub25lIn0sImhlYWQiOnsibmFtZSI6Im5vbmUifX19XSxbOSwxMSwiRlxcYWxwaGEiLDEseyJsYWJlbF9wb3NpdGlvbiI6NDAsInNob3J0ZW4iOnsic291cmNlIjoyMCwidGFyZ2V0IjoyMH0sInN0eWxlIjp7ImJvZHkiOnsibmFtZSI6Im5vbmUifSwiaGVhZCI6eyJuYW1lIjoibm9uZSJ9fX1dXQ==
        \begin{tikzcd}
          Fx & {Fx'} & {Fx'} \\
          Fy & {Fy'} \\
          Fz & {Fz'} & {Fz'}
          \arrow[""{name=0, anchor=center, inner sep=0}, "Ff"{inner sep=.8ex}, "\shortmid"{marking}, from=1-1, to=1-2]
          \arrow["Fm"', from=1-1, to=2-1]
          \arrow["\shortmid"{marking}, equals, from=1-2, to=1-3]
          \arrow["{Fm'}", from=1-2, to=2-2]
          \arrow[""{name=1, anchor=center, inner sep=0}, "{F(m' \odot n')}", from=1-3, to=3-3]
          \arrow[""{name=2, anchor=center, inner sep=0}, "Fg"{inner sep=.8ex}, "\shortmid"{marking}, from=2-1, to=2-2]
          \arrow["Fn"', from=2-1, to=3-1]
          \arrow["{Fn'}", from=2-2, to=3-2]
          \arrow[""{name=3, anchor=center, inner sep=0}, "Fh"'{inner sep=.8ex}, "\shortmid"{marking}, from=3-1, to=3-2]
          \arrow["\shortmid"{marking}, equals, from=3-2, to=3-3]
          \arrow["{F\alpha}"{description, pos=0.4}, draw=none, from=0, to=2]
          \arrow["{F\beta}"{description}, shift left, draw=none, from=2, to=3]
          \arrow["{F^{m',n'}}"{description}, draw=none, from=2-2, to=1]
        \end{tikzcd}
        \quad=\quad
        % https://q.uiver.app/#q=WzAsNyxbMSwwLCJGeCJdLFswLDAsIkZ4Il0sWzAsMSwiRnkiXSxbMCwyLCJGeiJdLFsxLDIsIkZ6Il0sWzIsMCwiRngnIl0sWzIsMiwiRnonIl0sWzEsMiwiRm0iLDJdLFsyLDMsIkZuIiwyXSxbMSwwLCIiLDEseyJsZXZlbCI6Miwic3R5bGUiOnsiYm9keSI6eyJuYW1lIjoiYmFycmVkIn0sImhlYWQiOnsibmFtZSI6Im5vbmUifX19XSxbMyw0LCIiLDAseyJsZXZlbCI6Miwic3R5bGUiOnsiYm9keSI6eyJuYW1lIjoiYmFycmVkIn0sImhlYWQiOnsibmFtZSI6Im5vbmUifX19XSxbNSw2LCJGKG0nIFxcb2RvdCBuJykiXSxbMCw1LCJGZiIsMCx7InN0eWxlIjp7ImJvZHkiOnsibmFtZSI6ImJhcnJlZCJ9fX1dLFs0LDYsIkZoIiwyLHsic3R5bGUiOnsiYm9keSI6eyJuYW1lIjoiYmFycmVkIn19fV0sWzAsNCwiRihtIFxcb2RvdCBuKSIsMV0sWzIsMTQsIkZee20sbn0iLDEseyJsYWJlbF9wb3NpdGlvbiI6NDAsInNob3J0ZW4iOnsidGFyZ2V0IjoyMH0sInN0eWxlIjp7ImJvZHkiOnsibmFtZSI6Im5vbmUifSwiaGVhZCI6eyJuYW1lIjoibm9uZSJ9fX1dLFsxNCwxMSwiRihcXGFscGhhIFxcb2RvdCBcXGJldGEpIiwxLHsibGFiZWxfcG9zaXRpb24iOjYwLCJzaG9ydGVuIjp7InNvdXJjZSI6MjAsInRhcmdldCI6MjB9LCJzdHlsZSI6eyJib2R5Ijp7Im5hbWUiOiJub25lIn0sImhlYWQiOnsibmFtZSI6Im5vbmUifX19XV0=
        \begin{tikzcd}[column sep=large]
          Fx & Fx & {Fx'} \\
          Fy \\
          Fz & Fz & {Fz'}
          \arrow["\shortmid"{marking}, equals, from=1-1, to=1-2]
          \arrow["Fm"', from=1-1, to=2-1]
          \arrow["Ff"{inner sep=.8ex}, "\shortmid"{marking}, from=1-2, to=1-3]
          \arrow[""{name=0, anchor=center, inner sep=0}, "{F(m \odot n)}"{description}, from=1-2, to=3-2]
          \arrow[""{name=1, anchor=center, inner sep=0}, "{F(m' \odot n')}", from=1-3, to=3-3]
          \arrow["Fn"', from=2-1, to=3-1]
          \arrow["\shortmid"{marking}, equals, from=3-1, to=3-2]
          \arrow["Fh"'{inner sep=.8ex}, "\shortmid"{marking}, from=3-2, to=3-3]
          \arrow["{F(\alpha \odot \beta)}"{description, pos=0.6}, draw=none, from=0, to=1]
          \arrow["{F^{m,n}}"{description, pos=0.4}, draw=none, from=2-1, to=0]
        \end{tikzcd}.
      \end{equation*}
    \item Naturality of loose-to-tight identity comparisons: for each arrow
      $f: x \to y$ in $\dbl{D}$,
      \begin{equation*}
        % https://q.uiver.app/#q=WzAsNixbMSwwLCJGeSJdLFsxLDEsIkZ5Il0sWzIsMCwiRnkiXSxbMiwxLCJGeSJdLFswLDAsIkZ4Il0sWzAsMSwiRngiXSxbMCwxLCIxX3tGeX0iLDFdLFsyLDMsIkYoXFxtYXRocm17aWR9X3kpIl0sWzAsMiwiIiwyLHsibGV2ZWwiOjIsInN0eWxlIjp7ImJvZHkiOnsibmFtZSI6ImJhcnJlZCJ9LCJoZWFkIjp7Im5hbWUiOiJub25lIn19fV0sWzEsMywiIiwwLHsibGV2ZWwiOjIsInN0eWxlIjp7ImJvZHkiOnsibmFtZSI6ImJhcnJlZCJ9LCJoZWFkIjp7Im5hbWUiOiJub25lIn19fV0sWzQsMCwiRmYiLDAseyJzdHlsZSI6eyJib2R5Ijp7Im5hbWUiOiJiYXJyZWQifX19XSxbNSwxLCJGZiIsMix7InN0eWxlIjp7ImJvZHkiOnsibmFtZSI6ImJhcnJlZCJ9fX1dLFs0LDUsIjFfe0Z4fSIsMl0sWzYsNywiRl55IiwxLHsic2hvcnRlbiI6eyJzb3VyY2UiOjIwLCJ0YXJnZXQiOjIwfSwic3R5bGUiOnsiYm9keSI6eyJuYW1lIjoibm9uZSJ9LCJoZWFkIjp7Im5hbWUiOiJub25lIn19fV0sWzEwLDExLCIxX3tGZn0iLDEseyJzaG9ydGVuIjp7InNvdXJjZSI6MjAsInRhcmdldCI6MjB9LCJzdHlsZSI6eyJib2R5Ijp7Im5hbWUiOiJub25lIn0sImhlYWQiOnsibmFtZSI6Im5vbmUifX19XV0=
        \begin{tikzcd}
          Fx & Fy & Fy \\
          Fx & Fy & Fy
          \arrow[""{name=0, anchor=center, inner sep=0}, "Ff"{inner sep=.8ex}, "\shortmid"{marking}, from=1-1, to=1-2]
          \arrow["{1_{Fx}}"', from=1-1, to=2-1]
          \arrow["\shortmid"{marking}, equals, from=1-2, to=1-3]
          \arrow[""{name=1, anchor=center, inner sep=0}, "{1_{Fy}}"{description}, from=1-2, to=2-2]
          \arrow[""{name=2, anchor=center, inner sep=0}, "{F(\mathrm{id}_y)}", from=1-3, to=2-3]
          \arrow[""{name=3, anchor=center, inner sep=0}, "Ff"'{inner sep=.8ex}, "\shortmid"{marking}, from=2-1, to=2-2]
          \arrow["\shortmid"{marking}, equals, from=2-2, to=2-3]
          \arrow["{1_{Ff}}"{description}, draw=none, from=0, to=3]
          \arrow["{F^y}"{description}, draw=none, from=1, to=2]
        \end{tikzcd}
        \quad=\quad
        % https://q.uiver.app/#q=WzAsNixbMCwwLCJGeCJdLFswLDEsIkZ4Il0sWzEsMCwiRngiXSxbMSwxLCJGeCJdLFsyLDAsIkZ5Il0sWzIsMSwiRnkiXSxbMCwxLCIxX3tGeH0iLDJdLFsyLDMsIkYgXFxtYXRocm17aWR9X3giLDFdLFswLDIsIiIsMix7ImxldmVsIjoyLCJzdHlsZSI6eyJib2R5Ijp7Im5hbWUiOiJiYXJyZWQifSwiaGVhZCI6eyJuYW1lIjoibm9uZSJ9fX1dLFsxLDMsIiIsMCx7ImxldmVsIjoyLCJzdHlsZSI6eyJib2R5Ijp7Im5hbWUiOiJiYXJyZWQifSwiaGVhZCI6eyJuYW1lIjoibm9uZSJ9fX1dLFsyLDQsIkZmIiwwLHsic3R5bGUiOnsiYm9keSI6eyJuYW1lIjoiYmFycmVkIn19fV0sWzMsNSwiRmYiLDIseyJzdHlsZSI6eyJib2R5Ijp7Im5hbWUiOiJiYXJyZWQifX19XSxbNCw1LCJGKFxcbWF0aHJte2lkfV95KSJdLFs2LDcsIkZeeCIsMSx7InNob3J0ZW4iOnsic291cmNlIjoyMCwidGFyZ2V0IjoyMH0sInN0eWxlIjp7ImJvZHkiOnsibmFtZSI6Im5vbmUifSwiaGVhZCI6eyJuYW1lIjoibm9uZSJ9fX1dLFsxMCwxMSwiRlxcbWF0aHJte2lkfV9mIiwxLHsic2hvcnRlbiI6eyJzb3VyY2UiOjIwLCJ0YXJnZXQiOjIwfSwic3R5bGUiOnsiYm9keSI6eyJuYW1lIjoibm9uZSJ9LCJoZWFkIjp7Im5hbWUiOiJub25lIn19fV1d
        \begin{tikzcd}
          Fx & Fx & Fy \\
          Fx & Fx & Fy
          \arrow["\shortmid"{marking}, equals, from=1-1, to=1-2]
          \arrow[""{name=0, anchor=center, inner sep=0}, "{1_{Fx}}"', from=1-1, to=2-1]
          \arrow[""{name=1, anchor=center, inner sep=0}, "Ff"{inner sep=.8ex}, "\shortmid"{marking}, from=1-2, to=1-3]
          \arrow[""{name=2, anchor=center, inner sep=0}, "{F \mathrm{id}_x}"{description}, from=1-2, to=2-2]
          \arrow["{F(\mathrm{id}_y)}", from=1-3, to=2-3]
          \arrow["\shortmid"{marking}, equals, from=2-1, to=2-2]
          \arrow[""{name=3, anchor=center, inner sep=0}, "Ff"'{inner sep=.8ex}, "\shortmid"{marking}, from=2-2, to=2-3]
          \arrow["{F^x}"{description}, draw=none, from=0, to=2]
          \arrow["{F\mathrm{id}_f}"{description}, draw=none, from=1, to=3]
        \end{tikzcd}.
      \end{equation*}
    \item Associativity of tight-to-loose comparisons: for each triple of
      composable arrows $w \xto{f} x \xto{g} y \xto{h} z$ in $\dbl{D}$,
      \begin{equation*}
        % https://q.uiver.app/#q=WzAsOSxbMCwwLCJGdyJdLFsxLDAsIkZ4Il0sWzIsMCwiRnkiXSxbMywwLCJGeiJdLFswLDEsIkZ3Il0sWzIsMSwiRnkiXSxbMywxLCJGeiJdLFswLDIsIkZ3Il0sWzMsMiwiRnoiXSxbMCwxLCJGZiIsMCx7InN0eWxlIjp7ImJvZHkiOnsibmFtZSI6ImJhcnJlZCJ9fX1dLFsxLDIsIkZnIiwwLHsic3R5bGUiOnsiYm9keSI6eyJuYW1lIjoiYmFycmVkIn19fV0sWzIsMywiRmgiLDAseyJzdHlsZSI6eyJib2R5Ijp7Im5hbWUiOiJiYXJyZWQifX19XSxbMCw0LCIiLDIseyJsZXZlbCI6Miwic3R5bGUiOnsiaGVhZCI6eyJuYW1lIjoibm9uZSJ9fX1dLFsyLDUsIiIsMix7ImxldmVsIjoyLCJzdHlsZSI6eyJoZWFkIjp7Im5hbWUiOiJub25lIn19fV0sWzQsNSwiRihmIFxcY2RvdCBnKSIsMix7InN0eWxlIjp7ImJvZHkiOnsibmFtZSI6ImJhcnJlZCJ9fX1dLFszLDYsIiIsMCx7ImxldmVsIjoyLCJzdHlsZSI6eyJoZWFkIjp7Im5hbWUiOiJub25lIn19fV0sWzUsNiwiRmgiLDIseyJzdHlsZSI6eyJib2R5Ijp7Im5hbWUiOiJiYXJyZWQifX19XSxbNiw4LCIiLDEseyJsZXZlbCI6Miwic3R5bGUiOnsiaGVhZCI6eyJuYW1lIjoibm9uZSJ9fX1dLFs0LDcsIiIsMSx7ImxldmVsIjoyLCJzdHlsZSI6eyJoZWFkIjp7Im5hbWUiOiJub25lIn19fV0sWzcsOCwiRihmIFxcY2RvdCBnIFxcY2RvdCBoKSIsMix7InN0eWxlIjp7ImJvZHkiOnsibmFtZSI6ImJhcnJlZCJ9fX1dLFsxLDE0LCJGX3tmLGd9IiwxLHsic2hvcnRlbiI6eyJ0YXJnZXQiOjIwfSwic3R5bGUiOnsiYm9keSI6eyJuYW1lIjoibm9uZSJ9LCJoZWFkIjp7Im5hbWUiOiJub25lIn19fV0sWzExLDE2LCIxX3tGaH0iLDEseyJzaG9ydGVuIjp7InNvdXJjZSI6MjAsInRhcmdldCI6MjB9LCJzdHlsZSI6eyJib2R5Ijp7Im5hbWUiOiJub25lIn0sImhlYWQiOnsibmFtZSI6Im5vbmUifX19XSxbMTgsMTcsIkZfe2YgXFxjZG90IGcsIGh9IiwxLHsib2Zmc2V0IjoxLCJzaG9ydGVuIjp7InNvdXJjZSI6MjAsInRhcmdldCI6MjB9LCJzdHlsZSI6eyJib2R5Ijp7Im5hbWUiOiJub25lIn0sImhlYWQiOnsibmFtZSI6Im5vbmUifX19XV0=
        \begin{tikzcd}
          Fw & Fx & Fy & Fz \\
          Fw && Fy & Fz \\
          Fw &&& Fz
          \arrow["Ff"{inner sep=.8ex}, "\shortmid"{marking}, from=1-1, to=1-2]
          \arrow[equals, from=1-1, to=2-1]
          \arrow["Fg"{inner sep=.8ex}, "\shortmid"{marking}, from=1-2, to=1-3]
          \arrow[""{name=0, anchor=center, inner sep=0}, "Fh"{inner sep=.8ex}, "\shortmid"{marking}, from=1-3, to=1-4]
          \arrow[equals, from=1-3, to=2-3]
          \arrow[equals, from=1-4, to=2-4]
          \arrow[""{name=1, anchor=center, inner sep=0}, "{F(f \cdot g)}"'{inner sep=.8ex}, "\shortmid"{marking}, from=2-1, to=2-3]
          \arrow[""{name=2, anchor=center, inner sep=0}, equals, from=2-1, to=3-1]
          \arrow[""{name=3, anchor=center, inner sep=0}, "Fh"'{inner sep=.8ex}, "\shortmid"{marking}, from=2-3, to=2-4]
          \arrow[""{name=4, anchor=center, inner sep=0}, equals, from=2-4, to=3-4]
          \arrow["{F(f \cdot g \cdot h)}"'{inner sep=.8ex}, "\shortmid"{marking}, from=3-1, to=3-4]
          \arrow["{F_{f,g}}"{description}, draw=none, from=1-2, to=1]
          \arrow["{1_{Fh}}"{description}, draw=none, from=0, to=3]
          \arrow["{F_{f \cdot g, h}}"{description}, shift right, draw=none, from=2, to=4]
        \end{tikzcd}
        \quad=\quad
        % https://q.uiver.app/#q=WzAsOSxbMCwwLCJGdyJdLFsxLDAsIkZ4Il0sWzIsMCwiRnkiXSxbMywwLCJGeiJdLFswLDIsIkZ3Il0sWzMsMiwiRnoiXSxbMywxLCJGeiJdLFswLDEsIkZ3Il0sWzEsMSwiRngiXSxbMCwxLCJGZiIsMCx7InN0eWxlIjp7ImJvZHkiOnsibmFtZSI6ImJhcnJlZCJ9fX1dLFsxLDIsIkZnIiwwLHsic3R5bGUiOnsiYm9keSI6eyJuYW1lIjoiYmFycmVkIn19fV0sWzIsMywiRmgiLDAseyJzdHlsZSI6eyJib2R5Ijp7Im5hbWUiOiJiYXJyZWQifX19XSxbNCw1LCJGKGYgXFxjZG90IGcgXFxjZG90IGgpIiwyLHsic3R5bGUiOnsiYm9keSI6eyJuYW1lIjoiYmFycmVkIn19fV0sWzMsNiwiIiwwLHsibGV2ZWwiOjIsInN0eWxlIjp7ImhlYWQiOnsibmFtZSI6Im5vbmUifX19XSxbNiw1LCIiLDEseyJsZXZlbCI6Miwic3R5bGUiOnsiaGVhZCI6eyJuYW1lIjoibm9uZSJ9fX1dLFswLDcsIiIsMix7ImxldmVsIjoyLCJzdHlsZSI6eyJoZWFkIjp7Im5hbWUiOiJub25lIn19fV0sWzcsNCwiIiwxLHsibGV2ZWwiOjIsInN0eWxlIjp7ImhlYWQiOnsibmFtZSI6Im5vbmUifX19XSxbMSw4LCIiLDAseyJsZXZlbCI6Miwic3R5bGUiOnsiaGVhZCI6eyJuYW1lIjoibm9uZSJ9fX1dLFs3LDgsIkZmIiwyLHsic3R5bGUiOnsiYm9keSI6eyJuYW1lIjoiYmFycmVkIn19fV0sWzgsNiwiRihnIFxcY2RvdCBoKSIsMix7InN0eWxlIjp7ImJvZHkiOnsibmFtZSI6ImJhcnJlZCJ9fX1dLFsxNiwxNCwiRl97ZiwgZyBcXGNkb3QgaH0iLDEseyJvZmZzZXQiOjEsInNob3J0ZW4iOnsic291cmNlIjoyMCwidGFyZ2V0IjoyMH0sInN0eWxlIjp7ImJvZHkiOnsibmFtZSI6Im5vbmUifSwiaGVhZCI6eyJuYW1lIjoibm9uZSJ9fX1dLFs5LDE4LCIxX3tGZn0iLDEseyJzaG9ydGVuIjp7InNvdXJjZSI6MjAsInRhcmdldCI6MjB9LCJzdHlsZSI6eyJib2R5Ijp7Im5hbWUiOiJub25lIn0sImhlYWQiOnsibmFtZSI6Im5vbmUifX19XSxbMiwxOSwiRl97ZyxofSIsMSx7InNob3J0ZW4iOnsidGFyZ2V0IjoyMH0sInN0eWxlIjp7ImJvZHkiOnsibmFtZSI6Im5vbmUifSwiaGVhZCI6eyJuYW1lIjoibm9uZSJ9fX1dXQ==
        \begin{tikzcd}
          Fw & Fx & Fy & Fz \\
          Fw & Fx && Fz \\
          Fw &&& Fz
          \arrow[""{name=0, anchor=center, inner sep=0}, "Ff"{inner sep=.8ex}, "\shortmid"{marking}, from=1-1, to=1-2]
          \arrow[equals, from=1-1, to=2-1]
          \arrow["Fg"{inner sep=.8ex}, "\shortmid"{marking}, from=1-2, to=1-3]
          \arrow[equals, from=1-2, to=2-2]
          \arrow["Fh"{inner sep=.8ex}, "\shortmid"{marking}, from=1-3, to=1-4]
          \arrow[equals, from=1-4, to=2-4]
          \arrow[""{name=1, anchor=center, inner sep=0}, "Ff"'{inner sep=.8ex}, "\shortmid"{marking}, from=2-1, to=2-2]
          \arrow[""{name=2, anchor=center, inner sep=0}, equals, from=2-1, to=3-1]
          \arrow[""{name=3, anchor=center, inner sep=0}, "{F(g \cdot h)}"'{inner sep=.8ex}, "\shortmid"{marking}, from=2-2, to=2-4]
          \arrow[""{name=4, anchor=center, inner sep=0}, equals, from=2-4, to=3-4]
          \arrow["{F(f \cdot g \cdot h)}"'{inner sep=.8ex}, "\shortmid"{marking}, from=3-1, to=3-4]
          \arrow["{1_{Ff}}"{description}, draw=none, from=0, to=1]
          \arrow["{F_{g,h}}"{description}, draw=none, from=1-3, to=3]
          \arrow["{F_{f, g \cdot h}}"{description}, shift right, draw=none, from=2, to=4]
        \end{tikzcd}.
      \end{equation*}
    \item Unitality of tight-to-loose comparisons: for each arrow $f: x \to y$
      in $\dbl{D}$,
      \begin{equation*}
        % https://q.uiver.app/#q=WzAsOCxbMCwwLCJGeCJdLFsxLDAsIkZ4Il0sWzIsMCwiRnkiXSxbMCwxLCJGeCJdLFsxLDEsIkZ4Il0sWzIsMSwiRnkiXSxbMCwyLCJGeCJdLFsyLDIsIkZ5Il0sWzEsMiwiRmYiLDAseyJzdHlsZSI6eyJib2R5Ijp7Im5hbWUiOiJiYXJyZWQifX19XSxbMCwxLCJcXGlkX3tGeH0iLDAseyJzdHlsZSI6eyJib2R5Ijp7Im5hbWUiOiJiYXJyZWQifX19XSxbMCwzLCIiLDIseyJsZXZlbCI6Miwic3R5bGUiOnsiaGVhZCI6eyJuYW1lIjoibm9uZSJ9fX1dLFsxLDQsIiIsMix7ImxldmVsIjoyLCJzdHlsZSI6eyJoZWFkIjp7Im5hbWUiOiJub25lIn19fV0sWzMsNCwiRigxX3gpIiwyLHsic3R5bGUiOnsiYm9keSI6eyJuYW1lIjoiYmFycmVkIn19fV0sWzIsNSwiIiwwLHsibGV2ZWwiOjIsInN0eWxlIjp7ImhlYWQiOnsibmFtZSI6Im5vbmUifX19XSxbNCw1LCJGZiIsMix7InN0eWxlIjp7ImJvZHkiOnsibmFtZSI6ImJhcnJlZCJ9fX1dLFs1LDcsIiIsMSx7ImxldmVsIjoyLCJzdHlsZSI6eyJoZWFkIjp7Im5hbWUiOiJub25lIn19fV0sWzMsNiwiIiwxLHsibGV2ZWwiOjIsInN0eWxlIjp7ImhlYWQiOnsibmFtZSI6Im5vbmUifX19XSxbNiw3LCJGZiIsMix7InN0eWxlIjp7ImJvZHkiOnsibmFtZSI6ImJhcnJlZCJ9fX1dLFs4LDE0LCIxX3tGZn0iLDEseyJzaG9ydGVuIjp7InNvdXJjZSI6MjAsInRhcmdldCI6MjB9LCJzdHlsZSI6eyJib2R5Ijp7Im5hbWUiOiJub25lIn0sImhlYWQiOnsibmFtZSI6Im5vbmUifX19XSxbOSwxMiwiRl94IiwxLHsic2hvcnRlbiI6eyJzb3VyY2UiOjIwLCJ0YXJnZXQiOjIwfSwic3R5bGUiOnsiYm9keSI6eyJuYW1lIjoibm9uZSJ9LCJoZWFkIjp7Im5hbWUiOiJub25lIn19fV0sWzQsMTcsIkZfe3gsZn0iLDEseyJzaG9ydGVuIjp7InRhcmdldCI6MjB9LCJzdHlsZSI6eyJib2R5Ijp7Im5hbWUiOiJub25lIn0sImhlYWQiOnsibmFtZSI6Im5vbmUifX19XV0=
        \begin{tikzcd}
          Fx & Fx & Fy \\
          Fx & Fx & Fy \\
          Fx && Fy
          \arrow[""{name=0, anchor=center, inner sep=0}, "{\id_{Fx}}"{inner sep=.8ex}, "\shortmid"{marking}, from=1-1, to=1-2]
          \arrow[equals, from=1-1, to=2-1]
          \arrow[""{name=1, anchor=center, inner sep=0}, "Ff"{inner sep=.8ex}, "\shortmid"{marking}, from=1-2, to=1-3]
          \arrow[equals, from=1-2, to=2-2]
          \arrow[equals, from=1-3, to=2-3]
          \arrow[""{name=2, anchor=center, inner sep=0}, "{F(1_x)}"'{inner sep=.8ex}, "\shortmid"{marking}, from=2-1, to=2-2]
          \arrow[equals, from=2-1, to=3-1]
          \arrow[""{name=3, anchor=center, inner sep=0}, "Ff"'{inner sep=.8ex}, "\shortmid"{marking}, from=2-2, to=2-3]
          \arrow[equals, from=2-3, to=3-3]
          \arrow[""{name=4, anchor=center, inner sep=0}, "Ff"'{inner sep=.8ex}, "\shortmid"{marking}, from=3-1, to=3-3]
          \arrow["{F_x}"{description}, draw=none, from=0, to=2]
          \arrow["{1_{Ff}}"{description}, draw=none, from=1, to=3]
          \arrow["{F_{x,f}}"{description}, draw=none, from=2-2, to=4]
        \end{tikzcd}
        \quad=\quad
        1_{Ff}
        \quad=\quad
        % https://q.uiver.app/#q=WzAsOCxbMCwwLCJGeCJdLFsyLDAsIkZ5Il0sWzAsMSwiRngiXSxbMiwxLCJGeSJdLFswLDIsIkZ4Il0sWzIsMiwiRnkiXSxbMSwwLCJGeSJdLFsxLDEsIkZ5Il0sWzAsMiwiIiwyLHsibGV2ZWwiOjIsInN0eWxlIjp7ImhlYWQiOnsibmFtZSI6Im5vbmUifX19XSxbMSwzLCIiLDAseyJsZXZlbCI6Miwic3R5bGUiOnsiaGVhZCI6eyJuYW1lIjoibm9uZSJ9fX1dLFszLDUsIiIsMSx7ImxldmVsIjoyLCJzdHlsZSI6eyJoZWFkIjp7Im5hbWUiOiJub25lIn19fV0sWzIsNCwiIiwxLHsibGV2ZWwiOjIsInN0eWxlIjp7ImhlYWQiOnsibmFtZSI6Im5vbmUifX19XSxbNCw1LCJGZiIsMix7InN0eWxlIjp7ImJvZHkiOnsibmFtZSI6ImJhcnJlZCJ9fX1dLFswLDYsIkZmIiwwLHsic3R5bGUiOnsiYm9keSI6eyJuYW1lIjoiYmFycmVkIn19fV0sWzYsMSwiXFxpZF97Rnl9IiwwLHsic3R5bGUiOnsiYm9keSI6eyJuYW1lIjoiYmFycmVkIn19fV0sWzYsNywiIiwwLHsibGV2ZWwiOjIsInN0eWxlIjp7ImhlYWQiOnsibmFtZSI6Im5vbmUifX19XSxbMiw3LCJGZiIsMix7InN0eWxlIjp7ImJvZHkiOnsibmFtZSI6ImJhcnJlZCJ9fX1dLFs3LDMsIkYoMV95KSIsMix7InN0eWxlIjp7ImJvZHkiOnsibmFtZSI6ImJhcnJlZCJ9fX1dLFs3LDEyLCJGX3tmLHl9IiwxLHsic2hvcnRlbiI6eyJ0YXJnZXQiOjIwfSwic3R5bGUiOnsiYm9keSI6eyJuYW1lIjoibm9uZSJ9LCJoZWFkIjp7Im5hbWUiOiJub25lIn19fV0sWzE0LDE3LCJGX3kiLDEseyJzaG9ydGVuIjp7InNvdXJjZSI6MjAsInRhcmdldCI6MjB9LCJzdHlsZSI6eyJib2R5Ijp7Im5hbWUiOiJub25lIn0sImhlYWQiOnsibmFtZSI6Im5vbmUifX19XSxbMTMsMTYsIjFfe0ZmfSIsMSx7InNob3J0ZW4iOnsic291cmNlIjoyMCwidGFyZ2V0IjoyMH0sInN0eWxlIjp7ImJvZHkiOnsibmFtZSI6Im5vbmUifSwiaGVhZCI6eyJuYW1lIjoibm9uZSJ9fX1dXQ==
        \begin{tikzcd}
          Fx & Fy & Fy \\
          Fx & Fy & Fy \\
          Fx && Fy
          \arrow[""{name=0, anchor=center, inner sep=0}, "Ff"{inner sep=.8ex}, "\shortmid"{marking}, from=1-1, to=1-2]
          \arrow[equals, from=1-1, to=2-1]
          \arrow[""{name=1, anchor=center, inner sep=0}, "{\id_{Fy}}"{inner sep=.8ex}, "\shortmid"{marking}, from=1-2, to=1-3]
          \arrow[equals, from=1-2, to=2-2]
          \arrow[equals, from=1-3, to=2-3]
          \arrow[""{name=2, anchor=center, inner sep=0}, "Ff"'{inner sep=.8ex}, "\shortmid"{marking}, from=2-1, to=2-2]
          \arrow[equals, from=2-1, to=3-1]
          \arrow[""{name=3, anchor=center, inner sep=0}, "{F(1_y)}"'{inner sep=.8ex}, "\shortmid"{marking}, from=2-2, to=2-3]
          \arrow[equals, from=2-3, to=3-3]
          \arrow[""{name=4, anchor=center, inner sep=0}, "Ff"'{inner sep=.8ex}, "\shortmid"{marking}, from=3-1, to=3-3]
          \arrow["{1_{Ff}}"{description}, draw=none, from=0, to=2]
          \arrow["{F_y}"{description}, draw=none, from=1, to=3]
          \arrow["{F_{f,y}}"{description}, draw=none, from=2-2, to=4]
        \end{tikzcd}.
      \end{equation*}
    \item Associativity of loose-to-tight comparisons: for each triple of composable
      proarrows $w \xproto{m} x \xproto{n} y \xproto{p} z$ in $\dbl{D}$, we have
      \begin{equation*}
        \begin{tikzcd}[column sep=scriptsize]
          Fw & Fw & Fw & Fw \\
          Fx && Fw & Fw \\
          Fy & Fy \\
          Fz & Fz & Fz & Fz
          \arrow[""{name=0, anchor=center, inner sep=0}, "\shortmid"{marking}, equals, from=1-1, to=1-2]
          \arrow["Fm"', from=1-1, to=2-1]
          \arrow[""{name=1, anchor=center, inner sep=0}, "\shortmid"{marking}, equals, from=1-2, to=1-3]
          \arrow["{F(m \odot n)}", from=1-2, to=3-2]
          \arrow[""{name=2, anchor=center, inner sep=0}, "\shortmid"{marking}, equals, from=1-3, to=1-4]
          \arrow[equals, from=1-3, to=2-3]
          \arrow[equals, from=1-4, to=2-4]
          \arrow["Fn"', from=2-1, to=3-1]
          \arrow[""{name=3, anchor=center, inner sep=0}, "{F(1_w)}"'{inner sep=.8ex}, "\shortmid"{marking}, from=2-3, to=2-4]
          \arrow[from=2-3, to=4-3]
          \arrow["{F(m \odot (n \odot p))}", no head, from=2-4, to=4-4]
          \arrow[""{name=4, anchor=center, inner sep=0}, "\shortmid"{marking}, equals, from=3-1, to=3-2]
          \arrow["Fp"', from=3-1, to=4-1]
          \arrow["Fp", from=3-2, to=4-2]
          \arrow[""{name=5, anchor=center, inner sep=0}, "\shortmid"{marking}, equals, from=4-1, to=4-2]
          \arrow[""{name=6, anchor=center, inner sep=0}, "\shortmid"{marking}, equals, from=4-2, to=4-3]
          \arrow[""{name=7, anchor=center, inner sep=0}, "{F(1_z)}"'{inner sep=.8ex}, "\shortmid"{marking}, from=4-3, to=4-4]
          \arrow["{F^{m,n}}"{description}, draw=none, from=0, to=4]
          \arrow["{F^{m \odot n, p}}"{description}, draw=none, from=1, to=6]
          \arrow["{F_w}"{description}, draw=none, from=2, to=3]
          \arrow["{F\alpha_{m,n,p}}"{description}, draw=none, from=3, to=7]
          \arrow["{\id_{Fp}}"{description}, draw=none, from=4, to=5]
        \end{tikzcd}
        =
        \begin{tikzcd}[column sep=scriptsize]
          Fw & Fw & Fw & Fw \\
          Fx & Fx \\
          Fy && Fz & Fz \\
          Fz & Fz & Fz & Fz
          \arrow[""{name=0, anchor=center, inner sep=0}, "\shortmid"{marking}, equals, from=1-1, to=1-2]
          \arrow["Fm"', from=1-1, to=2-1]
          \arrow[""{name=1, anchor=center, inner sep=0}, "\shortmid"{marking}, equals, from=1-2, to=1-3]
          \arrow["Fm", from=1-2, to=2-2]
          \arrow[""{name=2, anchor=center, inner sep=0}, "\shortmid"{marking}, equals, from=1-3, to=1-4]
          \arrow[no head, from=1-3, to=3-3]
          \arrow["{F(m \odot (n \odot p))}", from=1-4, to=3-4]
          \arrow[""{name=3, anchor=center, inner sep=0}, "\shortmid"{marking}, equals, from=2-1, to=2-2]
          \arrow["Fn"', from=2-1, to=3-1]
          \arrow["{F(n \odot p)}", from=2-2, to=4-2]
          \arrow["Fp"', from=3-1, to=4-1]
          \arrow[""{name=4, anchor=center, inner sep=0}, "\shortmid"{marking}, equals, from=3-3, to=3-4]
          \arrow[equals, from=3-3, to=4-3]
          \arrow[equals, from=3-4, to=4-4]
          \arrow[""{name=5, anchor=center, inner sep=0}, "\shortmid"{marking}, equals, from=4-1, to=4-2]
          \arrow[""{name=6, anchor=center, inner sep=0}, "\shortmid"{marking}, equals, from=4-2, to=4-3]
          \arrow[""{name=7, anchor=center, inner sep=0}, "{F(1_z)}"', from=4-3, to=4-4]
          \arrow["{\id_{Fm}}"{description}, draw=none, from=0, to=3]
          \arrow["{F^{m, n \odot p}}"{description}, draw=none, from=1, to=6]
          \arrow["\id"{description}, draw=none, from=2, to=4]
          \arrow["{F^{n,p}}"{description}, draw=none, from=3, to=5]
          \arrow["{F_z}"{description}, draw=none, from=4, to=7]
        \end{tikzcd},
      \end{equation*}
      where the isomorphism $\alpha_{m,n,p}: (m \odot n) \odot p \xto{\cong} m \odot (n \odot p)$ is the
      associator in $\dbl{D}$.
    \item Unitality of loose-to-tight comparisons: for each proarrow $m: x \proto y$
      in $\dbl{D}$, we have
      \begin{equation*}
        \begin{tikzcd}
          Fx & Fx & Fx & Fx \\
          Fx & Fx & Fx & Fx \\
          Fy & Fy & Fy & Fy
          \arrow[""{name=0, anchor=center, inner sep=0}, "\shortmid"{marking}, equals, from=1-1, to=1-2]
          \arrow["{1_{Fx}}"', from=1-1, to=2-1]
          \arrow[""{name=1, anchor=center, inner sep=0}, "\shortmid"{marking}, equals, from=1-2, to=1-3]
          \arrow["{F\id_x}", from=1-2, to=2-2]
          \arrow[""{name=2, anchor=center, inner sep=0}, "\shortmid"{marking}, equals, from=1-3, to=1-4]
          \arrow[equals, from=1-3, to=2-3]
          \arrow[equals, from=1-4, to=2-4]
          \arrow[""{name=3, anchor=center, inner sep=0}, "\shortmid"{marking}, equals, from=2-1, to=2-2]
          \arrow["Fm"', from=2-1, to=3-1]
          \arrow["Fm", from=2-2, to=3-2]
          \arrow[""{name=4, anchor=center, inner sep=0}, "{F(1_x)}"{inner sep=.8ex}, "\shortmid"{marking}, from=2-3, to=2-4]
          \arrow[from=2-3, to=3-3]
          \arrow["Fm", from=2-4, to=3-4]
          \arrow[""{name=5, anchor=center, inner sep=0}, "\shortmid"{marking}, equals, from=3-1, to=3-2]
          \arrow[""{name=6, anchor=center, inner sep=0}, "\shortmid"{marking}, equals, from=3-2, to=3-3]
          \arrow[""{name=7, anchor=center, inner sep=0}, "{F(1_y)}"'{inner sep=.8ex}, "\shortmid"{marking}, from=3-3, to=3-4]
          \arrow["{F^x}"{description}, draw=none, from=0, to=3]
          \arrow["{F^{x,m}}"{description}, draw=none, from=1, to=6]
          \arrow["{F_x}"{description, pos=0.4}, draw=none, from=2, to=4]
          \arrow["{\id_{Fm}}"{description}, draw=none, from=3, to=5]
          \arrow["{F(\lambda_m)}"{description}, draw=none, from=4, to=7]
        \end{tikzcd}
        \quad=\quad
        % https://q.uiver.app/#q=WzAsNixbMCwwLCJGeCJdLFsxLDAsIkZ4Il0sWzAsMSwiRnkiXSxbMSwxLCJGeSJdLFswLDIsIkZ5Il0sWzEsMiwiRnkiXSxbMCwyLCJGbSIsMl0sWzEsMywiRm0iXSxbMCwxLCIiLDAseyJsZXZlbCI6Miwic3R5bGUiOnsiYm9keSI6eyJuYW1lIjoiYmFycmVkIn0sImhlYWQiOnsibmFtZSI6Im5vbmUifX19XSxbMiwzLCIiLDIseyJsZXZlbCI6Miwic3R5bGUiOnsiYm9keSI6eyJuYW1lIjoiYmFycmVkIn0sImhlYWQiOnsibmFtZSI6Im5vbmUifX19XSxbMyw1LCIiLDIseyJsZXZlbCI6Miwic3R5bGUiOnsiaGVhZCI6eyJuYW1lIjoibm9uZSJ9fX1dLFsyLDQsIiIsMix7ImxldmVsIjoyLCJzdHlsZSI6eyJoZWFkIjp7Im5hbWUiOiJub25lIn19fV0sWzQsNSwiRigxX3kpIiwyLHsic3R5bGUiOnsiYm9keSI6eyJuYW1lIjoiYmFycmVkIn19fV0sWzgsOSwiXFxpZF97Rm19IiwxLHsic2hvcnRlbiI6eyJzb3VyY2UiOjIwLCJ0YXJnZXQiOjIwfSwic3R5bGUiOnsiYm9keSI6eyJuYW1lIjoibm9uZSJ9LCJoZWFkIjp7Im5hbWUiOiJub25lIn19fV0sWzksMTIsIkZfeSIsMSx7InNob3J0ZW4iOnsic291cmNlIjoyMCwidGFyZ2V0IjoyMH0sInN0eWxlIjp7ImJvZHkiOnsibmFtZSI6Im5vbmUifSwiaGVhZCI6eyJuYW1lIjoibm9uZSJ9fX1dXQ==
        \begin{tikzcd}
          Fx & Fx \\
          Fy & Fy \\
          Fy & Fy
          \arrow[""{name=0, anchor=center, inner sep=0}, "\shortmid"{marking}, equals, from=1-1, to=1-2]
          \arrow["Fm"', from=1-1, to=2-1]
          \arrow["Fm", from=1-2, to=2-2]
          \arrow[""{name=1, anchor=center, inner sep=0}, "\shortmid"{marking}, equals, from=2-1, to=2-2]
          \arrow[equals, from=2-1, to=3-1]
          \arrow[equals, from=2-2, to=3-2]
          \arrow[""{name=2, anchor=center, inner sep=0}, "{F(1_y)}"'{inner sep=.8ex}, "\shortmid"{marking}, from=3-1, to=3-2]
          \arrow["{\id_{Fm}}"{description}, draw=none, from=0, to=1]
          \arrow["{F_y}"{description}, draw=none, from=1, to=2]
        \end{tikzcd},
      \end{equation*}
      where the isomorphism $\lambda_m: \id_x \odot m \xto{\cong} m$ is the left unitor in
      $\dbl{D}$, and
      \begin{equation*}
        \begin{tikzcd}
          Fx & Fx & Fx & Fx \\
          Fy & Fy & Fx & Fx \\
          Fy & Fy & Fy & Fy
          \arrow[""{name=0, anchor=center, inner sep=0}, "\shortmid"{marking}, equals, from=1-1, to=1-2]
          \arrow["Fm"', from=1-1, to=2-1]
          \arrow[""{name=1, anchor=center, inner sep=0}, "\shortmid"{marking}, equals, from=1-2, to=1-3]
          \arrow["Fm", from=1-2, to=2-2]
          \arrow[""{name=2, anchor=center, inner sep=0}, "\shortmid"{marking}, equals, from=1-3, to=1-4]
          \arrow[equals, from=1-3, to=2-3]
          \arrow[equals, from=1-4, to=2-4]
          \arrow[""{name=3, anchor=center, inner sep=0}, "\shortmid"{marking}, equals, from=2-1, to=2-2]
          \arrow["{1_{Fy}}"', from=2-1, to=3-1]
          \arrow["{F\id_y}", from=2-2, to=3-2]
          \arrow[""{name=4, anchor=center, inner sep=0}, "{F(1_x)}"{inner sep=.8ex}, "\shortmid"{marking}, from=2-3, to=2-4]
          \arrow[from=2-3, to=3-3]
          \arrow["Fm", from=2-4, to=3-4]
          \arrow[""{name=5, anchor=center, inner sep=0}, "\shortmid"{marking}, equals, from=3-1, to=3-2]
          \arrow[""{name=6, anchor=center, inner sep=0}, "\shortmid"{marking}, equals, from=3-2, to=3-3]
          \arrow[""{name=7, anchor=center, inner sep=0}, "{F(1_y)}"'{inner sep=.8ex}, "\shortmid"{marking}, from=3-3, to=3-4]
          \arrow["{\id_{Fm}}"{description}, draw=none, from=0, to=3]
          \arrow["{F^{m,y}}"{description}, draw=none, from=1, to=6]
          \arrow["{F_x}"{description, pos=0.4}, draw=none, from=2, to=4]
          \arrow["{F^y}"{description}, draw=none, from=3, to=5]
          \arrow["{F(\rho_m)}"{description}, draw=none, from=4, to=7]
        \end{tikzcd}
        \quad=\quad
        % https://q.uiver.app/#q=WzAsNixbMCwwLCJGeCJdLFsxLDAsIkZ4Il0sWzAsMSwiRnkiXSxbMSwxLCJGeSJdLFswLDIsIkZ5Il0sWzEsMiwiRnkiXSxbMCwyLCJGbSIsMl0sWzEsMywiRm0iXSxbMCwxLCIiLDAseyJsZXZlbCI6Miwic3R5bGUiOnsiYm9keSI6eyJuYW1lIjoiYmFycmVkIn0sImhlYWQiOnsibmFtZSI6Im5vbmUifX19XSxbMiwzLCIiLDIseyJsZXZlbCI6Miwic3R5bGUiOnsiYm9keSI6eyJuYW1lIjoiYmFycmVkIn0sImhlYWQiOnsibmFtZSI6Im5vbmUifX19XSxbMyw1LCIiLDIseyJsZXZlbCI6Miwic3R5bGUiOnsiaGVhZCI6eyJuYW1lIjoibm9uZSJ9fX1dLFsyLDQsIiIsMix7ImxldmVsIjoyLCJzdHlsZSI6eyJoZWFkIjp7Im5hbWUiOiJub25lIn19fV0sWzQsNSwiRigxX3kpIiwyLHsic3R5bGUiOnsiYm9keSI6eyJuYW1lIjoiYmFycmVkIn19fV0sWzgsOSwiXFxpZF97Rm19IiwxLHsic2hvcnRlbiI6eyJzb3VyY2UiOjIwLCJ0YXJnZXQiOjIwfSwic3R5bGUiOnsiYm9keSI6eyJuYW1lIjoibm9uZSJ9LCJoZWFkIjp7Im5hbWUiOiJub25lIn19fV0sWzksMTIsIkZfeSIsMSx7InNob3J0ZW4iOnsic291cmNlIjoyMCwidGFyZ2V0IjoyMH0sInN0eWxlIjp7ImJvZHkiOnsibmFtZSI6Im5vbmUifSwiaGVhZCI6eyJuYW1lIjoibm9uZSJ9fX1dXQ==
        \begin{tikzcd}
          Fx & Fx \\
          Fy & Fy \\
          Fy & Fy
          \arrow[""{name=0, anchor=center, inner sep=0}, "\shortmid"{marking}, equals, from=1-1, to=1-2]
          \arrow["Fm"', from=1-1, to=2-1]
          \arrow["Fm", from=1-2, to=2-2]
          \arrow[""{name=1, anchor=center, inner sep=0}, "\shortmid"{marking}, equals, from=2-1, to=2-2]
          \arrow[equals, from=2-1, to=3-1]
          \arrow[equals, from=2-2, to=3-2]
          \arrow[""{name=2, anchor=center, inner sep=0}, "{F(1_y)}"'{inner sep=.8ex}, "\shortmid"{marking}, from=3-1, to=3-2]
          \arrow["{\id_{Fm}}"{description}, draw=none, from=0, to=1]
          \arrow["{F_y}"{description}, draw=none, from=1, to=2]
        \end{tikzcd},
      \end{equation*}
      where the isomorphism $\rho_m: m \odot \id_y \xto{\cong} m$ is the right unitor in
      $\dbl{D}$.
    \end{itemize}

  A \define{twisted lax functor} is a twisted doubly lax functor
  $F: \dbl{D} \twistto \dbl{E}$ such that the loose-to-tight comparisons
  $F^{m,n}$ and $F^x$ have inverses
  \begin{equation*}
    % https://q.uiver.app/#q=WzAsNSxbMSwwLCJGeCJdLFsxLDEsIkZ5Il0sWzEsMiwiRnoiXSxbMCwwLCJGeCJdLFswLDIsIkZ6Il0sWzAsMSwiRm0iXSxbMSwyLCJGbiJdLFszLDQsIkYobSBcXG9kb3QgbikiLDJdLFswLDMsIiIsMix7ImxldmVsIjoyLCJzdHlsZSI6eyJib2R5Ijp7Im5hbWUiOiJiYXJyZWQifSwiaGVhZCI6eyJuYW1lIjoibm9uZSJ9fX1dLFsyLDQsIiIsMCx7ImxldmVsIjoyLCJzdHlsZSI6eyJib2R5Ijp7Im5hbWUiOiJiYXJyZWQifSwiaGVhZCI6eyJuYW1lIjoibm9uZSJ9fX1dLFsxLDcsIihGXnttLG59KV57LTF9IiwxLHsibGFiZWxfcG9zaXRpb24iOjQwLCJzaG9ydGVuIjp7InRhcmdldCI6MjB9LCJzdHlsZSI6eyJib2R5Ijp7Im5hbWUiOiJub25lIn0sImhlYWQiOnsibmFtZSI6Im5vbmUifX19XV0=
    \begin{tikzcd}
      Fx & Fx \\
      & Fy \\
      Fz & Fz
      \arrow[""{name=0, anchor=center, inner sep=0}, "{F(m \odot n)}"', from=1-1, to=3-1]
      \arrow["\shortmid"{marking}, equals, from=1-2, to=1-1]
      \arrow["Fm", from=1-2, to=2-2]
      \arrow["Fn", from=2-2, to=3-2]
      \arrow["\shortmid"{marking}, equals, from=3-2, to=3-1]
      \arrow["{(F^{m,n})^{-1}}"{description, pos=0.4}, draw=none, from=2-2, to=0]
    \end{tikzcd}
    \qquad\text{and}\qquad
    % https://q.uiver.app/#q=WzAsNCxbMSwwLCJGeCJdLFsxLDEsIkZ4Il0sWzAsMCwiRngiXSxbMCwxLCJGeCJdLFswLDEsIjFfe0Z4fSJdLFsyLDMsIkYoXFxtYXRocm17aWR9X3gpIiwyXSxbMCwyLCIiLDIseyJsZXZlbCI6Miwic3R5bGUiOnsiYm9keSI6eyJuYW1lIjoiYmFycmVkIn0sImhlYWQiOnsibmFtZSI6Im5vbmUifX19XSxbMSwzLCIiLDAseyJsZXZlbCI6Miwic3R5bGUiOnsiYm9keSI6eyJuYW1lIjoiYmFycmVkIn0sImhlYWQiOnsibmFtZSI6Im5vbmUifX19XSxbNCw1LCIoRl54KV57LTF9IiwxLHsic2hvcnRlbiI6eyJzb3VyY2UiOjIwLCJ0YXJnZXQiOjIwfSwic3R5bGUiOnsiYm9keSI6eyJuYW1lIjoibm9uZSJ9LCJoZWFkIjp7Im5hbWUiOiJub25lIn19fV1d
    \begin{tikzcd}
      Fx & Fx \\
      Fx & Fx
      \arrow[""{name=0, anchor=center, inner sep=0}, "{F(\mathrm{id}_x)}"', from=1-1, to=2-1]
      \arrow["\shortmid"{marking}, equals, from=1-2, to=1-1]
      \arrow[""{name=1, anchor=center, inner sep=0}, "{1_{Fx}}", from=1-2, to=2-2]
      \arrow["\shortmid"{marking}, equals, from=2-2, to=2-1]
      \arrow["{(F^x)^{-1}}"{description}, draw=none, from=1, to=0]
    \end{tikzcd}
  \end{equation*}
  with respect to loose composition in $\dbl{E}$, for all proarrows
  $x \xproto{m} y \xproto{n} z$ and all objects $x$ in $\dbl{D}$.

  A twisted lax or doubly lax functor $F: \dbl{D} \twistto \dbl{E}$ is
  \define{normal} if the tight-to-loose identity comparisons $F_x$ have inverses
  \begin{equation*}
    % https://q.uiver.app/#q=WzAsNCxbMCwxLCJGeCJdLFsxLDEsIkZ4Il0sWzAsMCwiRngiXSxbMSwwLCJGeCJdLFswLDEsIlxcbWF0aHJte2lkfV97Rnh9IiwyLHsic3R5bGUiOnsiYm9keSI6eyJuYW1lIjoiYmFycmVkIn19fV0sWzIsMywiRigxX3gpIiwwLHsic3R5bGUiOnsiYm9keSI6eyJuYW1lIjoiYmFycmVkIn19fV0sWzAsMiwiIiwwLHsibGV2ZWwiOjIsInN0eWxlIjp7ImhlYWQiOnsibmFtZSI6Im5vbmUifX19XSxbMSwzLCIiLDAseyJsZXZlbCI6Miwic3R5bGUiOnsiaGVhZCI6eyJuYW1lIjoibm9uZSJ9fX1dLFs0LDUsIkZfeF57LTF9IiwxLHsic2hvcnRlbiI6eyJzb3VyY2UiOjIwLCJ0YXJnZXQiOjIwfSwic3R5bGUiOnsiYm9keSI6eyJuYW1lIjoibm9uZSJ9LCJoZWFkIjp7Im5hbWUiOiJub25lIn19fV1d
    \begin{tikzcd}
      Fx & Fx \\
      Fx & Fx
      \arrow[""{name=0, anchor=center, inner sep=0}, "{F(1_x)}", "\shortmid"{marking}, from=1-1, to=1-2]
      \arrow[Rightarrow, no head, from=2-1, to=1-1]
      \arrow[""{name=1, anchor=center, inner sep=0}, "{\mathrm{id}_{Fx}}"', "\shortmid"{marking}, from=2-1, to=2-2]
      \arrow[Rightarrow, no head, from=2-2, to=1-2]
      \arrow["{F_x^{-1}}"{description}, draw=none, from=1, to=0]
    \end{tikzcd}
  \end{equation*}
  with respect to tight composition in $\dbl{E}$, for all objects $x$ in
  $\dbl{D}$. It is \define{unitary} if all tight-to-loose identity comparisons
  $F_x$ are tight identities. Finally, a twisted lax functor $F$ is
  \define{pseudo} or \define{strong} if the tight-to-loose comparisons $F_{f,g}$
  and $F_x$ have inverses
  \begin{equation*}
    % https://q.uiver.app/#q=WzAsNSxbMCwxLCJGeCJdLFsxLDEsIkZ5Il0sWzIsMSwiRnoiXSxbMCwwLCJGeCJdLFsyLDAsIkZ6Il0sWzMsNCwiRihmIFxcY2RvdCBnKSIsMCx7InN0eWxlIjp7ImJvZHkiOnsibmFtZSI6ImJhcnJlZCJ9fX1dLFswLDMsIiIsMCx7ImxldmVsIjoyLCJzdHlsZSI6eyJoZWFkIjp7Im5hbWUiOiJub25lIn19fV0sWzIsNCwiIiwyLHsibGV2ZWwiOjIsInN0eWxlIjp7ImhlYWQiOnsibmFtZSI6Im5vbmUifX19XSxbMCwxLCJGZiIsMix7InN0eWxlIjp7ImJvZHkiOnsibmFtZSI6ImJhcnJlZCJ9fX1dLFsxLDIsIkZnIiwyLHsic3R5bGUiOnsiYm9keSI6eyJuYW1lIjoiYmFycmVkIn19fV0sWzEsNSwiRl97ZixnfV57LTF9IiwxLHsibGFiZWxfcG9zaXRpb24iOjQwLCJzaG9ydGVuIjp7InRhcmdldCI6MjB9LCJzdHlsZSI6eyJib2R5Ijp7Im5hbWUiOiJub25lIn0sImhlYWQiOnsibmFtZSI6Im5vbmUifX19XV0=
    \begin{tikzcd}
      Fx && Fz \\
      Fx & Fy & Fz
      \arrow[""{name=0, anchor=center, inner sep=0}, "{F(f \cdot g)}", "\shortmid"{marking}, from=1-1, to=1-3]
      \arrow[Rightarrow, no head, from=2-1, to=1-1]
      \arrow["Ff"', "\shortmid"{marking}, from=2-1, to=2-2]
      \arrow["Fg"', "\shortmid"{marking}, from=2-2, to=2-3]
      \arrow[Rightarrow, no head, from=2-3, to=1-3]
      \arrow["{F_{f,g}^{-1}}"{description, pos=0.4}, draw=none, from=2-2, to=0]
    \end{tikzcd}
    \qquad\text{and}\qquad
    % https://q.uiver.app/#q=WzAsNCxbMCwxLCJGeCJdLFsxLDEsIkZ4Il0sWzAsMCwiRngiXSxbMSwwLCJGeCJdLFswLDEsIlxcbWF0aHJte2lkfV97Rnh9IiwyLHsic3R5bGUiOnsiYm9keSI6eyJuYW1lIjoiYmFycmVkIn19fV0sWzIsMywiRigxX3gpIiwwLHsic3R5bGUiOnsiYm9keSI6eyJuYW1lIjoiYmFycmVkIn19fV0sWzAsMiwiIiwwLHsibGV2ZWwiOjIsInN0eWxlIjp7ImhlYWQiOnsibmFtZSI6Im5vbmUifX19XSxbMSwzLCIiLDAseyJsZXZlbCI6Miwic3R5bGUiOnsiaGVhZCI6eyJuYW1lIjoibm9uZSJ9fX1dLFs0LDUsIkZfeF57LTF9IiwxLHsic2hvcnRlbiI6eyJzb3VyY2UiOjIwLCJ0YXJnZXQiOjIwfSwic3R5bGUiOnsiYm9keSI6eyJuYW1lIjoibm9uZSJ9LCJoZWFkIjp7Im5hbWUiOiJub25lIn19fV1d
    \begin{tikzcd}
      Fx & Fx \\
      Fx & Fx
      \arrow[""{name=0, anchor=center, inner sep=0}, "{F(1_x)}", "\shortmid"{marking}, from=1-1, to=1-2]
      \arrow[Rightarrow, no head, from=2-1, to=1-1]
      \arrow[""{name=1, anchor=center, inner sep=0}, "{\mathrm{id}_{Fx}}"', "\shortmid"{marking}, from=2-1, to=2-2]
      \arrow[Rightarrow, no head, from=2-2, to=1-2]
      \arrow["{F_x^{-1}}"{description}, draw=none, from=1, to=0]
    \end{tikzcd}
  \end{equation*}
  with respect to tight composition in $\dbl{E}$, for all arrows
  $x \xto{f} y \xto{g} z$ and all objects $x$ in $\dbl{D}$.
\end{definition}

Our terminology for twisted double functors reflects the philosophy of the
paper. We regard twisted lax functors, which are lax in one direction and strong
in the other, as the fundamental notion. They are the twisted counterpart to lax
double functors, which are lax in one direction and strict in the other. All of
our examples of twisted functors will be twisted normal lax. Nevertheless, when
it is simplest to do so, we will state our constructions in generality for
twisted doubly lax functors and specialize them later as needed.

\begin{remark}[Related concepts]
  A notion of ``double pseudofunctor'' between strict double categories has been
  defined by Shulman \cite[\S{6}]{shulman2011}. As the name suggests, such a map
  is pseudo in both directions. The definition is similar to the one above in
  having comparisons in both directions, but no twisting is involved. Another
  approach to twisted functors is to start with a notion of ``doubly weak double
  category'' \cite{verity1992,fairbanks2026}. Then transposition always makes
  sense, and a twisted functor could be defined as a doubly weak map from one
  doubly weak double category to the transpose of another. While such an
  approach is likely viable, from our point of view, much of the utility of
  (pseudo) double categories comes precisely from the asymmetry between tight
  and loose. Here, for example, natural transformations between twisted functors
  are defined to have components in the tight direction, yielding a vertical
  composition that strictly satisfies the category laws.
\end{remark}

In keeping with mathematical tradition, we first give a few rather trivial
examples of twisted double functors, which nevertheless serve to motivate the
development. Subsequent sections will furnish examples that utilize all of the data
present in a twisted double functor.

\begin{example}[Twisted point] \label{ex:twisted-point}
  Consider a \emph{twisted point}, namely, a normal lax functor
  $P\colon\dbl{1} \twistto \Prof$ out of the terminal double category. This
  amounts to a category $\cat{C} \coloneqq P(*)$ together with a functor
  $P(\id_*)\colon \cat{C} \to \cat{C}$ and a profunctor
  $P(1_*)\colon \cat{C} \proto \cat{C}$. Since $P$ is normal in both directions,
  $P(\id_*)$ is naturally isomorphic to the identity functor $1_{\cat{C}}$ and
  $P(1_*)$ is naturally isomorphic to the Yoneda bifunctor
  $\Hom_{\cat{C}}$. Moreover, the other comparison cells are all derivable from
  these isomorphisms. Thus, a twisted point $P\colon\dbl{1}\twistto\Prof$ is
  essentially just a category.
\end{example}

\begin{example}[Transposition]
  As expected, when $\dbl{D}$ is a strict double category and $\dbl{E}$ is any
  double category, a (normal) lax functor $F: \dbl{D}^\top \to \dbl{E}$ defines a
  twisted (normal) lax functor $F: \dbl{D} \twistto \dbl{E}$ making the same
  assignments and having the same tight-to-loose comparison cells. However, it
  should be noted that, even when $\dbl{D}$ is strict, an arbitrary twisted lax
  functor $\dbl{D} \twistto \dbl{E}$ is more general than a lax functor
  $\dbl{D}^\top \to \dbl{E}$ since the former is allowed to have invertible but
  nontrivial loose-to-tight comparison cells.
\end{example}

Profunctor-valued twisted functors generalize indexed categories, as the
following proposition and examples show. Given a category $\cat{C}$, let
$\Loose(\cat{C})$ be the double category derived from $\cat{C}$ by concentrating
its structure in the loose direction; that is, $\Loose(\cat{C})$ has as objects
and proarrows the objects and morphisms of $\cat{C}$ and has only trivial arrows
and cells. Given a double category $\dbl{D}$, let $\UTwoCat(\dbl{D})$ denote the
2-category underlying $\dbl{D}$, consisting of $\dbl{D}$'s objects, tight
morphisms, and special cells.

\begin{proposition} \label{prop:pseudo-functors-as-twisted-functors}
  Let $\cat{C}$ be a category and let $\dbl{D}$ be a double category. There is a
  bijection (respectively, an equivalence)
  \begin{equation*}
    \Loose(\cat{C}) \twistto \dbl{D}
    \qquad\leftrightsquigarrow\qquad
    \cat{C} \to \UTwoCat(\dbl{D})
  \end{equation*}
  between, on the left, twisted lax functors that are unitary (respectively,
  normal) and, on the right, pseudofunctors out of the discrete 2-category
  $\cat{C}$.
\end{proposition}
\begin{proof}
  Nearly immediate from the definitions. The loose-to-tight comparison cells on
  the left correspond to the comparison cells on the right, and the
  tight-to-loose comparison cells are trivial.
\end{proof}

\begin{example}[Indexed categories] \label{ex:indexed-cats-as-tw-funcs}
  Taking $\dbl{D} = \Prof$ in the preceding proposition, a twisted normal lax
  functor $\Loose(\cat{C}) \twistto \Prof$ corresponds to a pseudofunctor
  $F: \cat{C} \to \Cat$. The latter is an op-indexed category over $\cat{C}$, in
  which a transition functor $f_! \coloneqq Ff: Fa \to Fb$ is
  pseudofunctorially associated to each morphism $f: a \to b$ in $\cat{C}$.
\end{example}

\begin{example}[Presheaves] \label{ex:copresheaves-as-tw-funcs}
  Taking $\dbl{D} = \Span$ in the preceding proposition, a twisted normal lax
  functor $\Loose(\cat{C}) \twistto \Span$ corresponds to a functor
  $F: \cat{C} \to \Set$, noting that $\UTwoCat(\Span) = \Set$ is discrete as a
  2-category. The latter is an ordinary set-valued copresheaf on $\cat{C}$,
  having a transition function $f_! \coloneqq Ff: Fa \to Fb$ for each morphism
  $f: a \to b$ in $\cat{C}$.
\end{example}

\subsection{The twisted Hom double functor}

A key example of a lax double functor is the Hom double functor on a double
category $\dbl{D}$,
\begin{equation*}
  \Hom_{\dbl{D}}: \dbl{D}^\op \times \dbl{D} \to \Span,
\end{equation*}
introduced by Paré to formulate a Yoneda theory for double categories
\cite{pare2011}. On objects and arrows, the Hom double functor agrees with the
usual Hom functor
\begin{equation*}
  (\Hom_{\dbl{D}})_0 \coloneqq \Hom_{\dbl{D}_0}: \dbl{D}_0^\op \times \dbl{D}_0 \to \Set
\end{equation*}
on the underlying category objects and arrows. In view of the correspondence
between span-valued lax functors and profunctor-valued \emph{normal} lax
functors (see, for example, \cite[Corollary 2.17]{lambert2024}), the Hom double
functor is equivalently a normal lax functor
\begin{equation*}
  \Hom_{\dbl{D}}: \dbl{D}^\op \times \dbl{D} \to \Prof
\end{equation*}
that now sends a pair of objects $x$ and $y$ in $\dbl{D}$ to the category of
arrows $x \to y$ and loosely globular cells between them.

This Hom double functor extracts the sets (or categories) of \emph{ordinary} 
arrows in the represented double category. Since a double category has two 
kinds of morphisms, it stands to reason that there should be another kind of 
Hom double functor that extracts the categories of proarrows in a double 
category. But unless the double category is strict, such a functor must be 
twisted, and it indeed provides our most fundamental example of a twisted 
double functor.

\begin{construction}[Twisted Hom functor] \label{construction:twistedhomfunctor}
  The \define{twisted Hom functor} on a double category $\dbl{D}$ is the twisted
  unitary lax functor
  \begin{equation*}
    \dbl{D}(-,=): \dbl{D}^\co \times \dbl{D} \twistto \Prof,
  \end{equation*}
  defined by the assignments:
  \begin{itemize}
    \item for each pair of objects $x$ and $y$ in $\dbl{D}$, the
      \define{hom-category} $\dbl{D}(x,y)$ has as objects the proarrows
      $m: x \proto y$ in $\dbl{D}$ and as morphisms $m \to n$ the tightly globular 
      cells
      \begin{equation*}
        % https://q.uiver.app/#q=WzAsNCxbMCwwLCJ4Il0sWzEsMCwieSJdLFswLDEsIngiXSxbMSwxLCJ5Il0sWzAsMSwibSIsMCx7InN0eWxlIjp7ImJvZHkiOnsibmFtZSI6ImJhcnJlZCJ9fX1dLFsyLDMsIm4iLDIseyJzdHlsZSI6eyJib2R5Ijp7Im5hbWUiOiJiYXJyZWQifX19XSxbMCwyLCIiLDIseyJsZXZlbCI6Miwic3R5bGUiOnsiaGVhZCI6eyJuYW1lIjoibm9uZSJ9fX1dLFsxLDMsIiIsMCx7ImxldmVsIjoyLCJzdHlsZSI6eyJoZWFkIjp7Im5hbWUiOiJub25lIn19fV0sWzQsNSwiXFxhbHBoYSIsMSx7InNob3J0ZW4iOnsic291cmNlIjoyMCwidGFyZ2V0IjoyMH0sInN0eWxlIjp7ImJvZHkiOnsibmFtZSI6Im5vbmUifSwiaGVhZCI6eyJuYW1lIjoibm9uZSJ9fX1dXQ==
        \begin{tikzcd}[row sep=scriptsize]
          x & y \\
          x & y
          \arrow[""{name=0, anchor=center, inner sep=0}, "m", "\shortmid"{marking}, from=1-1, to=1-2]
          \arrow[Rightarrow, no head, from=1-1, to=2-1]
          \arrow[Rightarrow, no head, from=1-2, to=2-2]
          \arrow[""{name=1, anchor=center, inner sep=0}, "n"', "\shortmid"{marking}, from=2-1, to=2-2]
          \arrow["\alpha"{description}, draw=none, from=0, to=1]
        \end{tikzcd};
      \end{equation*}
    \item for each pair of arrows $f: x \to w$ and $g: y \to z$ in $\dbl{D}$,
      the profunctor between hom-categories
      \begin{equation*}
        \dbl{D}(f,g): \dbl{D}(x,y) \proto \dbl{D}(w,z)
      \end{equation*}
      has as heteromorphisms from one proarrow $m: x \proto y$ to another
      $n: w \proto z$, the cells in $\dbl{D}$ of form
      \begin{equation*}
        % https://q.uiver.app/#q=WzAsNCxbMCwwLCJ4Il0sWzEsMCwieSJdLFswLDEsInciXSxbMSwxLCJ6Il0sWzAsMSwibSIsMCx7InN0eWxlIjp7ImJvZHkiOnsibmFtZSI6ImJhcnJlZCJ9fX1dLFsyLDMsIm4iLDIseyJzdHlsZSI6eyJib2R5Ijp7Im5hbWUiOiJiYXJyZWQifX19XSxbMCwyLCJmIiwyXSxbMSwzLCJnIl0sWzQsNSwiXFxhbHBoYSIsMSx7InNob3J0ZW4iOnsic291cmNlIjoyMCwidGFyZ2V0IjoyMH0sInN0eWxlIjp7ImJvZHkiOnsibmFtZSI6Im5vbmUifSwiaGVhZCI6eyJuYW1lIjoibm9uZSJ9fX1dXQ==
        \begin{tikzcd}
          x & y \\
          w & z
          \arrow[""{name=0, anchor=center, inner sep=0}, "m", "\shortmid"{marking}, from=1-1, to=1-2]
          \arrow["f"', from=1-1, to=2-1]
          \arrow["g", from=1-2, to=2-2]
          \arrow[""{name=1, anchor=center, inner sep=0}, "n"', "\shortmid"{marking}, from=2-1, to=2-2]
          \arrow["\alpha"{description}, draw=none, from=0, to=1]
        \end{tikzcd},
      \end{equation*}
      on which tightly globular cells act by pre- and post-composition in 
      $\dbl{D}_1$;
    \item for each pair of proarrows $p: x' \proto x$ and $q: y \proto y'$ in
      $\dbl{D}$, the functor
      \begin{equation*}
        \dbl{D}(p,q): \dbl{D}(x,y) \to \dbl{D}(x',y')
      \end{equation*}
      sends a proarrow $m: x \proto y$ to the composite
      $p \odot m \odot q: x' \proto y'$ and sends a tightly globular cell
      \begin{equation*}
        % https://q.uiver.app/#q=WzAsNCxbMCwwLCJ4Il0sWzEsMCwieSJdLFswLDEsIngiXSxbMSwxLCJ5Il0sWzAsMSwibSIsMCx7InN0eWxlIjp7ImJvZHkiOnsibmFtZSI6ImJhcnJlZCJ9fX1dLFsyLDMsIm4iLDIseyJzdHlsZSI6eyJib2R5Ijp7Im5hbWUiOiJiYXJyZWQifX19XSxbMCwyLCIiLDIseyJsZXZlbCI6Miwic3R5bGUiOnsiaGVhZCI6eyJuYW1lIjoibm9uZSJ9fX1dLFsxLDMsIiIsMCx7ImxldmVsIjoyLCJzdHlsZSI6eyJoZWFkIjp7Im5hbWUiOiJub25lIn19fV0sWzQsNSwiXFxhbHBoYSIsMSx7InNob3J0ZW4iOnsic291cmNlIjoyMCwidGFyZ2V0IjoyMH0sInN0eWxlIjp7ImJvZHkiOnsibmFtZSI6Im5vbmUifSwiaGVhZCI6eyJuYW1lIjoibm9uZSJ9fX1dXQ==
        \begin{tikzcd}[row sep=scriptsize]
          x & y \\
          x & y
          \arrow[""{name=0, anchor=center, inner sep=0}, "m", "\shortmid"{marking}, from=1-1, to=1-2]
          \arrow[Rightarrow, no head, from=1-1, to=2-1]
          \arrow[Rightarrow, no head, from=1-2, to=2-2]
          \arrow[""{name=1, anchor=center, inner sep=0}, "n"', "\shortmid"{marking}, from=2-1, to=2-2]
          \arrow["\alpha"{description}, draw=none, from=0, to=1]
        \end{tikzcd}
        \quad\mapsto\quad
        % https://q.uiver.app/#q=WzAsOCxbMSwwLCJ4Il0sWzIsMCwieSJdLFsxLDEsIngiXSxbMiwxLCJ5Il0sWzAsMCwieCciXSxbMCwxLCJ4JyJdLFszLDEsInknIl0sWzMsMCwieSciXSxbMCwxLCJtIiwwLHsic3R5bGUiOnsiYm9keSI6eyJuYW1lIjoiYmFycmVkIn19fV0sWzIsMywibiIsMix7InN0eWxlIjp7ImJvZHkiOnsibmFtZSI6ImJhcnJlZCJ9fX1dLFswLDIsIiIsMix7ImxldmVsIjoyLCJzdHlsZSI6eyJoZWFkIjp7Im5hbWUiOiJub25lIn19fV0sWzEsMywiIiwwLHsibGV2ZWwiOjIsInN0eWxlIjp7ImhlYWQiOnsibmFtZSI6Im5vbmUifX19XSxbNCwwLCJwIiwwLHsic3R5bGUiOnsiYm9keSI6eyJuYW1lIjoiYmFycmVkIn19fV0sWzUsMiwicCIsMix7InN0eWxlIjp7ImJvZHkiOnsibmFtZSI6ImJhcnJlZCJ9fX1dLFs0LDUsIiIsMix7ImxldmVsIjoyLCJzdHlsZSI6eyJoZWFkIjp7Im5hbWUiOiJub25lIn19fV0sWzEsNywicSIsMCx7InN0eWxlIjp7ImJvZHkiOnsibmFtZSI6ImJhcnJlZCJ9fX1dLFs3LDYsIiIsMCx7ImxldmVsIjoyLCJzdHlsZSI6eyJoZWFkIjp7Im5hbWUiOiJub25lIn19fV0sWzMsNiwicSIsMix7InN0eWxlIjp7ImJvZHkiOnsibmFtZSI6ImJhcnJlZCJ9fX1dLFs4LDksIlxcYWxwaGEiLDEseyJzaG9ydGVuIjp7InNvdXJjZSI6MjAsInRhcmdldCI6MjB9LCJzdHlsZSI6eyJib2R5Ijp7Im5hbWUiOiJub25lIn0sImhlYWQiOnsibmFtZSI6Im5vbmUifX19XSxbMTIsMTMsIjFfcCIsMSx7InNob3J0ZW4iOnsic291cmNlIjoyMCwidGFyZ2V0IjoyMH0sInN0eWxlIjp7ImJvZHkiOnsibmFtZSI6Im5vbmUifSwiaGVhZCI6eyJuYW1lIjoibm9uZSJ9fX1dLFsxNSwxNywiMV9xIiwxLHsic2hvcnRlbiI6eyJzb3VyY2UiOjIwLCJ0YXJnZXQiOjIwfSwic3R5bGUiOnsiYm9keSI6eyJuYW1lIjoibm9uZSJ9LCJoZWFkIjp7Im5hbWUiOiJub25lIn19fV1d
        \begin{tikzcd}[row sep=scriptsize]
          {x'} & x & y & {y'} \\
          {x'} & x & y & {y'}
          \arrow[""{name=0, anchor=center, inner sep=0}, "p", "\shortmid"{marking}, from=1-1, to=1-2]
          \arrow[Rightarrow, no head, from=1-1, to=2-1]
          \arrow[""{name=1, anchor=center, inner sep=0}, "m", "\shortmid"{marking}, from=1-2, to=1-3]
          \arrow[Rightarrow, no head, from=1-2, to=2-2]
          \arrow[""{name=2, anchor=center, inner sep=0}, "q", "\shortmid"{marking}, from=1-3, to=1-4]
          \arrow[Rightarrow, no head, from=1-3, to=2-3]
          \arrow[Rightarrow, no head, from=1-4, to=2-4]
          \arrow[""{name=3, anchor=center, inner sep=0}, "p"', "\shortmid"{marking}, from=2-1, to=2-2]
          \arrow[""{name=4, anchor=center, inner sep=0}, "n"', "\shortmid"{marking}, from=2-2, to=2-3]
          \arrow[""{name=5, anchor=center, inner sep=0}, "q"', "\shortmid"{marking}, from=2-3, to=2-4]
          \arrow["{1_p}"{description}, draw=none, from=0, to=3]
          \arrow["\alpha"{description}, draw=none, from=1, to=4]
          \arrow["{1_q}"{description}, draw=none, from=2, to=5]
        \end{tikzcd}
      \end{equation*}
      as on the left to the composite on the right;
    \item for each pair of cells $\inlineCell{x'}{x}{w'}{w}{p}{r}{f'}{f}{\phi}$
      and $\inlineCell{y}{y'}{z}{z'}{q}{s}{g}{g'}{\psi}$ in $\dbl{D}$, the
      natural transformation
      \begin{equation*}
        % https://q.uiver.app/#q=WzAsNCxbMCwwLCJcXGRibHtEfSh4LHkpIl0sWzEsMCwiXFxkYmx7RH0odyx6KSJdLFswLDEsIlxcZGJse0R9KHgnLHknKSJdLFsxLDEsIlxcZGJse0R9KHcnLHonKSJdLFswLDEsIlxcZGJse0R9KGYsZykiLDAseyJzdHlsZSI6eyJib2R5Ijp7Im5hbWUiOiJiYXJyZWQifX19XSxbMCwyLCJcXGRibHtEfShwLHEpIiwyXSxbMSwzLCJcXGRibHtEfShyLHMpIl0sWzIsMywiXFxkYmx7RH0oZicsZycpIiwyLHsic3R5bGUiOnsiYm9keSI6eyJuYW1lIjoiYmFycmVkIn19fV0sWzQsNywiXFxkYmx7RH0oXFxwaGksXFxwc2kpIiwxLHsic2hvcnRlbiI6eyJzb3VyY2UiOjIwLCJ0YXJnZXQiOjIwfSwic3R5bGUiOnsiYm9keSI6eyJuYW1lIjoibm9uZSJ9LCJoZWFkIjp7Im5hbWUiOiJub25lIn19fV1d
        \begin{tikzcd}
          {\dbl{D}(x,y)} & {\dbl{D}(w,z)} \\
          {\dbl{D}(x',y')} & {\dbl{D}(w',z')}
          \arrow[""{name=0, anchor=center, inner sep=0}, "{\dbl{D}(f,g)}", "\shortmid"{marking}, from=1-1, to=1-2]
          \arrow["{\dbl{D}(p,q)}"', from=1-1, to=2-1]
          \arrow["{\dbl{D}(r,s)}", from=1-2, to=2-2]
          \arrow[""{name=1, anchor=center, inner sep=0}, "{\dbl{D}(f',g')}"', "\shortmid"{marking}, from=2-1, to=2-2]
          \arrow["{\dbl{D}(\phi,\psi)}"{description}, draw=none, from=0, to=1]
        \end{tikzcd}
      \end{equation*}
      has as its component at proarrows $m: x \proto y$ and $n: w \proto z$, the
      function
      \begin{equation*}
        \dbl{D}(\phi,\psi)_{m,n}: \dbl{D}(f,g)(m,n) \to
          \dbl{D}(f',g')(p \odot m \odot q, r \odot n \odot s)
      \end{equation*}
      acting on cells by loose pre- and post-composition in $\dbl{D}$:
      \begin{equation*}
        % https://q.uiver.app/#q=WzAsNCxbMCwwLCJ4Il0sWzEsMCwieSJdLFswLDEsInciXSxbMSwxLCJ6Il0sWzAsMSwibSIsMCx7InN0eWxlIjp7ImJvZHkiOnsibmFtZSI6ImJhcnJlZCJ9fX1dLFsyLDMsIm4iLDIseyJzdHlsZSI6eyJib2R5Ijp7Im5hbWUiOiJiYXJyZWQifX19XSxbMCwyLCJmIiwyXSxbMSwzLCJnIl0sWzQsNSwiXFxhbHBoYSIsMSx7InNob3J0ZW4iOnsic291cmNlIjoyMCwidGFyZ2V0IjoyMH0sInN0eWxlIjp7ImJvZHkiOnsibmFtZSI6Im5vbmUifSwiaGVhZCI6eyJuYW1lIjoibm9uZSJ9fX1dXQ==
        \begin{tikzcd}
          x & y \\
          w & z
          \arrow[""{name=0, anchor=center, inner sep=0}, "m", "\shortmid"{marking}, from=1-1, to=1-2]
          \arrow["f"', from=1-1, to=2-1]
          \arrow["g", from=1-2, to=2-2]
          \arrow[""{name=1, anchor=center, inner sep=0}, "n"', "\shortmid"{marking}, from=2-1, to=2-2]
          \arrow["\alpha"{description}, draw=none, from=0, to=1]
        \end{tikzcd}
        \quad\mapsto\quad
        % https://q.uiver.app/#q=WzAsOCxbMSwwLCJ4Il0sWzIsMCwieSJdLFsxLDEsInciXSxbMiwxLCJ6Il0sWzAsMCwieCciXSxbMCwxLCJ3JyJdLFszLDAsInknIl0sWzMsMSwieiciXSxbMCwxLCJtIiwwLHsic3R5bGUiOnsiYm9keSI6eyJuYW1lIjoiYmFycmVkIn19fV0sWzIsMywibiIsMix7InN0eWxlIjp7ImJvZHkiOnsibmFtZSI6ImJhcnJlZCJ9fX1dLFswLDIsImYiLDJdLFsxLDMsImciXSxbNCwwLCJwIiwwLHsic3R5bGUiOnsiYm9keSI6eyJuYW1lIjoiYmFycmVkIn19fV0sWzUsMiwiciIsMix7InN0eWxlIjp7ImJvZHkiOnsibmFtZSI6ImJhcnJlZCJ9fX1dLFs0LDUsImYnIiwyXSxbMyw3LCJzIiwyLHsic3R5bGUiOnsiYm9keSI6eyJuYW1lIjoiYmFycmVkIn19fV0sWzYsNywiZyciXSxbMSw2LCJxIiwwLHsic3R5bGUiOnsiYm9keSI6eyJuYW1lIjoiYmFycmVkIn19fV0sWzgsOSwiXFxhbHBoYSIsMSx7InNob3J0ZW4iOnsic291cmNlIjoyMCwidGFyZ2V0IjoyMH0sInN0eWxlIjp7ImJvZHkiOnsibmFtZSI6Im5vbmUifSwiaGVhZCI6eyJuYW1lIjoibm9uZSJ9fX1dLFsxMiwxMywiXFxwaGkiLDEseyJzaG9ydGVuIjp7InNvdXJjZSI6MjAsInRhcmdldCI6MjB9LCJzdHlsZSI6eyJib2R5Ijp7Im5hbWUiOiJub25lIn0sImhlYWQiOnsibmFtZSI6Im5vbmUifX19XSxbMTcsMTUsIlxccHNpIiwxLHsic2hvcnRlbiI6eyJzb3VyY2UiOjIwLCJ0YXJnZXQiOjIwfSwic3R5bGUiOnsiYm9keSI6eyJuYW1lIjoibm9uZSJ9LCJoZWFkIjp7Im5hbWUiOiJub25lIn19fV1d
        \begin{tikzcd}
          {x'} & x & y & {y'} \\
          {w'} & w & z & {z'}
          \arrow[""{name=0, anchor=center, inner sep=0}, "p", "\shortmid"{marking}, from=1-1, to=1-2]
          \arrow["{f'}"', from=1-1, to=2-1]
          \arrow[""{name=1, anchor=center, inner sep=0}, "m", "\shortmid"{marking}, from=1-2, to=1-3]
          \arrow["f"', from=1-2, to=2-2]
          \arrow[""{name=2, anchor=center, inner sep=0}, "q", "\shortmid"{marking}, from=1-3, to=1-4]
          \arrow["g", from=1-3, to=2-3]
          \arrow["{g'}", from=1-4, to=2-4]
          \arrow[""{name=3, anchor=center, inner sep=0}, "r"', "\shortmid"{marking}, from=2-1, to=2-2]
          \arrow[""{name=4, anchor=center, inner sep=0}, "n"', "\shortmid"{marking}, from=2-2, to=2-3]
          \arrow[""{name=5, anchor=center, inner sep=0}, "s"', "\shortmid"{marking}, from=2-3, to=2-4]
          \arrow["\phi"{description}, draw=none, from=0, to=3]
          \arrow["\alpha"{description}, draw=none, from=1, to=4]
          \arrow["\psi"{description}, draw=none, from=2, to=5]
        \end{tikzcd}.
      \end{equation*}
  \end{itemize}
  The comparison cells of the twisted Hom functor are defined as follows.
  \begin{itemize}
    \item For pairs of composable arrows $x \xto{f} y \xto{g} z$ and
      $x' \xto{f'} y' \xto{g'} z'$ in $\dbl{D}$, the natural transformation
      \begin{equation*}
        % https://q.uiver.app/#q=WzAsNSxbMCwwLCJcXGRibHtEfSh4LHgnKSJdLFsxLDAsIlxcZGJse0R9KHkseScpIl0sWzIsMCwiXFxkYmx7RH0oeix6JykiXSxbMCwxLCJcXGRibHtEfSh4LHgnKSJdLFsyLDEsIlxcZGJse0R9KHoseicpIl0sWzAsMSwiXFxkYmx7RH0oZixmJykiLDAseyJzdHlsZSI6eyJib2R5Ijp7Im5hbWUiOiJiYXJyZWQifX19XSxbMSwyLCJcXGRibHtEfShnLGcnKSIsMCx7InN0eWxlIjp7ImJvZHkiOnsibmFtZSI6ImJhcnJlZCJ9fX1dLFswLDMsIiIsMix7ImxldmVsIjoyLCJzdHlsZSI6eyJoZWFkIjp7Im5hbWUiOiJub25lIn19fV0sWzIsNCwiIiwwLHsibGV2ZWwiOjIsInN0eWxlIjp7ImhlYWQiOnsibmFtZSI6Im5vbmUifX19XSxbMyw0LCJcXGRibHtEfShmIFxcY2RvdCBnLCBmJyBcXGNkb3QgZycpIiwyLHsic3R5bGUiOnsiYm9keSI6eyJuYW1lIjoiYmFycmVkIn19fV0sWzEsOSwiXFxkYmx7RH1feyhmLGYnKSwoZyxnJyl9IiwxLHsic2hvcnRlbiI6eyJ0YXJnZXQiOjIwfSwic3R5bGUiOnsiYm9keSI6eyJuYW1lIjoibm9uZSJ9LCJoZWFkIjp7Im5hbWUiOiJub25lIn19fV1d
        \begin{tikzcd}
          {\dbl{D}(x,x')} & {\dbl{D}(y,y')} & {\dbl{D}(z,z')} \\
          {\dbl{D}(x,x')} && {\dbl{D}(z,z')}
          \arrow["{\dbl{D}(f,f')}", "\shortmid"{marking}, from=1-1, to=1-2]
          \arrow[Rightarrow, no head, from=1-1, to=2-1]
          \arrow["{\dbl{D}(g,g')}", "\shortmid"{marking}, from=1-2, to=1-3]
          \arrow[Rightarrow, no head, from=1-3, to=2-3]
          \arrow[""{name=0, anchor=center, inner sep=0}, "{\dbl{D}(f \cdot g, f' \cdot g')}"', "\shortmid"{marking}, from=2-1, to=2-3]
          \arrow["{\dbl{D}_{(f,f'),(g,g')}}"{description}, draw=none, from=1-2, to=0]
        \end{tikzcd}
      \end{equation*}
      has as its components the functions that act by composing cells in
      $\dbl{D}_1$:
      \begin{equation*}
        \begin{gathered}
          % https://q.uiver.app/#q=WzAsNCxbMCwwLCJ4Il0sWzEsMCwieCciXSxbMCwxLCJ5Il0sWzEsMSwieSciXSxbMiwzLCJuIiwyLHsic3R5bGUiOnsiYm9keSI6eyJuYW1lIjoiYmFycmVkIn19fV0sWzAsMiwiZiIsMl0sWzEsMywiZiciXSxbMCwxLCJtIiwwLHsic3R5bGUiOnsiYm9keSI6eyJuYW1lIjoiYmFycmVkIn19fV0sWzcsNCwiXFxhbHBoYSIsMSx7InNob3J0ZW4iOnsic291cmNlIjoyMCwidGFyZ2V0IjoyMH0sInN0eWxlIjp7ImJvZHkiOnsibmFtZSI6Im5vbmUifSwiaGVhZCI6eyJuYW1lIjoibm9uZSJ9fX1dXQ==
          \begin{tikzcd}
            x & {x'} \\
            y & {y'}
            \arrow[""{name=0, anchor=center, inner sep=0}, "m", "\shortmid"{marking}, from=1-1, to=1-2]
            \arrow["f"', from=1-1, to=2-1]
            \arrow["{f'}", from=1-2, to=2-2]
            \arrow[""{name=1, anchor=center, inner sep=0}, "n"', "\shortmid"{marking}, from=2-1, to=2-2]
            \arrow["\alpha"{description}, draw=none, from=0, to=1]
          \end{tikzcd}
          \\
          % https://q.uiver.app/#q=WzAsNCxbMCwwLCJ5Il0sWzEsMCwieSciXSxbMCwxLCJ6Il0sWzEsMSwieiciXSxbMCwxLCJuIiwwLHsic3R5bGUiOnsiYm9keSI6eyJuYW1lIjoiYmFycmVkIn19fV0sWzAsMiwiZyIsMl0sWzEsMywiZyciXSxbMiwzLCJwIiwyLHsic3R5bGUiOnsiYm9keSI6eyJuYW1lIjoiYmFycmVkIn19fV0sWzQsNywiXFxiZXRhIiwxLHsic2hvcnRlbiI6eyJzb3VyY2UiOjIwLCJ0YXJnZXQiOjIwfSwic3R5bGUiOnsiYm9keSI6eyJuYW1lIjoibm9uZSJ9LCJoZWFkIjp7Im5hbWUiOiJub25lIn19fV1d
          \begin{tikzcd}
            y & {y'} \\
            z & {z'}
            \arrow[""{name=0, anchor=center, inner sep=0}, "n", "\shortmid"{marking}, from=1-1, to=1-2]
            \arrow["g"', from=1-1, to=2-1]
            \arrow["{g'}", from=1-2, to=2-2]
            \arrow[""{name=1, anchor=center, inner sep=0}, "p"', "\shortmid"{marking}, from=2-1, to=2-2]
            \arrow["\beta"{description}, draw=none, from=0, to=1]
          \end{tikzcd}
        \end{gathered}
        \quad\mapsto\quad
        % https://q.uiver.app/#q=WzAsNixbMCwwLCJ4Il0sWzAsMSwieSJdLFswLDIsInoiXSxbMSwwLCJ4JyJdLFsxLDEsInknIl0sWzEsMiwieiciXSxbMCwxLCJmIiwyXSxbMSwyLCJnIiwyXSxbMyw0LCJmJyJdLFs0LDUsImcnIl0sWzAsMywibSIsMCx7InN0eWxlIjp7ImJvZHkiOnsibmFtZSI6ImJhcnJlZCJ9fX1dLFsyLDUsInAiLDIseyJzdHlsZSI6eyJib2R5Ijp7Im5hbWUiOiJiYXJyZWQifX19XSxbMSw0LCJuIiwwLHsic3R5bGUiOnsiYm9keSI6eyJuYW1lIjoiYmFycmVkIn19fV0sWzEwLDEyLCJcXGFscGhhIiwxLHsic2hvcnRlbiI6eyJzb3VyY2UiOjIwLCJ0YXJnZXQiOjIwfSwic3R5bGUiOnsiYm9keSI6eyJuYW1lIjoibm9uZSJ9LCJoZWFkIjp7Im5hbWUiOiJub25lIn19fV0sWzEyLDExLCJcXGJldGEiLDEseyJzaG9ydGVuIjp7InNvdXJjZSI6MjAsInRhcmdldCI6MjB9LCJzdHlsZSI6eyJib2R5Ijp7Im5hbWUiOiJub25lIn0sImhlYWQiOnsibmFtZSI6Im5vbmUifX19XV0=
        \begin{tikzcd}
          x & {x'} \\
          y & {y'} \\
          z & {z'}
          \arrow[""{name=0, anchor=center, inner sep=0}, "m", "\shortmid"{marking}, from=1-1, to=1-2]
          \arrow["f"', from=1-1, to=2-1]
          \arrow["{f'}", from=1-2, to=2-2]
          \arrow[""{name=1, anchor=center, inner sep=0}, "n", "\shortmid"{marking}, from=2-1, to=2-2]
          \arrow["g"', from=2-1, to=3-1]
          \arrow["{g'}", from=2-2, to=3-2]
          \arrow[""{name=2, anchor=center, inner sep=0}, "p"', "\shortmid"{marking}, from=3-1, to=3-2]
          \arrow["\alpha"{description}, draw=none, from=0, to=1]
          \arrow["\beta"{description}, draw=none, from=1, to=2]
        \end{tikzcd}.
      \end{equation*}
    \item For a pair of objects $x$ and $x'$ in $\dbl{D}$, the natural
      isomorphism
      \begin{equation*}
        % https://q.uiver.app/#q=WzAsNCxbMCwwLCJcXGRibHtEfSh4LHgnKSJdLFsxLDAsIlxcZGJse0R9KHgseCcpIl0sWzAsMSwiXFxkYmx7RH0oeCx4JykiXSxbMSwxLCJcXGRibHtEfSh4LHgnKSJdLFswLDEsIlxcaWRfe1xcZGJse0R9KHgseCcpfSIsMCx7InN0eWxlIjp7ImJvZHkiOnsibmFtZSI6ImJhcnJlZCJ9fX1dLFswLDIsIiIsMix7ImxldmVsIjoyLCJzdHlsZSI6eyJoZWFkIjp7Im5hbWUiOiJub25lIn19fV0sWzEsMywiIiwwLHsibGV2ZWwiOjIsInN0eWxlIjp7ImhlYWQiOnsibmFtZSI6Im5vbmUifX19XSxbMiwzLCJcXGRibHtEfSgxX3gsMV97eCd9KSIsMix7InN0eWxlIjp7ImJvZHkiOnsibmFtZSI6ImJhcnJlZCJ9fX1dLFs0LDcsIlxcZGJse0R9X3soeCx4Jyl9IiwxLHsic2hvcnRlbiI6eyJzb3VyY2UiOjIwLCJ0YXJnZXQiOjIwfSwic3R5bGUiOnsiYm9keSI6eyJuYW1lIjoibm9uZSJ9LCJoZWFkIjp7Im5hbWUiOiJub25lIn19fV1d
        \begin{tikzcd}
          {\dbl{D}(x,x')} & {\dbl{D}(x,x')} \\
          {\dbl{D}(x,x')} & {\dbl{D}(x,x')}
          \arrow[""{name=0, anchor=center, inner sep=0}, "{\id_{\dbl{D}(x,x')}}", "\shortmid"{marking}, from=1-1, to=1-2]
          \arrow[Rightarrow, no head, from=1-1, to=2-1]
          \arrow[Rightarrow, no head, from=1-2, to=2-2]
          \arrow[""{name=1, anchor=center, inner sep=0}, "{\dbl{D}(1_x,1_{x'})}"', "\shortmid"{marking}, from=2-1, to=2-2]
          \arrow["{\dbl{D}_{(x,x')}}"{description}, draw=none, from=0, to=1]
        \end{tikzcd}
      \end{equation*}
      is simply the identity transformation at the Hom functor on the category
      $\dbl{D}(x,x')$.
    \item For pairs of composable proarrows $x'' \xproto{p'} x' \xproto{p} x$
      and $y \xproto{q} y' \xproto{q'} y''$ in $\dbl{D}$, the natural
      isomorphism
      \begin{equation*}
        % https://q.uiver.app/#q=WzAsNSxbMCwwLCJcXGRibHtEfSh4LHkpIl0sWzEsMCwiXFxkYmx7RH0oeCx5KSJdLFswLDEsIlxcZGJse0R9KHgnLHknKSJdLFswLDIsIlxcZGJse0R9KHgnJyx5JycpIl0sWzEsMiwiXFxkYmx7RH0oeCcnLHknJykiXSxbMCwxLCIiLDAseyJsZXZlbCI6Miwic3R5bGUiOnsiYm9keSI6eyJuYW1lIjoiYmFycmVkIn0sImhlYWQiOnsibmFtZSI6Im5vbmUifX19XSxbMCwyLCJcXGRibHtEfShwLHEpIiwyXSxbMiwzLCJcXGRibHtEfShwJyxxJykiLDJdLFsxLDQsIlxcZGJse0R9KHAnIFxcb2RvdCBwLHEgXFxvZG90IHEnKSJdLFszLDQsIiIsMix7ImxldmVsIjoyLCJzdHlsZSI6eyJib2R5Ijp7Im5hbWUiOiJiYXJyZWQifSwiaGVhZCI6eyJuYW1lIjoibm9uZSJ9fX1dLFs1LDksIlxcZGJse0R9XnsocCxxKSwocCcscScpfSIsMSx7Im9mZnNldCI6LTMsInNob3J0ZW4iOnsic291cmNlIjoyMCwidGFyZ2V0IjoyMH0sInN0eWxlIjp7ImJvZHkiOnsibmFtZSI6Im5vbmUifSwiaGVhZCI6eyJuYW1lIjoibm9uZSJ9fX1dXQ==
        \begin{tikzcd}
          {\dbl{D}(x,y)} & {\dbl{D}(x,y)} \\
          {\dbl{D}(x',y')} \\
          {\dbl{D}(x'',y'')} & {\dbl{D}(x'',y'')}
          \arrow[""{name=0, anchor=center, inner sep=0}, "\shortmid"{marking}, equals, from=1-1, to=1-2]
          \arrow["{\dbl{D}(p,q)}"', from=1-1, to=2-1]
          \arrow["{\dbl{D}(p' \odot p,q \odot q')}", from=1-2, to=3-2]
          \arrow["{\dbl{D}(p',q')}"', from=2-1, to=3-1]
          \arrow[""{name=1, anchor=center, inner sep=0}, "\shortmid"{marking}, equals, from=3-1, to=3-2]
          \arrow["{\dbl{D}^{(p,q),(p',q')}}"{description}, shift left=3, draw=none, from=0, to=1]
        \end{tikzcd}
      \end{equation*}
      has as its component at a proarrow $m: x \proto y$, the tightly globular
      isomorphism
      \begin{equation*}
        p' \odot (p \odot m \odot q) \odot q' \xto{\cong}
        (p' \odot p) \odot m \odot (q \odot q')
      \end{equation*}
      determined by the associators of $\dbl{D}$.
    \item For pairs of objects $x$ and $y$ in $\dbl{D}$, the natural isomorphism
      \begin{equation*}
        % https://q.uiver.app/#q=WzAsNCxbMCwwLCJcXGRibHtEfSh4LHkpIl0sWzEsMCwiXFxkYmx7RH0oeCx5KSJdLFswLDEsIlxcZGJse0R9KHgseSkiXSxbMSwxLCJcXGRibHtEfSh4LHkpIl0sWzAsMSwiIiwwLHsibGV2ZWwiOjIsInN0eWxlIjp7ImJvZHkiOnsibmFtZSI6ImJhcnJlZCJ9LCJoZWFkIjp7Im5hbWUiOiJub25lIn19fV0sWzAsMiwiMV97XFxkYmx7RH0oeCx5KX0iLDJdLFsyLDMsIiIsMix7ImxldmVsIjoyLCJzdHlsZSI6eyJib2R5Ijp7Im5hbWUiOiJiYXJyZWQifSwiaGVhZCI6eyJuYW1lIjoibm9uZSJ9fX1dLFsxLDMsIlxcZGJse0R9KFxcaWRfeCwgXFxpZF95KSJdLFs0LDYsIlxcZGJse0R9XnsoeCx5KX0iLDEseyJzaG9ydGVuIjp7InNvdXJjZSI6MjAsInRhcmdldCI6MjB9LCJzdHlsZSI6eyJib2R5Ijp7Im5hbWUiOiJub25lIn0sImhlYWQiOnsibmFtZSI6Im5vbmUifX19XV0=
        \begin{tikzcd}
          {\dbl{D}(x,y)} & {\dbl{D}(x,y)} \\
          {\dbl{D}(x,y)} & {\dbl{D}(x,y)}
          \arrow[""{name=0, anchor=center, inner sep=0}, "\shortmid"{marking}, equals, from=1-1, to=1-2]
          \arrow["{1_{\dbl{D}(x,y)}}"', from=1-1, to=2-1]
          \arrow["{\dbl{D}(\id_x, \id_y)}", from=1-2, to=2-2]
          \arrow[""{name=1, anchor=center, inner sep=0}, "\shortmid"{marking}, equals, from=2-1, to=2-2]
          \arrow["{\dbl{D}^{(x,y)}}"{description}, draw=none, from=0, to=1]
        \end{tikzcd}
      \end{equation*}
      has as its component at a proarrow $m: x \proto y$, the globular
      isomorphism
      \begin{equation*}
        m \xto{\cong} {\id_x} \odot {m} \odot {\id_y}
      \end{equation*}
      determined by the unitors of $\dbl{D}$. \qedhere
  \end{itemize}
\end{construction}

\begin{lemma}[Well-definedness of twisted Hom]
\label{lem:well-definedness}
  For any double category $\dbl{D}$, the twisted Hom $\dbl{D}(-, =)$ is well
  defined as a twisted lax functor and is unitary.
\end{lemma}
\begin{proof}
  Just as the data of the twisted Hom functor $\dbl{D}(-, =)$ repackages the
  data of the double category $\dbl{D}$, its axioms repackage the laws of a
  double category. That each hom-category $\dbl{D}(x,y)$ is indeed a category,
  and each $\dbl{D}(f,g)$ is a profunctor, uses the tight associativity and
  unitality of $\dbl{D}$; that each $\dbl{D}(p,q)$ is a functor and each
  $\dbl{D}(\phi,\psi)$ a natural transformation use special cases of the
  interchange law for $\dbl{D}$.

  Naturality of the tight-to-loose composition comparisons also follows from the
  interchange law: it says that given a composable grid
  \begin{equation*}
    \begin{dblArray}{ccc}
      \phi_1 & \alpha & \psi_1 \\
      \phi_2 & \beta & \psi_2
    \end{dblArray}
  \end{equation*}
  of cells in $\dbl{D}$, composing vertically then horizontally or horizontally
  then vertically yields the same result. Naturality of the tight-to-loose
  identity comparisons holds trivially because these comparisons are tight
  identities. Naturality of the loose-to-tight composition and identity
  comparisons follows from the naturality of the associators and unitors of
  $\dbl{D}$, respectively.

  Finally, associativity and unitality of the tight-to-loose comparisons uses
  tight associativity and unitality of $\dbl{D}$, whereas associativity and
  unitality of the loose-to-tight comparisons follows from the coherence laws
  for the associators and unitors of $\dbl{D}$.
\end{proof}

\subsection{Pre-composition with double functors}

A twisted double functor can be pre-composed with an ordinary double functor,
provided that the twisted functor is \emph{normal}. This construction will yield
a plethora of important examples, including the twisted representable functors.

\begin{construction}[Pre-composing onto a twisted functor]
  \label{construction:twisted-pre-composite}
  Let $F\colon \dbl{C}\to\dbl{D}$ be a lax double functor and let
  $G\colon \dbl{D} \twistto\dbl{E}$ be a twisted doubly lax functor that is normal.
  The \define{composite} of $F$ and $G$ is the twisted doubly lax functor
  \begin{equation*}
    GF \coloneqq G \circ F\colon \dbl{C} \twistto \dbl{E}
  \end{equation*}
  that maps objects, arrows, proarrows, and cells all in the evident way by
  applying $F$ and then $G$. Its comparison cells are defined as follows.
  \begin{itemize}
    \item For composable arrows $x \xto{f} y \xto{g} z$ in $\dbl{C}$, the
      composition comparison is
      \begin{equation*}
        % https://q.uiver.app/#q=WzAsNSxbMCwwLCJHRngiXSxbMSwwLCJHRnkiXSxbMiwwLCJHRnoiXSxbMiwxLCJHRnoiXSxbMCwxLCJHRngiXSxbMCwxLCJHRmYiLDAseyJzdHlsZSI6eyJib2R5Ijp7Im5hbWUiOiJiYXJyZWQifX19XSxbMSwyLCJHRmciLDAseyJzdHlsZSI6eyJib2R5Ijp7Im5hbWUiOiJiYXJyZWQifX19XSxbMiwzLCIiLDAseyJsZXZlbCI6Miwic3R5bGUiOnsiaGVhZCI6eyJuYW1lIjoibm9uZSJ9fX1dLFswLDQsIiIsMix7ImxldmVsIjoyLCJzdHlsZSI6eyJoZWFkIjp7Im5hbWUiOiJub25lIn19fV0sWzQsMywiR0YoZiBcXGNkb3QgZykiLDIseyJzdHlsZSI6eyJib2R5Ijp7Im5hbWUiOiJiYXJyZWQifX19XSxbMSw5LCIoR0YpX3tmLGd9IiwxLHsic2hvcnRlbiI6eyJ0YXJnZXQiOjIwfSwic3R5bGUiOnsiYm9keSI6eyJuYW1lIjoibm9uZSJ9LCJoZWFkIjp7Im5hbWUiOiJub25lIn19fV1d
        \begin{tikzcd}
          GFx & GFy & GFz \\
          GFx && GFz
          \arrow["GFf", "\shortmid"{marking}, from=1-1, to=1-2]
          \arrow[Rightarrow, no head, from=1-1, to=2-1]
          \arrow["GFg", "\shortmid"{marking}, from=1-2, to=1-3]
          \arrow[Rightarrow, no head, from=1-3, to=2-3]
          \arrow[""{name=0, anchor=center, inner sep=0}, "{GF(f \cdot g)}"', "\shortmid"{marking}, from=2-1, to=2-3]
          \arrow["{(GF)_{f,g}}"{description}, draw=none, from=1-2, to=0]
        \end{tikzcd}
        \quad\coloneqq\quad
        % https://q.uiver.app/#q=WzAsNSxbMCwwLCJHRngiXSxbMSwwLCJHRnkiXSxbMiwwLCJHRnoiXSxbMiwxLCJHRnoiXSxbMCwxLCJHRngiXSxbMCwxLCJHRmYiLDAseyJzdHlsZSI6eyJib2R5Ijp7Im5hbWUiOiJiYXJyZWQifX19XSxbMSwyLCJHRmciLDAseyJzdHlsZSI6eyJib2R5Ijp7Im5hbWUiOiJiYXJyZWQifX19XSxbMiwzLCIiLDAseyJsZXZlbCI6Miwic3R5bGUiOnsiaGVhZCI6eyJuYW1lIjoibm9uZSJ9fX1dLFswLDQsIiIsMix7ImxldmVsIjoyLCJzdHlsZSI6eyJoZWFkIjp7Im5hbWUiOiJub25lIn19fV0sWzQsMywiRyhGZiBcXGNkb3QgRmcpIiwyLHsic3R5bGUiOnsiYm9keSI6eyJuYW1lIjoiYmFycmVkIn19fV0sWzEsOSwiR197RmYsRmd9IiwxLHsic2hvcnRlbiI6eyJ0YXJnZXQiOjIwfSwic3R5bGUiOnsiYm9keSI6eyJuYW1lIjoibm9uZSJ9LCJoZWFkIjp7Im5hbWUiOiJub25lIn19fV1d
        \begin{tikzcd}
          GFx & GFy & GFz \\
          GFx && GFz
          \arrow["GFf", "\shortmid"{marking}, from=1-1, to=1-2]
          \arrow[Rightarrow, no head, from=1-1, to=2-1]
          \arrow["GFg", "\shortmid"{marking}, from=1-2, to=1-3]
          \arrow[Rightarrow, no head, from=1-3, to=2-3]
          \arrow[""{name=0, anchor=center, inner sep=0}, "{G(Ff \cdot Fg)}"', "\shortmid"{marking}, from=2-1, to=2-3]
          \arrow["{G_{Ff,Fg}}"{description}, draw=none, from=1-2, to=0]
        \end{tikzcd},
      \end{equation*}
      and for each object $x \in \dbl{C}$, the identity comparison is
      \begin{equation*}
        % https://q.uiver.app/#q=WzAsNCxbMCwwLCJHRngiXSxbMSwwLCJHRngiXSxbMSwxLCJHRngiXSxbMCwxLCJHRngiXSxbMCwxLCJcXGlkX3tHRnh9IiwwLHsic3R5bGUiOnsiYm9keSI6eyJuYW1lIjoiYmFycmVkIn19fV0sWzEsMiwiIiwyLHsibGV2ZWwiOjIsInN0eWxlIjp7ImhlYWQiOnsibmFtZSI6Im5vbmUifX19XSxbMCwzLCIiLDAseyJsZXZlbCI6Miwic3R5bGUiOnsiaGVhZCI6eyJuYW1lIjoibm9uZSJ9fX1dLFszLDIsIkdGMV94IiwyLHsic3R5bGUiOnsiYm9keSI6eyJuYW1lIjoiYmFycmVkIn19fV0sWzQsNywiR197Rnh9IiwxLHsic2hvcnRlbiI6eyJzb3VyY2UiOjIwLCJ0YXJnZXQiOjIwfSwic3R5bGUiOnsiYm9keSI6eyJuYW1lIjoibm9uZSJ9LCJoZWFkIjp7Im5hbWUiOiJub25lIn19fV1d
        \begin{tikzcd}
          GFx & GFx \\
          GFx & GFx
          \arrow[""{name=0, anchor=center, inner sep=0}, "{\id_{GFx}}", "\shortmid"{marking}, from=1-1, to=1-2]
          \arrow[Rightarrow, no head, from=1-1, to=2-1]
          \arrow[Rightarrow, no head, from=1-2, to=2-2]
          \arrow[""{name=1, anchor=center, inner sep=0}, "{GF1_x}"', "\shortmid"{marking}, from=2-1, to=2-2]
          \arrow["{(GF)_{x}}"{description}, draw=none, from=0, to=1]
        \end{tikzcd}
        \quad\coloneqq\quad
        % https://q.uiver.app/#q=WzAsNCxbMCwwLCJHRngiXSxbMSwwLCJHRngiXSxbMSwxLCJHRngiXSxbMCwxLCJHRngiXSxbMCwxLCJcXGlkX3tHRnh9IiwwLHsic3R5bGUiOnsiYm9keSI6eyJuYW1lIjoiYmFycmVkIn19fV0sWzEsMiwiIiwyLHsibGV2ZWwiOjIsInN0eWxlIjp7ImhlYWQiOnsibmFtZSI6Im5vbmUifX19XSxbMCwzLCIiLDAseyJsZXZlbCI6Miwic3R5bGUiOnsiaGVhZCI6eyJuYW1lIjoibm9uZSJ9fX1dLFszLDIsIkdGMV94IiwyLHsic3R5bGUiOnsiYm9keSI6eyJuYW1lIjoiYmFycmVkIn19fV0sWzQsNywiKEdGKV97eH0iLDEseyJzaG9ydGVuIjp7InNvdXJjZSI6MjAsInRhcmdldCI6MjB9LCJzdHlsZSI6eyJib2R5Ijp7Im5hbWUiOiJub25lIn0sImhlYWQiOnsibmFtZSI6Im5vbmUifX19XV0=
        \begin{tikzcd}
          GFx & GFx \\
          GFx & GFx
          \arrow[""{name=0, anchor=center, inner sep=0}, "{\id_{GFx}}", "\shortmid"{marking}, from=1-1, to=1-2]
          \arrow[Rightarrow, no head, from=1-1, to=2-1]
          \arrow[Rightarrow, no head, from=1-2, to=2-2]
          \arrow[""{name=1, anchor=center, inner sep=0}, "{GF1_x}"', "\shortmid"{marking}, from=2-1, to=2-2]
          \arrow["{G_{Fx}}"{description}, draw=none, from=0, to=1]
        \end{tikzcd}.
      \end{equation*}
    \item For composable proarrows $x \xproto{m} y \xproto{n} z$ in $\dbl{C}$, the
      composition comparison is
      \begin{equation*}
        % https://q.uiver.app/#q=WzAsNSxbMCwwLCJHRngiXSxbMSwwLCJHRngiXSxbMSwyLCJHRnoiXSxbMCwxLCJHRnkiXSxbMCwyLCJHRnoiXSxbMCwxLCIiLDAseyJsZXZlbCI6Miwic3R5bGUiOnsiYm9keSI6eyJuYW1lIjoiYmFycmVkIn0sImhlYWQiOnsibmFtZSI6Im5vbmUifX19XSxbMSwyLCJHRihtXFxvZG90IG4pIl0sWzAsMywiR0ZtIiwyXSxbMyw0LCJHRm4iLDJdLFs0LDIsIiIsMix7ImxldmVsIjoyLCJzdHlsZSI6eyJib2R5Ijp7Im5hbWUiOiJiYXJyZWQifSwiaGVhZCI6eyJuYW1lIjoibm9uZSJ9fX1dLFs1LDksIihHRilee20sbn0iLDEseyJvZmZzZXQiOi0xLCJzaG9ydGVuIjp7InNvdXJjZSI6MjAsInRhcmdldCI6MjB9LCJzdHlsZSI6eyJib2R5Ijp7Im5hbWUiOiJub25lIn0sImhlYWQiOnsibmFtZSI6Im5vbmUifX19XV0=
        \begin{tikzcd}
          GFx & GFx \\
          GFy \\
          GFz & GFz
          \arrow[""{name=0, anchor=center, inner sep=0}, "\shortmid"{marking}, equals, from=1-1, to=1-2]
          \arrow["GFm"', from=1-1, to=2-1]
          \arrow["{GF(m\odot n)}", from=1-2, to=3-2]
          \arrow["GFn"', from=2-1, to=3-1]
          \arrow[""{name=1, anchor=center, inner sep=0}, "\shortmid"{marking}, equals, from=3-1, to=3-2]
          \arrow["{(GF)^{m,n}}"{description}, shift left, draw=none, from=0, to=1]
        \end{tikzcd}
        \quad\coloneqq\quad
        \begin{tikzcd}[row sep=scriptsize]
          GFx & GFx & GFx \\
          & GFx & GFx \\
          GFy \\
          & GFz & GFz \\
          GFz & GFz & GFz
          \arrow[""{name=0, anchor=center, inner sep=0}, "\shortmid"{marking}, equals, from=1-1, to=1-2]
          \arrow["GFm"', from=1-1, to=3-1]
          \arrow[""{name=1, anchor=center, inner sep=0}, "{\id_{GFx}}"{inner sep=.8ex}, "\shortmid"{marking}, from=1-2, to=1-3]
          \arrow[equals, from=1-2, to=2-2]
          \arrow[equals, from=1-3, to=2-3]
          \arrow[""{name=2, anchor=center, inner sep=0}, "{G1_{Fx}}"'{inner sep=.8ex}, "\shortmid"{marking}, from=2-2, to=2-3]
          \arrow[from=2-2, to=4-2]
          \arrow["{GF(m\odot n)}", from=2-3, to=4-3]
          \arrow["GFn"', from=3-1, to=5-1]
          \arrow[""{name=3, anchor=center, inner sep=0}, "{G1_{Fz}}"{inner sep=.8ex}, "\shortmid"{marking}, from=4-2, to=4-3]
          \arrow[equals, from=4-2, to=5-2]
          \arrow[equals, from=4-3, to=5-3]
          \arrow[""{name=4, anchor=center, inner sep=0}, "\shortmid"{marking}, equals, from=5-1, to=5-2]
          \arrow[""{name=5, anchor=center, inner sep=0}, "{\id_{GFz}}"'{inner sep=.8ex}, "\shortmid"{marking}, from=5-2, to=5-3]
          \arrow["{G^{Fm,Fn}}"{description}, shift left, draw=none, from=0, to=4]
          \arrow["{G_{Fx}}"{description}, draw=none, from=1, to=2]
          \arrow["{G(F_{m,n})}"{description}, draw=none, from=2, to=3]
          \arrow["{G_{Fz}^{-1}}"{description}, draw=none, from=3, to=5]
        \end{tikzcd},
      \end{equation*}
      and for each object $x$ in $\dbl{C}$, the identity comparison is
      \begin{equation*}
        % https://q.uiver.app/#q=WzAsNCxbMCwwLCJHRngiXSxbMSwwLCJHRngiXSxbMSwxLCJHRngiXSxbMCwxLCJHRngiXSxbMCwxLCIiLDAseyJsZXZlbCI6Miwic3R5bGUiOnsiYm9keSI6eyJuYW1lIjoiYmFycmVkIn0sImhlYWQiOnsibmFtZSI6Im5vbmUifX19XSxbMSwyLCJHRlxcaWRfeCJdLFszLDIsIiIsMix7ImxldmVsIjoyLCJzdHlsZSI6eyJib2R5Ijp7Im5hbWUiOiJiYXJyZWQifSwiaGVhZCI6eyJuYW1lIjoibm9uZSJ9fX1dLFswLDMsIjFfe0dGeH0iLDJdLFs0LDYsIihHRileeCIsMSx7InNob3J0ZW4iOnsic291cmNlIjoyMCwidGFyZ2V0IjoyMH0sInN0eWxlIjp7ImJvZHkiOnsibmFtZSI6Im5vbmUifSwiaGVhZCI6eyJuYW1lIjoibm9uZSJ9fX1dXQ==
        \begin{tikzcd}
          GFx & GFx \\
          GFx & GFx
          \arrow[""{name=0, anchor=center, inner sep=0}, "\shortmid"{marking}, equals, from=1-1, to=1-2]
          \arrow["{1_{GFx}}"', from=1-1, to=2-1]
          \arrow["{GF\id_x}", from=1-2, to=2-2]
          \arrow[""{name=1, anchor=center, inner sep=0}, "\shortmid"{marking}, equals, from=2-1, to=2-2]
          \arrow["{(GF)^x}"{description}, draw=none, from=0, to=1]
        \end{tikzcd}
        \quad\coloneqq\quad
        \begin{tikzcd}
          GFx & GFx & GFx \\
          & GFx & GFx \\
          & GFx & GFx \\
          GFx & GFx & GFx
          \arrow[""{name=0, anchor=center, inner sep=0}, "\shortmid"{marking}, equals, from=1-1, to=1-2]
          \arrow["{1_{GFx}}"', from=1-1, to=4-1]
          \arrow[""{name=1, anchor=center, inner sep=0}, "{\id_{GFx}}"{inner sep=.8ex}, "\shortmid"{marking}, from=1-2, to=1-3]
          \arrow[equals, from=1-2, to=2-2]
          \arrow[equals, from=1-3, to=2-3]
          \arrow[""{name=2, anchor=center, inner sep=0}, "{G1_{Fx}}"{inner sep=.8ex}, "\shortmid"{marking}, from=2-2, to=2-3]
          \arrow["{G\id_{Fx}}"{description}, from=2-2, to=3-2]
          \arrow["{GF\id_x}", from=2-3, to=3-3]
          \arrow[""{name=3, anchor=center, inner sep=0}, "{G1_{Fx}}"'{inner sep=.8ex}, "\shortmid"{marking}, from=3-2, to=3-3]
          \arrow[equals, from=3-2, to=4-2]
          \arrow[equals, from=3-3, to=4-3]
          \arrow[""{name=4, anchor=center, inner sep=0}, "\shortmid"{marking}, equals, from=4-1, to=4-2]
          \arrow[""{name=5, anchor=center, inner sep=0}, "{\id_{GFx}}"'{inner sep=.8ex}, "\shortmid"{marking}, from=4-2, to=4-3]
          \arrow["{G^{Fx}}"{description}, draw=none, from=0, to=4]
          \arrow["{G_{Fx}}"{description, pos=0.4}, draw=none, from=1, to=2]
          \arrow["{G(F_x)}"{description}, draw=none, from=2, to=3]
          \arrow["{G_{Fx}^{-1}}"{description, pos=0.6}, draw=none, from=3, to=5]
        \end{tikzcd}.
        \qedhere
      \end{equation*}
  \end{itemize}
\end{construction}

\begin{lemma}[Well-definedness of precomposition]
  The composite $GF$ in \cref{construction:twisted-pre-composite} is well defined
  as a twisted doubly lax functor and is normal.

  Moreover, when $F: \dbl{C} \to \dbl{D}$ is a (pseudo) double functor and
  $G: \dbl{D} \twistto \dbl{E}$ is a twisted normal lax functor,
  $GF: \dbl{C} \twistto \dbl{E}$ is again a twisted normal lax functor.
\end{lemma}
\begin{proof}
  The proof is routine, though it involves a lengthy calculation if carried out in
  detail. First, note that the comparison cells have the correct shape: the
  comparison $(GF)_{f,g}$ is well-defined because the lax functor $F$ is
  strictly functorial on arrows, hence the equation $F(f \cdot g) = Ff \cdot Fg$ holds.
  We check the naturality condition for the loose-to-tight composition
  comparisons of $GF$. Given loosely composable cells
  $\inlineCell{a}{b}{x}{y}{m}{p}{f}{g}{\alpha}$ and
  $\inlineCell{b}{c}{y}{z}{n}{q}{g}{h}{\beta}$, calculate that
  \begin{align*}
      \begin{dblArray}{cc}
        GF\alpha & \Block{2-1}{(GF)^{p,q}} \\
        GF\beta &
      \end{dblArray} &=
      \begin{dblArray}{ccc}
        GF\alpha & \Block{2-1}{G^{Fp,Fq}} & \Block{2-1}{G(F_{p,q})} \\
        GF\beta & &
      \end{dblArray} &
      &(\text{definition}) \\ &=
      \begin{dblArray}{ccc}
        \Block{3-1}{G^{Fm,Fn}} & 1_{GFf} & G_{Fx} \\
        & G(F\alpha \odot F\beta) & G(F_{p,q}) \\
        & 1_{GFh} & G_{Fz}^{-1}
      \end{dblArray} &
      &(\text{loose-to-tight naturality})\\ &=
      \begin{dblArray}{ccc}
        \Block{3-1}{G^{Fm,Fn}} & 1_{GFf} & G_{Fx} \\
        & G(F\alpha \odot F\beta) & G(F_{p,q}) \\
        & \Block{1-2}{G_{Fh,Fz}}
      \end{dblArray} &
      &(\text{unitality}) \\ &=
      \begin{dblArray}{ccc}
        \Block{3-1}{G^{Fm,Fn}} & \quad 1_{GFf}\quad  & \quad G_{Fx}\quad \\
        & \Block{1-2}{G_{Ff,Fx}} \\
        & \Block{1-2}{G((F\alpha\odot F\beta) \cdot F_{p,q})}
      \end{dblArray} &
      &(\text{tight-to-loose naturality}) \\ &=
      \begin{dblArray}{cc}
        G^{Fm,Fn} & G((F\alpha\odot F\beta) \cdot F_{p,q})
      \end{dblArray} &
      &(\text{unitality}) \\ &=
      \begin{dblArray}{cc}
        G^{Fm,Fn} & G(F_{m,n} \cdot F(\alpha\odot\beta))
      \end{dblArray} &
      &(\text{naturality}) \\ &=
      \begin{dblArray}{ccc}
        \Block{3-1}{G^{Fm,Fn}} & \quad G_{Fa}\quad  & \quad 1_{GFf}\quad \\
        & \Block{1-2}{G_{Fa,Ff}} \\
        & \Block{1-2}{G(F_{m,n} \cdot F(\alpha\odot\beta))}
      \end{dblArray} &
      &(\text{unitality}) \\ &=
      \begin{dblArray}{ccc}
        \Block{3-1}{G^{Fm,Fn}} & G_{Fa} & 1_{GFf} \\
        & G(F_{m,n}) & GF(\alpha\odot\beta) \\
        & \Block{1-2}{G_{Fc,Fh}}
      \end{dblArray} &
      &(\text{tight-to-loose naturality}) \\ &=
      \begin{dblArray}{ccc}
        \Block{3-1}{G^{Fm,Fn}} & G_{Fa} & 1_{GFf} \\
        & G(F_{m,n}) & GF(\alpha\odot\beta) \\
        & G_{Fc}^{-1} & 1_{GFh}
      \end{dblArray} &
      &(\text{unitality}) \\ &=
      \begin{dblArray}{cc}
        (GF)^{Fm,Fn} & GF(\alpha\odot\beta)
      \end{dblArray}. &
      &(\text{definition})
  \end{align*}
  The verification of naturality for loose-to-tight identity comparisons is
  similar. Meanwhile, the tight-to-loose comparisons of $GF$ satisfy the
  required naturality and coherence conditions since they simply instantiate
  those of $G$ at components reindexed by $F$.

  Moreover, when $F$ is pseudo and $G$ is twisted normal lax, the comparisons
  $(GF)^{m,n}$ and $(GF)^x$ are invertible by the following lemma.
\end{proof}

\begin{lemma} \label{lem:tight-to-loose-invertibility}
  Let $F: \dbl{D} \twistto \dbl{E}$ be a twisted doubly lax functor that is
  normal. If $\alpha: m \To n$ is a tightly globular, tightly invertible cell
  between proarrows $m, n: x \proto y$ in $\dbl{D}$, then
  $F_x \cdot F\alpha \cdot F_y^{-1}$ is a loosely invertible cell in $\dbl{E}$.
\end{lemma}
\begin{proof}
  The loose inverse of $F_x \cdot F(\alpha) \cdot F_y^{-1}$ is
  $F_x \cdot F(\alpha^{-1}) \cdot F_y^{-1}$, where $\alpha^{-1}: n \To m$ is the
  tight inverse of $\alpha$. Checking the two equations for a loose inverse uses
  the naturality and unitality of the tight-to-loose comparisons of $F$.
\end{proof}

Pre-composing a twisted Hom functor with an ordinary double functor yields many
more examples of twisted functors. The twisted representables, to be discussed
shortly, are particularly important. Another example is:

\begin{example}[Twisted companions and conjoints]
  \label{ex:twisted-companion-conjoint}
  If $F\colon \dbl{D} \to \dbl{E}$ is an ordinary double functor, then its composites
  in either slot with the twisted Hom functor $\dbl{E}(-,=)$ are twisted normal
  lax functors
  \begin{equation*}
    \dbl{E}(F(-),=) \colon \dbl{D}^{\co} \times \dbl{E} \twistto \Prof
    \qquad\text{and}\qquad
    \dbl{E}(-,F(=)) \colon \dbl{E}^{\co} \times \dbl{D} \twistto \Prof.
  \end{equation*}
  Constructing such a ``twisted companion'' and ``twisted conjoint'' of a double
  functor is the first step toward defining a twisted double adjunction, an idea
  we revisit in \cref{ex:twisted-adjunction}.
\end{example}

\begin{lemma}[Functorality of precomposition]
  \label{lem:twisted-precomposition-functorality}
  Given a double category $\dbl{E}$, there is a functor $|\Dbll|_1^\op \to \Set$
  sending each double category $\dbl{D}$ to the set of twisted normal doubly lax
  functors $\dbl{D} \twistto \dbl{E}$ and acting with lax double functors by
  precomposition (\cref{construction:twisted-pre-composite}).
\end{lemma}
\begin{proof}
  Functorality for the underlying mappings and the tight-to-loose comparisons is
  immediate, as only function applications and reindexing are involved. As for
  the loose-to-tight comparisons, if $\dbl{A} \xto{F} \dbl{B} \xto{G} \dbl{C}$
  are lax functors and $H: \dbl{C} \twistto \dbl{E}$ is a twisted functor, then
  the loose-to-tight composition comparisons of $(HG)F$ and $H(GF)$ at proarrows
  $x \xproto{m} y \xproto{n} z$ in $\dbl{A}$ are
  \begin{equation*}
    \begin{dblArray}{ccc}
      \Block{3-1}{H^{GFm, GFn}} & H_{GFx} & H_{GFx} \\
      & H(G_{Fm,Fn}) & HG(F_{m,n}) \\
      & H_{GFz}^{-1} & H_{GFz}^{-1}
    \end{dblArray}
    \quad=\quad
    \begin{dblArray}{cc}
      \Block{3-1}{H^{GFm, GFn}} & H_{GFx} \\
      & H((GF)_{m,n}) \\
      & H_{GFz}^{-1}
    \end{dblArray}\,,
  \end{equation*}
  and similarly for the loose-to-tight identity comparisons.
\end{proof}

\begin{remark}\label{rmk:pseudo-pre-composition}
Though in \cref{construction:twisted-pre-composite} we made use of the strict functoriality of the lax functor $F \colon \dbl{C} \to \dbl{D}$ on tights, we do not believe this to be truly essential. Indeed, if $F \colon \dbl{C} \to \dbl{D}$ is only a double \emph{pseudofunctor} (as in Definition 6.1 of \cite{shulman2011}), then we may still restrict normal twisted doubly lax functors $G\colon \dbl{D} \twistto\dbl{E}$ by $F$, defining the composition comparison for composable tight arrows $x \xto{f} y \xto{g} z$ by

      \begin{equation*}
        % https://q.uiver.app/#q=WzAsNSxbMCwwLCJHRngiXSxbMSwwLCJHRnkiXSxbMiwwLCJHRnoiXSxbMiwxLCJHRnoiXSxbMCwxLCJHRngiXSxbMCwxLCJHRmYiLDAseyJzdHlsZSI6eyJib2R5Ijp7Im5hbWUiOiJiYXJyZWQifX19XSxbMSwyLCJHRmciLDAseyJzdHlsZSI6eyJib2R5Ijp7Im5hbWUiOiJiYXJyZWQifX19XSxbMiwzLCIiLDAseyJsZXZlbCI6Miwic3R5bGUiOnsiaGVhZCI6eyJuYW1lIjoibm9uZSJ9fX1dLFswLDQsIiIsMix7ImxldmVsIjoyLCJzdHlsZSI6eyJoZWFkIjp7Im5hbWUiOiJub25lIn19fV0sWzQsMywiR0YoZiBcXGNkb3QgZykiLDIseyJzdHlsZSI6eyJib2R5Ijp7Im5hbWUiOiJiYXJyZWQifX19XSxbMSw5LCIoR0YpX3tmLGd9IiwxLHsic2hvcnRlbiI6eyJ0YXJnZXQiOjIwfSwic3R5bGUiOnsiYm9keSI6eyJuYW1lIjoibm9uZSJ9LCJoZWFkIjp7Im5hbWUiOiJub25lIn19fV1d
        \begin{tikzcd}
          GFx & GFy & GFz \\
          GFx && GFz
          \arrow["GFf", "\shortmid"{marking}, from=1-1, to=1-2]
          \arrow[Rightarrow, no head, from=1-1, to=2-1]
          \arrow["GFg", "\shortmid"{marking}, from=1-2, to=1-3]
          \arrow[Rightarrow, no head, from=1-3, to=2-3]
          \arrow[""{name=0, anchor=center, inner sep=0}, "{GF(f \cdot g)}"', "\shortmid"{marking}, from=2-1, to=2-3]
          \arrow["{(GF)_{f,g}}"{description}, draw=none, from=1-2, to=0]
        \end{tikzcd}
        \quad\coloneqq\quad
% https://q.uiver.app/#q=WzAsNyxbMCwwLCJHRngiXSxbMSwwLCJHRnkiXSxbMiwwLCJHRnoiXSxbMiwxLCJHRnoiXSxbMCwxLCJHRngiXSxbMCwyLCJHRngiXSxbMiwyLCJHRnoiXSxbMCwxLCJHRmYiLDAseyJzdHlsZSI6eyJib2R5Ijp7Im5hbWUiOiJiYXJyZWQifX19XSxbMSwyLCJHRmciLDAseyJzdHlsZSI6eyJib2R5Ijp7Im5hbWUiOiJiYXJyZWQifX19XSxbMiwzLCIiLDAseyJsZXZlbCI6Miwic3R5bGUiOnsiaGVhZCI6eyJuYW1lIjoibm9uZSJ9fX1dLFswLDQsIiIsMix7ImxldmVsIjoyLCJzdHlsZSI6eyJoZWFkIjp7Im5hbWUiOiJub25lIn19fV0sWzQsMywiRyhGZiBcXGNkb3QgRmcpIiwyLHsic3R5bGUiOnsiYm9keSI6eyJuYW1lIjoiYmFycmVkIn19fV0sWzUsNiwiR0YoZiBcXGNkb3QgZykiLDIseyJzdHlsZSI6eyJib2R5Ijp7Im5hbWUiOiJiYXJyZWQifX19XSxbNCw1LCIiLDAseyJsZXZlbCI6Miwic3R5bGUiOnsiaGVhZCI6eyJuYW1lIjoibm9uZSJ9fX1dLFszLDYsIiIsMCx7ImxldmVsIjoyLCJzdHlsZSI6eyJoZWFkIjp7Im5hbWUiOiJub25lIn19fV0sWzEsMTEsIkdfe0ZmLEZnfSIsMSx7InNob3J0ZW4iOnsidGFyZ2V0IjoyMH0sInN0eWxlIjp7ImJvZHkiOnsibmFtZSI6Im5vbmUifSwiaGVhZCI6eyJuYW1lIjoibm9uZSJ9fX1dLFsxMSwxMiwiRyhGXntmLGp9KSIsMSx7InNob3J0ZW4iOnsic291cmNlIjoyMCwidGFyZ2V0IjoyMH0sInN0eWxlIjp7ImJvZHkiOnsibmFtZSI6Im5vbmUifSwiaGVhZCI6eyJuYW1lIjoibm9uZSJ9fX1dXQ==
\begin{tikzcd}
	GFx & GFy & GFz \\
	GFx && GFz \\
	GFx && GFz
	\arrow["GFf"{inner sep=.8ex}, "\shortmid"{marking}, from=1-1, to=1-2]
	\arrow[equals, from=1-1, to=2-1]
	\arrow["GFg"{inner sep=.8ex}, "\shortmid"{marking}, from=1-2, to=1-3]
	\arrow[equals, from=1-3, to=2-3]
	\arrow[""{name=0, anchor=center, inner sep=0}, "{G(Ff \cdot Fg)}"'{inner sep=.8ex}, "\shortmid"{marking}, from=2-1, to=2-3]
	\arrow[equals, from=2-1, to=3-1]
	\arrow[equals, from=2-3, to=3-3]
	\arrow[""{name=1, anchor=center, inner sep=0}, "{GF(f \cdot g)}"'{inner sep=.8ex}, "\shortmid"{marking}, from=3-1, to=3-3]
	\arrow["{G_{Ff,Fg}}"{description}, draw=none, from=1-2, to=0]
	\arrow["{G(F^{f,j})}"{description}, draw=none, from=0, to=1]
\end{tikzcd}
      \end{equation*}
      and we must similarly pad out the identity comparison.

      The verification of the laws in this case follows along the lines of \cref{lem:well-definedness}, though made more plodding by the additional padding and shuffling of coherences.
\end{remark}

\subsection{Twisted representable functors}
\label{subsection:twisted-representable-functors}

Restricting a twisted Hom functor along an ordinary point in one argument 
yields an important class of examples including all of the familiar indexed 
categories.

\begin{construction}[Twisted representables]
  \label{construction:twistedrepresentables}
  Let $a$ be any object in a double category $\dbl{D}$, viewed as a double
  functor $a: \dbl{1} \to \dbl{D}$. Composing with the twisted Hom functor on
  $\dbl{D}$ using \cref{construction:twisted-pre-composite} yields the
  \define{twisted representable functors}
  \begin{equation*}
    \dbl{D}(a, -) : \dbl{D} \twistto \Prof
    \qquad\text{and}\qquad
    \dbl{D}(-, a): \dbl{D}^{\co} \twistto \Prof.
  \end{equation*}
  It is worth explicitly describing one of the twisted representables.
  The covariant twisted representable $\dbl{D}(a,-)$ makes the following
  assignments.
  \begin{itemize}
    \item For each object $x$, the hom-category $\dbl{D}(a,x)$ has proarrows
      $a \proto x$ as objects and tightly globular cells as morphisms.
    \item For each arrow $f: x \to y$, the profunctor
      $\dbl{D}(a,f): \dbl{D}(a,x) \proto \dbl{D}(a,y)$ between hom-categories
      has as heteromorphisms from $m$ to $n$ the cells in $\dbl{D}$ of the form
      \begin{equation*}
        % https://q.uiver.app/#q=WzAsNCxbMCwwLCJhIl0sWzAsMSwiYSJdLFsxLDAsIngiXSxbMSwxLCJ5Il0sWzAsMiwibSIsMCx7InN0eWxlIjp7ImJvZHkiOnsibmFtZSI6ImJhcnJlZCJ9fX1dLFsxLDMsIm4iLDIseyJzdHlsZSI6eyJib2R5Ijp7Im5hbWUiOiJiYXJyZWQifX19XSxbMiwzLCJmIl0sWzAsMSwiIiwyLHsibGV2ZWwiOjIsInN0eWxlIjp7ImhlYWQiOnsibmFtZSI6Im5vbmUifX19XSxbNCw1LCJcXGFscGhhIiwxLHsic2hvcnRlbiI6eyJzb3VyY2UiOjIwLCJ0YXJnZXQiOjIwfSwic3R5bGUiOnsiYm9keSI6eyJuYW1lIjoibm9uZSJ9LCJoZWFkIjp7Im5hbWUiOiJub25lIn19fV1d
        \begin{tikzcd}
          a & x \\
          a & y
          \arrow[""{name=0, anchor=center, inner sep=0}, "m", "\shortmid"{marking}, from=1-1, to=1-2]
          \arrow[Rightarrow, no head, from=1-1, to=2-1]
          \arrow["f", from=1-2, to=2-2]
          \arrow[""{name=1, anchor=center, inner sep=0}, "n"', "\shortmid"{marking}, from=2-1, to=2-2]
          \arrow["\alpha"{description}, draw=none, from=0, to=1]
        \end{tikzcd}.
      \end{equation*}
    \item For each proarrow $p: x \proto y$, the functor
      $\dbl{D}(a,p) = (-) \odot p: \dbl{D}(a,x) \to \dbl{D}(a,y)$ between
      hom-categories maps a proarrow $m: a \proto x$ to the composite
      $m \odot p: a \proto y$ and a tightly globular cell $\alpha: m \To n$ to the
      composite $\alpha \odot 1_p: m \odot p \To n \odot p$.
    \item For each cell $\inlineCell{x}{w}{y}{z}{p}{q}{f}{g}{\phi}$, the natural
      transformation
      \begin{equation*}
        % https://q.uiver.app/#q=WzAsNCxbMCwwLCJcXGRibHtEfShhLHgpIl0sWzEsMCwiXFxkYmx7RH0oYSx5KSJdLFswLDEsIlxcZGJse0R9KGEsdykiXSxbMSwxLCJcXGRibHtEfShhLHopIl0sWzAsMSwiXFxkYmx7RH0oYSxmKSIsMCx7InN0eWxlIjp7ImJvZHkiOnsibmFtZSI6ImJhcnJlZCJ9fX1dLFswLDIsIlxcZGJse0R9KGEscCkiLDJdLFsxLDMsIlxcZGJse0R9KGEscSkiXSxbMiwzLCJcXGRibHtEfShhLGcpIiwyLHsic3R5bGUiOnsiYm9keSI6eyJuYW1lIjoiYmFycmVkIn19fV0sWzQsNywiXFxkYmx7RH0oYSxcXHBoaSkiLDEseyJzaG9ydGVuIjp7InNvdXJjZSI6MjAsInRhcmdldCI6MjB9LCJzdHlsZSI6eyJib2R5Ijp7Im5hbWUiOiJub25lIn0sImhlYWQiOnsibmFtZSI6Im5vbmUifX19XV0=
        \begin{tikzcd}
          {\dbl{D}(a,x)} & {\dbl{D}(a,y)} \\
          {\dbl{D}(a,w)} & {\dbl{D}(a,z)}
          \arrow[""{name=0, anchor=center, inner sep=0}, "{\dbl{D}(a,f)}", "\shortmid"{marking}, from=1-1, to=1-2]
          \arrow["{\dbl{D}(a,p)}"', from=1-1, to=2-1]
          \arrow["{\dbl{D}(a,q)}", from=1-2, to=2-2]
          \arrow[""{name=1, anchor=center, inner sep=0}, "{\dbl{D}(a,g)}"', "\shortmid"{marking}, from=2-1, to=2-2]
          \arrow["{\dbl{D}(a,\phi)}"{description}, draw=none, from=0, to=1]
        \end{tikzcd}
      \end{equation*}
      has components acting by loose post-composition with $\phi$:
      \begin{equation*}
        % https://q.uiver.app/#q=WzAsNCxbMCwwLCJhIl0sWzAsMSwiYSJdLFsxLDAsIngiXSxbMSwxLCJ5Il0sWzAsMiwibSIsMCx7InN0eWxlIjp7ImJvZHkiOnsibmFtZSI6ImJhcnJlZCJ9fX1dLFsxLDMsIm4iLDIseyJzdHlsZSI6eyJib2R5Ijp7Im5hbWUiOiJiYXJyZWQifX19XSxbMiwzLCJmIl0sWzAsMSwiIiwyLHsibGV2ZWwiOjIsInN0eWxlIjp7ImhlYWQiOnsibmFtZSI6Im5vbmUifX19XSxbNCw1LCJcXGFscGhhIiwxLHsic2hvcnRlbiI6eyJzb3VyY2UiOjIwLCJ0YXJnZXQiOjIwfSwic3R5bGUiOnsiYm9keSI6eyJuYW1lIjoibm9uZSJ9LCJoZWFkIjp7Im5hbWUiOiJub25lIn19fV1d
        \begin{tikzcd}
          a & x \\
          a & y
          \arrow[""{name=0, anchor=center, inner sep=0}, "m", "\shortmid"{marking}, from=1-1, to=1-2]
          \arrow[Rightarrow, no head, from=1-1, to=2-1]
          \arrow["f", from=1-2, to=2-2]
          \arrow[""{name=1, anchor=center, inner sep=0}, "n"', "\shortmid"{marking}, from=2-1, to=2-2]
          \arrow["\alpha"{description}, draw=none, from=0, to=1]
        \end{tikzcd}
        \quad\mapsto\quad
        % https://q.uiver.app/#q=WzAsNixbMCwwLCJhIl0sWzAsMSwiYSJdLFsxLDAsIngiXSxbMSwxLCJ5Il0sWzIsMCwidyJdLFsyLDEsInoiXSxbMCwyLCJtIiwwLHsic3R5bGUiOnsiYm9keSI6eyJuYW1lIjoiYmFycmVkIn19fV0sWzEsMywibiIsMix7InN0eWxlIjp7ImJvZHkiOnsibmFtZSI6ImJhcnJlZCJ9fX1dLFsyLDMsImYiLDFdLFswLDEsIiIsMix7ImxldmVsIjoyLCJzdHlsZSI6eyJoZWFkIjp7Im5hbWUiOiJub25lIn19fV0sWzIsNCwicCIsMCx7InN0eWxlIjp7ImJvZHkiOnsibmFtZSI6ImJhcnJlZCJ9fX1dLFszLDUsInEiLDIseyJzdHlsZSI6eyJib2R5Ijp7Im5hbWUiOiJiYXJyZWQifX19XSxbNCw1LCJnIl0sWzYsNywiXFxhbHBoYSIsMSx7InNob3J0ZW4iOnsic291cmNlIjoyMCwidGFyZ2V0IjoyMH0sInN0eWxlIjp7ImJvZHkiOnsibmFtZSI6Im5vbmUifSwiaGVhZCI6eyJuYW1lIjoibm9uZSJ9fX1dLFsxMCwxMSwiXFxwaGkiLDEseyJzaG9ydGVuIjp7InNvdXJjZSI6MjAsInRhcmdldCI6MjB9LCJzdHlsZSI6eyJib2R5Ijp7Im5hbWUiOiJub25lIn0sImhlYWQiOnsibmFtZSI6Im5vbmUifX19XV0=
        \begin{tikzcd}
          a & x & w \\
          a & y & z
          \arrow[""{name=0, anchor=center, inner sep=0}, "m", "\shortmid"{marking}, from=1-1, to=1-2]
          \arrow[Rightarrow, no head, from=1-1, to=2-1]
          \arrow[""{name=1, anchor=center, inner sep=0}, "p", "\shortmid"{marking}, from=1-2, to=1-3]
          \arrow["f"{description}, from=1-2, to=2-2]
          \arrow["g", from=1-3, to=2-3]
          \arrow[""{name=2, anchor=center, inner sep=0}, "n"', "\shortmid"{marking}, from=2-1, to=2-2]
          \arrow[""{name=3, anchor=center, inner sep=0}, "q"', "\shortmid"{marking}, from=2-2, to=2-3]
          \arrow["\alpha"{description}, draw=none, from=0, to=2]
          \arrow["\phi"{description}, draw=none, from=1, to=3]
        \end{tikzcd}.
      \end{equation*}
  \end{itemize}
  In addition, for each pair of arrows $f: x \to y$ and $g: y \to z$, the
  composition comparison cell
  \begin{equation*}
    % https://q.uiver.app/#q=WzAsNSxbMCwwLCJcXGRibHtEfShhLHgpIl0sWzAsMSwiXFxkYmx7RH0oYSx4KSJdLFsxLDAsIlxcZGJse0R9KGEseSkiXSxbMiwwLCJcXGRibHtEfShhLHopIl0sWzIsMSwiXFxkYmx7RH0oYSx6KSJdLFsxLDQsIlxcZGJse0R9KGEsIGYgXFxjZG90IGcpIiwyLHsic3R5bGUiOnsiYm9keSI6eyJuYW1lIjoiYmFycmVkIn19fV0sWzAsMSwiIiwyLHsibGV2ZWwiOjIsInN0eWxlIjp7ImhlYWQiOnsibmFtZSI6Im5vbmUifX19XSxbMyw0LCIiLDAseyJsZXZlbCI6Miwic3R5bGUiOnsiaGVhZCI6eyJuYW1lIjoibm9uZSJ9fX1dLFswLDIsIlxcZGJse0R9KGEsZikiLDAseyJzdHlsZSI6eyJib2R5Ijp7Im5hbWUiOiJiYXJyZWQifX19XSxbMiwzLCJcXGRibHtEfShhLGcpIiwwLHsic3R5bGUiOnsiYm9keSI6eyJuYW1lIjoiYmFycmVkIn19fV0sWzIsNSwiXFxkYmx7RH0oYSwtKV97ZixnfSIsMSx7InNob3J0ZW4iOnsidGFyZ2V0IjoyMH0sInN0eWxlIjp7ImJvZHkiOnsibmFtZSI6Im5vbmUifSwiaGVhZCI6eyJuYW1lIjoibm9uZSJ9fX1dXQ==
    \begin{tikzcd}
      {\dbl{D}(a,x)} & {\dbl{D}(a,y)} & {\dbl{D}(a,z)} \\
      {\dbl{D}(a,x)} && {\dbl{D}(a,z)}
      \arrow["{\dbl{D}(a,f)}", "\shortmid"{marking}, from=1-1, to=1-2]
      \arrow[Rightarrow, no head, from=1-1, to=2-1]
      \arrow["{\dbl{D}(a,g)}", "\shortmid"{marking}, from=1-2, to=1-3]
      \arrow[Rightarrow, no head, from=1-3, to=2-3]
      \arrow[""{name=0, anchor=center, inner sep=0}, "{\dbl{D}(a, f \cdot g)}"', "\shortmid"{marking}, from=2-1, to=2-3]
      \arrow["{\dbl{D}(a,-)_{f,g}}"{description}, draw=none, from=1-2, to=0]
    \end{tikzcd}
  \end{equation*}
  has components acting by tight composition in $\dbl{D}$:
  \begin{equation*}
    \begin{gathered}
      % https://q.uiver.app/#q=WzAsNCxbMCwwLCJhIl0sWzAsMSwiYSJdLFsxLDAsIngiXSxbMSwxLCJ5Il0sWzAsMiwibSIsMCx7InN0eWxlIjp7ImJvZHkiOnsibmFtZSI6ImJhcnJlZCJ9fX1dLFsxLDMsIm4iLDIseyJzdHlsZSI6eyJib2R5Ijp7Im5hbWUiOiJiYXJyZWQifX19XSxbMiwzLCJmIl0sWzAsMSwiIiwyLHsibGV2ZWwiOjIsInN0eWxlIjp7ImhlYWQiOnsibmFtZSI6Im5vbmUifX19XSxbNCw1LCJcXGFscGhhIiwxLHsic2hvcnRlbiI6eyJzb3VyY2UiOjIwLCJ0YXJnZXQiOjIwfSwic3R5bGUiOnsiYm9keSI6eyJuYW1lIjoibm9uZSJ9LCJoZWFkIjp7Im5hbWUiOiJub25lIn19fV1d
      \begin{tikzcd}
        a & x \\
        a & y
        \arrow[""{name=0, anchor=center, inner sep=0}, "m", "\shortmid"{marking}, from=1-1, to=1-2]
        \arrow[Rightarrow, no head, from=1-1, to=2-1]
        \arrow["f", from=1-2, to=2-2]
        \arrow[""{name=1, anchor=center, inner sep=0}, "n"', "\shortmid"{marking}, from=2-1, to=2-2]
        \arrow["\alpha"{description}, draw=none, from=0, to=1]
      \end{tikzcd}
      \\
      % https://q.uiver.app/#q=WzAsNCxbMCwwLCJhIl0sWzAsMSwiYSJdLFsxLDAsInkiXSxbMSwxLCJ6Il0sWzAsMiwibiIsMCx7InN0eWxlIjp7ImJvZHkiOnsibmFtZSI6ImJhcnJlZCJ9fX1dLFsxLDMsInAiLDIseyJzdHlsZSI6eyJib2R5Ijp7Im5hbWUiOiJiYXJyZWQifX19XSxbMiwzLCJnIl0sWzAsMSwiIiwyLHsibGV2ZWwiOjIsInN0eWxlIjp7ImhlYWQiOnsibmFtZSI6Im5vbmUifX19XSxbNCw1LCJcXGJldGEiLDEseyJzaG9ydGVuIjp7InNvdXJjZSI6MjAsInRhcmdldCI6MjB9LCJzdHlsZSI6eyJib2R5Ijp7Im5hbWUiOiJub25lIn0sImhlYWQiOnsibmFtZSI6Im5vbmUifX19XV0=
      \begin{tikzcd}
        a & y \\
        a & z
        \arrow[""{name=0, anchor=center, inner sep=0}, "n", "\shortmid"{marking}, from=1-1, to=1-2]
        \arrow[Rightarrow, no head, from=1-1, to=2-1]
        \arrow["g", from=1-2, to=2-2]
        \arrow[""{name=1, anchor=center, inner sep=0}, "p"', "\shortmid"{marking}, from=2-1, to=2-2]
        \arrow["\beta"{description}, draw=none, from=0, to=1]
      \end{tikzcd}
    \end{gathered}
    \quad\mapsto\quad
    % https://q.uiver.app/#q=WzAsNixbMCwxLCJhIl0sWzAsMiwiYSJdLFsxLDEsInkiXSxbMSwyLCJ6Il0sWzAsMCwiYSJdLFsxLDAsIngiXSxbMCwyLCJuIiwwLHsic3R5bGUiOnsiYm9keSI6eyJuYW1lIjoiYmFycmVkIn19fV0sWzEsMywicCIsMix7InN0eWxlIjp7ImJvZHkiOnsibmFtZSI6ImJhcnJlZCJ9fX1dLFsyLDMsImciXSxbMCwxLCIiLDIseyJsZXZlbCI6Miwic3R5bGUiOnsiaGVhZCI6eyJuYW1lIjoibm9uZSJ9fX1dLFs0LDUsIm0iLDAseyJzdHlsZSI6eyJib2R5Ijp7Im5hbWUiOiJiYXJyZWQifX19XSxbNSwyLCJmIl0sWzQsMCwiIiwwLHsibGV2ZWwiOjIsInN0eWxlIjp7ImhlYWQiOnsibmFtZSI6Im5vbmUifX19XSxbNiw3LCJcXGJldGEiLDEseyJzaG9ydGVuIjp7InNvdXJjZSI6MjAsInRhcmdldCI6MjB9LCJzdHlsZSI6eyJib2R5Ijp7Im5hbWUiOiJub25lIn0sImhlYWQiOnsibmFtZSI6Im5vbmUifX19XSxbMTAsNiwiXFxhbHBoYSIsMSx7InNob3J0ZW4iOnsic291cmNlIjoyMCwidGFyZ2V0IjoyMH0sInN0eWxlIjp7ImJvZHkiOnsibmFtZSI6Im5vbmUifSwiaGVhZCI6eyJuYW1lIjoibm9uZSJ9fX1dXQ==
    \begin{tikzcd}
      a & x \\
      a & y \\
      a & z
      \arrow[""{name=0, anchor=center, inner sep=0}, "m", "\shortmid"{marking}, from=1-1, to=1-2]
      \arrow[Rightarrow, no head, from=1-1, to=2-1]
      \arrow["f", from=1-2, to=2-2]
      \arrow[""{name=1, anchor=center, inner sep=0}, "n", "\shortmid"{marking}, from=2-1, to=2-2]
      \arrow[Rightarrow, no head, from=2-1, to=3-1]
      \arrow["g", from=2-2, to=3-2]
      \arrow[""{name=2, anchor=center, inner sep=0}, "p"', "\shortmid"{marking}, from=3-1, to=3-2]
      \arrow["\alpha"{description}, draw=none, from=0, to=1]
      \arrow["\beta"{description}, draw=none, from=1, to=2]
    \end{tikzcd}.
  \end{equation*}
  The other comparison cells are either trivial or given by associators or
  unitors in $\dbl{D}$.
\end{construction}

\begin{example}
  Consider a special case of the previous \cref{construction:twistedrepresentables}.
  Let $\dbl{C}$ denote the \emph{horizontal categorification} of an ordinary
  monoidal category $(\cat{C},\otimes,I)$. That is, $\dbl{C}$ is the double
  category on a single object, say $\star$, with only an identity arrow, and
  whose proarrows are the objects of $\cat{C}$ and whose cells (all tightly globular)
  are the morphisms of $\cat{C}$. The loose unit and composition are the
  monoidal unit $I$ and monoidal tensor $\otimes$. Trivially $\dbl{C}$ is an
  equipment.

  The twisted representable
    \begin{equation*}
      \dbl{C}(\star,-)\colon\dbl{C}\twistto\Prof
    \end{equation*}
  has the action sending
    \begin{itemize}
      \item the single object $\star$ to the hom category $\dbl{C}(\star,\star)
      \cong \cat{C}$;
      \item the identity $1_\star$ to the endo-profunctor on $\dbl{C}(\star,\star)$ taking proarrows $x, y\colon \star\proto\star$ to the set of all cells $x\Rightarrow y$;
      \item any proarrow $\star\xproto{x}\star$ to the endo-functor
      suggested by
        \begin{equation*}
          \cat{C} \xto{-\otimes x} \cat{C} \qquad\qquad y\mapsto y\otimes x
        \end{equation*}
      and on arrows (that is, cells in $\dbl{C}$) the action is \textit{whiskering};
      \item and finally each cell $f\colon x\Rightarrow y$ is sent to the transformation $-\otimes f$ given by whiskering with $f$.
    \end{itemize}
  The comparison isos for the action endo-functors $-\otimes x$ have components coming from the monoidal structure on $\cat{C}$. These satisfy the required conditions by the axioms for a monoidal category. Similarly for units.
\end{example}

Twisted representables often give rise to bifibrations, or rather to their
indexed analogues. As the next example shows, this includes the canonical
self-indexing of a category with finite limits, corresponding to the codomain
bifibration \cite[\S{1.1}]{jacobs1999}.

\begin{example}[Self-indexing] \label{ex:self-indexing}
  Let $\cat{C}$ be a category with finite limits. The twisted representable
  functor
  \begin{equation*}
    \Span(\cat{C})(1, -): \Span(\cat{C}) \twistto \Prof
  \end{equation*}
  at any terminal object in $\cat{C}$ sends
  \begin{itemize}
    \item each object $a$ in $\cat{C}$ to the slice category over $a$, namely
      \begin{equation*}
        \Span(\cat{C})(1, a) = \cat{C}/a;
      \end{equation*}
    \item each span $m = (a \xfrom{f} s \xto{g} b)$ in $\cat{C}$ to the functor
      \begin{equation*}
        \Span(\cat{C})(1,m) = (-) \odot m: \cat{C}/a \to \cat{C}/b
      \end{equation*}
      acting on a morphism $\phi: x \to a$ by pulling back along $f$ and then
      post-composing with $g$:
      \begin{equation*}
        % https://q.uiver.app/#q=WzAsNixbMCwwLCJ4Il0sWzAsMSwiYSJdLFsxLDEsInMiXSxbMSwwLCJ4IFxcdGltZXNfYSBzIl0sWzIsMSwiYiJdLFsyLDAsInggXFx0aW1lc19hIHMiXSxbMCwxLCJcXHBoaSIsMl0sWzIsMSwiZiJdLFszLDIsImZeKiBcXHBoaSJdLFszLDBdLFsyLDQsImciLDJdLFszLDEsIiIsMCx7InN0eWxlIjp7Im5hbWUiOiJjb3JuZXIifX1dLFszLDUsIiIsMCx7ImxldmVsIjoyLCJzdHlsZSI6eyJoZWFkIjp7Im5hbWUiOiJub25lIn19fV0sWzUsNCwiZ18hIGZeKiBcXHBoaSAiXV0=
        \begin{tikzcd}
          x & {x \times_a s} & {x \times_a s} \\
          a & s & b
          \arrow["\phi"', from=1-1, to=2-1]
          \arrow[from=1-2, to=1-1]
          \arrow[Rightarrow, no head, from=1-2, to=1-3]
          \arrow["\lrcorner"{anchor=center, pos=0.125, rotate=-90}, draw=none, from=1-2, to=2-1]
          \arrow["{f^* \phi}", from=1-2, to=2-2]
          \arrow["{g_! f^* \phi }", from=1-3, to=2-3]
          \arrow["f", from=2-2, to=2-1]
          \arrow["g"', from=2-2, to=2-3]
        \end{tikzcd}.
      \end{equation*}
  \end{itemize}
  As the notation suggests, the functor at a general span
  $m = (a \xfrom{f} s \xto{g} b)$ decomposes as
  \begin{equation*}
    \cat{C}/a \xto{(-) \odot f^*} \cat{C}/s \xto{(-) \odot g_!} \cat{C}/b,
  \end{equation*}
  where the span $f^* = (a \xfrom{f} s \xto{1_s} s)$ is the conjoint of $f$ and
  $g_! = (s \xfrom{1_s} s \xto{g} b)$ is the companion of $g$ in
  $\Span(\cat{C})$. For any morphism $f$ in $\cat{C}$, the functors given by
  post-composition with $f_!$ and $f^*$ are adjoint. The twisted representable
  $\Span(\cat{C})(1,-)$ thus contains within it the self-indexing of $\cat{C}$,
  a pseudofunctor
  \begin{equation*}
    \Span(\cat{C})_0 = \cat{C} \to \AdjTwo
  \end{equation*}
  into the 2-category of adjunctions recalled below.
\end{example}

\begin{example}[Subobjects] \label{ex:subobjects}
  Let $\cat{C}$ be a regular category. We then have the double category
  $\Rel(\cat{C})$ of relations in $\cat{C}$. The twisted representable functor
    \begin{equation*}
      \Rel(\cat{C})(1, -): \Rel(\cat{C}) \twistto \Prof
    \end{equation*}
  at a terminal object in $\cat{C}$ sends
    \begin{itemize}
      \item each object $a$ in $\cat{C}$ to the category of subobjects of $a$,
      that is,
        \begin{equation*}
          \Rel(\cat{C})(1, a) \cong \Sub(a);
        \end{equation*}
      \item each relation $m = (\langle d,c\rangle\colon r\to a\times b)$
      in $\cat{C}$ to the functor
        \begin{equation*}
          \Rel(\cat{C})(1,m) : \Sub(a) \to \Sub(b)
        \end{equation*}
      that acts on a representative monomorphism $s: u \rightarrowtail a$ by
      first pulling back along $d$ and then taking the image of the pushforward
      along $c$:
        \begin{equation*}
          % https://q.uiver.app/#q=WzAsNCxbMCwwLCJ1XFx0aW1lc19hciJdLFsxLDAsInIiXSxbMSwxLCJhIl0sWzAsMSwidSJdLFswLDEsInBfMiJdLFsxLDIsImQiXSxbMCwzLCJwXzEiLDJdLFszLDIsInMiLDIseyJzdHlsZSI6eyJ0YWlsIjp7Im5hbWUiOiJtb25vIn19fV0sWzAsMiwiIiwxLHsic3R5bGUiOnsibmFtZSI6ImNvcm5lciJ9fV1d
          \begin{tikzcd}
            {u\times_ar} & r \\
            u & a
            \arrow["{p_2}", from=1-1, to=1-2]
            \arrow["{p_1}"', from=1-1, to=2-1]
            \arrow["\lrcorner"{anchor=center, pos=0.125}, draw=none, from=1-1, to=2-2]
            \arrow["d", from=1-2, to=2-2]
            \arrow["s"', tail, from=2-1, to=2-2]
          \end{tikzcd}
          \qquad\qquad
          % https://q.uiver.app/#q=WzAsNCxbMCwwLCJ1XFx0aW1lc19hciJdLFsxLDAsInIiXSxbMiwwLCJiIl0sWzEsMSwidCJdLFswLDEsInBfMiJdLFsxLDIsImMiXSxbMCwzLCIiLDIseyJzdHlsZSI6eyJoZWFkIjp7Im5hbWUiOiJlcGkifX19XSxbMywyLCJcXFJlbCgxLG0pKHMpIiwyLHsic3R5bGUiOnsidGFpbCI6eyJuYW1lIjoibW9ubyJ9fX1dXQ==
          \begin{tikzcd}
            {u\times_ar} & r & b \\
            & t
            \arrow["{p_2}", from=1-1, to=1-2]
            \arrow[two heads, from=1-1, to=2-2]
            \arrow["c", from=1-2, to=1-3]
            \arrow["{\Rel(1,m)(s)}"', tail, from=2-2, to=1-3]
          \end{tikzcd}.
        \end{equation*}
    \end{itemize}
  As with the previous example, this functor induced by any such relation $\langle d,c\rangle$ factors as
    \begin{equation*}
      \Sub(a) \xto{(-) \odot d^*} \Sub(r) \xto{(-) \odot c_!} \Sub(b),
    \end{equation*}
  where now $d^*$ and $c_!$ are the cograph and the graph of $d$ and $c$, respectively,
  that is, $d^* = \langle d,1\rangle$ and $c_!=\langle 1,c\rangle$. For despite
  the incorporation of images into the construction, the composition of relations
  given in this way is indeed associative up to isomorphism. Again companions
  and conjoints in $\Rel(\cat{C})$ are adjoint in the horizontal bicategory, so that
  the twisted representable $\Rel(\cat{C})(1,-)$ in a manner of speaking contains
  the indexed category associated to the \emph{subobject bifibration}, which we
  think of as a pseudofunctor $\Rel(\cat{C})_0 = \cat{C} \to \AdjTwo$ into
  the 2-category of adjunctions.

  Now, let $\Adj$ denote the double category of categories, functors and 
  adjunctions \cite[\S{3.1.5}]{grandis2019}. Let $\cat{E}$ denote any
  well-powered regular category. Define a twisted functor
    \begin{equation*}
      \Sub(-)\colon \Rel(\cat{E})\twistto \Adj
    \end{equation*}
  by the assignments sending
    \begin{itemize}
      \item an object $x$ to the subobject lattice $\Sub(x)$;
      \item a morphism $f\colon x\to y$ to the adjunction
        \begin{equation*}
          \exists_f\colon \Sub(x) \rightleftarrows \Sub(y) \colon f^{-1},
        \end{equation*}
        where $\exists_f$ and $f^{-1}$ are the image and inverse image functors;
      \item a relation $\langle d,c\rangle\colon r\to x\times y$ to the functor
        \begin{equation*}
          -\odot \langle d,c\rangle \colon \Sub(x)\to \Sub(y)
        \end{equation*}
        sending a subobject (representative) $m\colon s \rightarrowtail x$ to the subobject
        of $y$ represented by the composite relation $m\odot \langle d,c\rangle$
        on $1\times y\cong y$.
    \end{itemize}
  If $\cat{E}$ is, for example, a topos \cite{maclane1992}, then we also have a
  contravariant twisted functor
    \begin{equation*}
      \Sub(-)\colon \Rel(\cat{E})^\op \twistto \Adj
    \end{equation*}
  making the same assignments but sending a morphism $f\colon x \to y$ to the
  adjunction $f^*\dashv f_*$ where now $f_*\colon \Sub(x)\to \Sub(y)$ is the
  direct image functor.
\end{example}

\begin{example}[Copresheaves]
  Let $\ADJ$ denote the double category of (large) categories, functors and
  adjunctions. Size issues are at play here, so we insist on the special notation 
  $\ADJ$ rather than $\Adj$ to indicate that we may have to leave merely small 
  categories. There is a twisted functor
    \begin{equation*}
      [-,\Set]\colon\Prof \twistto \ADJ
    \end{equation*}
  making the assignments sending
    \begin{itemize}
      \item a category $\cat{C}$ to the (possibly large) category of copresheaves
      $[\cat{C},\Set]$;
      \item a functor $F\colon\cat{C}\to\cat{D}$ to the adjunction
        \begin{equation*}
          F_!\colon [\cat{C},\Set] \rightleftarrows [\cat{D},\Set]\colon F^*,
        \end{equation*}
      where $F_!$ is the left Kan extension of $F$ and $F^*$ is the usual
      pullback/substitution functor;
      \item a profunctor $P\colon \cat{C}\proto\cat{D}$ to the functor
        \begin{equation*}
        -\odot P\colon [\cat{C},\Set]\to [\cat{D},\Set] \qquad T\mapsto T\odot P,
        \end{equation*}
        where `$\odot$' denotes profunctor composition.
    \end{itemize}
  There are variations on this paradigm. For example, there is a twisted lax functor
    \begin{equation*}
      [-,\Set]\colon\Prof^{\op} \twistto \ADJ
    \end{equation*}
  making the same assignments except that functors $F$ are sent to the induced
  \emph{geometric morphism} $F^*\dashv F_*$ whose direct image is the
  \emph{right} Kan extension. Conventionally this points in the same direction as
  the direct image $F_*$, but in $\ADJ$, a proarrow points in the same direction as
  the left adjoint.
\end{example}

\subsection{Twisted functors on equipments}

Using a bit of nonstandard terminology, we say that a \define{bi-indexed
  category} over a category $\cat{C}$ is a pseudofunctor $\cat{C} \to \AdjTwo$
\cite[\mbox{Definition 2.1.28}]{michel2011}, where $\AdjTwo$ is the 2-category
of categories, adjoint pairs of functors, and conjugate pairs of natural
transformations. Bi-indexed categories correspond to bifibrations via the
Grothendieck construction in the same way that indexed categories do to
fibrations and op-indexed categories do to opfibrations. For a precise
statement, see \cite[\mbox{Theorem 2.1.37}]{michel2011}.

In the previous section, we saw that the twisted representable on the double
category of spans contains a bi-indexed category, namely the canonical
self-indexing. The emergence of a bi-indexed category is not a coincidence. In
this section, we will see that the crucial property of the double category of
spans making this happen is that it is an \emph{equipment}, meaning among other
equivalent definitions that it has all companions and conjoints
\cite[\mbox{Theorem 4.1}]{shulman2008}. We begin with a fact about how twisted
double functors interact with companions and conjoints.

\begin{proposition}[Twisting companions and conjoints]
  \label{prop:twisting-companions-conjoints}
  Let $F: \dbl{D} \twistto \dbl{E}$ be a twisted doubly normal lax functor. If
  $f$ is an arrow with a companion $f_!$ in $\dbl{D}$, then the arrow $F(f_!)$
  has a companion in $\dbl{E}$, namely $Ff$. Likewise, if $g$ is an arrow with a
  conjoint $g^*$ in $\dbl{D}$, then the arrow $F(g^*)$ has a conjoint in
  $\dbl{E}$, namely $Fg$.
\end{proposition}
\begin{proof}
  We prove the statement for companions; the case of conjoints is dual. Suppose
  that $f_!: x \proto y$ is a companion of an arrow $f: x \to y$ in $\dbl{D}$,
  with binding cells
  \begin{equation*}
    % https://q.uiver.app/#q=WzAsNCxbMCwwLCJ4Il0sWzAsMSwieCJdLFsxLDAsIngiXSxbMSwxLCJ5Il0sWzAsMSwiIiwwLHsibGV2ZWwiOjIsInN0eWxlIjp7ImhlYWQiOnsibmFtZSI6Im5vbmUifX19XSxbMCwyLCIiLDIseyJsZXZlbCI6Miwic3R5bGUiOnsiYm9keSI6eyJuYW1lIjoiYmFycmVkIn0sImhlYWQiOnsibmFtZSI6Im5vbmUifX19XSxbMiwzLCJmIl0sWzEsMywiZl8hIiwyLHsic3R5bGUiOnsiYm9keSI6eyJuYW1lIjoiYmFycmVkIn19fV0sWzUsNywiXFxldGEiLDEseyJzaG9ydGVuIjp7InNvdXJjZSI6MjAsInRhcmdldCI6MjB9LCJzdHlsZSI6eyJib2R5Ijp7Im5hbWUiOiJub25lIn0sImhlYWQiOnsibmFtZSI6Im5vbmUifX19XV0=
    \begin{tikzcd}
      x & x \\
      x & y
      \arrow[""{name=0, anchor=center, inner sep=0}, "\shortmid"{marking}, Rightarrow, no head, from=1-1, to=1-2]
      \arrow[Rightarrow, no head, from=1-1, to=2-1]
      \arrow["f", from=1-2, to=2-2]
      \arrow[""{name=1, anchor=center, inner sep=0}, "{f_!}"', "\shortmid"{marking}, from=2-1, to=2-2]
      \arrow["\eta"{description}, draw=none, from=0, to=1]
    \end{tikzcd}
    \qquad\text{and}\qquad
    % https://q.uiver.app/#q=WzAsNCxbMCwwLCJ4Il0sWzEsMCwieSJdLFswLDEsInkiXSxbMSwxLCJ5Il0sWzAsMSwiZl8hIiwwLHsic3R5bGUiOnsiYm9keSI6eyJuYW1lIjoiYmFycmVkIn19fV0sWzAsMiwiZiIsMl0sWzEsMywiIiwwLHsibGV2ZWwiOjIsInN0eWxlIjp7ImhlYWQiOnsibmFtZSI6Im5vbmUifX19XSxbMiwzLCIiLDIseyJsZXZlbCI6Miwic3R5bGUiOnsiYm9keSI6eyJuYW1lIjoiYmFycmVkIn0sImhlYWQiOnsibmFtZSI6Im5vbmUifX19XSxbNCw3LCJcXHZhcmVwc2lsb24iLDEseyJzaG9ydGVuIjp7InNvdXJjZSI6MjAsInRhcmdldCI6MjB9LCJzdHlsZSI6eyJib2R5Ijp7Im5hbWUiOiJub25lIn0sImhlYWQiOnsibmFtZSI6Im5vbmUifX19XV0=
    \begin{tikzcd}
      x & y \\
      y & y
      \arrow[""{name=0, anchor=center, inner sep=0}, "{f_!}", "\shortmid"{marking}, from=1-1, to=1-2]
      \arrow["f"', from=1-1, to=2-1]
      \arrow[Rightarrow, no head, from=1-2, to=2-2]
      \arrow[""{name=1, anchor=center, inner sep=0}, "\shortmid"{marking}, Rightarrow, no head, from=2-1, to=2-2]
      \arrow["\varepsilon"{description}, draw=none, from=0, to=1]
    \end{tikzcd}.
  \end{equation*}
  To exhibit $Ff: Fx \proto Fy$ as a companion of $F(f_!): Fx \to Fy$, we form the
  binding cells
  \begin{equation*}
    \eta' \coloneqq
    \begin{tikzcd}
      Fx & Fx & Fx \\
      Fx & Fx & Fx \\
      Fx & Fx & Fy
      \arrow[""{name=0, anchor=center, inner sep=0}, "\shortmid"{marking}, equals, from=1-1, to=1-2]
      \arrow[equals, from=1-1, to=2-1]
      \arrow[""{name=1, anchor=center, inner sep=0}, "\shortmid"{marking}, equals, from=1-2, to=1-3]
      \arrow[equals, from=1-2, to=2-2]
      \arrow[equals, from=1-3, to=2-3]
      \arrow[""{name=2, anchor=center, inner sep=0}, "\shortmid"{marking}, equals, from=2-1, to=2-2]
      \arrow[equals, from=2-1, to=3-1]
      \arrow[""{name=3, anchor=center, inner sep=0}, "{F1_x}"{inner sep=.8ex}, "\shortmid"{marking}, from=2-2, to=2-3]
      \arrow["{F\id_x}"{description}, from=2-2, to=3-2]
      \arrow["{F(f_!)}", from=2-3, to=3-3]
      \arrow[""{name=4, anchor=center, inner sep=0}, "\shortmid"{marking}, equals, from=3-1, to=3-2]
      \arrow[""{name=5, anchor=center, inner sep=0}, "Ff"'{inner sep=.8ex}, "\shortmid"{marking}, from=3-2, to=3-3]
      \arrow["{1=\id}"{description}, draw=none, from=0, to=2]
      \arrow["{F_x}"{description, pos=0.4}, draw=none, from=1, to=3]
      \arrow["{F^x}"{description}, draw=none, from=2, to=4]
      \arrow["{F\eta}"{description}, draw=none, from=3, to=5]
    \end{tikzcd}
    \qquad\text{and}\qquad
    \varepsilon' \coloneqq
    \begin{tikzcd}
      Fx & Fy & Fy \\
      Fy & Fy & Fy \\
      Fy & Fy & Fy
      \arrow[""{name=0, anchor=center, inner sep=0}, "Ff"{inner sep=.8ex}, "\shortmid"{marking}, from=1-1, to=1-2]
      \arrow["{F(f_!)}"', from=1-1, to=2-1]
      \arrow[""{name=1, anchor=center, inner sep=0}, "\shortmid"{marking}, equals, from=1-2, to=1-3]
      \arrow["{F\id_y}"{description}, from=1-2, to=2-2]
      \arrow[equals, from=1-3, to=2-3]
      \arrow[""{name=2, anchor=center, inner sep=0}, "{F1_y}"'{inner sep=.8ex}, "\shortmid"{marking}, from=2-1, to=2-2]
      \arrow[equals, from=2-1, to=3-1]
      \arrow[""{name=3, anchor=center, inner sep=0}, "\shortmid"{marking}, equals, from=2-2, to=2-3]
      \arrow[equals, from=2-2, to=3-2]
      \arrow[equals, from=2-3, to=3-3]
      \arrow[""{name=4, anchor=center, inner sep=0}, "\shortmid"{marking}, equals, from=3-1, to=3-2]
      \arrow[""{name=5, anchor=center, inner sep=0}, "\shortmid"{marking}, equals, from=3-2, to=3-3]
      \arrow["{F\varepsilon}"{description}, draw=none, from=0, to=2]
      \arrow["{(F^y)^{-1}}"{description}, draw=none, from=1, to=3]
      \arrow["{F_y^{-1}}"{description, pos=0.6}, draw=none, from=2, to=4]
      \arrow["{1=\id}"{description}, draw=none, from=3, to=5]
    \end{tikzcd}.
  \end{equation*}
  We must show that the cells $\eta'$ and $\varepsilon'$ satisfy the two companion
  identities. The identity $\eta' \cdot \varepsilon' = \id_{Ff_!}$ is proved by
  calculating
  \begin{equation*}
    \begin{dblArray}{c}
      \eta' \\
      \varepsilon'
    \end{dblArray}
    =
    \begin{dblArray}{ccc}
      \id_{1_{Fx}} & F_x & \Block{2-1}{\id_{Ff_!}} \\
      F^x & F\eta & \\
      \Block{2-1}{\id_{Ff_!}} & F\varepsilon & (F^y)^{-1} \\
      & F_y^{-1} & \id_{1_{Fx}}
    \end{dblArray}
    =
    \begin{dblArray}{ccc}
      \id_{1_{Fx}} & F_x & \id_{1_{Fx}} \\
      F^x & F\eta & \Block{2-1}{F^{f_!,\id_y}} \\
      \Block{2-1}{\id_{Ff_!}} & F\varepsilon & \\
      & F_y^{-1} & \id_{1_{Fx}}
    \end{dblArray}
    =
    \begin{dblArray}{ccc}
      \Block{1-2}{\id_{1_{Fx}}} & & F_x \\
      F^x & \Block{2-1}{F^{\id_x, f_!}} & \Block{2-1}{F(\eta \odot \varepsilon)} \\
      \id_{Ff_!} & & \\
      \Block{1-2}{\id_{1_{Fx}}} & & F_y^{-1}
    \end{dblArray}
    = \id_{Ff_!},
  \end{equation*}
  where the second equality uses the unitality of composition comparisons, the
  third equality uses the naturality of composition comparisons, and the final
  equality uses the unitality of the composition comparisons again, as well as
  the naturality of the unit comparisons to compute that
  \begin{equation*}
    F_x \cdot F(\eta \odot \varepsilon) \cdot F_y^{-1}
      = F_x \cdot F(1_{f_!}) \cdot F_y^{-1} = \id_{Ff_!}.
  \end{equation*}
  Transposing this calculation proves the other companion identity, namely
  $\eta' \odot \varepsilon' = 1_{Ff}$.
\end{proof}

\begin{corollary} \label{cor:twisting-companion-conjoint-adjunctions}
  Let $F: \dbl{D} \twistto \dbl{E}$ be a twisted doubly normal lax functor. If
  $f$ is an arrow in $\dbl{D}$ with a companion $f_!$ and a conjoint $f^*$, then
  there is an adjunction $F(f_!) \dashv F(f^*)$ between arrows in the 2-category
  underlying $\dbl{E}$.
\end{corollary}
\begin{proof}
  By \cref{prop:twisting-companions-conjoints}, the proarrow $Ff$ is a companion
  of $F(f_!)$ and also a conjoint of $F(f^*)$ in $\dbl{E}$. Thus, by the
  transpose of \cite[\mbox{Proposition 1.4}]{grandis2004} or \cite[\mbox{Lemma
    3.21}]{shulman2010}, the arrows $Ff_!$ and $Ff^*$ are adjoint in the
  underlying 2-category of $\dbl{E}$.
\end{proof}

For the next theorem, recall that double categories of squares and of
adjunctions can be defined relative to any 2-category. Given a 2-category
$\twocat{C}$, the \define{double category of squares} in $\twocat{C}$, denoted
$\Sq(\twocat{C})$, has the same objects and arrows as $\twocat{C}$, has as its
proarrows also the arrows of $\twocat{C}$, and has as its cells of the form on
the left
\begin{equation*}
  % https://q.uiver.app/#q=WzAsNCxbMCwwLCJ4Il0sWzEsMCwieSJdLFswLDEsInciXSxbMSwxLCJ6Il0sWzAsMiwiZiIsMl0sWzEsMywiZyJdLFswLDEsImgiLDAseyJzdHlsZSI6eyJib2R5Ijp7Im5hbWUiOiJiYXJyZWQifX19XSxbMiwzLCJrIiwyLHsic3R5bGUiOnsiYm9keSI6eyJuYW1lIjoiYmFycmVkIn19fV0sWzYsNywiXFxhbHBoYSIsMSx7InNob3J0ZW4iOnsic291cmNlIjoyMCwidGFyZ2V0IjoyMH0sInN0eWxlIjp7ImJvZHkiOnsibmFtZSI6Im5vbmUifSwiaGVhZCI6eyJuYW1lIjoibm9uZSJ9fX1dXQ==
  \begin{tikzcd}
    x & y \\
    w & z
    \arrow[""{name=0, anchor=center, inner sep=0}, "h", "\shortmid"{marking}, from=1-1, to=1-2]
    \arrow["f"', from=1-1, to=2-1]
    \arrow["g", from=1-2, to=2-2]
    \arrow[""{name=1, anchor=center, inner sep=0}, "k"', "\shortmid"{marking}, from=2-1, to=2-2]
    \arrow["\alpha"{description}, draw=none, from=0, to=1]
  \end{tikzcd}
  \qquad\leftrightsquigarrow\qquad
  % https://q.uiver.app/#q=WzAsNCxbMSwwLCJ4Il0sWzIsMSwieSJdLFswLDEsInciXSxbMSwyLCJ6Il0sWzAsMiwiZiIsMix7ImN1cnZlIjoxfV0sWzEsMywiZyIsMCx7ImN1cnZlIjotMX1dLFswLDEsImgiLDAseyJjdXJ2ZSI6LTF9XSxbMiwzLCJrIiwyLHsiY3VydmUiOjF9XSxbMiwxLCJcXGFscGhhIiwwLHsic2hvcnRlbiI6eyJzb3VyY2UiOjIwLCJ0YXJnZXQiOjIwfSwibGV2ZWwiOjJ9XV0=
  \begin{tikzcd}[sep=small]
    & x \\
    w && y \\
    & z
    \arrow["f"', curve={height=6pt}, from=1-2, to=2-1]
    \arrow["h", curve={height=-6pt}, from=1-2, to=2-3]
    \arrow["\alpha", shorten <=10pt, shorten >=10pt, Rightarrow, from=2-1, to=2-3]
    \arrow["k"', curve={height=6pt}, from=2-1, to=3-2]
    \arrow["g", curve={height=-6pt}, from=2-3, to=3-2]
  \end{tikzcd}
\end{equation*}
the 2-cells in $\twocat{C}$ of the form on the right
\cite[\S{3.1.4}]{grandis2019}. The \define{double category of adjunctions} in
$\twocat{C}$, denoted $\Adj(\twocat{C})$, again has the same objects and arrows
as $\twocat{C}$. Its proarrows $x \proto y$ are adjunctions $(f,g,\eta,\varepsilon): x \leftrightarrows y$
in $\twocat{C}$, and its cells are mate pairs of 2-cells in $\twocat{C}$
\cite[\S{3.1.5}]{grandis2019}.

\begin{theorem}[Twisted functors on equipments] \label{thm:twisted-functors-on-equipments}
  Let $F: \dbl{D} \twistto \dbl{E}$ be a twisted normal lax functor.
  \begin{enumerate}[(i)]
    \item When $\dbl{D}$ has companions, there is a twisted pseudo functor
      $\dbl{F}: \dbl{D} \twistto \Sq(\UTwoCat(\dbl{E}))$ that agrees with $F$ on
      objects and proarrows and acts on arrows as $\dbl{F}(f) \coloneqq F(f_!)$.
    \item When $\dbl{D}$ has conjoints, there is a twisted pseudo functor
      $\dbl{F}: \dbl{D}^\co \twistto \Sq(\UTwoCat(\dbl{E}^\co))$ that agrees
      with $F$ on objects and proarrows and acts on arrows as
      $\dbl{F}(f) \coloneqq F(f^*)$.
    \item When $\dbl{D}$ is an equipment, there is a twisted pseudo functor
      $\dbl{F}: \dbl{D} \twistto \Adj(\UTwoCat(\dbl{E}))$ that agrees with $F$
      on objects and proarrows and acts on arrows as
      $\dbl{F}(f) \coloneqq (F(f_!) \dashv F(f^*))$.
  \end{enumerate}
\end{theorem}
\begin{proof}
  To prove part (i), we use the construction given in
  \cref{prop:pseudo-functors-as-twisted-functors}. Suppose that $\dbl{D}$ has
  companions, and make a choice of them. First, we define
  $\dbl{F}: \dbl{D} \twistto \Sq(\UTwoCat(\dbl{E}))$ on objects and proarrows by
  $\dbl{F}(x) \coloneqq F(x)$ and $\dbl{F}(m) \coloneqq F(m)$; the
  loose-to-tight comparison cells are also borrowed from $F$. All this makes
  sense because the double categories $\dbl{E}$ and $\Sq(\UTwoCat(\dbl{E}))$
  have the same objects, arrows, and loosely globular cells. Next, on arrows,
  define $\dbl{F}(f) \coloneqq F(f_!)$, and on cells declare $\dbl{F}(\alpha)$ to be
  the composite given by the sliding calculus for companions
  \cite[\S{A.1}]{patterson2024}:
    \begin{equation*}
      % https://q.uiver.app/#q=WzAsNCxbMCwwLCJGeCJdLFsxLDAsIkZ3Il0sWzAsMSwiRnkiXSxbMSwxLCJGeiJdLFswLDIsIkZtIiwyXSxbMCwxLCJGZiIsMCx7InN0eWxlIjp7ImJvZHkiOnsibmFtZSI6ImJhcnJlZCJ9fX1dLFsxLDMsIkZuIl0sWzIsMywiRmciLDIseyJzdHlsZSI6eyJib2R5Ijp7Im5hbWUiOiJiYXJyZWQifX19XSxbNSw3LCJGXFxhbHBoYSIsMSx7InNob3J0ZW4iOnsic291cmNlIjoyMCwidGFyZ2V0IjoyMH0sInN0eWxlIjp7ImJvZHkiOnsibmFtZSI6Im5vbmUifSwiaGVhZCI6eyJuYW1lIjoibm9uZSJ9fX1dXQ==
      \begin{tikzcd}
        Fx & Fw \\
        Fy & Fz
        \arrow[""{name=0, anchor=center, inner sep=0}, "Ff", "\shortmid"{marking}, from=1-1, to=1-2]
        \arrow["Fm"', from=1-1, to=2-1]
        \arrow["Fn", from=1-2, to=2-2]
        \arrow[""{name=1, anchor=center, inner sep=0}, "Fg"', "\shortmid"{marking}, from=2-1, to=2-2]
        \arrow["{F\alpha}"{description}, draw=none, from=0, to=1]
      \end{tikzcd}
      \qquad\leftrightsquigarrow\qquad
      % https://q.uiver.app/#q=WzAsOCxbMCwzLCJGeiJdLFsxLDMsIkZ6Il0sWzAsMCwiRngiXSxbMSwwLCJGdyJdLFswLDEsIkZ5Il0sWzEsMSwiRnoiXSxbMCwyLCJGeiJdLFsxLDIsIkZ6Il0sWzAsMSwiIiwyLHsibGV2ZWwiOjIsInN0eWxlIjp7ImJvZHkiOnsibmFtZSI6ImJhcnJlZCJ9LCJoZWFkIjp7Im5hbWUiOiJub25lIn19fV0sWzIsMywiIiwxLHsibGV2ZWwiOjIsInN0eWxlIjp7ImJvZHkiOnsibmFtZSI6ImJhcnJlZCJ9LCJoZWFkIjp7Im5hbWUiOiJub25lIn19fV0sWzQsNSwiRmYiLDAseyJzdHlsZSI6eyJib2R5Ijp7Im5hbWUiOiJiYXJyZWQifX19XSxbMiw0LCIxIiwyXSxbMyw1LCJGZl8hIl0sWzQsNiwiRm0iLDJdLFs1LDcsIkZuIl0sWzYsNywiRmciLDIseyJzdHlsZSI6eyJib2R5Ijp7Im5hbWUiOiJiYXJyZWQifX19XSxbNiwwLCJGZ18hIiwyXSxbNywxLCIxIl0sWzE1LDgsIkZcXGVwc2lsb24iLDEseyJsYWJlbF9wb3NpdGlvbiI6NjAsInNob3J0ZW4iOnsic291cmNlIjoyMCwidGFyZ2V0IjoyMH0sInN0eWxlIjp7ImJvZHkiOnsibmFtZSI6Im5vbmUifSwiaGVhZCI6eyJuYW1lIjoibm9uZSJ9fX1dLFsxMCwxNSwiRlxcYWxwaGEiLDEseyJzaG9ydGVuIjp7InNvdXJjZSI6MjAsInRhcmdldCI6MjB9LCJzdHlsZSI6eyJib2R5Ijp7Im5hbWUiOiJub25lIn0sImhlYWQiOnsibmFtZSI6Im5vbmUifX19XSxbOSwxMCwiRlxcZXRhIiwxLHsibGFiZWxfcG9zaXRpb24iOjQwLCJzaG9ydGVuIjp7InNvdXJjZSI6MjAsInRhcmdldCI6MjB9LCJzdHlsZSI6eyJib2R5Ijp7Im5hbWUiOiJub25lIn0sImhlYWQiOnsibmFtZSI6Im5vbmUifX19XV0=
      \begin{tikzcd}
        Fx & Fx \\
        Fx & Fw \\
        Fy & Fz \\
        Fz & Fz
        \arrow[""{name=0, anchor=center, inner sep=0}, "\shortmid"{marking}, equals, from=1-1, to=1-2]
        \arrow["1"', from=1-1, to=2-1]
        \arrow["{Ff_!}", from=1-2, to=2-2]
        \arrow[""{name=1, anchor=center, inner sep=0}, "Ff"{inner sep=.8ex}, "\shortmid"{marking}, from=2-1, to=2-2]
        \arrow["Fm"', from=2-1, to=3-1]
        \arrow["Fn", from=2-2, to=3-2]
        \arrow[""{name=2, anchor=center, inner sep=0}, "Fg"'{inner sep=.8ex}, "\shortmid"{marking}, from=3-1, to=3-2]
        \arrow["{Fg_!}"', from=3-1, to=4-1]
        \arrow["1", from=3-2, to=4-2]
        \arrow[""{name=3, anchor=center, inner sep=0}, "\shortmid"{marking}, equals, from=4-1, to=4-2]
        \arrow["{F\eta}"{description, pos=0.4}, draw=none, from=0, to=1]
        \arrow["{F\alpha}"{description}, draw=none, from=1, to=2]
        \arrow["{F\epsilon}"{description, pos=0.6}, draw=none, from=2, to=3]
      \end{tikzcd}
      \quad\eqqcolon\quad
      % https://q.uiver.app/#q=WzAsNCxbMCwwLCJcXGRibHtGfXgiXSxbMSwwLCJcXGRibHtGfXciXSxbMCwxLCJcXGRibHtGfXkiXSxbMSwxLCJcXGRibHtGfXoiXSxbMCwyLCJcXGRibHtGfW0iLDJdLFsxLDMsIlxcZGJse0Z9biJdLFsyLDMsIlxcZGJse0Z9ZyIsMix7InN0eWxlIjp7ImJvZHkiOnsibmFtZSI6ImJhcnJlZCJ9fX1dLFswLDEsIlxcZGJse0Z9ZiIsMCx7InN0eWxlIjp7ImJvZHkiOnsibmFtZSI6ImJhcnJlZCJ9fX1dLFs3LDYsIlxcZGJse0Z9XFxhbHBoYSIsMSx7InNob3J0ZW4iOnsic291cmNlIjoyMCwidGFyZ2V0IjoyMH0sInN0eWxlIjp7ImJvZHkiOnsibmFtZSI6Im5vbmUifSwiaGVhZCI6eyJuYW1lIjoibm9uZSJ9fX1dXQ==
      \begin{tikzcd}
        {\dbl{F}x} & {\dbl{F}w} \\
        {\dbl{F}y} & {\dbl{F}z}
        \arrow[""{name=0, anchor=center, inner sep=0}, "{\dbl{F}f}", "\shortmid"{marking}, from=1-1, to=1-2]
        \arrow["{\dbl{F}m}"', from=1-1, to=2-1]
        \arrow["{\dbl{F}n}", from=1-2, to=2-2]
        \arrow[""{name=1, anchor=center, inner sep=0}, "{\dbl{F}g}"', "\shortmid"{marking}, from=2-1, to=2-2]
        \arrow["{\dbl{F}\alpha}"{description}, draw=none, from=0, to=1]
      \end{tikzcd}.
    \end{equation*}
  We have been a bit cavalier here in that the top and bottom cells, namely,
  the images of the companion unit and counit, need to be suitably padded by
  identity comparison cells. Since $F$ is assumed to be normal these cells are all
  invertible.

  In any case, $\dbl{F}$ is a well-defined twisted functor. For, by
  construction, its comparison cells in both directions satisfy the
  associativity and unitality axioms. The naturality axioms for the comparison
  cells follow from the functorality properties of sliding. Specifically,
  loose-to-tight naturality holds by \cite[\mbox{Lemma A.9}]{patterson2024} and
  tight-to-loose naturality holds by \cite[\mbox{Lemma A.10}]{patterson2024}.

  Part (ii) is dual to part (i). For ease of reference, we record that when
  $\dbl{D}$ has conjoints, the twisted functor
  $\dbl{F}: \dbl{D}^\co \twistto \Sq(\UTwoCat(\dbl{E}^\co))$ is defined on arrows by
  $\dbl{F} f \coloneqq F(f^*)$, hence satisfies the equation
  $(\dbl{F} f)^* = Ff$ in $\dbl{E}$ by
  \cref{prop:twisting-companions-conjoints}. It is defined on cells by
  \begin{equation*}
    % https://q.uiver.app/#q=WzAsNCxbMCwwLCJGeCJdLFsxLDAsIkZ3Il0sWzAsMSwiRnkiXSxbMSwxLCJGeiJdLFswLDIsIkZtIiwyXSxbMCwxLCJGZiIsMCx7InN0eWxlIjp7ImJvZHkiOnsibmFtZSI6ImJhcnJlZCJ9fX1dLFsxLDMsIkZuIl0sWzIsMywiRmciLDIseyJzdHlsZSI6eyJib2R5Ijp7Im5hbWUiOiJiYXJyZWQifX19XSxbNSw3LCJGXFxhbHBoYSIsMSx7InNob3J0ZW4iOnsic291cmNlIjoyMCwidGFyZ2V0IjoyMH0sInN0eWxlIjp7ImJvZHkiOnsibmFtZSI6Im5vbmUifSwiaGVhZCI6eyJuYW1lIjoibm9uZSJ9fX1dXQ==
    \begin{tikzcd}
      Fx & Fw \\
      Fy & Fz
      \arrow[""{name=0, anchor=center, inner sep=0}, "Ff", "\shortmid"{marking}, from=1-1, to=1-2]
      \arrow["Fm"', from=1-1, to=2-1]
      \arrow["Fn", from=1-2, to=2-2]
      \arrow[""{name=1, anchor=center, inner sep=0}, "Fg"', "\shortmid"{marking}, from=2-1, to=2-2]
      \arrow["{F\alpha}"{description}, draw=none, from=0, to=1]
    \end{tikzcd}
    \qquad\leftrightsquigarrow\qquad
    % https://q.uiver.app/#q=WzAsNixbMCwxLCJGeCJdLFsxLDAsIkZ3Il0sWzAsMiwiRnkiXSxbMSwxLCJGeiJdLFsxLDIsIkZ5Il0sWzAsMCwiRnciXSxbMCwyLCJGbSIsMl0sWzEsMywiRm4iXSxbMyw0LCJGKGdeKikiXSxbNSwwLCJGKGZeKikiLDJdLFs1LDEsIiIsMCx7ImxldmVsIjoyLCJzdHlsZSI6eyJib2R5Ijp7Im5hbWUiOiJiYXJyZWQifSwiaGVhZCI6eyJuYW1lIjoibm9uZSJ9fX1dLFsyLDQsIiIsMix7ImxldmVsIjoyLCJzdHlsZSI6eyJib2R5Ijp7Im5hbWUiOiJiYXJyZWQifSwiaGVhZCI6eyJuYW1lIjoibm9uZSJ9fX1dLFsxMCwxMSwiXFxEb3duYXJyb3ciLDEseyJzaG9ydGVuIjp7InNvdXJjZSI6MjAsInRhcmdldCI6MjB9LCJzdHlsZSI6eyJib2R5Ijp7Im5hbWUiOiJub25lIn0sImhlYWQiOnsibmFtZSI6Im5vbmUifX19XV0=
    \begin{tikzcd}[row sep=scriptsize]
      Fw & Fw \\
      Fx & Fz \\
      Fy & Fy
      \arrow[""{name=0, anchor=center, inner sep=0}, "\shortmid"{marking}, Rightarrow, no head, from=1-1, to=1-2]
      \arrow["{F(f^*)}"', from=1-1, to=2-1]
      \arrow["Fn", from=1-2, to=2-2]
      \arrow["Fm"', from=2-1, to=3-1]
      \arrow["{F(g^*)}", from=2-2, to=3-2]
      \arrow[""{name=1, anchor=center, inner sep=0}, "\shortmid"{marking}, Rightarrow, no head, from=3-1, to=3-2]
      \arrow["\Downarrow"{description}, draw=none, from=0, to=1]
    \end{tikzcd}
    \quad\eqqcolon\quad
    % https://q.uiver.app/#q=WzAsNCxbMCwwLCJcXGRibHtGfXciXSxbMSwwLCJcXGRibHtGfXgiXSxbMCwxLCJcXGRibHtGfXoiXSxbMSwxLCJcXGRibHtGfXkiXSxbMCwyLCJcXGRibHtGfW4iLDJdLFsxLDMsIlxcZGJse0Z9bSJdLFsyLDMsIlxcZGJse0Z9ZyIsMix7InN0eWxlIjp7ImJvZHkiOnsibmFtZSI6ImJhcnJlZCJ9fX1dLFswLDEsIlxcZGJse0Z9ZiIsMCx7InN0eWxlIjp7ImJvZHkiOnsibmFtZSI6ImJhcnJlZCJ9fX1dLFs3LDYsIlxcZGJse0Z9XFxhbHBoYSIsMSx7InNob3J0ZW4iOnsic291cmNlIjoyMCwidGFyZ2V0IjoyMH0sInN0eWxlIjp7ImJvZHkiOnsibmFtZSI6Im5vbmUifSwiaGVhZCI6eyJuYW1lIjoibm9uZSJ9fX1dXQ==
    \begin{tikzcd}
      {\dbl{F}w} & {\dbl{F}x} \\
      {\dbl{F}z} & {\dbl{F}y}
      \arrow[""{name=0, anchor=center, inner sep=0}, "{\dbl{F}f}", "\shortmid"{marking}, from=1-1, to=1-2]
      \arrow["{\dbl{F}n}"', from=1-1, to=2-1]
      \arrow["{\dbl{F}m}", from=1-2, to=2-2]
      \arrow[""{name=1, anchor=center, inner sep=0}, "{\dbl{F}g}"', "\shortmid"{marking}, from=2-1, to=2-2]
      \arrow["{\dbl{F}\alpha}"{description}, draw=none, from=0, to=1]
    \end{tikzcd},
  \end{equation*}
  hence satisfies the equation $(\dbl{F} \alpha)^* = F\alpha$ between cells in
  $\dbl{E}$.

  Finally, to prove part (iii), suppose that $\dbl{D}$ is an equipment. Using
  parts (i) and (ii), we can construct twisted functors
  $\dbl{F}_\bullet: \dbl{D} \twistto \Sq(\UTwoCat(\dbl{E}))$ and
  $\dbl{F}^\bullet: \dbl{D}^\co \twistto \Sq(\UTwoCat(\dbl{E}^\co))$ that agree with
  each other (and with $F$) on objects and proarrows. By
  \cref{cor:twisting-companion-conjoint-adjunctions}, for any arrow $f$ in
  $\dbl{D}$, we see that $\dbl{F}_\bullet(f) = F(f_!)$ and $\dbl{F}^\bullet(f) = F(f^*)$ are
  adjoint as arrows in $\dbl{E}$. Moreover, for any cell $\alpha$ in $\dbl{D}$, the
  cells $\dbl{F}_\bullet(\alpha)$ and $\dbl{F}^\bullet(\alpha)$ viewed in $\dbl{E}$ are none other
  than the mates of each other. The twisted functors $\dbl{F}_\bullet$ and $\dbl{F}^\bullet$
  thus pair into a twisted normal lax functor
  $\dbl{F}: \dbl{D} \twistto \Adj(\UTwoCat(\dbl{E}))$.
\end{proof}

Applying the theorem to twisted representables, we deduce:

\begin{corollary}[Twisted representables on equipments]
  When $\dbl{D}$ is an equipment, any twisted representable functor
  $\dbl{D}(a,-): \dbl{D} \twistto \Prof$ defines a twisted pseudo functor
  $\dbl{D} \twistto \Adj$ and, in particular, a bi-indexed category over
  $\dbl{D}_0$.
\end{corollary}

In several examples from \cref{subsection:twisted-representable-functors}, such
as the canonical self-indexing of a finitely complete category
(\cref{ex:self-indexing}) and the subobjects of a regular category
(\cref{ex:subobjects}), we saw that a twisted representable furnishes a
bi-indexed category. The corollary explains why this happens: the original
double category, here that of spans or relations, is an equipment.

Another way of repackaging the data of a twisted functor, when the domain
$\dbl{D}$ has conjoints, is to ``untwist'' a twisted functor
$\dbl{D} \twistto \dbl{E}$ into a double \emph{pseudo}functor
$\dbl{D} \to \Sq(\UTwoCat(\dbl{E}))^{\op}$. Not to be confused with a pseudo
double functor, a \emph{double pseudofunctor} is a kind of functor between
double categories that is pseudo in both directions \cite[Definition
6.1]{shulman2011}.

\begin{theorem}[Untwisting a twisted functor]\label{thm:untwisting}
  Let $F : \dbl{D} \twistto \dbl{E}$ be a twisted normal lax functor. When
  $\dbl{D}$ has conjoints, there is a double pseudofunctor
  $\dbl{F} : \dbl{D} \to \Sq(\UTwoCat(\dbl{E}))^{\op}$ that agrees with $F$ on
  objects and loose arrows and acts on tight arrows by $\dbl{F}(f) := F(f^{\ast})$.
\end{theorem}
\begin{proof}
  We define $\dbl{F}$ on objects, by $\dbl{F}(x) := Fx$; on proarrows, by
  $\dbl{F}(m) := Fm$ for $m : x \proto y$; and on arrows, by
  $\dbl{F}(f) := F(f^{\ast}) : \dbl{F}y \to \dbl{F} x$ for $f : x \to y$. On cells, we
  use the sliding calculus for conjoints \cite[\S{A.1}]{patterson2024}:
  \begin{equation*}
    % https://q.uiver.app/#q=WzAsNCxbMCwwLCJGeCJdLFsxLDAsIkZ3Il0sWzAsMSwiRnkiXSxbMSwxLCJGeiJdLFswLDIsIkZtIiwyXSxbMCwxLCJGZiIsMCx7InN0eWxlIjp7ImJvZHkiOnsibmFtZSI6ImJhcnJlZCJ9fX1dLFsxLDMsIkZuIl0sWzIsMywiRmciLDIseyJzdHlsZSI6eyJib2R5Ijp7Im5hbWUiOiJiYXJyZWQifX19XSxbNSw3LCJGXFxhbHBoYSIsMSx7InNob3J0ZW4iOnsic291cmNlIjoyMCwidGFyZ2V0IjoyMH0sInN0eWxlIjp7ImJvZHkiOnsibmFtZSI6Im5vbmUifSwiaGVhZCI6eyJuYW1lIjoibm9uZSJ9fX1dXQ==
    \begin{tikzcd}
      Fx & Fw \\
      Fy & Fz
      \arrow[""{name=0, anchor=center, inner sep=0}, "Ff", "\shortmid"{marking}, from=1-1, to=1-2]
      \arrow["Fm"', from=1-1, to=2-1]
      \arrow["Fn", from=1-2, to=2-2]
      \arrow[""{name=1, anchor=center, inner sep=0}, "Fg"', "\shortmid"{marking}, from=2-1, to=2-2]
      \arrow["{F\alpha}"{description}, draw=none, from=0, to=1]
    \end{tikzcd}
    \qquad\leftrightsquigarrow\qquad
    % https://q.uiver.app/#q=WzAsOCxbMCwzLCJGeSJdLFsxLDMsIkZ5Il0sWzAsMCwiRnciXSxbMSwwLCJGdyJdLFswLDEsIkZ4Il0sWzEsMSwiRnciXSxbMCwyLCJGeSJdLFsxLDIsIkZ6Il0sWzAsMSwiIiwyLHsibGV2ZWwiOjIsInN0eWxlIjp7ImJvZHkiOnsibmFtZSI6ImJhcnJlZCJ9LCJoZWFkIjp7Im5hbWUiOiJub25lIn19fV0sWzIsMywiIiwxLHsibGV2ZWwiOjIsInN0eWxlIjp7ImJvZHkiOnsibmFtZSI6ImJhcnJlZCJ9LCJoZWFkIjp7Im5hbWUiOiJub25lIn19fV0sWzQsNSwiRmYiLDAseyJzdHlsZSI6eyJib2R5Ijp7Im5hbWUiOiJiYXJyZWQifX19XSxbMiw0LCJGZl57XFxhc3R9IiwyXSxbMyw1LCIxIl0sWzQsNiwiRm0iLDJdLFs1LDcsIkZuIl0sWzYsNywiRmciLDIseyJzdHlsZSI6eyJib2R5Ijp7Im5hbWUiOiJiYXJyZWQifX19XSxbNiwwLCIxIiwyXSxbNywxLCJGZ157XFxhc3R9Il0sWzE1LDgsIkZcXGVwc2lsb24iLDEseyJsYWJlbF9wb3NpdGlvbiI6NjAsInNob3J0ZW4iOnsic291cmNlIjoyMCwidGFyZ2V0IjoyMH0sInN0eWxlIjp7ImJvZHkiOnsibmFtZSI6Im5vbmUifSwiaGVhZCI6eyJuYW1lIjoibm9uZSJ9fX1dLFsxMCwxNSwiRlxcYWxwaGEiLDEseyJzaG9ydGVuIjp7InNvdXJjZSI6MjAsInRhcmdldCI6MjB9LCJzdHlsZSI6eyJib2R5Ijp7Im5hbWUiOiJub25lIn0sImhlYWQiOnsibmFtZSI6Im5vbmUifX19XSxbOSwxMCwiRlxcZXRhIiwxLHsibGFiZWxfcG9zaXRpb24iOjQwLCJzaG9ydGVuIjp7InNvdXJjZSI6MjAsInRhcmdldCI6MjB9LCJzdHlsZSI6eyJib2R5Ijp7Im5hbWUiOiJub25lIn0sImhlYWQiOnsibmFtZSI6Im5vbmUifX19XV0=
    \begin{tikzcd}
      Fw & Fw \\
      Fx & Fw \\
      Fy & Fz \\
      Fy & Fy
      \arrow[""{name=0, anchor=center, inner sep=0}, "\shortmid"{marking}, equals, from=1-1, to=1-2]
      \arrow["{Ff^{\ast}}"', from=1-1, to=2-1]
      \arrow["1", from=1-2, to=2-2]
      \arrow[""{name=1, anchor=center, inner sep=0}, "Ff"{inner sep=.8ex}, "\shortmid"{marking}, from=2-1, to=2-2]
      \arrow["Fm"', from=2-1, to=3-1]
      \arrow["Fn", from=2-2, to=3-2]
      \arrow[""{name=2, anchor=center, inner sep=0}, "Fg"'{inner sep=.8ex}, "\shortmid"{marking}, from=3-1, to=3-2]
      \arrow["1"', from=3-1, to=4-1]
      \arrow["{Fg^{\ast}}", from=3-2, to=4-2]
      \arrow[""{name=3, anchor=center, inner sep=0}, "\shortmid"{marking}, equals, from=4-1, to=4-2]
      \arrow["{F\eta}"{description, pos=0.4}, draw=none, from=0, to=1]
      \arrow["{F\alpha}"{description}, draw=none, from=1, to=2]
      \arrow["{F\epsilon}"{description, pos=0.6}, draw=none, from=2, to=3]
    \end{tikzcd}
    \quad\eqqcolon\quad
    % https://q.uiver.app/#q=WzAsNCxbMCwwLCJcXGRibHtGfXgiXSxbMCwxLCJcXGRibHtGfXciXSxbMSwwLCJcXGRibHtGfXkiXSxbMSwxLCJcXGRibHtGfXoiXSxbMCwyLCJcXGRibHtGfW0iXSxbMSwzLCJcXGRibHtGfW4iLDJdLFsyLDMsIlxcZGJse0Z9ZyIsMCx7InN0eWxlIjp7InRhaWwiOnsibmFtZSI6ImFycm93aGVhZCJ9LCJoZWFkIjp7Im5hbWUiOiJub25lIn19fV0sWzAsMSwiXFxkYmx7Rn1mIiwyLHsic3R5bGUiOnsidGFpbCI6eyJuYW1lIjoiYXJyb3doZWFkIn0sImhlYWQiOnsibmFtZSI6Im5vbmUifX19XSxbNyw2LCJcXGRibHtGfVxcYWxwaGEiLDEseyJzaG9ydGVuIjp7InNvdXJjZSI6MjAsInRhcmdldCI6MjB9LCJzdHlsZSI6eyJib2R5Ijp7Im5hbWUiOiJub25lIn0sImhlYWQiOnsibmFtZSI6Im5vbmUifX19XV0=
    \begin{tikzcd}
      {\dbl{F}x} & {\dbl{F}y} \\
      {\dbl{F}w} & {\dbl{F}z}
      \arrow["{\dbl{F}m}", from=1-1, to=1-2]
      \arrow[""{name=0, anchor=center, inner sep=0}, "{\dbl{F}f}"', tail reversed, no head, from=1-1, to=2-1]
      \arrow[""{name=1, anchor=center, inner sep=0}, "{\dbl{F}g}", tail reversed, no head, from=1-2, to=2-2]
      \arrow["{\dbl{F}n}"', from=2-1, to=2-2]
      \arrow["{\dbl{F}\alpha}"{description}, draw=none, from=0, to=1]
    \end{tikzcd}.
  \end{equation*}
  As in \cref{thm:twisted-functors-on-equipments}, we need to pad out the images
  of the conjoint unit and counit with identity cells, which are isomorphisms by
  assumption. The comparison cells of $\dbl{F}$ are the loose-to-tight
  comparisons of $F$, suitably arranged; note that these appear both for tights
  and for looses, and therefore $\dbl{F}$ must be a double \emph{pseudo}functor.
  Nevertheless, the coherences required of this double pseudofunctor are
  essentially those of the twisted double functor, suitably ``slid'' using
  conjoints.
\end{proof}

Among twisted copresheaves on equipments we now single out those on span double
categories, for which we will establish a tighter correspondence with bi-indexed
categories. Specifically, we will show that twisted copresheaves
$\Span(\cat{C}) \twistto \Prof$ correspond exactly with bi-indexed categories
$\cat{C} \to \AdjTwo$ satisfying the Beck-Chevalley condition.

\begin{definition}[Beck-Chevalley]
  Let $\cat{C}$ be a category with pullbacks. A bi-indexed category over
  $\cat{C}$, sending a morphism $f$ in $\cat{C}$ to an adjunction
  $f_! \dashv f^*$, satisfies the \define{Beck-Chevalley condition} if for each
  pullback square in $\cat{C}$ as on the left below, the natural transformation
  arising as the mate of the image of the square under
  $(-)^*: \cat{C}^\op \to \Cat$ is an isomorphism:
  \begin{equation*}
    % https://q.uiver.app/#q=WzAsNCxbMCwxLCJ4Il0sWzEsMiwieiJdLFsyLDEsInkiXSxbMSwwLCJ3Il0sWzAsMSwiZiIsMl0sWzIsMSwiZyJdLFszLDAsImsiLDJdLFszLDIsImgiXSxbMywxLCIiLDEseyJzdHlsZSI6eyJuYW1lIjoiY29ybmVyIn19XV0=
    \begin{tikzcd}[sep=scriptsize]
      & w & \\
      x && y \\
      & z
      \arrow["k"', from=1-2, to=2-1]
      \arrow["h", from=1-2, to=2-3]
      \arrow["\lrcorner"{anchor=center, pos=0.125, rotate=-45}, draw=none, from=1-2, to=3-2]
      \arrow["f"', from=2-1, to=3-2]
      \arrow["g", from=2-3, to=3-2]
    \end{tikzcd}
    \qquad\leadsto\qquad
    k^* \cdot h_! \cong f_! \cdot g^*.
    \qedhere
  \end{equation*}
\end{definition}

Our argument will make use of a correspondence due to Siqueira
\cite{siqueira2025}, which is in essence a double (1-)categorical version of
Barwick's \emph{unfurling} construction \cite{barwick2017}. We begin by
recalling Siqueira's correspondence as encoded in Propositions 3.1 and 3.6 of
\cite{siqueira2025}.

\begin{theorem}[\cite{siqueira2025}] \label{thm:siqueira.correspondence}
  Let $\cat{C}$ be a category with pullbacks. Bi-indexed categories over
  $\cat{C}$ satisfying the Beck-Chevalley condition correspond to double
  pseudofunctors $\Span(\cat{C})^\op \to \Sq(\Cat)$ (see Definition 2.9 of
  \cite{siqueira2025}).
\end{theorem}
\begin{proof}
  We note that Proposition 3.1 of \cite{siqueira2025} makes no use of the
  monoidal structure and provides the construction of the double pseudofunctor
  $\Span(\cat{C})^\op \to \Sq(\Cat)$ corresponding to a bi-indexed category
  $E : \cat{C} \to \AdjTwo$ satisfying Beck-Chevalley. The action on tight
  morphisms takes $f$ to $f^*$, while the action on spans sends
  $x \xleftarrow{\ell} s \xto{r} y$ to $\ell^{\ast} \cdot r_{!}$. The
  pseudo-functoriality of composition encodes the Beck-Chevalley condition. The
  converse is provided by Proposition 3.6.
\end{proof}

We will now adapt this correspondence to twisted copresheaves. First, we
construct a straightforward twisted copresheaf on $\Sq(\Cat)$, including it into
$\Prof$.

\begin{construction}
  Define the twisted inclusion $\iota : \Sq(\Cat)^\op \twistto \Prof$ that sends
  every category to itself, every loose morphism to itself (considered as a
  tight morphism of $\Prof$), and every tight morphism $F : \cat{X} \leftarrow \cat{Y}$
  to its conjoint $F^{\ast} = \cat{X}(-,F(=)): \cat{X} \proto \cat{Y}$.
\end{construction}

We can now prove our correspondence with twisted copresheaves.

\begin{theorem}[Twisted copresheaves on spans]
  \label{thm:correspondence-bi-indexed-twisted-copresheaf}
  Let $\cat{C}$ be a category with pullbacks. Given a bi-indexed category
  $E : \cat{C} \to \AdjTwo$ satisfying the Beck-Chevalley condition, we may
  associate to it the twisted copresheaf $\Span(\cat{C}) \twistto \Prof$ which
  acts as $E$ on objects, sends a span $x \xfrom{\ell} s \xto{r} y$ to the functor
  $\ell^{\ast} \cdot r_!$, and sends a map $f : x \to y$ to the profunctor
  $E(x)(-, f^{\ast}(=))$.

  This assignment yields a correspondence between bi-indexed categories over
  $\cat{C}$ satisfying Beck-Chevalley and twisted copresheaves
  $\Span(\cat{C}) \twistto \Prof$.
\end{theorem}
\begin{proof}
  Given a bi-indexed category $E : \cat{C} \to \AdjTwo$ satisfying
  Beck-Chevalley, consider the corresponding double pseudofunctor
  $E : \Span(\cat{C})^\op \to \Sq(\Cat)$ by \cref{thm:siqueira.correspondence}. We
  may then pre-compose $\iota$ by $E^\op$ to get a twisted copresheaf
  $\iota \circ E^\op : \Span(\cat{C}) \twistto \Prof$ (adapting
  \cref{construction:twisted-pre-composite} along the lines of
  \cref{rmk:pseudo-pre-composition}).

  Conversely, a twisted copresheaf $E : \Span(\cat{C}) \twistto \Prof$ gives
  rise to a double pseudofunctor $\dbl{E} : \Span(\cat{C}) \to \Sq(\Cat)^\op$ by
  \cref{thm:untwisting}; this then corresponds to $E : \cat{C} \to \AdjTwo$ by
  \cref{thm:siqueira.correspondence}.
\end{proof}

\section{Loosely discrete double opfibrations}

The central concept of a loosely discrete opfibration of double categories is
introduced in this section as an appropriate adaptation of the notion of an
internal discrete opfibration suited for the 2-category of categories, functors
and transformations. The important definition of a cleavage for such a loosely
discrete opfibration is the topic of \cref{subsection:cleavages}. The algebraic
axioms and the lifting properties of such cleavages are presented and studied in
detail. Lastly, the elements construction associated to a profunctor-valued
normal twisted functor is given in \cref{subsection:elements} and is shown to
give rise to such a loosely discrete opfibration.

\subsection{Definition and examples}

Recall that a functor $P\colon \cat{E}\to\cat{B}$ is a \define{discrete opfibration} if
in the diagram
\begin{equation*}
    % https://q.uiver.app/#q=WzAsNSxbMSwxLCJcXGNhdHtFfV8wXFx0aW1lc197XFxjYXR7Qn1fMH1cXGNhdHtCfV8xIl0sWzIsMSwiXFxjYXR7Qn1fMSJdLFsyLDIsIlxcY2F0e0J9XzAiXSxbMSwyLCJcXGNhdHtFfV8wIl0sWzAsMCwiXFxjYXR7RX1fMSJdLFswLDEsIlxccGlfMiJdLFsxLDIsImRfMCJdLFswLDMsIlxccGlfMSIsMl0sWzMsMiwiUF8wIiwyXSxbNCwxLCJQXzEiLDAseyJjdXJ2ZSI6LTJ9XSxbNCwzLCJkXzAiLDIseyJjdXJ2ZSI6Mn1dLFs0LDAsIlxcbGFuZ2xlIGRfMCxQXzFcXHJhbmdsZSIsMSx7InN0eWxlIjp7ImJvZHkiOnsibmFtZSI6ImRhc2hlZCJ9fX1dLFswLDIsIiIsMix7InN0eWxlIjp7Im5hbWUiOiJjb3JuZXIifX1dXQ==
    \begin{tikzcd}
      {\cat{E}_1} \\
      & {\cat{E}_0\times_{\cat{B}_0}\cat{B}_1} & {\cat{B}_1} \\
      & {\cat{E}_0} & {\cat{B}_0}
      \arrow["{\langle d_0,P_1\rangle}"{description}, dashed, from=1-1, to=2-2]
      \arrow["{P_1}", curve={height=-12pt}, from=1-1, to=2-3]
      \arrow["{d_0}"', curve={height=12pt}, from=1-1, to=3-2]
      \arrow["{\pi_2}", from=2-2, to=2-3]
      \arrow["{\pi_1}"', from=2-2, to=3-2]
      \arrow["\lrcorner"{anchor=center, pos=0.125}, draw=none, from=2-2, to=3-3]
      \arrow["{d_0}", from=2-3, to=3-3]
      \arrow["{P_0}"', from=3-2, to=3-3]
    \end{tikzcd},
\end{equation*}
the canonical function given by the dashed arrow is a bijection, i.e., if the
square forming the outside of the diagram presents the pullback of the domain
function of $\cat{B}$ along $P_0$. This standard description of a discrete
opfibration makes sense in any category with enough finite limits. However,
naively instantiating it in the category of categories to define a
double-categorical notion of discrete opfibration would yield too strict a
definition. We would be asking that the induced dashed arrow be an isomorphism
of categories, rarely the appropriate notion of sameness for categories. In the
following definition, we weaken the isomorphism to an equivalence.

\begin{definition}[Loosely discrete opfibration] \label{def:loosely-discrete-opfibration}
  A double functor $P\colon \dbl{E}\to\dbl{B}$ is a \textbf{loosely discrete 
  double opfibration} if in the diagram
  \begin{equation*}
      % https://q.uiver.app/#q=WzAsNSxbMSwxLCJcXGRibHtFfV8wXFx0aW1lc197XFxkYmx7Qn1fMH1cXGRibHtCfV8xIl0sWzIsMSwiXFxkYmx7Qn1fMSJdLFsyLDIsIlxcZGJse0J9XzAiXSxbMSwyLCJcXGRibHtFfV8wIl0sWzAsMCwiXFxkYmx7RX1fMSJdLFswLDEsIlxccGlfMiJdLFsxLDIsIlxcc3JjIl0sWzAsMywiXFxwaV8xIiwyXSxbMywyLCJQXzAiLDJdLFs0LDEsIlBfMSIsMCx7ImN1cnZlIjotMn1dLFs0LDMsIlxcc3JjIiwyLHsiY3VydmUiOjJ9XSxbNCwwLCJcXGxhbmdsZVxcc3JjLFBfMVxccmFuZ2xlIiwxLHsic3R5bGUiOnsiYm9keSI6eyJuYW1lIjoiZGFzaGVkIn19fV0sWzAsMiwiIiwyLHsic3R5bGUiOnsibmFtZSI6ImNvcm5lciJ9fV1d
      \begin{tikzcd}
        {\dbl{E}_1} \\
        & {\dbl{E}_0\times_{\dbl{B}_0}\dbl{B}_1} & {\dbl{B}_1} \\
        & {\dbl{E}_0} & {\dbl{B}_0}
        \arrow["{\langle\src,P_1\rangle}"{description}, dashed, from=1-1, to=2-2]
        \arrow["{P_1}", curve={height=-12pt}, from=1-1, to=2-3]
        \arrow["\src"', curve={height=12pt}, from=1-1, to=3-2]
        \arrow["{\pi_2}", from=2-2, to=2-3]
        \arrow["{\pi_1}"', from=2-2, to=3-2]
        \arrow["\lrcorner"{anchor=center, pos=0.125}, draw=none, from=2-2, to=3-3]
        \arrow["\src", from=2-3, to=3-3]
        \arrow["{P_0}"', from=3-2, to=3-3]
      \end{tikzcd},
  \end{equation*}
  the canonical functor given by the dashed arrow is a (weak) equivalence of
  categories.
\end{definition}

In the definition, by a \emph{(weak) equivalence} of categories, we mean a
functor that is fully faithful and essentially surjective on objects.

\begin{remark}[Loose discreteness] \label{rmk:on-loose-discreteness}
  A loosely discrete opfibration is ``locally discrete in the loose direction.'' We
  will explain what this means in steps. First, a double category is
  \define{loosely discrete} if it is equivalent to the double category $\Tight(\cat{D})$
  on some category $\cat{D}$. Such a double category has only trivial
  proarrows and cells. The double category $\Tight(\cat{D})$ is of course
  loosely discrete in this sense. Next, by ``locally,'' we mean ``for each
    fiber.'' Given a double functor $P: \dbl{E} \to \dbl{B}$, the \define{fiber} over
  an object $x\in\dbl{B}$, denoted $\dbl{E}_x$, is the double category given by
  pulling back along the point determined by $x$:
  \begin{equation*}
      % https://q.uiver.app/#q=WzAsNCxbMCwwLCJcXGRibHtFfV94Il0sWzEsMCwiXFxkYmx7RX0iXSxbMSwxLCJcXGRibHtCfSJdLFswLDEsIlxcZGJsezF9Il0sWzAsMV0sWzEsMiwiUCJdLFswLDNdLFszLDIsIngiLDJdLFswLDIsIiIsMSx7InN0eWxlIjp7Im5hbWUiOiJjb3JuZXIifX1dXQ==
      \begin{tikzcd}
        {\dbl{E}_x} & {\dbl{E}} \\
        {\dbl{1}} & {\dbl{B}}
        \arrow[from=1-1, to=1-2]
        \arrow[from=1-1, to=2-1]
        \arrow["\lrcorner"{anchor=center, pos=0.125}, draw=none, from=1-1, to=2-2]
        \arrow["P", from=1-2, to=2-2]
        \arrow["x"', from=2-1, to=2-2]
      \end{tikzcd}.
  \end{equation*}
  So $\dbl{E}_x$ is the double category consisting of those objects, arrows,
  proarrows, and cells of $\dbl{E}$ mapped by $P$ onto the object $x$ in
  $\dbl{B}$, its associated identity arrow, identity proarrow, and identity
  cell. Thus, in saying that the double functor $P$ is \define{locally
    discrete in the loose direction}, we mean that each fiber $\dbl{E}_x$ is
  loosely discrete as a double category. This is now proved more precisely in
  the next two results.
\end{remark}

\begin{lemma}
  A double functor $!: \dbl{E}\to\dbl{1}$ is a loosely discrete opfibration if, and
  only if, $\dbl{E}$ is loosely discrete as a double category.
\end{lemma}
\begin{proof}
  In this case, the diagram in the definition takes the form
    \begin{equation*}
      % https://q.uiver.app/#q=WzAsNSxbMSwxLCJcXGRibHtFfV8wXFx0aW1lc197XFxkYmx7MX1fMH1cXGRibHsxfV8xIl0sWzIsMSwiXFxkYmx7MX1fMSJdLFsyLDIsIlxcZGJsezF9XzAiXSxbMSwyLCJcXGRibHtFfV8wIl0sWzAsMCwiXFxkYmx7RX1fMSJdLFswLDEsIlxccGlfMiJdLFsxLDIsIlxcc3JjIl0sWzAsMywiXFxwaV8xIiwyXSxbMywyLCJQXzAiLDJdLFs0LDEsIlBfMSIsMCx7ImN1cnZlIjotMn1dLFs0LDMsIlxcc3JjIiwyLHsiY3VydmUiOjJ9XSxbNCwwLCJcXGxhbmdsZVxcc3JjLFBfMVxccmFuZ2xlIiwxLHsic3R5bGUiOnsiYm9keSI6eyJuYW1lIjoiZGFzaGVkIn19fV0sWzAsMiwiIiwyLHsic3R5bGUiOnsibmFtZSI6ImNvcm5lciJ9fV1d
      \begin{tikzcd}
        {\dbl{E}_1} \\
        & {\dbl{E}_0\times_{\dbl{1}_0}\dbl{1}_1} & {\dbl{1}_1} \\
        & {\dbl{E}_0} & {\dbl{1}_0}
        \arrow["{\langle\src,P_1\rangle}"{description}, dashed, from=1-1, to=2-2]
        \arrow["{P_1}", curve={height=-12pt}, from=1-1, to=2-3]
        \arrow["\src"', curve={height=12pt}, from=1-1, to=3-2]
        \arrow["{\pi_2}", from=2-2, to=2-3]
        \arrow["{\pi_1}"', from=2-2, to=3-2]
        \arrow["\lrcorner"{anchor=center, pos=0.125}, draw=none, from=2-2, to=3-3]
        \arrow["\src", from=2-3, to=3-3]
        \arrow["{P_0}"', from=3-2, to=3-3]
      \end{tikzcd}
    \end{equation*}
  But since the corner object of the pullback is equally well represented by 
  $\dbl{E}_0$ itself, the canonical functor is an equivalence
  if and only if the functor $\src\colon \dbl{E}_1 \to \dbl{E}_0$
  is an equivalence, which happens if and only if $\dbl{E}$ is equivalent to the
  tight double category $\Tight(\dbl{E}_0)$.
\end{proof}

\begin{corollary}
  If $P\colon\dbl{E}\to\dbl{B}$ is a loosely discrete opfibration, then for each object
  $x\in\dbl{B}$, the double-categorical fiber $\dbl{E}_x$ is loosely discrete in
  the sense that, up to isomorphism, its only proarrows and cells are identities.
\end{corollary}
\begin{proof}
  The morphism $\src\colon(\dbl{E}_x)_1\to(\dbl{E}_x)_0$ is always essentially 
  surjective because $\dbl{E}$ has loose identities. That it is also fully 
  faithful is a result of the fact that the only proarrows in $\dbl{E}_x$ are 
  up to isomorphism loose identities and the canonical functor 
  $\langle\src,P_1\rangle$ in \cref{def:loosely-discrete-opfibration} is assumed 
  to be fully faithful.
\end{proof}

Loosely discrete opfibrations over a fixed base assemble into a category.

\begin{definition}[Category of loosely discrete opfibrations]
  \label{def:cat-of-opfibrations}
  A \define{morphism} between loosely discrete double opfibrations
  $P\colon\dbl{E}\to \dbl{B}$ and $Q\colon\dbl{F}\to\dbl{B}$ is a double functor
  $F\colon \dbl{E} \to\dbl{F}$ commuting strictly with the projections into
  $\dbl{B}$. Denote by $\cat{DOpf}(\dbl{B})$ the resulting category of loosely
  discrete double opfibrations over $\dbl{B}$.
\end{definition}

We now give a first few examples of loosely discrete opfibrations.

\begin{example}[Tightly discrete double categories] \label{ex:disc-opfs-are-loosely-disc-opfs}
  Let $P\colon\dbl{E}\to\dbl{B}$ be a double functor between strict double categories
  that are both trivial in the tight direction. That is,
  $\dbl{E} = \Loose(\cat{E})$ and $\dbl{B} = \Loose(\cat{B})$ concentrate the
  morphisms of ordinary categories $\cat{E}$ and $\cat{B}$ in the loose
  direction. In this case, every category in the diagram
  \begin{equation*}
      % https://q.uiver.app/#q=WzAsNSxbMSwxLCJcXGRibHtFfV8wXFx0aW1lc197XFxkYmx7Qn1fMH1cXGRibHtCfV8xIl0sWzIsMSwiXFxkYmx7Qn1fMSJdLFsyLDIsIlxcZGJse0J9XzAiXSxbMSwyLCJcXGRibHtFfV8wIl0sWzAsMCwiXFxkYmx7RX1fMSJdLFswLDEsIlxccGlfMiJdLFsxLDIsIlxcc3JjIl0sWzAsMywiXFxwaV8xIiwyXSxbMywyLCJQXzAiLDJdLFs0LDEsIlBfMSIsMCx7ImN1cnZlIjotMn1dLFs0LDMsIlxcc3JjIiwyLHsiY3VydmUiOjJ9XSxbNCwwLCJcXGxhbmdsZVxcc3JjLFBfMVxccmFuZ2xlIiwxLHsic3R5bGUiOnsiYm9keSI6eyJuYW1lIjoiZGFzaGVkIn19fV0sWzAsMiwiIiwyLHsic3R5bGUiOnsibmFtZSI6ImNvcm5lciJ9fV1d
      \begin{tikzcd}
        {\dbl{E}_1} \\
        & {\dbl{E}_0\times_{\dbl{B}_0}\dbl{B}_1} & {\dbl{B}_1} \\
        & {\dbl{E}_0} & {\dbl{B}_0}
        \arrow["{\langle\src,P_1\rangle}"{description}, dashed, from=1-1, to=2-2]
        \arrow["{P_1}", curve={height=-12pt}, from=1-1, to=2-3]
        \arrow["\src"', curve={height=12pt}, from=1-1, to=3-2]
        \arrow["{\pi_2}", from=2-2, to=2-3]
        \arrow["{\pi_1}"', from=2-2, to=3-2]
        \arrow["\lrcorner"{anchor=center, pos=0.125}, draw=none, from=2-2, to=3-3]
        \arrow["\src", from=2-3, to=3-3]
        \arrow["{P_0}"', from=3-2, to=3-3]
      \end{tikzcd}
  \end{equation*}
  is just a set, and to ask that the dashed morphism be an equivalence is just
  to say it is an isomorphism, i.e., a bijection. Thus, a loosely discrete
  opfibration between tightly discrete double categories is a discrete
  opfibration in the usual sense. Conversely, any discrete opfibration
  $P\colon\cat{E}\to\cat{B}$ can be viewed as a loosely discrete opfibration in this
  way. Thus, as expected, loosely discrete opfibrations generalize the standard
  notion of discrete opfibration for categories.
\end{example}

\begin{example}
  The present notion of a loosely discrete opfibration was in a sense 
  anticipated in \cite[Proposition 2.17]{lambert2021} where given a 
  \emph{discrete double fibration} $P\colon\dbl{E}\to\dbl{B}$ of \emph{strict} 
  double categories, the resulting transpose 
  $P^\top\colon\dbl{E}^\top\to\dbl{B}^\top$ is a loosely discrete opfibration in 
  the present sense. Thus, transposing a loosely discrete opfibration of strict 
  double categories likewise results in a discrete double fibration. Accordingly,
  the transpose of any loosely discrete opfibration of strict double categories 
  is also a double fibration \cite{cruttwell2022} in a rather trivial sense. In this
  way we think of the present theory as transpose to, or orthogonal to, former 
  work on fibrations between double categories.
\end{example}

The following construction offers a more interesting example of a loosely
discrete double opfibration, generalizing the familiar situation that the
projection from a coslice category is a discrete opfibration.

\begin{construction}[Loose coslice] \label{construction:double-coslice}
  Let $\dbl{B}$ be a double category and fix an object $a\in\dbl{B}$. The
  \textbf{loose double coslice under} $a$ is the double category 
  $a\sslash\dbl{B}$ having
    \begin{itemize}[noitemsep]
      \item as objects, an object $x$ together with a proarrow $m\colon a\proto x$;
      \item as arrows $(x,m) \to (z,n)$, an arrow $f: x \to z$ together with a left
        semi-tightly-globular cell
        \begin{equation*}
          % https://q.uiver.app/#q=WzAsNCxbMCwwLCJhIl0sWzEsMCwieCJdLFsxLDEsInoiXSxbMCwxLCJhIl0sWzAsMSwibSIsMCx7InN0eWxlIjp7ImJvZHkiOnsibmFtZSI6ImJhcnJlZCJ9fX1dLFsxLDIsImYiXSxbMCwzLCIiLDIseyJsZXZlbCI6Miwic3R5bGUiOnsiaGVhZCI6eyJuYW1lIjoibm9uZSJ9fX1dLFszLDIsIm4iLDIseyJzdHlsZSI6eyJib2R5Ijp7Im5hbWUiOiJiYXJyZWQifX19XSxbNCw3LCJcXHBoaSIsMSx7InNob3J0ZW4iOnsic291cmNlIjoyMCwidGFyZ2V0IjoyMH0sInN0eWxlIjp7ImJvZHkiOnsibmFtZSI6Im5vbmUifSwiaGVhZCI6eyJuYW1lIjoibm9uZSJ9fX1dXQ==
          \begin{tikzcd}
            a & x \\
            a & z
            \arrow[""{name=0, anchor=center, inner sep=0}, "m"{inner sep=.8ex}, "\shortmid"{marking}, from=1-1, to=1-2]
            \arrow[equals, from=1-1, to=2-1]
            \arrow["f", from=1-2, to=2-2]
            \arrow[""{name=1, anchor=center, inner sep=0}, "n"'{inner sep=.8ex}, "\shortmid"{marking}, from=2-1, to=2-2]
            \arrow["\phi"{description}, draw=none, from=0, to=1]
          \end{tikzcd};
        \end{equation*}
      \item as proarrows $(x,m) \proto (y,m')$, a proarrow $p: x \proto y$ together
        with a tightly globular isocell
        \begin{equation*}
          % https://q.uiver.app/#q=WzAsNSxbMCwwLCJhIl0sWzEsMCwieCJdLFsyLDAsInkiXSxbMCwxLCJhIl0sWzIsMSwieSJdLFswLDEsIm0iLDAseyJzdHlsZSI6eyJib2R5Ijp7Im5hbWUiOiJiYXJyZWQifX19XSxbMSwyLCJwIiwwLHsic3R5bGUiOnsiYm9keSI6eyJuYW1lIjoiYmFycmVkIn19fV0sWzAsMywiIiwyLHsibGV2ZWwiOjIsInN0eWxlIjp7ImhlYWQiOnsibmFtZSI6Im5vbmUifX19XSxbMyw0LCJtJyIsMix7InN0eWxlIjp7ImJvZHkiOnsibmFtZSI6ImJhcnJlZCJ9fX1dLFsyLDQsIiIsMCx7ImxldmVsIjoyLCJzdHlsZSI6eyJoZWFkIjp7Im5hbWUiOiJub25lIn19fV0sWzEsOCwiXFxyaG8iLDEseyJsYWJlbF9wb3NpdGlvbiI6NDAsInNob3J0ZW4iOnsidGFyZ2V0IjoyMH0sInN0eWxlIjp7ImJvZHkiOnsibmFtZSI6Im5vbmUifSwiaGVhZCI6eyJuYW1lIjoibm9uZSJ9fX1dXQ==
          \begin{tikzcd}
            a & x & y \\
            a && y
            \arrow["m"{inner sep=.8ex}, "\shortmid"{marking}, from=1-1, to=1-2]
            \arrow[equals, from=1-1, to=2-1]
            \arrow["p"{inner sep=.8ex}, "\shortmid"{marking}, from=1-2, to=1-3]
            \arrow[equals, from=1-3, to=2-3]
            \arrow[""{name=0, anchor=center, inner sep=0}, "{m'}"'{inner sep=.8ex}, "\shortmid"{marking}, from=2-1, to=2-3]
            \arrow["\rho"{description, pos=0.4}, draw=none, from=1-2, to=0]
          \end{tikzcd};
        \end{equation*}
      \item as cells
        $\inlineCell{(x,m)}{(y,m')}{(z,n)}{(w,n')}{(p,\rho)}{(q,\omega)}{(f,\phi)}{(g,\gamma)}{}$,
        a cell $\inlineCell{x}{y}{z}{w}{p}{q}{f}{g}{\alpha}$ in $\dbl{B}$ for
        which there is an equality
        \begin{equation*}
          % https://q.uiver.app/#q=WzAsOCxbMSwwLCJ4Il0sWzIsMCwieSJdLFsyLDEsInciXSxbMSwxLCJ6Il0sWzAsMCwiYSJdLFswLDEsImEiXSxbMCwyLCJhIl0sWzIsMiwidyJdLFswLDEsInAiLDAseyJzdHlsZSI6eyJib2R5Ijp7Im5hbWUiOiJiYXJyZWQifX19XSxbMSwyLCJnIl0sWzAsMywiZiIsMV0sWzMsMiwicSIsMix7InN0eWxlIjp7ImJvZHkiOnsibmFtZSI6ImJhcnJlZCJ9fX1dLFs0LDAsIm0iLDAseyJzdHlsZSI6eyJib2R5Ijp7Im5hbWUiOiJiYXJyZWQifX19XSxbNCw1LCIiLDEseyJsZXZlbCI6Miwic3R5bGUiOnsiaGVhZCI6eyJuYW1lIjoibm9uZSJ9fX1dLFs1LDMsIm4iLDIseyJzdHlsZSI6eyJib2R5Ijp7Im5hbWUiOiJiYXJyZWQifX19XSxbNSw2LCIiLDAseyJsZXZlbCI6Miwic3R5bGUiOnsiaGVhZCI6eyJuYW1lIjoibm9uZSJ9fX1dLFs2LDcsIm4nIiwyLHsic3R5bGUiOnsiYm9keSI6eyJuYW1lIjoiYmFycmVkIn19fV0sWzIsNywiIiwwLHsibGV2ZWwiOjIsInN0eWxlIjp7ImhlYWQiOnsibmFtZSI6Im5vbmUifX19XSxbOCwxMSwiXFxhbHBoYSIsMSx7InNob3J0ZW4iOnsic291cmNlIjoyMCwidGFyZ2V0IjoyMH0sInN0eWxlIjp7ImJvZHkiOnsibmFtZSI6Im5vbmUifSwiaGVhZCI6eyJuYW1lIjoibm9uZSJ9fX1dLFsxMiwxNCwiXFxwaGkiLDEseyJzaG9ydGVuIjp7InNvdXJjZSI6MjAsInRhcmdldCI6MjB9LCJzdHlsZSI6eyJib2R5Ijp7Im5hbWUiOiJub25lIn0sImhlYWQiOnsibmFtZSI6Im5vbmUifX19XSxbMywxNiwiXFxvbWVnYSIsMSx7ImxhYmVsX3Bvc2l0aW9uIjo0MCwic2hvcnRlbiI6eyJ0YXJnZXQiOjIwfSwic3R5bGUiOnsiYm9keSI6eyJuYW1lIjoibm9uZSJ9LCJoZWFkIjp7Im5hbWUiOiJub25lIn19fV1d
          \begin{tikzcd}
            a & x & y \\
            a & z & w \\
            a && w
            \arrow[""{name=0, anchor=center, inner sep=0}, "m"{inner sep=.8ex}, "\shortmid"{marking}, from=1-1, to=1-2]
            \arrow[equals, from=1-1, to=2-1]
            \arrow[""{name=1, anchor=center, inner sep=0}, "p"{inner sep=.8ex}, "\shortmid"{marking}, from=1-2, to=1-3]
            \arrow["f"{description}, from=1-2, to=2-2]
            \arrow["g", from=1-3, to=2-3]
            \arrow[""{name=2, anchor=center, inner sep=0}, "n"'{inner sep=.8ex}, "\shortmid"{marking}, from=2-1, to=2-2]
            \arrow[equals, from=2-1, to=3-1]
            \arrow[""{name=3, anchor=center, inner sep=0}, "q"'{inner sep=.8ex}, "\shortmid"{marking}, from=2-2, to=2-3]
            \arrow[equals, from=2-3, to=3-3]
            \arrow[""{name=4, anchor=center, inner sep=0}, "{n'}"'{inner sep=.8ex}, "\shortmid"{marking}, from=3-1, to=3-3]
            \arrow["\phi"{description}, draw=none, from=0, to=2]
            \arrow["\alpha"{description}, draw=none, from=1, to=3]
            \arrow["\omega"{description, pos=0.4}, draw=none, from=2-2, to=4]
          \end{tikzcd}
          \qquad=\qquad
          % https://q.uiver.app/#q=WzAsNyxbMiwxLCJ5Il0sWzIsMiwidyJdLFswLDEsImEiXSxbMCwyLCJhIl0sWzAsMCwiYSJdLFsyLDAsInkiXSxbMSwwLCJ4Il0sWzAsMSwiZyJdLFsyLDAsIm0nIiwwLHsic3R5bGUiOnsiYm9keSI6eyJuYW1lIjoiYmFycmVkIn19fV0sWzIsMywiIiwxLHsibGV2ZWwiOjIsInN0eWxlIjp7ImhlYWQiOnsibmFtZSI6Im5vbmUifX19XSxbMywxLCJuJyIsMix7InN0eWxlIjp7ImJvZHkiOnsibmFtZSI6ImJhcnJlZCJ9fX1dLFs0LDIsIiIsMCx7ImxldmVsIjoyLCJzdHlsZSI6eyJoZWFkIjp7Im5hbWUiOiJub25lIn19fV0sWzUsMCwiIiwyLHsibGV2ZWwiOjIsInN0eWxlIjp7ImhlYWQiOnsibmFtZSI6Im5vbmUifX19XSxbNCw2LCJtIiwwLHsic3R5bGUiOnsiYm9keSI6eyJuYW1lIjoiYmFycmVkIn19fV0sWzYsNSwicCIsMCx7InN0eWxlIjp7ImJvZHkiOnsibmFtZSI6ImJhcnJlZCJ9fX1dLFs4LDEwLCJcXGdhbW1hIiwxLHsic2hvcnRlbiI6eyJzb3VyY2UiOjIwLCJ0YXJnZXQiOjIwfSwic3R5bGUiOnsiYm9keSI6eyJuYW1lIjoibm9uZSJ9LCJoZWFkIjp7Im5hbWUiOiJub25lIn19fV0sWzYsOCwiXFxyaG8iLDEseyJsYWJlbF9wb3NpdGlvbiI6MzAsInNob3J0ZW4iOnsidGFyZ2V0IjoyMH0sInN0eWxlIjp7ImJvZHkiOnsibmFtZSI6Im5vbmUifSwiaGVhZCI6eyJuYW1lIjoibm9uZSJ9fX1dXQ==
          \begin{tikzcd}
            a & x & y \\
            a && y \\
            a && w
            \arrow["m"{inner sep=.8ex}, "\shortmid"{marking}, from=1-1, to=1-2]
            \arrow[equals, from=1-1, to=2-1]
            \arrow["p"{inner sep=.8ex}, "\shortmid"{marking}, from=1-2, to=1-3]
            \arrow[equals, from=1-3, to=2-3]
            \arrow[""{name=0, anchor=center, inner sep=0}, "{m'}"{inner sep=.8ex}, "\shortmid"{marking}, from=2-1, to=2-3]
            \arrow[equals, from=2-1, to=3-1]
            \arrow["g", from=2-3, to=3-3]
            \arrow[""{name=1, anchor=center, inner sep=0}, "{n'}"'{inner sep=.8ex}, "\shortmid"{marking}, from=3-1, to=3-3]
            \arrow["\rho"{description, pos=0.3}, draw=none, from=1-2, to=0]
            \arrow["\gamma"{description}, draw=none, from=0, to=1]
          \end{tikzcd}.
        \end{equation*}
    \end{itemize}
  Composition of arrows is by tight composition of cells in
  $\dbl{B}$, while composition of proarrows and cells is by loose composition of
  cells in $\dbl{B}$.

  Moreover, there is a projection double functor
  \begin{equation*}
      J\colon a\sslash\dbl{B}\to \dbl{B}
  \end{equation*}
  that extracts the first component on objects and morphisms and is the identity on cells.
  This projection is a loosely discrete double opfibration: the commutative square
  \begin{equation*}
      % https://q.uiver.app/#q=WzAsNCxbMCwwLCIoYVxcc3NsYXNoXFxkYmx7Qn0pXzEiXSxbMSwwLCJcXGRibHtCfV8xIl0sWzEsMSwiXFxkYmx7Qn1fMCJdLFswLDEsIihhXFxzc2xhc2hcXGRibHtCfSlfMCJdLFswLDEsIkpfMSJdLFsxLDIsIlxcc3JjIl0sWzAsMywiXFxzcmMiLDJdLFszLDIsIkpfMCIsMl1d
      \begin{tikzcd}
        {(a\sslash\dbl{B})_1} & {\dbl{B}_1} \\
        {(a\sslash\dbl{B})_0} & {\dbl{B}_0}
        \arrow["{J_1}", from=1-1, to=1-2]
        \arrow["\src"', from=1-1, to=2-1]
        \arrow["\src", from=1-2, to=2-2]
        \arrow["{J_0}"', from=2-1, to=2-2]
      \end{tikzcd}
  \end{equation*}
  presents the pullback of $\src: \dbl{B}_1 \to \dbl{B}_0$ along $J_0$ up to
  equivalence, due to the definition of proarrows including a tightly globular isocell.
  In fact, there is a canonical choice of pseudo-inverse to the induced functor
  to the pullback, taking the pair $(a \xproto{m} x, x \xproto{p} y)$ to the
  proarrow $(p, 1_{m \odot p}): (x, m) \proto (y, m \odot p)$ in
  $a\sslash\dbl{B}$ and likewise on morphisms. We denote this functor by $\sigma$ as
  in the diagram
  \begin{equation*}
        % https://q.uiver.app/#q=WzAsNSxbMCwwLCIoYVxcc3NsYXNoXFxkYmx7Qn0pXzEiXSxbMiwxLCJcXGRibHtCfV8xIl0sWzIsMiwiXFxkYmx7Qn1fMCJdLFsxLDIsIihhXFxzc2xhc2hcXGRibHtCfSlfMCJdLFsxLDEsIihhXFxzc2xhc2hcXGRibHtCfSlfMFxcdGltZXNfe1xcZGJse0J9XzB9XFxkYmx7Qn1fMSJdLFswLDEsIkpfMSIsMCx7ImN1cnZlIjotM31dLFsxLDIsIlxcc3JjIl0sWzAsMywiXFxzcmMiLDIseyJjdXJ2ZSI6NH1dLFszLDIsIkpfMCIsMl0sWzQsM10sWzQsMV0sWzAsNCwiXFxsYW5nbGVcXHNyYyxKXzFcXHJhbmdsZSIsMix7Im9mZnNldCI6MX1dLFs0LDAsIlxcc2lnbWEiLDIseyJvZmZzZXQiOjIsInN0eWxlIjp7ImJvZHkiOnsibmFtZSI6ImRhc2hlZCJ9fX1dXQ==
        \begin{tikzcd}
          {(a\sslash\dbl{B})_1} \\
          & {(a\sslash\dbl{B})_0\times_{\dbl{B}_0}\dbl{B}_1} & {\dbl{B}_1} \\
          & {(a\sslash\dbl{B})_0} & {\dbl{B}_0}
          \arrow["{\langle\src,J_1\rangle}"', shift right, from=1-1, to=2-2]
          \arrow["{J_1}", curve={height=-18pt}, from=1-1, to=2-3]
          \arrow["\src"', curve={height=24pt}, from=1-1, to=3-2]
          \arrow["\sigma"', shift right=2, dashed, from=2-2, to=1-1]
          \arrow[from=2-2, to=2-3]
          \arrow[from=2-2, to=3-2]
          \arrow["\src", from=2-3, to=3-3]
          \arrow["{J_0}"', from=3-2, to=3-3]
        \end{tikzcd}.
  \end{equation*}
  In this example, $\sigma$ is the \emph{left-adjoint-right-inverse} of an adjoint
  equivalence of categories. Note that the second component of $\sigma$ is not
  just an isomorphism, but in fact the identity cell. This is typical, as we
  shall see below.
\end{construction}

\subsection{Cloven discrete opfibrations} \label{subsection:cleavages}

To motivate this section, we recall more of the theory of discrete opfibrations
internal to a category, for which see \cite[\S{2.1}]{johnstone1977},
\cite[\S{B2.5}]{johnstone2002}, or \cite[\S{V.7}]{maclane1992}. As already
noted, the definition of a discrete opfibration makes sense in any category with
enough finite limits. Let $P\colon\cat{E}\to\cat{B}$ be such a discrete opfibration in
a category $\cat{C}$ with pullbacks. We stress that $\cat{E}$ and $\cat{B}$ are
categories internal to $\cat{C}$. Now, the base category $\cat{B}$, with
object-of-objects $b_0$ and object-of-morphisms $b_1$, determines a monad on
$\cat{C}/b_0$ that acts on a morphism $f\colon a\to b_0$ by first pulling back along
$d_0$ and then composing the second projection with $d_1$, as depicted below:
\begin{equation*}
  % https://q.uiver.app/#q=WzAsNSxbMCwwLCJhXFx0aW1lc197Yl8wfWJfMSJdLFsxLDAsImJfMSJdLFsxLDEsImJfMCJdLFswLDEsImEiXSxbMiwwLCJiXzAiXSxbMCwzLCJcXHBpXzEiLDJdLFszLDIsImYiLDJdLFswLDEsIlxccGlfMiJdLFsxLDIsImRfMCJdLFsxLDQsImRfMSJdLFswLDIsIiIsMSx7InN0eWxlIjp7Im5hbWUiOiJjb3JuZXIifX1dXQ==
  \begin{tikzcd}
    {a\times_{b_0}b_1} & {b_1} & {b_0} \\
    a & {b_0}
    \arrow["{\pi_2}", from=1-1, to=1-2]
    \arrow["{\pi_1}"', from=1-1, to=2-1]
    \arrow["\lrcorner"{anchor=center, pos=0.125}, draw=none, from=1-1, to=2-2]
    \arrow["{d_1}", from=1-2, to=1-3]
    \arrow["{d_0}", from=1-2, to=2-2]
    \arrow["f"', from=2-1, to=2-2]
  \end{tikzcd}.
\end{equation*}
In other words, this operation is base change along $d_0$ followed by
pushforward along $d_1$:
\begin{equation*}
  % https://q.uiver.app/#q=WzAsMyxbMCwwLCJcXGNhdHtDfS9iXzAiXSxbMSwwLCJcXGNhdHtDfS9iXzEiXSxbMiwwLCJcXGNhdHtDfS9iXzAiXSxbMCwxLCJkXzBeKiJdLFsxLDIsIihkXzEpXyEiXV0=
  \begin{tikzcd}
    {\cat{C}/b_0} & {\cat{C}/b_1} & {\cat{C}/b_0}
    \arrow["{d_0^*}", from=1-1, to=1-2]
    \arrow["{(d_1)_!}", from=1-2, to=1-3]
  \end{tikzcd}.
\end{equation*}
We refer to this monad as the \define{pull-push monad associated to} $\cat{B}$.
The relevant result is that the forgetful functor from internal discrete
opfibrations to $\cat{C}/b_0$ is \emph{monadic}: the category of internal
discrete opfibrations is equivalent, even isomorphic, to the category of
algebras for the pull-push monad associated to $\cat{B}$; see \cite[Proposition
2.21]{johnstone1977} or \cite[Theorem V.7.2]{maclane1992}.

The thesis of this section is that loosely discrete opfibrations have a similar
algebraic characterization. However, due to the weakening of the standard
internalized definition, more care is required. First, we construct the
candidate pseudo-monad.

\begin{construction}[Pull-push monad] \label{construction:pull-push-monad}
  Let $\dbl{B}$ be a double category, denoting its loose source and target maps
  by $\src,\tgt\colon \dbl{B}_1\rightrightarrows\dbl{B}_0$. Let $\Cat/\dbl{B}_0$
  denote the strict slice 2-category of $\Cat$ over $\dbl{B}_0$. Define a
  2-functor
    \begin{equation*}
      \Cat/\dbl{B}_0\to\Cat/\dbl{B}_0
    \end{equation*} 
  on objects by sending a functor $F\colon\cat{A}\to\dbl{B}_0$ to the top composite
  $\tgt \circ \pi_2$ in the diagram
    \begin{equation*}
      % https://q.uiver.app/#q=WzAsNSxbMCwwLCJcXGNhdHtBfVxcdGltZXNfe1xcZGJse0J9XzB9XFxkYmx7Qn1fMSJdLFsxLDAsIlxcZGJse0J9XzEiXSxbMSwxLCJcXGRibHtCfV8wIl0sWzAsMSwiXFxjYXR7QX0iXSxbMiwwLCJcXGRibHtCfV8wIl0sWzAsMywiXFxwaV8xIiwyXSxbMywyLCJGIiwyXSxbMCwxLCJcXHBpXzIiXSxbMSwyLCJcXHNyYyJdLFsxLDQsIlxcdGd0Il0sWzAsMiwiIiwxLHsic3R5bGUiOnsibmFtZSI6ImNvcm5lciJ9fV1d
      \begin{tikzcd}
        {\cat{A}\times_{\dbl{B}_0}\dbl{B}_1} & {\dbl{B}_1} & {\dbl{B}_0} \\
        {\cat{A}} & {\dbl{B}_0}
        \arrow["{\pi_2}", from=1-1, to=1-2]
        \arrow["{\pi_1}"', from=1-1, to=2-1]
        \arrow["\lrcorner"{anchor=center, pos=0.125}, draw=none, from=1-1, to=2-2]
        \arrow["\tgt", from=1-2, to=1-3]
        \arrow["\src", from=1-2, to=2-2]
        \arrow["F"', from=2-1, to=2-2]
      \end{tikzcd}
    \end{equation*}
  and on morphisms and cells of $\Cat/\dbl{B}_0$ using the
  universal property of the (strict) pullback. The assignment is thus
  given by pulling back along the source functor and then pushing forward along the
  target functor. This 2-functor underlies a pseudo-monad whose multiplication and unit
  2-natural transformations have components given by the top morphisms in the
  diagrams
    \begin{equation*}
      % https://q.uiver.app/#q=WzAsMyxbMCwwLCJcXGNhdHtBfVxcdGltZXNfe1xcZGJse0J9XzB9XFxkYmx7Qn1fMVxcdGltZXNfe1xcZGJse0J9XzB9XFxkYmx7Qn1fMSJdLFsyLDAsIlxcY2F0e0F9XFx0aW1lc197XFxkYmx7Qn1fMH1cXGRibHtCfV8xIl0sWzEsMSwiXFxkYmx7Qn1fMCJdLFswLDEsIjFcXHRpbWVzKC1cXG9kb3QtKSJdLFswLDIsIlxcdGd0IFxccGlfMiIsMl0sWzEsMiwiXFx0Z3QgXFxwaV8yIl1d
      \begin{tikzcd}
        {\cat{A}\times_{\dbl{B}_0}\dbl{B}_1\times_{\dbl{B}_0}\dbl{B}_1} && {\cat{A}\times_{\dbl{B}_0}\dbl{B}_1} \\
        & {\dbl{B}_0}
        \arrow["{1\times(-\odot-)}", from=1-1, to=1-3]
        \arrow["{\tgt \pi_3}"', from=1-1, to=2-2]
        \arrow["{\tgt \pi_2}", from=1-3, to=2-2]
      \end{tikzcd}
      \qquad\text{and}\qquad
      % https://q.uiver.app/#q=WzAsMyxbMCwwLCJcXGNhdHtBfSJdLFsyLDAsIlxcY2F0e0F9XFx0aW1lc197XFxkYmx7Qn1fMH1cXGRibHtCfV8xIl0sWzEsMSwiXFxkYmx7Qn1fMCJdLFswLDEsIlxcbGFuZ2xlMSxcXGlkIEZcXHJhbmdsZSJdLFswLDIsIkYiLDJdLFsxLDIsIlxcdGd0IFxccGlfMiJdXQ==
      \begin{tikzcd}
        {\cat{A}} && {\cat{A}\times_{\dbl{B}_0}\dbl{B}_1} \\
        & {\dbl{B}_0}
        \arrow["{\langle1,\id F\rangle}", from=1-1, to=1-3]
        \arrow["F"', from=1-1, to=2-2]
        \arrow["{\tgt \pi_2}", from=1-3, to=2-2]
      \end{tikzcd}.
    \end{equation*}
  The associators and unitors of the pseudo-monad are of the form
    \begin{equation*}
      % https://q.uiver.app/#q=WzAsNCxbMCwwLCJcXGNhdHtBfVxcdGltZXNfe1xcZGJse0J9XzB9XFxkYmx7Qn1fMVxcdGltZXNfe1xcZGJse0J9XzB9XFxkYmx7Qn1fMVxcdGltZXNfe1xcZGJse0J9XzB9XFxkYmx7Qn1fMSJdLFsxLDAsIlxcY2F0e0F9XFx0aW1lc197XFxkYmx7Qn1fMH1cXGRibHtCfV8xXFx0aW1lc197XFxkYmx7Qn1fMH1cXGRibHtCfV8xIl0sWzEsMSwiXFxjYXR7QX1cXHRpbWVzX3tcXGRibHtCfV8wfVxcZGJse0J9XzEiXSxbMCwxLCJcXGNhdHtBfVxcdGltZXNfe1xcZGJse0J9XzB9XFxkYmx7Qn1fMVxcdGltZXNfe1xcZGJse0J9XzB9XFxkYmx7Qn1fMSJdLFswLDEsIjFcXHRpbWVzXFxvZG90XFx0aW1lczEiXSxbMSwyLCIxXFx0aW1lc1xcb2RvdCJdLFswLDMsIjFcXHRpbWVzMVxcdGltZXNcXG9kb3QiLDJdLFszLDIsIjFcXHRpbWVzXFxvZG90IiwyXSxbMSwzLCIxXFx0aW1lc1xcbWF0aGZyYWt7YX0iLDEseyJzdHlsZSI6eyJib2R5Ijp7Im5hbWUiOiJub25lIn0sImhlYWQiOnsibmFtZSI6Im5vbmUifX19XV0=
      \begin{tikzcd}[column sep=large]
        {\cat{A}\times_{\dbl{B}_0}\dbl{B}_1\times_{\dbl{B}_0}\dbl{B}_1\times_{\dbl{B}_0}\dbl{B}_1} & {\cat{A}\times_{\dbl{B}_0}\dbl{B}_1\times_{\dbl{B}_0}\dbl{B}_1} \\
        {\cat{A}\times_{\dbl{B}_0}\dbl{B}_1\times_{\dbl{B}_0}\dbl{B}_1} & {\cat{A}\times_{\dbl{B}_0}\dbl{B}_1}
        \arrow["{1\times\odot\times1}", from=1-1, to=1-2]
        \arrow["{1\times1\times\odot}"', from=1-1, to=2-1]
        \arrow["{1\times\mathfrak{a}}"{description}, draw=none, from=1-2, to=2-1]
        \arrow["{1\times\odot}", from=1-2, to=2-2]
        \arrow["{1\times\odot}"', from=2-1, to=2-2]
      \end{tikzcd}
    \end{equation*}
  and 
    \begin{equation*}
      % https://q.uiver.app/#q=WzAsMyxbMCwwLCJcXGNhdHtBfVxcdGltZXNfe1xcZGJse0J9XzB9XFxkYmx7Qn1fMSJdLFsxLDEsIlxcY2F0e0F9XFx0aW1lc197XFxkYmx7Qn1fMH1cXGRibHtCfV8xIl0sWzEsMCwiXFxjYXR7QX1cXHRpbWVzX3tcXGRibHtCfV8wfVxcZGJse0J9XzFcXHRpbWVzX3tcXGRibHtCfV8wfVxcZGJse0J9XzEiXSxbMCwxLCIiLDIseyJsZXZlbCI6Miwic3R5bGUiOnsiaGVhZCI6eyJuYW1lIjoibm9uZSJ9fX1dLFswLDIsIlxcbGFuZ2xlIDEsXFxpZCBcXHRndCBcXHBpXzIgXFxyYW5nbGUiXSxbMiwxLCIxXFx0aW1lc1xcb2RvdCJdLFsyLDMsIjFcXHRpbWVzXFxyaG8iLDEseyJzaG9ydGVuIjp7InRhcmdldCI6MjB9LCJzdHlsZSI6eyJib2R5Ijp7Im5hbWUiOiJub25lIn0sImhlYWQiOnsibmFtZSI6Im5vbmUifX19XV0=
      \begin{tikzcd}[sep=large]
        {\cat{A}\times_{\dbl{B}_0}\dbl{B}_1} & {\cat{A}\times_{\dbl{B}_0}\dbl{B}_1\times_{\dbl{B}_0}\dbl{B}_1} \\
        & {\cat{A}\times_{\dbl{B}_0}\dbl{B}_1}
        \arrow["{\langle 1,\id \tgt \pi_2 \rangle}", from=1-1, to=1-2]
        \arrow[""{name=0, anchor=center, inner sep=0}, equals, from=1-1, to=2-2]
        \arrow["{1\times\odot}", from=1-2, to=2-2]
        \arrow["{1\times\rho}"{description}, draw=none, from=1-2, to=0]
      \end{tikzcd}
      \qquad
      % https://q.uiver.app/#q=WzAsMyxbMCwwLCJcXGNhdHtBfVxcdGltZXNfe1xcZGJse0J9XzB9XFxkYmx7Qn1fMSJdLFsxLDEsIlxcY2F0e0F9XFx0aW1lc197XFxkYmx7Qn1fMH1cXGRibHtCfV8xIl0sWzEsMCwiXFxjYXR7QX1cXHRpbWVzX3tcXGRibHtCfV8wfVxcZGJse0J9XzFcXHRpbWVzX3tcXGRibHtCfV8wfVxcZGJse0J9XzEiXSxbMCwxLCIiLDIseyJsZXZlbCI6Miwic3R5bGUiOnsiaGVhZCI6eyJuYW1lIjoibm9uZSJ9fX1dLFswLDIsIlxcbGFuZ2xlIDEsXFxpZCBGIFxccmFuZ2xlIFxcdGltZXMgMSJdLFsyLDEsIjFcXHRpbWVzXFxvZG90Il0sWzIsMywiMVxcdGltZXNcXGxhbWJkYSIsMSx7InNob3J0ZW4iOnsidGFyZ2V0IjoyMH0sInN0eWxlIjp7ImJvZHkiOnsibmFtZSI6Im5vbmUifSwiaGVhZCI6eyJuYW1lIjoibm9uZSJ9fX1dXQ==
      \begin{tikzcd}[sep=large]
        {\cat{A}\times_{\dbl{B}_0}\dbl{B}_1} & {\cat{A}\times_{\dbl{B}_0}\dbl{B}_1\times_{\dbl{B}_0}\dbl{B}_1} \\
        & {\cat{A}\times_{\dbl{B}_0}\dbl{B}_1}
        \arrow["{\langle 1,\id F \rangle \times 1}", from=1-1, to=1-2]
        \arrow[""{name=0, anchor=center, inner sep=0}, equals, from=1-1, to=2-2]
        \arrow["{1\times\odot}", from=1-2, to=2-2]
        \arrow["{1\times\lambda}"{description}, draw=none, from=1-2, to=0]
      \end{tikzcd},
    \end{equation*}
  built from the associators and unitors of $\dbl{B}$. That the components are appropriately natural and that the accompanying
  transformations satisfy the various conditions of a pseudo-monad on 
  $\Cat/\dbl{B}_0$ is straightforward to check using the universal properties of 
  the pullback defining the object assignment. Throughout we shall refer to such 
  a pseudo-monad, determined by the double category $\dbl{B}$, as the \define{pull-push
  monad associated to} $\dbl{B}$.
\end{construction}

\begin{remark}[Pseudo-algebras]
  As pull-push along source and then target defines a pseudo-monad, there
  is a 2-category of pseudo-algebras of interest. Such an algebra is a
  functor $F\colon\cat{A}\to\dbl{B}_0$ together with a structure map and two 
  cells
    \begin{equation*}
      % https://q.uiver.app/#q=WzAsMyxbMCwwLCJcXGNhdHtBfVxcdGltZXNfe1xcZGJse0J9XzB9XFxkYmx7Qn1fMSJdLFsxLDEsIlxcZGJse0J9XzAiXSxbMiwwLCJcXGNhdHtBfSJdLFswLDEsIlxcdGd0LlxccGlfMiIsMl0sWzIsMSwiRiJdLFswLDIsIkwiXV0=
      \begin{tikzcd}
        {\cat{A}\times_{\dbl{B}_0}\dbl{B}_1} && {\cat{A}} \\
        & {\dbl{B}_0}
        \arrow["L", from=1-1, to=1-3]
        \arrow["{\tgt.\pi_2}"', from=1-1, to=2-2]
        \arrow["F", from=1-3, to=2-2]
      \end{tikzcd}
      \qquad
      % https://q.uiver.app/#q=WzAsNCxbMSwwLCJcXGNhdHtBfVxcdGltZXNfe1xcZGJse0J9XzB9XFxkYmx7Qn1fMSJdLFsxLDEsIlxcY2F0e0F9Il0sWzAsMSwiXFxjYXR7QX1cXHRpbWVzX3tcXGRibHtCfV8wfVxcZGJse0J9XzEiXSxbMCwwLCJcXGNhdHtBfVxcdGltZXNfe1xcZGJse0J9XzB9XFxkYmx7Qn1fMVxcdGltZXNfe1xcZGJse0J9XzB9XFxkYmx7Qn1fMSJdLFswLDEsIkwiXSxbMiwxLCJMIiwyXSxbMywwLCJMXFx0aW1lcyAxIl0sWzMsMiwiMSBcXHRpbWVzICgtXFxvZG90LSkiLDJdLFs2LDUsIlxcbXUiLDEseyJzaG9ydGVuIjp7InNvdXJjZSI6MjAsInRhcmdldCI6MjB9LCJzdHlsZSI6eyJib2R5Ijp7Im5hbWUiOiJub25lIn0sImhlYWQiOnsibmFtZSI6Im5vbmUifX19XV0=
      \begin{tikzcd}
        {\cat{A}\times_{\dbl{B}_0}\dbl{B}_1\times_{\dbl{B}_0}\dbl{B}_1} & {\cat{A}\times_{\dbl{B}_0}\dbl{B}_1} \\
        {\cat{A}\times_{\dbl{B}_0}\dbl{B}_1} & {\cat{A}}
        \arrow[""{name=0, anchor=center, inner sep=0}, "{L\times 1}", from=1-1, to=1-2]
        \arrow["{1 \times (-\odot-)}"', from=1-1, to=2-1]
        \arrow["L", from=1-2, to=2-2]
        \arrow[""{name=1, anchor=center, inner sep=0}, "L"', from=2-1, to=2-2]
        \arrow["\mu"{description}, draw=none, from=0, to=1]
      \end{tikzcd}
      \qquad
      % https://q.uiver.app/#q=WzAsMyxbMSwwLCJcXGNhdHtBfVxcdGltZXNfe1xcZGJse0J9XzB9XFxkYmx7Qn1fMSJdLFsxLDEsIlxcY2F0e0F9Il0sWzAsMCwiXFxjYXR7QX0iXSxbMCwxLCJMIl0sWzIsMCwiXFxsYW5nbGUgMSxcXGlkIEYgXFxyYW5nbGUiXSxbMiwxLCIiLDIseyJsZXZlbCI6Miwic3R5bGUiOnsiaGVhZCI6eyJuYW1lIjoibm9uZSJ9fX1dLFswLDUsIlxca2FwcGEiLDEseyJzaG9ydGVuIjp7InRhcmdldCI6MjB9LCJzdHlsZSI6eyJib2R5Ijp7Im5hbWUiOiJub25lIn0sImhlYWQiOnsibmFtZSI6Im5vbmUifX19XV0=
      \begin{tikzcd}
        {\cat{A}} & {\cat{A}\times_{\dbl{B}_0}\dbl{B}_1} \\
        & {\cat{A}}
        \arrow["{\langle 1,\id F \rangle}", from=1-1, to=1-2]
        \arrow[""{name=0, anchor=center, inner sep=0}, equals, from=1-1, to=2-2]
        \arrow["L", from=1-2, to=2-2]
        \arrow["\kappa"{description}, draw=none, from=1-2, to=0]
      \end{tikzcd}
    \end{equation*}
  satisfying the appropriate axioms (\cite[\S1]{lack2002} and
  \cref{remark:cleavage-pseudo-algebra-relations}). The associativity axiom
  in particular we think of as a \emph{codescent condition} as described below. 
  Our thesis will be that loosely discrete opfibrations equipped with cleavages are 
  closely related to such algebras. It might seem strange that cleavages are 
  required but this is exactly the situation for ordinary opfibrations.
\end{remark}

Let $P\colon \cat{E}\to\cat{B}$ be an opfibration of categories equipped with a
cleavage $\sigma$ associating to every pair $(a,f)$ a chosen opcartesian arrow
$\sigma(a,f) \colon a \to f_!(a)$ of $\cat{E}$. An elegant categorical description of the
cleavage is as a functor $\sigma$, as shown below
  \begin{equation*}
    % https://q.uiver.app/#q=WzAsNSxbMCwwLCJcXGNhdHtFfV57XFxjYXQgMn0iXSxbMiwxLCJcXGNhdHtCfV57XFxjYXQgMn0iXSxbMiwyLCJcXGNhdHtCfSJdLFsxLDIsIlxcY2F0e0V9Il0sWzEsMSwiXFxjYXR7RX1cXHRpbWVzX3tcXGNhdHtCfX1cXGNhdHtCfV57XFxjYXQgMn0iXSxbMCwxLCJQXntcXGNhdCAyfSIsMCx7ImN1cnZlIjotM31dLFsxLDIsImRfMCJdLFswLDMsImRfMCIsMix7ImN1cnZlIjo1fV0sWzMsMiwiUCIsMl0sWzQsMywiXFxwaV8xIiwyXSxbNCwxLCJcXHBpXzIiXSxbMCw0LCJcXGxhbmdsZSBkXzAsIFBee1xcY2F0IDJ9IFxccmFuZ2xlIiwyLHsibGFiZWxfcG9zaXRpb24iOjcwfV0sWzQsMCwiXFxzaWdtYSIsMix7Im9mZnNldCI6Miwic3R5bGUiOnsiYm9keSI6eyJuYW1lIjoiZGFzaGVkIn19fV1d
    \begin{tikzcd}
      {\cat{E}^{\cat 2}} && \\
      & {\cat{E}\times_{\cat{B}}\cat{B}^{\cat 2}} & {\cat{B}^{\cat 2}} \\
      & {\cat{E}} & {\cat{B}}
      \arrow["{\langle d_0, P^{\cat 2} \rangle}"'{pos=0.7}, from=1-1, to=2-2]
      \arrow["{P^{\cat 2}}", curve={height=-18pt}, from=1-1, to=2-3]
      \arrow["{d_0}"', curve={height=30pt}, from=1-1, to=3-2]
      \arrow["\sigma"', shift right=2, dashed, from=2-2, to=1-1]
      \arrow["{\pi_2}", from=2-2, to=2-3]
      \arrow["{\pi_1}"', from=2-2, to=3-2]
      \arrow["{d_0}", from=2-3, to=3-3]
      \arrow["P"', from=3-2, to=3-3]
    \end{tikzcd},
  \end{equation*}
that is a left-adjoint-right-inverse (LARI) to the pair functor into the
pullback; see \cite[\mbox{Proposition 3.11}]{gray1966} for the dual statement
for fibrations. Now, the diagram recalls that defining discrete
opfibrations with some changes, namely, that the internal categories are formed by the
cotensors $d_0,d_1\colon \cat{B}^{\cat 2}\rightrightarrows \cat{B}$ and $d_0,
d_1\colon \cat{E}^{\cat 2}\rightrightarrows \cat{E}$ in place of the generic
internal categories and that we specifically ask for $\sigma$ to be a LARI.
Additionally, the cotensor $\cat{B}^{\cat 2}$ induces a 2-monad on the slice
$\Cat/\cat{B}$ following the familiar pattern:
  \begin{equation*}
    % https://q.uiver.app/#q=WzAsMyxbMCwwLCJcXGNhdC9cXGNhdHtCfSJdLFsxLDAsIlxcY2F0L1xcY2F0e0J9XntcXGNhdCAyfSJdLFsyLDAsIlxcY2F0L1xcY2F0e0J9Il0sWzAsMSwiZF8wXioiXSxbMSwyLCIoZF8xKV8hIl1d
    \begin{tikzcd}
      {\Cat/\cat{B}} & {\Cat/\cat{B}^{\cat 2}} & {\Cat/\cat{B}}
      \arrow["{d_0^*}", from=1-1, to=1-2]
      \arrow["{(d_1)_!}", from=1-2, to=1-3]
    \end{tikzcd}.
  \end{equation*}
That is, a functor $F: \cat{A} \to \cat{B}$ is sent to the top composite in the diagram
  \begin{equation*}
    % https://q.uiver.app/#q=WzAsNSxbMCwwLCJcXGNhdHtBfVxcdGltZXNfe1xcY2F0e0J9fVxcY2F0e0J9XntcXGNhdCAyfSJdLFsxLDAsIlxcY2F0e0J9XntcXGNhdCAyfSJdLFsxLDEsIlxcY2F0e0J9Il0sWzAsMSwiXFxjYXR7QX0iXSxbMiwwLCJcXGNhdHtCfSJdLFswLDMsIlxccGlfMSIsMl0sWzMsMiwiRiIsMl0sWzAsMSwiXFxwaV8yIl0sWzEsMiwiZF8wIl0sWzEsNCwiZF8xIl0sWzAsMiwiIiwxLHsic3R5bGUiOnsibmFtZSI6ImNvcm5lciJ9fV1d
    \begin{tikzcd}
      {\cat{A}\times_{\cat{B}}\cat{B}^{\cat 2}} & {\cat{B}^{\cat 2}} & {\cat{B}} \\
      {\cat{A}} & {\cat{B}}
      \arrow["{\pi_2}", from=1-1, to=1-2]
      \arrow["{\pi_1}"', from=1-1, to=2-1]
      \arrow["\lrcorner"{anchor=center, pos=0.125}, draw=none, from=1-1, to=2-2]
      \arrow["{d_1}", from=1-2, to=1-3]
      \arrow["{d_0}", from=1-2, to=2-2]
      \arrow["F"', from=2-1, to=2-2]
    \end{tikzcd}.
  \end{equation*}
The point is that it is precisely the cleavage $\sigma$ that makes the opfibration into a
(pseudo!) algebra for this 2-monad. For it specifies without ambiguity the value
of the action functor
  \begin{equation*}
    % https://q.uiver.app/#q=WzAsMyxbMCwwLCJcXGNhdHtFfVxcdGltZXNfe1xcY2F0e0J9fVxcY2F0e0J9XntcXGNhdCAyfSJdLFsxLDAsIlxcY2F0e0V9XntcXGNhdCAyfSJdLFsyLDAsIlxcY2F0e0V9Il0sWzAsMSwiXFxzaWdtYSJdLFsxLDIsImRfMSJdXQ==
    \begin{tikzcd}
      {\cat{E}\times_{\cat{B}}\cat{B}^{\cat 2}} & {\cat{E}^{\cat 2}} & {\cat{E}}
      \arrow["\sigma", from=1-1, to=1-2]
      \arrow["{d_1}", from=1-2, to=1-3]
    \end{tikzcd},
  \end{equation*}
and the normalized pseudo-functoriality of the cleavage is precisely the 
associativity and unitality of the pseudo-algebra structure.
Note that both cotensors are the underlying spans of simple double categories, 
namely, the \emph{quintet constructions} on the 1-categories $\cat{B}$ and 
$\cat{E}$.

In our generalization of a cleavage for an arbitrary base double category
$\dbl{B}$, we ask for a LARI split adjoint equivalence rather than a mere
adjunction to stay within the setting of discrete opfibrations. Notice that the
pseudo-algebra structure is now axiomatized since the unique algebra structure
present for ordinary discrete opfibrations is lost in the double-categorical
weakening.

\begin{definition}[Cleavage] \label{def:cleavage-for-disc-opfibration}
  A \textbf{cleavage} for a loosely discrete opfibration $P\colon\dbl{E}\to
  \dbl{B}$ is a functor
    \begin{equation*}
      \sigma: \dbl{E}_0 \times_{\dbl{B}_0} \dbl{B}_1 \to \dbl{E}_1
    \end{equation*}
  that is a left-adjoint-right-inverse to the canonically induced functor $\langle\src,P_1\rangle$ below
    \begin{equation*}
      % https://q.uiver.app/#q=WzAsNSxbMCwwLCJcXGRibHtFfV8xIl0sWzIsMSwiXFxkYmx7Qn1fMSJdLFsyLDIsIlxcZGJse0J9XzAiXSxbMSwyLCJcXGRibHtFfV8wIl0sWzEsMSwiXFxkYmx7RX1fMFxcdGltZXNfe1xcZGJse0J9XzB9XFxkYmx7Qn1fMSJdLFswLDEsIlBfMSIsMCx7ImN1cnZlIjotM31dLFsxLDIsIlxcc3JjIl0sWzAsMywiXFxzcmMiLDIseyJjdXJ2ZSI6NX1dLFszLDIsIlBfMCIsMl0sWzQsMywiXFxwaV8xIiwyXSxbNCwxLCJcXHBpXzIiXSxbMCw0LCJcXGxhbmdsZVxcc3JjLFBfMVxccmFuZ2xlIiwyLHsibGFiZWxfcG9zaXRpb24iOjYwfV0sWzQsMCwiXFxzaWdtYSIsMix7Im9mZnNldCI6Miwic3R5bGUiOnsiYm9keSI6eyJuYW1lIjoiZGFzaGVkIn19fV1d
      \begin{tikzcd}
        {\dbl{E}_1} \\
        & {\dbl{E}_0\times_{\dbl{B}_0}\dbl{B}_1} & {\dbl{B}_1} \\
        & {\dbl{E}_0} & {\dbl{B}_0}
        \arrow["{\langle\src,P_1\rangle}"'{pos=0.6}, from=1-1, to=2-2]
        \arrow["{P_1}", curve={height=-18pt}, from=1-1, to=2-3]
        \arrow["\src"', curve={height=30pt}, from=1-1, to=3-2]
        \arrow["\sigma"', shift right=2, dashed, from=2-2, to=1-1]
        \arrow["{\pi_2}", from=2-2, to=2-3]
        \arrow["{\pi_1}"', from=2-2, to=3-2]
        \arrow["\src", from=2-3, to=3-3]
        \arrow["{P_0}"', from=3-2, to=3-3]
      \end{tikzcd},
    \end{equation*}
  such that the composite $\tgt \sigma: \dbl{E}_0 \times_{\dbl{B}_0} \dbl{B}_1 \to \dbl{E}_0$
  is the structure map making $P_0: \dbl{E}_0 \to \dbl{B}_0$ into a
  pseudo-algebra for the pull-push monad associated to
  $\dbl{B}$ (\cref{construction:pull-push-monad}).
  Additionally, $\sigma$ should satisfy the following coherence conditions:
    \begin{itemize}
      \item \emph{Associativity}: the cell
        \begin{equation*}
          % https://q.uiver.app/#q=WzAsNyxbMCwxLCJcXGRibHtFfV8xXnsoMil9Il0sWzEsMCwiXFxkYmx7RX1fMFxcdGltZXNfe1xcZGJse0J9XzB9XFxkYmx7Qn1fMV57KDIpfSJdLFsyLDEsIlxcZGJse0V9XzBcXHRpbWVzX3tcXGRibHtCfV8wfVxcZGJse0J9XzEiXSxbMSwyLCJcXGRibHtFfV8xIl0sWzMsMSwiXFxkYmx7RX1fMCJdLFsyLDAsIlxcZGJse0V9XzBcXHRpbWVzX3tcXGRibHtCfV8wfVxcZGJse0J9XzEiXSxbMiwyLCJcXGRibHtFfV8wIl0sWzAsMSwiXFxsYW5nbGVcXHNyYyxQXzFcXHJhbmdsZVxcdGltZXMgUCJdLFsxLDIsIjFcXHRpbWVzXFxvZG90IiwyXSxbMCwzLCJcXG9kb3QiLDJdLFszLDIsIlxcbGFuZ2xlXFxzcmMsUF8xXFxyYW5nbGUiXSxbMiw0LCJcXHRndCBcXHNpZ21hIiwyXSxbMSw1LCJcXHRndCBcXHNpZ21hIFxcdGltZXMgMSJdLFs1LDQsIlxcdGd0IFxcc2lnbWEiXSxbMyw2LCJcXHRndCIsMl0sWzYsNCwiMSIsMl0sWzEsMywiXFxtYXRoZnJha3thfSIsMSx7InN0eWxlIjp7ImJvZHkiOnsibmFtZSI6Im5vbmUifSwiaGVhZCI6eyJuYW1lIjoibm9uZSJ9fX1dLFs1LDIsIlxcbXUiLDEseyJzdHlsZSI6eyJib2R5Ijp7Im5hbWUiOiJub25lIn0sImhlYWQiOnsibmFtZSI6Im5vbmUifX19XSxbMiw2LCJcXHRndCBcXGVwc2lsb24iLDEseyJzdHlsZSI6eyJib2R5Ijp7Im5hbWUiOiJub25lIn0sImhlYWQiOnsibmFtZSI6Im5vbmUifX19XV0=
          \begin{tikzcd}
            & {\dbl{E}_0\times_{\dbl{B}_0}\dbl{B}_1^{(2)}} & {\dbl{E}_0\times_{\dbl{B}_0}\dbl{B}_1} & \\
            {\dbl{E}_1^{(2)}} && {\dbl{E}_0\times_{\dbl{B}_0}\dbl{B}_1} & {\dbl{E}_0} \\
            & {\dbl{E}_1} & {\dbl{E}_0}
            \arrow["{\tgt \sigma \times 1}", from=1-2, to=1-3]
            \arrow["{1\times\odot}"', from=1-2, to=2-3]
            \arrow["{\mathfrak{a}}"{description}, draw=none, from=1-2, to=3-2]
            \arrow["\mu"{description}, draw=none, from=1-3, to=2-3]
            \arrow["{\tgt \sigma}", from=1-3, to=2-4]
            \arrow["{\langle\src,P_1\rangle\times P}", from=2-1, to=1-2]
            \arrow["\odot"', from=2-1, to=3-2]
            \arrow["{\tgt \sigma}"', from=2-3, to=2-4]
            \arrow["{\tgt \epsilon}"{description}, draw=none, from=2-3, to=3-3]
            \arrow["{\langle\src,P_1\rangle}", from=3-2, to=2-3]
            \arrow["\tgt"', from=3-2, to=3-3]
            \arrow["1"', from=3-3, to=2-4]
          \end{tikzcd}
        \end{equation*}
      is equal to 
        \begin{equation*}
          % https://q.uiver.app/#q=WzAsNyxbMCwxLCJcXGRibHtFfV8xXnsoMil9Il0sWzEsMCwiXFxkYmx7RX1fMFxcdGltZXNfe1xcZGJse0J9XzB9XFxkYmx7Qn1fMV57KDIpfSJdLFsxLDIsIlxcZGJse0V9XzEiXSxbMywxLCJcXGRibHtFfV8wIl0sWzIsMCwiXFxkYmx7RX1fMFxcdGltZXNfe1xcZGJse0J9XzB9XFxkYmx7Qn1fMSJdLFsyLDIsIlxcZGJse0V9XzAiXSxbMSwxLCJcXGRibHtFfV8xIl0sWzAsMSwiXFxsYW5nbGVcXHNyYyxQXzFcXHJhbmdsZVxcdGltZXMgUCJdLFswLDIsIlxcb2RvdCIsMl0sWzEsNCwiXFx0Z3RcXHNpZ21hXFx0aW1lcyAxIl0sWzQsMywiXFx0Z3QgXFxzaWdtYSJdLFsyLDUsIlxcdGd0IiwyXSxbNSwzLCIxIiwyXSxbNiw0LCJcXGxhbmdsZVxcc3JjLFBfMVxccmFuZ2xlIiwyXSxbNiw1LCJcXHRndCJdLFswLDYsIlxccGlfMiIsMl0sWzQsNSwiXFx0Z3QgXFxlcHNpbG9uIiwxLHsic3R5bGUiOnsiYm9keSI6eyJuYW1lIjoibm9uZSJ9LCJoZWFkIjp7Im5hbWUiOiJub25lIn19fV0sWzEsNiwiXFx0Z3RcXGVwc2lsb25cXHRpbWVzIDEiLDEseyJzdHlsZSI6eyJib2R5Ijp7Im5hbWUiOiJub25lIn0sImhlYWQiOnsibmFtZSI6Im5vbmUifX19XV0=
          \begin{tikzcd}
            & {\dbl{E}_0\times_{\dbl{B}_0}\dbl{B}_1^{(2)}} & {\dbl{E}_0\times_{\dbl{B}_0}\dbl{B}_1} & \\
            {\dbl{E}_1^{(2)}} & {\dbl{E}_1} && {\dbl{E}_0} \\
            & {\dbl{E}_1} & {\dbl{E}_0}
            \arrow["{\tgt\sigma\times 1}", from=1-2, to=1-3]
            \arrow["{\tgt\epsilon\times 1}"{description}, draw=none, from=1-2, to=2-2]
            \arrow["{\tgt \sigma}", from=1-3, to=2-4]
            \arrow["{\tgt \epsilon}"{description}, draw=none, from=1-3, to=3-3]
            \arrow["{\langle\src,P_1\rangle\times P}", from=2-1, to=1-2]
            \arrow["{\pi_2}"', from=2-1, to=2-2]
            \arrow["\odot"', from=2-1, to=3-2]
            \arrow["{\langle\src,P_1\rangle}"', from=2-2, to=1-3]
            \arrow["\tgt", from=2-2, to=3-3]
            \arrow["\tgt"', from=3-2, to=3-3]
            \arrow["1"', from=3-3, to=2-4]
          \end{tikzcd}
        \end{equation*}
      where $\mu$ denotes the algebra structure isocell and $\epsilon$ is the counit of the adjunction; and 
      \item \emph{Unitality}:
        \begin{equation*}
          % https://q.uiver.app/#q=WzAsNixbMCwxLCJcXGRibHtFfV8wIl0sWzEsMCwiXFxkYmx7RX1fMCJdLFsyLDEsIlxcZGJse0V9XzBcXHRpbWVzX3tcXGRibHtCfV8wfVxcZGJse0J9XzEiXSxbMSwyLCJcXGRibHtFfV8xIl0sWzMsMSwiXFxkYmx7RX1fMCJdLFsyLDIsIlxcZGJse0V9XzAiXSxbMCwxLCIxIl0sWzEsMiwiXFxsYW5nbGUxLFxcaWQgUFxccmFuZ2xlIiwyXSxbMCwzLCJcXGlkIiwyXSxbMywyLCJcXGxhbmdsZVxcc3JjLFBfMVxccmFuZ2xlIl0sWzIsNCwiXFx0Z3RcXHNpZ21hIl0sWzEsNCwiMSIsMCx7ImN1cnZlIjotMn1dLFszLDUsIlxcdGd0IiwyXSxbNSw0LCIxIiwyXSxbMiw1LCJcXHRndFxcZXBzaWxvbiIsMSx7InN0eWxlIjp7ImJvZHkiOnsibmFtZSI6Im5vbmUifSwiaGVhZCI6eyJuYW1lIjoibm9uZSJ9fX1dLFsxLDMsIlBfeygtKX0iLDEseyJzdHlsZSI6eyJib2R5Ijp7Im5hbWUiOiJub25lIn0sImhlYWQiOnsibmFtZSI6Im5vbmUifX19XSxbMTEsMiwiXFxrYXBwYSIsMSx7InNob3J0ZW4iOnsic291cmNlIjoyMH0sInN0eWxlIjp7ImJvZHkiOnsibmFtZSI6Im5vbmUifSwiaGVhZCI6eyJuYW1lIjoibm9uZSJ9fX1dXQ==
          \begin{tikzcd}
            & {\dbl{E}_0} && \\
            {\dbl{E}_0} && {\dbl{E}_0\times_{\dbl{B}_0}\dbl{B}_1} & {\dbl{E}_0} \\
            & {\dbl{E}_1} & {\dbl{E}_0}
            \arrow["{\langle1,\id P\rangle}"', from=1-2, to=2-3]
            \arrow[""{name=0, anchor=center, inner sep=0}, "1", curve={height=-12pt}, from=1-2, to=2-4]
            \arrow["{P_{(-)}}"{description}, draw=none, from=1-2, to=3-2]
            \arrow["1", from=2-1, to=1-2]
            \arrow["\id"', from=2-1, to=3-2]
            \arrow["{\tgt\sigma}", from=2-3, to=2-4]
            \arrow["{\tgt\epsilon}"{description}, draw=none, from=2-3, to=3-3]
            \arrow["{\langle\src,P_1\rangle}", from=3-2, to=2-3]
            \arrow["\tgt"', from=3-2, to=3-3]
            \arrow["1"', from=3-3, to=2-4]
            \arrow["\kappa"{description}, draw=none, from=0, to=2-3]
          \end{tikzcd}
          \qquad=\qquad
          1_{\tgt\id}
        \end{equation*}
      where $\kappa$ is the algebra structure unit isocell. 
    \end{itemize}
  A loosely discrete opfibration equipped with such a cleavage is said to be 
  \textbf{cloven}. 
\end{definition}

\begin{example}
  The projection $a\sslash\dbl{B}\to\dbl{B}$ from the double coslice is thus a 
  cloven loosely discrete opfibration.
\end{example}

\begin{remark}[Category of cloven loosely discrete opfibrations]
  A \define{morphism} between cloven loosely discrete opfibrations is simply a
  morphism as in \cref{def:cat-of-opfibrations}, that is, we do not ask that
  such morphisms be cleavage-preserving. The reason will become apparent in
  subsequent developments. In essence, the elements construction from
  \cref{subsection:elements} below can be extended to transformations between
  twisted functors (\cref{def:naturaltransformation}), but the induced morphisms
  of opfibrations are not cleavage-preserving
  (\cref{rmk:on-cleavage-preservation}).
\end{remark}

There are two facets of the definition to explore in the remainder of the subsection. 
First, we would like to understand the consequences of the adjunction equations and 
how this yields a lifting property for chosen opcartesian proarrows. Secondly, we 
would like a clearer unpacking of the algebra conditions. We study the adjunction 
equations first.

\begin{lemma} \label{lemma:counit-equations}
  For any cloven loosely discrete opfibration $P\colon\dbl{E}\to\dbl{B}$ with 
  cleavage $\sigma$, the counit $\epsilon$ of the accompanying adjunction 
  satisfies the two equations:
    \begin{enumerate}
      \item $P\epsilon_u=1$ for any proarrow $u\colon a\proto b$ of $\dbl{E}$.
      \item $\epsilon_{\sigma(a,m)} = 1$ for any object $a\in\dbl{E}$ over $x$ and
        $m\colon x\proto y$ of $\dbl{B}$.
    \end{enumerate}
  Consequently, targets of the components of $\epsilon$ are $P$-vertical and 
  targets of cleavage-indexed components of $\epsilon$ are tight identities in 
  $\dbl{E}$.
\end{lemma}
\begin{proof}
  That the cleavage $\sigma$ is left-adjoint-right-inverse to $\langle\src,P_1\rangle$ 
  implies that two equations hold, namely,
    \begin{enumerate}
      \item $\src(\sigma(a,m)) = a$
      \item $P(\sigma(a,m)) = m$,
    \end{enumerate}
  as one would expect for any such cleavage. That is, $\sigma(-,=)$ commutes with the 
  two legs of the spans in the diagram of \cref{def:cleavage-for-disc-opfibration}.
  In particular, the second equation says that $\sigma(a,m)$ is above the
  proarrow $m$ which it is supposed to lift to $\dbl{E}$. The first equation 
  says that the lift has the correct source. That $\sigma$ is a left adjoint in 
  an adjoint equivalence says that every proarrow $u\colon a\proto b$ in $\dbl{E}$ is
  left-sesquiglobularly isomorphic to the chosen proarrow $\sigma(a,Pu)$ via the
  corresponding component of the counit:
    \begin{equation*}
      % https://q.uiver.app/#q=WzAsNCxbMCwwLCJhIl0sWzEsMCwibV8hKGEpIl0sWzEsMSwiYiJdLFswLDEsImEiXSxbMCwxLCJcXHNpZ21hKGEsUHUpIiwwLHsic3R5bGUiOnsiYm9keSI6eyJuYW1lIjoiYmFycmVkIn19fV0sWzEsMiwiXFx0Z3RcXGVwc2lsb25fdSJdLFswLDMsIiIsMix7ImxldmVsIjoyLCJzdHlsZSI6eyJoZWFkIjp7Im5hbWUiOiJub25lIn19fV0sWzMsMiwidSIsMix7InN0eWxlIjp7ImJvZHkiOnsibmFtZSI6ImJhcnJlZCJ9fX1dLFs0LDcsIlxcZXBzaWxvbl91IiwxLHsic2hvcnRlbiI6eyJzb3VyY2UiOjIwLCJ0YXJnZXQiOjIwfSwic3R5bGUiOnsiYm9keSI6eyJuYW1lIjoibm9uZSJ9LCJoZWFkIjp7Im5hbWUiOiJub25lIn19fV1d
      \begin{tikzcd}
        a & {m_!(a)} \\
        a & b
        \arrow[""{name=0, anchor=center, inner sep=0}, "{\sigma(a,Pu)}"{inner sep=.8ex}, "\shortmid"{marking}, from=1-1, to=1-2]
        \arrow[equals, from=1-1, to=2-1]
        \arrow["{\tgt\epsilon_u}", from=1-2, to=2-2]
        \arrow[""{name=1, anchor=center, inner sep=0}, "u"'{inner sep=.8ex}, "\shortmid"{marking}, from=2-1, to=2-2]
        \arrow["{\epsilon_u}"{description}, draw=none, from=0, to=1]
      \end{tikzcd}.
    \end{equation*}
  Note that each such cell is invertible; in particular, the target morphism is 
  invertible. As for the triangle equations, on the one hand, the composite of 
  $\epsilon$ followed by the pairing of $\src$ and $P_1$ should be the identity:
    \begin{equation*}
      % https://q.uiver.app/#q=WzAsNCxbMCwwLCJcXGRibHtFfV8xIl0sWzEsMCwiXFxkYmx7RX1fMFxcdGltZXNfe1xcZGJse0J9XzB9XFxkYmx7Qn1fMSJdLFsxLDEsIlxcZGJse0V9XzEiXSxbMiwxLCJcXGRibHtFfV8wXFx0aW1lc197XFxkYmx7Qn1fMH1cXGRibHtCfV8xIl0sWzAsMSwiXFxsYW5nbGVcXHNyYyxQXzFcXHJhbmdsZSJdLFsxLDIsIlxcc2lnbWEiXSxbMiwzLCJcXGxhbmdsZVxcc3JjLFBfMVxccmFuZ2xlIiwyXSxbMSwzLCIiLDAseyJjdXJ2ZSI6LTEsImxldmVsIjoyLCJzdHlsZSI6eyJoZWFkIjp7Im5hbWUiOiJub25lIn19fV0sWzAsMiwiIiwwLHsiY3VydmUiOjEsImxldmVsIjoyLCJzdHlsZSI6eyJoZWFkIjp7Im5hbWUiOiJub25lIn19fV0sWzEsOCwiXFxlcHNpbG9uIiwwLHsibGFiZWxfcG9zaXRpb24iOjMwLCJzaG9ydGVuIjp7InNvdXJjZSI6MTAsInRhcmdldCI6MzB9fV1d
      \begin{tikzcd}
        {\dbl{E}_1} & {\dbl{E}_0\times_{\dbl{B}_0}\dbl{B}_1} \\
        & {\dbl{E}_1} & {\dbl{E}_0\times_{\dbl{B}_0}\dbl{B}_1}
        \arrow["{\langle\src,P_1\rangle}", from=1-1, to=1-2]
        \arrow[""{name=0, anchor=center, inner sep=0}, curve={height=6pt}, equals, from=1-1, to=2-2]
        \arrow["\sigma", from=1-2, to=2-2]
        \arrow[curve={height=-6pt}, equals, from=1-2, to=2-3]
        \arrow["{\langle\src,P_1\rangle}"', from=2-2, to=2-3]
        \arrow["\epsilon"{pos=0.3}, between={0.1}{0.7}, Rightarrow, from=1-2, to=0]
      \end{tikzcd} \qquad = \qquad 1_{\langle\src,P_1\rangle}.
    \end{equation*}
  And on the other hand, the cleavage $\sigma$ followed by $\epsilon$ should be 
  the identity:
    \begin{equation*}
      % https://q.uiver.app/#q=WzAsNCxbMCwxLCJcXGRibHtFfV8xIl0sWzEsMSwiXFxkYmx7RX1fMFxcdGltZXNfe1xcZGJse0J9XzB9XFxkYmx7Qn1fMSJdLFsxLDIsIlxcZGJse0V9XzEiXSxbMCwwLCJcXGRibHtFfV8wXFx0aW1lc197XFxkYmx7Qn1fMH1cXGRibHtCfV8xIl0sWzAsMSwiXFxsYW5nbGVcXHNyYyxQXzFcXHJhbmdsZSJdLFsxLDIsIlxcc2lnbWEiXSxbMCwyLCIiLDAseyJjdXJ2ZSI6MSwibGV2ZWwiOjIsInN0eWxlIjp7ImhlYWQiOnsibmFtZSI6Im5vbmUifX19XSxbMywwLCJcXHNpZ21hIiwyXSxbMywxLCIiLDAseyJjdXJ2ZSI6LTEsImxldmVsIjoyLCJzdHlsZSI6eyJoZWFkIjp7Im5hbWUiOiJub25lIn19fV0sWzEsNiwiXFxlcHNpbG9uIiwwLHsibGFiZWxfcG9zaXRpb24iOjMwLCJzaG9ydGVuIjp7InNvdXJjZSI6MTAsInRhcmdldCI6MzB9fV1d
      \begin{tikzcd}
        {\dbl{E}_0\times_{\dbl{B}_0}\dbl{B}_1} \\
        {\dbl{E}_1} & {\dbl{E}_0\times_{\dbl{B}_0}\dbl{B}_1} \\
        & {\dbl{E}_1}
        \arrow["\sigma"', from=1-1, to=2-1]
        \arrow[curve={height=-6pt}, equals, from=1-1, to=2-2]
        \arrow["{\langle\src,P_1\rangle}", from=2-1, to=2-2]
        \arrow[""{name=0, anchor=center, inner sep=0}, curve={height=6pt}, equals, from=2-1, to=3-2]
        \arrow["\sigma", from=2-2, to=3-2]
        \arrow["\epsilon"{pos=0.3}, between={0.1}{0.7}, Rightarrow, from=2-2, to=0]
      \end{tikzcd}\qquad = \qquad 1_\sigma.
    \end{equation*}
  Unpacking the former, we have that each component $\epsilon_u$ as above sits 
  over an identity cell in $\dbl{B}$ via $P$ and so in particular 
  $P(\tgt\epsilon_u)$ is a tight identity arrow of $\dbl{B}$. That is, 
  $\tgt\epsilon_u$ belongs to the fiber of $P$ over $b$. In the
  traditional language of fibrations, this is to say that $\tgt\epsilon_u$ is a 
  \emph{vertical} arrow relative to $P$. Likewise, the second equation shows 
  that the cleavage-indexed components of $\epsilon$ are exactly identity cells 
  of $\dbl{E}$ itself in that $\epsilon_{\sigma(a,m)} = 1$ and thus also the 
  targets of such cells are tight identities of $\dbl{E}$. 
\end{proof}

As a consequence, we can see that each chosen proarrow is \emph{opcartesian} in 
the following sense. This property will be used without comment throughout 
\cref{subsection:pseudo-inverse-construction}.

\begin{lemma}[Lifting property] \label{lemma:lifting-property}
  Any niche in $\dbl{E}$ of the form
    \begin{equation*}
      % https://q.uiver.app/#q=WzAsNCxbMCwwLCJhIl0sWzEsMCwibV8hKGEpIl0sWzAsMSwiYiJdLFsxLDEsImMiXSxbMCwxLCJcXHNpZ21hKGEsbSkiLDAseyJzdHlsZSI6eyJib2R5Ijp7Im5hbWUiOiJiYXJyZWQifX19XSxbMiwzLCJ1IiwyLHsic3R5bGUiOnsiYm9keSI6eyJuYW1lIjoiYmFycmVkIn19fV0sWzAsMiwicyIsMl1d
      \begin{tikzcd}
        a & {m_!(a)} \\
        b & c
        \arrow["{\sigma(a,m)}"{inner sep=.8ex}, "\shortmid"{marking}, from=1-1, to=1-2]
        \arrow["s"', from=1-1, to=2-1]
        \arrow["u"'{inner sep=.8ex}, "\shortmid"{marking}, from=2-1, to=2-2]
      \end{tikzcd}
    \end{equation*}
  over a cell of $\dbl{B}$ has a unique $P$-lift to a cell in $\dbl{E}$:
    \begin{equation*}
      % https://q.uiver.app/#q=WzAsNCxbMCwwLCJhIl0sWzEsMCwibV8hKGEpIl0sWzAsMSwiYiJdLFsxLDEsImMiXSxbMCwxLCJcXHNpZ21hKGEsbSkiLDAseyJzdHlsZSI6eyJib2R5Ijp7Im5hbWUiOiJiYXJyZWQifX19XSxbMiwzLCJ1IiwyLHsic3R5bGUiOnsiYm9keSI6eyJuYW1lIjoiYmFycmVkIn19fV0sWzAsMiwicyIsMl0sWzEsMywiIiwwLHsic3R5bGUiOnsiYm9keSI6eyJuYW1lIjoiZGFzaGVkIn19fV0sWzQsNSwiXFxleGlzdHNcXCwhIiwxLHsic2hvcnRlbiI6eyJzb3VyY2UiOjIwLCJ0YXJnZXQiOjIwfSwic3R5bGUiOnsiYm9keSI6eyJuYW1lIjoibm9uZSJ9LCJoZWFkIjp7Im5hbWUiOiJub25lIn19fV1d
      \begin{tikzcd}
        a & {m_!(a)} \\
        b & c
        \arrow[""{name=0, anchor=center, inner sep=0}, "{\sigma(a,m)}"{inner sep=.8ex}, "\shortmid"{marking}, from=1-1, to=1-2]
        \arrow["s"', from=1-1, to=2-1]
        \arrow[dashed, from=1-2, to=2-2]
        \arrow[""{name=1, anchor=center, inner sep=0}, "u"'{inner sep=.8ex}, "\shortmid"{marking}, from=2-1, to=2-2]
        \arrow["{\exists\,!}"{description}, draw=none, from=0, to=1]
      \end{tikzcd}
      \qquad\xmapsto{P}\qquad
      % https://q.uiver.app/#q=WzAsNCxbMCwwLCJ4Il0sWzEsMCwieSJdLFswLDEsIlBiIl0sWzEsMSwiUGMiXSxbMCwxLCJtIiwwLHsic3R5bGUiOnsiYm9keSI6eyJuYW1lIjoiYmFycmVkIn19fV0sWzIsMywiUHUiLDIseyJzdHlsZSI6eyJib2R5Ijp7Im5hbWUiOiJiYXJyZWQifX19XSxbMCwyLCJQcyIsMl0sWzEsMywiZyJdLFs0LDUsIlxcdGhldGEiLDEseyJzaG9ydGVuIjp7InNvdXJjZSI6MjAsInRhcmdldCI6MjB9LCJzdHlsZSI6eyJib2R5Ijp7Im5hbWUiOiJub25lIn0sImhlYWQiOnsibmFtZSI6Im5vbmUifX19XV0=
      \begin{tikzcd}
        x & y \\
        Pb & Pc
        \arrow[""{name=0, anchor=center, inner sep=0}, "m"{inner sep=.8ex}, "\shortmid"{marking}, from=1-1, to=1-2]
        \arrow["Ps"', from=1-1, to=2-1]
        \arrow["g", from=1-2, to=2-2]
        \arrow[""{name=1, anchor=center, inner sep=0}, "Pu"'{inner sep=.8ex}, "\shortmid"{marking}, from=2-1, to=2-2]
        \arrow["\theta"{description}, draw=none, from=0, to=1]
      \end{tikzcd}.
    \end{equation*}
\end{lemma}
\begin{proof}
  The lifting property is proved by a more detailed examination of the LARI
  equivalence, furnishing isomorphisms
    \begin{equation*}
      \dbl{E}_1(\sigma(a,m),u)\cong(\dbl{E}_0\times_{\dbl{B}_0}\dbl{B}_1)((a,m),(b,Pu)),
    \end{equation*}
  natural in objects $a$ over $x$, proarrows $m: x \proto y$ in $\dbl{B}$, and proarrows $u: b \proto c$ in $\dbl{E}$.
  Starting with a morphism on the right, that is, an arrow $s$ of $\dbl{E}$ and a cell $\theta$ of $\dbl{B}$ of the form
    \begin{equation*}
      % https://q.uiver.app/#q=WzAsMixbMCwwLCJhIl0sWzAsMSwiYiJdLFswLDEsInMiXV0=
      \begin{tikzcd}
        a \\
        b
        \arrow["s", from=1-1, to=2-1]
      \end{tikzcd}
      \qquad\qquad
      % https://q.uiver.app/#q=WzAsNCxbMCwwLCJ4Il0sWzEsMCwieSJdLFsxLDEsIlBjIl0sWzAsMSwiUGIiXSxbMSwyLCJnIl0sWzAsMywiZiIsMl0sWzMsMiwiUHUiLDIseyJzdHlsZSI6eyJib2R5Ijp7Im5hbWUiOiJiYXJyZWQifX19XSxbMCwxLCJtIiwwLHsic3R5bGUiOnsiYm9keSI6eyJuYW1lIjoiYmFycmVkIn19fV0sWzcsNiwiXFx0aGV0YSIsMSx7InNob3J0ZW4iOnsic291cmNlIjoyMCwidGFyZ2V0IjoyMH0sInN0eWxlIjp7ImJvZHkiOnsibmFtZSI6Im5vbmUifSwiaGVhZCI6eyJuYW1lIjoibm9uZSJ9fX1dXQ==
      \begin{tikzcd}
        x & y \\
        Pb & Pc
        \arrow[""{name=0, anchor=center, inner sep=0}, "m"{inner sep=.8ex}, "\shortmid"{marking}, from=1-1, to=1-2]
        \arrow["f"', from=1-1, to=2-1]
        \arrow["g", from=1-2, to=2-2]
        \arrow[""{name=1, anchor=center, inner sep=0}, "Pu"'{inner sep=.8ex}, "\shortmid"{marking}, from=2-1, to=2-2]
        \arrow["\theta"{description}, draw=none, from=0, to=1]
      \end{tikzcd},
    \end{equation*}
  we apply the morphism underlying the bijection, namely, first apply $\sigma$ and then post-compose with the counit of the adjunction as in:
    \begin{equation*}
      % https://q.uiver.app/#q=WzAsNixbMCwwLCJhIl0sWzEsMCwibV8hKGEpIl0sWzEsMSwiUCh1KV8hKGIpIl0sWzAsMSwiYiJdLFswLDIsImIiXSxbMSwyLCJjIl0sWzAsMSwiXFxzaWdtYShhLG0pIiwwLHsic3R5bGUiOnsiYm9keSI6eyJuYW1lIjoiYmFycmVkIn19fV0sWzEsMiwiXFx0aGV0YV8hKHMpIl0sWzAsMywicyIsMl0sWzMsMiwiXFxzaWdtYShiLFB1KSIsMix7InN0eWxlIjp7ImJvZHkiOnsibmFtZSI6ImJhcnJlZCJ9fX1dLFszLDQsIiIsMix7ImxldmVsIjoyLCJzdHlsZSI6eyJoZWFkIjp7Im5hbWUiOiJub25lIn19fV0sWzQsNSwidSIsMix7InN0eWxlIjp7ImJvZHkiOnsibmFtZSI6ImJhcnJlZCJ9fX1dLFsyLDUsIlxcdGd0XFxlcHNpbG9uX3UiXSxbNiw5LCJcXHNpZ21hKHMsXFx0aGV0YSkiLDEseyJzaG9ydGVuIjp7InNvdXJjZSI6MjAsInRhcmdldCI6MjB9LCJzdHlsZSI6eyJib2R5Ijp7Im5hbWUiOiJub25lIn0sImhlYWQiOnsibmFtZSI6Im5vbmUifX19XSxbOSwxMSwiXFxlcHNpbG9uX3UiLDEseyJzaG9ydGVuIjp7InNvdXJjZSI6MjAsInRhcmdldCI6MjB9LCJzdHlsZSI6eyJib2R5Ijp7Im5hbWUiOiJub25lIn0sImhlYWQiOnsibmFtZSI6Im5vbmUifX19XV0=
      \begin{tikzcd}
        a & {m_!(a)} \\
        b & {P(u)_!(b)} \\
        b & c
        \arrow[""{name=0, anchor=center, inner sep=0}, "{\sigma(a,m)}"{inner sep=.8ex}, "\shortmid"{marking}, from=1-1, to=1-2]
        \arrow["s"', from=1-1, to=2-1]
        \arrow["{\theta_!(s)}", from=1-2, to=2-2]
        \arrow[""{name=1, anchor=center, inner sep=0}, "{\sigma(b,Pu)}"'{inner sep=.8ex}, "\shortmid"{marking}, from=2-1, to=2-2]
        \arrow[equals, from=2-1, to=3-1]
        \arrow["{\tgt\epsilon_u}", from=2-2, to=3-2]
        \arrow[""{name=2, anchor=center, inner sep=0}, "u"'{inner sep=.8ex}, "\shortmid"{marking}, from=3-1, to=3-2]
        \arrow["{\sigma(s,\theta)}"{description}, draw=none, from=0, to=1]
        \arrow["{\epsilon_u}"{description}, draw=none, from=1, to=2]
      \end{tikzcd}.
  \end{equation*}
  This cell fills the niche. It is both a lift of $\theta$ and the unique one via
  $P$, as the inverse correspondence of the bijection above simply applies $P$
  since the unit of the adjunction is an identity.
\end{proof}

\begin{remark}
  The usual opcartesian lifting property for an ordinary opfibration is that of the existence of a unique morphism completing a triangle in the over category and above that in the base category by pre-composition with the chosen opcartesian arrow:
    \begin{equation*}
      % https://q.uiver.app/#q=WzAsMyxbMCwxLCJhIl0sWzIsMSwiZl8hKGEpIl0sWzIsMCwiYyJdLFswLDEsIlxcc2lnbWEoYSxmKSIsMl0sWzAsMiwibCJdLFsxLDIsIlxcaGF0IGciLDIseyJzdHlsZSI6eyJib2R5Ijp7Im5hbWUiOiJkYXNoZWQifX19XV0=
      \begin{tikzcd}
        && c \\
        a && {f_!(a)}
        \arrow["l", from=2-1, to=1-3]
        \arrow["{\sigma(a,f)}"', from=2-1, to=2-3]
        \arrow["{\hat g}"', dashed, from=2-3, to=1-3]
      \end{tikzcd} \qquad \xmapsto{P} \qquad 
      % https://q.uiver.app/#q=WzAsMyxbMCwxLCJ4Il0sWzIsMSwieSJdLFsyLDAsInoiXSxbMCwxLCJmIiwyXSxbMCwyLCJQKGwpIl0sWzEsMiwiZyIsMl1d
      \begin{tikzcd}
        && z \\
        x && y
        \arrow["{P(l)}", from=2-1, to=1-3]
        \arrow["f"', from=2-1, to=2-3]
        \arrow["g"', from=2-3, to=1-3]
      \end{tikzcd}
    \end{equation*}
  This property (and its higher-dimensional analogues) is used in the construction
  of (higher-)category representations via the filling of squares rather than 
  triangles:
    \begin{equation*}
      % https://q.uiver.app/#q=WzAsNCxbMCwxLCJhIl0sWzIsMSwiZl8hKGEpIl0sWzIsMCwiYyJdLFswLDAsImQiXSxbMCwxLCJcXHNpZ21hKGEsZikiLDJdLFsxLDIsIlxcaGF0IGciLDIseyJzdHlsZSI6eyJib2R5Ijp7Im5hbWUiOiJkYXNoZWQifX19XSxbMywyLCJsIl0sWzAsMywiciJdXQ==
      \begin{tikzcd}
        d && c \\
        a && {f_!(a)}
        \arrow["l", from=1-1, to=1-3]
        \arrow["r", from=2-1, to=1-1]
        \arrow["{\sigma(a,f)}"', from=2-1, to=2-3]
        \arrow["{\hat g}"', dashed, from=2-3, to=1-3]
      \end{tikzcd} \qquad \xmapsto{P} \qquad
      % https://q.uiver.app/#q=WzAsNCxbMCwxLCJ4Il0sWzIsMSwieSJdLFsyLDAsInoiXSxbMCwwLCJ3Il0sWzAsMSwiZiIsMl0sWzEsMiwiZyIsMl0sWzMsMiwiUChsKSJdLFswLDMsIlAocikiXV0=
      \begin{tikzcd}
        w && z \\
        x && y
        \arrow["{P(l)}", from=1-1, to=1-3]
        \arrow["{P(r)}", from=2-1, to=1-1]
        \arrow["f"', from=2-1, to=2-3]
        \arrow["g"', from=2-3, to=1-3]
      \end{tikzcd}
    \end{equation*}
  In this way, we see that the lifting property of \cref{lemma:lifting-property} 
  is an analogous \emph{square-filling} property that we have simply chosen to 
  present in a different orientation. And indeed the pseudo-inverse construction 
  (\cref{subsection:pseudo-inverse-construction}) uses precisely this unique lifting
  condition via square filling to give transition functors across fiber categories 
  in a manner analogous to the classical construction.
\end{remark}

To conclude the subsection, we now give a statement and unpacking of the algebra
conditions in the definition of a cleavage.

\begin{remark}[Pseudo-algebra equations]
  \label{remark:cleavage-pseudo-algebra-relations}
  The pseudo algebra $\dbl{E}\ltimes\dbl{B}$ from 
  \cref{def:cleavage-for-disc-opfibration} above can be described in more detail 
  in the following way. There are two categories
    \begin{enumerate}
      \item objects and arrows: $(\dbl{E}\ltimes\dbl{B})_0:=\dbl{E}_0$
      \item proarrows and cells: $(\dbl{E}\ltimes\dbl{B})_1:=\dbl{E}_0
      \times_{\dbl{B}_0}\dbl{B}_1$
    \end{enumerate}
  where the latter category is the corner of the pullback of
  $P_0: \dbl{E}_0 \to \dbl{B}_0$ along $\src: \dbl{B}_1 \to \dbl{B}_0$.
  The structure morphism is the composite $\tgt\sigma\colon\dbl{E}_0
  \times_{\dbl{B}_0}\dbl{B}_1\to\dbl{E}_0$. The comparison cells are then
    \begin{equation*}
      % https://q.uiver.app/#q=WzAsNCxbMSwwLCJcXGRibHtFfV8wXFx0aW1lc197XFxkYmx7Qn1fMH1cXGRibHtCfV8xIl0sWzEsMSwiXFxkYmx7RX1fMCJdLFswLDEsIlxcZGJse0V9XzBcXHRpbWVzX3tcXGRibHtCfV8wfVxcZGJse0J9XzEiXSxbMCwwLCJcXGRibHtFfV8wXFx0aW1lc197XFxkYmx7Qn1fMH1cXGRibHtCfV8xXFx0aW1lc197XFxkYmx7Qn1fMH1cXGRibHtCfV8xIl0sWzAsMSwiXFx0Z3QuXFxzaWdtYSJdLFsyLDEsIlxcdGd0Llxcc2lnbWEiLDJdLFszLDAsIlxcdGd0Llxcc2lnbWFcXHRpbWVzIDEiXSxbMywyLCIxXFx0aW1lcy1cXG9kb3QtIiwyXSxbNiw1LCJcXG11IiwxLHsic2hvcnRlbiI6eyJzb3VyY2UiOjIwLCJ0YXJnZXQiOjIwfSwic3R5bGUiOnsiYm9keSI6eyJuYW1lIjoibm9uZSJ9LCJoZWFkIjp7Im5hbWUiOiJub25lIn19fV1d
      \begin{tikzcd}
        {\dbl{E}_0\times_{\dbl{B}_0}\dbl{B}_1\times_{\dbl{B}_0}\dbl{B}_1} & {\dbl{E}_0\times_{\dbl{B}_0}\dbl{B}_1} \\
        {\dbl{E}_0\times_{\dbl{B}_0}\dbl{B}_1} & {\dbl{E}_0}
        \arrow[""{name=0, anchor=center, inner sep=0}, "{\tgt.\sigma\times 1}", from=1-1, to=1-2]
        \arrow["{1\times-\odot-}"', from=1-1, to=2-1]
        \arrow["{\tgt.\sigma}", from=1-2, to=2-2]
        \arrow[""{name=1, anchor=center, inner sep=0}, "{\tgt.\sigma}"', from=2-1, to=2-2]
        \arrow["\mu"{description}, draw=none, from=0, to=1]
      \end{tikzcd}\qquad\qquad
      % https://q.uiver.app/#q=WzAsMyxbMSwwLCJcXGRibHtFfV8wXFx0aW1lc197XFxkYmx7Qn1fMH1cXGRibHtCfV8xIl0sWzEsMSwiXFxkYmx7RX1fMCJdLFswLDAsIlxcZGJse0V9XzAiXSxbMCwxLCJcXHRndC5cXHNpZ21hIl0sWzIsMCwiXFxsYW5nbGUgMSxcXGlkLlBfMFxccmFuZ2xlIl0sWzIsMSwiIiwyLHsibGV2ZWwiOjIsInN0eWxlIjp7ImhlYWQiOnsibmFtZSI6Im5vbmUifX19XSxbMCw1LCJcXGthcHBhIiwxLHsic2hvcnRlbiI6eyJ0YXJnZXQiOjIwfSwic3R5bGUiOnsiYm9keSI6eyJuYW1lIjoibm9uZSJ9LCJoZWFkIjp7Im5hbWUiOiJub25lIn19fV1d
      \begin{tikzcd}
        {\dbl{E}_0} & {\dbl{E}_0\times_{\dbl{B}_0}\dbl{B}_1} \\
        & {\dbl{E}_0}
        \arrow["{\langle 1,\id.P_0\rangle}", from=1-1, to=1-2]
        \arrow[""{name=0, anchor=center, inner sep=0}, equals, from=1-1, to=2-2]
        \arrow["{\tgt.\sigma}", from=1-2, to=2-2]
        \arrow["\kappa"{description}, draw=none, from=1-2, to=0]
      \end{tikzcd}.
    \end{equation*}
  A component of $\mu$ in particular is then an isomorphism of $\dbl{E}$ of the 
  form 
    \begin{equation*}
      \mu_{a,m,n}\colon n_!m_!(a)\xrightarrow{\cong}(m\odot n)_!(a).
    \end{equation*}
  To state the associativity condition for $\mu$, we denote the corner objects 
  of the two- and three-fold actions of the pull-push monad on $P_0$ by 
  $\dbl{E}_0\times_{\dbl{B}_0}\dbl{B}_1^{(2)}$ and $\dbl{E}_0\times_{\dbl{B}_0}
  \dbl{B}_1^{(3)}$ respectively. Thus, diagrammatically, the associativity 
  condition for $\mu$ is the statement that the composite diagram
    \begin{equation*}
      % https://q.uiver.app/#q=WzAsNyxbMCwwLCJcXGRibHtFfV8wXFx0aW1lc197XFxkYmx7Qn1fMH1cXGRibHtCfV8xXnsoMyl9Il0sWzEsMSwiXFxkYmx7RX1fMFxcdGltZXNfe1xcZGJse0J9XzB9XFxkYmx7Qn1fMV57KDIpfSJdLFsxLDAsIlxcZGJse0V9XzBcXHRpbWVzX3tcXGRibHtCfV8wfVxcZGJse0J9XzFeeygyKX0iXSxbMiwxLCJcXGRibHtFfV8wXFx0aW1lc197XFxkYmx7Qn1fMH1cXGRibHtCfV8xIl0sWzAsMSwiXFxkYmx7RX1fMFxcdGltZXNfe1xcZGJse0J9XzB9XFxkYmx7Qn1fMV57KDIpfSJdLFsxLDIsIlxcZGJse0V9XzBcXHRpbWVzX3tcXGRibHtCfV8wfVxcZGJse0J9XzEiXSxbMiwyLCJcXGRibHtFfV8wIl0sWzAsMSwiMVxcdGltZXNcXG9kb3RcXHRpbWVzIDEiLDFdLFswLDIsIlxcdGd0Llxcc2lnbWFcXHRpbWVzIDEiXSxbMiwzLCJcXHRndC5cXHNpZ21hIFxcdGltZXMgMSJdLFswLDQsIjFcXHRpbWVzIDFcXHRpbWVzIFxcb2RvdCIsMl0sWzQsNSwiMVxcdGltZXMgXFxvZG90IiwyXSxbMSw1LCIxXFx0aW1lc1xcb2RvdCJdLFsxLDMsIjFcXHRpbWVzIFxcdGd0Llxcc2lnbWEiXSxbMyw2LCJcXHRndC5cXHNpZ21hIl0sWzUsNiwiXFx0Z3QuXFxzaWdtYSIsMl0sWzMsNSwiXFxtdSIsMSx7InN0eWxlIjp7ImJvZHkiOnsibmFtZSI6Im5vbmUifSwiaGVhZCI6eyJuYW1lIjoibm9uZSJ9fX1dLFsyLDEsIlxcbXVcXHRpbWVzIDEiLDEseyJzdHlsZSI6eyJib2R5Ijp7Im5hbWUiOiJub25lIn0sImhlYWQiOnsibmFtZSI6Im5vbmUifX19XSxbNywxMSwiMVxcdGltZXNcXG1hdGhmcmFrIGEiLDEseyJzaG9ydGVuIjp7InNvdXJjZSI6MjAsInRhcmdldCI6MjB9LCJzdHlsZSI6eyJib2R5Ijp7Im5hbWUiOiJub25lIn0sImhlYWQiOnsibmFtZSI6Im5vbmUifX19XV0=
      \begin{tikzcd}
        {\dbl{E}_0\times_{\dbl{B}_0}\dbl{B}_1^{(3)}} & {\dbl{E}_0\times_{\dbl{B}_0}\dbl{B}_1^{(2)}} \\
        {\dbl{E}_0\times_{\dbl{B}_0}\dbl{B}_1^{(2)}} & {\dbl{E}_0\times_{\dbl{B}_0}\dbl{B}_1^{(2)}} & {\dbl{E}_0\times_{\dbl{B}_0}\dbl{B}_1} \\
        & {\dbl{E}_0\times_{\dbl{B}_0}\dbl{B}_1} & {\dbl{E}_0}
        \arrow["{\tgt.\sigma\times 1}", from=1-1, to=1-2]
        \arrow["{1\times 1\times \odot}"', from=1-1, to=2-1]
        \arrow[""{name=0, anchor=center, inner sep=0}, "{1\times\odot\times 1}"{description}, from=1-1, to=2-2]
        \arrow["{\mu\times 1}"{description}, draw=none, from=1-2, to=2-2]
        \arrow["{\tgt.\sigma \times 1}", from=1-2, to=2-3]
        \arrow[""{name=1, anchor=center, inner sep=0}, "{1\times \odot}"', from=2-1, to=3-2]
        \arrow["{1\times \tgt.\sigma}", from=2-2, to=2-3]
        \arrow["{1\times\odot}", from=2-2, to=3-2]
        \arrow["\mu"{description}, draw=none, from=2-3, to=3-2]
        \arrow["{\tgt.\sigma}", from=2-3, to=3-3]
        \arrow["{\tgt.\sigma}"', from=3-2, to=3-3]
        \arrow["{1\times\mathfrak a}"{description}, draw=none, from=0, to=1]
      \end{tikzcd}
    \end{equation*}
  is equal to the composite diagram
    \begin{equation*}
      % https://q.uiver.app/#q=WzAsNyxbMCwwLCJcXGRibHtFfV8wXFx0aW1lc197XFxkYmx7Qn1fMH1cXGRibHtCfV8xXnsoMyl9Il0sWzEsMSwiXFxkYmx7RX1fMFxcdGltZXNfe1xcZGJse0J9XzB9XFxkYmx7Qn1fMSJdLFsxLDAsIlxcZGJse0V9XzBcXHRpbWVzX3tcXGRibHtCfV8wfVxcZGJse0J9XzFeeygyKX0iXSxbMiwxLCJcXGRibHtFfV8wXFx0aW1lc197XFxkYmx7Qn1fMH1cXGRibHtCfV8xIl0sWzAsMSwiXFxkYmx7RX1fMFxcdGltZXNfe1xcZGJse0J9XzB9XFxkYmx7Qn1fMV57KDIpfSJdLFsxLDIsIlxcZGJse0V9XzBcXHRpbWVzX3tcXGRibHtCfV8wfVxcZGJse0J9XzEiXSxbMiwyLCJcXGRibHtFfV8wIl0sWzAsMiwiXFx0Z3QuXFxzaWdtYVxcdGltZXMgMSJdLFsyLDMsIlxcdGd0Llxcc2lnbWEgXFx0aW1lcyAxIl0sWzAsNCwiMVxcdGltZXMgMVxcdGltZXMgXFxvZG90IiwyXSxbNCw1LCIxXFx0aW1lcyBcXG9kb3QiLDJdLFszLDYsIlxcdGd0Llxcc2lnbWEiXSxbNSw2LCJcXHRndC5cXHNpZ21hIiwyXSxbNCwxLCJcXHRndC5cXHNpZ21hIFxcdGltZXMgMSJdLFsyLDEsIjFcXHRpbWVzIFxcb2RvdCIsMl0sWzEsNiwiXFx0Z3QuXFxzaWdtYSIsMV0sWzEsNSwiXFxtdSIsMSx7InN0eWxlIjp7ImJvZHkiOnsibmFtZSI6Im5vbmUifSwiaGVhZCI6eyJuYW1lIjoibm9uZSJ9fX1dLFs4LDE1LCJcXG11IiwxLHsic2hvcnRlbiI6eyJzb3VyY2UiOjIwLCJ0YXJnZXQiOjIwfSwic3R5bGUiOnsiYm9keSI6eyJuYW1lIjoibm9uZSJ9LCJoZWFkIjp7Im5hbWUiOiJub25lIn19fV1d
      \begin{tikzcd}
        {\dbl{E}_0\times_{\dbl{B}_0}\dbl{B}_1^{(3)}} & {\dbl{E}_0\times_{\dbl{B}_0}\dbl{B}_1^{(2)}} \\
        {\dbl{E}_0\times_{\dbl{B}_0}\dbl{B}_1^{(2)}} & {\dbl{E}_0\times_{\dbl{B}_0}\dbl{B}_1} & {\dbl{E}_0\times_{\dbl{B}_0}\dbl{B}_1} \\
        & {\dbl{E}_0\times_{\dbl{B}_0}\dbl{B}_1} & {\dbl{E}_0}
        \arrow["{\tgt.\sigma\times 1}", from=1-1, to=1-2]
        \arrow["{1\times 1\times \odot}"', from=1-1, to=2-1]
        \arrow["{1\times \odot}"', from=1-2, to=2-2]
        \arrow[""{name=0, anchor=center, inner sep=0}, "{\tgt.\sigma \times 1}", from=1-2, to=2-3]
        \arrow["{\tgt.\sigma \times 1}", from=2-1, to=2-2]
        \arrow["{1\times \odot}"', from=2-1, to=3-2]
        \arrow["\mu"{description}, draw=none, from=2-2, to=3-2]
        \arrow[""{name=1, anchor=center, inner sep=0}, "{\tgt.\sigma}"{description}, from=2-2, to=3-3]
        \arrow["{\tgt.\sigma}", from=2-3, to=3-3]
        \arrow["{\tgt.\sigma}"', from=3-2, to=3-3]
        \arrow["\mu"{description}, draw=none, from=0, to=1]
      \end{tikzcd}
    \end{equation*}
  where $1\times\mathfrak a$ is the structure cell for the monad given by the 
  appropriate component of the associator for loose composition in $\dbl{B}$.
  On a component, the associativity condition for $\mu$ is thus the 
  commutativity of the pentagonal figure
    \begin{equation*}
      % https://q.uiver.app/#q=WzAsNSxbMCwwLCJwXyFuXyFtXyEoYSkiXSxbMiwwLCJwXyEobVxcb2RvdCBuKV8hKGEpIl0sWzAsMSwiKG5cXG9kb3QgcClfIVxcb2RvdCBtXyEoYSkiXSxbMiwxLCIoKG1cXG9kb3QgbilcXG9kb3QgcClfIShhKSJdLFsxLDIsIihtXFxvZG90IChuXFxvZG90IHApKV8hKGEpIl0sWzAsMSwicF8hKFxcbXVfe2EsbSxufSkiXSxbMCwyLCJcXG11X3ttXyEoYSksbixwfSIsMl0sWzEsMywiXFxtdV97YSxtXFxvZG90IG4scH0iLDJdLFsyLDQsIlxcbXVfe2EsbSxuXFxvZG90IHB9IiwyXSxbMyw0LCJcXGNvbmciXV0=
      \begin{tikzcd}
        {p_!n_!m_!(a)} && {p_!(m\odot n)_!(a)} \\
        {(n\odot p)_!\odot m_!(a)} && {((m\odot n)\odot p)_!(a)} \\
        & {(m\odot (n\odot p))_!(a)}
        \arrow["{p_!(\mu_{a,m,n})}", from=1-1, to=1-3]
        \arrow["{\mu_{m_!(a),n,p}}"', from=1-1, to=2-1]
        \arrow["{\mu_{a,m\odot n,p}}", from=1-3, to=2-3]
        \arrow["{\mu_{a,m,n\odot p}}"', from=2-1, to=3-2]
        \arrow["\cong", from=2-3, to=3-2]
      \end{tikzcd}
    \end{equation*}
  hence our view that this is a codescent condition. Clockwise is given by 
  the left side of the preceding cell diagram and counterclockwise by the right.
  There likewise are two unit conditions. In diagrams these are 
    \begin{equation*}
      % https://q.uiver.app/#q=WzAsNixbMCwwLCJcXGRibHtFfV8wXFx0aW1lc197XFxkYmx7Qn1fMH1cXGRibHtCfV8xIl0sWzEsMCwiXFxkYmx7RX1fMCJdLFsxLDEsIlxcZGJse0V9XzBcXHRpbWVzX3tcXGRibHtCfV8wfVxcZGJse0J9XzEiXSxbMCwxLCJcXGRibHtFfV8wXFx0aW1lc197XFxkYmx7Qn1fMH1cXGRibHtCfV8xXFx0aW1lc197XFxkYmx7Qn1fMH1cXGRibHtCfV8xIl0sWzIsMSwiXFxkYmx7RX1fMCJdLFsxLDIsIlxcZGJse0V9XzBcXHRpbWVzX3tcXGRibHtCfV8wfVxcZGJse0J9XzEiXSxbMCwxLCJcXHRndC5cXHNpZ21hIl0sWzEsMiwiXFxsYW5nbGUxLFxcaWQuUFxccmFuZ2xlIiwyXSxbMCwzLCJcXGxhbmdsZSAxLFxcaWQuUFxccmFuZ2xlIiwyXSxbMywyLCJcXHRndC5cXHNpZ21hXFx0aW1lcyAxIl0sWzIsNCwiXFx0Z3QuXFxzaWdtYSJdLFszLDUsIjFcXHRpbWVzXFxvZG90IiwyXSxbNSw0LCJcXHRndC5cXHNpZ21hIiwyXSxbMSw0LCIiLDEseyJjdXJ2ZSI6LTIsImxldmVsIjoyLCJzdHlsZSI6eyJoZWFkIjp7Im5hbWUiOiJub25lIn19fV0sWzIsNSwiXFxtdSIsMSx7InN0eWxlIjp7ImJvZHkiOnsibmFtZSI6Im5vbmUifSwiaGVhZCI6eyJuYW1lIjoibm9uZSJ9fX1dLFsxMywyLCJcXGthcHBhIiwxLHsic2hvcnRlbiI6eyJzb3VyY2UiOjIwfSwic3R5bGUiOnsiYm9keSI6eyJuYW1lIjoibm9uZSJ9LCJoZWFkIjp7Im5hbWUiOiJub25lIn19fV1d
      \begin{tikzcd}
        {\dbl{E}_0\times_{\dbl{B}_0}\dbl{B}_1} & {\dbl{E}_0} \\
        {\dbl{E}_0\times_{\dbl{B}_0}\dbl{B}_1\times_{\dbl{B}_0}\dbl{B}_1} & {\dbl{E}_0\times_{\dbl{B}_0}\dbl{B}_1} & {\dbl{E}_0} \\
        & {\dbl{E}_0\times_{\dbl{B}_0}\dbl{B}_1}
        \arrow["{\tgt.\sigma}", from=1-1, to=1-2]
        \arrow["{\langle 1,\id.P\rangle}"', from=1-1, to=2-1]
        \arrow["{\langle1,\id.P\rangle}"', from=1-2, to=2-2]
        \arrow[""{name=0, anchor=center, inner sep=0}, curve={height=-12pt}, equals, from=1-2, to=2-3]
        \arrow["{\tgt.\sigma\times 1}", from=2-1, to=2-2]
        \arrow["{1\times\odot}"', from=2-1, to=3-2]
        \arrow["{\tgt.\sigma}", from=2-2, to=2-3]
        \arrow["\mu"{description}, draw=none, from=2-2, to=3-2]
        \arrow["{\tgt.\sigma}"', from=3-2, to=2-3]
        \arrow["\kappa"{description}, draw=none, from=0, to=2-2]
      \end{tikzcd}\;=
      % https://q.uiver.app/#q=WzAsNCxbMCwwLCJcXGRibHtFfV8wXFx0aW1lc197XFxkYmx7Qn1fMH1cXGRibHtCfV8xIl0sWzAsMSwiXFxkYmx7RX1fMFxcdGltZXNfe1xcZGJse0J9XzB9XFxkYmx7Qn1fMVxcdGltZXNfe1xcZGJse0J9XzB9XFxkYmx7Qn1fMSJdLFsyLDEsIlxcZGJse0V9XzAiXSxbMSwxLCJcXGRibHtFfV8wXFx0aW1lc197XFxkYmx7Qn1fMH1cXGRibHtCfV8xIl0sWzAsMSwiXFxsYW5nbGUgMSxcXGlkLlBcXHJhbmdsZSIsMl0sWzEsMywiMVxcdGltZXNcXG9kb3QiLDJdLFszLDIsIlxcdGd0Llxcc2lnbWEiLDJdLFswLDMsIiIsMix7ImN1cnZlIjotMSwibGV2ZWwiOjIsInN0eWxlIjp7ImhlYWQiOnsibmFtZSI6Im5vbmUifX19XSxbNywxLCJcXHJobyIsMSx7InNob3J0ZW4iOnsic291cmNlIjoyMH0sInN0eWxlIjp7ImJvZHkiOnsibmFtZSI6Im5vbmUifSwiaGVhZCI6eyJuYW1lIjoibm9uZSJ9fX1dXQ==
      \begin{tikzcd}
        {\dbl{E}_0\times_{\dbl{B}_0}\dbl{B}_1} \\
        {\dbl{E}_0\times_{\dbl{B}_0}\dbl{B}_1\times_{\dbl{B}_0}\dbl{B}_1} & {\dbl{E}_0\times_{\dbl{B}_0}\dbl{B}_1} & {\dbl{E}_0}
        \arrow["{\langle 1,\id.P\rangle}"', from=1-1, to=2-1]
        \arrow[""{name=0, anchor=center, inner sep=0}, curve={height=-6pt}, equals, from=1-1, to=2-2]
        \arrow["{1\times\odot}"', from=2-1, to=2-2]
        \arrow["{\tgt.\sigma}"', from=2-2, to=2-3]
        \arrow["\rho"{description}, draw=none, from=0, to=2-1]
      \end{tikzcd}
    \end{equation*}
  and 
    \begin{equation*}
      % https://q.uiver.app/#q=WzAsNSxbMCwxLCJcXGRibHtFfV8wXFx0aW1lc197XFxkYmx7Qn1fMH1cXGRibHtCfV8xIl0sWzEsMSwiXFxkYmx7RX1fMFxcdGltZXNfe1xcZGJse0J9XzB9XFxkYmx7Qn1fMVxcdGltZXNfe1xcZGJse0J9XzB9XFxkYmx7Qn1fMSJdLFsyLDIsIlxcZGJse0V9XzBcXHRpbWVzX3tcXGRibHtCfV8wfVxcZGJse0J9XzEiXSxbMiwwLCJcXGRibHtFfV8wXFx0aW1lc197XFxkYmx7Qn1fMH1cXGRibHtCfV8xIl0sWzMsMSwiXFxkYmx7RX1fMCJdLFswLDEsIlxcbGFuZ2xlIDEsXFxpZC5QXFxyYW5nbGVcXHRpbWVzIDEiLDJdLFsxLDIsIjFcXHRpbWVzXFxvZG90Il0sWzEsMywiXFx0Z3QuXFxzaWdtYVxcdGltZXMgMSIsMl0sWzMsNCwiXFx0Z3QuXFxzaWdtYSJdLFsyLDQsIlxcdGd0Llxcc2lnbWEiLDJdLFswLDMsIiIsMSx7ImN1cnZlIjotMiwibGV2ZWwiOjIsInN0eWxlIjp7ImhlYWQiOnsibmFtZSI6Im5vbmUifX19XSxbMCwyLCIiLDEseyJjdXJ2ZSI6MiwibGV2ZWwiOjIsInN0eWxlIjp7ImhlYWQiOnsibmFtZSI6Im5vbmUifX19XSxbMywyLCJcXG11IiwxLHsic3R5bGUiOnsiYm9keSI6eyJuYW1lIjoibm9uZSJ9LCJoZWFkIjp7Im5hbWUiOiJub25lIn19fV0sWzEwLDEsIlxca2FwcGFcXHRpbWVzIDEiLDEseyJzaG9ydGVuIjp7InNvdXJjZSI6MjB9LCJzdHlsZSI6eyJib2R5Ijp7Im5hbWUiOiJub25lIn0sImhlYWQiOnsibmFtZSI6Im5vbmUifX19XSxbMSwxMSwiXFxsYW1iZGEiLDEseyJzaG9ydGVuIjp7InRhcmdldCI6MjB9LCJzdHlsZSI6eyJib2R5Ijp7Im5hbWUiOiJub25lIn0sImhlYWQiOnsibmFtZSI6Im5vbmUifX19XV0=
      \begin{tikzcd}
        && {\dbl{E}_0\times_{\dbl{B}_0}\dbl{B}_1} \\
        {\dbl{E}_0\times_{\dbl{B}_0}\dbl{B}_1} & {\dbl{E}_0\times_{\dbl{B}_0}\dbl{B}_1\times_{\dbl{B}_0}\dbl{B}_1} && {\dbl{E}_0} \\
        && {\dbl{E}_0\times_{\dbl{B}_0}\dbl{B}_1}
        \arrow["{\tgt.\sigma}", from=1-3, to=2-4]
        \arrow["\mu"{description}, draw=none, from=1-3, to=3-3]
        \arrow[""{name=0, anchor=center, inner sep=0}, curve={height=-12pt}, equals, from=2-1, to=1-3]
        \arrow["{\langle 1,\id.P\rangle\times 1}"', from=2-1, to=2-2]
        \arrow[""{name=1, anchor=center, inner sep=0}, curve={height=12pt}, equals, from=2-1, to=3-3]
        \arrow["{\tgt.\sigma\times 1}"', from=2-2, to=1-3]
        \arrow["{1\times\odot}", from=2-2, to=3-3]
        \arrow["{\tgt.\sigma}"', from=3-3, to=2-4]
        \arrow["{\kappa\times 1}"{description}, draw=none, from=0, to=2-2]
        \arrow["\lambda"{description}, draw=none, from=2-2, to=1]
      \end{tikzcd} \qquad = \qquad 1_\mu
    \end{equation*}
  which the reader can unwind by chasing an object through each side of each 
  equality and comparing component arrows.
\end{remark}

\subsection{Elements construction}
\label{subsection:elements}

Famously, a copresheaf gives rise to a discrete opfibration via the category of
elements construction. Here we show that a twisted copresheaf, or
profunctor-valued twisted normal lax functor, gives rise to a cloven loosely
discrete opfibration via an elements construction. In fact, although our
examples of twisted functors will be normal, even unitary, that assumption is
not needed for the construction.

\begin{construction}[Elements construction] \label{construction:weak-elements}
  Let $F\colon \dbl{B}\twistto\Prof$ be a twisted lax functor.
  The \textbf{double category of elements} associated to $F$, denoted
  by $\Elt(F)$, has 
    \begin{enumerate}
      \item as objects, an object $x$ of $\dbl{B}$ together with an object $a$ of
        $Fx$;
      \item as arrows $(x,a)\to (z,c)$, an arrow $f\colon x\to z$ of $\dbl{B}$ together
        with a heteromorphism $\bar f$ in $Ff(a,c)$, where $Ff\colon Fx\proto Fz$ is
        the profunctor assigned to $f$ by $F$;
      \item as proarrows $(x,a)\proto (y,b)$, a proarrow $m\colon x\proto y$ of $\dbl{B}$
        together with an isomorphism $\overline m\colon Fm(a)\xrightarrow{\cong} b$ of
        $Fy$, where $Fm: Fx \to Fy$ is the functor assigned to $m$ by $F$;
      \item as cells
        \begin{equation*}
          % https://q.uiver.app/#q=WzAsNCxbMCwwLCIoeCxhKSJdLFsxLDAsIih5LGIpIl0sWzAsMSwiKHosYykiXSxbMSwxLCIodyxkKSJdLFswLDEsIihtLFxcb3ZlcmxpbmUgbSkiLDAseyJzdHlsZSI6eyJib2R5Ijp7Im5hbWUiOiJiYXJyZWQifX19XSxbMCwyLCIoZixcXG92ZXJsaW5lIGYpIiwyXSxbMiwzLCIobixcXG92ZXJsaW5lIG4pIiwyLHsic3R5bGUiOnsiYm9keSI6eyJuYW1lIjoiYmFycmVkIn19fV0sWzEsMywiKGcsXFxvdmVybGluZSBnKSJdLFs0LDYsIlxcYWxwaGEiLDEseyJzaG9ydGVuIjp7InNvdXJjZSI6MjAsInRhcmdldCI6MjB9LCJzdHlsZSI6eyJib2R5Ijp7Im5hbWUiOiJub25lIn0sImhlYWQiOnsibmFtZSI6Im5vbmUifX19XV0=
          \begin{tikzcd}
            {(x,a)} & {(y,b)} \\
            {(z,c)} & {(w,d)}
            \arrow[""{name=0, anchor=center, inner sep=0}, "{(m,\overline m)}", "\shortmid"{marking}, from=1-1, to=1-2]
            \arrow["{(f,\overline f)}"', from=1-1, to=2-1]
            \arrow["{(g,\overline g)}", from=1-2, to=2-2]
            \arrow[""{name=1, anchor=center, inner sep=0}, "{(n,\overline n)}"', "\shortmid"{marking}, from=2-1, to=2-2]
            \arrow["\alpha"{description}, draw=none, from=0, to=1]
          \end{tikzcd},
        \end{equation*}
      those cells of $\dbl{B}$ 
        \begin{equation*}
          % https://q.uiver.app/#q=WzAsNCxbMCwwLCJ4Il0sWzEsMCwieSJdLFswLDEsInoiXSxbMSwxLCJ3Il0sWzAsMSwibSIsMCx7InN0eWxlIjp7ImJvZHkiOnsibmFtZSI6ImJhcnJlZCJ9fX1dLFswLDIsImYiLDJdLFsyLDMsIm4iLDIseyJzdHlsZSI6eyJib2R5Ijp7Im5hbWUiOiJiYXJyZWQifX19XSxbMSwzLCJnIl0sWzQsNiwiXFxhbHBoYSIsMSx7InNob3J0ZW4iOnsic291cmNlIjoyMCwidGFyZ2V0IjoyMH0sInN0eWxlIjp7ImJvZHkiOnsibmFtZSI6Im5vbmUifSwiaGVhZCI6eyJuYW1lIjoibm9uZSJ9fX1dXQ==
          \begin{tikzcd}
            x & y \\
            z & w
            \arrow[""{name=0, anchor=center, inner sep=0}, "m", "\shortmid"{marking}, from=1-1, to=1-2]
            \arrow["f"', from=1-1, to=2-1]
            \arrow["g", from=1-2, to=2-2]
            \arrow[""{name=1, anchor=center, inner sep=0}, "n"', "\shortmid"{marking}, from=2-1, to=2-2]
            \arrow["\alpha"{description}, draw=none, from=0, to=1]
          \end{tikzcd}
        \end{equation*}
      such that the equation in $\Prof$ holds:
        \begin{equation*}
          % https://q.uiver.app/#q=WzAsMTAsWzAsMCwiMSJdLFsxLDAsIjEiXSxbMSwxLCJGeiJdLFswLDEsIkZ4Il0sWzEsMiwiRnciXSxbMCwyLCJGeSJdLFsyLDAsIjEiXSxbMiwyLCJGdyJdLFswLDMsIkZ5Il0sWzIsMywiRnciXSxbMCwxLCJcXGlkIiwwLHsic3R5bGUiOnsiYm9keSI6eyJuYW1lIjoiYmFycmVkIn19fV0sWzEsMiwiYyJdLFswLDMsImEiLDJdLFszLDIsIkZmIiwwLHsic3R5bGUiOnsiYm9keSI6eyJuYW1lIjoiYmFycmVkIn19fV0sWzIsNCwiRm4iXSxbMyw1LCJGbSIsMl0sWzUsNCwiRmciLDIseyJzdHlsZSI6eyJib2R5Ijp7Im5hbWUiOiJiYXJyZWQifX19XSxbMSw2LCJcXGlkIiwwLHsic3R5bGUiOnsiYm9keSI6eyJuYW1lIjoiYmFycmVkIn19fV0sWzYsNywiZCJdLFs0LDcsIlxcaWQiLDIseyJzdHlsZSI6eyJib2R5Ijp7Im5hbWUiOiJiYXJyZWQifX19XSxbNSw4LCIiLDAseyJsZXZlbCI6Miwic3R5bGUiOnsiaGVhZCI6eyJuYW1lIjoibm9uZSJ9fX1dLFs4LDksIkZnIiwyLHsic3R5bGUiOnsiYm9keSI6eyJuYW1lIjoiYmFycmVkIn19fV0sWzcsOSwiIiwwLHsibGV2ZWwiOjIsInN0eWxlIjp7ImhlYWQiOnsibmFtZSI6Im5vbmUifX19XSxbMTMsMTYsIkZcXGFscGhhIiwxLHsic2hvcnRlbiI6eyJzb3VyY2UiOjIwLCJ0YXJnZXQiOjIwfSwic3R5bGUiOnsiYm9keSI6eyJuYW1lIjoibm9uZSJ9LCJoZWFkIjp7Im5hbWUiOiJub25lIn19fV0sWzEwLDEzLCJcXG92ZXJsaW5lIGYiLDEseyJsYWJlbF9wb3NpdGlvbiI6NDAsInNob3J0ZW4iOnsic291cmNlIjoyMCwidGFyZ2V0IjoyMH0sInN0eWxlIjp7ImJvZHkiOnsibmFtZSI6Im5vbmUifSwiaGVhZCI6eyJuYW1lIjoibm9uZSJ9fX1dLFsxNywxOSwiXFxvdmVybGluZSBuIiwxLHsic2hvcnRlbiI6eyJzb3VyY2UiOjIwLCJ0YXJnZXQiOjIwfSwic3R5bGUiOnsiYm9keSI6eyJuYW1lIjoibm9uZSJ9LCJoZWFkIjp7Im5hbWUiOiJub25lIn19fV0sWzQsMjEsIlxcY29uZyIsMSx7InNob3J0ZW4iOnsidGFyZ2V0IjoyMH0sInN0eWxlIjp7ImJvZHkiOnsibmFtZSI6Im5vbmUifSwiaGVhZCI6eyJuYW1lIjoibm9uZSJ9fX1dXQ==
          \begin{tikzcd}
            1 & 1 & 1 \\
            Fx & Fz \\
            Fy & Fw & Fw \\
            Fy && Fw
            \arrow[""{name=0, anchor=center, inner sep=0}, "\id"{inner sep=.8ex}, "\shortmid"{marking}, from=1-1, to=1-2]
            \arrow["a"', from=1-1, to=2-1]
            \arrow[""{name=1, anchor=center, inner sep=0}, "\id"{inner sep=.8ex}, "\shortmid"{marking}, from=1-2, to=1-3]
            \arrow["c", from=1-2, to=2-2]
            \arrow["d", from=1-3, to=3-3]
            \arrow[""{name=2, anchor=center, inner sep=0}, "Ff"{inner sep=.8ex}, "\shortmid"{marking}, from=2-1, to=2-2]
            \arrow["Fm"', from=2-1, to=3-1]
            \arrow["Fn", from=2-2, to=3-2]
            \arrow[""{name=3, anchor=center, inner sep=0}, "Fg"'{inner sep=.8ex}, "\shortmid"{marking}, from=3-1, to=3-2]
            \arrow[equals, from=3-1, to=4-1]
            \arrow[""{name=4, anchor=center, inner sep=0}, "\id"'{inner sep=.8ex}, "\shortmid"{marking}, from=3-2, to=3-3]
            \arrow[equals, from=3-3, to=4-3]
            \arrow[""{name=5, anchor=center, inner sep=0}, "Fg"'{inner sep=.8ex}, "\shortmid"{marking}, from=4-1, to=4-3]
            \arrow["{\overline f}"{description, pos=0.4}, draw=none, from=0, to=2]
            \arrow["{\overline n}"{description}, draw=none, from=1, to=4]
            \arrow["{F\alpha}"{description}, draw=none, from=2, to=3]
            \arrow["\cong"{description}, draw=none, from=3-2, to=5]
          \end{tikzcd}
          \quad=\quad
          % https://q.uiver.app/#q=WzAsOSxbMCwwLCIxIl0sWzEsMCwiMSJdLFswLDEsIkZ4Il0sWzEsMiwiRnkiXSxbMCwyLCJGeSJdLFsyLDAsIjEiXSxbMiwyLCJGdyJdLFswLDMsIkZ5Il0sWzIsMywiRnciXSxbMCwxLCJcXGlkIiwwLHsic3R5bGUiOnsiYm9keSI6eyJuYW1lIjoiYmFycmVkIn19fV0sWzAsMiwiYSIsMl0sWzIsNCwiRm0iLDJdLFs0LDMsIlxcaWQiLDIseyJzdHlsZSI6eyJib2R5Ijp7Im5hbWUiOiJiYXJyZWQifX19XSxbMSw1LCJcXGlkIiwwLHsic3R5bGUiOnsiYm9keSI6eyJuYW1lIjoiYmFycmVkIn19fV0sWzUsNiwiZCJdLFszLDYsIkZnIiwyLHsic3R5bGUiOnsiYm9keSI6eyJuYW1lIjoiYmFycmVkIn19fV0sWzcsOCwiRmciLDIseyJzdHlsZSI6eyJib2R5Ijp7Im5hbWUiOiJiYXJyZWQifX19XSxbMSwzLCJiIl0sWzYsOCwiIiwwLHsibGV2ZWwiOjIsInN0eWxlIjp7ImhlYWQiOnsibmFtZSI6Im5vbmUifX19XSxbNCw3LCIiLDAseyJsZXZlbCI6Miwic3R5bGUiOnsiaGVhZCI6eyJuYW1lIjoibm9uZSJ9fX1dLFsxMywxNSwiXFxvdmVybGluZSBnIiwxLHsic2hvcnRlbiI6eyJzb3VyY2UiOjIwLCJ0YXJnZXQiOjIwfSwic3R5bGUiOnsiYm9keSI6eyJuYW1lIjoibm9uZSJ9LCJoZWFkIjp7Im5hbWUiOiJub25lIn19fV0sWzMsMTYsIlxcY29uZyIsMSx7InNob3J0ZW4iOnsidGFyZ2V0IjoyMH0sInN0eWxlIjp7ImJvZHkiOnsibmFtZSI6Im5vbmUifSwiaGVhZCI6eyJuYW1lIjoibm9uZSJ9fX1dLFs5LDEyLCJcXG92ZXJsaW5lIG0iLDEseyJzaG9ydGVuIjp7InNvdXJjZSI6MjAsInRhcmdldCI6MjB9LCJzdHlsZSI6eyJib2R5Ijp7Im5hbWUiOiJub25lIn0sImhlYWQiOnsibmFtZSI6Im5vbmUifX19XV0=
          \begin{tikzcd}
            1 & 1 & 1 \\
            Fx \\
            Fy & Fy & Fw \\
            Fy && Fw
            \arrow[""{name=0, anchor=center, inner sep=0}, "\id"{inner sep=.8ex}, "\shortmid"{marking}, from=1-1, to=1-2]
            \arrow["a"', from=1-1, to=2-1]
            \arrow[""{name=1, anchor=center, inner sep=0}, "\id"{inner sep=.8ex}, "\shortmid"{marking}, from=1-2, to=1-3]
            \arrow["b", from=1-2, to=3-2]
            \arrow["d", from=1-3, to=3-3]
            \arrow["Fm"', from=2-1, to=3-1]
            \arrow[""{name=2, anchor=center, inner sep=0}, "\id"'{inner sep=.8ex}, "\shortmid"{marking}, from=3-1, to=3-2]
            \arrow[equals, from=3-1, to=4-1]
            \arrow[""{name=3, anchor=center, inner sep=0}, "Fg"'{inner sep=.8ex}, "\shortmid"{marking}, from=3-2, to=3-3]
            \arrow[equals, from=3-3, to=4-3]
            \arrow[""{name=4, anchor=center, inner sep=0}, "Fg"'{inner sep=.8ex}, "\shortmid"{marking}, from=4-1, to=4-3]
            \arrow["{\overline m}"{description}, draw=none, from=0, to=2]
            \arrow["{\overline g}"{description}, draw=none, from=1, to=3]
            \arrow["\cong"{description}, draw=none, from=3-2, to=4]
          \end{tikzcd}.
        \end{equation*}
      \end{enumerate}
  Composition and identities in the tight direction are given as follows.
    \begin{itemize}
      \item For composable arrows $(f,\overline f)\colon (x,a)\to(y,b)$ and
        $(g,\overline g)\colon (y,b)\to (z,c)$, define the composite to be the pair
          \begin{equation*}
            (gf,F_{f,g}(\overline f,\overline g))\colon (x,a)\to (z,c),
          \end{equation*}
        where $F_{f,g}(\overline f,\overline g)$ is the component of the
        tight-to-loose composition comparison of $F$ applied to $\overline f$
        and $\overline g$.
      \item For an object $(x,a)$, define the identity arrow to be
          \begin{equation*}
            (1_x,F_x(1_a))\colon (x,a)\to (x,a),
          \end{equation*}
        where $F_x(1_a)$ is the component of the tight-to-loose identity
        comparison of $F$ applied to the identity morphism $1_a$ in $Fx$.
    \end{itemize}
  Composition and identities in the loose direction are given in the
  usual category-of-elements style.
    \begin{itemize}
      \item For composable proarrows $(m,\overline m)\colon (x,a)\proto (y,b)$ and
        $(n,\overline n)\colon (y,b)\proto (z,c)$, define the loose composite 
        to be the pair
          \begin{equation*}
            (m\odot n, \overline{m\odot n})\colon (x,a)\proto (z,c) 
          \end{equation*}
        whose second component is the usual elements-composite
          \begin{equation*}
            % https://q.uiver.app/#q=WzAsNCxbMCwwLCIobVxcb2RvdCBuKV8hKGEpIl0sWzEsMCwibl8hbV8hKGEpIl0sWzIsMCwibl8hKGIpIl0sWzMsMCwiYyJdLFswLDEsIlxcY29uZyJdLFsxLDIsIm5fIShcXG92ZXJsaW5lIG0pIl0sWzIsMywiXFxvdmVybGluZSBuIl1d
            \begin{tikzcd}
              {(m\odot n)_!(a)} & {n_!m_!(a)} & {n_!(b)} & c
              \arrow["\cong", from=1-1, to=1-2]
              \arrow["{n_!(\overline m)}", from=1-2, to=1-3]
              \arrow["{\overline n}", from=1-3, to=1-4]
            \end{tikzcd},
          \end{equation*}
        where the first isomorphism is the inverse of a component of the
        loose-to-tight composition comparison $F^{m,n}$.
      \item The loose identity on $(x,a)$ is the pair
        $(\id_x, \overline{1_a})\colon (x,a)\proto (x,a)$ whose second component
        $(\id_x)_!(a) \xto{\cong} a$ is the inverse of a component of the
        loose-to-tight identity comparison $F^x$.
    \end{itemize}
  Finally, composition of cells is by cell composition in the base double
  category $\dbl{B}$.

  It is tedious but straightforward to check that this construction results in a
  genuine double category $\Elt(F)$ and the evident projection
  $P_F\colon \Elt(F)\to\dbl{B}$ is a strict double functor.
\end{construction}

\begin{remark}[Cells in elements construction] \label{rmk:elements-cell-well-definition}
  The cells of $\Elt(F)$ for a twisted functor $F\colon\dbl{B}\twistto\Prof$ are thus
  cells $\alpha$ of the base double category $\dbl{B}$ satisfying
  appropriate well-typing and well-definition conditions. The latter is an
  equation, phrased in diagrams above. It is useful for computations to unpack
  this in terms of components. To this end, suppose we have such a cell:
    \begin{equation*}
          % https://q.uiver.app/#q=WzAsNCxbMCwwLCIoeCxhKSJdLFsxLDAsIih5LGIpIl0sWzAsMSwiKHosYykiXSxbMSwxLCIodyxkKSJdLFswLDEsIihtLFxcb3ZlcmxpbmUgbSkiLDAseyJzdHlsZSI6eyJib2R5Ijp7Im5hbWUiOiJiYXJyZWQifX19XSxbMCwyLCIoZixcXG92ZXJsaW5lIGYpIiwyXSxbMiwzLCIobixcXG92ZXJsaW5lIG4pIiwyLHsic3R5bGUiOnsiYm9keSI6eyJuYW1lIjoiYmFycmVkIn19fV0sWzEsMywiKGcsXFxvdmVybGluZSBnKSJdLFs0LDYsIlxcYWxwaGEiLDEseyJzaG9ydGVuIjp7InNvdXJjZSI6MjAsInRhcmdldCI6MjB9LCJzdHlsZSI6eyJib2R5Ijp7Im5hbWUiOiJub25lIn0sImhlYWQiOnsibmFtZSI6Im5vbmUifX19XV0=
          \begin{tikzcd}
            {(x,a)} & {(y,b)} \\
            {(z,c)} & {(w,d)}
            \arrow[""{name=0, anchor=center, inner sep=0}, "{(m,\overline m)}", "\shortmid"{marking}, from=1-1, to=1-2]
            \arrow["{(f,\overline f)}"', from=1-1, to=2-1]
            \arrow["{(g,\overline g)}", from=1-2, to=2-2]
            \arrow[""{name=1, anchor=center, inner sep=0}, "{(n,\overline n)}"', "\shortmid"{marking}, from=2-1, to=2-2]
            \arrow["\alpha"{description}, draw=none, from=0, to=1]
          \end{tikzcd}.
        \end{equation*}
  There are three significant mappings involved in the well-definition diagram,
  namely, a component of $F(\alpha)$, and the set functions $Fg(1,\overline{n})$ 
  and $Fg(\overline{m},1)$, related by 
      \begin{equation*}
        % https://q.uiver.app/#q=WzAsNCxbMCwwLCJGZihhLGMpIl0sWzAsMSwiRmcoRm0oYSksRm4oYykpIl0sWzIsMSwiRmcoRm0oYSksZCkiXSxbMiwwLCJGZyhiLGQpIl0sWzAsMSwiRihcXGFscGhhKV97YSxjfSIsMl0sWzEsMiwiRmcoMSxcXG92ZXJsaW5lIG4pIiwyXSxbMywyLCJGZyhcXG92ZXJsaW5lIG0sIDEpIl1d
        \begin{tikzcd}
          {Ff(a,c)} && {Fg(b,d)} \\
          {Fg(Fm(a),Fn(c))} && {Fg(Fm(a),d)}
          \arrow["{F(\alpha)_{a,c}}"', from=1-1, to=2-1]
          \arrow["{Fg(\overline m, 1)}", from=1-3, to=2-3]
          \arrow["{Fg(1,\overline n)}"', from=2-1, to=2-3]
        \end{tikzcd}.
      \end{equation*}
  Since $\overline f$ is an element of $Ff(a,c)$ and $\overline g$ is an element 
  of $Fg(b,d)$, we can run these elements around each side and the well-definition 
  condition is that the results are equal:
      \begin{equation*}
        Fg(1,\overline n)(F(\alpha)_{a,c}(\overline f)) = Fg(\overline m, 1)(\overline g)
      \end{equation*}
  in $Fg(Fm(a),d)$. A reasonable shorthand for this equation is thus
  $\overline n\, F(\alpha)(\overline f) = \overline g\, \overline m$, or,
  as a commutative square:
  \begin{equation*}
    % https://q.uiver.app/#q=WzAsNCxbMCwwLCJGbShhKSJdLFsxLDAsIkZuKGMpIl0sWzEsMSwiZCJdLFswLDEsImIiXSxbMCwxLCJGXFxhbHBoYShcXGJhciBmKSIsMCx7InN0eWxlIjp7ImJvZHkiOnsibmFtZSI6ImRhc2hlZCJ9fX1dLFsxLDIsIlxcYmFyIG4iXSxbMCwzLCJcXGJhciBtIiwyXSxbMywyLCJcXGJhciBnIiwyLHsic3R5bGUiOnsiYm9keSI6eyJuYW1lIjoiZGFzaGVkIn19fV1d
    \begin{tikzcd}
      {Fm(a)} & {Fn(c)} \\
      b & d
      \arrow["{F\alpha(\bar f)}", dashed, from=1-1, to=1-2]
      \arrow["{\bar m}"', from=1-1, to=2-1]
      \arrow["{\bar n}", from=1-2, to=2-2]
      \arrow["{\bar g}"', dashed, from=2-1, to=2-2]
    \end{tikzcd}. \qedhere
  \end{equation*}
\end{remark}

\begin{lemma}
  For any twisted normal lax functor $T\colon\dbl{B}\twistto\Prof$, the
  projection $J\colon\Elt(T)\to\dbl{B}$ is a loosely discrete double 
  opfibration in the sense of \cref{def:loosely-discrete-opfibration}. 
  Additionally, it is canonically cloven in the sense of 
  \cref{def:cleavage-for-disc-opfibration}.
\end{lemma}
\begin{proof}
  The commutative square 
    \begin{equation*}
      % https://q.uiver.app/#q=WzAsNCxbMCwwLCJcXEVsdChUKV8xIl0sWzEsMCwiXFxkYmx7Qn1fMSJdLFsxLDEsIlxcZGJse0J9XzAiXSxbMCwxLCJcXEVsdChUKV8wIl0sWzAsMSwiSl8xIl0sWzEsMiwiXFxzcmMiXSxbMCwzLCJcXHNyYyIsMl0sWzMsMiwiSl8wIiwyXV0=
      \begin{tikzcd}
        {\Elt(T)_1} & {\dbl{B}_1} \\
        {\Elt(T)_0} & {\dbl{B}_0}
        \arrow["{J_1}", from=1-1, to=1-2]
        \arrow["\src"', from=1-1, to=2-1]
        \arrow["\src", from=1-2, to=2-2]
        \arrow["{J_0}"', from=2-1, to=2-2]
      \end{tikzcd}
    \end{equation*}
  up to equivalence of categories presents the pullback of $\src$ along $J_0$ 
  by construction. Associated to a proarrow $m\colon x\proto y$ with $a\in Tx$ 
  is the proarrow $(m,1_{Tm(a)})\colon(x,a)\proto(y,Tm(a))$ of $\Elt(T)$. This is the
  object assignment of the cleavage $\sigma$. It is extended to arrows in the 
  evident way and functorial by the assumptions on $T$.
\end{proof}

\begin{example}[Elements construction of twisted point]
  Let $P\colon \dbl{1}\twistto\Prof$ denote a twisted point as in
  \cref{ex:twisted-point}. The double functor $\Elt(P)\to\dbl{1}$ is thus a 
  loosely discrete opfibration and the double category $\Elt(P)$ is loosely 
  discrete by the general considerations above. It can moreover be seen directly 
  that $\Elt(P)_0\cong P(*)$ and that there is a suitably structured 
  equivalence $\src\colon \Elt(P)_1 \to \Elt(P)_0$ using the fact that 
  $P(\id_*)$ is isomorphic to the identity functor on $P(*)$. 
\end{example}

\begin{example}[Elements construction for indexed categories]
  The elements construction applied to $F\colon\dbl{E}\twistto\Prof$ built from
  an ordinary op-indexed category $F\colon \cat{E}\to\Cat$ as in
  \cref{ex:indexed-cats-as-tw-funcs} results in a strict double category 
  $\Elt(F)$ and a strict double functor $P\colon \Elt(F)\to\dbl{E}$ whose underlying
  1-functor in the loose direction is a cloven opfibration. That is, level 0
  of $P^\top$ is an ordinary cloven opfibration.
\end{example}

\begin{example}[Elements construction for presheaves]
  Likewise if $F\colon\dbl{E}\twistto\Span$ is built from an ordinary 
  copresheaf as in \cref{ex:copresheaves-as-tw-funcs}, then the resulting 
  double functor $P\colon\Elt(F)\to\dbl{E}$ has $(P^\top)_0$ an ordinary 
  \emph{discrete} opfibration.
\end{example}

\begin{example}[Elements construction of twisted representable]
  The elements construction applied to a representable twisted copresheaf 
  (\cref{construction:twistedrepresentables})
    \begin{equation*}
      \dbl{B}(a,-)\colon\dbl{B}\twistto\Prof
    \end{equation*}
  is precisely the loose double coslice under $a$
  (\cref{construction:double-coslice}) in the sense that there is an evident 
  isomorphism of double categories 
    \begin{equation*}
      \Elt(\dbl{B}(a,-)) \cong a\sslash\dbl{B}
    \end{equation*}
  that commutes with the projections to $\dbl{B}$. That is, the elements 
  construction applied to the twisted representable at $a$ is the double coslice
  under $a$ up to isomorphism in $\cat{DOpf}(\dbl{B})$, the category of loosely
  discrete double opfibrations.
\end{example}

\section{2-categories of twisted double functors}

Like ordinary double functors, twisted double functors between a fixed pair of
double categories form the objects of a 2-category. The 1- and 2-morphisms in
this 2-category are natural transformations and modifications, suitably adjusted
to twisted functors. That is familiar enough, but in order to produce cells
between our key examples of twisted functors, the twisted Homs, we will need to
generalize our notion of modification to allow a nontrivial loose boundary.
These \emph{relative modifications} will be the 2-morphisms in a larger
2-category of twisted functors, in which the domain double category is allowed
to vary while the codomain is held fixed.

\subsection{Natural transformations}

As with transformations between ordinary double functors
\cite[\S{1.5}]{pare2011}, the well-behaved notion of natural transformation
between twisted double functors has components in the \emph{tight} direction of
the target double category.

\begin{definition}[Natural transformation] \label{def:naturaltransformation}
  A \define{lax natural transformation} $\alpha: F \To G$ between twisted doubly
  lax functors $F,G: \dbl{D} \twistto \dbl{E}$ consists of:
  \begin{itemize}
    \item for each object $x$ in $\dbl{D}$, the \define{component} of $\alpha$
      at $x$, an arrow $\alpha_x: Fx \to Gx$ in $\dbl{E}$;
    \item for each arrow $f: x \to y$ in $\dbl{D}$, the \define{component} of
      $\alpha$ at $f$, a cell in $\dbl{E}$ of the form
      \begin{equation*}
        % https://q.uiver.app/#q=WzAsNCxbMCwwLCJGeCJdLFsxLDAsIkZ5Il0sWzAsMSwiR3giXSxbMSwxLCJHeSJdLFswLDEsIkZmIiwwLHsic3R5bGUiOnsiYm9keSI6eyJuYW1lIjoiYmFycmVkIn19fV0sWzAsMiwiXFxhbHBoYV94IiwyXSxbMiwzLCJHZiIsMix7InN0eWxlIjp7ImJvZHkiOnsibmFtZSI6ImJhcnJlZCJ9fX1dLFsxLDMsIlxcYWxwaGFfeSJdLFs0LDYsIlxcYWxwaGFfZiIsMSx7InNob3J0ZW4iOnsic291cmNlIjoyMCwidGFyZ2V0IjoyMH0sInN0eWxlIjp7ImJvZHkiOnsibmFtZSI6Im5vbmUifSwiaGVhZCI6eyJuYW1lIjoibm9uZSJ9fX1dXQ==
        \begin{tikzcd}
          Fx & Fy \\
          Gx & Gy
          \arrow[""{name=0, anchor=center, inner sep=0}, "Ff", "\shortmid"{marking}, from=1-1, to=1-2]
          \arrow["{\alpha_x}"', from=1-1, to=2-1]
          \arrow["{\alpha_y}", from=1-2, to=2-2]
          \arrow[""{name=1, anchor=center, inner sep=0}, "Gf"', "\shortmid"{marking}, from=2-1, to=2-2]
          \arrow["{\alpha_f}"{description}, draw=none, from=0, to=1]
        \end{tikzcd}
      \end{equation*}
    \item for each proarrow $m: x \proto y$ in $\dbl{D}$, the \define{naturality
      comparison} of $\alpha$ at $m$, a cell in $\dbl{E}$ of the form
      \begin{equation*}
        % https://q.uiver.app/#q=WzAsNixbMCwwLCJGeCJdLFswLDEsIkd4Il0sWzAsMiwiR3kiXSxbMSwwLCJGeCJdLFsxLDEsIkZ5Il0sWzEsMiwiR3kiXSxbMCwxLCJcXGFscGhhX3giLDJdLFsxLDIsIkdtIiwyXSxbMyw0LCJGbSJdLFs0LDUsIlxcYWxwaGFfeSJdLFswLDMsIiIsMix7ImxldmVsIjoyLCJzdHlsZSI6eyJib2R5Ijp7Im5hbWUiOiJiYXJyZWQifSwiaGVhZCI6eyJuYW1lIjoibm9uZSJ9fX1dLFsyLDUsIiIsMCx7ImxldmVsIjoyLCJzdHlsZSI6eyJib2R5Ijp7Im5hbWUiOiJiYXJyZWQifSwiaGVhZCI6eyJuYW1lIjoibm9uZSJ9fX1dLFsxMCwxMSwiXFxhbHBoYV9tIiwxLHsibGV2ZWwiOjEsInN0eWxlIjp7ImJvZHkiOnsibmFtZSI6Im5vbmUifSwiaGVhZCI6eyJuYW1lIjoibm9uZSJ9fX1dXQ==
        \begin{tikzcd}[row sep=scriptsize]
          Fx & Fx \\
          Gx & Fy \\
          Gy & Gy
          \arrow[""{name=0, anchor=center, inner sep=0}, "\shortmid"{marking}, Rightarrow, no head, from=1-1, to=1-2]
          \arrow["{\alpha_x}"', from=1-1, to=2-1]
          \arrow["Fm", from=1-2, to=2-2]
          \arrow["Gm"', from=2-1, to=3-1]
          \arrow["{\alpha_y}", from=2-2, to=3-2]
          \arrow[""{name=1, anchor=center, inner sep=0}, "\shortmid"{marking}, Rightarrow, no head, from=3-1, to=3-2]
          \arrow["{\alpha_m}"{description}, draw=none, from=0, to=1]
        \end{tikzcd}.
      \end{equation*}
  \end{itemize}
  The following axioms must be satisfied.
  \begin{itemize}
    \item Naturality with respect to cells: for each cell
      $\stdInlineCell{\gamma}$ in $\dbl{D}$,
      \begin{equation*}
        % https://q.uiver.app/#q=WzAsOSxbMCwwLCJGeCJdLFsxLDAsIkZ3Il0sWzAsMSwiR3giXSxbMSwxLCJHdyJdLFswLDIsIkd5Il0sWzEsMiwiR3oiXSxbMiwwLCJGdyJdLFsyLDIsIkd6Il0sWzIsMSwiRnoiXSxbMCwxLCJGZiIsMCx7InN0eWxlIjp7ImJvZHkiOnsibmFtZSI6ImJhcnJlZCJ9fX1dLFswLDIsIlxcYWxwaGFfeCIsMl0sWzIsMywiR2YiLDAseyJzdHlsZSI6eyJib2R5Ijp7Im5hbWUiOiJiYXJyZWQifX19XSxbMSwzLCJcXGFscGhhX3ciXSxbMiw0LCJHbSIsMl0sWzMsNSwiR24iXSxbNCw1LCJHZyIsMix7InN0eWxlIjp7ImJvZHkiOnsibmFtZSI6ImJhcnJlZCJ9fX1dLFsxLDYsIiIsMSx7ImxldmVsIjoyLCJzdHlsZSI6eyJib2R5Ijp7Im5hbWUiOiJiYXJyZWQifSwiaGVhZCI6eyJuYW1lIjoibm9uZSJ9fX1dLFs1LDcsIiIsMSx7ImxldmVsIjoyLCJzdHlsZSI6eyJib2R5Ijp7Im5hbWUiOiJiYXJyZWQifSwiaGVhZCI6eyJuYW1lIjoibm9uZSJ9fX1dLFs2LDgsIkZuIl0sWzgsNywiXFxhbHBoYV96Il0sWzksMTEsIlxcYWxwaGFfZiIsMSx7ImxhYmVsX3Bvc2l0aW9uIjo0MCwic2hvcnRlbiI6eyJzb3VyY2UiOjIwLCJ0YXJnZXQiOjIwfSwic3R5bGUiOnsiYm9keSI6eyJuYW1lIjoibm9uZSJ9LCJoZWFkIjp7Im5hbWUiOiJub25lIn19fV0sWzExLDE1LCJHXFxnYW1tYSIsMSx7InNob3J0ZW4iOnsic291cmNlIjoyMCwidGFyZ2V0IjoyMH0sInN0eWxlIjp7ImJvZHkiOnsibmFtZSI6Im5vbmUifSwiaGVhZCI6eyJuYW1lIjoibm9uZSJ9fX1dLFsxNiwxNywiXFxhbHBoYV9uIiwxLHsibGV2ZWwiOjEsInN0eWxlIjp7ImJvZHkiOnsibmFtZSI6Im5vbmUifSwiaGVhZCI6eyJuYW1lIjoibm9uZSJ9fX1dXQ==
        \begin{tikzcd}
          Fx & Fw & Fw \\
          Gx & Gw & Fz \\
          Gy & Gz & Gz
          \arrow[""{name=0, anchor=center, inner sep=0}, "Ff", "\shortmid"{marking}, from=1-1, to=1-2]
          \arrow["{\alpha_x}"', from=1-1, to=2-1]
          \arrow[""{name=1, anchor=center, inner sep=0}, "\shortmid"{marking}, Rightarrow, no head, from=1-2, to=1-3]
          \arrow["{\alpha_w}", from=1-2, to=2-2]
          \arrow["Fn", from=1-3, to=2-3]
          \arrow[""{name=2, anchor=center, inner sep=0}, "Gf", "\shortmid"{marking}, from=2-1, to=2-2]
          \arrow["Gm"', from=2-1, to=3-1]
          \arrow["Gn", from=2-2, to=3-2]
          \arrow["{\alpha_z}", from=2-3, to=3-3]
          \arrow[""{name=3, anchor=center, inner sep=0}, "Gg"', "\shortmid"{marking}, from=3-1, to=3-2]
          \arrow[""{name=4, anchor=center, inner sep=0}, "\shortmid"{marking}, Rightarrow, no head, from=3-2, to=3-3]
          \arrow["{\alpha_f}"{description, pos=0.4}, draw=none, from=0, to=2]
          \arrow["{\alpha_n}"{description}, draw=none, from=1, to=4]
          \arrow["{G\gamma}"{description}, draw=none, from=2, to=3]
        \end{tikzcd}
        \quad=\quad
        % https://q.uiver.app/#q=WzAsOSxbMSwwLCJGeCJdLFsxLDEsIkZ5Il0sWzIsMCwiRnciXSxbMiwxLCJGeiJdLFsxLDIsIkd5Il0sWzIsMiwiR3oiXSxbMCwwLCJGeCJdLFswLDEsIkd4Il0sWzAsMiwiR3kiXSxbMCwxLCJGbSIsMl0sWzAsMiwiRmYiLDAseyJzdHlsZSI6eyJib2R5Ijp7Im5hbWUiOiJiYXJyZWQifX19XSxbMiwzLCJGbiJdLFsxLDMsIkZnIiwwLHsic3R5bGUiOnsiYm9keSI6eyJuYW1lIjoiYmFycmVkIn19fV0sWzEsNCwiXFxhbHBoYV95IiwyXSxbNCw1LCJHZyIsMix7InN0eWxlIjp7ImJvZHkiOnsibmFtZSI6ImJhcnJlZCJ9fX1dLFszLDUsIlxcYWxwaGFfeiJdLFs2LDcsIlxcYWxwaGFfeCIsMl0sWzcsOCwiR20iLDJdLFs2LDAsIiIsMCx7ImxldmVsIjoyLCJzdHlsZSI6eyJib2R5Ijp7Im5hbWUiOiJiYXJyZWQifSwiaGVhZCI6eyJuYW1lIjoibm9uZSJ9fX1dLFs4LDQsIiIsMCx7ImxldmVsIjoyLCJzdHlsZSI6eyJib2R5Ijp7Im5hbWUiOiJiYXJyZWQifSwiaGVhZCI6eyJuYW1lIjoibm9uZSJ9fX1dLFsxMCwxMiwiRlxcZ2FtbWEiLDEseyJsYWJlbF9wb3NpdGlvbiI6NDAsInNob3J0ZW4iOnsic291cmNlIjoyMCwidGFyZ2V0IjoyMH0sInN0eWxlIjp7ImJvZHkiOnsibmFtZSI6Im5vbmUifSwiaGVhZCI6eyJuYW1lIjoibm9uZSJ9fX1dLFsxMiwxNCwiXFxhbHBoYV9nIiwxLHsic2hvcnRlbiI6eyJzb3VyY2UiOjIwLCJ0YXJnZXQiOjIwfSwic3R5bGUiOnsiYm9keSI6eyJuYW1lIjoibm9uZSJ9LCJoZWFkIjp7Im5hbWUiOiJub25lIn19fV0sWzE4LDE5LCJcXGFscGhhX20iLDEseyJsZXZlbCI6MSwic3R5bGUiOnsiYm9keSI6eyJuYW1lIjoibm9uZSJ9LCJoZWFkIjp7Im5hbWUiOiJub25lIn19fV1d
        \begin{tikzcd}
          Fx & Fx & Fw \\
          Gx & Fy & Fz \\
          Gy & Gy & Gz
          \arrow[""{name=0, anchor=center, inner sep=0}, "\shortmid"{marking}, Rightarrow, no head, from=1-1, to=1-2]
          \arrow["{\alpha_x}"', from=1-1, to=2-1]
          \arrow[""{name=1, anchor=center, inner sep=0}, "Ff", "\shortmid"{marking}, from=1-2, to=1-3]
          \arrow["Fm"', from=1-2, to=2-2]
          \arrow["Fn", from=1-3, to=2-3]
          \arrow["Gm"', from=2-1, to=3-1]
          \arrow[""{name=2, anchor=center, inner sep=0}, "Fg", "\shortmid"{marking}, from=2-2, to=2-3]
          \arrow["{\alpha_y}"', from=2-2, to=3-2]
          \arrow["{\alpha_z}", from=2-3, to=3-3]
          \arrow[""{name=3, anchor=center, inner sep=0}, "\shortmid"{marking}, Rightarrow, no head, from=3-1, to=3-2]
          \arrow[""{name=4, anchor=center, inner sep=0}, "Gg"', "\shortmid"{marking}, from=3-2, to=3-3]
          \arrow["{\alpha_m}"{description}, draw=none, from=0, to=3]
          \arrow["{F\gamma}"{description, pos=0.4}, draw=none, from=1, to=2]
          \arrow["{\alpha_g}"{description}, draw=none, from=2, to=4]
        \end{tikzcd};
      \end{equation*}
    \item Functorality of component cells: for each pair of arrows
      $x \xto{f} y \xto{g} z$ in $\dbl{D}$,
      \begin{equation*}
        % https://q.uiver.app/#q=WzAsOCxbMCwwLCJGeCJdLFsxLDAsIkZ5Il0sWzAsMSwiR3giXSxbMSwxLCJHeSJdLFsyLDAsIkZ6Il0sWzIsMSwiR3oiXSxbMCwyLCJHeCJdLFsyLDIsIkd6Il0sWzAsMSwiRmYiLDAseyJzdHlsZSI6eyJib2R5Ijp7Im5hbWUiOiJiYXJyZWQifX19XSxbMCwyLCJcXGFscGhhX3giLDJdLFsyLDMsIkdmIiwyLHsic3R5bGUiOnsiYm9keSI6eyJuYW1lIjoiYmFycmVkIn19fV0sWzEsMywiXFxhbHBoYV95IiwxXSxbMSw0LCJGZyIsMCx7InN0eWxlIjp7ImJvZHkiOnsibmFtZSI6ImJhcnJlZCJ9fX1dLFszLDUsIkdnIiwyLHsic3R5bGUiOnsiYm9keSI6eyJuYW1lIjoiYmFycmVkIn19fV0sWzQsNSwiXFxhbHBoYV96Il0sWzYsNywiRyhmIFxcY2RvdCBnKSIsMix7InN0eWxlIjp7ImJvZHkiOnsibmFtZSI6ImJhcnJlZCJ9fX1dLFsyLDYsIiIsMCx7ImxldmVsIjoyLCJzdHlsZSI6eyJoZWFkIjp7Im5hbWUiOiJub25lIn19fV0sWzUsNywiIiwwLHsibGV2ZWwiOjIsInN0eWxlIjp7ImhlYWQiOnsibmFtZSI6Im5vbmUifX19XSxbOCwxMCwiXFxhbHBoYV9mIiwxLHsic2hvcnRlbiI6eyJzb3VyY2UiOjIwLCJ0YXJnZXQiOjIwfSwic3R5bGUiOnsiYm9keSI6eyJuYW1lIjoibm9uZSJ9LCJoZWFkIjp7Im5hbWUiOiJub25lIn19fV0sWzEyLDEzLCJcXGFscGhhX2ciLDEseyJzaG9ydGVuIjp7InNvdXJjZSI6MjAsInRhcmdldCI6MjB9LCJzdHlsZSI6eyJib2R5Ijp7Im5hbWUiOiJub25lIn0sImhlYWQiOnsibmFtZSI6Im5vbmUifX19XSxbMywxNSwiR197ZixnfSIsMSx7ImxhYmVsX3Bvc2l0aW9uIjo0MCwic2hvcnRlbiI6eyJ0YXJnZXQiOjIwfSwic3R5bGUiOnsiYm9keSI6eyJuYW1lIjoibm9uZSJ9LCJoZWFkIjp7Im5hbWUiOiJub25lIn19fV1d
        \begin{tikzcd}
          Fx & Fy & Fz \\
          Gx & Gy & Gz \\
          Gx && Gz
          \arrow[""{name=0, anchor=center, inner sep=0}, "Ff", "\shortmid"{marking}, from=1-1, to=1-2]
          \arrow["{\alpha_x}"', from=1-1, to=2-1]
          \arrow[""{name=1, anchor=center, inner sep=0}, "Fg", "\shortmid"{marking}, from=1-2, to=1-3]
          \arrow["{\alpha_y}"{description}, from=1-2, to=2-2]
          \arrow["{\alpha_z}", from=1-3, to=2-3]
          \arrow[""{name=2, anchor=center, inner sep=0}, "Gf"', "\shortmid"{marking}, from=2-1, to=2-2]
          \arrow[Rightarrow, no head, from=2-1, to=3-1]
          \arrow[""{name=3, anchor=center, inner sep=0}, "Gg"', "\shortmid"{marking}, from=2-2, to=2-3]
          \arrow[Rightarrow, no head, from=2-3, to=3-3]
          \arrow[""{name=4, anchor=center, inner sep=0}, "{G(f \cdot g)}"', "\shortmid"{marking}, from=3-1, to=3-3]
          \arrow["{\alpha_f}"{description}, draw=none, from=0, to=2]
          \arrow["{\alpha_g}"{description}, draw=none, from=1, to=3]
          \arrow["{G_{f,g}}"{description, pos=0.4}, draw=none, from=2-2, to=4]
        \end{tikzcd}
        \quad=\quad
        % https://q.uiver.app/#q=WzAsNyxbMCwwLCJGeCJdLFsyLDAsIkZ6Il0sWzEsMCwiRnkiXSxbMCwxLCJGeCJdLFsyLDEsIkZ6Il0sWzAsMiwiR3giXSxbMiwyLCJHeiJdLFswLDIsIkZmIiwwLHsic3R5bGUiOnsiYm9keSI6eyJuYW1lIjoiYmFycmVkIn19fV0sWzIsMSwiRmciLDAseyJzdHlsZSI6eyJib2R5Ijp7Im5hbWUiOiJiYXJyZWQifX19XSxbMCwzLCIiLDAseyJsZXZlbCI6Miwic3R5bGUiOnsiaGVhZCI6eyJuYW1lIjoibm9uZSJ9fX1dLFsxLDQsIiIsMix7ImxldmVsIjoyLCJzdHlsZSI6eyJoZWFkIjp7Im5hbWUiOiJub25lIn19fV0sWzMsNCwiRihmIFxcY2RvdCBnKSIsMV0sWzUsNiwiRyhmXFxjZG90IGcpIiwyLHsic3R5bGUiOnsiYm9keSI6eyJuYW1lIjoiYmFycmVkIn19fV0sWzQsNiwiXFxhbHBoYV96Il0sWzMsNSwiXFxhbHBoYV94IiwyXSxbMiwxMSwiRl97ZixnfSIsMSx7ImxhYmVsX3Bvc2l0aW9uIjo0MCwic2hvcnRlbiI6eyJ0YXJnZXQiOjIwfSwic3R5bGUiOnsiYm9keSI6eyJuYW1lIjoibm9uZSJ9LCJoZWFkIjp7Im5hbWUiOiJub25lIn19fV0sWzExLDEyLCJcXGFscGhhX3tmIFxcY2RvdCBnfSIsMSx7InNob3J0ZW4iOnsic291cmNlIjoyMCwidGFyZ2V0IjoyMH0sInN0eWxlIjp7ImJvZHkiOnsibmFtZSI6Im5vbmUifSwiaGVhZCI6eyJuYW1lIjoibm9uZSJ9fX1dXQ==
        \begin{tikzcd}
          Fx & Fy & Fz \\
          Fx && Fz \\
          Gx && Gz
          \arrow["Ff", "\shortmid"{marking}, from=1-1, to=1-2]
          \arrow[Rightarrow, no head, from=1-1, to=2-1]
          \arrow["Fg", "\shortmid"{marking}, from=1-2, to=1-3]
          \arrow[Rightarrow, no head, from=1-3, to=2-3]
          \arrow[""{name=0, anchor=center, inner sep=0}, "{F(f \cdot g)}"{description}, from=2-1, to=2-3]
          \arrow["{\alpha_x}"', from=2-1, to=3-1]
          \arrow["{\alpha_z}", from=2-3, to=3-3]
          \arrow[""{name=1, anchor=center, inner sep=0}, "{G(f\cdot g)}"', "\shortmid"{marking}, from=3-1, to=3-3]
          \arrow["{F_{f,g}}"{description, pos=0.4}, draw=none, from=1-2, to=0]
          \arrow["{\alpha_{f \cdot g}}"{description}, draw=none, from=0, to=1]
        \end{tikzcd},
      \end{equation*}
      and for each object $x$ in $\dbl{D}$,
      \begin{equation*}
        % https://q.uiver.app/#q=WzAsNixbMCwyLCJHeCJdLFsxLDIsIkd4Il0sWzAsMSwiR3giXSxbMSwxLCJHeCJdLFswLDAsIkZ4Il0sWzEsMCwiRngiXSxbMiwwLCIiLDEseyJsZXZlbCI6Miwic3R5bGUiOnsiaGVhZCI6eyJuYW1lIjoibm9uZSJ9fX1dLFszLDEsIiIsMSx7ImxldmVsIjoyLCJzdHlsZSI6eyJoZWFkIjp7Im5hbWUiOiJub25lIn19fV0sWzAsMSwiRyAxX3giLDIseyJzdHlsZSI6eyJib2R5Ijp7Im5hbWUiOiJiYXJyZWQifX19XSxbMiwzLCJcXG1hdGhybXtpZH1fe0d4fSIsMV0sWzQsNSwiXFxtYXRocm17aWR9X3tGeH0iLDAseyJzdHlsZSI6eyJib2R5Ijp7Im5hbWUiOiJiYXJyZWQifX19XSxbNSwzLCJcXGFscGhhX3giXSxbNCwyLCJcXGFscGhhX3giLDJdLFs5LDgsIkdfeCIsMSx7InNob3J0ZW4iOnsic291cmNlIjoyMCwidGFyZ2V0IjoyMH0sInN0eWxlIjp7ImJvZHkiOnsibmFtZSI6Im5vbmUifSwiaGVhZCI6eyJuYW1lIjoibm9uZSJ9fX1dLFsxMCw5LCJcXG1hdGhybXtpZH1fe1xcYWxwaGFfeH0iLDEseyJzaG9ydGVuIjp7InNvdXJjZSI6MjAsInRhcmdldCI6MjB9LCJzdHlsZSI6eyJib2R5Ijp7Im5hbWUiOiJub25lIn0sImhlYWQiOnsibmFtZSI6Im5vbmUifX19XV0=
        \begin{tikzcd}
          Fx & Fx \\
          Gx & Gx \\
          Gx & Gx
          \arrow[""{name=0, anchor=center, inner sep=0}, "{\mathrm{id}_{Fx}}", "\shortmid"{marking}, from=1-1, to=1-2]
          \arrow["{\alpha_x}"', from=1-1, to=2-1]
          \arrow["{\alpha_x}", from=1-2, to=2-2]
          \arrow[""{name=1, anchor=center, inner sep=0}, "{\mathrm{id}_{Gx}}"{description}, from=2-1, to=2-2]
          \arrow[Rightarrow, no head, from=2-1, to=3-1]
          \arrow[Rightarrow, no head, from=2-2, to=3-2]
          \arrow[""{name=2, anchor=center, inner sep=0}, "{G 1_x}"', "\shortmid"{marking}, from=3-1, to=3-2]
          \arrow["{\mathrm{id}_{\alpha_x}}"{description}, draw=none, from=0, to=1]
          \arrow["{G_x}"{description}, draw=none, from=1, to=2]
        \end{tikzcd}
        \quad=\quad
        % https://q.uiver.app/#q=WzAsNixbMSwxLCJGeCJdLFswLDEsIkZ4Il0sWzAsMCwiRngiXSxbMSwwLCJGeCJdLFswLDIsIkd4Il0sWzEsMiwiR3giXSxbMSwwLCJGIDFfeCIsMV0sWzIsMywiXFxtYXRocm17aWR9X3tGeH0iLDAseyJzdHlsZSI6eyJib2R5Ijp7Im5hbWUiOiJiYXJyZWQifX19XSxbMiwxLCIiLDAseyJsZXZlbCI6Miwic3R5bGUiOnsiaGVhZCI6eyJuYW1lIjoibm9uZSJ9fX1dLFszLDAsIiIsMix7ImxldmVsIjoyLCJzdHlsZSI6eyJoZWFkIjp7Im5hbWUiOiJub25lIn19fV0sWzAsNSwiXFxhbHBoYV94Il0sWzQsNSwiRyAxX3giLDIseyJzdHlsZSI6eyJib2R5Ijp7Im5hbWUiOiJiYXJyZWQifX19XSxbMSw0LCJcXGFscGhhX3giLDJdLFs3LDYsIkZfeCIsMSx7InNob3J0ZW4iOnsic291cmNlIjoyMCwidGFyZ2V0IjoyMH0sInN0eWxlIjp7ImJvZHkiOnsibmFtZSI6Im5vbmUifSwiaGVhZCI6eyJuYW1lIjoibm9uZSJ9fX1dLFs2LDExLCJcXGFscGhhX3sxX3h9IiwxLHsic2hvcnRlbiI6eyJzb3VyY2UiOjIwLCJ0YXJnZXQiOjIwfSwic3R5bGUiOnsiYm9keSI6eyJuYW1lIjoibm9uZSJ9LCJoZWFkIjp7Im5hbWUiOiJub25lIn19fV1d
        \begin{tikzcd}
          Fx & Fx \\
          Fx & Fx \\
          Gx & Gx
          \arrow[""{name=0, anchor=center, inner sep=0}, "{\mathrm{id}_{Fx}}", "\shortmid"{marking}, from=1-1, to=1-2]
          \arrow[Rightarrow, no head, from=1-1, to=2-1]
          \arrow[Rightarrow, no head, from=1-2, to=2-2]
          \arrow[""{name=1, anchor=center, inner sep=0}, "{F 1_x}"{description}, from=2-1, to=2-2]
          \arrow["{\alpha_x}"', from=2-1, to=3-1]
          \arrow["{\alpha_x}", from=2-2, to=3-2]
          \arrow[""{name=2, anchor=center, inner sep=0}, "{G 1_x}"', "\shortmid"{marking}, from=3-1, to=3-2]
          \arrow["{F_x}"{description}, draw=none, from=0, to=1]
          \arrow["{\alpha_{1_x}}"{description}, draw=none, from=1, to=2]
        \end{tikzcd};
      \end{equation*}
    \item Coherence of naturality comparisons: for each pair of proarrows
      $x \xproto{m} y \xproto{n} z$ in $\dbl{D}$,
      \begin{equation*}
        \begin{tikzcd}
          Fx & Fx & Fx & Fx \\
          Gx & Fy & Fy \\
          Gy & Gy & Fz & Fz \\
          Gz & Gz & Gz & Gz
          \arrow[""{name=0, anchor=center, inner sep=0}, "\shortmid"{marking}, equals, from=1-1, to=1-2]
          \arrow["{\alpha_x}"', from=1-1, to=2-1]
          \arrow[""{name=1, anchor=center, inner sep=0}, "\shortmid"{marking}, equals, from=1-2, to=1-3]
          \arrow["Fm"', from=1-2, to=2-2]
          \arrow[""{name=2, anchor=center, inner sep=0}, "\shortmid"{marking}, equals, from=1-3, to=1-4]
          \arrow["Fm", from=1-3, to=2-3]
          \arrow["{F(m \odot n)}", from=1-4, to=3-4]
          \arrow["Gm"', from=2-1, to=3-1]
          \arrow[""{name=3, anchor=center, inner sep=0}, "\shortmid"{marking}, equals, from=2-2, to=2-3]
          \arrow["{\alpha_y}"', from=2-2, to=3-2]
          \arrow["Fn", from=2-3, to=3-3]
          \arrow[""{name=4, anchor=center, inner sep=0}, "\shortmid"{marking}, equals, from=3-1, to=3-2]
          \arrow["Gn"', from=3-1, to=4-1]
          \arrow["Gn", from=3-2, to=4-2]
          \arrow[""{name=5, anchor=center, inner sep=0}, "\shortmid"{marking}, equals, from=3-3, to=3-4]
          \arrow["{\alpha_z}"', from=3-3, to=4-3]
          \arrow["{\alpha_z}", from=3-4, to=4-4]
          \arrow[""{name=6, anchor=center, inner sep=0}, "\shortmid"{marking}, equals, from=4-1, to=4-2]
          \arrow[""{name=7, anchor=center, inner sep=0}, "\shortmid"{marking}, equals, from=4-2, to=4-3]
          \arrow[""{name=8, anchor=center, inner sep=0}, "\shortmid"{marking}, equals, from=4-3, to=4-4]
          \arrow["{\alpha_m}"{description}, draw=none, from=0, to=4]
          \arrow["{\id_{Fm}}"{description}, draw=none, from=1, to=3]
          \arrow["{F^{m,n}}"{description}, draw=none, from=2, to=5]
          \arrow["{\alpha_n}"{description}, draw=none, from=3, to=7]
          \arrow["{\id_{Gn}}"{description}, draw=none, from=4, to=6]
          \arrow["{\id_{\alpha_z}}"{description}, draw=none, from=5, to=8]
        \end{tikzcd}
        \quad=\quad
        % https://q.uiver.app/#q=WzAsMTAsWzEsMCwiRngiXSxbMiwwLCJGeCJdLFsxLDEsIkd4Il0sWzIsMiwiRnoiXSxbMSwzLCJHeiJdLFsyLDMsIkd6Il0sWzAsMCwiRngiXSxbMCwxLCJHeCJdLFswLDIsIkd5Il0sWzAsMywiR3oiXSxbMCwxLCIiLDAseyJsZXZlbCI6Miwic3R5bGUiOnsiYm9keSI6eyJuYW1lIjoiYmFycmVkIn0sImhlYWQiOnsibmFtZSI6Im5vbmUifX19XSxbMCwyLCJcXGFscGhhX3giXSxbMSwzLCJGKG0gXFxvZG90IG4pIl0sWzMsNSwiXFxhbHBoYV96Il0sWzIsNCwiRyhtIFxcb2RvdCBuKSIsMV0sWzQsNSwiIiwyLHsibGV2ZWwiOjIsInN0eWxlIjp7ImJvZHkiOnsibmFtZSI6ImJhcnJlZCJ9LCJoZWFkIjp7Im5hbWUiOiJub25lIn19fV0sWzYsNywiXFxhbHBoYV94IiwyXSxbNiwwLCIiLDAseyJsZXZlbCI6Miwic3R5bGUiOnsiYm9keSI6eyJuYW1lIjoiYmFycmVkIn0sImhlYWQiOnsibmFtZSI6Im5vbmUifX19XSxbNyw4LCJHbSIsMl0sWzgsOSwiR24iLDJdLFs5LDQsIiIsMix7ImxldmVsIjoyLCJzdHlsZSI6eyJib2R5Ijp7Im5hbWUiOiJiYXJyZWQifSwiaGVhZCI6eyJuYW1lIjoibm9uZSJ9fX1dLFs3LDIsIiIsMix7ImxldmVsIjoyLCJzdHlsZSI6eyJib2R5Ijp7Im5hbWUiOiJiYXJyZWQifSwiaGVhZCI6eyJuYW1lIjoibm9uZSJ9fX1dLFsxMCwxNSwiXFxhbHBoYV97bSBcXG9kb3Qgbn0iLDEseyJzaG9ydGVuIjp7InNvdXJjZSI6MjAsInRhcmdldCI6MjB9LCJzdHlsZSI6eyJib2R5Ijp7Im5hbWUiOiJub25lIn0sImhlYWQiOnsibmFtZSI6Im5vbmUifX19XSxbMTcsMjEsIlxcaWRfe1xcYWxwaGFfeH0iLDEseyJzaG9ydGVuIjp7InNvdXJjZSI6MjAsInRhcmdldCI6MjB9LCJzdHlsZSI6eyJib2R5Ijp7Im5hbWUiOiJub25lIn0sImhlYWQiOnsibmFtZSI6Im5vbmUifX19XSxbMjEsMjAsIkdee20sbn0iLDEseyJvZmZzZXQiOjEsInNob3J0ZW4iOnsic291cmNlIjoyMCwidGFyZ2V0IjoyMH0sInN0eWxlIjp7ImJvZHkiOnsibmFtZSI6Im5vbmUifSwiaGVhZCI6eyJuYW1lIjoibm9uZSJ9fX1dXQ==
        \begin{tikzcd}
          Fx & Fx & Fx \\
          Gx & Gx \\
          Gy && Fz \\
          Gz & Gz & Gz
          \arrow[""{name=0, anchor=center, inner sep=0}, "\shortmid"{marking}, equals, from=1-1, to=1-2]
          \arrow["{\alpha_x}"', from=1-1, to=2-1]
          \arrow[""{name=1, anchor=center, inner sep=0}, "\shortmid"{marking}, equals, from=1-2, to=1-3]
          \arrow["{\alpha_x}", from=1-2, to=2-2]
          \arrow["{F(m \odot n)}", from=1-3, to=3-3]
          \arrow[""{name=2, anchor=center, inner sep=0}, "\shortmid"{marking}, equals, from=2-1, to=2-2]
          \arrow["Gm"', from=2-1, to=3-1]
          \arrow["{G(m \odot n)}"{description}, from=2-2, to=4-2]
          \arrow["Gn"', from=3-1, to=4-1]
          \arrow["{\alpha_z}", from=3-3, to=4-3]
          \arrow[""{name=3, anchor=center, inner sep=0}, "\shortmid"{marking}, equals, from=4-1, to=4-2]
          \arrow[""{name=4, anchor=center, inner sep=0}, "\shortmid"{marking}, equals, from=4-2, to=4-3]
          \arrow["{\id_{\alpha_x}}"{description}, draw=none, from=0, to=2]
          \arrow["{\alpha_{m \odot n}}"{description}, draw=none, from=1, to=4]
          \arrow["{G^{m,n}}"{description}, shift right, draw=none, from=2, to=3]
        \end{tikzcd},
      \end{equation*}
      and for each object $x$ in $\dbl{D}$,
      \begin{equation*}
        % https://q.uiver.app/#q=WzAsNixbMCwwLCJGeCJdLFsxLDAsIkZ4Il0sWzAsMSwiRngiXSxbMSwxLCJGeCJdLFsxLDIsIkd4Il0sWzAsMiwiR3giXSxbMCwxLCIiLDAseyJsZXZlbCI6Miwic3R5bGUiOnsiYm9keSI6eyJuYW1lIjoiYmFycmVkIn0sImhlYWQiOnsibmFtZSI6Im5vbmUifX19XSxbMiwzLCIiLDAseyJsZXZlbCI6Miwic3R5bGUiOnsiYm9keSI6eyJuYW1lIjoiYmFycmVkIn0sImhlYWQiOnsibmFtZSI6Im5vbmUifX19XSxbMCwyLCIxX3tGeH0iLDJdLFsxLDMsIkZcXGlkX3giXSxbMyw0LCJcXGFscGhhX3giXSxbMiw1LCJcXGFscGhhX3giLDJdLFs1LDQsIiIsMix7ImxldmVsIjoyLCJzdHlsZSI6eyJib2R5Ijp7Im5hbWUiOiJiYXJyZWQifSwiaGVhZCI6eyJuYW1lIjoibm9uZSJ9fX1dLFs2LDcsIkZeeCIsMSx7InNob3J0ZW4iOnsic291cmNlIjoyMCwidGFyZ2V0IjoyMH0sInN0eWxlIjp7ImJvZHkiOnsibmFtZSI6Im5vbmUifSwiaGVhZCI6eyJuYW1lIjoibm9uZSJ9fX1dLFs3LDEyLCJcXGlkX3tcXGFscGhhX3h9IiwxLHsic2hvcnRlbiI6eyJzb3VyY2UiOjIwLCJ0YXJnZXQiOjIwfSwic3R5bGUiOnsiYm9keSI6eyJuYW1lIjoibm9uZSJ9LCJoZWFkIjp7Im5hbWUiOiJub25lIn19fV1d
        \begin{tikzcd}
          Fx & Fx \\
          Fx & Fx \\
          Gx & Gx
          \arrow[""{name=0, anchor=center, inner sep=0}, "\shortmid"{marking}, equals, from=1-1, to=1-2]
          \arrow["{1_{Fx}}"', from=1-1, to=2-1]
          \arrow["{F\id_x}", from=1-2, to=2-2]
          \arrow[""{name=1, anchor=center, inner sep=0}, "\shortmid"{marking}, equals, from=2-1, to=2-2]
          \arrow["{\alpha_x}"', from=2-1, to=3-1]
          \arrow["{\alpha_x}", from=2-2, to=3-2]
          \arrow[""{name=2, anchor=center, inner sep=0}, "\shortmid"{marking}, equals, from=3-1, to=3-2]
          \arrow["{F^x}"{description}, draw=none, from=0, to=1]
          \arrow["{\id_{\alpha_x}}"{description}, draw=none, from=1, to=2]
        \end{tikzcd}
        \quad=\quad
        % https://q.uiver.app/#q=WzAsOSxbMSwwLCJGeCJdLFsxLDEsIkd4Il0sWzEsMiwiR3giXSxbMiwwLCJGeCJdLFsyLDEsIkZ4Il0sWzIsMiwiR3giXSxbMCwxLCJHeCJdLFswLDIsIkd4Il0sWzAsMCwiRngiXSxbMSwyLCJHIFxcaWRfeCJdLFswLDEsIlxcYWxwaGFfeCJdLFszLDQsIkZcXGlkX3giXSxbNCw1LCJcXGFscGhhX3giXSxbMCwzLCIiLDAseyJsZXZlbCI6Miwic3R5bGUiOnsiYm9keSI6eyJuYW1lIjoiYmFycmVkIn0sImhlYWQiOnsibmFtZSI6Im5vbmUifX19XSxbMiw1LCIiLDIseyJsZXZlbCI6Miwic3R5bGUiOnsiYm9keSI6eyJuYW1lIjoiYmFycmVkIn0sImhlYWQiOnsibmFtZSI6Im5vbmUifX19XSxbNiw3LCIxX3tHeH0iLDJdLFs2LDEsIiIsMix7ImxldmVsIjoyLCJzdHlsZSI6eyJib2R5Ijp7Im5hbWUiOiJiYXJyZWQifSwiaGVhZCI6eyJuYW1lIjoibm9uZSJ9fX1dLFs3LDIsIiIsMSx7ImxldmVsIjoyLCJzdHlsZSI6eyJib2R5Ijp7Im5hbWUiOiJiYXJyZWQifSwiaGVhZCI6eyJuYW1lIjoibm9uZSJ9fX1dLFs4LDYsIlxcYWxwaGFfeCIsMl0sWzgsMCwiIiwyLHsibGV2ZWwiOjIsInN0eWxlIjp7ImJvZHkiOnsibmFtZSI6ImJhcnJlZCJ9LCJoZWFkIjp7Im5hbWUiOiJub25lIn19fV0sWzEzLDE0LCJcXGFscGhhX3tcXGlkX3h9IiwxLHsic2hvcnRlbiI6eyJzb3VyY2UiOjIwLCJ0YXJnZXQiOjIwfSwic3R5bGUiOnsiYm9keSI6eyJuYW1lIjoibm9uZSJ9LCJoZWFkIjp7Im5hbWUiOiJub25lIn19fV0sWzE2LDE3LCJHXngiLDEseyJzaG9ydGVuIjp7InNvdXJjZSI6MjAsInRhcmdldCI6MjB9LCJzdHlsZSI6eyJib2R5Ijp7Im5hbWUiOiJub25lIn0sImhlYWQiOnsibmFtZSI6Im5vbmUifX19XSxbMTksMTYsIlxcaWRfe1xcYWxwaGFfeH0iLDEseyJzaG9ydGVuIjp7InNvdXJjZSI6MjAsInRhcmdldCI6MjB9LCJzdHlsZSI6eyJib2R5Ijp7Im5hbWUiOiJub25lIn0sImhlYWQiOnsibmFtZSI6Im5vbmUifX19XV0=
        \begin{tikzcd}
          Fx & Fx & Fx \\
          Gx & Gx & Fx \\
          Gx & Gx & Gx
          \arrow[""{name=0, anchor=center, inner sep=0}, "\shortmid"{marking}, equals, from=1-1, to=1-2]
          \arrow["{\alpha_x}"', from=1-1, to=2-1]
          \arrow[""{name=1, anchor=center, inner sep=0}, "\shortmid"{marking}, equals, from=1-2, to=1-3]
          \arrow["{\alpha_x}", from=1-2, to=2-2]
          \arrow["{F\id_x}", from=1-3, to=2-3]
          \arrow[""{name=2, anchor=center, inner sep=0}, "\shortmid"{marking}, equals, from=2-1, to=2-2]
          \arrow["{1_{Gx}}"', from=2-1, to=3-1]
          \arrow["{G \id_x}", from=2-2, to=3-2]
          \arrow["{\alpha_x}", from=2-3, to=3-3]
          \arrow[""{name=3, anchor=center, inner sep=0}, "\shortmid"{marking}, equals, from=3-1, to=3-2]
          \arrow[""{name=4, anchor=center, inner sep=0}, "\shortmid"{marking}, equals, from=3-2, to=3-3]
          \arrow["{\id_{\alpha_x}}"{description}, draw=none, from=0, to=2]
          \arrow["{\alpha_{\id_x}}"{description}, draw=none, from=1, to=4]
          \arrow["{G^x}"{description}, draw=none, from=2, to=3]
        \end{tikzcd}.
      \end{equation*}
  \end{itemize}
  A lax natural transformation is \define{pseudo} or \define{strong} if each
  naturality comparison has an inverse with respect to loose composition, and it
  is \define{strict} if each naturality comparison is a loose identity.
\end{definition}

The prime example of a twisted functor is the twisted Hom functor on a double
category (\cref{construction:twistedhomfunctor}). Moving up a level, any
ordinary double functor induces a natural transformation between the
corresponding twisted Homs, once the codomain Hom is restricted along the
functor. This construction serves not only as a first example of a natural
transformation between twisted functors, but as the next step toward embedding
the 2-category of double categories into a 2-category of twisted functors.

\begin{construction}[Transformation between twisted Homs]
  \label{construction:double-functor-to-twisted-transformation}
  A lax double functor $F: \dbl{D} \to \dbl{E}$ induces a lax natural
  transformation
  \begin{equation*}
    \tilde F: \dbl{D}(-, =) \To \dbl{E}(F(-), F(=)):
      \dbl{D}^\co \times \dbl{D} \twistto \Prof
  \end{equation*}
  constructed as follows:
  \begin{itemize}
    \item The component at a pair of objects $x$ and $y$ in $\dbl{D}$ is the functor
      \begin{equation*}
        \tilde F_{x,y}: \dbl{D}(x,y) \to \dbl{E}(Fx, Fy)
      \end{equation*}
      that sends a proarrow $m: x \proto y$ in $\dbl{D}$ to $Fm: Fx \proto Fy$
      and a tightly globular cell $\alpha: m \To n$ to $F\alpha: Fm \To Fn$.
    \item The component at a pair of arrows $f: x \to w$ and $g: y \to z$ in $\dbl{D}$
      is the natural transformation
      \begin{equation*}
        % https://q.uiver.app/#q=WzAsNCxbMCwwLCJcXGRibHtEfSh4LHkpIl0sWzEsMCwiXFxkYmx7RH0odyx6KSJdLFswLDEsIlxcZGJse0V9KEZ4LEZ5KSJdLFsxLDEsIlxcZGJse0V9KEZ3LEZ6KSJdLFswLDEsIlxcZGJse0R9KGYsZykiLDAseyJzdHlsZSI6eyJib2R5Ijp7Im5hbWUiOiJiYXJyZWQifX19XSxbMCwyLCJcXHRpbGRlIEZfe3gseX0iLDJdLFsxLDMsIlxcdGlsZGUgRl97dyx6fSJdLFsyLDMsIlxcZGJse0V9KEZmLEZnKSIsMix7InN0eWxlIjp7ImJvZHkiOnsibmFtZSI6ImJhcnJlZCJ9fX1dLFs0LDcsIlxcdGlsZGUgRl97ZixnfSIsMSx7InNob3J0ZW4iOnsic291cmNlIjoyMCwidGFyZ2V0IjoyMH0sInN0eWxlIjp7ImJvZHkiOnsibmFtZSI6Im5vbmUifSwiaGVhZCI6eyJuYW1lIjoibm9uZSJ9fX1dXQ==
        \begin{tikzcd}
          {\dbl{D}(x,y)} & {\dbl{D}(w,z)} \\
          {\dbl{E}(Fx,Fy)} & {\dbl{E}(Fw,Fz)}
          \arrow[""{name=0, anchor=center, inner sep=0}, "{\dbl{D}(f,g)}"{inner sep=.8ex}, "\shortmid"{marking}, from=1-1, to=1-2]
          \arrow["{\tilde F_{x,y}}"', from=1-1, to=2-1]
          \arrow["{\tilde F_{w,z}}", from=1-2, to=2-2]
          \arrow[""{name=1, anchor=center, inner sep=0}, "{\dbl{E}(Ff,Fg)}"'{inner sep=.8ex}, "\shortmid"{marking}, from=2-1, to=2-2]
          \arrow["{\tilde F_{f,g}}"{description}, draw=none, from=0, to=1]
        \end{tikzcd}
      \end{equation*}
      that sends a cell $\stdInlineCell{\alpha}$ in $\dbl{D}$ to
      $\inlineCell{Fx}{Fy}{Fw}{Fz}{Fm}{Fn}{Ff}{Fg}{F\alpha}$.
    \item The naturality comparison at a pair of proarrows $p: x' \proto x$ and
      $q: y \proto y'$ in $\dbl{D}$ is the natural transformation
      \begin{equation*}
        % https://q.uiver.app/#q=WzAsNixbMCwwLCJcXGRibHtEfSh4LHkpIl0sWzAsMSwiXFxkYmx7RX0oRngsRnkpIl0sWzAsMiwiXFxkYmx7RX0oRngnLEZ5JykiXSxbMSwyLCJcXGRibHtFfShGeCcsRnknKSJdLFsxLDAsIlxcZGJse0R9KHgseSkiXSxbMSwxLCJcXGRibHtEfSh4Jyx5JykiXSxbMCwxLCJcXHRpbGRlIEZfe3gseX0iLDJdLFsxLDIsIlxcZGJse0V9KEZwLEZxKSIsMl0sWzIsMywiIiwyLHsibGV2ZWwiOjIsInN0eWxlIjp7ImJvZHkiOnsibmFtZSI6ImJhcnJlZCJ9LCJoZWFkIjp7Im5hbWUiOiJub25lIn19fV0sWzQsNSwiXFxkYmx7RH0ocCxxKSJdLFs1LDMsIlxcdGlsZGUgRl97eCcseSd9Il0sWzAsNCwiIiwwLHsibGV2ZWwiOjIsInN0eWxlIjp7ImJvZHkiOnsibmFtZSI6ImJhcnJlZCJ9LCJoZWFkIjp7Im5hbWUiOiJub25lIn19fV0sWzExLDgsIlxcdGlsZGUgRl97cCxxfSIsMSx7InNob3J0ZW4iOnsic291cmNlIjoyMCwidGFyZ2V0IjoyMH0sInN0eWxlIjp7ImJvZHkiOnsibmFtZSI6Im5vbmUifSwiaGVhZCI6eyJuYW1lIjoibm9uZSJ9fX1dXQ==
        \begin{tikzcd}
          {\dbl{D}(x,y)} & {\dbl{D}(x,y)} \\
          {\dbl{E}(Fx,Fy)} & {\dbl{D}(x',y')} \\
          {\dbl{E}(Fx',Fy')} & {\dbl{E}(Fx',Fy')}
          \arrow[""{name=0, anchor=center, inner sep=0}, "\shortmid"{marking}, equals, from=1-1, to=1-2]
          \arrow["{\tilde F_{x,y}}"', from=1-1, to=2-1]
          \arrow["{\dbl{D}(p,q)}", from=1-2, to=2-2]
          \arrow["{\dbl{E}(Fp,Fq)}"', from=2-1, to=3-1]
          \arrow["{\tilde F_{x',y'}}", from=2-2, to=3-2]
          \arrow[""{name=1, anchor=center, inner sep=0}, "\shortmid"{marking}, equals, from=3-1, to=3-2]
          \arrow["{\tilde F_{p,q}}"{description}, draw=none, from=0, to=1]
        \end{tikzcd}
      \end{equation*}
      whose component at a proarrow $m: x \proto y$ in $\dbl{D}$ is the globular
      cell in $\dbl{E}$
      \begin{equation*}
        (\tilde F_{p,q})_m \coloneqq F_{p,m,q}:
          Fp \odot Fm \odot Fq \to F(p \odot m \odot q)  
      \end{equation*}
      given by composition comparisons of $F$.
  \end{itemize}
  Moreover, when $F$ is a pseudo (resp.\ strict) double functor, $\tilde F$ is a
  pseudo (resp.\ strict) natural transformation between twisted functors.
\end{construction}

\begin{lemma}
  The lax natural transformation $\tilde F$ between twisted functors induced by
  a lax double functor $F$ is well-defined.
\end{lemma}
\begin{proof}
  A special case of \cref{lem:module-to-twisted-transformation} below.
\end{proof}

Having made this construction, one naturally asks whether it can be generalized
to allow restricting the two slots of the codomain twisted Hom along
\emph{different} double functors. The input data suited to this generalization
is a module between double functors. Recall that a \define{module}
$M: F \proTo G$ between lax double functors $F,G: \dbl{D} \to \dbl{E}$ is a lax
double functor $M: \dbl{D} \times \Loose \to \dbl{E}$ such that $M(-, 0) = F$ and
$M(-, 1) = G$. For a complete unpacking of this definition, see the original
source \cite[\mbox{Definition 3.2}]{pare2011} or alternatively
\cite[\mbox{Definition 9.1}]{lambert2024}.

\begin{construction}[Transformation between twisted Homs, generalized]
  \label{construction:module-to-twisted-transformation}
  Let $M: F \proTo G$ be a module between lax double functors
  $F,G: \dbl{D} \to \dbl{E}$. A lax natural transformation
  \begin{equation*}
    \tilde M: \dbl{D}(-, =) \To \dbl{E}(F(-), G(=)): \dbl{D}^\co \times \dbl{D} \twistto \Prof
  \end{equation*}
  is constructed as follows:
  \begin{itemize}
    \item The component at a pair of objects $x$ and $y$ in $\dbl{D}$ is the
      functor
      \begin{equation*}
        \tilde M_{x,y}: \dbl{D}(x,y) \to \dbl{E}(Fx, Gy)
      \end{equation*}
      that sends a proarrow $m: x \proto y$ in $\dbl{D}$ to its image
      $Mm: Fx \proto Gy$ and a tightly globular cell $\alpha: m \To n$ to 
      $M\alpha: Mm \To Mn$.
    \item The component at a pair of arrows $f: x \to w$ and $g: y \to z$ in
      $\dbl{D}$ is the natural transformation
      \begin{equation*}
        % https://q.uiver.app/#q=WzAsNCxbMCwwLCJcXGRibHtEfSh4LHkpIl0sWzEsMCwiXFxkYmx7RH0odyx6KSJdLFswLDEsIlxcZGJse0V9KEZ4LEd5KSJdLFsxLDEsIlxcZGJse0V9KEZ3LEd6KSJdLFswLDEsIlxcZGJse0R9KGYsZykiLDAseyJzdHlsZSI6eyJib2R5Ijp7Im5hbWUiOiJiYXJyZWQifX19XSxbMCwyLCJcXHRpbGRlIE1fe3gseX0iLDJdLFsxLDMsIlxcdGlsZGUgTV97dyx6fSJdLFsyLDMsIlxcZGJse0V9KEZmLEdnKSIsMix7InN0eWxlIjp7ImJvZHkiOnsibmFtZSI6ImJhcnJlZCJ9fX1dLFs0LDcsIlxcdGlsZGUgTV97ZixnfSIsMSx7InNob3J0ZW4iOnsic291cmNlIjoyMCwidGFyZ2V0IjoyMH0sInN0eWxlIjp7ImJvZHkiOnsibmFtZSI6Im5vbmUifSwiaGVhZCI6eyJuYW1lIjoibm9uZSJ9fX1dXQ==
        \begin{tikzcd}
          {\dbl{D}(x,y)} & {\dbl{D}(w,z)} \\
          {\dbl{E}(Fx,Gy)} & {\dbl{E}(Fw,Gz)}
          \arrow[""{name=0, anchor=center, inner sep=0}, "{\dbl{D}(f,g)}"{inner sep=.8ex}, "\shortmid"{marking}, from=1-1, to=1-2]
          \arrow["{\tilde M_{x,y}}"', from=1-1, to=2-1]
          \arrow["{\tilde M_{w,z}}", from=1-2, to=2-2]
          \arrow[""{name=1, anchor=center, inner sep=0}, "{\dbl{E}(Ff,Gg)}"'{inner sep=.8ex}, "\shortmid"{marking}, from=2-1, to=2-2]
          \arrow["{\tilde M_{f,g}}"{description}, draw=none, from=0, to=1]
        \end{tikzcd}
      \end{equation*}
      that sends a cell $\stdInlineCell{\alpha}$ in $\dbl{D}$ to its image
      $\inlineCell{Fx}{Gy}{Fw}{Gz}{Mm}{Mn}{Ff}{Gg}{M\alpha}$.
    \item The naturality comparison at a pair of proarrows $p: x' \proto x$ and
      $q: y \proto y'$ in $\dbl{D}$ is the natural transformation
      \begin{equation*}
        % https://q.uiver.app/#q=WzAsNixbMCwwLCJcXGRibHtEfSh4LHkpIl0sWzAsMSwiXFxkYmx7RX0oRngsR3kpIl0sWzAsMiwiXFxkYmx7RX0oRngnLEd5JykiXSxbMSwyLCJcXGRibHtFfShGeCcsR3knKSJdLFsxLDAsIlxcZGJse0R9KHgseSkiXSxbMSwxLCJcXGRibHtEfSh4Jyx5JykiXSxbMCwxLCJcXHRpbGRlIE1fe3gseX0iLDJdLFsxLDIsIlxcZGJse0V9KEZwLEdxKSIsMl0sWzIsMywiIiwyLHsibGV2ZWwiOjIsInN0eWxlIjp7ImJvZHkiOnsibmFtZSI6ImJhcnJlZCJ9LCJoZWFkIjp7Im5hbWUiOiJub25lIn19fV0sWzQsNSwiXFxkYmx7RH0ocCxxKSJdLFs1LDMsIlxcdGlsZGUgTV97eCcseSd9Il0sWzAsNCwiIiwwLHsibGV2ZWwiOjIsInN0eWxlIjp7ImJvZHkiOnsibmFtZSI6ImJhcnJlZCJ9LCJoZWFkIjp7Im5hbWUiOiJub25lIn19fV0sWzExLDgsIlxcdGlsZGUgTV97cCxxfSIsMSx7InNob3J0ZW4iOnsic291cmNlIjoyMCwidGFyZ2V0IjoyMH0sInN0eWxlIjp7ImJvZHkiOnsibmFtZSI6Im5vbmUifSwiaGVhZCI6eyJuYW1lIjoibm9uZSJ9fX1dXQ==
        \begin{tikzcd}
          {\dbl{D}(x,y)} & {\dbl{D}(x,y)} \\
          {\dbl{E}(Fx,Gy)} & {\dbl{D}(x',y')} \\
          {\dbl{E}(Fx',Gy')} & {\dbl{E}(Fx',Gy')}
          \arrow[""{name=0, anchor=center, inner sep=0}, "\shortmid"{marking}, equals, from=1-1, to=1-2]
          \arrow["{\tilde M_{x,y}}"', from=1-1, to=2-1]
          \arrow["{\dbl{D}(p,q)}", from=1-2, to=2-2]
          \arrow["{\dbl{E}(Fp,Gq)}"', from=2-1, to=3-1]
          \arrow["{\tilde M_{x',y'}}", from=2-2, to=3-2]
          \arrow[""{name=1, anchor=center, inner sep=0}, "\shortmid"{marking}, equals, from=3-1, to=3-2]
          \arrow["{\tilde M_{p,q}}"{description}, draw=none, from=0, to=1]
        \end{tikzcd}
      \end{equation*}
      whose component at a proarrow $m: x \proto y$ in $\dbl{D}$ is the tightly 
      globular cell in $\dbl{E}$
      \begin{equation*}
        (\tilde M_{p,q})_m \coloneqq M_{p,m,q}: Fp \odot Mm \odot Gq \to M(p \odot m \odot q)
      \end{equation*}
      given by applying the left and right actions of the module $M$ (in either
      order, equivalent by the associativity/compatibility axiom for the
      actions). \qedhere
  \end{itemize}
\end{construction}

This construction recovers the previous one by taking the identity module
$\id_F: F \proTo F$ on a lax double functor $F$. Specializing in a different
direction, when $F$ and $G$ are \emph{pseudo} double functors, a loose
transformation $\tau: F \proTo G$ is equivalently a pseudo double functor
$\tau: \dbl{D} \times \Loose \to \dbl{E}$ such that $\tau(-, 0) = F$ and $\tau(-, 1) = G$. In
this case, the induced natural transformation between twisted functors is pseudo
too.

\begin{lemma} \label{lem:module-to-twisted-transformation}
  The lax natural transformation $\tilde M$ between twisted functors induced by
  a module $M$ between lax double functors is well-defined.
\end{lemma}
\begin{proof}
  We follow \cite[Definition 9.1]{lambert2024} in naming the axioms of a module.

  That the component arrows and cells of $\tilde M$ \emph{are} functors and
  natural transformations, respectively, as well as the functorality of the
  component cells, all follows from the functorality on cells of the module $M$.
  That the naturality comparisons of $\tilde M$ \emph{are} natural
  transformations, as well as the naturality of $\tilde M$ with respect to
  cells, follows from the naturality of the left and right actions of $M$.
  Finally, the coherence of the naturality comparisons of $\tilde M$ follows
  from the associativity and unitality of the left/right actions of $M$.
\end{proof}

Having seen a few examples of natural transformations between twisted functors,
we return to the general formalism, exhibiting the twisted functor categories
whose morphisms are natural transformations.

\begin{lemma}[Category of twisted functors] \label{lem:cat-twisted-lax}
  Given a pair of double categories $\dbl{D}$ and $\dbl{E}$, there is a category
  of twisted doubly lax functors $\dbl{D} \twistto \dbl{E}$ and natural
  transformations between them, for any choice of lax, pseudo, or strict natural
  transformations.

  In this category, the composite $\alpha \cdot \beta: F \To H$ of natural
  transformations $\alpha: F \To G$ and $\beta: G \To H$ has components
  $(\alpha \cdot \beta)_x \coloneqq \alpha_x \cdot \beta_x$ for each $x$ in
  $\dbl{D}$ and
  \begin{equation} \label{eq:composite-natural-transformation-components}
    (\alpha \cdot \beta)_f
    \quad\coloneqq\quad
    % https://q.uiver.app/#q=WzAsNixbMCwwLCJGeCJdLFsxLDAsIkZ5Il0sWzAsMSwiR3giXSxbMSwxLCJHeSJdLFswLDIsIkh4Il0sWzEsMiwiSHkiXSxbMCwxLCJGZiIsMCx7InN0eWxlIjp7ImJvZHkiOnsibmFtZSI6ImJhcnJlZCJ9fX1dLFswLDIsIlxcYWxwaGFfeCIsMl0sWzIsMywiR2YiLDAseyJzdHlsZSI6eyJib2R5Ijp7Im5hbWUiOiJiYXJyZWQifX19XSxbMSwzLCJcXGFscGhhX3kiXSxbMiw0LCJcXGJldGFfeCIsMl0sWzMsNSwiXFxiZXRhX3kiXSxbNCw1LCJIZiIsMix7InN0eWxlIjp7ImJvZHkiOnsibmFtZSI6ImJhcnJlZCJ9fX1dLFs2LDgsIlxcYWxwaGFfZiIsMSx7ImxhYmVsX3Bvc2l0aW9uIjo0MCwic2hvcnRlbiI6eyJzb3VyY2UiOjIwLCJ0YXJnZXQiOjIwfSwic3R5bGUiOnsiYm9keSI6eyJuYW1lIjoibm9uZSJ9LCJoZWFkIjp7Im5hbWUiOiJub25lIn19fV0sWzgsMTIsIlxcYmV0YV9mIiwxLHsic2hvcnRlbiI6eyJzb3VyY2UiOjIwLCJ0YXJnZXQiOjIwfSwic3R5bGUiOnsiYm9keSI6eyJuYW1lIjoibm9uZSJ9LCJoZWFkIjp7Im5hbWUiOiJub25lIn19fV1d
    \begin{tikzcd}
      Fx & Fy \\
      Gx & Gy \\
      Hx & Hy
      \arrow[""{name=0, anchor=center, inner sep=0}, "Ff", "\shortmid"{marking}, from=1-1, to=1-2]
      \arrow["{\alpha_x}"', from=1-1, to=2-1]
      \arrow["{\alpha_y}", from=1-2, to=2-2]
      \arrow[""{name=1, anchor=center, inner sep=0}, "Gf", "\shortmid"{marking}, from=2-1, to=2-2]
      \arrow["{\beta_x}"', from=2-1, to=3-1]
      \arrow["{\beta_y}", from=2-2, to=3-2]
      \arrow[""{name=2, anchor=center, inner sep=0}, "Hf"', "\shortmid"{marking}, from=3-1, to=3-2]
      \arrow["{\alpha_f}"{description, pos=0.4}, draw=none, from=0, to=1]
      \arrow["{\beta_f}"{description}, draw=none, from=1, to=2]
    \end{tikzcd}
  \end{equation}
  for each arrow $f: x \to y$ in $\dbl{D}$, and naturality comparisons
  \begin{equation} \label{eq:composite-natural-transformation-comparisons}
    (\alpha \cdot \beta)_m
    \quad\coloneqq\quad
    \begin{tikzcd}[row sep=scriptsize]
      Fx && Fx \\
      Fx & Fx & Fx \\
      Gx & Gx & Fy \\
      Hx & Gy & Gy \\
      Hy & Hy & Hy \\
      Hy && Hy
      \arrow[""{name=0, anchor=center, inner sep=0}, "\shortmid"{marking}, equals, from=1-1, to=1-3]
      \arrow[equals, from=1-1, to=2-1]
      \arrow[equals, from=1-3, to=2-3]
      \arrow[""{name=1, anchor=center, inner sep=0}, "\shortmid"{marking}, equals, from=2-1, to=2-2]
      \arrow["{\alpha_x}"', from=2-1, to=3-1]
      \arrow[""{name=2, anchor=center, inner sep=0}, "\shortmid"{marking}, equals, from=2-2, to=2-3]
      \arrow["{\alpha_x}", from=2-2, to=3-2]
      \arrow["Fm", from=2-3, to=3-3]
      \arrow[""{name=3, anchor=center, inner sep=0}, "\shortmid"{marking}, equals, from=3-1, to=3-2]
      \arrow["{\beta_x}"', from=3-1, to=4-1]
      \arrow["Gm"', from=3-2, to=4-2]
      \arrow["{\alpha_y}", from=3-3, to=4-3]
      \arrow["Hm"', from=4-1, to=5-1]
      \arrow[""{name=4, anchor=center, inner sep=0}, "\shortmid"{marking}, equals, from=4-2, to=4-3]
      \arrow["{\beta_y}"', from=4-2, to=5-2]
      \arrow["{\beta_y}", from=4-3, to=5-3]
      \arrow[""{name=5, anchor=center, inner sep=0}, "\shortmid"{marking}, equals, from=5-1, to=5-2]
      \arrow[equals, from=5-1, to=6-1]
      \arrow[""{name=6, anchor=center, inner sep=0}, "\shortmid"{marking}, equals, from=5-2, to=5-3]
      \arrow[equals, from=5-3, to=6-3]
      \arrow[""{name=7, anchor=center, inner sep=0}, "\shortmid"{marking}, equals, from=6-1, to=6-3]
      \arrow["\cong"{description, pos=0.6}, draw=none, from=0, to=2-2]
      \arrow["{\id_{\alpha_x}}"{description}, draw=none, from=1, to=3]
      \arrow["{\alpha_m}"{description}, draw=none, from=2, to=4]
      \arrow["{\beta_m}"{description}, draw=none, from=3, to=5]
      \arrow["{\id_{\beta_y}}"{description}, draw=none, from=4, to=6]
      \arrow["\cong"{description, pos=0.4}, draw=none, from=5-2, to=7]
    \end{tikzcd}
  \end{equation}
  for each proarrow $m: x \proto y$ in $\dbl{D}$. The identity transformation is
  the strict natural transformation whose components are identities in
  $\dbl{E}_0$ and $\dbl{E}_1$.
\end{lemma}
\begin{proof}[Proof sketch]
  It is straightforward to verify that the composite and identity
  transformations are well defined. The associativity and unitality laws hold by
  the corresponding laws in both directions of the codomain double category
  $\dbl{E}$.
\end{proof}

\begin{remark} \label{rmk:on-cleavage-preservation}
  Now that we have a category of twisted functors, one might expect that the 
  elements correspondence (\cref{construction:weak-elements}) extends to a 
  functorial assignment into a category of loosely discrete opfibrations 
  (\cref{def:cat-of-opfibrations}). Although it is possible, we intend to leave 
  the full account of this correspondence to a further study. We pause here, 
  however, to explain why the category of loosely discrete opfibrations was chosen 
  to have morphisms that are not cleavage-preserving. This is, of course, a 
  standard move for fibration theories inasmuch as cleavages are \emph{extra structure} (\cite[\S B1.3]{johnstone2002}). 
  But now that we have a definition of morphism of twisted functors in hand, we 
  can see a technical obstacle as well. Given a lax transformation 
  $\alpha\colon F\Rightarrow G$ as above, an elements-style functor 
    \begin{equation*}
      \Elt(\alpha)\colon\Elt(F)\to\Elt(G)
    \end{equation*}  
  over the base would send $(x,a)\mapsto (x,\alpha_x(a))$ using the correct 
  component of the transformation. An assignment on proarrows would look like 
    \begin{equation*}
      (m,\overline m)\colon (x,a)\proto (y,b) \quad \mapsto \quad (m,\ast) \colon (x,\alpha_x(a))\proto (y,\alpha_y(b))
    \end{equation*}
  where the component arrow $\ast$ is defined as the composite 
    \begin{equation*}
      % https://q.uiver.app/#q=WzAsMyxbMCwwLCJHbVxcYWxwaGFfeChhKSJdLFswLDEsIlxcYWxwaGFfeUZtKGEpIl0sWzEsMSwiXFxhbHBoYV95KGIpIl0sWzAsMSwiKFxcYWxwaGFfbSlfe2EsYX0oMV9hKSIsMl0sWzEsMiwiXFxhbHBoYV95KFxcb3ZlcmxpbmUgbSkiLDJdLFswLDIsIlxcYXN0IiwwLHsic3R5bGUiOnsiYm9keSI6eyJuYW1lIjoiZGFzaGVkIn19fV1d
      \begin{tikzcd}
        {Gm\alpha_x(a)} & \\
        {\alpha_yFm(a)} & {\alpha_y(b)}
        \arrow["{(\alpha_m)_{a,a}(1_a)}"', from=1-1, to=2-1]
        \arrow["\ast", dashed, from=1-1, to=2-2]
        \arrow["{\alpha_y(\overline m)}"', from=2-1, to=2-2]
      \end{tikzcd}
    \end{equation*}
  where the vertically displayed arrow is the $(a,a)$-component of $\alpha_m$, 
  the naturality comparison, evaluated at the identity arrow $1_a$. This arrow 
  $\ast$ is invertible if $\alpha$ is pseudo but not generally an identity. The 
  canonical cleavage for $\Elt(F)$ has identity arrows in each second component, 
  as in $(m,1)$ with $\overline m =1$. So, what this remark shows is that 
  $\Elt(\alpha)$, under any reasonable definition, does not preserve the $1$ in 
  the second argument if $\alpha$ is lax or pseudo. Inasmuch as either should be 
  the basic notion of transformation between twisted functors, we therefore do 
  not require our morphisms of opfibrations preserve cleavages, lest the elements 
  construction land in the wrong category.
\end{remark}

Pre-composing an ordinary double functor onto a twisted functor
(\cref{construction:twisted-pre-composite}) extends to pre-whiskering a natural
transformation between twisted functors.

\begin{construction}[Pre-whiskering] \label{construction:twisted-pre-whiskering}
  Given a lax double functor $F : \dbl{C} \to \dbl{D}$ and a (lax, pseudo, or
  strict) natural transformation $\alpha : G \To H$ between twisted normal doubly lax
  functors $G,H : \dbl{D} \twistto \dbl{E}$, the \define{pre-whiskering}
  \begin{equation*}
    \alpha F : GF \To HF: \dbl{C} \twistto \dbl{E}
  \end{equation*}
  is the (correspondingly lax, pseudo, or strict) natural transformation
  constructed as follows:
  \begin{itemize}
    \item For each object $x \in \dbl{C}$, the component at $x$ is
      $(\alpha F)_{x} := \alpha_{F x} : GFx \to HFx$.
    \item For each arrow $f : x \to y$ in $\dbl{C}$, the component at $f$ is
      \begin{equation*}
        % https://q.uiver.app/#q=WzAsNCxbMCwwLCJHRngiXSxbMCwxLCJIRngiXSxbMSwwLCJHRnkiXSxbMSwxLCJIRnkiXSxbMCwxLCIoXFxhbHBoYSBGKV94IiwyXSxbMCwyLCJHRmYiLDAseyJzdHlsZSI6eyJib2R5Ijp7Im5hbWUiOiJiYXJyZWQifX19XSxbMiwzLCIoXFxhbHBoYSBGKV95Il0sWzEsMywiSEZmIiwyLHsic3R5bGUiOnsiYm9keSI6eyJuYW1lIjoiYmFycmVkIn19fV0sWzUsNywiKFxcYWxwaGEgRilfZiIsMSx7InNob3J0ZW4iOnsic291cmNlIjoyMCwidGFyZ2V0IjoyMH0sInN0eWxlIjp7ImJvZHkiOnsibmFtZSI6Im5vbmUifSwiaGVhZCI6eyJuYW1lIjoibm9uZSJ9fX1dXQ==
        \begin{tikzcd}
          GFx & GFy \\
          HFx & HFy
          \arrow[""{name=0, anchor=center, inner sep=0}, "GFf"{inner sep=.8ex}, "\shortmid"{marking}, from=1-1, to=1-2]
          \arrow["{(\alpha F)_x}"', from=1-1, to=2-1]
          \arrow["{(\alpha F)_y}", from=1-2, to=2-2]
          \arrow[""{name=1, anchor=center, inner sep=0}, "HFf"'{inner sep=.8ex}, "\shortmid"{marking}, from=2-1, to=2-2]
          \arrow["{(\alpha F)_f}"{description}, draw=none, from=0, to=1]
        \end{tikzcd}
        \quad\coloneqq\quad
        % https://q.uiver.app/#q=WzAsNCxbMCwwLCJHRngiXSxbMCwxLCJIRngiXSxbMSwwLCJHRnkiXSxbMSwxLCJIRnkiXSxbMCwxLCJcXGFscGhhX3tGeH0iLDJdLFswLDIsIkdGZiIsMCx7InN0eWxlIjp7ImJvZHkiOnsibmFtZSI6ImJhcnJlZCJ9fX1dLFsyLDMsIlxcYWxwaGFfe0Z5fSJdLFsxLDMsIkhGZiIsMix7InN0eWxlIjp7ImJvZHkiOnsibmFtZSI6ImJhcnJlZCJ9fX1dLFs1LDcsIlxcYWxwaGFfe0ZmfSIsMSx7InNob3J0ZW4iOnsic291cmNlIjoyMCwidGFyZ2V0IjoyMH0sInN0eWxlIjp7ImJvZHkiOnsibmFtZSI6Im5vbmUifSwiaGVhZCI6eyJuYW1lIjoibm9uZSJ9fX1dXQ==
        \begin{tikzcd}
          GFx & GFy \\
          HFx & HFy
          \arrow[""{name=0, anchor=center, inner sep=0}, "GFf"{inner sep=.8ex}, "\shortmid"{marking}, from=1-1, to=1-2]
          \arrow["{\alpha_{Fx}}"', from=1-1, to=2-1]
          \arrow["{\alpha_{Fy}}", from=1-2, to=2-2]
          \arrow[""{name=1, anchor=center, inner sep=0}, "HFf"'{inner sep=.8ex}, "\shortmid"{marking}, from=2-1, to=2-2]
          \arrow["{\alpha_{Ff}}"{description}, draw=none, from=0, to=1]
        \end{tikzcd}.
      \end{equation*}
    \item For each proarrow $m : x \proto y$ in $\dbl{C}$, the naturality comparison
      at $m$ is $(\alpha F)_{m} := \alpha_{Fm}$.
  \end{itemize}
\end{construction}

\begin{lemma}[Functorality of pre-whiskering] \label{lemma:pre-composition-functor}
  Let $F : \dbl{C} \to \dbl{D}$ be a lax double functor. Then the assignment
  \begin{equation*}
    (-)F: G \mapsto GF, \quad \alpha \mapsto \alpha F
  \end{equation*}
  sending a twisted normal doubly lax functor $G : \dbl{D} \twistto \dbl{E}$ to
  its composite $GF : \dbl{C} \twistto \dbl{E}$
  (\cref{construction:twisted-pre-composite}) and a natural transformation
  $\alpha: G \To H$ to its whiskering $\alpha F: GF \To HF$
  (\cref{construction:twisted-pre-whiskering}) is functorial.

  Moreover, given a double category $\dbl{E}$, there is a functor
  $|\Dbll|_1^{\op} \to \Cat$ sending a double category $\dbl{D}$ to the category
  of twisted normal doubly lax functors $\dbl{D} \twistto \dbl{E}$ and lax
  natural transformations between them, and acting with lax double functors by
  precomposition and pre-whiskering.
\end{lemma}
\begin{proof}
  We have already constructed a functor $|\Dbll|_1^{\op} \to \Set$ in
  \cref{lem:twisted-precomposition-functorality}. Upgrading to a category-valued
  presheaf is straightforward since composite and identity transformations are
  defined componentwise and pre-whiskering simply restricts the components. For
  example, for each arrow $f$ in the domain double category, we have
  \begin{equation*}
    ((\alpha \cdot \beta) F)_f
      = (\alpha \cdot \beta)_{Ff}
      = \alpha_{Ff} \cdot \beta_{Ff}
      = (\alpha F)_f \cdot (\beta F)_f
      = (\alpha F \cdot \beta F)_f
  \end{equation*}
  and
  \begin{equation*}
    (\alpha (G \circ F))_f
      = \alpha_{G(F(f))}
      = (\alpha G)_{Ff}
      = ((\alpha G) F)_f. \qedhere
  \end{equation*}
\end{proof}

\subsection{Modifications}

Modifications, like transformations, are defined to have components in the tight
direction of the target double category.

\begin{definition}[Modification] \label{def:modification-of-lax-transfs}
  Let $F,G: \dbl{D} \twistto \dbl{E}$ be twisted doubly lax functors and let
  $\alpha, \beta: F \To G$ be lax natural transformations between them. A
  \define{modification} $\mu: \alpha \Tto \beta$ consists of, for each object $x$ in
  $\dbl{D}$, its \define{component} at $x$, a cell
  \begin{equation*}
    % https://q.uiver.app/#q=WzAsNCxbMCwwLCJGeCJdLFswLDEsIkd4Il0sWzEsMCwiRngiXSxbMSwxLCJHeCJdLFswLDEsIlxcYWxwaGFfeCIsMl0sWzIsMywiXFxiZXRhX3giXSxbMCwyLCIiLDAseyJsZXZlbCI6Miwic3R5bGUiOnsiYm9keSI6eyJuYW1lIjoiYmFycmVkIn0sImhlYWQiOnsibmFtZSI6Im5vbmUifX19XSxbMSwzLCIiLDIseyJsZXZlbCI6Miwic3R5bGUiOnsiYm9keSI6eyJuYW1lIjoiYmFycmVkIn0sImhlYWQiOnsibmFtZSI6Im5vbmUifX19XSxbNiw3LCJcXG11X3giLDEseyJzaG9ydGVuIjp7InNvdXJjZSI6MjAsInRhcmdldCI6MjB9LCJzdHlsZSI6eyJib2R5Ijp7Im5hbWUiOiJub25lIn0sImhlYWQiOnsibmFtZSI6Im5vbmUifX19XV0=
    \begin{tikzcd}
      Fx & Fx \\
      Gx & Gx
      \arrow[""{name=0, anchor=center, inner sep=0}, "\shortmid"{marking}, Rightarrow, no head, from=1-1, to=1-2]
      \arrow["{\alpha_x}"', from=1-1, to=2-1]
      \arrow["{\beta_x}", from=1-2, to=2-2]
      \arrow[""{name=1, anchor=center, inner sep=0}, "\shortmid"{marking}, Rightarrow, no head, from=2-1, to=2-2]
      \arrow["{\mu_x}"{description}, draw=none, from=0, to=1]
    \end{tikzcd}
  \end{equation*}
  in $\dbl{E}$, satisfying two axioms.
  \begin{itemize}
    \item Tight-to-loose equivariance: for each arrow $f: x \to y$ in $\dbl{D}$,
      \begin{equation*}
        % https://q.uiver.app/#q=WzAsNixbMSwwLCJGeSJdLFsxLDEsIkd5Il0sWzIsMCwiRnkiXSxbMiwxLCJHeSJdLFswLDAsIkZ4Il0sWzAsMSwiR3kiXSxbMCwxLCJcXGFscGhhX3kiLDFdLFsyLDMsIlxcYmV0YV95Il0sWzAsMiwiIiwwLHsibGV2ZWwiOjIsInN0eWxlIjp7ImJvZHkiOnsibmFtZSI6ImJhcnJlZCJ9LCJoZWFkIjp7Im5hbWUiOiJub25lIn19fV0sWzEsMywiIiwyLHsibGV2ZWwiOjIsInN0eWxlIjp7ImJvZHkiOnsibmFtZSI6ImJhcnJlZCJ9LCJoZWFkIjp7Im5hbWUiOiJub25lIn19fV0sWzQsMCwiRmYiLDAseyJzdHlsZSI6eyJib2R5Ijp7Im5hbWUiOiJiYXJyZWQifX19XSxbNCw1LCJcXGFscGhhX3giLDJdLFs1LDEsIkdmIiwyLHsic3R5bGUiOnsiYm9keSI6eyJuYW1lIjoiYmFycmVkIn19fV0sWzgsOSwiXFxtdV95IiwxLHsic2hvcnRlbiI6eyJzb3VyY2UiOjIwLCJ0YXJnZXQiOjIwfSwic3R5bGUiOnsiYm9keSI6eyJuYW1lIjoibm9uZSJ9LCJoZWFkIjp7Im5hbWUiOiJub25lIn19fV0sWzEwLDEyLCJcXGFscGhhX2YiLDEseyJzaG9ydGVuIjp7InNvdXJjZSI6MjAsInRhcmdldCI6MjB9LCJzdHlsZSI6eyJib2R5Ijp7Im5hbWUiOiJub25lIn0sImhlYWQiOnsibmFtZSI6Im5vbmUifX19XV0=
        \begin{tikzcd}
          Fx & Fy & Fy \\
          Gy & Gy & Gy
          \arrow[""{name=0, anchor=center, inner sep=0}, "Ff", "\shortmid"{marking}, from=1-1, to=1-2]
          \arrow["{\alpha_x}"', from=1-1, to=2-1]
          \arrow[""{name=1, anchor=center, inner sep=0}, "\shortmid"{marking}, Rightarrow, no head, from=1-2, to=1-3]
          \arrow["{\alpha_y}"{description}, from=1-2, to=2-2]
          \arrow["{\beta_y}", from=1-3, to=2-3]
          \arrow[""{name=2, anchor=center, inner sep=0}, "Gf"', "\shortmid"{marking}, from=2-1, to=2-2]
          \arrow[""{name=3, anchor=center, inner sep=0}, "\shortmid"{marking}, Rightarrow, no head, from=2-2, to=2-3]
          \arrow["{\alpha_f}"{description}, draw=none, from=0, to=2]
          \arrow["{\mu_y}"{description}, draw=none, from=1, to=3]
        \end{tikzcd}
        \quad=\quad
        % https://q.uiver.app/#q=WzAsNixbMCwwLCJGeCJdLFswLDEsIkd4Il0sWzEsMCwiRngiXSxbMSwxLCJHeCJdLFsyLDAsIkZ5Il0sWzIsMSwiR3kiXSxbMCwxLCJcXGFscGhhX3giLDJdLFsyLDMsIlxcYmV0YV94IiwxXSxbMCwyLCIiLDAseyJsZXZlbCI6Miwic3R5bGUiOnsiYm9keSI6eyJuYW1lIjoiYmFycmVkIn0sImhlYWQiOnsibmFtZSI6Im5vbmUifX19XSxbMSwzLCIiLDIseyJsZXZlbCI6Miwic3R5bGUiOnsiYm9keSI6eyJuYW1lIjoiYmFycmVkIn0sImhlYWQiOnsibmFtZSI6Im5vbmUifX19XSxbMiw0LCJGZiIsMCx7InN0eWxlIjp7ImJvZHkiOnsibmFtZSI6ImJhcnJlZCJ9fX1dLFs0LDUsIlxcYmV0YV95Il0sWzMsNSwiR2YiLDIseyJzdHlsZSI6eyJib2R5Ijp7Im5hbWUiOiJiYXJyZWQifX19XSxbOCw5LCJcXG11X3giLDEseyJzaG9ydGVuIjp7InNvdXJjZSI6MjAsInRhcmdldCI6MjB9LCJzdHlsZSI6eyJib2R5Ijp7Im5hbWUiOiJub25lIn0sImhlYWQiOnsibmFtZSI6Im5vbmUifX19XSxbMTAsMTIsIlxcYmV0YV9mIiwxLHsic2hvcnRlbiI6eyJzb3VyY2UiOjIwLCJ0YXJnZXQiOjIwfSwic3R5bGUiOnsiYm9keSI6eyJuYW1lIjoibm9uZSJ9LCJoZWFkIjp7Im5hbWUiOiJub25lIn19fV1d
        \begin{tikzcd}
          Fx & Fx & Fy \\
          Gx & Gx & Gy
          \arrow[""{name=0, anchor=center, inner sep=0}, "\shortmid"{marking}, Rightarrow, no head, from=1-1, to=1-2]
          \arrow["{\alpha_x}"', from=1-1, to=2-1]
          \arrow[""{name=1, anchor=center, inner sep=0}, "Ff", "\shortmid"{marking}, from=1-2, to=1-3]
          \arrow["{\beta_x}"{description}, from=1-2, to=2-2]
          \arrow["{\beta_y}", from=1-3, to=2-3]
          \arrow[""{name=2, anchor=center, inner sep=0}, "\shortmid"{marking}, Rightarrow, no head, from=2-1, to=2-2]
          \arrow[""{name=3, anchor=center, inner sep=0}, "Gf"', "\shortmid"{marking}, from=2-2, to=2-3]
          \arrow["{\mu_x}"{description}, draw=none, from=0, to=2]
          \arrow["{\beta_f}"{description}, draw=none, from=1, to=3]
        \end{tikzcd}.
      \end{equation*}
    \item Loose-to-tight equivariance: for each proarrow $m: x \proto y$ in $\dbl{D}$,
      \begin{equation*}
        % https://q.uiver.app/#q=WzAsOSxbMSwxLCJGeSJdLFsxLDIsIkd5Il0sWzIsMSwiRnkiXSxbMiwyLCJHeSJdLFsxLDAsIkZ4Il0sWzIsMCwiRngiXSxbMCwwLCJGeCJdLFswLDEsIkd4Il0sWzAsMiwiR3kiXSxbMCwxLCJcXGFscGhhX3kiLDJdLFsyLDMsIlxcYmV0YV95Il0sWzAsMiwiIiwwLHsibGV2ZWwiOjIsInN0eWxlIjp7ImJvZHkiOnsibmFtZSI6ImJhcnJlZCJ9LCJoZWFkIjp7Im5hbWUiOiJub25lIn19fV0sWzEsMywiIiwyLHsibGV2ZWwiOjIsInN0eWxlIjp7ImJvZHkiOnsibmFtZSI6ImJhcnJlZCJ9LCJoZWFkIjp7Im5hbWUiOiJub25lIn19fV0sWzUsMiwiRm0iXSxbNCwwLCJGbSIsMl0sWzQsNSwiIiwwLHsibGV2ZWwiOjIsInN0eWxlIjp7ImJvZHkiOnsibmFtZSI6ImJhcnJlZCJ9LCJoZWFkIjp7Im5hbWUiOiJub25lIn19fV0sWzYsNywiXFxhbHBoYV94IiwyXSxbNyw4LCJHbSIsMl0sWzYsNCwiIiwyLHsibGV2ZWwiOjIsInN0eWxlIjp7ImJvZHkiOnsibmFtZSI6ImJhcnJlZCJ9LCJoZWFkIjp7Im5hbWUiOiJub25lIn19fV0sWzgsMSwiIiwyLHsibGV2ZWwiOjIsInN0eWxlIjp7ImJvZHkiOnsibmFtZSI6ImJhcnJlZCJ9LCJoZWFkIjp7Im5hbWUiOiJub25lIn19fV0sWzExLDEyLCJcXG11X3kiLDEseyJzaG9ydGVuIjp7InNvdXJjZSI6MjAsInRhcmdldCI6MjB9LCJzdHlsZSI6eyJib2R5Ijp7Im5hbWUiOiJub25lIn0sImhlYWQiOnsibmFtZSI6Im5vbmUifX19XSxbMTgsMTksIlxcYWxwaGFfbSIsMSx7InNob3J0ZW4iOnsic291cmNlIjoyMCwidGFyZ2V0IjoyMH0sInN0eWxlIjp7ImJvZHkiOnsibmFtZSI6Im5vbmUifSwiaGVhZCI6eyJuYW1lIjoibm9uZSJ9fX1dLFsxNSwxMSwiXFxpZF97Rm19IiwxLHsic2hvcnRlbiI6eyJzb3VyY2UiOjIwLCJ0YXJnZXQiOjIwfSwic3R5bGUiOnsiYm9keSI6eyJuYW1lIjoibm9uZSJ9LCJoZWFkIjp7Im5hbWUiOiJub25lIn19fV1d
        \begin{tikzcd}
          Fx & Fx & Fx \\
          Gx & Fy & Fy \\
          Gy & Gy & Gy
          \arrow[""{name=0, anchor=center, inner sep=0}, "\shortmid"{marking}, Rightarrow, no head, from=1-1, to=1-2]
          \arrow["{\alpha_x}"', from=1-1, to=2-1]
          \arrow[""{name=1, anchor=center, inner sep=0}, "\shortmid"{marking}, Rightarrow, no head, from=1-2, to=1-3]
          \arrow["Fm"', from=1-2, to=2-2]
          \arrow["Fm", from=1-3, to=2-3]
          \arrow["Gm"', from=2-1, to=3-1]
          \arrow[""{name=2, anchor=center, inner sep=0}, "\shortmid"{marking}, Rightarrow, no head, from=2-2, to=2-3]
          \arrow["{\alpha_y}"', from=2-2, to=3-2]
          \arrow["{\beta_y}", from=2-3, to=3-3]
          \arrow[""{name=3, anchor=center, inner sep=0}, "\shortmid"{marking}, Rightarrow, no head, from=3-1, to=3-2]
          \arrow[""{name=4, anchor=center, inner sep=0}, "\shortmid"{marking}, Rightarrow, no head, from=3-2, to=3-3]
          \arrow["{\alpha_m}"{description}, draw=none, from=0, to=3]
          \arrow["{\id_{Fm}}"{description}, draw=none, from=1, to=2]
          \arrow["{\mu_y}"{description}, draw=none, from=2, to=4]
        \end{tikzcd}
        \quad=\quad
        % https://q.uiver.app/#q=WzAsOSxbMCwwLCJGeCJdLFswLDEsIkd4Il0sWzEsMCwiRngiXSxbMSwxLCJHeCJdLFsyLDAsIkZ4Il0sWzEsMiwiR3kiXSxbMiwxLCJGeSJdLFsyLDIsIkd5Il0sWzAsMiwiR3kiXSxbMCwxLCJcXGFscGhhX3giLDJdLFsyLDMsIlxcYmV0YV94Il0sWzAsMiwiIiwwLHsibGV2ZWwiOjIsInN0eWxlIjp7ImJvZHkiOnsibmFtZSI6ImJhcnJlZCJ9LCJoZWFkIjp7Im5hbWUiOiJub25lIn19fV0sWzEsMywiIiwyLHsibGV2ZWwiOjIsInN0eWxlIjp7ImJvZHkiOnsibmFtZSI6ImJhcnJlZCJ9LCJoZWFkIjp7Im5hbWUiOiJub25lIn19fV0sWzMsNSwiR20iXSxbNCw2LCJGbSJdLFsyLDQsIiIsMCx7ImxldmVsIjoyLCJzdHlsZSI6eyJib2R5Ijp7Im5hbWUiOiJiYXJyZWQifSwiaGVhZCI6eyJuYW1lIjoibm9uZSJ9fX1dLFs2LDcsIlxcYmV0YV95Il0sWzUsNywiIiwxLHsibGV2ZWwiOjIsInN0eWxlIjp7ImJvZHkiOnsibmFtZSI6ImJhcnJlZCJ9LCJoZWFkIjp7Im5hbWUiOiJub25lIn19fV0sWzEsOCwiR20iLDJdLFs4LDUsIiIsMix7ImxldmVsIjoyLCJzdHlsZSI6eyJib2R5Ijp7Im5hbWUiOiJiYXJyZWQifSwiaGVhZCI6eyJuYW1lIjoibm9uZSJ9fX1dLFsxMSwxMiwiXFxtdV94IiwxLHsic2hvcnRlbiI6eyJzb3VyY2UiOjIwLCJ0YXJnZXQiOjIwfSwic3R5bGUiOnsiYm9keSI6eyJuYW1lIjoibm9uZSJ9LCJoZWFkIjp7Im5hbWUiOiJub25lIn19fV0sWzE1LDE3LCJcXGJldGFfbSIsMSx7InNob3J0ZW4iOnsic291cmNlIjoyMCwidGFyZ2V0IjoyMH0sInN0eWxlIjp7ImJvZHkiOnsibmFtZSI6Im5vbmUifSwiaGVhZCI6eyJuYW1lIjoibm9uZSJ9fX1dLFsxMiwxOSwiXFxpZF97R219IiwxLHsic2hvcnRlbiI6eyJzb3VyY2UiOjIwLCJ0YXJnZXQiOjIwfSwic3R5bGUiOnsiYm9keSI6eyJuYW1lIjoibm9uZSJ9LCJoZWFkIjp7Im5hbWUiOiJub25lIn19fV1d
        \begin{tikzcd}
          Fx & Fx & Fx \\
          Gx & Gx & Fy \\
          Gy & Gy & Gy
          \arrow[""{name=0, anchor=center, inner sep=0}, "\shortmid"{marking}, Rightarrow, no head, from=1-1, to=1-2]
          \arrow["{\alpha_x}"', from=1-1, to=2-1]
          \arrow[""{name=1, anchor=center, inner sep=0}, "\shortmid"{marking}, Rightarrow, no head, from=1-2, to=1-3]
          \arrow["{\beta_x}", from=1-2, to=2-2]
          \arrow["Fm", from=1-3, to=2-3]
          \arrow[""{name=2, anchor=center, inner sep=0}, "\shortmid"{marking}, Rightarrow, no head, from=2-1, to=2-2]
          \arrow["Gm"', from=2-1, to=3-1]
          \arrow["Gm", from=2-2, to=3-2]
          \arrow["{\beta_y}", from=2-3, to=3-3]
          \arrow[""{name=3, anchor=center, inner sep=0}, "\shortmid"{marking}, Rightarrow, no head, from=3-1, to=3-2]
          \arrow[""{name=4, anchor=center, inner sep=0}, "\shortmid"{marking}, Rightarrow, no head, from=3-2, to=3-3]
          \arrow["{\mu_x}"{description}, draw=none, from=0, to=2]
          \arrow["{\beta_m}"{description}, draw=none, from=1, to=4]
          \arrow["{\id_{Gm}}"{description}, draw=none, from=2, to=3]
        \end{tikzcd}.
        \qedhere
      \end{equation*}
  \end{itemize}
\end{definition}

The second axiom is the familiar \emph{modification axiom} or \emph{cylinder
  condition} from 2-category theory \cite[\S{4.4}]{johnson2021}. The first axiom
is specific to double categories and has a counterpart in the existing notion of
modification for non-twisted double functors \cite[\S{3.2.7}]{grandis2019}.

\begin{lemma}[2-category of twisted functors] \label{lem:2cat-twisted-functors}
  Given a pair of double categories $\dbl{D}$ and $\dbl{E}$, there is a
  2-category of twisted doubly lax functors $\dbl{D} \twistto \dbl{E}$, natural
  transformations, and modifications, for any choice of lax, pseudo, or strict
  transformations.

  In this 2-category, the vertical composite $\mu \cdot \nu: \alpha \Tto \gamma$
  of modifications $\mu: \alpha \Tto \beta$ and $\nu: \beta \Tto \gamma$ has
  components
  \begin{equation*}
    (\mu \cdot \nu)_x
    \quad\coloneqq\quad
    \begin{tikzcd}[row sep=scriptsize]
      Fx && Fx \\
      Fx & Fx & Fx \\
      Gx & Gx & Gx \\
      Gx && Gx
      \arrow[""{name=0, anchor=center, inner sep=0}, "\shortmid"{marking}, equals, from=1-1, to=1-3]
      \arrow[equals, from=1-1, to=2-1]
      \arrow[equals, from=1-3, to=2-3]
      \arrow[""{name=1, anchor=center, inner sep=0}, "\shortmid"{marking}, equals, from=2-1, to=2-2]
      \arrow["{\alpha_x}"', from=2-1, to=3-1]
      \arrow[""{name=2, anchor=center, inner sep=0}, "\shortmid"{marking}, equals, from=2-2, to=2-3]
      \arrow["{\beta_x}"{description}, from=2-2, to=3-2]
      \arrow["{\gamma_x}", from=2-3, to=3-3]
      \arrow[""{name=3, anchor=center, inner sep=0}, "\shortmid"{marking}, equals, from=3-1, to=3-2]
      \arrow[equals, from=3-1, to=4-1]
      \arrow[""{name=4, anchor=center, inner sep=0}, "\shortmid"{marking}, equals, from=3-2, to=3-3]
      \arrow[equals, from=3-3, to=4-3]
      \arrow[""{name=5, anchor=center, inner sep=0}, "\shortmid"{marking}, equals, from=4-1, to=4-3]
      \arrow["\cong"{description, pos=0.6}, draw=none, from=0, to=2-2]
      \arrow["{\mu_x}"{description}, draw=none, from=1, to=3]
      \arrow["{\nu_x}"{description}, draw=none, from=2, to=4]
      \arrow["\cong"{description, pos=0.4}, draw=none, from=3-2, to=5]
    \end{tikzcd}
  \end{equation*}
  for each $x \in \dbl{D}$. The horizontal composite
  $\mu * \nu: \alpha \cdot \gamma \Tto \beta \cdot \delta$ of modifications
  \begin{equation*}
    % https://q.uiver.app/#q=WzAsMyxbMCwwLCJGIl0sWzEsMCwiRyJdLFsyLDAsIkgiXSxbMCwxLCJcXGFscGhhIiwwLHsiY3VydmUiOi0yfV0sWzAsMSwiXFxiZXRhIiwyLHsiY3VydmUiOjJ9XSxbMSwyLCJcXGdhbW1hIiwwLHsiY3VydmUiOi0yfV0sWzEsMiwiXFxkZWx0YSIsMix7ImN1cnZlIjoyfV0sWzMsNCwiXFxtdSIsMix7InNob3J0ZW4iOnsic291cmNlIjozMCwidGFyZ2V0IjozMH19XSxbNSw2LCJcXG51IiwyLHsic2hvcnRlbiI6eyJzb3VyY2UiOjMwLCJ0YXJnZXQiOjMwfX1dXQ==
    \begin{tikzcd}
      F & G & H
      \arrow[""{name=0, anchor=center, inner sep=0}, "\alpha", curve={height=-12pt}, from=1-1, to=1-2]
      \arrow[""{name=1, anchor=center, inner sep=0}, "\beta"', curve={height=12pt}, from=1-1, to=1-2]
      \arrow[""{name=2, anchor=center, inner sep=0}, "\gamma", curve={height=-12pt}, from=1-2, to=1-3]
      \arrow[""{name=3, anchor=center, inner sep=0}, "\delta"', curve={height=12pt}, from=1-2, to=1-3]
      \arrow["\mu"', shorten <=5pt, shorten >=5pt, Rightarrow, from=0, to=1]
      \arrow["\nu"', shorten <=5pt, shorten >=5pt, Rightarrow, from=2, to=3]
    \end{tikzcd}
  \end{equation*}
  has components
  \begin{equation*}
    (\mu * \nu)_x
    \quad\coloneqq\quad
    % https://q.uiver.app/#q=WzAsNixbMCwwLCJGeCJdLFswLDEsIkd4Il0sWzAsMiwiSHgiXSxbMSwwLCJGeCJdLFsxLDEsIkd4Il0sWzEsMiwiSHgiXSxbMCwxLCJcXGFscGhhX3giLDJdLFsxLDIsIlxcZ2FtbWFfeCIsMl0sWzMsNCwiXFxiZXRhX3giXSxbMCwzLCIiLDAseyJsZXZlbCI6Miwic3R5bGUiOnsiYm9keSI6eyJuYW1lIjoiYmFycmVkIn0sImhlYWQiOnsibmFtZSI6Im5vbmUifX19XSxbMSw0LCIiLDAseyJsZXZlbCI6Miwic3R5bGUiOnsiYm9keSI6eyJuYW1lIjoiYmFycmVkIn0sImhlYWQiOnsibmFtZSI6Im5vbmUifX19XSxbNCw1LCJcXGRlbHRhX3giXSxbMiw1LCIiLDIseyJsZXZlbCI6Miwic3R5bGUiOnsiYm9keSI6eyJuYW1lIjoiYmFycmVkIn0sImhlYWQiOnsibmFtZSI6Im5vbmUifX19XSxbOSwxMCwiXFxtdV94IiwxLHsic2hvcnRlbiI6eyJzb3VyY2UiOjIwLCJ0YXJnZXQiOjIwfSwic3R5bGUiOnsiYm9keSI6eyJuYW1lIjoibm9uZSJ9LCJoZWFkIjp7Im5hbWUiOiJub25lIn19fV0sWzEwLDEyLCJcXG51X3giLDEseyJzaG9ydGVuIjp7InNvdXJjZSI6MjAsInRhcmdldCI6MjB9LCJzdHlsZSI6eyJib2R5Ijp7Im5hbWUiOiJub25lIn0sImhlYWQiOnsibmFtZSI6Im5vbmUifX19XV0=
    \begin{tikzcd}
      Fx & Fx \\
      Gx & Gx \\
      Hx & Hx
      \arrow[""{name=0, anchor=center, inner sep=0}, "\shortmid"{marking}, Rightarrow, no head, from=1-1, to=1-2]
      \arrow["{\alpha_x}"', from=1-1, to=2-1]
      \arrow["{\beta_x}", from=1-2, to=2-2]
      \arrow[""{name=1, anchor=center, inner sep=0}, "\shortmid"{marking}, Rightarrow, no head, from=2-1, to=2-2]
      \arrow["{\gamma_x}"', from=2-1, to=3-1]
      \arrow["{\delta_x}", from=2-2, to=3-2]
      \arrow[""{name=2, anchor=center, inner sep=0}, "\shortmid"{marking}, Rightarrow, no head, from=3-1, to=3-2]
      \arrow["{\mu_x}"{description}, draw=none, from=0, to=1]
      \arrow["{\nu_x}"{description}, draw=none, from=1, to=2]
    \end{tikzcd}
  \end{equation*}
  for each $x \in \dbl{D}$.
\end{lemma}

\begin{remark}[Invertibility]
  Note that a modification $\mu\colon \alpha\Rrightarrow\beta$ is invertible with respect to
  vertical composition in the 2-category of twisted functors if and only if each
  component $\mu_x$ is loosely invertible in the receiving double category
  $\dbl{E}$. For on the one hand, if $\mu$ is invertible in this sense, then
  certainly each component is; on the other hand, if each component is
  invertible, then these are easily seen to be the components of a modification
  $\mu^{-1}\colon \beta\Rrightarrow\alpha$ satisfying the two appropriate axioms above that is
  inverse to $\mu$ with respect to vertical composition. We can thus think of
  invertible modifications as being \emph{invertible componentwise}.
\end{remark}

With a 2-category of twisted functors now in hand, we can propose a notion of
twisted adjunction between double categories.

\begin{example}[Twisted adjunctions] \label{ex:twisted-adjunction}
  Continuing \cref{ex:twisted-companion-conjoint}, suppose that
  $F: \dbl{D} \leftrightarrows \dbl{E}: G$ are ordinary double functors. A \define{twisted
    double adjunction}, making $F$ be twisted left adjoint to $G$ and $G$ be
  twisted right adjoint to $F$, is an equivalence
  \begin{equation*}
    \dbl{E}(F(-),=) \simeq \dbl{D}(-,G(=))
  \end{equation*}
  in the 2-category of twisted double functors
  $\dbl{D}^\co \times \dbl{E} \twistto \Prof$, strong natural transformations,
  and modifications.

  Among other data and coherence equations, a twisted double adjunction gives a
  correspondence between proarrows
  \begin{equation*}
    F(x) \proto y \quad\leftrightarrows\quad x \proto G(y)
  \end{equation*}
  for each object $x \in \dbl{D}$ and $y \in \dbl{E}$, and, compatibly, a
  correspondence between cells
  \begin{equation*}
    % https://q.uiver.app/#q=WzAsNCxbMCwwLCJGKHgpIl0sWzEsMCwieSJdLFswLDEsIkYodykiXSxbMSwxLCJ6Il0sWzAsMiwiRihmKSIsMl0sWzEsMywiZyJdLFswLDEsIiIsMSx7InN0eWxlIjp7ImJvZHkiOnsibmFtZSI6ImJhcnJlZCJ9fX1dLFsyLDMsIiIsMSx7InN0eWxlIjp7ImJvZHkiOnsibmFtZSI6ImJhcnJlZCJ9fX1dLFs2LDcsIiIsMSx7InNob3J0ZW4iOnsic291cmNlIjozMCwidGFyZ2V0IjozMH19XV0=
    \begin{tikzcd}
      {F(x)} & y \\
      {F(w)} & z
      \arrow[""{name=0, anchor=center, inner sep=0}, "\shortmid"{marking}, from=1-1, to=1-2]
      \arrow["{F(f)}"', from=1-1, to=2-1]
      \arrow["g", from=1-2, to=2-2]
      \arrow[""{name=1, anchor=center, inner sep=0}, "\shortmid"{marking}, from=2-1, to=2-2]
      \arrow[between={0.3}{0.7}, Rightarrow, from=0, to=1]
    \end{tikzcd}
    \quad\leftrightarrows\quad
    % https://q.uiver.app/#q=WzAsNCxbMCwwLCJ4Il0sWzEsMCwiRyh5KSJdLFswLDEsInciXSxbMSwxLCJHKHopIl0sWzAsMSwiIiwzLHsic3R5bGUiOnsiYm9keSI6eyJuYW1lIjoiYmFycmVkIn19fV0sWzAsMiwiZiIsMl0sWzEsMywiRyhnKSJdLFsyLDMsIiIsMix7InN0eWxlIjp7ImJvZHkiOnsibmFtZSI6ImJhcnJlZCJ9fX1dLFs0LDcsIiIsMyx7InNob3J0ZW4iOnsic291cmNlIjozMCwidGFyZ2V0IjozMH19XV0=
    \begin{tikzcd}
      x & {G(y)} \\
      w & {G(z)}
      \arrow[""{name=0, anchor=center, inner sep=0}, "\shortmid"{marking}, from=1-1, to=1-2]
      \arrow["f"', from=1-1, to=2-1]
      \arrow["{G(g)}", from=1-2, to=2-2]
      \arrow[""{name=1, anchor=center, inner sep=0}, "\shortmid"{marking}, from=2-1, to=2-2]
      \arrow[between={0.3}{0.7}, Rightarrow, from=0, to=1]
    \end{tikzcd}
  \end{equation*}
  for each arrow $f: x \to w$ in $\dbl{D}$ and $g: y \to z$ in $\dbl{E}$. While
  twisted adjunctions and their application to loose compactness was our
  original motivation to study twisted functors, this definition is considered
  provisional, and a full accounting is deferred to future work.
\end{example}

\begin{construction}[Pre-whiskering] \label{construction:twisted-pre-2-whiskering}
  Let $F: \dbl{C} \to \dbl{D}$ be a lax double functor and let $\alpha, \beta : G \To H$ be
  lax natural transformations between twisted normal doubly lax functors
  $G,H : \dbl{D} \twistto \dbl{E}$. Given a modification $\mu : \alpha \Tto \beta$, the
  \define{pre-whiskering} $\mu F : \alpha F \Rrightarrow \beta F$ is the modification with
  components
  \begin{equation*}
    % https://q.uiver.app/#q=WzAsNCxbMCwwLCJHRngiXSxbMCwxLCJIRngiXSxbMSwwLCJHRngiXSxbMSwxLCJIR3giXSxbMCwxLCIoXFxhbHBoYSBGKV94IiwyXSxbMiwzLCIoXFxiZXRhIEYpX3giXSxbMCwyLCIiLDAseyJsZXZlbCI6Miwic3R5bGUiOnsiYm9keSI6eyJuYW1lIjoiYmFycmVkIn0sImhlYWQiOnsibmFtZSI6Im5vbmUifX19XSxbMSwzLCIiLDIseyJsZXZlbCI6Miwic3R5bGUiOnsiYm9keSI6eyJuYW1lIjoiYmFycmVkIn0sImhlYWQiOnsibmFtZSI6Im5vbmUifX19XSxbNiw3LCIoXFxtdSBGKV94IiwxLHsic2hvcnRlbiI6eyJzb3VyY2UiOjIwLCJ0YXJnZXQiOjIwfSwic3R5bGUiOnsiYm9keSI6eyJuYW1lIjoibm9uZSJ9LCJoZWFkIjp7Im5hbWUiOiJub25lIn19fV1d
    \begin{tikzcd}
      GFx & GFx \\
      HFx & HGx
      \arrow[""{name=0, anchor=center, inner sep=0}, "\shortmid"{marking}, equals, from=1-1, to=1-2]
      \arrow["{(\alpha F)_x}"', from=1-1, to=2-1]
      \arrow["{(\beta F)_x}", from=1-2, to=2-2]
      \arrow[""{name=1, anchor=center, inner sep=0}, "\shortmid"{marking}, equals, from=2-1, to=2-2]
      \arrow["{(\mu F)_x}"{description}, draw=none, from=0, to=1]
    \end{tikzcd}
    \quad\coloneqq\quad
    % https://q.uiver.app/#q=WzAsNCxbMCwwLCJHRngiXSxbMCwxLCJIRngiXSxbMSwwLCJHRngiXSxbMSwxLCJIR3giXSxbMCwxLCJcXGFscGhhX3tGeH0iLDJdLFsyLDMsIlxcYmV0YV97Rnh9Il0sWzAsMiwiIiwwLHsibGV2ZWwiOjIsInN0eWxlIjp7ImJvZHkiOnsibmFtZSI6ImJhcnJlZCJ9LCJoZWFkIjp7Im5hbWUiOiJub25lIn19fV0sWzEsMywiIiwyLHsibGV2ZWwiOjIsInN0eWxlIjp7ImJvZHkiOnsibmFtZSI6ImJhcnJlZCJ9LCJoZWFkIjp7Im5hbWUiOiJub25lIn19fV0sWzYsNywiXFxtdV97Rnh9IiwxLHsic2hvcnRlbiI6eyJzb3VyY2UiOjIwLCJ0YXJnZXQiOjIwfSwic3R5bGUiOnsiYm9keSI6eyJuYW1lIjoibm9uZSJ9LCJoZWFkIjp7Im5hbWUiOiJub25lIn19fV1d
    \begin{tikzcd}
      GFx & GFx \\
      HFx & HGx
      \arrow[""{name=0, anchor=center, inner sep=0}, "\shortmid"{marking}, equals, from=1-1, to=1-2]
      \arrow["{\alpha_{Fx}}"', from=1-1, to=2-1]
      \arrow["{\beta_{Fx}}", from=1-2, to=2-2]
      \arrow[""{name=1, anchor=center, inner sep=0}, "\shortmid"{marking}, equals, from=2-1, to=2-2]
      \arrow["{\mu_{Fx}}"{description}, draw=none, from=0, to=1]
    \end{tikzcd}
  \end{equation*}
  for each object $x \in \dbl{C}$.
\end{construction}

\begin{lemma}[2-functorality of pre-whiskering]
  Let $F : \dbl{C} \to \dbl{D}$ be a lax double functor. Then the functor $(-)F$
  of \cref{lemma:pre-composition-functor} sending a twisted normal doubly lax
  functor $G : \dbl{D} \twistto \dbl{E}$ to its composite
  $GF : \dbl{C} \twistto \dbl{E}$ extends to a 2-functor by prewhiskering
  modifications (\cref{construction:twisted-pre-2-whiskering}).
\end{lemma}
\begin{proof}
  This is again straightforward since the components of the prewhiskering $\mu F$
  are just the components of $\mu$ restricted along $F$.
\end{proof}

\subsection{Relative modifications}

While modifications have their uses, such as defining a notion of equivalence
between twisted functors (\cref{ex:twisted-adjunction}), modifications have not
made an appearance in our running example of twisted Homs. The problem is easy
to see: since a pair of double functors $F,G: \dbl{D} \to \dbl{E}$ induce natural
transformations between twisted functors
\begin{equation*}
  \tilde F: \dbl{D}(-, =) \To \dbl{E}(F(-), F(=))
  \qquad\text{and}\qquad
  \tilde G: \dbl{D}(-, =) \To \dbl{E}(G(-), G(=)),
\end{equation*}
one would expect a natural transformation $\alpha: F \To G$ to induce a
modification between $\tilde F$ and $\tilde G$, but that does not make sense
because $\tilde F$ and $\tilde G$ have different codomains. The problem is
solved by generalizing our notion of modification.

\begin{definition}[Relative modification]
\label{defn:relative.modification}
  Let $F: \dbl{C} \twistto \dbl{E}$ and $G: \dbl{D} \twistto \dbl{E}$ be twisted
  normal doubly lax functors, and let $(A, \alpha), (B, \beta): F \To G$ consist of lax
  functors $A, B: \dbl{C} \to \dbl{D}$ and lax natural transformations
  $\alpha: F \To GA$ and $\beta: F \To GB$. Given a natural transformation $\tau: A \To B$,
  a \define{modification} $\mu: \alpha \Tto \beta$ \define{relative to} $\tau$, or
  \define{$\tau$-relative modification} for short, consists of, for each object
  $x$ in $\dbl{C}$, a cell
  \begin{equation*}
    % https://q.uiver.app/#q=WzAsNCxbMCwwLCJGeCJdLFsxLDAsIkZ4Il0sWzAsMSwiR0F4Il0sWzEsMSwiR0J4Il0sWzAsMSwiIiwwLHsibGV2ZWwiOjIsInN0eWxlIjp7ImJvZHkiOnsibmFtZSI6ImJhcnJlZCJ9LCJoZWFkIjp7Im5hbWUiOiJub25lIn19fV0sWzAsMiwiXFxhbHBoYV94IiwyXSxbMSwzLCJcXGJldGFfeCJdLFsyLDMsIkdcXHRhdV94IiwyLHsic3R5bGUiOnsiYm9keSI6eyJuYW1lIjoiYmFycmVkIn19fV0sWzQsNywiXFxtdV94IiwxLHsic2hvcnRlbiI6eyJzb3VyY2UiOjIwLCJ0YXJnZXQiOjIwfSwic3R5bGUiOnsiYm9keSI6eyJuYW1lIjoibm9uZSJ9LCJoZWFkIjp7Im5hbWUiOiJub25lIn19fV1d
    \begin{tikzcd}
      Fx & Fx \\
      GAx & GBx
      \arrow[""{name=0, anchor=center, inner sep=0}, "\shortmid"{marking}, Rightarrow, no head, from=1-1, to=1-2]
      \arrow["{\alpha_x}"', from=1-1, to=2-1]
      \arrow["{\beta_x}", from=1-2, to=2-2]
      \arrow[""{name=1, anchor=center, inner sep=0}, "{G\tau_x}"', "\shortmid"{marking}, from=2-1, to=2-2]
      \arrow["{\mu_x}"{description}, draw=none, from=0, to=1]
    \end{tikzcd}
  \end{equation*}
  in $\dbl{E}$, called its \define{component} at $x$. The following two laws
  must be satisfied.
  \begin{itemize}
    \item Tight-to-loose equivariance: for each arrow $f : x \to y$ in $\dbl{C}$,
    \[
      % https://q.uiver.app/#q=WzAsOCxbMSwwLCJGeSJdLFsxLDEsIkdBeSJdLFsyLDAsIkZ5Il0sWzIsMSwiR0J5Il0sWzAsMCwiRngiXSxbMCwxLCJHQXgiXSxbMCwyLCJHQXgiXSxbMiwyLCJHQnkiXSxbMCwxLCJcXGFscGhhX3kiLDFdLFsyLDMsIlxcYmV0YV95Il0sWzAsMiwiIiwwLHsibGV2ZWwiOjIsInN0eWxlIjp7ImJvZHkiOnsibmFtZSI6ImJhcnJlZCJ9LCJoZWFkIjp7Im5hbWUiOiJub25lIn19fV0sWzEsMywiR1xcdGF1X3kiLDIseyJzdHlsZSI6eyJib2R5Ijp7Im5hbWUiOiJiYXJyZWQifX19XSxbNCwwLCJGZiIsMCx7InN0eWxlIjp7ImJvZHkiOnsibmFtZSI6ImJhcnJlZCJ9fX1dLFs0LDUsIlxcYWxwaGFfeCIsMl0sWzUsMSwiR0FmIiwyLHsic3R5bGUiOnsiYm9keSI6eyJuYW1lIjoiYmFycmVkIn19fV0sWzYsNywiRyhBZiBcXGNkb3QgXFx0YXVfeSkiLDIseyJzdHlsZSI6eyJib2R5Ijp7Im5hbWUiOiJiYXJyZWQifX19XSxbNSw2LCIiLDIseyJsZXZlbCI6Miwic3R5bGUiOnsiaGVhZCI6eyJuYW1lIjoibm9uZSJ9fX1dLFszLDcsIiIsMix7ImxldmVsIjoyLCJzdHlsZSI6eyJoZWFkIjp7Im5hbWUiOiJub25lIn19fV0sWzEwLDExLCJcXG11X3kiLDEseyJzaG9ydGVuIjp7InNvdXJjZSI6MjAsInRhcmdldCI6MjB9LCJzdHlsZSI6eyJib2R5Ijp7Im5hbWUiOiJub25lIn0sImhlYWQiOnsibmFtZSI6Im5vbmUifX19XSxbMTIsMTQsIlxcYWxwaGFfZiIsMSx7InNob3J0ZW4iOnsic291cmNlIjoyMCwidGFyZ2V0IjoyMH0sInN0eWxlIjp7ImJvZHkiOnsibmFtZSI6Im5vbmUifSwiaGVhZCI6eyJuYW1lIjoibm9uZSJ9fX1dLFsxLDE1LCJHX3tBZixcXHRhdV95fSIsMSx7InNob3J0ZW4iOnsidGFyZ2V0IjoyMH0sInN0eWxlIjp7ImJvZHkiOnsibmFtZSI6Im5vbmUifSwiaGVhZCI6eyJuYW1lIjoibm9uZSJ9fX1dXQ==
      \begin{tikzcd}
        Fx & Fy & Fy \\
        GAx & GAy & GBy \\
        GAx && GBy
        \arrow[""{name=0, anchor=center, inner sep=0}, "Ff"{inner sep=.8ex}, "\shortmid"{marking}, from=1-1, to=1-2]
        \arrow["{\alpha_x}"', from=1-1, to=2-1]
        \arrow[""{name=1, anchor=center, inner sep=0}, "\shortmid"{marking}, equals, from=1-2, to=1-3]
        \arrow["{\alpha_y}"{description}, from=1-2, to=2-2]
        \arrow["{\beta_y}", from=1-3, to=2-3]
        \arrow[""{name=2, anchor=center, inner sep=0}, "GAf"'{inner sep=.8ex}, "\shortmid"{marking}, from=2-1, to=2-2]
        \arrow[equals, from=2-1, to=3-1]
        \arrow[""{name=3, anchor=center, inner sep=0}, "{G\tau_y}"'{inner sep=.8ex}, "\shortmid"{marking}, from=2-2, to=2-3]
        \arrow[equals, from=2-3, to=3-3]
        \arrow[""{name=4, anchor=center, inner sep=0}, "{G(Af \cdot \tau_y)}"'{inner sep=.8ex}, "\shortmid"{marking}, from=3-1, to=3-3]
        \arrow["{\alpha_f}"{description}, draw=none, from=0, to=2]
        \arrow["{\mu_y}"{description}, draw=none, from=1, to=3]
        \arrow["{G_{Af,\tau_y}}"{description}, draw=none, from=2-2, to=4]
      \end{tikzcd}
      \quad=\quad
      % https://q.uiver.app/#q=WzAsOCxbMCwwLCJGeCJdLFswLDEsIkdBeCJdLFsxLDAsIkZ4Il0sWzEsMSwiR0J4Il0sWzIsMCwiRnkiXSxbMiwxLCJHQnkiXSxbMCwyLCJHQXgiXSxbMiwyLCJHQnkiXSxbMCwxLCJcXGFscGhhX3giLDJdLFsyLDMsIlxcYmV0YV94IiwxXSxbMCwyLCIiLDAseyJsZXZlbCI6Miwic3R5bGUiOnsiYm9keSI6eyJuYW1lIjoiYmFycmVkIn0sImhlYWQiOnsibmFtZSI6Im5vbmUifX19XSxbMSwzLCJHXFx0YXVfeCIsMix7InN0eWxlIjp7ImJvZHkiOnsibmFtZSI6ImJhcnJlZCJ9fX1dLFsyLDQsIkZmIiwwLHsic3R5bGUiOnsiYm9keSI6eyJuYW1lIjoiYmFycmVkIn19fV0sWzQsNSwiXFxiZXRhX3kiXSxbMyw1LCJHQmYiLDIseyJzdHlsZSI6eyJib2R5Ijp7Im5hbWUiOiJiYXJyZWQifX19XSxbNiw3LCJHKFxcdGF1X3ggXFxjZG90IEJmKSIsMix7InN0eWxlIjp7ImJvZHkiOnsibmFtZSI6ImJhcnJlZCJ9fX1dLFsxLDYsIiIsMix7ImxldmVsIjoyLCJzdHlsZSI6eyJoZWFkIjp7Im5hbWUiOiJub25lIn19fV0sWzUsNywiIiwyLHsibGV2ZWwiOjIsInN0eWxlIjp7ImhlYWQiOnsibmFtZSI6Im5vbmUifX19XSxbMTAsMTEsIlxcbXVfeCIsMSx7InNob3J0ZW4iOnsic291cmNlIjoyMCwidGFyZ2V0IjoyMH0sInN0eWxlIjp7ImJvZHkiOnsibmFtZSI6Im5vbmUifSwiaGVhZCI6eyJuYW1lIjoibm9uZSJ9fX1dLFsxMiwxNCwiXFxiZXRhX2YiLDEseyJzaG9ydGVuIjp7InNvdXJjZSI6MjAsInRhcmdldCI6MjB9LCJzdHlsZSI6eyJib2R5Ijp7Im5hbWUiOiJub25lIn0sImhlYWQiOnsibmFtZSI6Im5vbmUifX19XSxbMywxNSwiR197XFx0YXVfeCxCZn0iLDEseyJsZXZlbCI6MSwic3R5bGUiOnsiYm9keSI6eyJuYW1lIjoibm9uZSJ9LCJoZWFkIjp7Im5hbWUiOiJub25lIn19fV1d
      \begin{tikzcd}
        Fx & Fx & Fy \\
        GAx & GBx & GBy \\
        GAx && GBy
        \arrow[""{name=0, anchor=center, inner sep=0}, "\shortmid"{marking}, equals, from=1-1, to=1-2]
        \arrow["{\alpha_x}"', from=1-1, to=2-1]
        \arrow[""{name=1, anchor=center, inner sep=0}, "Ff"{inner sep=.8ex}, "\shortmid"{marking}, from=1-2, to=1-3]
        \arrow["{\beta_x}"{description}, from=1-2, to=2-2]
        \arrow["{\beta_y}", from=1-3, to=2-3]
        \arrow[""{name=2, anchor=center, inner sep=0}, "{G\tau_x}"'{inner sep=.8ex}, "\shortmid"{marking}, from=2-1, to=2-2]
        \arrow[equals, from=2-1, to=3-1]
        \arrow[""{name=3, anchor=center, inner sep=0}, "GBf"'{inner sep=.8ex}, "\shortmid"{marking}, from=2-2, to=2-3]
        \arrow[equals, from=2-3, to=3-3]
        \arrow[""{name=4, anchor=center, inner sep=0}, "{G(\tau_x \cdot Bf)}"'{inner sep=.8ex}, "\shortmid"{marking}, from=3-1, to=3-3]
        \arrow["{\mu_x}"{description}, draw=none, from=0, to=2]
        \arrow["{\beta_f}"{description}, draw=none, from=1, to=3]
        \arrow["{G_{\tau_x,Bf}}"{description}, draw=none, from=2-2, to=4]
      \end{tikzcd},
    \]
    which is a well-typed equation by the naturality of $\tau$.
    \item Loose-to-tight equivariance: for each proarrow $m : x \proto y$ in $\dbl{C}$,
    \[
      % https://q.uiver.app/#q=WzAsOSxbMSwxLCJGeSJdLFsxLDIsIkdBeSJdLFsyLDEsIkZ5Il0sWzIsMiwiR0J5Il0sWzEsMCwiRngiXSxbMiwwLCJGeCJdLFswLDAsIkZ4Il0sWzAsMSwiR0F4Il0sWzAsMiwiR0F5Il0sWzAsMSwiXFxhbHBoYV95IiwyXSxbMiwzLCJcXGJldGFfeSJdLFswLDIsIiIsMCx7ImxldmVsIjoyLCJzdHlsZSI6eyJib2R5Ijp7Im5hbWUiOiJiYXJyZWQifSwiaGVhZCI6eyJuYW1lIjoibm9uZSJ9fX1dLFsxLDMsIkdcXHRhdV95IiwyLHsic3R5bGUiOnsiYm9keSI6eyJuYW1lIjoiYmFycmVkIn19fV0sWzUsMiwiRm0iXSxbNCwwLCJGbSIsMl0sWzQsNSwiIiwwLHsibGV2ZWwiOjIsInN0eWxlIjp7ImJvZHkiOnsibmFtZSI6ImJhcnJlZCJ9LCJoZWFkIjp7Im5hbWUiOiJub25lIn19fV0sWzYsNywiXFxhbHBoYV94IiwyXSxbNyw4LCJHQW0iLDJdLFs2LDQsIiIsMix7ImxldmVsIjoyLCJzdHlsZSI6eyJib2R5Ijp7Im5hbWUiOiJiYXJyZWQifSwiaGVhZCI6eyJuYW1lIjoibm9uZSJ9fX1dLFs4LDEsIiIsMix7ImxldmVsIjoyLCJzdHlsZSI6eyJib2R5Ijp7Im5hbWUiOiJiYXJyZWQifSwiaGVhZCI6eyJuYW1lIjoibm9uZSJ9fX1dLFsxMSwxMiwiXFxtdV95IiwxLHsic2hvcnRlbiI6eyJzb3VyY2UiOjIwLCJ0YXJnZXQiOjIwfSwic3R5bGUiOnsiYm9keSI6eyJuYW1lIjoibm9uZSJ9LCJoZWFkIjp7Im5hbWUiOiJub25lIn19fV0sWzE4LDE5LCJcXGFscGhhX20iLDEseyJzaG9ydGVuIjp7InNvdXJjZSI6MjAsInRhcmdldCI6MjB9LCJzdHlsZSI6eyJib2R5Ijp7Im5hbWUiOiJub25lIn0sImhlYWQiOnsibmFtZSI6Im5vbmUifX19XSxbMTUsMTEsIlxcaWRfe0ZtfSIsMSx7InNob3J0ZW4iOnsic291cmNlIjoyMCwidGFyZ2V0IjoyMH0sInN0eWxlIjp7ImJvZHkiOnsibmFtZSI6Im5vbmUifSwiaGVhZCI6eyJuYW1lIjoibm9uZSJ9fX1dXQ==
      \begin{tikzcd}
        Fx & Fx & Fx \\
        GAx & Fy & Fy \\
        GAy & GAy & GBy
        \arrow[""{name=0, anchor=center, inner sep=0}, "\shortmid"{marking}, Rightarrow, no head, from=1-1, to=1-2]
        \arrow["{\alpha_x}"', from=1-1, to=2-1]
        \arrow[""{name=1, anchor=center, inner sep=0}, "\shortmid"{marking}, Rightarrow, no head, from=1-2, to=1-3]
        \arrow["Fm"', from=1-2, to=2-2]
        \arrow["Fm", from=1-3, to=2-3]
        \arrow["GAm"', from=2-1, to=3-1]
        \arrow[""{name=2, anchor=center, inner sep=0}, "\shortmid"{marking}, Rightarrow, no head, from=2-2, to=2-3]
        \arrow["{\alpha_y}"', from=2-2, to=3-2]
        \arrow["{\beta_y}", from=2-3, to=3-3]
        \arrow[""{name=3, anchor=center, inner sep=0}, "\shortmid"{marking}, Rightarrow, no head, from=3-1, to=3-2]
        \arrow[""{name=4, anchor=center, inner sep=0}, "{G\tau_y}"', "\shortmid"{marking}, from=3-2, to=3-3]
        \arrow["{\alpha_m}"{description}, draw=none, from=0, to=3]
        \arrow["{\id_{Fm}}"{description}, draw=none, from=1, to=2]
        \arrow["{\mu_y}"{description}, draw=none, from=2, to=4]
      \end{tikzcd}
      \quad = \quad
      % https://q.uiver.app/#q=WzAsOSxbMCwwLCJGeCJdLFswLDEsIkdBeCJdLFsxLDAsIkZ4Il0sWzEsMSwiR0J4Il0sWzIsMCwiRngiXSxbMSwyLCJHQnkiXSxbMiwxLCJGeSJdLFsyLDIsIkd5Il0sWzAsMiwiR0F5Il0sWzAsMSwiXFxhbHBoYV94IiwyXSxbMiwzLCJcXGJldGFfeCJdLFswLDIsIiIsMCx7ImxldmVsIjoyLCJzdHlsZSI6eyJib2R5Ijp7Im5hbWUiOiJiYXJyZWQifSwiaGVhZCI6eyJuYW1lIjoibm9uZSJ9fX1dLFsxLDMsIkdcXHRhdV94IiwyLHsic3R5bGUiOnsiYm9keSI6eyJuYW1lIjoiYmFycmVkIn19fV0sWzMsNSwiR0JtIl0sWzQsNiwiRm0iXSxbMiw0LCIiLDAseyJsZXZlbCI6Miwic3R5bGUiOnsiYm9keSI6eyJuYW1lIjoiYmFycmVkIn0sImhlYWQiOnsibmFtZSI6Im5vbmUifX19XSxbNiw3LCJcXGJldGFfeSJdLFs1LDcsIiIsMSx7ImxldmVsIjoyLCJzdHlsZSI6eyJib2R5Ijp7Im5hbWUiOiJiYXJyZWQifSwiaGVhZCI6eyJuYW1lIjoibm9uZSJ9fX1dLFsxLDgsIkdBbSIsMl0sWzgsNSwiR1xcdGF1X3kiLDIseyJzdHlsZSI6eyJib2R5Ijp7Im5hbWUiOiJiYXJyZWQifX19XSxbMTEsMTIsIlxcbXVfeCIsMSx7InNob3J0ZW4iOnsic291cmNlIjoyMCwidGFyZ2V0IjoyMH0sInN0eWxlIjp7ImJvZHkiOnsibmFtZSI6Im5vbmUifSwiaGVhZCI6eyJuYW1lIjoibm9uZSJ9fX1dLFsxNSwxNywiXFxiZXRhX20iLDEseyJzaG9ydGVuIjp7InNvdXJjZSI6MjAsInRhcmdldCI6MjB9LCJzdHlsZSI6eyJib2R5Ijp7Im5hbWUiOiJub25lIn0sImhlYWQiOnsibmFtZSI6Im5vbmUifX19XSxbMTIsMTksIkdcXHRhdV9tIiwxLHsic2hvcnRlbiI6eyJzb3VyY2UiOjIwLCJ0YXJnZXQiOjIwfSwic3R5bGUiOnsiYm9keSI6eyJuYW1lIjoibm9uZSJ9LCJoZWFkIjp7Im5hbWUiOiJub25lIn19fV1d
      \begin{tikzcd}
        Fx & Fx & Fx \\
        GAx & GBx & Fy \\
        GAy & GBy & GBy
        \arrow[""{name=0, anchor=center, inner sep=0}, "\shortmid"{marking}, Rightarrow, no head, from=1-1, to=1-2]
        \arrow["{\alpha_x}"', from=1-1, to=2-1]
        \arrow[""{name=1, anchor=center, inner sep=0}, "\shortmid"{marking}, Rightarrow, no head, from=1-2, to=1-3]
        \arrow["{\beta_x}", from=1-2, to=2-2]
        \arrow["Fm", from=1-3, to=2-3]
        \arrow[""{name=2, anchor=center, inner sep=0}, "{G\tau_x}"', "\shortmid"{marking}, from=2-1, to=2-2]
        \arrow["GAm"', from=2-1, to=3-1]
        \arrow["GBm", from=2-2, to=3-2]
        \arrow["{\beta_y}", from=2-3, to=3-3]
        \arrow[""{name=3, anchor=center, inner sep=0}, "{G\tau_y}"', "\shortmid"{marking}, from=3-1, to=3-2]
        \arrow[""{name=4, anchor=center, inner sep=0}, "\shortmid"{marking}, Rightarrow, no head, from=3-2, to=3-3]
        \arrow["{\mu_x}"{description}, draw=none, from=0, to=2]
        \arrow["{\beta_m}"{description}, draw=none, from=1, to=4]
        \arrow["{G\tau_m}"{description}, draw=none, from=2, to=3]
      \end{tikzcd}.
      \qedhere
    \]
  \end{itemize}
\end{definition}

The anticipated 2-cell between twisted Homs is now straightforward to construct.

\begin{construction}[Relative modification between twisted Homs]
  \label{construction:transformation-to-relative-modification}
  A natural transformation $\alpha: F \To G$ between lax double functors
  $F, G: \dbl{D} \to \dbl{E}$ induces a relative modification
  \begin{equation*}
    (\alpha^\co \times \alpha, \tilde \alpha):
      (F^\co \times F, \tilde F) \Tto (G^\co \times G, \tilde G):
      \dbl{D}(-, =) \To \dbl{E}(-, =).
  \end{equation*}
  Its component at a pair of objects $x$ and $y$ in $\dbl{D}$ is the natural
  transformation
  \begin{equation*}
    % https://q.uiver.app/#q=WzAsNCxbMCwwLCJcXGRibHtEfSh4LHkpIl0sWzEsMCwiXFxkYmx7RH0oeCx5KSJdLFswLDEsIlxcZGJse0V9KEZ4LEZ5KSJdLFsxLDEsIlxcZGJse0V9KEd4LEd5KSJdLFswLDEsIiIsMCx7ImxldmVsIjoyLCJzdHlsZSI6eyJib2R5Ijp7Im5hbWUiOiJiYXJyZWQifSwiaGVhZCI6eyJuYW1lIjoibm9uZSJ9fX1dLFswLDIsIlxcdGlsZGUgRl97eCx5fSIsMl0sWzEsMywiXFx0aWxkZSBHX3t4LHl9Il0sWzIsMywiXFxkYmx7RX0oXFxhbHBoYV94LFxcYWxwaGFfeSkiLDIseyJzdHlsZSI6eyJib2R5Ijp7Im5hbWUiOiJiYXJyZWQifX19XSxbNCw3LCJcXHRpbGRlIFxcYWxwaGFfe3gseX0iLDEseyJzaG9ydGVuIjp7InNvdXJjZSI6MjAsInRhcmdldCI6MjB9LCJzdHlsZSI6eyJib2R5Ijp7Im5hbWUiOiJub25lIn0sImhlYWQiOnsibmFtZSI6Im5vbmUifX19XV0=
    \begin{tikzcd}
      {\dbl{D}(x,y)} & {\dbl{D}(x,y)} \\
      {\dbl{E}(Fx,Fy)} & {\dbl{E}(Gx,Gy)}
      \arrow[""{name=0, anchor=center, inner sep=0}, "\shortmid"{marking}, equals, from=1-1, to=1-2]
      \arrow["{\tilde F_{x,y}}"', from=1-1, to=2-1]
      \arrow["{\tilde G_{x,y}}", from=1-2, to=2-2]
      \arrow[""{name=1, anchor=center, inner sep=0}, "{\dbl{E}(\alpha_x,\alpha_y)}"'{inner sep=.8ex}, "\shortmid"{marking}, from=2-1, to=2-2]
      \arrow["{\tilde \alpha_{x,y}}"{description}, draw=none, from=0, to=1]
    \end{tikzcd}
  \end{equation*}
  that acts on tightly globular cells in $\dbl{D}(x,y)$ as
  \begin{equation*}
    % https://q.uiver.app/#q=WzAsNCxbMCwwLCJ4Il0sWzEsMCwieSJdLFswLDEsIngiXSxbMSwxLCJ5Il0sWzAsMSwibSIsMCx7InN0eWxlIjp7ImJvZHkiOnsibmFtZSI6ImJhcnJlZCJ9fX1dLFsxLDMsIiIsMCx7ImxldmVsIjoyLCJzdHlsZSI6eyJoZWFkIjp7Im5hbWUiOiJub25lIn19fV0sWzAsMiwiIiwyLHsibGV2ZWwiOjIsInN0eWxlIjp7ImhlYWQiOnsibmFtZSI6Im5vbmUifX19XSxbMiwzLCJuIiwyLHsic3R5bGUiOnsiYm9keSI6eyJuYW1lIjoiYmFycmVkIn19fV0sWzQsNywiXFxnYW1tYSIsMSx7InNob3J0ZW4iOnsic291cmNlIjoyMCwidGFyZ2V0IjoyMH0sInN0eWxlIjp7ImJvZHkiOnsibmFtZSI6Im5vbmUifSwiaGVhZCI6eyJuYW1lIjoibm9uZSJ9fX1dXQ==
    \begin{tikzcd}
      x & y \\
      x & y
      \arrow[""{name=0, anchor=center, inner sep=0}, "m"{inner sep=.8ex}, "\shortmid"{marking}, from=1-1, to=1-2]
      \arrow[equals, from=1-1, to=2-1]
      \arrow[equals, from=1-2, to=2-2]
      \arrow[""{name=1, anchor=center, inner sep=0}, "n"'{inner sep=.8ex}, "\shortmid"{marking}, from=2-1, to=2-2]
      \arrow["\gamma"{description}, draw=none, from=0, to=1]
    \end{tikzcd}
    \qquad\mapsto\qquad
    % https://q.uiver.app/#q=WzAsNixbMCwwLCJGeCJdLFsxLDAsIkZ5Il0sWzAsMSwiRngiXSxbMSwxLCJGeSJdLFswLDIsIkd4Il0sWzEsMiwiR3kiXSxbMCwxLCJGbSIsMCx7InN0eWxlIjp7ImJvZHkiOnsibmFtZSI6ImJhcnJlZCJ9fX1dLFsxLDMsIiIsMCx7ImxldmVsIjoyLCJzdHlsZSI6eyJoZWFkIjp7Im5hbWUiOiJub25lIn19fV0sWzAsMiwiIiwyLHsibGV2ZWwiOjIsInN0eWxlIjp7ImhlYWQiOnsibmFtZSI6Im5vbmUifX19XSxbMiwzLCJGbiIsMix7InN0eWxlIjp7ImJvZHkiOnsibmFtZSI6ImJhcnJlZCJ9fX1dLFs0LDUsIkduIiwyLHsic3R5bGUiOnsiYm9keSI6eyJuYW1lIjoiYmFycmVkIn19fV0sWzIsNCwiXFxhbHBoYV94IiwyXSxbMyw1LCJcXGFscGhhX3kiXSxbNiw5LCJGXFxnYW1tYSIsMSx7InNob3J0ZW4iOnsic291cmNlIjoyMCwidGFyZ2V0IjoyMH0sInN0eWxlIjp7ImJvZHkiOnsibmFtZSI6Im5vbmUifSwiaGVhZCI6eyJuYW1lIjoibm9uZSJ9fX1dLFs5LDEwLCJcXGFscGhhX24iLDEseyJsYWJlbF9wb3NpdGlvbiI6NjAsInNob3J0ZW4iOnsic291cmNlIjoyMCwidGFyZ2V0IjoyMH0sInN0eWxlIjp7ImJvZHkiOnsibmFtZSI6Im5vbmUifSwiaGVhZCI6eyJuYW1lIjoibm9uZSJ9fX1dXQ==
    \begin{tikzcd}
      Fx & Fy \\
      Fx & Fy \\
      Gx & Gy
      \arrow[""{name=0, anchor=center, inner sep=0}, "Fm"{inner sep=.8ex}, "\shortmid"{marking}, from=1-1, to=1-2]
      \arrow[equals, from=1-1, to=2-1]
      \arrow[equals, from=1-2, to=2-2]
      \arrow[""{name=1, anchor=center, inner sep=0}, "Fn"'{inner sep=.8ex}, "\shortmid"{marking}, from=2-1, to=2-2]
      \arrow["{\alpha_x}"', from=2-1, to=3-1]
      \arrow["{\alpha_y}", from=2-2, to=3-2]
      \arrow[""{name=2, anchor=center, inner sep=0}, "Gn"'{inner sep=.8ex}, "\shortmid"{marking}, from=3-1, to=3-2]
      \arrow["{F\gamma}"{description}, draw=none, from=0, to=1]
      \arrow["{\alpha_n}"{description, pos=0.6}, draw=none, from=1, to=2]
    \end{tikzcd}
    =
    % https://q.uiver.app/#q=WzAsNixbMCwxLCJHeCJdLFsxLDEsIkd5Il0sWzAsMiwiR3giXSxbMSwyLCJHeSJdLFswLDAsIkZ4Il0sWzEsMCwiRnkiXSxbMCwxLCJHbSIsMCx7InN0eWxlIjp7ImJvZHkiOnsibmFtZSI6ImJhcnJlZCJ9fX1dLFsxLDMsIiIsMCx7ImxldmVsIjoyLCJzdHlsZSI6eyJoZWFkIjp7Im5hbWUiOiJub25lIn19fV0sWzAsMiwiIiwyLHsibGV2ZWwiOjIsInN0eWxlIjp7ImhlYWQiOnsibmFtZSI6Im5vbmUifX19XSxbMiwzLCJHbiIsMix7InN0eWxlIjp7ImJvZHkiOnsibmFtZSI6ImJhcnJlZCJ9fX1dLFs1LDEsIlxcYWxwaGFfeSJdLFs0LDAsIlxcYWxwaGFfeCIsMl0sWzQsNSwiRm0iLDAseyJzdHlsZSI6eyJib2R5Ijp7Im5hbWUiOiJiYXJyZWQifX19XSxbNiw5LCJHXFxnYW1tYSIsMSx7InNob3J0ZW4iOnsic291cmNlIjoyMCwidGFyZ2V0IjoyMH0sInN0eWxlIjp7ImJvZHkiOnsibmFtZSI6Im5vbmUifSwiaGVhZCI6eyJuYW1lIjoibm9uZSJ9fX1dLFsxMiw2LCJcXGFscGhhX20iLDEseyJsYWJlbF9wb3NpdGlvbiI6NDAsInNob3J0ZW4iOnsic291cmNlIjoyMCwidGFyZ2V0IjoyMH0sInN0eWxlIjp7ImJvZHkiOnsibmFtZSI6Im5vbmUifSwiaGVhZCI6eyJuYW1lIjoibm9uZSJ9fX1dXQ==
    \begin{tikzcd}
      Fx & Fy \\
      Gx & Gy \\
      Gx & Gy
      \arrow[""{name=0, anchor=center, inner sep=0}, "Fm"{inner sep=.8ex}, "\shortmid"{marking}, from=1-1, to=1-2]
      \arrow["{\alpha_x}"', from=1-1, to=2-1]
      \arrow["{\alpha_y}", from=1-2, to=2-2]
      \arrow[""{name=1, anchor=center, inner sep=0}, "Gm"{inner sep=.8ex}, "\shortmid"{marking}, from=2-1, to=2-2]
      \arrow[equals, from=2-1, to=3-1]
      \arrow[equals, from=2-2, to=3-2]
      \arrow[""{name=2, anchor=center, inner sep=0}, "Gn"'{inner sep=.8ex}, "\shortmid"{marking}, from=3-1, to=3-2]
      \arrow["{\alpha_m}"{description, pos=0.4}, draw=none, from=0, to=1]
      \arrow["{G\gamma}"{description}, draw=none, from=1, to=2]
    \end{tikzcd},
  \end{equation*}
  where the equality on the right uses the naturality of $\alpha$ with respect
  to cells.
\end{construction}

\begin{lemma}
  The relative modification $\tilde \alpha$ induced by a natural transformation
  $\alpha$ between lax double functors is well-defined.
\end{lemma}
\begin{proof}
  A special case of \cref{lem:modulation-to-relative-modification} below.
\end{proof}

If the construction of a transformation between twisted Homs generalizes from
lax double functors to modules, then the construction of a relative modification
should generalize from tight transformations between lax functors to morphisms
between modules, which are called \emph{modulations}. A \define{modulation}
between modules $M,N: \dbl{D} \times \Loose \to \dbl{E}$ can be succinctly
defined as a (tight) natural transformation $\mu: M \To N$. For an unpacking of
this definition, see the original source \cite[\mbox{Definition 3.3}]{pare2011}
or \cite[\mbox{Definition 9.7}]{lambert2024}.

\begin{construction}[Relative modification between twisted Homs, generalized]
  \label{construction:modulation-to-relative-modification}
  Let $F,G,H,K: \dbl{D} \to \dbl{E}$ be lax double functors, let
  $\alpha: F \To H$ and $\beta: G \To K$ be natural transformations, and let
  $M: F \proTo G$ and $N: H \proTo K$ be modules. Then a modulation
  \begin{equation*}
    % https://q.uiver.app/#q=WzAsNCxbMCwwLCJGIl0sWzEsMCwiRyJdLFswLDEsIkgiXSxbMSwxLCJLIl0sWzAsMSwiTSIsMCx7InN0eWxlIjp7ImJvZHkiOnsibmFtZSI6ImJhcnJlZCJ9fX1dLFswLDIsIlxcYWxwaGEiLDJdLFsxLDMsIlxcYmV0YSJdLFsyLDMsIk4iLDIseyJzdHlsZSI6eyJib2R5Ijp7Im5hbWUiOiJiYXJyZWQifX19XSxbNCw3LCJcXG11IiwxLHsic2hvcnRlbiI6eyJzb3VyY2UiOjIwLCJ0YXJnZXQiOjIwfSwic3R5bGUiOnsiYm9keSI6eyJuYW1lIjoibm9uZSJ9LCJoZWFkIjp7Im5hbWUiOiJub25lIn19fV1d
    \begin{tikzcd}
      F & G \\
      H & K
      \arrow[""{name=0, anchor=center, inner sep=0}, "M"{inner sep=.8ex}, "\shortmid"{marking}, from=1-1, to=1-2]
      \arrow["\alpha"', from=1-1, to=2-1]
      \arrow["\beta", from=1-2, to=2-2]
      \arrow[""{name=1, anchor=center, inner sep=0}, "N"'{inner sep=.8ex}, "\shortmid"{marking}, from=2-1, to=2-2]
      \arrow["\mu"{description}, draw=none, from=0, to=1]
    \end{tikzcd}
  \end{equation*}
  induces a relative modification
  \begin{equation*}
    (\alpha^\co \times \beta, \tilde \mu):
      (F^\co \times G, \tilde M) \Tto (H^\co \times K, \tilde N):
      \dbl{D}(-, =) \To \dbl{E}(-, =).
  \end{equation*}
  Its component at a pair of objects $x$ and $y$ in $\dbl{D}$ is the natural
  transformation
  \begin{equation*}
    % https://q.uiver.app/#q=WzAsNCxbMCwwLCJcXGRibHtEfSh4LHkpIl0sWzEsMCwiXFxkYmx7RH0oeCx5KSJdLFswLDEsIlxcZGJse0V9KEZ4LEd5KSJdLFsxLDEsIlxcZGJse0V9KEh4LEt5KSJdLFswLDEsIiIsMCx7ImxldmVsIjoyLCJzdHlsZSI6eyJib2R5Ijp7Im5hbWUiOiJiYXJyZWQifSwiaGVhZCI6eyJuYW1lIjoibm9uZSJ9fX1dLFswLDIsIlxcdGlsZGUgTV97eCx5fSIsMl0sWzEsMywiXFx0aWxkZSBOX3t4LHl9Il0sWzIsMywiXFxkYmx7RX0oXFxhbHBoYV94LFxcYmV0YV95KSIsMix7InN0eWxlIjp7ImJvZHkiOnsibmFtZSI6ImJhcnJlZCJ9fX1dLFs0LDcsIlxcdGlsZGUgXFxtdV97eCx5fSIsMSx7InNob3J0ZW4iOnsic291cmNlIjoyMCwidGFyZ2V0IjoyMH0sInN0eWxlIjp7ImJvZHkiOnsibmFtZSI6Im5vbmUifSwiaGVhZCI6eyJuYW1lIjoibm9uZSJ9fX1dXQ==
    \begin{tikzcd}
      {\dbl{D}(x,y)} & {\dbl{D}(x,y)} \\
      {\dbl{E}(Fx,Gy)} & {\dbl{E}(Hx,Ky)}
      \arrow[""{name=0, anchor=center, inner sep=0}, "\shortmid"{marking}, equals, from=1-1, to=1-2]
      \arrow["{\tilde M_{x,y}}"', from=1-1, to=2-1]
      \arrow["{\tilde N_{x,y}}", from=1-2, to=2-2]
      \arrow[""{name=1, anchor=center, inner sep=0}, "{\dbl{E}(\alpha_x,\beta_y)}"'{inner sep=.8ex}, "\shortmid"{marking}, from=2-1, to=2-2]
      \arrow["{\tilde \mu_{x,y}}"{description}, draw=none, from=0, to=1]
    \end{tikzcd}
  \end{equation*}
  that acts on tightly globular cells in $\dbl{D}(x,y)$ as
  \begin{equation*}
    % https://q.uiver.app/#q=WzAsNCxbMCwwLCJ4Il0sWzEsMCwieSJdLFswLDEsIngiXSxbMSwxLCJ5Il0sWzAsMSwibSIsMCx7InN0eWxlIjp7ImJvZHkiOnsibmFtZSI6ImJhcnJlZCJ9fX1dLFsxLDMsIiIsMCx7ImxldmVsIjoyLCJzdHlsZSI6eyJoZWFkIjp7Im5hbWUiOiJub25lIn19fV0sWzAsMiwiIiwyLHsibGV2ZWwiOjIsInN0eWxlIjp7ImhlYWQiOnsibmFtZSI6Im5vbmUifX19XSxbMiwzLCJuIiwyLHsic3R5bGUiOnsiYm9keSI6eyJuYW1lIjoiYmFycmVkIn19fV0sWzQsNywiXFxnYW1tYSIsMSx7InNob3J0ZW4iOnsic291cmNlIjoyMCwidGFyZ2V0IjoyMH0sInN0eWxlIjp7ImJvZHkiOnsibmFtZSI6Im5vbmUifSwiaGVhZCI6eyJuYW1lIjoibm9uZSJ9fX1dXQ==
    \begin{tikzcd}
      x & y \\
      x & y
      \arrow[""{name=0, anchor=center, inner sep=0}, "m"{inner sep=.8ex}, "\shortmid"{marking}, from=1-1, to=1-2]
      \arrow[equals, from=1-1, to=2-1]
      \arrow[equals, from=1-2, to=2-2]
      \arrow[""{name=1, anchor=center, inner sep=0}, "n"'{inner sep=.8ex}, "\shortmid"{marking}, from=2-1, to=2-2]
      \arrow["\gamma"{description}, draw=none, from=0, to=1]
    \end{tikzcd}
    \qquad\mapsto\qquad
    % https://q.uiver.app/#q=WzAsNixbMCwwLCJGeCJdLFsxLDAsIkd5Il0sWzAsMSwiRngiXSxbMSwxLCJHeSJdLFswLDIsIkh4Il0sWzEsMiwiS3kiXSxbMCwxLCJNbSIsMCx7InN0eWxlIjp7ImJvZHkiOnsibmFtZSI6ImJhcnJlZCJ9fX1dLFsxLDMsIiIsMCx7ImxldmVsIjoyLCJzdHlsZSI6eyJoZWFkIjp7Im5hbWUiOiJub25lIn19fV0sWzAsMiwiIiwyLHsibGV2ZWwiOjIsInN0eWxlIjp7ImhlYWQiOnsibmFtZSI6Im5vbmUifX19XSxbMiwzLCJNbiIsMix7InN0eWxlIjp7ImJvZHkiOnsibmFtZSI6ImJhcnJlZCJ9fX1dLFs0LDUsIk5uIiwyLHsic3R5bGUiOnsiYm9keSI6eyJuYW1lIjoiYmFycmVkIn19fV0sWzIsNCwiXFxhbHBoYV94IiwyXSxbMyw1LCJcXGJldGFfeSJdLFs2LDksIk1cXGdhbW1hIiwxLHsic2hvcnRlbiI6eyJzb3VyY2UiOjIwLCJ0YXJnZXQiOjIwfSwic3R5bGUiOnsiYm9keSI6eyJuYW1lIjoibm9uZSJ9LCJoZWFkIjp7Im5hbWUiOiJub25lIn19fV0sWzksMTAsIlxcbXVfbiIsMSx7ImxhYmVsX3Bvc2l0aW9uIjo2MCwic2hvcnRlbiI6eyJzb3VyY2UiOjIwLCJ0YXJnZXQiOjIwfSwic3R5bGUiOnsiYm9keSI6eyJuYW1lIjoibm9uZSJ9LCJoZWFkIjp7Im5hbWUiOiJub25lIn19fV1d
    \begin{tikzcd}
      Fx & Gy \\
      Fx & Gy \\
      Hx & Ky
      \arrow[""{name=0, anchor=center, inner sep=0}, "Mm"{inner sep=.8ex}, "\shortmid"{marking}, from=1-1, to=1-2]
      \arrow[equals, from=1-1, to=2-1]
      \arrow[equals, from=1-2, to=2-2]
      \arrow[""{name=1, anchor=center, inner sep=0}, "Mn"'{inner sep=.8ex}, "\shortmid"{marking}, from=2-1, to=2-2]
      \arrow["{\alpha_x}"', from=2-1, to=3-1]
      \arrow["{\beta_y}", from=2-2, to=3-2]
      \arrow[""{name=2, anchor=center, inner sep=0}, "Nn"'{inner sep=.8ex}, "\shortmid"{marking}, from=3-1, to=3-2]
      \arrow["{M\gamma}"{description}, draw=none, from=0, to=1]
      \arrow["{\mu_n}"{description, pos=0.6}, draw=none, from=1, to=2]
    \end{tikzcd}
    =
    % https://q.uiver.app/#q=WzAsNixbMCwxLCJIeCJdLFsxLDEsIkt5Il0sWzAsMiwiSHgiXSxbMSwyLCJLeSJdLFswLDAsIkZ4Il0sWzEsMCwiR3kiXSxbMCwxLCJObSIsMCx7InN0eWxlIjp7ImJvZHkiOnsibmFtZSI6ImJhcnJlZCJ9fX1dLFsxLDMsIiIsMCx7ImxldmVsIjoyLCJzdHlsZSI6eyJoZWFkIjp7Im5hbWUiOiJub25lIn19fV0sWzAsMiwiIiwyLHsibGV2ZWwiOjIsInN0eWxlIjp7ImhlYWQiOnsibmFtZSI6Im5vbmUifX19XSxbMiwzLCJObiIsMix7InN0eWxlIjp7ImJvZHkiOnsibmFtZSI6ImJhcnJlZCJ9fX1dLFs1LDEsIlxcYmV0YV95Il0sWzQsMCwiXFxhbHBoYV94IiwyXSxbNCw1LCJNbSIsMCx7InN0eWxlIjp7ImJvZHkiOnsibmFtZSI6ImJhcnJlZCJ9fX1dLFsxMiw2LCJcXG11X20iLDEseyJsYWJlbF9wb3NpdGlvbiI6NDAsInNob3J0ZW4iOnsic291cmNlIjoyMCwidGFyZ2V0IjoyMH0sInN0eWxlIjp7ImJvZHkiOnsibmFtZSI6Im5vbmUifSwiaGVhZCI6eyJuYW1lIjoibm9uZSJ9fX1dLFs2LDksIk5cXGdhbW1hIiwxLHsic2hvcnRlbiI6eyJzb3VyY2UiOjIwLCJ0YXJnZXQiOjIwfSwic3R5bGUiOnsiYm9keSI6eyJuYW1lIjoibm9uZSJ9LCJoZWFkIjp7Im5hbWUiOiJub25lIn19fV1d
    \begin{tikzcd}
      Fx & Gy \\
      Hx & Ky \\
      Hx & Ky
      \arrow[""{name=0, anchor=center, inner sep=0}, "Mm"{inner sep=.8ex}, "\shortmid"{marking}, from=1-1, to=1-2]
      \arrow["{\alpha_x}"', from=1-1, to=2-1]
      \arrow["{\beta_y}", from=1-2, to=2-2]
      \arrow[""{name=1, anchor=center, inner sep=0}, "Nm"{inner sep=.8ex}, "\shortmid"{marking}, from=2-1, to=2-2]
      \arrow[equals, from=2-1, to=3-1]
      \arrow[equals, from=2-2, to=3-2]
      \arrow[""{name=2, anchor=center, inner sep=0}, "Nn"'{inner sep=.8ex}, "\shortmid"{marking}, from=3-1, to=3-2]
      \arrow["{\mu_m}"{description, pos=0.4}, draw=none, from=0, to=1]
      \arrow["{N\gamma}"{description}, draw=none, from=1, to=2]
    \end{tikzcd},
  \end{equation*}
  where the equality on the right uses the naturality of the modulation $\mu$.
\end{construction}

This construction recovers the previous one by taking the modulation to be the
loose identity $\id_\alpha: \id_F \Tto \id_G$ on a natural transformation
$\alpha: F \To G$ between lax double functors.

\begin{lemma} \label{lem:modulation-to-relative-modification}
  The modification $\tilde \mu$ relative to $\alpha^\co \times \beta$ induced by
  a modulation $\mu$ is well-defined.
\end{lemma}
\begin{proof}
  We use the terminology of \cite[\mbox{Definition 9.7}]{lambert2024} to name
  the modulation axioms.

  The tight-to-loose equivariance axiom of $\tilde \mu$ is exactly the naturality
  axiom of the modulation $\mu$ (holding for arbitrary cells in $\dbl{D}$, not
  just globular ones as in the construction), while the loose-to-tight
  equivariance of $\tilde \mu$ combines the left and right equivariance axioms of
  $\mu$ into one equation
  \begin{equation*}
    % https://q.uiver.app/#q=WzAsOCxbMCwwLCJGeCciXSxbMSwwLCJGeCJdLFsyLDAsIkd5Il0sWzMsMCwiR3knIl0sWzAsMSwiRngnIl0sWzMsMSwiR3knIl0sWzAsMiwiSHgnIl0sWzMsMiwiS3knIl0sWzAsMSwiRnAiLDAseyJzdHlsZSI6eyJib2R5Ijp7Im5hbWUiOiJiYXJyZWQifX19XSxbMSwyLCJNbSIsMCx7InN0eWxlIjp7ImJvZHkiOnsibmFtZSI6ImJhcnJlZCJ9fX1dLFsyLDMsIkdxIiwwLHsic3R5bGUiOnsiYm9keSI6eyJuYW1lIjoiYmFycmVkIn19fV0sWzQsNSwiTShwIFxcb2RvdCBtIFxcb2RvdCBxKSIsMix7InN0eWxlIjp7ImJvZHkiOnsibmFtZSI6ImJhcnJlZCJ9fX1dLFswLDQsIiIsMix7ImxldmVsIjoyLCJzdHlsZSI6eyJoZWFkIjp7Im5hbWUiOiJub25lIn19fV0sWzMsNSwiIiwwLHsibGV2ZWwiOjIsInN0eWxlIjp7ImhlYWQiOnsibmFtZSI6Im5vbmUifX19XSxbNCw2LCJcXGFscGhhX3t4J30iLDJdLFs1LDcsIlxcYmV0YV97eSd9Il0sWzYsNywiTihwIFxcb2RvdCBtIFxcb2RvdCBxKSIsMix7InN0eWxlIjp7ImJvZHkiOnsibmFtZSI6ImJhcnJlZCJ9fX1dLFs5LDExLCJNX3twLG0scX0iLDEseyJzaG9ydGVuIjp7InNvdXJjZSI6MjAsInRhcmdldCI6MjB9LCJzdHlsZSI6eyJib2R5Ijp7Im5hbWUiOiJub25lIn0sImhlYWQiOnsibmFtZSI6Im5vbmUifX19XSxbMTEsMTYsIlxcbXVfe3AgXFxvZG90IG0gXFxvZG90IHF9IiwxLHsibGFiZWxfcG9zaXRpb24iOjYwLCJzaG9ydGVuIjp7InNvdXJjZSI6MjAsInRhcmdldCI6MjB9LCJzdHlsZSI6eyJib2R5Ijp7Im5hbWUiOiJub25lIn0sImhlYWQiOnsibmFtZSI6Im5vbmUifX19XV0=
    \begin{tikzcd}
      {Fx'} & Fx & Gy & {Gy'} \\
      {Fx'} &&& {Gy'} \\
      {Hx'} &&& {Ky'}
      \arrow["Fp"{inner sep=.8ex}, "\shortmid"{marking}, from=1-1, to=1-2]
      \arrow[equals, from=1-1, to=2-1]
      \arrow[""{name=0, anchor=center, inner sep=0}, "Mm"{inner sep=.8ex}, "\shortmid"{marking}, from=1-2, to=1-3]
      \arrow["Gq"{inner sep=.8ex}, "\shortmid"{marking}, from=1-3, to=1-4]
      \arrow[equals, from=1-4, to=2-4]
      \arrow[""{name=1, anchor=center, inner sep=0}, "{M(p \odot m \odot q)}"'{inner sep=.8ex}, "\shortmid"{marking}, from=2-1, to=2-4]
      \arrow["{\alpha_{x'}}"', from=2-1, to=3-1]
      \arrow["{\beta_{y'}}", from=2-4, to=3-4]
      \arrow[""{name=2, anchor=center, inner sep=0}, "{N(p \odot m \odot q)}"'{inner sep=.8ex}, "\shortmid"{marking}, from=3-1, to=3-4]
      \arrow["{M_{p,m,q}}"{description}, draw=none, from=0, to=1]
      \arrow["{\mu_{p \odot m \odot q}}"{description, pos=0.6}, draw=none, from=1, to=2]
    \end{tikzcd}
    \quad=\quad
    % https://q.uiver.app/#q=WzAsMTAsWzAsMCwiRngnIl0sWzEsMCwiRngiXSxbMiwwLCJHeSJdLFszLDAsIkd5JyJdLFswLDIsIkh4JyJdLFszLDIsIkt5JyJdLFswLDEsIkh4JyJdLFszLDEsIkt5JyJdLFsxLDEsIkh4Il0sWzIsMSwiS3kiXSxbMCwxLCJGcCIsMCx7InN0eWxlIjp7ImJvZHkiOnsibmFtZSI6ImJhcnJlZCJ9fX1dLFsxLDIsIk1tIiwwLHsic3R5bGUiOnsiYm9keSI6eyJuYW1lIjoiYmFycmVkIn19fV0sWzIsMywiR3EiLDAseyJzdHlsZSI6eyJib2R5Ijp7Im5hbWUiOiJiYXJyZWQifX19XSxbNCw1LCJOKHAgXFxvZG90IG0gXFxvZG90IHEpIiwyLHsic3R5bGUiOnsiYm9keSI6eyJuYW1lIjoiYmFycmVkIn19fV0sWzAsNiwiXFxhbHBoYV97eCd9IiwyXSxbNiw0LCIiLDIseyJsZXZlbCI6Miwic3R5bGUiOnsiaGVhZCI6eyJuYW1lIjoibm9uZSJ9fX1dLFs3LDUsIiIsMCx7ImxldmVsIjoyLCJzdHlsZSI6eyJoZWFkIjp7Im5hbWUiOiJub25lIn19fV0sWzMsNywiXFxiZXRhX3t5J30iXSxbNiw4LCJIcCIsMix7InN0eWxlIjp7ImJvZHkiOnsibmFtZSI6ImJhcnJlZCJ9fX1dLFsxLDgsIlxcYWxwaGFfeCIsMV0sWzIsOSwiXFxiZXRhX3kiLDFdLFs4LDksIk5tIiwyLHsic3R5bGUiOnsiYm9keSI6eyJuYW1lIjoiYmFycmVkIn19fV0sWzksNywiS3EiLDIseyJzdHlsZSI6eyJib2R5Ijp7Im5hbWUiOiJiYXJyZWQifX19XSxbMTAsMTgsIlxcYWxwaGFfcCIsMSx7InNob3J0ZW4iOnsic291cmNlIjoyMCwidGFyZ2V0IjoyMH0sInN0eWxlIjp7ImJvZHkiOnsibmFtZSI6Im5vbmUifSwiaGVhZCI6eyJuYW1lIjoibm9uZSJ9fX1dLFsxMSwyMSwiXFxtdV9tIiwxLHsic2hvcnRlbiI6eyJzb3VyY2UiOjIwLCJ0YXJnZXQiOjIwfSwic3R5bGUiOnsiYm9keSI6eyJuYW1lIjoibm9uZSJ9LCJoZWFkIjp7Im5hbWUiOiJub25lIn19fV0sWzEyLDIyLCJcXGJldGFfcSIsMSx7InNob3J0ZW4iOnsic291cmNlIjoyMCwidGFyZ2V0IjoyMH0sInN0eWxlIjp7ImJvZHkiOnsibmFtZSI6Im5vbmUifSwiaGVhZCI6eyJuYW1lIjoibm9uZSJ9fX1dLFsyMSwxMywiTl97cCxtLHF9IiwxLHsibGFiZWxfcG9zaXRpb24iOjYwLCJzaG9ydGVuIjp7InNvdXJjZSI6MjAsInRhcmdldCI6MjB9LCJzdHlsZSI6eyJib2R5Ijp7Im5hbWUiOiJub25lIn0sImhlYWQiOnsibmFtZSI6Im5vbmUifX19XV0=
    \begin{tikzcd}
      {Fx'} & Fx & Gy & {Gy'} \\
      {Hx'} & Hx & Ky & {Ky'} \\
      {Hx'} &&& {Ky'}
      \arrow[""{name=0, anchor=center, inner sep=0}, "Fp"{inner sep=.8ex}, "\shortmid"{marking}, from=1-1, to=1-2]
      \arrow["{\alpha_{x'}}"', from=1-1, to=2-1]
      \arrow[""{name=1, anchor=center, inner sep=0}, "Mm"{inner sep=.8ex}, "\shortmid"{marking}, from=1-2, to=1-3]
      \arrow["{\alpha_x}"{description}, from=1-2, to=2-2]
      \arrow[""{name=2, anchor=center, inner sep=0}, "Gq"{inner sep=.8ex}, "\shortmid"{marking}, from=1-3, to=1-4]
      \arrow["{\beta_y}"{description}, from=1-3, to=2-3]
      \arrow["{\beta_{y'}}", from=1-4, to=2-4]
      \arrow[""{name=3, anchor=center, inner sep=0}, "Hp"'{inner sep=.8ex}, "\shortmid"{marking}, from=2-1, to=2-2]
      \arrow[equals, from=2-1, to=3-1]
      \arrow[""{name=4, anchor=center, inner sep=0}, "Nm"'{inner sep=.8ex}, "\shortmid"{marking}, from=2-2, to=2-3]
      \arrow[""{name=5, anchor=center, inner sep=0}, "Kq"'{inner sep=.8ex}, "\shortmid"{marking}, from=2-3, to=2-4]
      \arrow[equals, from=2-4, to=3-4]
      \arrow[""{name=6, anchor=center, inner sep=0}, "{N(p \odot m \odot q)}"'{inner sep=.8ex}, "\shortmid"{marking}, from=3-1, to=3-4]
      \arrow["{\alpha_p}"{description}, draw=none, from=0, to=3]
      \arrow["{\mu_m}"{description}, draw=none, from=1, to=4]
      \arrow["{\beta_q}"{description}, draw=none, from=2, to=5]
      \arrow["{N_{p,m,q}}"{description, pos=0.6}, draw=none, from=4, to=6]
    \end{tikzcd},
  \end{equation*}
  holding for every composable triple of proarrows $p: x' \proto x$,
  $m: x \proto y$, and $q: y \proto y'$ in $\dbl{D}$.
\end{proof}

\subsection{Twisted copresheaves}

Having relativized the notion of modification, we can now construct a 2-category
of twisted functors valued in a fixed double category. When the target double
category is $\Prof$, we will obtain a 2-category into which the 2-category of
double categories embeds.

\begin{construction}[2-category of twisted copresheaves]
  \label{construction:twisted-lax-2-cat}
  Let $\dbl{E}$ be a double category. There is a 2-category $\Twnll_{\dbl{E}}$
  having
  \begin{itemize}
    \item as objects, a double category $\dbl{A}$ together with a twisted normal
      doubly lax functor $F: \dbl{A} \twistto \dbl{E}$;
    \item as 1-morphisms
      $(F : \dbl{A} \twistto \dbl{E}) \to (G : \dbl{B} \twistto \dbl{E})$, a lax
      double functor $A : \dbl{A} \to \dbl{B}$ together with a natural
      transformation $\alpha : F \To GA$ (using
      \cref{construction:twisted-pre-composite}), for any choice of lax, pseudo,
      or strict natural transformations;
    \item as 2-morphisms $(A, \alpha) \To (B, \beta)$, a (tight) natural
      transformation $\tau : A \Rightarrow B$ together with a $\tau$-relative
      modification $\mu : \alpha \Tto \beta$.
  \end{itemize}
  Composition of 1-morphisms is defined by
  $(A, \alpha) \cdot (B, \beta) \coloneqq (A \cdot B, \alpha \cdot \beta A)$, i.e., by the pasting
  \begin{equation*}
    % https://q.uiver.app/#q=WzAsNCxbMCwwLCJcXGRibHtBfSJdLFsxLDAsIlxcZGJse0J9Il0sWzIsMCwiXFxkYmx7Q30iXSxbMSwxLCJcXGRibHtFfSJdLFswLDEsIkEiXSxbMSwyLCJCIl0sWzEsMywiRyIsMV0sWzAsMywiRiIsMl0sWzIsMywiSCJdLFs3LDYsIlxcYWxwaGEiLDAseyJzaG9ydGVuIjp7InNvdXJjZSI6MjAsInRhcmdldCI6MjB9fV0sWzYsOCwiXFxiZXRhIiwwLHsic2hvcnRlbiI6eyJzb3VyY2UiOjIwLCJ0YXJnZXQiOjIwfX1dXQ==
    \begin{tikzcd}[column sep=large]
      {\dbl{A}} & {\dbl{B}} & {\dbl{C}} \\
      & {\dbl{E}}
      \arrow["A", from=1-1, to=1-2]
      \arrow[""{name=0, anchor=center, inner sep=0}, "F"', from=1-1, to=2-2]
      \arrow["B", from=1-2, to=1-3]
      \arrow[""{name=1, anchor=center, inner sep=0}, "G"{description}, from=1-2, to=2-2]
      \arrow[""{name=2, anchor=center, inner sep=0}, "H", from=1-3, to=2-2]
      \arrow["\alpha", between={0.2}{0.8}, Rightarrow, from=0, to=1]
      \arrow["\beta", between={0.2}{0.8}, Rightarrow, from=1, to=2]
    \end{tikzcd}
  \end{equation*}
  where \cref{construction:twisted-pre-whiskering} is used to form the
  whiskering $\beta A$. The identity 1-morphism on $F: \dbl{A} \twistto \dbl{E}$ is
  $(1_{\dbl{A}}, 1_F)$, where $1_F$ is the identity transformation.

  Composition of 2-morphisms is defined as follows. The vertical composite of
  relative modifications $(\tau, \mu) : (A, \alpha) \To (B, \beta)$ and
  $(\sigma, \nu) : (B, \beta) \To (C, \gamma)$ is $(\tau, \mu) \cdot (\sigma, \nu) \coloneqq (\tau \cdot \sigma, \mu \cdot \nu)$,
  where $\mu \cdot \nu$ has components
  \begin{equation} \label{eq:vertical-composite-relative-modification}
    (\mu \cdot \nu)_x
    \quad\coloneqq\quad
    \begin{tikzcd}
      Fx && Fx \\
      Fx & Fx & Fx \\
      GAx & GBx & GCx \\
      GAx && GCx
      \arrow[""{name=0, anchor=center, inner sep=0}, "\shortmid"{marking}, equals, from=1-1, to=1-3]
      \arrow[equals, from=1-1, to=2-1]
      \arrow[equals, from=1-3, to=2-3]
      \arrow[""{name=1, anchor=center, inner sep=0}, "\shortmid"{marking}, equals, from=2-1, to=2-2]
      \arrow["{\alpha_x}"', from=2-1, to=3-1]
      \arrow[""{name=2, anchor=center, inner sep=0}, "\shortmid"{marking}, equals, from=2-2, to=2-3]
      \arrow["{\beta_x}"{description}, from=2-2, to=3-2]
      \arrow["{\gamma_x}", from=2-3, to=3-3]
      \arrow[""{name=3, anchor=center, inner sep=0}, "{G\tau_x}"'{inner sep=.8ex}, "\shortmid"{marking}, from=3-1, to=3-2]
      \arrow[equals, from=3-1, to=4-1]
      \arrow[""{name=4, anchor=center, inner sep=0}, "{G\sigma_x}"'{inner sep=.8ex}, "\shortmid"{marking}, from=3-2, to=3-3]
      \arrow[equals, from=3-3, to=4-3]
      \arrow[""{name=5, anchor=center, inner sep=0}, "{G(\tau_x \cdot \sigma_x)}"', from=4-1, to=4-3]
      \arrow["\cong"{description, pos=0.6}, draw=none, from=0, to=2-2]
      \arrow["{\mu_x}"{description}, draw=none, from=1, to=3]
      \arrow["{\nu_x}"{description}, draw=none, from=2, to=4]
      \arrow["{G_{\tau_x,\sigma_x}}"{description}, draw=none, from=3-2, to=5]
    \end{tikzcd},
    \qquad x \in \dbl{A}.
  \end{equation}
  It remains to define horizontal composition. The pre-whiskering of a relative
  modification $(\tau, \mu): (B, \beta) \To (C, \gamma): (\dbl{B}, G) \to (\dbl{C}, H)$ by a
  morphism $(A, \alpha): (\dbl{A}, F) \to (\dbl{B}, G)$ is the relative modification
  $(A * \tau, \alpha * \mu)$, where $\alpha * \mu$ has components
  \begin{equation} \label{eq:pre-whisker-relative-modification}
    (\alpha * \mu)_x
    \quad\coloneqq\quad
    % https://q.uiver.app/#q=WzAsNixbMCwwLCJGeCJdLFsxLDAsIkZ4Il0sWzAsMSwiR0F4Il0sWzEsMSwiR0F4Il0sWzAsMiwiSEJBeCJdLFsxLDIsIkhDQXgiXSxbMCwyLCJcXGFscGhhX3giLDJdLFsxLDMsIlxcYWxwaGFfeCJdLFswLDEsIiIsMCx7ImxldmVsIjoyLCJzdHlsZSI6eyJib2R5Ijp7Im5hbWUiOiJiYXJyZWQifSwiaGVhZCI6eyJuYW1lIjoibm9uZSJ9fX1dLFsyLDMsIiIsMix7ImxldmVsIjoyLCJzdHlsZSI6eyJib2R5Ijp7Im5hbWUiOiJiYXJyZWQifSwiaGVhZCI6eyJuYW1lIjoibm9uZSJ9fX1dLFsyLDQsIlxcYmV0YV97QXh9IiwyXSxbMyw1LCJcXGdhbW1hX3tBeH0iXSxbNCw1LCJIXFx0YXVfe0F4fSIsMix7InN0eWxlIjp7ImJvZHkiOnsibmFtZSI6ImJhcnJlZCJ9fX1dLFs4LDksIlxcaWRfe1xcYWxwaGFfeH0iLDEseyJzaG9ydGVuIjp7InNvdXJjZSI6MjAsInRhcmdldCI6MjB9LCJzdHlsZSI6eyJib2R5Ijp7Im5hbWUiOiJub25lIn0sImhlYWQiOnsibmFtZSI6Im5vbmUifX19XSxbOSwxMiwiXFxtdV97QXh9IiwxLHsic2hvcnRlbiI6eyJzb3VyY2UiOjIwLCJ0YXJnZXQiOjIwfSwic3R5bGUiOnsiYm9keSI6eyJuYW1lIjoibm9uZSJ9LCJoZWFkIjp7Im5hbWUiOiJub25lIn19fV1d
    \begin{tikzcd}
      Fx & Fx \\
      GAx & GAx \\
      HBAx & HCAx
      \arrow[""{name=0, anchor=center, inner sep=0}, "\shortmid"{marking}, equals, from=1-1, to=1-2]
      \arrow["{\alpha_x}"', from=1-1, to=2-1]
      \arrow["{\alpha_x}", from=1-2, to=2-2]
      \arrow[""{name=1, anchor=center, inner sep=0}, "\shortmid"{marking}, equals, from=2-1, to=2-2]
      \arrow["{\beta_{Ax}}"', from=2-1, to=3-1]
      \arrow["{\gamma_{Ax}}", from=2-2, to=3-2]
      \arrow[""{name=2, anchor=center, inner sep=0}, "{H\tau_{Ax}}"'{inner sep=.8ex}, "\shortmid"{marking}, from=3-1, to=3-2]
      \arrow["{\id_{\alpha_x}}"{description}, draw=none, from=0, to=1]
      \arrow["{\mu_{Ax}}"{description}, draw=none, from=1, to=2]
    \end{tikzcd},
    \qquad x \in \dbl{A}.
  \end{equation}
  The post-whiskering of $(\tau, \mu)$ by a morphism
  $(D, \delta): (\dbl{C}, H) \to (\dbl{D}, K)$ is the relative modification
  $(\tau * D, \mu * \delta)$, where $\mu * \delta$ has components
  \begin{equation} \label{eq:post-whisker-relative-modification}
    (\mu * \delta)_x
    \quad\coloneqq\quad
    % https://q.uiver.app/#q=WzAsNixbMCwwLCJHeCJdLFswLDEsIkhCeCJdLFsxLDAsIkd4Il0sWzEsMSwiSEN4Il0sWzAsMiwiS0RCeCJdLFsxLDIsIktEQ3giXSxbMCwxLCJcXGJldGFfeCIsMl0sWzIsMywiXFxnYW1tYV94Il0sWzAsMiwiIiwwLHsibGV2ZWwiOjIsInN0eWxlIjp7ImJvZHkiOnsibmFtZSI6ImJhcnJlZCJ9LCJoZWFkIjp7Im5hbWUiOiJub25lIn19fV0sWzEsMywiSFxcdGF1X3giLDIseyJzdHlsZSI6eyJib2R5Ijp7Im5hbWUiOiJiYXJyZWQifX19XSxbMSw0LCJcXGRlbHRhX3tCeH0iLDJdLFszLDUsIlxcZGVsdGFfe0N4fSJdLFs0LDUsIktEXFx0YXVfeCIsMix7InN0eWxlIjp7ImJvZHkiOnsibmFtZSI6ImJhcnJlZCJ9fX1dLFs4LDksIlxcbXVfeCIsMSx7InNob3J0ZW4iOnsic291cmNlIjoyMCwidGFyZ2V0IjoyMH0sInN0eWxlIjp7ImJvZHkiOnsibmFtZSI6Im5vbmUifSwiaGVhZCI6eyJuYW1lIjoibm9uZSJ9fX1dLFs5LDEyLCJcXGRlbHRhX3tcXHRhdV94fSIsMSx7ImxhYmVsX3Bvc2l0aW9uIjo2MCwic2hvcnRlbiI6eyJzb3VyY2UiOjIwLCJ0YXJnZXQiOjIwfSwic3R5bGUiOnsiYm9keSI6eyJuYW1lIjoibm9uZSJ9LCJoZWFkIjp7Im5hbWUiOiJub25lIn19fV1d
    \begin{tikzcd}
      Gx & Gx \\
      HBx & HCx \\
      KDBx & KDCx
      \arrow[""{name=0, anchor=center, inner sep=0}, "\shortmid"{marking}, equals, from=1-1, to=1-2]
      \arrow["{\beta_x}"', from=1-1, to=2-1]
      \arrow["{\gamma_x}", from=1-2, to=2-2]
      \arrow[""{name=1, anchor=center, inner sep=0}, "{H\tau_x}"'{inner sep=.8ex}, "\shortmid"{marking}, from=2-1, to=2-2]
      \arrow["{\delta_{Bx}}"', from=2-1, to=3-1]
      \arrow["{\delta_{Cx}}", from=2-2, to=3-2]
      \arrow[""{name=2, anchor=center, inner sep=0}, "{KD\tau_x}"'{inner sep=.8ex}, "\shortmid"{marking}, from=3-1, to=3-2]
      \arrow["{\mu_x}"{description}, draw=none, from=0, to=1]
      \arrow["{\delta_{\tau_x}}"{description, pos=0.6}, draw=none, from=1, to=2]
    \end{tikzcd},
    \qquad x \in \dbl{B}.
  \end{equation}
  This suffices to determine the horizontal composition of the 2-category (but
  see the proof of \cref{lem:twisted-lax-2-cat} below for the general formula).

  We write $\Twnll$ for the important special case $\Twnll_{\Prof}$ in which the
  twisted functors are valued in profunctors.
\end{construction}

\begin{lemma} \label{lem:twisted-lax-2-cat}
  For any double category $\dbl{E}$, the 2-category $\Twnll_{\dbl{E}}$ in
  \cref{construction:twisted-lax-2-cat} is well-defined.
\end{lemma}
\begin{proof}
  The first component of the claimed 2-category is none other than the
  2-category $\Dbll$ of double categories, lax double functors, and natural
  transformations. Moreover, the underlying 1-category of the claimed 2-category
  is the Grothendieck construction of the indexed category
  $|\Dbll|_1^{\op} \to \Cat$ from \cref{lemma:pre-composition-functor}, so is, in
  particular, a well-defined category.

  We must check the remaining data and axioms concerning 2-morphisms, for which
  it is convenient to use the description of a 2-category as a sesquicategory
  satisfying an interchange law \mbox{\cite[\S{2.3.2}]{grandis2019}}. First, we
  must see that the whiskerings are well-defined relative modifications.
  Tight-to-loose equivariance of the pre- and post-whiskered relative
  modifications is proved using loose unitality in the codomain double category
  and functorality of the component cells of the \mbox{transformation $\delta$},
  respectively, along with the tight-to-loose equivariance of the original
  relative modification. Loose-to-tight equivariance is also straightforward to
  check.

  As for the 2-category axioms, associativity and unitality of vertical
  composition follow from the associativity and unitality of the tight-to-loose
  comparisons of the codomain twisted functor. Associativity and unitality of
  pre- and post-whiskering follow directly from how composition and identities
  are defined for natural transformations (\cref{lem:cat-twisted-lax}).
  Distributivity of pre- and post-whiskering is again, respectively, loose
  unitality in the codomain double category and functorality of the component
  cells of a transformation between twisted functors.

  Finally, to prove the reduced interchange axiom, consider a pair of relative
  modifications that are horizontally composable:
  \begin{equation*}
    % https://q.uiver.app/#q=WzAsMyxbMCwwLCIoXFxkYmx7QX0sRikiXSxbMSwwLCIoXFxkYmx7Qn0sRykiXSxbMiwwLCIoXFxkYmx7Q30sSCkiXSxbMCwxLCIoQSxcXGFscGhhKSIsMCx7ImN1cnZlIjotM31dLFswLDEsIihCLFxcYmV0YSkiLDIseyJjdXJ2ZSI6M31dLFsxLDIsIihDLFxcZ2FtbWEpIiwwLHsiY3VydmUiOi0zfV0sWzEsMiwiKEQsXFxkZWx0YSkiLDIseyJjdXJ2ZSI6M31dLFszLDQsIihcXHRhdSxcXG11KSIsMSx7InNob3J0ZW4iOnsic291cmNlIjoyMCwidGFyZ2V0IjoyMH19XSxbNSw2LCIoXFxzaWdtYSxcXG51KSIsMSx7InNob3J0ZW4iOnsic291cmNlIjoyMCwidGFyZ2V0IjoyMH19XV0=
    \begin{tikzcd}
      {(\dbl{A},F)} & {(\dbl{B},G)} & {(\dbl{C},H)}
      \arrow[""{name=0, anchor=center, inner sep=0}, "{(A,\alpha)}", curve={height=-18pt}, from=1-1, to=1-2]
      \arrow[""{name=1, anchor=center, inner sep=0}, "{(B,\beta)}"', curve={height=18pt}, from=1-1, to=1-2]
      \arrow[""{name=2, anchor=center, inner sep=0}, "{(C,\gamma)}", curve={height=-18pt}, from=1-2, to=1-3]
      \arrow[""{name=3, anchor=center, inner sep=0}, "{(D,\delta)}"', curve={height=18pt}, from=1-2, to=1-3]
      \arrow["{(\tau,\mu)}"{description}, between={0.2}{0.8}, Rightarrow, from=0, to=1]
      \arrow["{(\sigma,\nu)}"{description}, between={0.2}{0.8}, Rightarrow, from=2, to=3]
    \end{tikzcd}.
  \end{equation*}
  We must show that, for each object $x \in \dbl{A}$, the equation holds:
  \begin{equation*}
    \begin{tikzcd}
      Fx && Fx \\
      Fx & Fx & Fx \\
      GAx & GBx & GBx \\
      HCAx & HCBx & HDBx \\
      HCAx && HDBx
      \arrow[""{name=0, anchor=center, inner sep=0}, "\shortmid"{marking}, equals, from=1-1, to=1-3]
      \arrow[equals, from=1-1, to=2-1]
      \arrow[equals, from=1-3, to=2-3]
      \arrow[""{name=1, anchor=center, inner sep=0}, "\shortmid"{marking}, equals, from=2-1, to=2-2]
      \arrow["{\alpha_x}"', from=2-1, to=3-1]
      \arrow[""{name=2, anchor=center, inner sep=0}, "\shortmid"{marking}, equals, from=2-2, to=2-3]
      \arrow["{\beta_x}"{description}, from=2-2, to=3-2]
      \arrow["{\beta_x}", from=2-3, to=3-3]
      \arrow[""{name=3, anchor=center, inner sep=0}, "{G\tau_x}"{inner sep=.8ex}, "\shortmid"{marking}, from=3-1, to=3-2]
      \arrow["{\gamma_{Ax}}"', from=3-1, to=4-1]
      \arrow[""{name=4, anchor=center, inner sep=0}, "\shortmid"{marking}, equals, from=3-2, to=3-3]
      \arrow["{\gamma_{Bx}}"{description}, from=3-2, to=4-2]
      \arrow["{\delta_{Bx}}", from=3-3, to=4-3]
      \arrow[""{name=5, anchor=center, inner sep=0}, "{HC\tau_x}"'{inner sep=.8ex}, "\shortmid"{marking}, from=4-1, to=4-2]
      \arrow[equals, from=4-1, to=5-1]
      \arrow[""{name=6, anchor=center, inner sep=0}, "{H\sigma_{Bx}}"'{inner sep=.8ex}, "\shortmid"{marking}, from=4-2, to=4-3]
      \arrow[equals, from=4-3, to=5-3]
      \arrow[""{name=7, anchor=center, inner sep=0}, "{H(C\tau_x \cdot \sigma_{Bx})}"'{inner sep=.8ex}, "\shortmid"{marking}, from=5-1, to=5-3]
      \arrow["\cong"{description, pos=0.6}, draw=none, from=0, to=2-2]
      \arrow["{\mu_x}"{description, pos=0.4}, draw=none, from=1, to=3]
      \arrow["{\id_{\beta_x}}"{description}, draw=none, from=2, to=4]
      \arrow["{\gamma_{\tau_x}}"{description}, draw=none, from=3, to=5]
      \arrow["{\nu_{Bx}}"{description}, draw=none, from=4, to=6]
      \arrow["{H_{C\tau_x, \sigma_{Bx}}}"{description}, draw=none, from=4-2, to=7]
    \end{tikzcd}
    \quad=\quad
    \begin{tikzcd}
      Fx && Fx \\
      Fx & Fx & Fx \\
      GAx & GAx & GBx \\
      HCAx & HDAx & HDBx \\
      HCAx && HDBx
      \arrow[""{name=0, anchor=center, inner sep=0}, "\shortmid"{marking}, equals, from=1-1, to=1-3]
      \arrow[equals, from=1-1, to=2-1]
      \arrow[equals, from=1-3, to=2-3]
      \arrow[""{name=1, anchor=center, inner sep=0}, "\shortmid"{marking}, equals, from=2-1, to=2-2]
      \arrow["{\alpha_x}"', from=2-1, to=3-1]
      \arrow[""{name=2, anchor=center, inner sep=0}, "\shortmid"{marking}, equals, from=2-2, to=2-3]
      \arrow["{\alpha_x}"{description}, from=2-2, to=3-2]
      \arrow["{\beta_x}", from=2-3, to=3-3]
      \arrow[""{name=3, anchor=center, inner sep=0}, "\shortmid"{marking}, equals, from=3-1, to=3-2]
      \arrow["{\gamma_{Ax}}"', from=3-1, to=4-1]
      \arrow[""{name=4, anchor=center, inner sep=0}, "{G\tau_x}"{inner sep=.8ex}, "\shortmid"{marking}, from=3-2, to=3-3]
      \arrow["{\delta_{Ax}}"{description}, from=3-2, to=4-2]
      \arrow["{\delta_{Bx}}", from=3-3, to=4-3]
      \arrow[""{name=5, anchor=center, inner sep=0}, "{H \sigma_{Ax}}"'{inner sep=.8ex}, "\shortmid"{marking}, from=4-1, to=4-2]
      \arrow[equals, from=4-1, to=5-1]
      \arrow[""{name=6, anchor=center, inner sep=0}, "{HD\tau_x}"'{inner sep=.8ex}, "\shortmid"{marking}, from=4-2, to=4-3]
      \arrow[equals, from=4-3, to=5-3]
      \arrow[""{name=7, anchor=center, inner sep=0}, "{H(\sigma_{Ax} \cdot D\tau_x)}"'{inner sep=.8ex}, "\shortmid"{marking}, from=5-1, to=5-3]
      \arrow["\cong"{description, pos=0.6}, draw=none, from=0, to=2-2]
      \arrow["{\id_{\alpha_x}}"{description}, draw=none, from=1, to=3]
      \arrow["{\mu_x}"{description, pos=0.4}, draw=none, from=2, to=4]
      \arrow["{\nu_{Ax}}"{description}, draw=none, from=3, to=5]
      \arrow["{\delta_{\tau_x}}"{description}, draw=none, from=4, to=6]
      \arrow["{H_{\sigma_{Ax}, D\tau_x}}"{description}, draw=none, from=4-2, to=7]
    \end{tikzcd}.
  \end{equation*}
  Note that the equation is well-typed since
  $C\tau_x \cdot \sigma_{Bx} = (\tau * \sigma)_x = \sigma_{Ax} \cdot D\tau_x$ is the component of the horizontal
  composite $\tau * \sigma$ at $x$. The equation itself is proved by applying
  $\id_{\alpha_x} \odot \mu_x \cong \mu_x \odot \id_{\beta_x}$ and the tight-to-loose equivariance of the
  relative modification $(\sigma,\nu)$ at the arrow $\tau_x: Ax \to Bx$. Conclude that
  horizontal composition in the 2-category is well-defined by
  $(\tau, \mu) * (\sigma, \nu) \coloneqq (\tau * \sigma, \mu * \nu)$, where $\mu * \nu$ has component at $x$
  given by either side of the above equation.
\end{proof}

\begin{proposition}[Embedding into twisted copresheaves]
\label{prop:embedding-into-twisted-copresheaves}
  There is a 2-functor
  \begin{equation*}
    \Dbll \to \Twnll,
  \end{equation*}
  that sends
  \begin{itemize}[noitemsep]
    \item a double category $\dbl{D}$ to its twisted Hom $\dbl{D}(-, =)$
      (\cref{construction:twistedhomfunctor})
    \item a lax double functor $F$ to the lax natural transformation
      $(F^\co \times F, \tilde F)$
      (\cref{construction:double-functor-to-twisted-transformation})
    \item a natural transformation $\alpha$ to the relative modification
      $(\alpha^\co \times \alpha, \tilde \alpha)$
      (\cref{construction:transformation-to-relative-modification}).
  \end{itemize}
\end{proposition}
\begin{proof}[Proof sketch]
  As indicated in the statement, we have already constructed the assignment of
  this 2-functor; what must be proved is its 2-functorality. For 1-functorality,
  notice that, given composable lax double functors $F$ and $G$, the components
  \eqref{eq:composite-natural-transformation-components} of the composite
  $((F \cdot G)^\co \times (F \cdot G), \tilde F \cdot {\tilde G} F)$ just apply $F$ and then $G$
  to proarrows and cells in the domain, while the naturality comparisons
  \eqref{eq:composite-natural-transformation-comparisons} of the composite
  transformation are given by the composition comparisons of the composite
  functor $F \cdot G$ (cf.\ \cite[Equation 3.63]{grandis2019}). Next, given
  vertically composable natural transformations $\alpha$ and $\beta$ between double
  functors, the composite relative modification
  $((\alpha \cdot \beta)^\co \times (\alpha \cdot \beta), \tilde\alpha \cdot \tilde\beta)$, defined by
  \cref{eq:vertical-composite-relative-modification}, acts on globular cells by
  composing with the cell components of $\alpha$ and then $\beta$. Finally,
  pre-whiskering \eqref{eq:pre-whisker-relative-modification} and
  post-whiskering \eqref{eq:post-whisker-relative-modification} of the relative
  modification $(\alpha^\co \times \alpha, \tilde \alpha)$ agree with the componentwise pre- and
  post-whiskering of a natural transformation $\alpha$.
\end{proof}

\begin{remark}[Meaning of embedding]
  The 2-functor $\Dbll \to \Twnll$ just constructed is the two-sided version of
  a twisted Yoneda embedding, which raises the question of in what sense it
  \emph{is} an embedding. That this 2-functor is \emph{locally} an embedding,
  i.e., that its hom-functors are injective on objects and faithful, is
  immediate, since 1- and 2-morphisms in the domain 2-category $\Dbll$ appear
  directly in the first component of the defining assignment. Inasmuch as the
  twisted Hom $\dbl{D}(-,=)$ just repackages the data of the double category
  $\dbl{D}$, the 2-functor should be essentially injective on objects, but
  injectivity in the strictest sense apparently fails due to the two-sided way
  that the loose-to-tight comparisons assemble the associators and unitors of
  $\dbl{D}$. A careful analysis is too far a digression to undertake here.
\end{remark}

We have seen that the input to construct a transformation between twisted Homs
generalizes from double functors to modules, and the input to construct a
relative modification generalizes from natural transformations to modulations.
The embedding of double categories into twisted copresheaves ought to generalize
accordingly, but to make sense of that, we must know how to compose modules.
Fortunately, the relevant notion of composition of modules is not that along a
shared double functor, which is fraught with difficulties \cite{pare2013}. The
notion appropriate here is much simpler, yet to our knowledge does not appear in
the literature.

\begin{construction}[2-category of double categories, modules, and modulations]
\label{constr:modulations}
  Denote by $\DblMod$ the co-Kleisli 2-category of the product 2-comonad
  $(-) \times \Loose$ on the 2-category $\Dbll$, where $\Loose$ is the walking
  loose morphism.
\end{construction}

So the 2-category $\DblMod$ has:
\begin{itemize}[noitemsep]
  \item as objects, double categories;
  \item as 1-morphisms from $\dbl{D}$ to $\dbl{E}$, lax double functors
    $M: \dbl{D} \times \Loose \to \dbl{E}$, i.e., a pair of lax double functors
    $F,G: \dbl{D} \to \dbl{E}$ together with a module $M: F \proTo G$;
  \item as 2-morphisms from $M$ to $N$, natural transformations $\mu: M \To N$,
    i.e., modulations between modules.
\end{itemize}
Modules $M: \dbl{C} \times \Loose \to \dbl{D}$ and
$N: \dbl{D} \times \Loose \to \dbl{E}$ compose via the formula
\begin{equation*}
  \dbl{C} \times \Loose \xto{1_{\dbl{C}} \times \Delta_{\Loose}}
  \dbl{C} \times \Loose^2 \xto{M \times 1_{\Loose}}
  \dbl{D} \times \Loose \xto{N}
  \dbl{E}.
\end{equation*}

\begin{proposition}[Embedding into twisted copresheaves, generalized]
  There is a 2-functor
  \begin{equation*}
    \DblMod \to \Twnll
  \end{equation*}
  that sends
  \begin{itemize}[noitemsep]
    \item a double category $\dbl{D}$ to its twisted Hom $\dbl{D}(-, =)$
      (\cref{construction:twistedhomfunctor})
      \item a pair of lax double functors $F,G$ together with a module $M: F \proTo G$ to
      the lax natural transformation $(F^\co \times G, \tilde M)$
      (\cref{construction:module-to-twisted-transformation})
    \item a pair of natural transformations $\alpha,\beta$ together with a
      modulation $\mu: \alpha \proTo \beta$ to the relative modification
      $(\alpha^\co \times \beta, \tilde \mu)$
      (\cref{construction:modulation-to-relative-modification}).
  \end{itemize}
\end{proposition}

The proof of this proposition is similar to the previous one; it also follows as a special case of the upcoming \cref{prop:double-barrel-to-twisted-copresheaf}.

\section{Double barrels and collages}
\label{sec:double-barrel}

In this section, we will compare twisted bimodules valued in $\Prof$ with \emph{double barrels}: double functors into the walking proarrow. Recall that, in \cref{constr:modulations}, we defined  the 2-category $\DblMod$ to be the co-Kleisli 2-category of the product 2-comonad $(-) \times \Loose$ on $\Dbll$. The corresponding 2-category of coalgebras is the slice $\Dbll \downarrow \Loose$ over the walking proarrow $\Loose$. Lax double functors $M : \dbl{D} \to \Loose$ (which must, in fact, be strict functors because $\Loose$ has no non-trivial squares) were called \emph{double barrels} in \cite{brown2025}, since they are the double-categorical analogue of Joyal's ``barrel'' presentation of profunctors as functors into the walking arrow. To restate:

\begin{definition}[Double barrel]
 A \define{double barrel} is a double functor $M : \dbl{M} \to \Loose$. The 2-category $\lModl$ of double barrels, lax morphisms, and tight transformations is defined to be the \emph{strict} 2-categorical slice.

 The \define{source} of a double barrel $M : \dbl{M} \to \Loose$ is the inverse image $\barrelsrc{\dbl{M}} \coloneqq M^{{-1}}(0)$, and the \define{target} of a double barrel is the inverse image $\barreltgt{\dbl{M}} \coloneqq M^{{-1}}(1)$. We denote the inclusions $\barrelsrc{\dbl{M}} \hookrightarrow \dbl{M}$ and $\barreltgt{\dbl{M}} \hookrightarrow \dbl{M}$ by $\inl$ and $\inr$ respectively. We write $M : \barrelsrc{\dbl{M}} \lbimodto \barreltgt{\dbl{M}}$ to say that $M$ is a double barrel from its source to its target. Taking the source and target of a double barrel gives a 2-functor $\lModl \to \Dbll \times \Dbll$.
\end{definition}

\subsection{From double barrels to twisted bimodules}

To go from a double barrel to a twisted bimodule is easy: we can restrict its twisted Hom functor.

\begin{construction}[Twisted bimodule of sections]
\label{construction:twisted-bimod-sections}
Given a double barrel $M : \dbl{M} \to \Loose$ from $\dbl{D}$ to $\dbl{E}$, we define its \define{twisted bimodule of sections} to be the restriction
\begin{equation*}
 \twsection{M} \coloneqq \dbl{M}(\inl(-), \inr (=)) : \dbl{D}^{\co} \times \dbl{E} \twistto \Prof
\end{equation*}
of the twisted Hom functor of $\dbl{M}$ (\cref{construction:twistedhomfunctor}).
\end{construction}

This construction extends functorially in much the same way that \cref{construction:module-to-twisted-transformation} was defined.

\begin{construction}[Lax transformation from morphism of double barrels]
\label{construction:loose-bimod-morphism}
Let $F : M \to N$ be a lax morphism of double barrels; that is, a commuting triangle
\[
% https://q.uiver.app/#q=WzAsMyxbMCwwLCJcXGRibHtNfSJdLFsyLDAsIlxcZGJse059Il0sWzEsMSwiXFxMb29zZSJdLFswLDIsIk0iLDJdLFsxLDIsIk4iXSxbMCwxLCJGIl1d
\begin{tikzcd}
	{\dbl{M}} && {\dbl{N}} \\
	& \Loose
	\arrow["F", from=1-1, to=1-3]
	\arrow["M"', from=1-1, to=2-2]
	\arrow["N", from=1-3, to=2-2]
\end{tikzcd}.
\]
We define its associated lax natural transformation
\begin{equation*}
  \twsection{F} : \twsection{M}(-,=) \To \twsection{N}(\barrelsrc{F}(-), \barreltgt{F}(=)) : \barrelsrc{\dbl{M}}^\co \times \barreltgt{\dbl{M}} \twistto \Prof
\end{equation*}
to be the transformation $\tilde F : \dbl{M}(-, =) \To \dbl{N}(F(-), F(=))$ of \cref{construction:double-functor-to-twisted-transformation}, restricted on the left by $\inl$ and on the right by $\inr$ in both the domain and codomain.
Explicitly:
\begin{itemize}
        \item The component at a pair of objects $x$ in $\barrelsrc{\dbl{M}}$ and $y$ in $\barreltgt{\dbl{M}}$ is the functor
        \begin{equation*}
          \tilde F_{x, y} : \dbl{M}(\inl(x), \inr(y)) \to \dbl{N}(\inl(\barrelsrc{F}(x)), \inr(\barreltgt{F}(y)))
        \end{equation*}
        sending a proarrow $m : x \proto y$ in $\dbl{M}$ to $Fm : Fx \proto Fy$ in $\dbl{N}$, and similarly for globular squares. 
        \item The component at a pair of arrows $f : x \to w$ in $\barrelsrc{\dbl{M}}$ and $g : y \to z$ in $\barreltgt{\dbl{M}}$ is the transformation
        % https://q.uiver.app/#q=WzAsNCxbMCwwLCJcXGRibHtNfSh4LHkpIl0sWzEsMCwiXFxkYmx7TX0odyx6KSJdLFswLDEsIlxcZGJse059KEZ4LEZ5KSJdLFsxLDEsIlxcZGJseyBOfShGdyxGeikiXSxbMCwxLCJcXGRibHtNfShmLGcpIiwwLHsic3R5bGUiOnsiYm9keSI6eyJuYW1lIjoiYmFycmVkIn19fV0sWzAsMiwiRl97eCx5fSIsMl0sWzEsMywiRl97dyx6fSJdLFsyLDMsIlxcZGJse059KEZmLEZnKSIsMix7InN0eWxlIjp7ImJvZHkiOnsibmFtZSI6ImJhcnJlZCJ9fX1dLFs0LDcsIkZfe2YsZ30iLDEseyJzaG9ydGVuIjp7InNvdXJjZSI6MjAsInRhcmdldCI6MjB9LCJzdHlsZSI6eyJib2R5Ijp7Im5hbWUiOiJub25lIn0sImhlYWQiOnsibmFtZSI6Im5vbmUifX19XV0=
\begin{tikzcd}
	{\dbl{M}(x,y)} & {\dbl{M}(w,z)} \\
	{\dbl{N}(Fx,Fy)} & {\dbl{ N}(Fw,Fz)}
	\arrow[""{name=0, anchor=center, inner sep=0}, "{\dbl{M}(f,g)}"{inner sep=.8ex}, "\shortmid"{marking}, from=1-1, to=1-2]
	\arrow["{\tilde F_{x,y}}"', from=1-1, to=2-1]
	\arrow["{\tilde F_{w,z}}", from=1-2, to=2-2]
	\arrow[""{name=1, anchor=center, inner sep=0}, "{\dbl{N}(Ff,Fg)}"'{inner sep=.8ex}, "\shortmid"{marking}, from=2-1, to=2-2]
	\arrow["{\tilde F_{f,g}}"{description}, draw=none, from=0, to=1]
\end{tikzcd}
  that sends a cell $\stdInlineCell{\alpha}$ in $\dbl{M}$ to its image $F\alpha$.
        \item The naturality comparison at a pair of proarrows $p: x' \proto x$ in $\barrelsrc{\dbl{M}}$ and $q : y \proto y'$ in $\barreltgt{\dbl{M}}$ is the transformation
        % https://q.uiver.app/#q=WzAsNixbMCwwLCJcXGRibHtNfSh4LHkpIl0sWzAsMSwiXFxkYmx7Tn0oRngsRnkpIl0sWzAsMiwiXFxkYmx7Tn0oRngnLEZ5JykiXSxbMSwyLCJcXGRibHtOfShGeCcsR3knKSJdLFsxLDAsIlxcZGJse019KHgseSkiXSxbMSwxLCJcXGRibHtNfSh4Jyx5JykiXSxbMCwxLCJGX3t4LHl9IiwyXSxbMSwyLCJcXGRibHtOfShGcCxGcSkiLDJdLFsyLDMsIiIsMix7ImxldmVsIjoyLCJzdHlsZSI6eyJib2R5Ijp7Im5hbWUiOiJiYXJyZWQifSwiaGVhZCI6eyJuYW1lIjoibm9uZSJ9fX1dLFs0LDUsIlxcZGJse019KHAscSkiXSxbNSwzLCJGX3t4Jyx5J30iXSxbMCw0LCIiLDAseyJsZXZlbCI6Miwic3R5bGUiOnsiYm9keSI6eyJuYW1lIjoiYmFycmVkIn0sImhlYWQiOnsibmFtZSI6Im5vbmUifX19XSxbMTEsOCwiRl97cCxxfSIsMSx7InNob3J0ZW4iOnsic291cmNlIjoyMCwidGFyZ2V0IjoyMH0sInN0eWxlIjp7ImJvZHkiOnsibmFtZSI6Im5vbmUifSwiaGVhZCI6eyJuYW1lIjoibm9uZSJ9fX1dXQ==
\begin{tikzcd}
	{\dbl{M}(x,y)} & {\dbl{M}(x,y)} \\
	{\dbl{N}(Fx,Fy)} & {\dbl{M}(x',y')} \\
	{\dbl{N}(Fx',Fy')} & {\dbl{N}(Fx',Gy')}
	\arrow[""{name=0, anchor=center, inner sep=0}, "\shortmid"{marking}, equals, from=1-1, to=1-2]
	\arrow["{\tilde F_{x,y}}"', from=1-1, to=2-1]
	\arrow["{\dbl{M}(p,q)}", from=1-2, to=2-2]
	\arrow["{\dbl{N}(Fp,Fq)}"', from=2-1, to=3-1]
	\arrow["{\tilde F_{x',y'}}", from=2-2, to=3-2]
	\arrow[""{name=1, anchor=center, inner sep=0}, "\shortmid"{marking}, equals, from=3-1, to=3-2]
	\arrow["{\tilde F_{p,q}}"{description}, draw=none, from=0, to=1]
\end{tikzcd}
whose component at a proarrow $m : x \proto y$ in $\dbl{M}$ is the tightly globular 
cell in $\dbl{N}$
      \begin{equation*}
        (\tilde F_{p,q})_m \coloneqq F_{p,m,q}: Fp \odot Fm \odot Fq \to F(p \odot m \odot q)
      \end{equation*}
      given by the laxators of $F$ (in either
      order, equivalent by the associativity axiom for laxators).
\end{itemize}
\end{construction}

Note that \cref{construction:loose-bimod-morphism} does indeed generalize \cref{construction:module-to-twisted-transformation} when applied to a lax double functor $F : \dbl{D} \times \Loose \to \dbl{E} \times \Loose$ in the slice over $\Loose$, recalling that such a functor is determined by its first component, a module $\pi_1 \circ F : \dbl{D} \times \Loose \to \dbl{E}$. (The second component of $F$ must be the product projection for it to live in the strict slice over $\Loose$.)

Finally, we may generalize \cref{construction:modulation-to-relative-modification} to double barrels as follows.
\begin{construction}[Relative modification from 2-cell of double barrels]
\label{construction:relative-modification-from-barrel}
Let $\alpha : F \To G : \dbl{M} \to \dbl{N}$ be a 2-cell in $\lModl$ between lax functors $F$ and $G$ between double barrels $M : \dbl{M} \to \Loose$ and $N : \dbl{N} \to \Loose$. Then $\alpha$ induces a relative modification
\begin{equation*}
\twsection{\alpha} : \twsection{F} \Tto \twsection{G} : \twsection{M}(-, =) \To \twsection{N}(-,=),
\end{equation*}
given by restricting $(\alpha^\co \times \alpha, \tilde \alpha)$ of \cref{construction:transformation-to-relative-modification} along $\inl$ and $\inr$ respectively.

Explicitly, its component at a pair of objects $x$ in $\barrelsrc{\dbl{M}}$ and $y$ in $\barreltgt{\dbl{M}}$ is the natural transformation
  \begin{equation*}
    % https://q.uiver.app/#q=WzAsNCxbMCwwLCJcXGRibHtEfSh4LHkpIl0sWzEsMCwiXFxkYmx7RH0oeCx5KSJdLFswLDEsIlxcZGJse0V9KEZ4LEZ5KSJdLFsxLDEsIlxcZGJse0V9KEd4LEd5KSJdLFswLDEsIiIsMCx7ImxldmVsIjoyLCJzdHlsZSI6eyJib2R5Ijp7Im5hbWUiOiJiYXJyZWQifSwiaGVhZCI6eyJuYW1lIjoibm9uZSJ9fX1dLFswLDIsIlxcdGlsZGUgRl97eCx5fSIsMl0sWzEsMywiXFx0aWxkZSBHX3t4LHl9Il0sWzIsMywiXFxkYmx7RX0oXFxhbHBoYV94LFxcYWxwaGFfeSkiLDIseyJzdHlsZSI6eyJib2R5Ijp7Im5hbWUiOiJiYXJyZWQifX19XSxbNCw3LCJcXHRpbGRlIFxcYWxwaGFfe3gseX0iLDEseyJzaG9ydGVuIjp7InNvdXJjZSI6MjAsInRhcmdldCI6MjB9LCJzdHlsZSI6eyJib2R5Ijp7Im5hbWUiOiJub25lIn0sImhlYWQiOnsibmFtZSI6Im5vbmUifX19XV0=
    \begin{tikzcd}
      {\dbl{M}(x,y)} & {\dbl{M}(x,y)} \\
      {\dbl{N}(Fx,Fy)} & {\dbl{N}(Gx,Gy)}
      \arrow[""{name=0, anchor=center, inner sep=0}, "\shortmid"{marking}, equals, from=1-1, to=1-2]
      \arrow["{\tilde F_{x,y}}"', from=1-1, to=2-1]
      \arrow["{\tilde G_{x,y}}", from=1-2, to=2-2]
      \arrow[""{name=1, anchor=center, inner sep=0}, "{\dbl{E}(\alpha_x,\alpha_y)}"'{inner sep=.8ex}, "\shortmid"{marking}, from=2-1, to=2-2]
      \arrow["{\tilde \alpha_{x,y}}"{description}, draw=none, from=0, to=1]
    \end{tikzcd}
  \end{equation*}
  that acts on tightly globular cells in $\dbl{M}(x,y)$ as
  \begin{equation*}
    % https://q.uiver.app/#q=WzAsNCxbMCwwLCJ4Il0sWzEsMCwieSJdLFswLDEsIngiXSxbMSwxLCJ5Il0sWzAsMSwibSIsMCx7InN0eWxlIjp7ImJvZHkiOnsibmFtZSI6ImJhcnJlZCJ9fX1dLFsxLDMsIiIsMCx7ImxldmVsIjoyLCJzdHlsZSI6eyJoZWFkIjp7Im5hbWUiOiJub25lIn19fV0sWzAsMiwiIiwyLHsibGV2ZWwiOjIsInN0eWxlIjp7ImhlYWQiOnsibmFtZSI6Im5vbmUifX19XSxbMiwzLCJuIiwyLHsic3R5bGUiOnsiYm9keSI6eyJuYW1lIjoiYmFycmVkIn19fV0sWzQsNywiXFxnYW1tYSIsMSx7InNob3J0ZW4iOnsic291cmNlIjoyMCwidGFyZ2V0IjoyMH0sInN0eWxlIjp7ImJvZHkiOnsibmFtZSI6Im5vbmUifSwiaGVhZCI6eyJuYW1lIjoibm9uZSJ9fX1dXQ==
    \begin{tikzcd}
      x & y \\
      x & y
      \arrow[""{name=0, anchor=center, inner sep=0}, "m"{inner sep=.8ex}, "\shortmid"{marking}, from=1-1, to=1-2]
      \arrow[equals, from=1-1, to=2-1]
      \arrow[equals, from=1-2, to=2-2]
      \arrow[""{name=1, anchor=center, inner sep=0}, "n"'{inner sep=.8ex}, "\shortmid"{marking}, from=2-1, to=2-2]
      \arrow["\gamma"{description}, draw=none, from=0, to=1]
    \end{tikzcd}
    \qquad\mapsto\qquad
    % https://q.uiver.app/#q=WzAsNixbMCwwLCJGeCJdLFsxLDAsIkZ5Il0sWzAsMSwiRngiXSxbMSwxLCJGeSJdLFswLDIsIkd4Il0sWzEsMiwiR3kiXSxbMCwxLCJGbSIsMCx7InN0eWxlIjp7ImJvZHkiOnsibmFtZSI6ImJhcnJlZCJ9fX1dLFsxLDMsIiIsMCx7ImxldmVsIjoyLCJzdHlsZSI6eyJoZWFkIjp7Im5hbWUiOiJub25lIn19fV0sWzAsMiwiIiwyLHsibGV2ZWwiOjIsInN0eWxlIjp7ImhlYWQiOnsibmFtZSI6Im5vbmUifX19XSxbMiwzLCJGbiIsMix7InN0eWxlIjp7ImJvZHkiOnsibmFtZSI6ImJhcnJlZCJ9fX1dLFs0LDUsIkduIiwyLHsic3R5bGUiOnsiYm9keSI6eyJuYW1lIjoiYmFycmVkIn19fV0sWzIsNCwiXFxhbHBoYV94IiwyXSxbMyw1LCJcXGFscGhhX3kiXSxbNiw5LCJGXFxnYW1tYSIsMSx7InNob3J0ZW4iOnsic291cmNlIjoyMCwidGFyZ2V0IjoyMH0sInN0eWxlIjp7ImJvZHkiOnsibmFtZSI6Im5vbmUifSwiaGVhZCI6eyJuYW1lIjoibm9uZSJ9fX1dLFs5LDEwLCJcXGFscGhhX24iLDEseyJsYWJlbF9wb3NpdGlvbiI6NjAsInNob3J0ZW4iOnsic291cmNlIjoyMCwidGFyZ2V0IjoyMH0sInN0eWxlIjp7ImJvZHkiOnsibmFtZSI6Im5vbmUifSwiaGVhZCI6eyJuYW1lIjoibm9uZSJ9fX1dXQ==
    \begin{tikzcd}
      Fx & Fy \\
      Fx & Fy \\
      Gx & Gy
      \arrow[""{name=0, anchor=center, inner sep=0}, "Fm"{inner sep=.8ex}, "\shortmid"{marking}, from=1-1, to=1-2]
      \arrow[equals, from=1-1, to=2-1]
      \arrow[equals, from=1-2, to=2-2]
      \arrow[""{name=1, anchor=center, inner sep=0}, "Fn"'{inner sep=.8ex}, "\shortmid"{marking}, from=2-1, to=2-2]
      \arrow["{\alpha_x}"', from=2-1, to=3-1]
      \arrow["{\alpha_y}", from=2-2, to=3-2]
      \arrow[""{name=2, anchor=center, inner sep=0}, "Gn"'{inner sep=.8ex}, "\shortmid"{marking}, from=3-1, to=3-2]
      \arrow["{F\gamma}"{description}, draw=none, from=0, to=1]
      \arrow["{\alpha_n}"{description, pos=0.6}, draw=none, from=1, to=2]
    \end{tikzcd}
    =
    % https://q.uiver.app/#q=WzAsNixbMCwxLCJHeCJdLFsxLDEsIkd5Il0sWzAsMiwiR3giXSxbMSwyLCJHeSJdLFswLDAsIkZ4Il0sWzEsMCwiRnkiXSxbMCwxLCJHbSIsMCx7InN0eWxlIjp7ImJvZHkiOnsibmFtZSI6ImJhcnJlZCJ9fX1dLFsxLDMsIiIsMCx7ImxldmVsIjoyLCJzdHlsZSI6eyJoZWFkIjp7Im5hbWUiOiJub25lIn19fV0sWzAsMiwiIiwyLHsibGV2ZWwiOjIsInN0eWxlIjp7ImhlYWQiOnsibmFtZSI6Im5vbmUifX19XSxbMiwzLCJHbiIsMix7InN0eWxlIjp7ImJvZHkiOnsibmFtZSI6ImJhcnJlZCJ9fX1dLFs1LDEsIlxcYWxwaGFfeSJdLFs0LDAsIlxcYWxwaGFfeCIsMl0sWzQsNSwiRm0iLDAseyJzdHlsZSI6eyJib2R5Ijp7Im5hbWUiOiJiYXJyZWQifX19XSxbNiw5LCJHXFxnYW1tYSIsMSx7InNob3J0ZW4iOnsic291cmNlIjoyMCwidGFyZ2V0IjoyMH0sInN0eWxlIjp7ImJvZHkiOnsibmFtZSI6Im5vbmUifSwiaGVhZCI6eyJuYW1lIjoibm9uZSJ9fX1dLFsxMiw2LCJcXGFscGhhX20iLDEseyJsYWJlbF9wb3NpdGlvbiI6NDAsInNob3J0ZW4iOnsic291cmNlIjoyMCwidGFyZ2V0IjoyMH0sInN0eWxlIjp7ImJvZHkiOnsibmFtZSI6Im5vbmUifSwiaGVhZCI6eyJuYW1lIjoibm9uZSJ9fX1dXQ==
    \begin{tikzcd}
      Fx & Fy \\
      Gx & Gy \\
      Gx & Gy
      \arrow[""{name=0, anchor=center, inner sep=0}, "Fm"{inner sep=.8ex}, "\shortmid"{marking}, from=1-1, to=1-2]
      \arrow["{\alpha_x}"', from=1-1, to=2-1]
      \arrow["{\alpha_y}", from=1-2, to=2-2]
      \arrow[""{name=1, anchor=center, inner sep=0}, "Gm"{inner sep=.8ex}, "\shortmid"{marking}, from=2-1, to=2-2]
      \arrow[equals, from=2-1, to=3-1]
      \arrow[equals, from=2-2, to=3-2]
      \arrow[""{name=2, anchor=center, inner sep=0}, "Gn"'{inner sep=.8ex}, "\shortmid"{marking}, from=3-1, to=3-2]
      \arrow["{\alpha_m}"{description, pos=0.4}, draw=none, from=0, to=1]
      \arrow["{G\gamma}"{description}, draw=none, from=1, to=2]
    \end{tikzcd},
  \end{equation*}
  where the equality on the right uses the naturality of $\alpha$ with respect
  to cells.
\end{construction}

\begin{proposition}[Double barrels to twisted copresheaves]
\label{prop:double-barrel-to-twisted-copresheaf}
There is a 2-functor
\begin{equation*}
\twsection{-}: \lModl \to \Twnll
\end{equation*}
that sends
\begin{itemize}[noitemsep]
  \item a double barrel $M : \dbl{M} \to \Loose$ to its twisted bimodule of sections $\twsection{M}$ (\cref{construction:twisted-bimod-sections})
  \item a lax morphism $F : \dbl{M} \to \dbl{N}$ of double barrels (satisfying $M = N \circ F$) to the lax natural transformation $\twsection{F}$ (\cref{construction:loose-bimod-morphism})
  \item a 2-cell $\alpha : F \To G$ in $\lModl$ to the relative modification $\twsection{\alpha}$ (\cref{construction:relative-modification-from-barrel}).
\end{itemize}
\end{proposition}
\begin{proof}[Proof sketch]
 By \cref{prop:embedding-into-twisted-copresheaves}, the assignment $\dbl{M} \mapsto \dbl{M}(-,=)$ is 2-functorial; all that needs to be checked is that this 2-functor restricts appropriately to $\twsection{M}$, since $\twsection{M} \coloneqq \dbl{M}(\inl(-), \inr(=))$ is the restriction along the inclusions $ \inl : \barrelsrc{\dbl{M}} \hookrightarrow \dbl{M}$ and $\inr : \barreltgt{\dbl{M}} \hookrightarrow \dbl{M}$. This follows straightforwardly from the strict commutativity over $\Loose$.
\end{proof}

It is clear from its construction that we may see the 2-functor $\twsection{-}$
as landing not just in the 2-category of twisted copresheaves, but in the
following 2-category of \emph{twisted bimodules}.

\begin{definition}[Twisted bimodule] \label{def:twisted-bimodule}
Let $\dbl{D}$ and $\dbl{E}$ be double categories. A \define{twisted bimodule} $M : \dbl{D} \twbimodto \dbl{E}$ is a twisted normal lax functor $M : \dbl{D}^{\co} \times \dbl{E} \twistto \Prof$.

We define the 2-category $\TwModl$ of twisted bimodules to be the following pullback of 2-categories:
\[
% https://q.uiver.app/#q=WzAsNCxbMSwwLCJcXFR3bmxsIl0sWzEsMSwiXFxEYmxsIl0sWzAsMSwiXFxEYmxsIFxcdGltZXMgXFxEYmxsIl0sWzAsMCwiXFxidWxsZXQiXSxbMiwxLCIoLSlee1xcbWF0aHNme2NvfX0gXFx0aW1lcyAoLSkiLDJdLFswLDFdLFszLDJdLFszLDBdLFszLDEsIiIsMSx7InN0eWxlIjp7Im5hbWUiOiJjb3JuZXIifX1dXQ==
\begin{tikzcd}
	\TwModl & \Twnll \\
	{\Dbll \times \Dbll} & \Dbll
	\arrow[from=1-1, to=1-2]
	\arrow[from=1-1, to=2-1]
	\arrow["\lrcorner"{anchor=center, pos=0.125}, draw=none, from=1-1, to=2-2]
	\arrow[from=1-2, to=2-2]
	\arrow["{(-)^{\mathsf{co}} \times (-)}"', from=2-1, to=2-2]
\end{tikzcd},
\]
where $\Twnll$ was defined in \cref{construction:twisted-lax-2-cat}. In other words,
the objects of $\TwModl$ consist of a pair of double categories $\dbl{D}$ and $\dbl{E}$ along with a twisted bimodule $\dbl{D} \twbimodto \dbl{E}$ between them.
\end{definition}

\subsection{Collage of a twisted bimodule}
\label{subsection:collage}

In this section, we will give an inverse to \cref{prop:double-barrel-to-twisted-copresheaf} given by the \emph{collage} construction.

\begin{theorem}[Equivalence of twisted bimodules and double barrels]
\label{thm:equivalence-of-twisted-bimod-dbl-barrel}
There is a 2-functor $\collage{-} : \TwModl \to \lModl$ sending a twisted bimodule to its \emph{collage}. 
\end{theorem}

The collage $\collage{M}$ of a twisted bimodule $M : \dbl{D}^{\co} \times \dbl{E} \twistto \Prof$ is a double category equipped with a labeling $\collage{M} \to \Loose$ in the walking loose arrow for which $\barrelsrc{\collage{M}} = \dbl{D}$ and $\barreltgt{\collage{M}} = \dbl{E}$. The idea is that there is a loose morphism $m : x \proto y$ in $\collage{M}$ from $x \in \dbl{D}$ to $y \in \dbl{E}$ for every $m \in M(x, y)$, and composition is given by the action of $M$ as a twisted double functor. We now spell this out in detail.

\begin{construction}[Collage of a twisted bimodule]\label{construction:collage-twisted-bimod}
  Let $M : \dbl{D} \twbimodto \dbl{E}$ be a twisted bimodule. We define its \define{collage} $\collage{M}$ as follows:
  \begin{enumerate}
          \item The tight category $\collage{M}_{0}$ is the coproduct $\dbl{D}_{0} + \dbl{E}_{0}$ of the tight categories of $\dbl{D}$ and $\dbl{E}$. Write $\inl : \dbl{D}_{0} \to \collage{M}_{0}$ and $\inr : \dbl{E}_{0} \to \collage{M}_{0}$ for the respective inclusions.
    \item For objects $x$ and $y$, we define the set $\collage{M}_{1}(x, y)$ of loose morphisms $x \proto y$ to be loose morphisms of $\dbl{D}$ (or $\dbl{E}$) if both $x$ and $y$ are objects of $\dbl{D}$ (respectively, $\dbl{E}$), or objects of $M(x, y)$ if $x \in \dbl{D}$ and $y \in \dbl{E}$. Explicitly:
          \begin{enumerate}
                  \item $\collage{M}_{1}(\inl x, \inl y) \coloneqq \dbl{D}_{1}(x, y)$
                  \item $\collage{M}_{1}(\inr x, \inr y) \coloneqq \dbl{E}_{1}(x, y)$
                  \item $\collage{M}_{1}(\inl x, \inr y) \coloneqq M(x, y)_{0}$
                  \item $\collage{M}_{1}(\inr x, \inl y) \coloneqq \emptyset$.
          \end{enumerate}
    \item For tight morphisms $f$ and $g$, we define the set $\collage{M}_{1}(f, g)$ of squares similarly:
          \begin{enumerate}
                  \item $\collage{M}_{1}(\inl f, \inl g) \coloneqq \dbl{D}_{1}(f, g)$
                  \item $\collage{M}_{1}(\inr f, \inr g) \coloneqq \dbl{E}_{1}(f, g)$
                  \item $\collage{M}_{1}(\inl f, \inr g) \coloneqq \sum_{(m, n)} M(f, g)(m, n)$
                  \item $\collage{M}_{1}(\inr f, \inl g) \coloneqq \emptyset$.
          \end{enumerate}
    \item Tight composition is given by composition in $\dbl{D}$, $\dbl{E}$, or by the tight comparison cells $M_{(f_{0}, g_{0}),(f_{1}, g_{1})}$. Associativity and unitality follow by the same rules for the tight comparison cells of $M$.
    \item Composition of loose morphisms $\collage{M}_{1}(x, y) \times \collage{M}_{1}(y, z) \to \collage{M}_{1}(x, z)$ is either computed by loose composition in $\dbl{D}$ or $\dbl{E}$, or given by the action of $M$. Explicitly:
          \begin{enumerate}
                  \item  $\collage{M}_{1}(\inl x, \inl y) \times \collage{M}_{1}(\inl y, \inl z) \to \collage{M}_{1}(\inl x, \inl z)$ is composition in $\dbl{D}$.
                  \item  $\collage{M}_{1}(\inr x, \inr y) \times \collage{M}_{1}(\inr y, \inr z) \to \collage{M}_{1}(\inr x, \inr z)$ is composition in $\dbl{E}$.
                  \item  $\collage{M}_{1}(\inl x, \inl y) \times \collage{M}_{1}(\inl y, \inr z) \to \collage{M}_{1}(\inl x, \inr z)$ is given by $(m, n) \mapsto M(m, \id)(n)$.
                  \item  $\collage{M}_{1}(\inl x, \inr y) \times \collage{M}_{1}(\inr y, \inr z) \to \collage{M}_{1}(\inl x, \inr z)$ is given by $(m, n) \mapsto M(\id, n)(m)$.
          \end{enumerate}
      \item Loose composition of squares $\collage{M}_{1}(f, g) \times \collage{M}_{1}(g, h) \to \collage{M}_{1}(f, h)$ is defined similarly:
          \begin{enumerate}
                  \item  $\collage{M}_{1}(\inl f, \inl g) \times \collage{M}_{1}(\inl g, \inl h) \to \collage{M}_{1}(\inl f, \inl h)$ is composition in $\dbl{D}$.
                  \item  $\collage{M}_{1}(\inr f, \inr g) \times \collage{M}_{1}(\inr g, \inr h) \to \collage{M}_{1}(\inr f, \inr h)$ is composition in $\dbl{E}$.
                  \item  $\collage{M}_{1}(\inl f, \inl g) \times \collage{M}_{1}(\inl g, \inr h) \to \collage{M}_{1}(\inl f, \inr h)$ is given by $(\alpha, \beta) \mapsto M(\alpha, \id)(\beta)$.
                  \item  $\collage{M}_{1}(\inl f, \inr g) \times \collage{M}_{1}(\inr g, \inr h) \to \collage{M}_{1}(\inl f, \inr h)$ is given by $(\alpha, \beta) \mapsto M(\id, \beta)(\alpha)$.
          \end{enumerate}
    \item The associators and unitors for loose composition are defined using the loose-to-tight comparison cells of the twisted bimodule:
          \begin{enumerate}
            \item Either the associators are taken in $\dbl{D}$ or $\dbl{E}$, or we are in one of three situations:
                  \begin{enumerate}
                    \item We have loose morphisms $x \xproto{m} y \xproto{n} z$ in $\dbl{D}$ and $p \in M(z, w)$, in which case we take the component at $p$ of the natural transformation
                      \[\begin{tikzcd}
                          {M(z, w)} & {M(z,w)} & {M(z,w)} \\
                          & {M(z, w)} & {M(z, w)} \\
                          {M(y, w)} \\
                          & {M(x, w)} & {M(x, w)} \\
                          {M(x, w)} & {M(x,w)} & {M(x,w)}
                          \arrow[""{name=0, anchor=center, inner sep=0}, "\shortmid"{marking}, equals, from=1-1, to=1-2]
                          \arrow["{M(n, \id)}"', from=1-1, to=3-1]
                          \arrow[""{name=1, anchor=center, inner sep=0}, "\shortmid"{marking}, equals, from=1-2, to=1-3]
                          \arrow[equals, from=1-2, to=2-2]
                          \arrow[equals, from=1-3, to=2-3]
                          \arrow[""{name=2, anchor=center, inner sep=0}, "{M(1,1)}"'{inner sep=.8ex}, "\shortmid"{marking}, from=2-2, to=2-3]
                          \arrow[""{name=3, anchor=center, inner sep=0}, "{M(m \odot n, \id \odot \id)}"{description}, from=2-2, to=4-2]
                          \arrow[""{name=4, anchor=center, inner sep=0}, "{M(m \odot n, \id)}", from=2-3, to=4-3]
                          \arrow["{M(m, \id)}"', from=3-1, to=5-1]
                          \arrow[""{name=5, anchor=center, inner sep=0}, "{M(1,1)}"{inner sep=.8ex}, "\shortmid"{marking}, from=4-2, to=4-3]
                          \arrow[equals, from=4-2, to=5-2]
                          \arrow[equals, from=4-3, to=5-3]
                          \arrow[""{name=6, anchor=center, inner sep=0}, "\shortmid"{marking}, equals, from=5-1, to=5-2]
                          \arrow[""{name=7, anchor=center, inner sep=0}, "\shortmid"{marking}, equals, from=5-2, to=5-3]
                          \arrow["{M^{(n,\id),(m,\id)}}"{description, pos=0.4}, draw=none, from=0, to=6]
                          \arrow["{M_{(z,w)}}"{description}, draw=none, from=1, to=2]
                          \arrow["{M(\cong)}"{description, pos=0.7}, draw=none, from=3, to=4]
                          \arrow["{M_{(x,w)}^{-1}}"{description}, draw=none, from=5, to=7]
                        \end{tikzcd}\]
                      where in the middle right we are using the unitor isomorphism $\id \odot \id \cong \id$ in $\dbl{E}$.
                    \item We have $p \in M(x, y)$ and $y \xproto{m} z \xproto{n} w$ in $\dbl{E}$, in which case we do essentially the same thing but on the other side.
                    \item We have $m: x \proto y$ in $\dbl{D}$, $p \in M(y, z)$, and $n: z \proto w$ in $\dbl{E}$, in which case we take the component at $p$ of the natural transformation
                      \[\begin{tikzcd}
                          {M(y, z)} & {M(y,z)} & {M(y,z)} \\
                          & {M(y, z)} & {M(y, z)} \\
                          {M(x, z)} \\
                          & {M(x, w)} & {M(x, w)} \\
                          {M(x, w)} & {M(x,w)} & {M(x,w)}
                          \arrow[""{name=0, anchor=center, inner sep=0}, "\shortmid"{marking}, equals, from=1-1, to=1-2]
                          \arrow["{M(m, \id)}"', from=1-1, to=3-1]
                          \arrow[""{name=1, anchor=center, inner sep=0}, "\shortmid"{marking}, equals, from=1-2, to=1-3]
                          \arrow[equals, from=1-2, to=2-2]
                          \arrow[equals, from=1-3, to=2-3]
                          \arrow[""{name=2, anchor=center, inner sep=0}, "{M(1,1)}"'{inner sep=.8ex}, "\shortmid"{marking}, from=2-2, to=2-3]
                          \arrow[""{name=3, anchor=center, inner sep=0}, "{M(m \odot \id, \id \odot n)}"{description}, from=2-2, to=4-2]
                          \arrow[""{name=4, anchor=center, inner sep=0}, "{M(m, n)}", from=2-3, to=4-3]
                          \arrow["{M(\id, n)}"', from=3-1, to=5-1]
                          \arrow[""{name=5, anchor=center, inner sep=0}, "{M(1,1)}"{inner sep=.8ex}, "\shortmid"{marking}, from=4-2, to=4-3]
                          \arrow[equals, from=4-2, to=5-2]
                          \arrow[equals, from=4-3, to=5-3]
                          \arrow[""{name=6, anchor=center, inner sep=0}, "\shortmid"{marking}, equals, from=5-1, to=5-2]
                          \arrow[""{name=7, anchor=center, inner sep=0}, "\shortmid"{marking}, equals, from=5-2, to=5-3]
                          \arrow["{M^{(m,\id),(\id,n)}}"{description, pos=0.4}, draw=none, from=0, to=6]
                          \arrow["{M_{(y,z)}}"{description}, draw=none, from=1, to=2]
                          \arrow["{M(\cong)}"{description, pos=0.7}, draw=none, from=3, to=4]
                          \arrow["{M_{(x,w)}^{-1}}"{description}, draw=none, from=5, to=7]
                        \end{tikzcd}\]
                      where in the middle right we use the pair of unitor isomorphisms $(m \odot \id \cong m, \id \odot n \cong n)$ in $\dbl{D}$ and $\dbl{E}$.
                  \end{enumerate}
            \item Either the unitors are taken in $\dbl{D}$ or $\dbl{E}$, or we have $m \in M(x, y)$ and are composing it on the left or right by a loose identity; in either case, we use the loose-to-tight comparison cell $M^{(x,y)} : 1_{M(x,y)} \To M(\id_x,\id_y)$.
          \end{enumerate}
    \item The coherence equations for $\collage{M}$ follow by those for $\dbl{D}$ or $\dbl{E}$, or by those for $M$.
  \end{enumerate}
We denote the evident inclusions by $\inl : \dbl{D} \hookrightarrow \collage{M}$ and $\inr : \dbl{E} \hookrightarrow \collage{M}$ respectively.
\end{construction}

\begin{remark}
  In the following, it will be convenient to abuse notation and identify $m \in M(x,y)$ with $M_{(x,y)}(1_m) \in M(1_x,1_y)$. This is most easily justified if the twisted bimodule is strictly unitary in the tight-to-loose direction, but it is straightforward to accommodate even in the normal case.
\end{remark}

We can now define the labeling of a collage in the walking loose arrow. Together with its labeling, the collage becomes a double barrel; we expect that this gives an inverse to the $\twsection{-}$ construction (see \cref{conjecture:sec.collage}).
\begin{construction}[Labeling] \label{construction:collage-labeling}
  Let $M : \dbl{D} \twbimodto \dbl{E}$ be a twisted bimodule. We define the \define{labeling} of its collage $\ell_{M} : \collage{M} \to \Loose$ by sending all of $\dbl{D}$ to $0$ (or its identities) and $\dbl{E}$ to $1$ (or its identities) and the rest to the walking loose arrow itself (or its identity).
\end{construction}

Next, we upgrade the collage construction to a functor.

\begin{construction}[Collage construction on morphisms] \label{construction:collage-transformation}
Let $\kappa : M(-,=) \To N(F^{\co}(-), G(=))$ be a morphism of twisted bimodules (1-cell in $\TwModl$), where $M: \dbl{D} \twbimodto \dbl{E}$ and $N : \dbl{A} \twbimodto \dbl{B}$ are twisted bimodules and $F : \dbl{D} \to \dbl{A}$ and $G : \dbl{E} \to \dbl{B}$ are lax functors. We construct a lax functor $\collage{\kappa} : \collage{M} \to \collage{N}$ as follows:
          \begin{itemize}[noitemsep]
            \item $\collage{\kappa}$ acts as $F$ on objects, tights, looses, and squares of $\dbl{D}$, and similarly acts as $G$ on elements of $\dbl{E}$.
            \item For $x \in \dbl{D}$ and $y \in \dbl{E}$ and $m : x \proto y$ in $\collage{M}$ (so, $m \in M(x,y)$), we define \[\collage{\kappa}(m) \coloneqq \kappa_{(x,y)}(m) \in N(Fx,Gy).\]
            \item For a cell $\stdInlineCell{\beta}$ with $f \in \dbl{D}$ and $g \in \dbl{E}$ (so, $\beta \in M(f,g)(m,n)$) we define \[\collage{\kappa}(\beta) \coloneqq \kappa_{(f,g)}(\beta) \in N(Ff,Gg)(\kappa_{(x,y)}(m),\kappa_{(w,z)}(n)).\]
            \item For $x \xproto{m} y \xproto{n} z$ in $\collage{M}$, we define the compositor $\collage{\kappa}_{m, n}$ as follows:
                  \begin{enumerate}
                          \item If $m$ and $n$ are in $\dbl{D}$ or $\dbl{E}$, we define $\collage{\kappa}$ to be the compositor of $F$ or $G$ respectively.
                          \item If $m$ is in $\dbl{D}$ and $z \in \dbl{E}$, we define the compositor $\collage{\kappa}_{m,n} : \collage{\kappa}(m) \odot \collage{\kappa}(n) \To \collage{\kappa}(m \odot n)$ to be the component of the following cell at $n \in M(y, z)$:
              
                          \[
                           % https://q.uiver.app/#q=WzAsOSxbMSwwLCJNKHkseikiXSxbMSwxLCJOKEZ5LEd6KSJdLFsxLDIsIk4oRngsR3opIl0sWzIsMCwiTSh5LHopIl0sWzIsMSwiTSh4LHopIl0sWzIsMiwiTihGeCxHeikiXSxbMCwxLCJOKEZ5LEd6KSJdLFswLDIsIk4oRnksR3opIl0sWzAsMCwiTSh5LCB6KSJdLFswLDEsIlxca2FwcGFfeyh5LHopfSIsMl0sWzEsMiwiTihGbSxHXFxpZF96KSJdLFszLDQsIk0obSxcXGlkX3opIl0sWzQsNSwiXFxrYXBwYV97KHgseil9Il0sWzAsMywiIiwwLHsibGV2ZWwiOjIsInN0eWxlIjp7ImJvZHkiOnsibmFtZSI6ImJhcnJlZCJ9LCJoZWFkIjp7Im5hbWUiOiJub25lIn19fV0sWzIsNSwiIiwyLHsibGV2ZWwiOjIsInN0eWxlIjp7ImJvZHkiOnsibmFtZSI6ImJhcnJlZCJ9LCJoZWFkIjp7Im5hbWUiOiJub25lIn19fV0sWzYsNywiTihGbSwgXFxpZF97R3p9KSIsMl0sWzYsMSwiIiwwLHsibGV2ZWwiOjIsInN0eWxlIjp7ImJvZHkiOnsibmFtZSI6ImJhcnJlZCJ9LCJoZWFkIjp7Im5hbWUiOiJub25lIn19fV0sWzcsMiwiIiwwLHsibGV2ZWwiOjIsInN0eWxlIjp7ImJvZHkiOnsibmFtZSI6ImJhcnJlZCJ9LCJoZWFkIjp7Im5hbWUiOiJub25lIn19fV0sWzgsMCwiIiwxLHsibGV2ZWwiOjIsInN0eWxlIjp7ImJvZHkiOnsibmFtZSI6ImJhcnJlZCJ9LCJoZWFkIjp7Im5hbWUiOiJub25lIn19fV0sWzgsNiwiXFxrYXBwYV97KHkseil9IiwyXSxbMTMsMTQsIlxca2FwcGFfeyhtLFxcaWRfen0pIiwxLHsic2hvcnRlbiI6eyJzb3VyY2UiOjIwLCJ0YXJnZXQiOjIwfSwic3R5bGUiOnsiYm9keSI6eyJuYW1lIjoibm9uZSJ9LCJoZWFkIjp7Im5hbWUiOiJub25lIn19fV0sWzE2LDE3LCJOKFxcaWRfe0ZtfSwgR196KSIsMSx7InNob3J0ZW4iOnsic291cmNlIjoyMCwidGFyZ2V0IjoyMH0sInN0eWxlIjp7ImJvZHkiOnsibmFtZSI6Im5vbmUifSwiaGVhZCI6eyJuYW1lIjoibm9uZSJ9fX1dXQ==
\begin{tikzcd}
	{M(y, z)} & {M(y,z)} & {M(y,z)} \\
	{N(Fy,Gz)} & {N(Fy,Gz)} & {M(x,z)} \\
	{N(Fy,Gz)} & {N(Fx,Gz)} & {N(Fx,Gz)}
	\arrow["\shortmid"{marking}, equals, from=1-1, to=1-2]
	\arrow["{\kappa_{(y,z)}}"', from=1-1, to=2-1]
	\arrow[""{name=0, anchor=center, inner sep=0}, "\shortmid"{marking}, equals, from=1-2, to=1-3]
	\arrow["{\kappa_{(y,z)}}"', from=1-2, to=2-2]
	\arrow["{M(m,\id_z)}", from=1-3, to=2-3]
	\arrow[""{name=1, anchor=center, inner sep=0}, "\shortmid"{marking}, equals, from=2-1, to=2-2]
	\arrow["{N(Fm, \id_{Gz})}"', from=2-1, to=3-1]
	\arrow["{N(Fm,G\id_z)}", from=2-2, to=3-2]
	\arrow["{\kappa_{(x,z)}}", from=2-3, to=3-3]
	\arrow[""{name=2, anchor=center, inner sep=0}, "\shortmid"{marking}, equals, from=3-1, to=3-2]
	\arrow[""{name=3, anchor=center, inner sep=0}, "\shortmid"{marking}, equals, from=3-2, to=3-3]
	\arrow["{\kappa_{(m,\id_z)}}"{description}, draw=none, from=0, to=3]
	\arrow["{N(\id_{Fm}, G_z)}"{description}, draw=none, from=1, to=2]
\end{tikzcd} 
                          \]
                          \item If $n$ is in $\dbl{E}$ and $x \in \dbl{D}$, we define the compositor $\collage{\kappa}_{m,n} : \collage{\kappa}(m) \odot \collage{\kappa}(n) \To \collage{\kappa}(m \odot n)$ to be the component of the following cell at $m \in M(x, y)$:
                          \[
% https://q.uiver.app/#q=WzAsOSxbMSwwLCJNKHgseSkiXSxbMSwxLCJOKEZ4LEd5KSJdLFsxLDIsIk4oRngsR3opIl0sWzIsMCwiTSh4LCB5KSJdLFsyLDEsIk0oeCx6KSJdLFsyLDIsIk4oRngsR3opIl0sWzAsMSwiTihGeCxHeSkiXSxbMCwyLCJOKEZ4LEd6KSJdLFswLDAsIk0oeCwgeSkiXSxbMCwxLCJcXGthcHBhX3soeSx6KX0iLDJdLFsxLDIsIk4oRlxcaWRfeCxHbikiXSxbMyw0LCJNKFxcaWRfeCxuKSJdLFs0LDUsIlxca2FwcGFfeyh4LHopfSJdLFswLDMsIiIsMCx7ImxldmVsIjoyLCJzdHlsZSI6eyJib2R5Ijp7Im5hbWUiOiJiYXJyZWQifSwiaGVhZCI6eyJuYW1lIjoibm9uZSJ9fX1dLFsyLDUsIiIsMix7ImxldmVsIjoyLCJzdHlsZSI6eyJib2R5Ijp7Im5hbWUiOiJiYXJyZWQifSwiaGVhZCI6eyJuYW1lIjoibm9uZSJ9fX1dLFs2LDcsIk4oXFxpZF97Rnh9LCBHbikiLDJdLFs2LDEsIiIsMCx7ImxldmVsIjoyLCJzdHlsZSI6eyJib2R5Ijp7Im5hbWUiOiJiYXJyZWQifSwiaGVhZCI6eyJuYW1lIjoibm9uZSJ9fX1dLFs3LDIsIiIsMCx7ImxldmVsIjoyLCJzdHlsZSI6eyJib2R5Ijp7Im5hbWUiOiJiYXJyZWQifSwiaGVhZCI6eyJuYW1lIjoibm9uZSJ9fX1dLFs4LDAsIiIsMSx7ImxldmVsIjoyLCJzdHlsZSI6eyJib2R5Ijp7Im5hbWUiOiJiYXJyZWQifSwiaGVhZCI6eyJuYW1lIjoibm9uZSJ9fX1dLFs4LDYsIlxca2FwcGFfeyh5LHopfSIsMl0sWzEzLDE0LCJcXGthcHBhX3soXFxpZF94LCBufSkiLDEseyJzaG9ydGVuIjp7InNvdXJjZSI6MjAsInRhcmdldCI6MjB9LCJzdHlsZSI6eyJib2R5Ijp7Im5hbWUiOiJub25lIn0sImhlYWQiOnsibmFtZSI6Im5vbmUifX19XSxbMTYsMTcsIk4oRl94LFxcaWRfe0dufSkiLDEseyJzaG9ydGVuIjp7InNvdXJjZSI6MjAsInRhcmdldCI6MjB9LCJzdHlsZSI6eyJib2R5Ijp7Im5hbWUiOiJub25lIn0sImhlYWQiOnsibmFtZSI6Im5vbmUifX19XV0=
\begin{tikzcd}
	{M(x, y)} & {M(x,y)} & {M(x, y)} \\
	{N(Fx,Gy)} & {N(Fx,Gy)} & {M(x,z)} \\
	{N(Fx,Gz)} & {N(Fx,Gz)} & {N(Fx,Gz)}
	\arrow["\shortmid"{marking}, equals, from=1-1, to=1-2]
	\arrow["{\kappa_{(x,y)}}"', from=1-1, to=2-1]
	\arrow[""{name=0, anchor=center, inner sep=0}, "\shortmid"{marking}, equals, from=1-2, to=1-3]
	\arrow["{\kappa_{(x,y)}}"', from=1-2, to=2-2]
	\arrow["{M(\id_x,n)}", from=1-3, to=2-3]
	\arrow[""{name=1, anchor=center, inner sep=0}, "\shortmid"{marking}, equals, from=2-1, to=2-2]
	\arrow["{N(\id_{Fx}, Gn)}"', from=2-1, to=3-1]
	\arrow["{N(F\id_x,Gn)}", from=2-2, to=3-2]
	\arrow["{\kappa_{(x,z)}}", from=2-3, to=3-3]
	\arrow[""{name=2, anchor=center, inner sep=0}, "\shortmid"{marking}, equals, from=3-1, to=3-2]
	\arrow[""{name=3, anchor=center, inner sep=0}, "\shortmid"{marking}, equals, from=3-2, to=3-3]
	\arrow["{\kappa_{(\id_x, n)}}"{description}, draw=none, from=0, to=3]
	\arrow["{N(F_x,\id_{Gn})}"{description}, draw=none, from=1, to=2]
\end{tikzcd}
                          \]
                  \end{enumerate}
          \end{itemize}
\end{construction}

\begin{lemma}
 \cref{construction:collage-transformation} gives a well defined lax double functor.
\end{lemma}
\begin{proof}
  The equations we need follow from the corresponding laws (\cref{def:naturaltransformation}) defining the lax transformation $\kappa$.

  The four laws governing a lax double functor, as they appear in \cite[\mbox{Definition 3.5.1}]{grandis2019}, follow smoothly once the definitions involved have been fully expanded. We address the case of \emph{Naturality of composition comparisons} in full; for the others we simply describe the key law from \cref{def:naturaltransformation}.
  \begin{enumerate}
    \item (\emph{Naturality of identity comparisons}) This follows by the appropriate conditions for $F$ and $G$, since by definition $\collage{\kappa}$ acts like $F$ or $G$ on $\dbl{D}$ and $\dbl{E}$ respectively.
    \item (\emph{Naturality of composition comparisons}) Let $\beta_{1}$ and $\beta_{2}$ be loosely composable cells in $\collage{M}$. We need to show that
          \[
\begin{tikzcd}
	\collage{\kappa}x & \collage{\kappa}y & \collage{\kappa}z && \collage{\kappa}x & \collage{\kappa}y & \collage{\kappa}z \\
	\collage{\kappa}x && \collage{\kappa}z & {=} & {\collage{\kappa}x'} & {\collage{\kappa}y'} & {\collage{\kappa}z'} \\
	{\collage{\kappa}x'} && \collage{\kappa}z && {\collage{\kappa}x'} && {\collage{\kappa}x'}
	\arrow["{\collage{\kappa}(m)}", "\shortmid"{marking}, from=1-1, to=1-2]
	\arrow[Rightarrow, no head, from=1-1, to=2-1]
	\arrow["{\collage{\kappa}(n)}", "\shortmid"{marking}, from=1-2, to=1-3]
	\arrow[Rightarrow, no head, from=1-3, to=2-3]
	\arrow[""{name=0, anchor=center, inner sep=0}, "{\collage{\kappa}(m)}", "\shortmid"{marking}, from=1-5, to=1-6]
	\arrow[from=1-5, to=2-5]
	\arrow[""{name=1, anchor=center, inner sep=0}, "{\collage{\kappa}(n)}", "\shortmid"{marking}, from=1-6, to=1-7]
	\arrow[from=1-6, to=2-6]
	\arrow[from=1-7, to=2-7]
	\arrow[""{name=2, anchor=center, inner sep=0}, "{\collage{\kappa}(m \odot n)}"', "\shortmid"{marking}, from=2-1, to=2-3]
	\arrow[from=2-1, to=3-1]
	\arrow[from=2-3, to=3-3]
	\arrow[""{name=3, anchor=center, inner sep=0}, "{\collage{\kappa}(m')}"', "\shortmid"{marking}, from=2-5, to=2-6]
	\arrow[Rightarrow, no head, from=2-5, to=3-5]
	\arrow[""{name=4, anchor=center, inner sep=0}, "{\collage{\kappa}(n')}"', "\shortmid"{marking}, from=2-6, to=2-7]
	\arrow[Rightarrow, no head, from=2-7, to=3-7]
	\arrow[""{name=5, anchor=center, inner sep=0}, "{\collage{\kappa}(m' \odot n')}"', "\shortmid"{marking}, from=3-1, to=3-3]
	\arrow[""{name=6, anchor=center, inner sep=0}, "{\collage{\kappa}(m' \odot n')}"', "\shortmid"{marking}, from=3-5, to=3-7]
	\arrow["{\collage{\kappa}_{m,n}}"{description}, draw=none, from=1-2, to=2]
	\arrow["{\collage{\kappa}(\beta_1)}"{description}, draw=none, from=0, to=3]
	\arrow["{\collage{\kappa}(\beta_2)}"{description}, draw=none, from=1, to=4]
	\arrow["{\collage{\kappa}(\beta_1\mid \beta_2)}"{description}, draw=none, from=2, to=5]
	\arrow["{\collage{\kappa}(m',n')}"{description}, draw=none, from=2-6, to=6]
\end{tikzcd}
          \]
          We must work by case analysis. If both $\beta_{1}$ and $\beta_{2}$ are in $\dbl{D}$ or $\dbl{E}$, then $C \coloneqq \collage{\kappa}$ acts as $F$ or $G$ and this is the tight naturality of laxators. Now consider the case that $\beta_{2}$ is in $\dbl{E}$ and $\beta_{1}$ is in $M(f, g)(m, m')$ (the other case will work similarly). Expanding the definitions, we see that we need to prove
\[
% https://q.uiver.app/#q=WzAsNixbMCwwLCJGeCJdLFszLDAsIkd6Il0sWzAsMSwiRngiXSxbMywxLCJHeiJdLFswLDIsIkZ4JyJdLFszLDIsIkd6Il0sWzAsMiwiIiwwLHsibGV2ZWwiOjIsInN0eWxlIjp7ImhlYWQiOnsibmFtZSI6Im5vbmUifX19XSxbMiw0XSxbMCwxLCJOKFxcaWQsIG4pKFxcYWxwaGFfeyh4LCB5KX0obSkpIiwwLHsic3R5bGUiOnsiYm9keSI6eyJuYW1lIjoiYmFycmVkIn19fV0sWzIsMywiXFxrYXBwYXsoeCx6KX1NKFxcaWQsbikobSkiLDIseyJzdHlsZSI6eyJib2R5Ijp7Im5hbWUiOiJiYXJyZWQifX19XSxbNCw1LCJcXGFscGhhX3soeCx6KX1NKFxcaWQsIG4nKShtJykiLDIseyJzdHlsZSI6eyJib2R5Ijp7Im5hbWUiOiJiYXJyZWQifX19XSxbMyw1XSxbMSwzLCIiLDAseyJsZXZlbCI6Miwic3R5bGUiOnsiaGVhZCI6eyJuYW1lIjoibm9uZSJ9fX1dLFs5LDEwLCJcXGthcHBhX3soZixoKX1NKFxcaWQsIFxcYmV0YV8yKShcXGJldGFfMSkiLDEseyJzaG9ydGVuIjp7InNvdXJjZSI6MjAsInRhcmdldCI6MjB9LCJzdHlsZSI6eyJib2R5Ijp7Im5hbWUiOiJub25lIn0sImhlYWQiOnsibmFtZSI6Im5vbmUifX19XSxbOCw5LCJOKEZfeCwgXFxpZF97R259KShcXGthcHBhX3soeCx5KX0obSkpO1xca2FwcGFfeyhcXGlkLCBuKX0obSkiLDEseyJzaG9ydGVuIjp7InRhcmdldCI6MjB9LCJzdHlsZSI6eyJib2R5Ijp7Im5hbWUiOiJub25lIn0sImhlYWQiOnsibmFtZSI6Im5vbmUifX19XV0=
\begin{tikzcd}[column sep = huge]
	Fx &&& Gz \\
	Fx &&& Gz \\
	{Fx'} &&& Gz
	\arrow[""{name=0, anchor=center, inner sep=0}, "{N(\id, n)(\alpha_{(x, y)}(m))}"{inner sep=.8ex}, "\shortmid"{marking}, from=1-1, to=1-4]
	\arrow[equals, from=1-1, to=2-1]
	\arrow[equals, from=1-4, to=2-4]
	\arrow[""{name=1, anchor=center, inner sep=0}, "{\kappa{(x,z)}M(\id,n)(m)}"'{inner sep=.8ex}, "\shortmid"{marking}, from=2-1, to=2-4]
	\arrow[from=2-1, to=3-1]
	\arrow[from=2-4, to=3-4]
	\arrow[""{name=2, anchor=center, inner sep=0}, "{\alpha_{(x,z)}M(\id, n')(m')}"'{inner sep=.8ex}, "\shortmid"{marking}, from=3-1, to=3-4]
	\arrow["{N(F_x, \id_{Gn})(\kappa_{(x,y)}(m));\kappa_{(\id, n)}(m)}"{description}, draw=none, from=0, to=1]
	\arrow["{\kappa_{(f,h)}M(\id, \beta_2)(\beta_1)}"{description}, draw=none, from=1, to=2]
\end{tikzcd}
\quad = \quad
  % https://q.uiver.app/#q=WzAsNixbMCwwLCJGeCJdLFszLDAsIkd6Il0sWzAsMSwiRngnIl0sWzMsMSwiR3onIl0sWzAsMiwiRngnIl0sWzMsMiwiR3gnIl0sWzAsMl0sWzIsNCwiIiwwLHsibGV2ZWwiOjIsInN0eWxlIjp7ImhlYWQiOnsibmFtZSI6Im5vbmUifX19XSxbNCw1LCJcXGthcHBhX3soeCx6KX1NKFxcaWQsIG4nKShtJykiLDIseyJzdHlsZSI6eyJib2R5Ijp7Im5hbWUiOiJiYXJyZWQifX19XSxbMiwzLCJOKFxcaWQsIG4nKShcXGthcHBhX3soeCcseScpfShtJykpIiwyLHsic3R5bGUiOnsiYm9keSI6eyJuYW1lIjoiYmFycmVkIn19fV0sWzAsMSwiTihcXGlkLCBuKShcXGthcHBhX3soeCx5KX0obSkpIiwwLHsic3R5bGUiOnsiYm9keSI6eyJuYW1lIjoiYmFycmVkIn19fV0sWzEsM10sWzMsNSwiIiwyLHsibGV2ZWwiOjIsInN0eWxlIjp7ImhlYWQiOnsibmFtZSI6Im5vbmUifX19XSxbMTAsOSwiTihcXGlkLCBcXGJldGFfMikoXFxrYXBwYV97KGYsZyl9KFxcYmV0YV8xKSkiLDEseyJzaG9ydGVuIjp7InNvdXJjZSI6MjAsInRhcmdldCI6MjB9LCJzdHlsZSI6eyJib2R5Ijp7Im5hbWUiOiJub25lIn0sImhlYWQiOnsibmFtZSI6Im5vbmUifX19XSxbOSw4LCJOKEZfe3gnfSwgXFxpZF97R24nfSkoXFxrYXBwYV97KHgnLHknKX0obSkpO1xca2FwcGFfeyhcXGlkLCBuJyl9KG0nKSIsMSx7InNob3J0ZW4iOnsidGFyZ2V0IjoyMH0sInN0eWxlIjp7ImJvZHkiOnsibmFtZSI6Im5vbmUifSwiaGVhZCI6eyJuYW1lIjoibm9uZSJ9fX1dXQ==
\begin{tikzcd}[column sep = huge]
	Fx &&& Gz \\
	{Fx'} &&& {Gz'} \\
	{Fx'} &&& {Gx'}
	\arrow[""{name=0, anchor=center, inner sep=0}, "{N(\id, n)(\kappa_{(x,y)}(m))}"{inner sep=.8ex}, "\shortmid"{marking}, from=1-1, to=1-4]
	\arrow[from=1-1, to=2-1]
	\arrow[from=1-4, to=2-4]
	\arrow[""{name=1, anchor=center, inner sep=0}, "{N(\id, n')(\kappa_{(x',y')}(m'))}"'{inner sep=.8ex}, "\shortmid"{marking}, from=2-1, to=2-4]
	\arrow[equals, from=2-1, to=3-1]
	\arrow[equals, from=2-4, to=3-4]
	\arrow[""{name=2, anchor=center, inner sep=0}, "{\kappa_{(x,z)}M(\id, n')(m')}"'{inner sep=.8ex}, "\shortmid"{marking}, from=3-1, to=3-4]
	\arrow["{N(\id, \beta_2)(\kappa_{(f,g)}(\beta_1))}"{description}, draw=none, from=0, to=1]
	\arrow["{N(F_{x'}, \id_{Gn'})(\kappa_{(x',y')}(m));\kappa_{(\id, n')}(m')}"{description}, draw=none, from=1, to=2]
\end{tikzcd}
\]

          We take the following equation and apply it at $\beta_1 \in M(f, g)$:
          \[
\begin{tikzcd}
	& {M(x,y)} & {M(x,y)} & {M(x,y)} & {M(x',y')} \\
	& {N(Fx, Gy)} & {N(Fx,Gy)} & {M(x,z)} & {M(x', z')} \\
	& {N(Fx, Gz)} & {N(Fx, Gy)} & {N(Fx, Gz)} & {N(Fx', Gz')} \\
	& {M(x,y)} & {M(x, y)} & {M(x',y')} & {M(x',y')} \\
	{=} & {N(Fx, Gy)} & {N(Fx,Gy)} & {N(Fx', Gy')} & {M(x', z')} \\
	& {N(Fx,Gz)} & {N(Fx, Gy)} & {N(Fx', Gz')} & {N(Fx', Gz')} \\
	& {M(x,y)} & {M(x', y')} & {M(x',y')} & {M(x',y')} \\
	{=} & {N(Fx, Gy)} & {N(Fx',Gy')} & {N(Fx', Gy')} & {M(x', z')} \\
	& {N(Fx,Gz)} & {N(Fx', Gy')} & {N(Fx', Gz')} & {N(Fx', Gz')}
	\arrow["\shortmid"{marking}, equals, from=1-2, to=1-3]
	\arrow["{\kappa_{(x, y)}}"', from=1-2, to=2-2]
	\arrow[""{name=0, anchor=center, inner sep=0}, "\shortmid"{marking}, equals, from=1-3, to=1-4]
	\arrow["{\kappa_{(x,y)}}"', from=1-3, to=2-3]
	\arrow[""{name=1, anchor=center, inner sep=0}, "{M(f, g)}"{inner sep=.8ex}, "\shortmid"{marking}, from=1-4, to=1-5]
	\arrow[from=1-4, to=2-4]
	\arrow["{M(\id_{x'}, n')}", from=1-5, to=2-5]
	\arrow[""{name=2, anchor=center, inner sep=0}, "\shortmid"{marking}, equals, from=2-2, to=2-3]
	\arrow["{N(\id_{Fx}, Gn)}"', from=2-2, to=3-2]
	\arrow["{N(\id_x, n)}", from=2-3, to=3-3]
	\arrow[""{name=3, anchor=center, inner sep=0}, "{M(f, h)}"{inner sep=.8ex}, "\shortmid"{marking}, from=2-4, to=2-5]
	\arrow[from=2-4, to=3-4]
	\arrow["{\kappa_{(x', y')}}", from=2-5, to=3-5]
	\arrow[""{name=4, anchor=center, inner sep=0}, "\shortmid"{marking}, equals, from=3-2, to=3-3]
	\arrow[""{name=5, anchor=center, inner sep=0}, "\shortmid"{marking}, equals, from=3-3, to=3-4]
	\arrow[""{name=6, anchor=center, inner sep=0}, "{N(Ff, Gh)}"'{inner sep=.8ex}, "\shortmid"{marking}, from=3-4, to=3-5]
	\arrow["\shortmid"{marking}, equals, from=4-2, to=4-3]
	\arrow["{\kappa_{(x,y)}}"', from=4-2, to=5-2]
	\arrow[""{name=7, anchor=center, inner sep=0}, "{M(f, g)}"{inner sep=.8ex}, "\shortmid"{marking}, from=4-3, to=4-4]
	\arrow["{\kappa_{(x, y)}}"', from=4-3, to=5-3]
	\arrow[""{name=8, anchor=center, inner sep=0}, "\shortmid"{marking}, equals, from=4-4, to=4-5]
	\arrow["{\kappa_{x',y'}}", from=4-4, to=5-4]
	\arrow["{M(\id_{x'}, n')}", from=4-5, to=5-5]
	\arrow[""{name=9, anchor=center, inner sep=0}, "\shortmid"{marking}, equals, from=5-2, to=5-3]
	\arrow[from=5-2, to=6-2]
	\arrow[""{name=10, anchor=center, inner sep=0}, "{N(Ff, Gg)}"{inner sep=.8ex}, "\shortmid"{marking}, from=5-3, to=5-4]
	\arrow["{N(\id_x, n)}"{description}, from=5-3, to=6-3]
	\arrow["{N(\id_{x'}, n')}", from=5-4, to=6-4]
	\arrow["{\kappa_{(x', y')}}", from=5-5, to=6-5]
	\arrow[""{name=11, anchor=center, inner sep=0}, "\shortmid"{marking}, equals, from=6-2, to=6-3]
	\arrow[""{name=12, anchor=center, inner sep=0}, "{N(Ff, Gh)}"'{inner sep=.8ex}, "\shortmid"{marking}, from=6-3, to=6-4]
	\arrow[""{name=13, anchor=center, inner sep=0}, "\shortmid"{marking}, equals, from=6-4, to=6-5]
	\arrow[""{name=14, anchor=center, inner sep=0}, "\shortmid"{marking}, from=7-2, to=7-3]
	\arrow[from=7-2, to=8-2]
	\arrow["\shortmid"{marking}, equals, from=7-3, to=7-4]
	\arrow[from=7-3, to=8-3]
	\arrow[""{name=15, anchor=center, inner sep=0}, "\shortmid"{marking}, equals, from=7-4, to=7-5]
	\arrow[from=7-4, to=8-4]
	\arrow[from=7-5, to=8-5]
	\arrow[""{name=16, anchor=center, inner sep=0}, "\shortmid"{marking}, from=8-2, to=8-3]
	\arrow[from=8-2, to=9-2]
	\arrow[""{name=17, anchor=center, inner sep=0}, "\shortmid"{marking}, equals, from=8-3, to=8-4]
	\arrow[from=8-3, to=9-3]
	\arrow[from=8-4, to=9-4]
	\arrow[from=8-5, to=9-5]
	\arrow[""{name=18, anchor=center, inner sep=0}, "\shortmid"{marking}, from=9-2, to=9-3]
	\arrow[""{name=19, anchor=center, inner sep=0}, "\shortmid"{marking}, equals, from=9-3, to=9-4]
	\arrow[""{name=20, anchor=center, inner sep=0}, "\shortmid"{marking}, equals, from=9-4, to=9-5]
	\arrow["{\kappa_{\id_x, n}}"{description}, draw=none, from=0, to=5]
	\arrow["{M(\id_f, \beta_2)}"{description}, draw=none, from=1, to=3]
	\arrow["{N(F_x, \id_{Gn})}"{description}, draw=none, from=2, to=4]
	\arrow["{\kappa_{f, h}}"{description}, draw=none, from=3, to=6]
	\arrow["{\kappa_{f,g}}"{description}, draw=none, from=7, to=10]
	\arrow["{\kappa_{\id_{x'}, n'}}"{description}, draw=none, from=8, to=13]
	\arrow["{N(F_x, \id_{Gn})}"{description}, draw=none, from=9, to=11]
	\arrow["{N(\id_f, \beta_2)}"{description}, draw=none, from=10, to=12]
	\arrow["{\kappa_{f, g}}"{description}, draw=none, from=14, to=16]
	\arrow["{\kappa_{\id_{x'},n'}}"{description}, draw=none, from=15, to=20]
	\arrow["{N(\id_f, \beta_2)}"{description}, draw=none, from=16, to=18]
	\arrow["{N(F_{x'}, \id_{Gn'})}"{description}, draw=none, from=17, to=19]
\end{tikzcd}
          \]

The first equation follows by \emph{naturality with respect to cells} for $\kappa$, taking $\gamma \coloneqq (\id_{f}, \beta_{2})$. The second follows by tight naturality of unitors for $F$.

    \item (\emph{Coherence with unitors}) The unitors of $\collage{\kappa}$ are those of $F$ and $G$. The cases $m : x \proto y$ in $\dbl{D}$ or $\dbl{E}$ are handled by $F$ and $G$ respectively. Suppose, then, that $m \in M(x, y)$ and consider $m \odot \id_{y}$ without loss of generality. This then follows from \emph{coherence of naturality comparisons} for objects.

%%          N_{(x,y)}(m) = (\id_{{\kappa_{(x, y)}(m)}} \odot G_{y}) ; \kappa_{(\id, G\id_{y})} ; \kappa_{(\id, \id)}(M_{(x, y)}(m))
% 

    \item (\emph{Coherence with associators}) Consider $x \xproto{m} y \xproto{n} z \xproto{p} w$. There are three interesting cases:
          \begin{enumerate}
                  \item Suppose that $m$ and $n$ are in $\dbl{D}$ and $p \in M(z, w)$; this follows from \emph{coherence of naturality comparisons} for composable pairs, taking $m = (m, \id)$ and $n = (n, \id)$.
               %%   \id \odot \kappa_{(n, \id)}(p) ; \kappa_{(m, \id)}(M(n, \id)(p)) ; M_{(m, \id),(n,\id)}(p) = N_{(m, \id),(n, \id)}(\kappa_{(z,w)}(p)) ; \kappa_{(m \odot n, \id)}(\id_{p}).
                  \item If $m \in M(x, y)$ and $n$ and $p$ are in $\dbl{E}$, then coherence with associators follows from \emph{coherence of naturality comparisons} for composable pairs, taking $m = (\id, n)$ and $n = (\id, p)$.
                  \item If $n \in M(y, z)$, then coherence with associators follows from \emph{functoriality of component cells}, taking $\gamma = (\id_{m}, \id_{p})$.
          \end{enumerate}
  \end{enumerate}
\end{proof}

Finally, we extend the collage construction to a 2-functor.

\begin{construction}[Collage construction on cells] \label{construction:collage-modification}
  Let $\kappa : M \To N(F_0, G_0)$ and $\lambda : M \To N(F_1, G_1)$ be 1-cells in $\TwModl$ between $M : \dbl{D} \twbimodto \dbl{E}$ and $N$, and let $\mu : \kappa \Rrightarrow \lambda$ be a 2-cell between them, namely a modification relative to $\alpha^{\co}\times \beta$ where $\alpha : F_0 \To F_1$ and $\beta : G_0 \To G_1$. We construct a tight transformation $\collage{\mu} : \collage{\kappa} \To \collage{\lambda}$ as follows:
  \begin{itemize}
          \item For each object $x \in \collage{M}$, we define $\collage{\mu}_{x}$ to be $\alpha_x$ if $x \in \dbl{D}$ and $\beta_{x}$ if $x \in \dbl{E}$.
          \item For each loose morphism $m : x \proto y$ in $\collage{M}$, if $m$ is in $\dbl{D}$ or $\dbl{E}$ we define $\collage{\mu}(m)$ to be $\alpha_m$ or $\beta_m$ respectively, and if $x \in \dbl{D}$ and $y \in \dbl{E}$ we define $\collage{\mu}(m)$ to be $\mu_{(x,y)}(m)$, that is, the image of $m$ under the transformation
      \begin{equation*}
        % https://q.uiver.app/#q=WzAsNCxbMCwwLCJNKHgseSkiXSxbMCwxLCJOKEZfMCB4LCBHXzAgeSkiXSxbMSwwLCJNKHgseSkiXSxbMSwxLCJOKEZfMSB4LCBHXzEgeSkiXSxbMCwxLCJcXGthcHBhX3soeCx5KX0iLDJdLFswLDIsIiIsMCx7ImxldmVsIjoyLCJzdHlsZSI6eyJib2R5Ijp7Im5hbWUiOiJiYXJyZWQifSwiaGVhZCI6eyJuYW1lIjoibm9uZSJ9fX1dLFsyLDMsIlxcbnVfeyh4LHkpfSJdLFsxLDMsIk4oXFxhbHBoYV94LCBcXGJldGFfeSkiLDIseyJzdHlsZSI6eyJib2R5Ijp7Im5hbWUiOiJiYXJyZWQifX19XSxbNSw3LCJcXG11X3soeCx5KX0iLDEseyJzaG9ydGVuIjp7InNvdXJjZSI6MjAsInRhcmdldCI6MjB9LCJzdHlsZSI6eyJib2R5Ijp7Im5hbWUiOiJub25lIn0sImhlYWQiOnsibmFtZSI6Im5vbmUifX19XV0=
        \begin{tikzcd}
          {M(x,y)} & {M(x,y)} \\
          {N(F_0 x, G_0 y)} & {N(F_1 x, G_1 y)}
          \arrow[""{name=0, anchor=center, inner sep=0}, "\shortmid"{marking}, equals, from=1-1, to=1-2]
          \arrow["{\kappa_{(x,y)}}"', from=1-1, to=2-1]
          \arrow["{\lambda_{(x,y)}}", from=1-2, to=2-2]
          \arrow[""{name=1, anchor=center, inner sep=0}, "{N(\alpha_x, \beta_y)}"'{inner sep=.8ex}, "\shortmid"{marking}, from=2-1, to=2-2]
          \arrow["{\mu_{(x,y)}}"{description}, draw=none, from=0, to=1]
        \end{tikzcd}. \qedhere
      \end{equation*}
  \end{itemize}
\end{construction}

\begin{lemma}
 \cref{construction:collage-modification} gives a well defined (tight) natural transformation between lax double functors.
\end{lemma}
\begin{proof}
  We check the conditions of \cite[\S~3.5.4]{grandis2019}:
  \begin{enumerate}
    \item (\emph{Naturality}) Consider a cell $\stdInlineCell{\gamma}$ in $\collage{M}$. If $\gamma$ is in $\dbl{D}$ or $\dbl{E}$, naturality holds by the corresponding laws for $\alpha$ or $\beta$. Now suppose that $\gamma \in M(f, g)(m, n)$, seeking to show that
          \[\collage{\kappa}(\gamma) ; \collage{\mu}(n) = \collage{\mu}(m) ; \collage{\lambda}(\gamma)\]
    Unfolding definitions, we have $\kappa_{(f, g)}(\gamma) ; \mu_{(w, z)}(n) = \mu_{(x, y)}(m) ; \lambda_{(f, g)}(\gamma)$; this follows by tight-to-loose equivariance of $\mu$.
    \item (\emph{Coherence with loose identities}) Since identities must be either in $\dbl{D}$ or $\dbl{E}$, this follows by naturality of identity comparisons in $\alpha$ or $\beta$ respectively.
    
    \item (\emph{Coherence with loose composition}) Suppose that $x \xproto{m} y \xproto{n} z$ in $\collage{M}$; there are two interesting cases to check which will be symmetric, so we may suppose that $m \in M(x, y)$ and $n \in \dbl{E}$. We need to show that
          \[\collage{\kappa}_{m,n} ; \collage{\mu}(m \odot n) = (\collage{\mu}(m) \odot \collage{\mu}(n)) ; \collage{\lambda}_{m, n}.\]
          Expanding definitions this becomes
          \begin{align*}
          &N(F_{0,x}, \id_{G_0n})(\kappa_{(x,y)}(m));\kappa_{(\id, n)}(m) ; \mu_{(x,z)}(m \odot n)  \\ =\, &(\mu_{(x, y)}(m) \odot \mu_{(y,z)}(n)) ;  N(F_{1,x}, \id_{G_1n})(\lambda_{(x,y)}(m));\lambda_{(\id, n)}(m).
          \end{align*}
          We may take the component of the following equation at $m$:
          \[
\begin{tikzcd}
	& {M(x,y)} & {M(x,y)} & {M(x,y)} & {M(x,y)} \\
	& {N(F_0x, G_0y)} & {N(F_0x, G_0y)} & {M(x,z)} & {M(x,z)} \\
	& {N(F_0x,G_0z)} & {N(F_0x, G_0z) } & {N(F_0 x, G_0z)} & {N(F_1x, G_1z)} \\
	& {M(x,y)} & {M(x,y)} & {M(x,y)} & {M(x,y)} \\
	{=} & {N(F_0x, G_0y)} & {N(F_0x, G_0y)} & {N(F_1x, G_1y)} & {M(x,z)} \\
	& {N(F_0x,G_0z)} & {N(F_0x, G_0z) } & {N(F_1x, G_1z)} & {N(F_1x, G_1z)n} \\
	& {M(x,y)} & {M(x,y)} & {M(x,y)} & {M(x,y)} \\
	{=} & {N(F_0x, G_0y)} & {N(F_1x, G_1y)} & {N(F_1x, G_1y)} & {M(x,z)} \\
	& {N(F_0x,G_0z)} & {N(F_1x, G_1z)} & {N(F_1x, G_1z)} & {N(F_1x, G_1z)}
	\arrow["\shortmid"{marking}, equals, from=1-2, to=1-3]
	\arrow["{\kappa_{(x,y)}}"', from=1-2, to=2-2]
	\arrow[""{name=0, anchor=center, inner sep=0}, "\shortmid"{marking}, equals, from=1-3, to=1-4]
	\arrow["{\kappa_{(x,y)}}"', from=1-3, to=2-3]
	\arrow["\shortmid"{marking}, equals, from=1-4, to=1-5]
	\arrow["{M(\id, n)}", from=1-4, to=2-4]
	\arrow["{M(\id, n)}", from=1-5, to=2-5]
	\arrow[""{name=1, anchor=center, inner sep=0}, "\shortmid"{marking}, equals, from=2-2, to=2-3]
	\arrow["{N(\id_{F_0x}, G_0n)}"', from=2-2, to=3-2]
	\arrow["{N(F_0\id, G_0n)}", from=2-3, to=3-3]
	\arrow[""{name=2, anchor=center, inner sep=0}, "\shortmid"{marking}, equals, from=2-4, to=2-5]
	\arrow["{\kappa_{(x,z)}}", from=2-4, to=3-4]
	\arrow["{\lambda_{(x,z)}}", from=2-5, to=3-5]
	\arrow[""{name=3, anchor=center, inner sep=0}, "\shortmid"{marking}, equals, from=3-2, to=3-3]
	\arrow[""{name=4, anchor=center, inner sep=0}, "\shortmid"{marking}, equals, from=3-3, to=3-4]
	\arrow[""{name=5, anchor=center, inner sep=0}, "{N(\alpha_x, \beta_z)}"'{inner sep=.8ex}, "\shortmid"{marking}, from=3-4, to=3-5]
	\arrow["\shortmid"{marking}, equals, from=4-2, to=4-3]
	\arrow["{\kappa_{(x,y)}}"', from=4-2, to=5-2]
	\arrow[""{name=6, anchor=center, inner sep=0}, "\shortmid"{marking}, equals, from=4-3, to=4-4]
	\arrow["{\kappa_{(x,y)}}"', from=4-3, to=5-3]
	\arrow[""{name=7, anchor=center, inner sep=0}, "\shortmid"{marking}, equals, from=4-4, to=4-5]
	\arrow["{\lambda_{(x,y)}}", from=4-4, to=5-4]
	\arrow["{M(\id, n)}", from=4-5, to=5-5]
	\arrow[""{name=8, anchor=center, inner sep=0}, "\shortmid"{marking}, equals, from=5-2, to=5-3]
	\arrow["{N(\id_{F_0x}, G_0n)}"', from=5-2, to=6-2]
	\arrow[""{name=9, anchor=center, inner sep=0}, "\shortmid"{marking}, from=5-3, to=5-4]
	\arrow[from=5-3, to=6-3]
	\arrow["{N(F\id, Gn)}", from=5-4, to=6-4]
	\arrow["{\lambda_{(x,z)}}", from=5-5, to=6-5]
	\arrow[""{name=10, anchor=center, inner sep=0}, "\shortmid"{marking}, equals, from=6-2, to=6-3]
	\arrow[""{name=11, anchor=center, inner sep=0}, "{N(\alpha_x, \beta_z)}"'{inner sep=.8ex}, "\shortmid"{marking}, from=6-3, to=6-4]
	\arrow[""{name=12, anchor=center, inner sep=0}, "\shortmid"{marking}, equals, from=6-4, to=6-5]
	\arrow[""{name=13, anchor=center, inner sep=0}, "\shortmid"{marking}, equals, from=7-2, to=7-3]
	\arrow["{\kappa_{(x,y)}}"', from=7-2, to=8-2]
	\arrow["\shortmid"{marking}, equals, from=7-3, to=7-4]
	\arrow["{\lambda_{(x,y)}}", from=7-3, to=8-3]
	\arrow[""{name=14, anchor=center, inner sep=0}, "\shortmid"{marking}, equals, from=7-4, to=7-5]
	\arrow["{\lambda_{(x,y)}}"', from=7-4, to=8-4]
	\arrow["{M(\id, n)}", from=7-5, to=8-5]
	\arrow[""{name=15, anchor=center, inner sep=0}, "\shortmid"{marking}, from=8-2, to=8-3]
	\arrow["{N(\id_{F_0x}, G_0n)}"', from=8-2, to=9-2]
	\arrow[""{name=16, anchor=center, inner sep=0}, "\shortmid"{marking}, equals, from=8-3, to=8-4]
	\arrow[from=8-3, to=9-3]
	\arrow[from=8-4, to=9-4]
	\arrow["{\lambda_{(x,z)}}", from=8-5, to=9-5]
	\arrow[""{name=17, anchor=center, inner sep=0}, "{N(\alpha_x, \beta_z)}"{description}, "\shortmid"{marking}, from=9-2, to=9-3]
	\arrow[""{name=18, anchor=center, inner sep=0}, "\shortmid"{marking}, equals, from=9-3, to=9-4]
	\arrow[""{name=19, anchor=center, inner sep=0}, "\shortmid"{marking}, equals, from=9-4, to=9-5]
	\arrow["{\kappa_{(\id, n)}}"{description}, draw=none, from=0, to=4]
	\arrow["{N(F_{0,x}, \id_{G_0n})}"{description}, draw=none, from=1, to=3]
	\arrow["{\mu_{(x,z)}}"{description}, draw=none, from=2, to=5]
	\arrow["{\mu_{(x,y)}}"{description}, draw=none, from=6, to=9]
	\arrow["{\lambda_{(\id, n)}}"{description}, draw=none, from=7, to=12]
	\arrow["{N(F_{0,x}, \id_{G_0n})}"{description}, draw=none, from=8, to=10]
	\arrow["{N(\alpha_{\id}, \beta_{n})}"{description}, draw=none, from=9, to=11]
	\arrow["{\mu_{(x,y)}}"{description}, draw=none, from=13, to=15]
	\arrow["{\lambda(\id,n)}"{description}, draw=none, from=14, to=19]
	\arrow["{N(\alpha_{\id}, \beta_n)}"{description}, draw=none, from=15, to=17]
	\arrow["{N(F_{1,x}, \id_{G_1n})}"{description}, draw=none, from=16, to=18]
\end{tikzcd}
          \]
          The first equality follows by loose-to-tight equivariance of $\mu$, and the second by \emph{coherence with loose identities} of $\alpha$ as a tight transformation.
  \end{enumerate}
\end{proof}

\begin{lemma}[2-functorality of collage] \label{theorem:collage-2functor}
  The collage construction extends to a 2-functor $\collage{-} : \TwModl \to \lModl$.
\end{lemma}
\begin{proof}
  We check 2-functoriality for morphisms and cells.
  \begin{itemize}
    \item Suppose $\kappa : M \To N$ and $\lambda : N \To P$ are lax transformations of twisted bimodules, looking to compare $\collage{\kappa \cdot \lambda}$ with $\collage{\kappa} ; \collage{\lambda}$. We check componentwise:
\begin{enumerate}
         \item For $m : x \proto y$ in $\collage{M}$, the only nontrivial case to check is where $m \in M(x, y)$, in which case we have $\collage{\kappa \cdot \lambda}(m) = (\kappa \cdot \lambda)_{(x, y)}(m) = \lambda_{(x, y)}(\kappa_{(x, y)}(m))$, which is evidently $(\collage{\kappa}; \collage{\lambda})(m)$.
  \item For $\stdInlineCell{\gamma}$ in $\collage{M}$, the only nontrivial case to check is where $\gamma \in M(f, g)(m, n)$, in which case both sides reduce to $\lambda_{(f, g)}(\kappa_{(f, g)}(\gamma))$.
        \item The compositors of $\collage{\kappa \cdot \lambda}$ are given (in the generic nontrivial case $m \in M(x, y)$, $n : y \proto z$ in $\dbl{E}$) by applying the cell on the left of the equation below to $m \in M(x, y)$, while the compositors of the composite $\collage{\kappa} \cdot \collage{\lambda}$ are given by applying the cell on the right to $m \in M(x, y)$.

        \[
        \begin{array}{c}
        \begin{dblArray}{wc{4.4cm}wc{3.4cm}wc{3.4cm}wc{1.8cm}}
          1_{\kappa_{(x,y)}} & 1_{\kappa_{(x,y)}} & 1_{\kappa_{(x,y)}} & \Block{2-1}{\kappa_{(\id, n)}} \\
          1_{\lambda_{(F_0x,G_0y)}} & N(F_{0,x}, \id_{G_0n}) & \Block{2-1}{\lambda_{(F_0\id, G_0n)}} & \\
          P(F_{1,F_0x}, \id_{G_1G_0n}) & 1_{\lambda_{(F_0x,G_0z)}} & & 1_{\lambda_{(F_0x,G_0z)}}
        \end{dblArray}
        \\[2.5ex]
        =
        \\[2.5ex]
        \begin{dblArray}{wc{4.4cm}wc{3.4cm}wc{3.4cm}wc{1.8cm}}
          1_{\kappa_{(x,y)}} & 1_{\kappa_{(x,y)}} & 1_{\kappa_{(x,y)}} & \Block{2-1}{\kappa_{(\id, n)}} \\
          1_{\lambda_{(F_0x,G_0y)}} & \Block{2-1}{\lambda_{(\id, G_0n)}} & N(F_{0,x}, \id_{G_0n}) & \\
          P(F_{1,F_0x}, \id_{G_1G_0n}) & & 1_{\lambda_{(F_0x,G_0z)}} & 1_{\lambda_{(F_0x,G_0z)}}
        \end{dblArray}
        \end{array}
        \]

These two cells are equal by naturality with respect to cells for $\lambda$, applied to the unitor of $F_0$.
        \end{enumerate}

    \item  Suppose $\mu : \kappa \Rrightarrow \lambda$ and $\nu : \lambda \Rrightarrow \gamma$, looking to compare $\collage{\mu} \cdot \collage{\nu}$ with $\collage{\mu \cdot \nu}$. The only nontrivial case to consider is $m \in M(x, y)$, in which case $\collage{\mu \cdot \nu}(m) = (\mu \cdot \nu)_{(x, y)}(m) = \mu_{(x, y)}(m) \cdot \nu_{(x, y)}(m)$, and this is $(\collage{\mu} \cdot \collage{\nu})(m)$.
    \item Suppose $\mu : \kappa \Rrightarrow \lambda$ and $\nu : \gamma \Rrightarrow \delta$, looking to compare $\collage{\mu \ast \nu}$ with $\collage{\mu} \ast \collage{\nu}$. In the only nontrivial case of $m \in M(x, y)$, we have $\collage{\mu \ast \nu}(m) = (\mu \ast \nu)_{(x, y)}(m) = \gamma_{(x, y)}(\mu_{(x, y)}(m)) \cdot \nu_{(x, y)}(\gamma_{(x, y)}(m))$, and again this is $( \collage{\mu}  \ast \collage{\nu})(m)$.
  \end{itemize}
\end{proof}

We expect that the collage construction $\collage{-}$ is an inverse to the sections construction $\twsection{-}$. For the most part, this should be a matter of expanding definitions. For any twisted bimodule $M : \dbl{D} \twbimodto \dbl{E}$, we have $\barrelsrc{\collage{M}} \cong \dbl{D}$ and $\barreltgt{\collage{M}} \cong \dbl{E}$ by construction. Elements of $\twsection{\collage{M}}(x, y)$ are, by definition, those that map to the walking loose arrow under the labeling $\ell_{M}$, which are exactly the elements of $M(x, y)$ by the construction of the collage.

The tricky part is seeing how the comparison cells line up, since the associators and unitors of $\collage{M}$ are defined by one-sided isomorphs of the comparisons of $M$. To note the difficulty, recall that the comparison cells of the twisted Hom functor were defined by 5-ary associators and left-and-right unitors. Showing that, when round-tripping, these are recovered should amount to a general coherence principle for twisted bimodules.

We expect the rest of the data to correspond as well, modulo these subtle difficulties with coherence.

\begin{conjecture}\label{conjecture:sec.collage}
The 2-functor $\twsection{-} : \lModl \to \TwModl$ is an equivalence of 2-categories, with inverse given by the \emph{collage} construction (\cref{construction:collage-twisted-bimod}).
\end{conjecture}

\section{Opfibrations, algebras, and twisted functors}

The object of this section, having defined transformations and modifications of
twisted functors, is to begin the 2-categorical study of, on the one hand, the
relationship between twisted functors and loosely discrete opfibrations; and on
the other hand, the relationship between loosely discrete opfibrations and
algebras of pull-push monads, hence also between twisted functors and algebras.
We will show that every cloven loosely discrete opfibration canonically gives
rise to a twisted functor on the base double category
(\cref{subsection:pseudo-inverse-construction}). This construction is
essentially inverse to the elements construction of a twisted functor, so we
think of cloven loosely discrete opfibrations as corresponding to twisted
functors. Formally, we construct the components, in both directions, of an
equivalence between the two constructions
(\cref{subsection:completing-elements-correspondence}). This is the object-part
of a representation theorem for loosely discrete opfibrations, following the
pattern established for ordinary fibrations, (discrete) 2-fibrations, and double
fibrations, as discussed in the introduction to the paper. However, we leave the
complete 2-categorical development to a future study.

\subsection{Pseudo-inverse construction}
\label{subsection:pseudo-inverse-construction}

In this subsection, we associate to each cloven loosely discrete opfibration
$P\colon \dbl{E}\to \dbl{B}$ a \emph{unitary} profunctor-valued twisted functor
$\Tw(P)\colon\dbl{B}\twistto\Prof$ on the base double category. So, throughout fix a
loosely discrete opfibration $P\colon \dbl{E}\to \dbl{B}$ equipped with a cleavage
$\sigma: \dbl{E}_0 \times_{\dbl{B}_0} \dbl{B}_1 \to \dbl{E}_1$. Throughout this 
subsection, we shall use the following notation. 

\begin{definition}[Fibers] \label{definition:fiber_over_an_object}
  The \textbf{fiber} above an object $x\in\dbl{B}$ is the category $\dbl{E}_x$
  consisting of objects $a\in\dbl{E}$ such that $Pa=x$ and arrows $u\colon a\to b$ 
  of $\dbl{E}$ such that $Pu=1_x$.
\end{definition}

Indeed this notation involves an abuse inasmuch as it conflicts with that of
\cref{rmk:on-loose-discreteness}. However, that immediately above denotes 
essentially the underlying 1-category of the pullback definition in the remark. 
Additionally, the double-categorical fibers in the loosely discrete case are just 
that: loosely discrete. That is, the loose structure is in this sense trivial, 
and so to view these fibers as 1-categories is not to lose essential information.

\begin{construction}[Profibers] \label{construction:pro_fibers}
  Fix an arrow $f\colon x\to y$ of $\dbl{B}$. Define a profunctor by
    \begin{equation*}
      \dbl{E}_f\colon\dbl{E}_x^{\op}\times\dbl{E}_y\to\Set,
      \qquad
      (a,b) \mapsto \lbrace u\colon a\to b\mid Pu=f\rbrace.
    \end{equation*}
  This is functorial since the arrows of $\dbl{E}_x$ and $\dbl{E}_y$ are over 
  identity arrows on $x$ and $y$ and since $P$ is strictly functorial on arrows 
  of $\dbl{E}$. Given a composable arrow $g\colon y\to z$ of $\dbl{B}$, there is an 
  associated comparison cell  
    \begin{equation*}
      % https://q.uiver.app/#q=WzAsNSxbMCwwLCJcXGRibHtFfV94Il0sWzEsMCwiXFxkYmx7RX1feSJdLFsyLDAsIlxcZGJse0V9X3oiXSxbMiwxLCJcXGRibHtFfV96Il0sWzAsMSwiXFxkYmx7RX1feCJdLFswLDEsIlxcZGJse0V9X2YiLDAseyJzdHlsZSI6eyJib2R5Ijp7Im5hbWUiOiJiYXJyZWQifX19XSxbMSwyLCJcXGRibHtFfV9nIiwwLHsic3R5bGUiOnsiYm9keSI6eyJuYW1lIjoiYmFycmVkIn19fV0sWzIsMywiIiwwLHsibGV2ZWwiOjIsInN0eWxlIjp7ImhlYWQiOnsibmFtZSI6Im5vbmUifX19XSxbMCw0LCIiLDIseyJsZXZlbCI6Miwic3R5bGUiOnsiaGVhZCI6eyJuYW1lIjoibm9uZSJ9fX1dLFs0LDMsIlxcZGJse0V9X3tnZn0iLDIseyJzdHlsZSI6eyJib2R5Ijp7Im5hbWUiOiJiYXJyZWQifX19XSxbMSw5LCJcXGRibHtFfV97ZixnfSIsMSx7ImxhYmVsX3Bvc2l0aW9uIjo0MCwic2hvcnRlbiI6eyJ0YXJnZXQiOjIwfSwic3R5bGUiOnsiYm9keSI6eyJuYW1lIjoibm9uZSJ9LCJoZWFkIjp7Im5hbWUiOiJub25lIn19fV1d
      \begin{tikzcd}
        {\dbl{E}_x} & {\dbl{E}_y} & {\dbl{E}_z} \\
        {\dbl{E}_x} && {\dbl{E}_z}
        \arrow["{\dbl{E}_f}", "\shortmid"{marking}, from=1-1, to=1-2]
        \arrow[equals, from=1-1, to=2-1]
        \arrow["{\dbl{E}_g}", "\shortmid"{marking}, from=1-2, to=1-3]
        \arrow[equals, from=1-3, to=2-3]
        \arrow[""{name=0, anchor=center, inner sep=0}, "{\dbl{E}_{gf}}"', "\shortmid"{marking}, from=2-1, to=2-3]
        \arrow["{\dbl{E}_{f,g}}"{description, pos=0.4}, draw=none, from=1-2, to=0]
      \end{tikzcd}
    \end{equation*}
  given in the following way. The profunctor composite $\dbl{E}_f\odot
  \dbl{E}_g$ has the usual formula
    \begin{equation*}
      (\dbl{E}_f\odot \dbl{E}_g)(a,c) = \left(
        \sum_{b\in\dbl{E}_y}[\dbl{E}_f(a,b) \times \dbl{E}_g(b,c)]
      \right)/\approx,
    \end{equation*}
  where the relation $\approx$ is given by 
    \begin{equation*}
      (u,v)\approx (u',v') \text{ if, and only if, } gu=u' \text{ and } v=v'g
      \text{ for some } g\colon b\to b'.
    \end{equation*}
  Thus, the component map 
    \begin{equation*}
      (\dbl{E}_{f,g})_{a,c}\colon (\dbl{E}_f\odot \dbl{E}_g)(a,c) \to
      \dbl{E}_{gf}(a,c)
    \end{equation*} 
  is evidently given by tight composition $(u,v)\mapsto vu$ and it is easy to see
  that this is well-defined. Each such $\dbl{E}_{f,g}$ is natural since 
  composition of ordinary arrows in $\dbl{E}$ is strictly associative. Note that 
  such a comparison cell is not invertible in general. As for unitor 
  comparisons, in the special case where $f$ is $1_x\colon x\to x$ in $\dbl{B}$, 
  the associated profunctor $\dbl{E}_{1_x}$ is in fact the hom-functor 
  $\dbl{E}_x(-,=)$. Accordingly, the resulting twisted functor will be 
  \emph{unitary}.
\end{construction}

\begin{construction}[Transition functors]
  \label{construction:transition_functors}
  Fix a proarrow $m\colon x\proto y$ of $\dbl{B}$. Define a functor 
    \begin{equation*}
      m_!\colon \dbl{E}_x\to\dbl{E}_y  
    \end{equation*}
  by, for each $a\in\dbl{E}_x$, taking $m_!(a)$ to be the target of the proarrow
  $\sigma(a,m)$ in $\dbl{E}_1$ above $m$:
    \begin{equation*}
      (\sigma(a,m): a \proto m_!(a))
      \qquad\xmapsto{P}\qquad
      (m: x \proto y).
    \end{equation*}
  Think of this as the action of $m$ on $a$. Define the arrow assignment by
  taking an arrow $u\colon a\to b$ of $\dbl{E}_x$ to $m_!(u)$, the target of the cell
  $\sigma(u,1_m)$ over the tight identity on $m$:
    \begin{equation*}
      % https://q.uiver.app/#q=WzAsNCxbMCwwLCJhIl0sWzAsMSwiYiJdLFsxLDAsIm1fIShhKSJdLFsxLDEsIm1fIShiKSJdLFswLDEsInUiLDJdLFswLDIsIlxcc2lnbWEoYSxtKSIsMCx7InN0eWxlIjp7ImJvZHkiOnsibmFtZSI6ImJhcnJlZCJ9fX1dLFsxLDMsIlxcc2lnbWEoYixtKSIsMix7InN0eWxlIjp7ImJvZHkiOnsibmFtZSI6ImJhcnJlZCJ9fX1dLFsyLDMsIm1fISh1KSIsMCx7InN0eWxlIjp7ImJvZHkiOnsibmFtZSI6ImRhc2hlZCJ9fX1dLFs1LDYsIlxcc2lnbWEodSwgMV9tKSIsMSx7InNob3J0ZW4iOnsic291cmNlIjoyMCwidGFyZ2V0IjoyMH0sInN0eWxlIjp7ImJvZHkiOnsibmFtZSI6Im5vbmUifSwiaGVhZCI6eyJuYW1lIjoibm9uZSJ9fX1dXQ==
      \begin{tikzcd}
        a & {m_!(a)} \\
        b & {m_!(b)}
        \arrow[""{name=0, anchor=center, inner sep=0}, "{\sigma(a,m)}"{inner sep=.8ex}, "\shortmid"{marking}, from=1-1, to=1-2]
        \arrow["u"', from=1-1, to=2-1]
        \arrow["{m_!(u)}", dashed, from=1-2, to=2-2]
        \arrow[""{name=1, anchor=center, inner sep=0}, "{\sigma(b,m)}"'{inner sep=.8ex}, "\shortmid"{marking}, from=2-1, to=2-2]
        \arrow["{\sigma(u, 1_m)}"{description}, draw=none, from=0, to=1]
      \end{tikzcd}
      \qquad\xmapsto{P}\qquad
      % https://q.uiver.app/#q=WzAsNCxbMCwwLCJ4Il0sWzAsMSwieCJdLFsxLDAsInkiXSxbMSwxLCJ5Il0sWzAsMSwiMV94IiwyXSxbMCwyLCJtIiwwLHsic3R5bGUiOnsiYm9keSI6eyJuYW1lIjoiYmFycmVkIn19fV0sWzEsMywibSIsMix7InN0eWxlIjp7ImJvZHkiOnsibmFtZSI6ImJhcnJlZCJ9fX1dLFsyLDMsIjFfeSJdLFs1LDYsIjFfbSIsMSx7InNob3J0ZW4iOnsic291cmNlIjoyMCwidGFyZ2V0IjoyMH0sInN0eWxlIjp7ImJvZHkiOnsibmFtZSI6Im5vbmUifSwiaGVhZCI6eyJuYW1lIjoibm9uZSJ9fX1dXQ==
      \begin{tikzcd}
        x & y \\
        x & y
        \arrow[""{name=0, anchor=center, inner sep=0}, "m", "\shortmid"{marking}, from=1-1, to=1-2]
        \arrow["{1_x}"', from=1-1, to=2-1]
        \arrow["{1_y}", from=1-2, to=2-2]
        \arrow[""{name=1, anchor=center, inner sep=0}, "m"', "\shortmid"{marking}, from=2-1, to=2-2]
        \arrow["{1_m}"{description}, draw=none, from=0, to=1]
      \end{tikzcd}.
    \end{equation*} 
  This assignment is functorial because $\sigma$ is a functor. Suppose now that
  $m\colon x\proto y$ and $n\colon y\proto z$ are composable proarrows of 
  $\dbl{B}$. Define a comparison
    \begin{equation*}
      % https://q.uiver.app/#q=WzAsNSxbMCwwLCJcXGRibHtFfV94Il0sWzAsMSwiXFxkYmx7RX1feSJdLFswLDIsIlxcZGJse0V9X3oiXSxbMSwyLCJcXGRibHtFfV96Il0sWzEsMCwiXFxkYmx7RX1feCJdLFswLDEsIm1fISIsMl0sWzEsMiwibl8hIiwyXSxbMiwzLCJcXGlkIiwyLHsic3R5bGUiOnsiYm9keSI6eyJuYW1lIjoiYmFycmVkIn19fV0sWzAsNCwiXFxpZCIsMCx7InN0eWxlIjp7ImJvZHkiOnsibmFtZSI6ImJhcnJlZCJ9fX1dLFs0LDMsIihtXFxvZG90IG4pXyEiXSxbOCw3LCJcXGRibHtFfV57bSxufSIsMSx7InNob3J0ZW4iOnsic291cmNlIjoyMCwidGFyZ2V0IjoyMH0sInN0eWxlIjp7ImJvZHkiOnsibmFtZSI6Im5vbmUifSwiaGVhZCI6eyJuYW1lIjoibm9uZSJ9fX1dXQ==
      \begin{tikzcd}
        {\dbl{E}_x} & {\dbl{E}_x} \\
        {\dbl{E}_y} \\
        {\dbl{E}_z} & {\dbl{E}_z}
        \arrow[""{name=0, anchor=center, inner sep=0}, "\id"{inner sep=.8ex}, "\shortmid"{marking}, from=1-1, to=1-2]
        \arrow["{m_!}"', from=1-1, to=2-1]
        \arrow["{(m\odot n)_!}", from=1-2, to=3-2]
        \arrow["{n_!}"', from=2-1, to=3-1]
        \arrow[""{name=1, anchor=center, inner sep=0}, "\id"'{inner sep=.8ex}, "\shortmid"{marking}, from=3-1, to=3-2]
        \arrow["{\dbl{E}^{m,n}}"{description}, draw=none, from=0, to=1]
      \end{tikzcd}
      \qquad\leftrightsquigarrow\qquad
      % https://q.uiver.app/#q=WzAsMyxbMCwwLCJcXGRibHtFfV94XntcXG9wfVxcdGltZXMgXFxkYmx7RX1feCJdLFswLDIsIlxcZGJse0V9X3pee1xcb3B9XFx0aW1lcyBcXGRibHtFfV96Il0sWzIsMSwiXFxTZXQiXSxbMCwxLCIobl8hbV8hKV5cXG9wXFx0aW1lcyAobVxcb2RvdCBuKV8hIiwyXSxbMSwyLCJcXGRibHtFfV96KC0sPSkiLDJdLFswLDIsIlxcZGJse0V9X3goLSw9KSJdLFs1LDQsIlxcZGJse0V9XnttLG59IiwyLHsib2Zmc2V0Ijo1LCJzaG9ydGVuIjp7InNvdXJjZSI6MzAsInRhcmdldCI6MzB9fV1d
      \begin{tikzcd}
        {\dbl{E}_x^{\op}\times \dbl{E}_x} && \\
        && \Set \\
        {\dbl{E}_z^{\op}\times \dbl{E}_z}
        \arrow[""{name=0, anchor=center, inner sep=0}, "{\dbl{E}_x(-,=)}", from=1-1, to=2-3]
        \arrow["{(n_!m_!)^\op\times (m\odot n)_!}"', from=1-1, to=3-1]
        \arrow[""{name=1, anchor=center, inner sep=0}, "{\dbl{E}_z(-,=)}"', from=3-1, to=2-3]
        \arrow["{\dbl{E}^{m,n}}"', shift right=5, between={0.3}{0.7}, Rightarrow, from=0, to=1]
      \end{tikzcd}
    \end{equation*}
  whose component at an arrow $u\colon a\to b$ of $\dbl{E}_x$ is taken to be 
  the composite of dashed arrows in the diagram:
    \begin{equation*}
      % https://q.uiver.app/#q=WzAsNyxbMCwwLCJhIl0sWzEsMCwibV8hKGEpIl0sWzIsMCwibl8hbV8hKGEpIl0sWzAsMSwiYSJdLFsyLDEsIihtXFxvZG90IG4pXyEoYSkiXSxbMCwyLCJiIl0sWzIsMiwiKG1cXG9kb3QgbilfIShiKSJdLFswLDEsIlxcc2lnbWEoYSxtKSIsMCx7InN0eWxlIjp7ImJvZHkiOnsibmFtZSI6ImJhcnJlZCJ9fX1dLFsxLDIsIlxcc2lnbWEobV8hKGEpLG4pIiwwLHsic3R5bGUiOnsiYm9keSI6eyJuYW1lIjoiYmFycmVkIn19fV0sWzMsNCwiXFxzaWdtYShhLG1cXG9kb3QgbikiLDAseyJzdHlsZSI6eyJib2R5Ijp7Im5hbWUiOiJiYXJyZWQifX19XSxbMiw0LCJcXG11X3thLG0sbn0iLDAseyJzdHlsZSI6eyJib2R5Ijp7Im5hbWUiOiJkYXNoZWQifX19XSxbMyw1LCJ1IiwyXSxbNSw2LCJcXHNpZ21hKGIsbVxcb2RvdCBuKSIsMix7InN0eWxlIjp7ImJvZHkiOnsibmFtZSI6ImJhcnJlZCJ9fX1dLFs0LDYsIihtXFxvZG90IG4pXyEodSkiLDAseyJzdHlsZSI6eyJib2R5Ijp7Im5hbWUiOiJkYXNoZWQifX19XSxbOSwxMiwiXFxzaWdtYSh1LDFfe21cXG9kb3Qgbn0pIiwxLHsic2hvcnRlbiI6eyJzb3VyY2UiOjIwLCJ0YXJnZXQiOjIwfSwic3R5bGUiOnsiYm9keSI6eyJuYW1lIjoibm9uZSJ9LCJoZWFkIjp7Im5hbWUiOiJub25lIn19fV1d
      \begin{tikzcd}
        a & {m_!(a)} & {n_!m_!(a)} \\
        a && {(m\odot n)_!(a)} \\
        b && {(m\odot n)_!(b)}
        \arrow["{\sigma(a,m)}"{inner sep=.8ex}, "\shortmid"{marking}, from=1-1, to=1-2]
        \arrow["{\sigma(m_!(a),n)}"{inner sep=.8ex}, "\shortmid"{marking}, from=1-2, to=1-3]
        \arrow["{\mu_{a,m,n}}", dashed, from=1-3, to=2-3]
        \arrow[""{name=0, anchor=center, inner sep=0}, "{\sigma(a,m\odot n)}"{inner sep=.8ex}, "\shortmid"{marking}, from=2-1, to=2-3]
        \arrow["u"', from=2-1, to=3-1]
        \arrow["{(m\odot n)_!(u)}", dashed, from=2-3, to=3-3]
        \arrow[""{name=1, anchor=center, inner sep=0}, "{\sigma(b,m\odot n)}"'{inner sep=.8ex}, "\shortmid"{marking}, from=3-1, to=3-3]
        \arrow["{\sigma(u,1_{m\odot n})}"{description}, draw=none, from=0, to=1]
      \end{tikzcd}
    \end{equation*}
  where $\mu$ is given by the pseudo-algebra structure assumed for $\sigma$ as 
  described in \cref{remark:cleavage-pseudo-algebra-relations} and the bottom 
  arrow is just the target of $\sigma(u,1_{m\odot n})$. Note that there is no 
  cell coming with $\mu$. In any case, such a choice is suitably natural and the 
  components at identity arrows $1_a\colon a\to a$ are thus just the invertible 
  arrow $\mu_{a,m,n}$. Thus, the resulting cell $\dbl{E}_{m,n}$ is loosely 
  invertible. These cells make the assignment $m\mapsto m_!$ pseudo-functorial in 
  the sense that there is an equality
    \begin{equation*}
      % https://q.uiver.app/#q=WzAsOSxbMCwwLCJcXGRibHtFfV94Il0sWzAsMSwiXFxkYmx7RX1feSJdLFswLDIsIlxcZGJse0V9X3oiXSxbMSwyLCJcXGRibHtFfV96Il0sWzEsMCwiXFxkYmx7RX1feCJdLFswLDMsIlxcZGJse0V9X3ciXSxbMSwzLCJcXGRibHtFfV93Il0sWzIsMCwiXFxkYmx7RX1feCJdLFsyLDMsIlxcZGJse0V9X3ciXSxbMCwxLCJtXyEiLDJdLFsxLDIsIm5fISIsMl0sWzIsMywiXFxpZCIsMCx7InN0eWxlIjp7ImJvZHkiOnsibmFtZSI6ImJhcnJlZCJ9fX1dLFswLDQsIlxcaWQiLDAseyJzdHlsZSI6eyJib2R5Ijp7Im5hbWUiOiJiYXJyZWQifX19XSxbNCwzLCIobVxcb2RvdCBuKV8hIl0sWzIsNSwicF8hIiwyXSxbNSw2LCJcXGlkIiwyLHsic3R5bGUiOnsiYm9keSI6eyJuYW1lIjoiYmFycmVkIn19fV0sWzMsNiwicF8hIl0sWzQsNywiXFxpZCIsMCx7InN0eWxlIjp7ImJvZHkiOnsibmFtZSI6ImJhcnJlZCJ9fX1dLFs3LDgsIigobVxcb2RvdCBuKSBcXG9kb3QgcClfISJdLFs2LDgsIlxcaWQiLDIseyJzdHlsZSI6eyJib2R5Ijp7Im5hbWUiOiJiYXJyZWQifX19XSxbMTcsMTksIlxcZGJse0V9XnttXFxvZG90IG4sIHB9IiwxLHsic2hvcnRlbiI6eyJzb3VyY2UiOjIwLCJ0YXJnZXQiOjIwfSwic3R5bGUiOnsiYm9keSI6eyJuYW1lIjoibm9uZSJ9LCJoZWFkIjp7Im5hbWUiOiJub25lIn19fV0sWzEyLDExLCJcXGRibHtFfV57bSxufSIsMSx7InNob3J0ZW4iOnsic291cmNlIjoyMCwidGFyZ2V0IjoyMH0sInN0eWxlIjp7ImJvZHkiOnsibmFtZSI6Im5vbmUifSwiaGVhZCI6eyJuYW1lIjoibm9uZSJ9fX1dLFsxMSwxNSwiXFxpZF97cF8hfSIsMSx7InNob3J0ZW4iOnsic291cmNlIjoyMCwidGFyZ2V0IjoyMH0sInN0eWxlIjp7ImJvZHkiOnsibmFtZSI6Im5vbmUifSwiaGVhZCI6eyJuYW1lIjoibm9uZSJ9fX1dXQ==
      \begin{tikzcd}
        {\dbl{E}_x} & {\dbl{E}_x} & {\dbl{E}_x} \\
        {\dbl{E}_y} \\
        {\dbl{E}_z} & {\dbl{E}_z} \\
        {\dbl{E}_w} & {\dbl{E}_w} & {\dbl{E}_w}
        \arrow[""{name=0, anchor=center, inner sep=0}, "\id"{inner sep=.8ex}, "\shortmid"{marking}, from=1-1, to=1-2]
        \arrow["{m_!}"', from=1-1, to=2-1]
        \arrow[""{name=1, anchor=center, inner sep=0}, "\id"{inner sep=.8ex}, "\shortmid"{marking}, from=1-2, to=1-3]
        \arrow["{(m\odot n)_!}", from=1-2, to=3-2]
        \arrow["{((m\odot n) \odot p)_!}", from=1-3, to=4-3]
        \arrow["{n_!}"', from=2-1, to=3-1]
        \arrow[""{name=2, anchor=center, inner sep=0}, "\id"{inner sep=.8ex}, "\shortmid"{marking}, from=3-1, to=3-2]
        \arrow["{p_!}"', from=3-1, to=4-1]
        \arrow["{p_!}", from=3-2, to=4-2]
        \arrow[""{name=3, anchor=center, inner sep=0}, "\id"'{inner sep=.8ex}, "\shortmid"{marking}, from=4-1, to=4-2]
        \arrow[""{name=4, anchor=center, inner sep=0}, "\id"'{inner sep=.8ex}, "\shortmid"{marking}, from=4-2, to=4-3]
        \arrow["{\dbl{E}^{m,n}}"{description}, draw=none, from=0, to=2]
        \arrow["{\dbl{E}^{m\odot n, p}}"{description}, draw=none, from=1, to=4]
        \arrow["{\id_{p_!}}"{description}, draw=none, from=2, to=3]
      \end{tikzcd}
      \qquad=\qquad
      % https://q.uiver.app/#q=WzAsOSxbMCwwLCJcXGRibHtFfV94Il0sWzAsMSwiXFxkYmx7RX1feSJdLFswLDIsIlxcZGJse0V9X3oiXSxbMSwwLCJcXGRibHtFfV94Il0sWzAsMywiXFxkYmx7RX1fdyJdLFsxLDMsIlxcZGJse0V9X3ciXSxbMiwwLCJcXGRibHtFfV94Il0sWzIsMywiXFxkYmx7RX1fdyJdLFsxLDEsIlxcZGJse0V9X3kiXSxbMCwxLCJtXyEiLDJdLFsxLDIsIm5fISIsMl0sWzAsMywiXFxpZCIsMCx7InN0eWxlIjp7ImJvZHkiOnsibmFtZSI6ImJhcnJlZCJ9fX1dLFsyLDQsInBfISIsMl0sWzQsNSwiXFxpZCIsMix7InN0eWxlIjp7ImJvZHkiOnsibmFtZSI6ImJhcnJlZCJ9fX1dLFszLDYsIlxcaWQiLDAseyJzdHlsZSI6eyJib2R5Ijp7Im5hbWUiOiJiYXJyZWQifX19XSxbNiw3LCIobVxcb2RvdCAoblxcb2RvdCBwKSlfISJdLFs1LDcsIlxcaWQiLDIseyJzdHlsZSI6eyJib2R5Ijp7Im5hbWUiOiJiYXJyZWQifX19XSxbMyw4LCJtXyEiXSxbMSw4LCJcXGlkIiwyLHsic3R5bGUiOnsiYm9keSI6eyJuYW1lIjoiYmFycmVkIn19fV0sWzgsNSwiKG5cXG9kb3QgcClfISJdLFsxNCwxNiwiXFxkYmx7RX1ee20sblxcb2RvdCBwfSIsMSx7InNob3J0ZW4iOnsic291cmNlIjoyMCwidGFyZ2V0IjoyMH0sInN0eWxlIjp7ImJvZHkiOnsibmFtZSI6Im5vbmUifSwiaGVhZCI6eyJuYW1lIjoibm9uZSJ9fX1dLFsxMSwxOCwiXFxpZF97bV8hfSIsMSx7InNob3J0ZW4iOnsic291cmNlIjoyMCwidGFyZ2V0IjoyMH0sInN0eWxlIjp7ImJvZHkiOnsibmFtZSI6Im5vbmUifSwiaGVhZCI6eyJuYW1lIjoibm9uZSJ9fX1dLFsxOCwxMywiXFxkYmx7RX1ee24scH0iLDEseyJzaG9ydGVuIjp7InNvdXJjZSI6MjAsInRhcmdldCI6MjB9LCJzdHlsZSI6eyJib2R5Ijp7Im5hbWUiOiJub25lIn0sImhlYWQiOnsibmFtZSI6Im5vbmUifX19XV0=
      \begin{tikzcd}
        {\dbl{E}_x} & {\dbl{E}_x} & {\dbl{E}_x} \\
        {\dbl{E}_y} & {\dbl{E}_y} \\
        {\dbl{E}_z} \\
        {\dbl{E}_w} & {\dbl{E}_w} & {\dbl{E}_w}
        \arrow[""{name=0, anchor=center, inner sep=0}, "\id"{inner sep=.8ex}, "\shortmid"{marking}, from=1-1, to=1-2]
        \arrow["{m_!}"', from=1-1, to=2-1]
        \arrow[""{name=1, anchor=center, inner sep=0}, "\id"{inner sep=.8ex}, "\shortmid"{marking}, from=1-2, to=1-3]
        \arrow["{m_!}", from=1-2, to=2-2]
        \arrow["{(m\odot (n\odot p))_!}", from=1-3, to=4-3]
        \arrow[""{name=2, anchor=center, inner sep=0}, "\id"'{inner sep=.8ex}, "\shortmid"{marking}, from=2-1, to=2-2]
        \arrow["{n_!}"', from=2-1, to=3-1]
        \arrow["{(n\odot p)_!}", from=2-2, to=4-2]
        \arrow["{p_!}"', from=3-1, to=4-1]
        \arrow[""{name=3, anchor=center, inner sep=0}, "\id"'{inner sep=.8ex}, "\shortmid"{marking}, from=4-1, to=4-2]
        \arrow[""{name=4, anchor=center, inner sep=0}, "\id"'{inner sep=.8ex}, "\shortmid"{marking}, from=4-2, to=4-3]
        \arrow["{\id_{m_!}}"{description}, draw=none, from=0, to=2]
        \arrow["{\dbl{E}^{m,n\odot p}}"{description}, draw=none, from=1, to=4]
        \arrow["{\dbl{E}^{n,p}}"{description}, draw=none, from=2, to=3]
      \end{tikzcd}
    \end{equation*}
  as can be seen by chasing a component through each side of the diagram and 
  inserting the suppressed associativity isomorphisms. The essential aspects of 
  the argument are the naturality of $\mu$ combined with the assumed 
  associativity condition spelled out in 
  \cref{remark:cleavage-pseudo-algebra-relations}. Unit comparisons 
    \begin{equation*}
      % https://q.uiver.app/#q=WzAsNCxbMCwwLCJcXGRibHtFfV94Il0sWzEsMCwiXFxkYmx7RX1feCJdLFsxLDEsIlxcZGJse0V9X3giXSxbMCwxLCJcXGRibHtFfV94Il0sWzAsMSwiXFxpZCIsMCx7InN0eWxlIjp7ImJvZHkiOnsibmFtZSI6ImJhcnJlZCJ9fX1dLFsxLDIsIlxcZGJse0V9X3tcXGlkX3h9Il0sWzAsMywiMSIsMl0sWzMsMiwiXFxpZCIsMix7InN0eWxlIjp7ImJvZHkiOnsibmFtZSI6ImJhcnJlZCJ9fX1dLFs0LDcsIlxcZGJse0V9XngiLDEseyJzaG9ydGVuIjp7InNvdXJjZSI6MjAsInRhcmdldCI6MjB9LCJzdHlsZSI6eyJib2R5Ijp7Im5hbWUiOiJub25lIn0sImhlYWQiOnsibmFtZSI6Im5vbmUifX19XV0=
      \begin{tikzcd}
        {\dbl{E}_x} & {\dbl{E}_x} \\
        {\dbl{E}_x} & {\dbl{E}_x}
        \arrow[""{name=0, anchor=center, inner sep=0}, "\id"{inner sep=.8ex}, "\shortmid"{marking}, from=1-1, to=1-2]
        \arrow["1"', from=1-1, to=2-1]
        \arrow["{\dbl{E}_{\id_x}}", from=1-2, to=2-2]
        \arrow[""{name=1, anchor=center, inner sep=0}, "\id"'{inner sep=.8ex}, "\shortmid"{marking}, from=2-1, to=2-2]
        \arrow["{\dbl{E}^x}"{description}, draw=none, from=0, to=1]
      \end{tikzcd}
    \end{equation*}
  are given in a similar fashion using an appropriate component of the unitor 
  $\kappa$ coming with the pseudo-algebra structure. That is, for a morphism 
  $u\colon a\to b$ in $\dbl{E}_x$, take the component to be the composite
    \begin{equation*}
      % https://q.uiver.app/#q=WzAsMyxbMCwwLCJhIl0sWzEsMCwiYiJdLFsyLDAsIihcXGlkX3gpXyEoYikiXSxbMCwxLCJ1Il0sWzEsMiwiXFxrYXBwYV9iIl1d
\begin{tikzcd}
	a & b & {(\id_x)_!(b)}
	\arrow["u", from=1-1, to=1-2]
	\arrow["{\kappa_b}", from=1-2, to=1-3]
\end{tikzcd}
    \end{equation*}
  This is suitably natural. There are two loose unitor conditions that must 
  be satisfied. These are derived straightforwardly by the naturality of 
  $\kappa$ and by the two unit conditions on $\kappa$ described at the end of 
  \cref{remark:cleavage-pseudo-algebra-relations}.
\end{construction}

\begin{construction}[Cell assignment]
  \label{construction:cell_assignment_pseudo_inverse}
  The goal is now to associate to each cell $\alpha$ of $\dbl{B}$ as on the left, a
  cell in $\Prof$ as on the right
    \begin{equation*}
      % https://q.uiver.app/#q=WzAsNCxbMCwwLCJ4Il0sWzEsMCwieSJdLFsxLDEsInciXSxbMCwxLCJ6Il0sWzAsMSwibSIsMCx7InN0eWxlIjp7ImJvZHkiOnsibmFtZSI6ImJhcnJlZCJ9fX1dLFsxLDIsImciXSxbMCwzLCJmIiwyXSxbMywyLCJuIiwyLHsic3R5bGUiOnsiYm9keSI6eyJuYW1lIjoiYmFycmVkIn19fV0sWzQsNywiXFxhbHBoYSIsMSx7InNob3J0ZW4iOnsic291cmNlIjoyMCwidGFyZ2V0IjoyMH0sInN0eWxlIjp7ImJvZHkiOnsibmFtZSI6Im5vbmUifSwiaGVhZCI6eyJuYW1lIjoibm9uZSJ9fX1dXQ==
      \begin{tikzcd}
        x & y \\
        z & w
        \arrow[""{name=0, anchor=center, inner sep=0}, "m", "\shortmid"{marking}, from=1-1, to=1-2]
        \arrow["f"', from=1-1, to=2-1]
        \arrow["g", from=1-2, to=2-2]
        \arrow[""{name=1, anchor=center, inner sep=0}, "n"', "\shortmid"{marking}, from=2-1, to=2-2]
        \arrow["\alpha"{description}, draw=none, from=0, to=1]
      \end{tikzcd} \qquad \leadsto \qquad
      % https://q.uiver.app/#q=WzAsNCxbMCwwLCJcXGRibHtFfV94Il0sWzEsMCwiXFxkYmx7RX1feiJdLFsxLDEsIlxcZGJse0V9X3ciXSxbMCwxLCJcXGRibHtFfV95Il0sWzEsMiwibl8hIl0sWzAsMywibV8hIiwyXSxbMywyLCJcXGRibHtFfV9nIiwyLHsic3R5bGUiOnsiYm9keSI6eyJuYW1lIjoiYmFycmVkIn19fV0sWzAsMSwiXFxkYmx7RX1fZiIsMCx7InN0eWxlIjp7ImJvZHkiOnsibmFtZSI6ImJhcnJlZCJ9fX1dLFs3LDYsIlxcYWxwaGFfISIsMSx7InNob3J0ZW4iOnsic291cmNlIjoyMCwidGFyZ2V0IjoyMH0sInN0eWxlIjp7ImJvZHkiOnsibmFtZSI6Im5vbmUifSwiaGVhZCI6eyJuYW1lIjoibm9uZSJ9fX1dXQ==
      \begin{tikzcd}
        {\dbl{E}_x} & {\dbl{E}_z} \\
        {\dbl{E}_y} & {\dbl{E}_w}
        \arrow[""{name=0, anchor=center, inner sep=0}, "{\dbl{E}_f}"{inner sep=.8ex}, "\shortmid"{marking}, from=1-1, to=1-2]
        \arrow["{m_!}"', from=1-1, to=2-1]
        \arrow["{n_!}", from=1-2, to=2-2]
        \arrow[""{name=1, anchor=center, inner sep=0}, "{\dbl{E}_g}"'{inner sep=.8ex}, "\shortmid"{marking}, from=2-1, to=2-2]
        \arrow["{\alpha_!}"{description}, draw=none, from=0, to=1]
      \end{tikzcd},
    \end{equation*}
  where the functors and profunctors are those of Constructions
  \ref{construction:transition_functors} and \ref{construction:pro_fibers}
  above. Define the component function of $\alpha_!$ at a pair $(a,b)$ in 
  $\dbl{E}_x^{\op}\times\dbl{E}_z$ to be the function on the left whose
  action on an arrow $u\colon a\to b$ is the dashed arrow on the right:
    \begin{equation*}
      % https://q.uiver.app/#q=WzAsNCxbMCwwLCJhIl0sWzEsMCwibV8hKGEpIl0sWzEsMSwibl8hKGIpIl0sWzAsMSwiYiJdLFswLDEsIlxcc2lnbWEobSxhKSIsMCx7InN0eWxlIjp7ImJvZHkiOnsibmFtZSI6ImJhcnJlZCJ9fX1dLFsxLDIsIihcXGFscGhhXyEpX3thLGJ9KHUpIiwwLHsic3R5bGUiOnsiYm9keSI6eyJuYW1lIjoiZGFzaGVkIn19fV0sWzAsMywidSIsMl0sWzMsMiwiXFxzaWdtYShuLGIpIiwyLHsic3R5bGUiOnsiYm9keSI6eyJuYW1lIjoiYmFycmVkIn19fV0sWzQsNywiXFxzaWdtYSh1LFxcYWxwaGEpIiwxLHsic2hvcnRlbiI6eyJzb3VyY2UiOjIwLCJ0YXJnZXQiOjIwfSwic3R5bGUiOnsiYm9keSI6eyJuYW1lIjoibm9uZSJ9LCJoZWFkIjp7Im5hbWUiOiJub25lIn19fV1d
      \begin{tikzcd}
        a & {m_!(a)} \\
        b & {n_!(b)}
        \arrow[""{name=0, anchor=center, inner sep=0}, "{\sigma(m,a)}"{inner sep=.8ex}, "\shortmid"{marking}, from=1-1, to=1-2]
        \arrow["u"', from=1-1, to=2-1]
        \arrow["{(\alpha_!)_{a,b}(u)}", dashed, from=1-2, to=2-2]
        \arrow[""{name=1, anchor=center, inner sep=0}, "{\sigma(n,b)}"'{inner sep=.8ex}, "\shortmid"{marking}, from=2-1, to=2-2]
        \arrow["{\sigma(u,\alpha)}"{description}, draw=none, from=0, to=1]
      \end{tikzcd}
      \qquad\text{given by}\qquad
      (u: a \to b) \quad\mapsto\quad
      % https://q.uiver.app/#q=WzAsNCxbMCwwLCJhIl0sWzEsMCwibV8hKGEpIl0sWzEsMSwibl8hKGIpIl0sWzAsMSwiYiJdLFswLDEsIlxcc2lnbWEobSxhKSIsMCx7InN0eWxlIjp7ImJvZHkiOnsibmFtZSI6ImJhcnJlZCJ9fX1dLFsxLDIsIihcXGFscGhhXyEpX3thLGJ9KHUpIiwwLHsic3R5bGUiOnsiYm9keSI6eyJuYW1lIjoiZGFzaGVkIn19fV0sWzAsMywidSIsMl0sWzMsMiwiXFxzaWdtYShuLGIpIiwyLHsic3R5bGUiOnsiYm9keSI6eyJuYW1lIjoiYmFycmVkIn19fV0sWzQsNywiXFxzaWdtYSh1LFxcYWxwaGEpIiwxLHsic2hvcnRlbiI6eyJzb3VyY2UiOjIwLCJ0YXJnZXQiOjIwfSwic3R5bGUiOnsiYm9keSI6eyJuYW1lIjoibm9uZSJ9LCJoZWFkIjp7Im5hbWUiOiJub25lIn19fV1d
      \begin{tikzcd}
        a & {m_!(a)} \\
        b & {n_!(b)}
        \arrow[""{name=0, anchor=center, inner sep=0}, "{\sigma(m,a)}"{inner sep=.8ex}, "\shortmid"{marking}, from=1-1, to=1-2]
        \arrow["u"', from=1-1, to=2-1]
        \arrow["{(\alpha_!)_{a,b}(u)}", dashed, from=1-2, to=2-2]
        \arrow[""{name=1, anchor=center, inner sep=0}, "{\sigma(n,b)}"'{inner sep=.8ex}, "\shortmid"{marking}, from=2-1, to=2-2]
        \arrow["{\sigma(u,\alpha)}"{description}, draw=none, from=0, to=1]
      \end{tikzcd}.
    \end{equation*}
  That is, each component acts by taking the target arrow of the cell 
  $\sigma(u,\alpha)$. These components are natural in $(a,b)$ owing to the 
  definitions and the fact that such lifts are unique.
\end{construction}

\begin{proposition}[Twisted functor from loosely discrete opfibration]
  \label{prop:twisted-functor-assoc-to-fibration}
  Given a cloven loosely discrete opfibration $P: \dbl{E} \to \dbl{B}$, the data
  of the previous constructions assemble into a \emph{unitary} twisted functor
  $\Tw(P): \dbl{B}\twistto\Prof$ with assignments
    \begin{equation*}
      \Tw(P)(x) \coloneqq \dbl{E}_x, \quad
      \Tw(P)(f) \coloneqq \dbl{E}_f, \quad
      \Tw(P)(m) \coloneqq m_!, \quad
      \Tw(P)(\alpha) \coloneqq \alpha_!
    \end{equation*}
  on objects, arrows, proarrows, and cells, respectively. The composition
  and identity comparisons are also defined in the constructions above.
\end{proposition}
\begin{proof}
  We shall verify the two naturality conditions of 
  \cref{defn:twisted-double-functor}. First for tightly composable cells as 
  at the left below, we should have an equality of composite cells as on the 
  right: 
    \begin{equation*}
      % https://q.uiver.app/#q=WzAsNixbMCwwLCJ4Il0sWzEsMCwieSJdLFsxLDEsInknIl0sWzAsMSwieCciXSxbMCwyLCJ4JyciXSxbMSwyLCJ5JyciXSxbMCwxLCJtIiwwLHsic3R5bGUiOnsiYm9keSI6eyJuYW1lIjoiYmFycmVkIn19fV0sWzEsMiwiZyJdLFswLDMsImYiLDJdLFszLDIsIm4iLDIseyJzdHlsZSI6eyJib2R5Ijp7Im5hbWUiOiJiYXJyZWQifX19XSxbMyw0LCJoIiwyXSxbNCw1LCJwIiwyLHsic3R5bGUiOnsiYm9keSI6eyJuYW1lIjoiYmFycmVkIn19fV0sWzIsNSwiayJdLFs5LDExLCJcXGJldGEiLDEseyJzaG9ydGVuIjp7InNvdXJjZSI6MjAsInRhcmdldCI6MjB9LCJzdHlsZSI6eyJib2R5Ijp7Im5hbWUiOiJub25lIn0sImhlYWQiOnsibmFtZSI6Im5vbmUifX19XSxbNiw5LCJcXGFscGhhIiwxLHsic2hvcnRlbiI6eyJzb3VyY2UiOjIwLCJ0YXJnZXQiOjIwfSwic3R5bGUiOnsiYm9keSI6eyJuYW1lIjoibm9uZSJ9LCJoZWFkIjp7Im5hbWUiOiJub25lIn19fV1d
      \begin{tikzcd}
        x & y \\
        {x'} & {y'} \\
        {x''} & {y''}
        \arrow[""{name=0, anchor=center, inner sep=0}, "m", "\shortmid"{marking}, from=1-1, to=1-2]
        \arrow["f"', from=1-1, to=2-1]
        \arrow["g", from=1-2, to=2-2]
        \arrow[""{name=1, anchor=center, inner sep=0}, "n"', "\shortmid"{marking}, from=2-1, to=2-2]
        \arrow["h"', from=2-1, to=3-1]
        \arrow["k", from=2-2, to=3-2]
        \arrow[""{name=2, anchor=center, inner sep=0}, "p"', "\shortmid"{marking}, from=3-1, to=3-2]
        \arrow["\alpha"{description}, draw=none, from=0, to=1]
        \arrow["\beta"{description}, draw=none, from=1, to=2]
      \end{tikzcd}
      \qquad\leadsto\qquad
      % https://q.uiver.app/#q=WzAsOCxbMCwwLCJcXGRibHtFfV94Il0sWzEsMCwiXFxkYmx7RX1fe3gnfSJdLFsxLDEsIlxcZGJse0V9X3t5J30iXSxbMCwxLCJcXGRibHtFfV95Il0sWzIsMCwiXFxkYmx7RX1fe3gnJ30iXSxbMiwxLCJcXGRibHtFfV97eScnfSJdLFswLDIsIlxcZGJse0V9X3kiXSxbMiwyLCJcXGRibHtFfV97eScnfSJdLFswLDEsIlxcZGJse0V9X2YiLDAseyJzdHlsZSI6eyJib2R5Ijp7Im5hbWUiOiJiYXJyZWQifX19XSxbMSwyLCJuXyEiLDFdLFswLDMsIm1fISIsMl0sWzMsMiwiXFxkYmx7RX1fZyIsMix7InN0eWxlIjp7ImJvZHkiOnsibmFtZSI6ImJhcnJlZCJ9fX1dLFsyLDUsIlxcZGJse0V9X2siLDIseyJzdHlsZSI6eyJib2R5Ijp7Im5hbWUiOiJiYXJyZWQifX19XSxbNCw1LCJwXyEiXSxbMyw2LCIiLDIseyJsZXZlbCI6Miwic3R5bGUiOnsiaGVhZCI6eyJuYW1lIjoibm9uZSJ9fX1dLFs2LDcsIlxcZGJse0V9X3trZ30iLDIseyJzdHlsZSI6eyJib2R5Ijp7Im5hbWUiOiJiYXJyZWQifX19XSxbNSw3LCIiLDAseyJsZXZlbCI6Miwic3R5bGUiOnsiaGVhZCI6eyJuYW1lIjoibm9uZSJ9fX1dLFsxLDQsIlxcZGJse0V9X2giLDAseyJzdHlsZSI6eyJib2R5Ijp7Im5hbWUiOiJiYXJyZWQifX19XSxbOCwxMSwiXFxhbHBoYV8hIiwxLHsic2hvcnRlbiI6eyJzb3VyY2UiOjIwLCJ0YXJnZXQiOjIwfSwic3R5bGUiOnsiYm9keSI6eyJuYW1lIjoibm9uZSJ9LCJoZWFkIjp7Im5hbWUiOiJub25lIn19fV0sWzIsMTUsIlxcZGJse0V9X3tnLGt9IiwxLHsibGFiZWxfcG9zaXRpb24iOjQwLCJzaG9ydGVuIjp7InRhcmdldCI6MjB9LCJzdHlsZSI6eyJib2R5Ijp7Im5hbWUiOiJub25lIn0sImhlYWQiOnsibmFtZSI6Im5vbmUifX19XSxbMTcsMTIsIlxcYmV0YV8hIiwxLHsic2hvcnRlbiI6eyJzb3VyY2UiOjIwLCJ0YXJnZXQiOjIwfSwic3R5bGUiOnsiYm9keSI6eyJuYW1lIjoibm9uZSJ9LCJoZWFkIjp7Im5hbWUiOiJub25lIn19fV1d
      \begin{tikzcd}
        {\dbl{E}_x} & {\dbl{E}_{x'}} & {\dbl{E}_{x''}} \\
        {\dbl{E}_y} & {\dbl{E}_{y'}} & {\dbl{E}_{y''}} \\
        {\dbl{E}_y} && {\dbl{E}_{y''}}
        \arrow[""{name=0, anchor=center, inner sep=0}, "{\dbl{E}_f}"{inner sep=.8ex}, "\shortmid"{marking}, from=1-1, to=1-2]
        \arrow["{m_!}"', from=1-1, to=2-1]
        \arrow[""{name=1, anchor=center, inner sep=0}, "{\dbl{E}_h}"{inner sep=.8ex}, "\shortmid"{marking}, from=1-2, to=1-3]
        \arrow["{n_!}"{description}, from=1-2, to=2-2]
        \arrow["{p_!}", from=1-3, to=2-3]
        \arrow[""{name=2, anchor=center, inner sep=0}, "{\dbl{E}_g}"'{inner sep=.8ex}, "\shortmid"{marking}, from=2-1, to=2-2]
        \arrow[equals, from=2-1, to=3-1]
        \arrow[""{name=3, anchor=center, inner sep=0}, "{\dbl{E}_k}"'{inner sep=.8ex}, "\shortmid"{marking}, from=2-2, to=2-3]
        \arrow[equals, from=2-3, to=3-3]
        \arrow[""{name=4, anchor=center, inner sep=0}, "{\dbl{E}_{kg}}"'{inner sep=.8ex}, "\shortmid"{marking}, from=3-1, to=3-3]
        \arrow["{\alpha_!}"{description}, draw=none, from=0, to=2]
        \arrow["{\beta_!}"{description}, draw=none, from=1, to=3]
        \arrow["{\dbl{E}_{g,k}}"{description, pos=0.4}, draw=none, from=2-2, to=4]
      \end{tikzcd}
      \;\; = \;\;
      % https://q.uiver.app/#q=WzAsNyxbMCwwLCJcXGRibHtFfV94Il0sWzEsMCwiXFxkYmx7RX1fe3gnfSJdLFswLDEsIlxcZGJse0V9X3kiXSxbMiwwLCJcXGRibHtFfV97eCcnfSJdLFsyLDEsIlxcZGJse0V9X3t5Jyd9Il0sWzAsMiwiXFxkYmx7RX1feSJdLFsyLDIsIlxcZGJse0V9X3t5Jyd9Il0sWzAsMSwiXFxkYmx7RX1fZiIsMCx7InN0eWxlIjp7ImJvZHkiOnsibmFtZSI6ImJhcnJlZCJ9fX1dLFswLDIsIiIsMix7ImxldmVsIjoyLCJzdHlsZSI6eyJoZWFkIjp7Im5hbWUiOiJub25lIn19fV0sWzEsMywiXFxkYmx7RX1fZyIsMCx7InN0eWxlIjp7ImJvZHkiOnsibmFtZSI6ImJhcnJlZCJ9fX1dLFszLDQsIiIsMCx7ImxldmVsIjoyLCJzdHlsZSI6eyJoZWFkIjp7Im5hbWUiOiJub25lIn19fV0sWzIsNSwibV8hIiwyXSxbNSw2LCJcXGRibHtFfV97a2d9IiwyLHsic3R5bGUiOnsiYm9keSI6eyJuYW1lIjoiYmFycmVkIn19fV0sWzQsNiwicF8hIl0sWzIsNCwiIiwxLHsic3R5bGUiOnsiYm9keSI6eyJuYW1lIjoiYmFycmVkIn19fV0sWzE0LDEyLCIoXFxiZXRhXFxhbHBoYSlfISIsMSx7InNob3J0ZW4iOnsic291cmNlIjoyMCwidGFyZ2V0IjoyMH0sInN0eWxlIjp7ImJvZHkiOnsibmFtZSI6Im5vbmUifSwiaGVhZCI6eyJuYW1lIjoibm9uZSJ9fX1dLFsxLDE0LCJcXGRibHtFfV97ZixofSIsMSx7ImxhYmVsX3Bvc2l0aW9uIjo0MCwic2hvcnRlbiI6eyJ0YXJnZXQiOjIwfSwic3R5bGUiOnsiYm9keSI6eyJuYW1lIjoibm9uZSJ9LCJoZWFkIjp7Im5hbWUiOiJub25lIn19fV1d
      \begin{tikzcd}
        {\dbl{E}_x} & {\dbl{E}_{x'}} & {\dbl{E}_{x''}} \\
        {\dbl{E}_y} && {\dbl{E}_{y''}} \\
        {\dbl{E}_y} && {\dbl{E}_{y''}}
        \arrow["{\dbl{E}_f}"{inner sep=.8ex}, "\shortmid"{marking}, from=1-1, to=1-2]
        \arrow[equals, from=1-1, to=2-1]
        \arrow["{\dbl{E}_g}"{inner sep=.8ex}, "\shortmid"{marking}, from=1-2, to=1-3]
        \arrow[equals, from=1-3, to=2-3]
        \arrow[""{name=0, anchor=center, inner sep=0}, "\shortmid"{marking}, from=2-1, to=2-3]
        \arrow["{m_!}"', from=2-1, to=3-1]
        \arrow["{p_!}", from=2-3, to=3-3]
        \arrow[""{name=1, anchor=center, inner sep=0}, "{\dbl{E}_{kg}}"'{inner sep=.8ex}, "\shortmid"{marking}, from=3-1, to=3-3]
        \arrow["{\dbl{E}_{f,h}}"{description, pos=0.4}, draw=none, from=1-2, to=0]
        \arrow["{(\beta\alpha)_!}"{description}, draw=none, from=0, to=1]
      \end{tikzcd}.
    \end{equation*}
  This is however easily seen to be the case. For chasing the action of an 
  arbitrary component through each diagram we arrive at a purported equality 
  between the dashed arrows on the right boundaries:
    \begin{equation*}
      % https://q.uiver.app/#q=WzAsNixbMCwwLCJhIl0sWzEsMCwibV8hKGEpIl0sWzEsMSwibl8hKGIpIl0sWzAsMSwiYiJdLFswLDIsImMiXSxbMSwyLCJwXyEoYykiXSxbMCwxLCJcXHNpZ21hKG0sYSkiLDAseyJzdHlsZSI6eyJib2R5Ijp7Im5hbWUiOiJiYXJyZWQifX19XSxbMSwyLCIiLDAseyJzdHlsZSI6eyJib2R5Ijp7Im5hbWUiOiJkYXNoZWQifX19XSxbMCwzLCJ1IiwyXSxbMyw0LCJ2IiwyXSxbNCw1LCJcXHNpZ21hKHAsYykiLDIseyJzdHlsZSI6eyJib2R5Ijp7Im5hbWUiOiJiYXJyZWQifX19XSxbMiw1LCIiLDAseyJzdHlsZSI6eyJib2R5Ijp7Im5hbWUiOiJkYXNoZWQifX19XSxbMywyLCJcXHNpZ21hKG4sYikiLDAseyJzdHlsZSI6eyJib2R5Ijp7Im5hbWUiOiJiYXJyZWQifX19XSxbMTIsMTAsIlxcc2lnbWEodixcXGJldGEpIiwxLHsic2hvcnRlbiI6eyJzb3VyY2UiOjIwLCJ0YXJnZXQiOjIwfSwic3R5bGUiOnsiYm9keSI6eyJuYW1lIjoibm9uZSJ9LCJoZWFkIjp7Im5hbWUiOiJub25lIn19fV0sWzYsMTIsIlxcc2lnbWEodSxcXGFscGhhKSIsMSx7ImxhYmVsX3Bvc2l0aW9uIjo0MCwic2hvcnRlbiI6eyJzb3VyY2UiOjIwLCJ0YXJnZXQiOjIwfSwic3R5bGUiOnsiYm9keSI6eyJuYW1lIjoibm9uZSJ9LCJoZWFkIjp7Im5hbWUiOiJub25lIn19fV1d
      \begin{tikzcd}
        a & {m_!(a)} \\
        b & {n_!(b)} \\
        c & {p_!(c)}
        \arrow[""{name=0, anchor=center, inner sep=0}, "{\sigma(m,a)}"{inner sep=.8ex}, "\shortmid"{marking}, from=1-1, to=1-2]
        \arrow["u"', from=1-1, to=2-1]
        \arrow[dashed, from=1-2, to=2-2]
        \arrow[""{name=1, anchor=center, inner sep=0}, "{\sigma(n,b)}"{inner sep=.8ex}, "\shortmid"{marking}, from=2-1, to=2-2]
        \arrow["v"', from=2-1, to=3-1]
        \arrow[dashed, from=2-2, to=3-2]
        \arrow[""{name=2, anchor=center, inner sep=0}, "{\sigma(p,c)}"'{inner sep=.8ex}, "\shortmid"{marking}, from=3-1, to=3-2]
        \arrow["{\sigma(u,\alpha)}"{description, pos=0.4}, draw=none, from=0, to=1]
        \arrow["{\sigma(v,\beta)}"{description}, draw=none, from=1, to=2]
      \end{tikzcd} \quad = \quad 
      % https://q.uiver.app/#q=WzAsNCxbMCwwLCJhIl0sWzEsMCwibV8hKGEpIl0sWzAsMSwiYyJdLFsxLDEsInBfIShjKSJdLFswLDEsIlxcc2lnbWEobSxhKSIsMCx7InN0eWxlIjp7ImJvZHkiOnsibmFtZSI6ImJhcnJlZCJ9fX1dLFsyLDMsIlxcc2lnbWEocCxjKSIsMix7InN0eWxlIjp7ImJvZHkiOnsibmFtZSI6ImJhcnJlZCJ9fX1dLFswLDIsInZ1IiwyXSxbMSwzLCIiLDAseyJzdHlsZSI6eyJib2R5Ijp7Im5hbWUiOiJkYXNoZWQifX19XSxbNCw1LCJcXHNpZ21hKHZ1LFxcYmV0YVxcYWxwaGEpIiwxLHsic2hvcnRlbiI6eyJzb3VyY2UiOjIwLCJ0YXJnZXQiOjIwfSwic3R5bGUiOnsiYm9keSI6eyJuYW1lIjoibm9uZSJ9LCJoZWFkIjp7Im5hbWUiOiJub25lIn19fV1d
      \begin{tikzcd}
        a & {m_!(a)} \\
        c & {p_!(c)}
        \arrow[""{name=0, anchor=center, inner sep=0}, "{\sigma(m,a)}"{inner sep=.8ex}, "\shortmid"{marking}, from=1-1, to=1-2]
        \arrow["vu"', from=1-1, to=2-1]
        \arrow[dashed, from=1-2, to=2-2]
        \arrow[""{name=1, anchor=center, inner sep=0}, "{\sigma(p,c)}"'{inner sep=.8ex}, "\shortmid"{marking}, from=2-1, to=2-2]
        \arrow["{\sigma(vu,\beta\alpha)}"{description}, draw=none, from=0, to=1]
      \end{tikzcd}.
    \end{equation*}
  But these must be equal owing to the fact that $\sigma$ is a functor. That is, 
  $\sigma$ preserves tight composition of cells and thus of targets in particular.
  Second, for loosely composable cells 
    \begin{equation*}
      % https://q.uiver.app/#q=WzAsNixbMCwwLCJ4Il0sWzEsMCwieSJdLFsyLDAsInoiXSxbMiwxLCJ6JyJdLFswLDEsIngnIl0sWzEsMSwieSciXSxbMCwxLCJtIiwwLHsic3R5bGUiOnsiYm9keSI6eyJuYW1lIjoiYmFycmVkIn19fV0sWzEsMiwibiIsMCx7InN0eWxlIjp7ImJvZHkiOnsibmFtZSI6ImJhcnJlZCJ9fX1dLFsyLDMsImgiXSxbMCw0LCJmIiwyXSxbNCw1LCJwIiwyLHsic3R5bGUiOnsiYm9keSI6eyJuYW1lIjoiYmFycmVkIn19fV0sWzUsMywicSIsMix7InN0eWxlIjp7ImJvZHkiOnsibmFtZSI6ImJhcnJlZCJ9fX1dLFsxLDUsImciLDFdLFs2LDEwLCJcXGFscGhhIiwxLHsic2hvcnRlbiI6eyJzb3VyY2UiOjIwLCJ0YXJnZXQiOjIwfSwic3R5bGUiOnsiYm9keSI6eyJuYW1lIjoibm9uZSJ9LCJoZWFkIjp7Im5hbWUiOiJub25lIn19fV0sWzcsMTEsIlxcYmV0YSIsMSx7InNob3J0ZW4iOnsic291cmNlIjoyMCwidGFyZ2V0IjoyMH0sInN0eWxlIjp7ImJvZHkiOnsibmFtZSI6Im5vbmUifSwiaGVhZCI6eyJuYW1lIjoibm9uZSJ9fX1dXQ==
      \begin{tikzcd}
        x & y & z \\
        {x'} & {y'} & {z'}
        \arrow[""{name=0, anchor=center, inner sep=0}, "m"{inner sep=.8ex}, "\shortmid"{marking}, from=1-1, to=1-2]
        \arrow["f"', from=1-1, to=2-1]
        \arrow[""{name=1, anchor=center, inner sep=0}, "n"{inner sep=.8ex}, "\shortmid"{marking}, from=1-2, to=1-3]
        \arrow["g"{description}, from=1-2, to=2-2]
        \arrow["h", from=1-3, to=2-3]
        \arrow[""{name=2, anchor=center, inner sep=0}, "p"'{inner sep=.8ex}, "\shortmid"{marking}, from=2-1, to=2-2]
        \arrow[""{name=3, anchor=center, inner sep=0}, "q"'{inner sep=.8ex}, "\shortmid"{marking}, from=2-2, to=2-3]
        \arrow["\alpha"{description}, draw=none, from=0, to=2]
        \arrow["\beta"{description}, draw=none, from=1, to=3]
      \end{tikzcd},
    \end{equation*}
  we should have a corresponding equality of cells:
    \begin{equation*}
      % https://q.uiver.app/#q=WzAsOCxbMCwwLCJcXGRibHtFfV94Il0sWzEsMCwiXFxkYmx7RX1fe3gnfSJdLFsxLDEsIlxcZGJse0V9X3t5J30iXSxbMCwxLCJcXGRibHtFfV95Il0sWzAsMiwiXFxkYmx7RX1feiJdLFsxLDIsIlxcZGJse0V9X3t6J30iXSxbMiwwLCJcXGRibHtFfV97eCd9Il0sWzIsMiwiXFxkYmx7RX1fe3onfSJdLFswLDEsIlxcZGJse0V9X2YiLDAseyJzdHlsZSI6eyJib2R5Ijp7Im5hbWUiOiJiYXJyZWQifX19XSxbMSwyLCJwXyEiXSxbMCwzLCJtXyEiLDJdLFszLDIsIiIsMCx7InN0eWxlIjp7ImJvZHkiOnsibmFtZSI6ImJhcnJlZCJ9fX1dLFszLDQsIm5fISIsMl0sWzQsNSwiXFxkYmx7RX1faCIsMix7InN0eWxlIjp7ImJvZHkiOnsibmFtZSI6ImJhcnJlZCJ9fX1dLFsyLDUsInFfISJdLFsxLDYsIiIsMCx7ImxldmVsIjoyLCJzdHlsZSI6eyJib2R5Ijp7Im5hbWUiOiJiYXJyZWQifSwiaGVhZCI6eyJuYW1lIjoibm9uZSJ9fX1dLFs2LDcsIihwXFxvZG90IHEpXyEiXSxbNSw3LCIiLDAseyJsZXZlbCI6Miwic3R5bGUiOnsiYm9keSI6eyJuYW1lIjoiYmFycmVkIn0sImhlYWQiOnsibmFtZSI6Im5vbmUifX19XSxbMTEsMTMsIlxcYmV0YV8hIiwxLHsic2hvcnRlbiI6eyJzb3VyY2UiOjIwLCJ0YXJnZXQiOjIwfSwic3R5bGUiOnsiYm9keSI6eyJuYW1lIjoibm9uZSJ9LCJoZWFkIjp7Im5hbWUiOiJub25lIn19fV0sWzgsMTEsIlxcYWxwaGFfISIsMSx7InNob3J0ZW4iOnsic291cmNlIjoyMCwidGFyZ2V0IjoyMH0sInN0eWxlIjp7ImJvZHkiOnsibmFtZSI6Im5vbmUifSwiaGVhZCI6eyJuYW1lIjoibm9uZSJ9fX1dLFsxNSwxNywiXFxkYmx7RX1ee3AscX0iLDEseyJzaG9ydGVuIjp7InNvdXJjZSI6MjAsInRhcmdldCI6MjB9LCJzdHlsZSI6eyJib2R5Ijp7Im5hbWUiOiJub25lIn0sImhlYWQiOnsibmFtZSI6Im5vbmUifX19XV0=
      \begin{tikzcd}
        {\dbl{E}_x} & {\dbl{E}_{x'}} & {\dbl{E}_{x'}} \\
        {\dbl{E}_y} & {\dbl{E}_{y'}} \\
        {\dbl{E}_z} & {\dbl{E}_{z'}} & {\dbl{E}_{z'}}
        \arrow[""{name=0, anchor=center, inner sep=0}, "{\dbl{E}_f}"{inner sep=.8ex}, "\shortmid"{marking}, from=1-1, to=1-2]
        \arrow["{m_!}"', from=1-1, to=2-1]
        \arrow[""{name=1, anchor=center, inner sep=0}, "\shortmid"{marking}, equals, from=1-2, to=1-3]
        \arrow["{p_!}", from=1-2, to=2-2]
        \arrow["{(p\odot q)_!}", from=1-3, to=3-3]
        \arrow[""{name=2, anchor=center, inner sep=0}, "\shortmid"{marking}, from=2-1, to=2-2]
        \arrow["{n_!}"', from=2-1, to=3-1]
        \arrow["{q_!}", from=2-2, to=3-2]
        \arrow[""{name=3, anchor=center, inner sep=0}, "{\dbl{E}_h}"'{inner sep=.8ex}, "\shortmid"{marking}, from=3-1, to=3-2]
        \arrow[""{name=4, anchor=center, inner sep=0}, "\shortmid"{marking}, equals, from=3-2, to=3-3]
        \arrow["{\alpha_!}"{description}, draw=none, from=0, to=2]
        \arrow["{\dbl{E}^{p,q}}"{description}, draw=none, from=1, to=4]
        \arrow["{\beta_!}"{description}, draw=none, from=2, to=3]
      \end{tikzcd}
      \quad=\quad
      % https://q.uiver.app/#q=WzAsNyxbMCwwLCJcXGRibHtFfV94Il0sWzEsMCwiXFxkYmx7RX1feCJdLFswLDEsIlxcZGJse0V9X3kiXSxbMCwyLCJcXGRibHtFfV96Il0sWzEsMiwiXFxkYmx7RX1feCJdLFsyLDAsIlxcZGJse0V9X3t4J30iXSxbMiwyLCJcXGRibHtFfV97eid9Il0sWzAsMSwiIiwwLHsibGV2ZWwiOjIsInN0eWxlIjp7ImJvZHkiOnsibmFtZSI6ImJhcnJlZCJ9LCJoZWFkIjp7Im5hbWUiOiJub25lIn19fV0sWzAsMiwibV8hIiwyXSxbMiwzLCJuXyEiLDJdLFszLDQsIiIsMix7ImxldmVsIjoyLCJzdHlsZSI6eyJib2R5Ijp7Im5hbWUiOiJiYXJyZWQifSwiaGVhZCI6eyJuYW1lIjoibm9uZSJ9fX1dLFsxLDUsIlxcZGJse0V9X2YiLDAseyJzdHlsZSI6eyJib2R5Ijp7Im5hbWUiOiJiYXJyZWQifX19XSxbNSw2LCIocFxcb2RvdCBxKV8hIl0sWzQsNiwiXFxkYmx7RX1faCIsMix7InN0eWxlIjp7ImJvZHkiOnsibmFtZSI6ImJhcnJlZCJ9fX1dLFsxLDRdLFsxMSwxMywiKFxcYWxwaGFcXG9kb3RcXGJldGEpXyEiLDEseyJzaG9ydGVuIjp7InNvdXJjZSI6MjAsInRhcmdldCI6MjB9LCJzdHlsZSI6eyJib2R5Ijp7Im5hbWUiOiJub25lIn0sImhlYWQiOnsibmFtZSI6Im5vbmUifX19XSxbNywxMCwiXFxkYmx7RX1ee20sbn0iLDEseyJzaG9ydGVuIjp7InNvdXJjZSI6MjAsInRhcmdldCI6MjB9LCJzdHlsZSI6eyJib2R5Ijp7Im5hbWUiOiJub25lIn0sImhlYWQiOnsibmFtZSI6Im5vbmUifX19XV0=
      \begin{tikzcd}
        {\dbl{E}_x} & {\dbl{E}_x} & {\dbl{E}_{x'}} \\
        {\dbl{E}_y} \\
        {\dbl{E}_z} & {\dbl{E}_x} & {\dbl{E}_{z'}}
        \arrow[""{name=0, anchor=center, inner sep=0}, "\shortmid"{marking}, equals, from=1-1, to=1-2]
        \arrow["{m_!}"', from=1-1, to=2-1]
        \arrow[""{name=1, anchor=center, inner sep=0}, "{\dbl{E}_f}"{inner sep=.8ex}, "\shortmid"{marking}, from=1-2, to=1-3]
        \arrow[from=1-2, to=3-2]
        \arrow["{(p\odot q)_!}", from=1-3, to=3-3]
        \arrow["{n_!}"', from=2-1, to=3-1]
        \arrow[""{name=2, anchor=center, inner sep=0}, "\shortmid"{marking}, equals, from=3-1, to=3-2]
        \arrow[""{name=3, anchor=center, inner sep=0}, "{\dbl{E}_h}"'{inner sep=.8ex}, "\shortmid"{marking}, from=3-2, to=3-3]
        \arrow["{\dbl{E}^{m,n}}"{description}, draw=none, from=0, to=2]
        \arrow["{(\alpha\odot\beta)_!}"{description}, draw=none, from=1, to=3]
      \end{tikzcd}.
    \end{equation*}
  Chasing an arrow $u\colon a\to b$ over $f$ via $P$ through the cell 
  composites on each side, we have a purported equality of cells
    \begin{equation*}
      % https://q.uiver.app/#q=WzAsOCxbMCwwLCJhIl0sWzEsMCwibV8hYSJdLFsyLDAsIm5fIW1fIWEiXSxbMiwxLCJxXyFwXyFiIl0sWzAsMSwiYiJdLFsxLDEsInBfIWIiXSxbMCwyLCJiIl0sWzIsMiwiKHBcXG9kb3QgcSlfIWIiXSxbMCwxLCJcXHNpZ21hKG0sYSkiLDAseyJzdHlsZSI6eyJib2R5Ijp7Im5hbWUiOiJiYXJyZWQifX19XSxbMSwyLCJcXHNpZ21hKG4sbV8hYSkiLDAseyJzdHlsZSI6eyJib2R5Ijp7Im5hbWUiOiJiYXJyZWQifX19XSxbMiwzLCIiLDAseyJzdHlsZSI6eyJib2R5Ijp7Im5hbWUiOiJkYXNoZWQifX19XSxbMCw0LCJ1IiwyXSxbNCw1LCJcXHNpZ21hKHAsYikiLDIseyJzdHlsZSI6eyJib2R5Ijp7Im5hbWUiOiJiYXJyZWQifX19XSxbNSwzLCJcXHNpZ21hKHEscF8hYikiLDIseyJzdHlsZSI6eyJib2R5Ijp7Im5hbWUiOiJiYXJyZWQifX19XSxbMSw1LCIiLDAseyJzdHlsZSI6eyJib2R5Ijp7Im5hbWUiOiJkYXNoZWQifX19XSxbNCw2LCIiLDEseyJsZXZlbCI6Miwic3R5bGUiOnsiaGVhZCI6eyJuYW1lIjoibm9uZSJ9fX1dLFs2LDcsIlxcc2lnbWEocFxcb2RvdCBxLGIpIiwyLHsic3R5bGUiOnsiYm9keSI6eyJuYW1lIjoiYmFycmVkIn19fV0sWzMsNywiIiwwLHsic3R5bGUiOnsiYm9keSI6eyJuYW1lIjoiZGFzaGVkIn19fV0sWzgsMTIsIlxcc2lnbWEodSxcXGFscGhhKSIsMSx7InNob3J0ZW4iOnsic291cmNlIjoyMCwidGFyZ2V0IjoyMH0sInN0eWxlIjp7ImJvZHkiOnsibmFtZSI6Im5vbmUifSwiaGVhZCI6eyJuYW1lIjoibm9uZSJ9fX1dLFs5LDEzLCJcXHNpZ21hKChcXGFscGhhXyEpX3UsXFxiZXRhKSIsMSx7Im9mZnNldCI6LTEsInNob3J0ZW4iOnsic291cmNlIjoyMCwidGFyZ2V0IjoyMH0sInN0eWxlIjp7ImJvZHkiOnsibmFtZSI6Im5vbmUifSwiaGVhZCI6eyJuYW1lIjoibm9uZSJ9fX1dLFs1LDE2LCJcXHNpZ21hKDFfYixQXnstMX0pIiwxLHsibGFiZWxfcG9zaXRpb24iOjQwLCJzaG9ydGVuIjp7InRhcmdldCI6MjB9LCJzdHlsZSI6eyJib2R5Ijp7Im5hbWUiOiJub25lIn0sImhlYWQiOnsibmFtZSI6Im5vbmUifX19XV0=
      \begin{tikzcd}
        a & {m_!a} & {n_!m_!a} \\
        b & {p_!b} & {q_!p_!b} \\
        b && {(p\odot q)_!b}
        \arrow[""{name=0, anchor=center, inner sep=0}, "{\sigma(m,a)}"{inner sep=.8ex}, "\shortmid"{marking}, from=1-1, to=1-2]
        \arrow["u"', from=1-1, to=2-1]
        \arrow[""{name=1, anchor=center, inner sep=0}, "{\sigma(n,m_!a)}"{inner sep=.8ex}, "\shortmid"{marking}, from=1-2, to=1-3]
        \arrow[dashed, from=1-2, to=2-2]
        \arrow[dashed, from=1-3, to=2-3]
        \arrow[""{name=2, anchor=center, inner sep=0}, "{\sigma(p,b)}"'{inner sep=.8ex}, "\shortmid"{marking}, from=2-1, to=2-2]
        \arrow[equals, from=2-1, to=3-1]
        \arrow[""{name=3, anchor=center, inner sep=0}, "{\sigma(q,p_!b)}"'{inner sep=.8ex}, "\shortmid"{marking}, from=2-2, to=2-3]
        \arrow[dashed, from=2-3, to=3-3]
        \arrow[""{name=4, anchor=center, inner sep=0}, "{\sigma(p\odot q,b)}"'{inner sep=.8ex}, "\shortmid"{marking}, from=3-1, to=3-3]
        \arrow["{\sigma(u,\alpha)}"{description}, draw=none, from=0, to=2]
        \arrow["{\sigma((\alpha_!)_u,\beta)}"{description}, shift left, draw=none, from=1, to=3]
        \arrow["{\sigma(1_b,P^{-1})}"{description, pos=0.4}, draw=none, from=2-2, to=4]
      \end{tikzcd} \quad = \quad 
      % https://q.uiver.app/#q=WzAsNyxbMCwwLCJhIl0sWzEsMCwibV8hYSJdLFsyLDAsIm5fIW1fIWEiXSxbMiwxLCIobVxcb2RvdCBuKV8hYiJdLFswLDEsImEiXSxbMCwyLCJiIl0sWzIsMiwiKHBcXG9kb3QgcSlfIWIiXSxbMCwxLCJcXHNpZ21hKG0sYSkiLDAseyJzdHlsZSI6eyJib2R5Ijp7Im5hbWUiOiJiYXJyZWQifX19XSxbMSwyLCJcXHNpZ21hKG4sbV8hYSkiLDAseyJzdHlsZSI6eyJib2R5Ijp7Im5hbWUiOiJiYXJyZWQifX19XSxbMiwzLCIiLDAseyJzdHlsZSI6eyJib2R5Ijp7Im5hbWUiOiJkYXNoZWQifX19XSxbMCw0LCIiLDIseyJsZXZlbCI6Miwic3R5bGUiOnsiaGVhZCI6eyJuYW1lIjoibm9uZSJ9fX1dLFs0LDUsInUiLDJdLFs1LDYsIlxcc2lnbWEocFxcb2RvdCBxLGIpIiwyLHsic3R5bGUiOnsiYm9keSI6eyJuYW1lIjoiYmFycmVkIn19fV0sWzMsNiwiIiwwLHsic3R5bGUiOnsiYm9keSI6eyJuYW1lIjoiZGFzaGVkIn19fV0sWzQsMywiXFxzaWdtYShtXFxvZG90IG4sYikiLDAseyJzdHlsZSI6eyJib2R5Ijp7Im5hbWUiOiJiYXJyZWQifX19XSxbMSwxNCwiXFxzaWdtYSgxX2EsUF57LTF9KSIsMSx7ImxhYmVsX3Bvc2l0aW9uIjozMCwic2hvcnRlbiI6eyJ0YXJnZXQiOjIwfSwic3R5bGUiOnsiYm9keSI6eyJuYW1lIjoibm9uZSJ9LCJoZWFkIjp7Im5hbWUiOiJub25lIn19fV0sWzE0LDEyLCJcXHNpZ21hKHUsXFxhbHBoYVxcb2RvdFxcYmV0YSkiLDEseyJzaG9ydGVuIjp7InNvdXJjZSI6MjAsInRhcmdldCI6MjB9LCJsZXZlbCI6MSwic3R5bGUiOnsiYm9keSI6eyJuYW1lIjoibm9uZSJ9LCJoZWFkIjp7Im5hbWUiOiJub25lIn19fV1d
      \begin{tikzcd}
        a & {m_!a} & {n_!m_!a} \\
        a && {(m\odot n)_!b} \\
        b && {(p\odot q)_!b}
        \arrow["{\sigma(m,a)}"{inner sep=.8ex}, "\shortmid"{marking}, from=1-1, to=1-2]
        \arrow[equals, from=1-1, to=2-1]
        \arrow["{\sigma(n,m_!a)}"{inner sep=.8ex}, "\shortmid"{marking}, from=1-2, to=1-3]
        \arrow[dashed, from=1-3, to=2-3]
        \arrow[""{name=0, anchor=center, inner sep=0}, "{\sigma(m\odot n,b)}"{inner sep=.8ex}, "\shortmid"{marking}, from=2-1, to=2-3]
        \arrow["u"', from=2-1, to=3-1]
        \arrow[dashed, from=2-3, to=3-3]
        \arrow[""{name=1, anchor=center, inner sep=0}, "{\sigma(p\odot q,b)}"'{inner sep=.8ex}, "\shortmid"{marking}, from=3-1, to=3-3]
        \arrow["{\sigma(1_a,P^{-1})}"{description, pos=0.3}, draw=none, from=1-2, to=0]
        \arrow["{\sigma(u,\alpha\odot\beta)}"{description}, draw=none, from=0, to=1]
      \end{tikzcd}
    \end{equation*}
  where each $P^{-1}$ denotes the inverse of an appropriate component of a laxator 
  comparison coming with $P$. But this equality holds since $P$ is a double functor 
  (and thus the comparison isomorphisms are natural) and such lifts are unique. For 
  the loose composition comparison isomorphisms coming with $P$ are natural with 
  respect to loose composition of cells of $\dbl{E}$. Thus, in particular the dashed 
  arrow composites on the right side of each boundary must be equal. The 
  remaining checks are straightforward.

  As for the claim that $\Tw(P)$ is unitary: this was essentially observed 
  already in \cref{construction:pro_fibers} since by definition
    \begin{equation*}
      \Tw(P)(1_x) = \dbl{E}_{1_x}\colon\dbl{E}_x^{\op}\times\dbl{E}_x\to\Set \qquad (a,b) \mapsto 
      \lbrace u\colon a\to b\mid Pu=1_x\rbrace
    \end{equation*}
  which is exactly the set of morphisms of $\dbl{E}$ above the identity on $x$, 
  which is how we defined the morphisms of the fiber category $\dbl{E}_x$. So, 
  the profunctor $\Tw(P)(1_x) = \dbl{E}_x(-,=)$ identically and the comparisons 
  can be taken to be identities.
\end{proof}

\subsection{Completing the elements correspondence}
\label{subsection:completing-elements-correspondence}

The goal of this subsection is to show that there is a correspondence
  \begin{equation*}
    \lbrace \text{cloven loosely discrete opfibrations}\rbrace \quad\leftrightarrow\quad\lbrace \text{normal profunctor-valued twisted functors} \rbrace.
  \end{equation*}
The elements (\cref{construction:weak-elements}) and twisted functor
(\cref{prop:twisted-functor-assoc-to-fibration}) constructions guarantee that we
have assignments in each direction, but in order to have a proper \emph{elements
  correspondence}, we need to see that these round-trip in some sense. Such a
round-tripping is also a prelude to a full 2-categorical account of the
relationship between loosely discrete opfibrations and twisted functors.

Start with a twisted normal lax functor $T\colon\dbl{B}\twistto\Prof$. The elements
construction applied to $T$ results in a loosely discrete opfibration
$P_T\colon\Elt(T)\to\dbl{B}$. Applying the \emph{pseudo-inverse} construction from
\cref{subsection:pseudo-inverse-construction} then yields a twisted functor
$\Tw(P_T)\colon\dbl{B}\twistto\Prof$. The following development will produce a natural
transformation
  \begin{equation*}
    \gamma\colon T\To \Tw(P_T)
  \end{equation*}
in the sense of \cref{def:naturaltransformation}, all of whose components are
invertible. In this sense, every twisted functor occurs as the twisted functor
canonically constructed out of a loosely discrete opfibration. Moreover,
since the constructed twisted functor is unitary, this implies that every normal
twisted functor is comparable to, even isomorphic to, a unitary twisted functor.

The following three constructions supply the data of this purported
transformation $\gamma$.

\begin{construction}[Arrow component]
  First, fix an object $x\in\dbl{B}$. Required is a functor
    \begin{equation*}
      \gamma_x\colon Tx\to\Tw(P_T)(x).
    \end{equation*}
  This is given on objects by sending $a\in Tx$ to the pair $(x,a)$, and on arrows in $Tx$ by
    \begin{equation*}
      (u\colon a\to b) \qquad \mapsto \qquad (1_x,T_x(u))\colon (x,a)\to (x,b),
    \end{equation*}
  that is, taking the arrow $u$ and moving it into $T1_x(a,b)$ via the identity
  comparison and then viewing it as a pair indexed by $1$ itself. The assignment $\gamma_x$
  is functorial by the unit coherence conditions for a twisted lax functor and 
  by the unit conditions in $\Prof$. Note that $\gamma_x$ has a genuine \emph{inverse}
  since $T$ itself is normal and thus $T_x$ has an inverse with
  respect to tight composition of cells.
\end{construction}

\begin{construction}[Cell component]
  Now fix an arrow $f\colon x\to z$ of $\dbl{B}$. Required is the further
  data of a cell 
    \begin{equation*}
      % https://q.uiver.app/#q=WzAsNCxbMCwwLCJUeCJdLFsxLDAsIlR6Il0sWzEsMSwiXFxUdyhQX1QpKHopIl0sWzAsMSwiXFxUdyhQX1QpKHgpIl0sWzAsMSwiVGYiLDAseyJzdHlsZSI6eyJib2R5Ijp7Im5hbWUiOiJiYXJyZWQifX19XSxbMSwyLCJcXGdhbW1hX3oiXSxbMCwzLCJcXGdhbW1hX3giLDJdLFszLDIsIlxcVHcoUF9UKShmKSIsMix7InN0eWxlIjp7ImJvZHkiOnsibmFtZSI6ImJhcnJlZCJ9fX1dLFs0LDcsIlxcZ2FtbWFfZiIsMSx7InNob3J0ZW4iOnsic291cmNlIjoyMCwidGFyZ2V0IjoyMH0sInN0eWxlIjp7ImJvZHkiOnsibmFtZSI6Im5vbmUifSwiaGVhZCI6eyJuYW1lIjoibm9uZSJ9fX1dXQ==
      \begin{tikzcd}[column sep=large]
        Tx & Tz \\
        {\Tw(P_T)(x)} & {\Tw(P_T)(z)}
        \arrow[""{name=0, anchor=center, inner sep=0}, "Tf"{inner sep=.8ex}, "\shortmid"{marking}, from=1-1, to=1-2]
        \arrow["{\gamma_x}"', from=1-1, to=2-1]
        \arrow["{\gamma_z}", from=1-2, to=2-2]
        \arrow[""{name=1, anchor=center, inner sep=0}, "{\Tw(P_T)(f)}"'{inner sep=.8ex}, "\shortmid"{marking}, from=2-1, to=2-2]
        \arrow["{\gamma_f}"{description}, draw=none, from=0, to=1]
      \end{tikzcd}
    \end{equation*}
  that is, a family of functions 
    \begin{equation*}
      (\gamma_f)_{a,c}\colon Tf(a,c) \to \Tw(P_T)(f)((x,a),(z,c))
    \end{equation*}
  natural in $a$ and $c$. Of course, the assignment is to take a heteromorphism 
  $u\in Tf(a,c)$ to the pair formed by viewing $f$ itself as an index:
    \begin{equation*}
      u\in Tf(a,c) \qquad\mapsto\qquad (f,u)\colon (x,a)\to (z,c).
    \end{equation*}
  Naturality does indeed hold and follows specifically by a double application of 
  unit coherence. Again each such component $(\gamma_f)_{a,c}$ is in fact 
  invertible owing to the fact that there is an evident inverse correspondence,
  namely, that of just projection via \emph{forgetting the indexing morphism}.
\end{construction}

\begin{construction}[Naturality comparisons]
  Let $m\colon x\proto y$ denote any proarrow of $\dbl{B}$. The associated 
  naturality comparison cell 
    \begin{equation*}
      % https://q.uiver.app/#q=WzAsNixbMCwwLCJUeCJdLFsxLDAsIlR4Il0sWzAsMSwiXFxUdyhQX1QpKHgpIl0sWzAsMiwiXFxUdyhQX1QpKHkpIl0sWzEsMiwiXFxUdyhQX1QpKHkpIl0sWzEsMSwiVHkiXSxbMCwxLCIiLDIseyJsZXZlbCI6Miwic3R5bGUiOnsiYm9keSI6eyJuYW1lIjoiYmFycmVkIn0sImhlYWQiOnsibmFtZSI6Im5vbmUifX19XSxbMCwyLCJcXGdhbW1hX3giLDJdLFsyLDMsIlxcVHcoUF9UKShtKSIsMl0sWzMsNCwiIiwwLHsibGV2ZWwiOjIsInN0eWxlIjp7ImJvZHkiOnsibmFtZSI6ImJhcnJlZCJ9LCJoZWFkIjp7Im5hbWUiOiJub25lIn19fV0sWzEsNSwiVG0iXSxbNSw0LCJcXGdhbW1hX3kiXSxbNiw5LCJcXGdhbW1hX20iLDEseyJzaG9ydGVuIjp7InNvdXJjZSI6MjAsInRhcmdldCI6MjB9LCJzdHlsZSI6eyJib2R5Ijp7Im5hbWUiOiJub25lIn0sImhlYWQiOnsibmFtZSI6Im5vbmUifX19XV0=
      \begin{tikzcd}
        Tx & Tx \\
        {\Tw(P_T)(x)} & Ty \\
        {\Tw(P_T)(y)} & {\Tw(P_T)(y)}
        \arrow[""{name=0, anchor=center, inner sep=0}, "\shortmid"{marking}, equals, from=1-1, to=1-2]
        \arrow["{\gamma_x}"', from=1-1, to=2-1]
        \arrow["Tm", from=1-2, to=2-2]
        \arrow["{\Tw(P_T)(m)}"', from=2-1, to=3-1]
        \arrow["{\gamma_y}", from=2-2, to=3-2]
        \arrow[""{name=1, anchor=center, inner sep=0}, "\shortmid"{marking}, equals, from=3-1, to=3-2]
        \arrow["{\gamma_m}"{description}, draw=none, from=0, to=1]
      \end{tikzcd}
    \end{equation*}
  is the loose identity. A component thus takes a morphism $u$ of $Tx$ to the pair
    \begin{equation*}
      (u\colon a\to b) \qquad\mapsto\qquad (y, Tm(u))\colon (y,Tm(a))\to (y,Tm(b)).
    \end{equation*}
  Note that this is well-typed since the source and target of the cell above have the 
  same action and are thus the same functor.
\end{construction}

Now we can show that every profunctor-valued twisted functor is the twisted functor 
canonically associated to some loosely discrete opfibration.

\begin{proposition}
  The data from the preceding constructions results in an invertible strict 
  natural transformation $\gamma\colon T\To\Tw(P_T)$.
\end{proposition}
\begin{proof}
  It has already been seen that the components of $\gamma$ are well-defined 
  and invertible. To see that $\gamma$ is a lax transformation, we need to verify the conditions of \cref{def:naturaltransformation}. There are three batches of such conditions: naturality, functoriality and coherence. All of these, however, follow by construction of $\gamma$, the definitions of the elements construction, the canonical twisted functor construction and the original properties of $T$. For example, cell naturality asserts that for each cell
    \begin{equation*}
      % https://q.uiver.app/#q=WzAsNCxbMCwwLCJ4Il0sWzEsMCwieSJdLFsxLDEsInciXSxbMCwxLCJ6Il0sWzAsMSwibSIsMCx7InN0eWxlIjp7ImJvZHkiOnsibmFtZSI6ImJhcnJlZCJ9fX1dLFsxLDIsImciXSxbMCwzLCJmIiwyXSxbMywyLCJuIiwyLHsic3R5bGUiOnsiYm9keSI6eyJuYW1lIjoiYmFycmVkIn19fV0sWzQsNywiXFxhbHBoYSIsMSx7InNob3J0ZW4iOnsic291cmNlIjoyMCwidGFyZ2V0IjoyMH0sInN0eWxlIjp7ImJvZHkiOnsibmFtZSI6Im5vbmUifSwiaGVhZCI6eyJuYW1lIjoibm9uZSJ9fX1dXQ==
      \begin{tikzcd}
        x & y \\
        z & w
        \arrow[""{name=0, anchor=center, inner sep=0}, "m"{inner sep=.8ex}, "\shortmid"{marking}, from=1-1, to=1-2]
        \arrow["f"', from=1-1, to=2-1]
        \arrow["g", from=1-2, to=2-2]
        \arrow[""{name=1, anchor=center, inner sep=0}, "n"'{inner sep=.8ex}, "\shortmid"{marking}, from=2-1, to=2-2]
        \arrow["\alpha"{description}, draw=none, from=0, to=1]
      \end{tikzcd}
    \end{equation*} 
  of $\dbl{B}$, we have an equality of composite diagrams 
    \begin{equation*}
      % https://q.uiver.app/#q=WzAsNixbMCwwLCJUeCJdLFsxLDAsIlR6Il0sWzEsMSwiXFxUdyhQKSh6KSJdLFswLDEsIlxcVHcoUCkoeCkiXSxbMCwyLCJcXFR3KFApKHkpIl0sWzEsMiwiXFxUdyhQKSh3KSJdLFswLDEsIlRmIiwwLHsic3R5bGUiOnsiYm9keSI6eyJuYW1lIjoiYmFycmVkIn19fV0sWzEsMiwiXFxnYW1tYV96Il0sWzAsMywiXFxnYW1tYV94IiwyXSxbMywyLCJcXFR3KFApKGYpIiwwLHsic3R5bGUiOnsiYm9keSI6eyJuYW1lIjoiYmFycmVkIn19fV0sWzMsNCwiXFxUdyhQKShtKSIsMl0sWzQsNSwiXFxUdyhQKShnKSIsMix7InN0eWxlIjp7ImJvZHkiOnsibmFtZSI6ImJhcnJlZCJ9fX1dLFsyLDUsIlxcVHcoUCkobikiLDFdLFs2LDksIlxcZ2FtbWFfZiIsMSx7ImxhYmVsX3Bvc2l0aW9uIjo0MCwic2hvcnRlbiI6eyJzb3VyY2UiOjIwLCJ0YXJnZXQiOjIwfSwic3R5bGUiOnsiYm9keSI6eyJuYW1lIjoibm9uZSJ9LCJoZWFkIjp7Im5hbWUiOiJub25lIn19fV0sWzksMTEsIlxcVHcoUCkoXFxhbHBoYSkiLDEseyJzaG9ydGVuIjp7InNvdXJjZSI6MjAsInRhcmdldCI6MjB9LCJzdHlsZSI6eyJib2R5Ijp7Im5hbWUiOiJub25lIn0sImhlYWQiOnsibmFtZSI6Im5vbmUifX19XV0=
      \begin{tikzcd}
        Tx & Tz \\
        {\Tw(P)(x)} & {\Tw(P)(z)} \\
        {\Tw(P)(y)} & {\Tw(P)(w)}
        \arrow[""{name=0, anchor=center, inner sep=0}, "Tf"{inner sep=.8ex}, "\shortmid"{marking}, from=1-1, to=1-2]
        \arrow["{\gamma_x}"', from=1-1, to=2-1]
        \arrow["{\gamma_z}", from=1-2, to=2-2]
        \arrow[""{name=1, anchor=center, inner sep=0}, "{\Tw(P)(f)}"{inner sep=.8ex}, "\shortmid"{marking}, from=2-1, to=2-2]
        \arrow["{\Tw(P)(m)}"', from=2-1, to=3-1]
        \arrow["{\Tw(P)(n)}"{description}, from=2-2, to=3-2]
        \arrow[""{name=2, anchor=center, inner sep=0}, "{\Tw(P)(g)}"'{inner sep=.8ex}, "\shortmid"{marking}, from=3-1, to=3-2]
        \arrow["{\gamma_f}"{description, pos=0.4}, draw=none, from=0, to=1]
        \arrow["{\Tw(P)(\alpha)}"{description}, draw=none, from=1, to=2]
      \end{tikzcd}\qquad = \qquad
      % https://q.uiver.app/#q=WzAsNixbMCwwLCJUeCJdLFswLDEsIlR5Il0sWzAsMiwiXFxUdyhQKSh5KSJdLFsxLDAsIlR6Il0sWzEsMSwiVHciXSxbMSwyLCJcXFR3KFApKHcpIl0sWzAsMSwiVG0iLDFdLFsxLDIsIlxcZ2FtbWFfeSIsMV0sWzAsMywiVGYiLDAseyJzdHlsZSI6eyJib2R5Ijp7Im5hbWUiOiJiYXJyZWQifX19XSxbMyw0LCJUbiJdLFsxLDQsIlRnIiwwLHsic3R5bGUiOnsiYm9keSI6eyJuYW1lIjoiYmFycmVkIn19fV0sWzIsNSwiXFxUdyhQKShnKSIsMix7InN0eWxlIjp7ImJvZHkiOnsibmFtZSI6ImJhcnJlZCJ9fX1dLFs0LDUsIlxcZ2FtbWFfdyJdLFsxMCwxMSwiXFxnYW1tYV9nIiwxLHsic2hvcnRlbiI6eyJzb3VyY2UiOjIwLCJ0YXJnZXQiOjIwfSwic3R5bGUiOnsiYm9keSI6eyJuYW1lIjoibm9uZSJ9LCJoZWFkIjp7Im5hbWUiOiJub25lIn19fV0sWzgsMTAsIlRcXGFscGhhIiwxLHsibGFiZWxfcG9zaXRpb24iOjQwLCJzaG9ydGVuIjp7InNvdXJjZSI6MjAsInRhcmdldCI6MjB9LCJzdHlsZSI6eyJib2R5Ijp7Im5hbWUiOiJub25lIn0sImhlYWQiOnsibmFtZSI6Im5vbmUifX19XV0=
      \begin{tikzcd}
        Tx & Tz \\
        Ty & Tw \\
        {\Tw(P)(y)} & {\Tw(P)(w)}
        \arrow[""{name=0, anchor=center, inner sep=0}, "Tf"{inner sep=.8ex}, "\shortmid"{marking}, from=1-1, to=1-2]
        \arrow["Tm"{description}, from=1-1, to=2-1]
        \arrow["Tn", from=1-2, to=2-2]
        \arrow[""{name=1, anchor=center, inner sep=0}, "Tg"{inner sep=.8ex}, "\shortmid"{marking}, from=2-1, to=2-2]
        \arrow["{\gamma_y}"{description}, from=2-1, to=3-1]
        \arrow["{\gamma_w}", from=2-2, to=3-2]
        \arrow[""{name=2, anchor=center, inner sep=0}, "{\Tw(P)(g)}"'{inner sep=.8ex}, "\shortmid"{marking}, from=3-1, to=3-2]
        \arrow["{T\alpha}"{description, pos=0.4}, draw=none, from=0, to=1]
        \arrow["{\gamma_g}"{description}, draw=none, from=1, to=2]
      \end{tikzcd}
    \end{equation*}
  with the naturality comparison cells suppressed since they are loose identities. 
  By the definition of the unit comparison cells in $\Prof$, chasing through each 
  of these composite diagrams amounts to chasing a heteromorphism $u\in Tf(a,c)$
  through the stacked squares on either side. The action in each case is exactly 
  $(g, T\alpha(u))$ and thus the composites are equal. The verifications for the 
  other batches of conditions are similar.

  As for invertibility: the pattern of definition of $\gamma$ suggests one for 
  a transformation in the reverse direction $\delta\colon \Tw(P_T)\Rightarrow T$
  utilizing the inverses of components of the unitor comparison morphisms 
  $(T^x)^{-1}$ as $x$ ranges over the objects of $\dbl{B}$. Put another way, the 
  two components of $\delta$ would be the literal inverse assignments suggested
  above in the construction of the components of $\gamma$.
\end{proof}

\begin{corollary}
  Every normal profunctor-valued twisted functor is canonically isomorphic to a
  \emph{unitary} profunctor-valued twisted functor.
\end{corollary}
\begin{proof}
  This directly follows from the previous result, the fact that $\Tw(P_T)$ is 
  unitary (\cref{prop:twisted-functor-assoc-to-fibration}), and the fact 
  that a transformation with invertible components is indeed an isomorphism in 
  the category of twisted functors and lax transformations.
\end{proof}

In the other direction, fixing a loosely discrete opfibration
$P\colon\dbl{E}\to\dbl{B}$ with a cleavage $\sigma$, we now seek a comparison
double functor from $P$ itself to the elements construction applied to the
twisted functor built from $P$:
  \begin{equation*}
    % https://q.uiver.app/#q=WzAsMyxbMCwwLCJcXGRibHtFfSJdLFsyLDAsIlxcRWx0KFxcVHcoUCkpIl0sWzEsMSwiXFxkYmx7Qn0iXSxbMCwxLCJSIl0sWzAsMiwiUCIsMl0sWzEsMiwiSiJdXQ==
    \begin{tikzcd}
      {\dbl{E}} && {\Elt(\Tw(P))} \\
      & {\dbl{B}}
      \arrow["R", from=1-1, to=1-3]
      \arrow["P"', from=1-1, to=2-2]
      \arrow["J", from=1-3, to=2-2]
    \end{tikzcd}.
  \end{equation*}
The rest of the subsection is dedicated to this task and to describing the sense
in which the resulting double functor $R$ is an equivalence.

\begin{construction}
  Fix a loosely discrete opfibration $P: \dbl{E} \to \dbl{B}$. Assume it is
  equipped with a cleavage
  $\sigma: \dbl{E}_0 \times_{\dbl{B}_0} \dbl{B}_1 \to \dbl{E}_1$
  (\cref{def:cleavage-for-disc-opfibration}). Here we propose assignments for
  what will amount to a double functor
    \begin{equation*}
      R\colon \dbl{E} \to \Elt(\Tw(P)).
    \end{equation*}
  On the underlying category of objects, these are:
    \begin{itemize}
      \item An object $a$ maps to $(Pa,a)$;
      \item An arrow $u\colon a\to b$ maps to $(Pu, u)\colon (Pa,a)\to (Pb,b)$.
    \end{itemize}
  These assignments are well-defined and functorial. The resulting functor
  $R_0\colon \dbl{E}_0 \to \Elt(\Tw(P))_0$ is an equivalence of categories as $\Tw(P)$
  acts on $\dbl{B}_0$ simply by taking preimages under $P$.

  For the assignments on proarrows and cells, recall that $\epsilon$ denotes the counit
  of the split adjoint equivalence defining the cleavage $\sigma$, whose
  component at a proarrow $t: a \proto b$ in $\dbl{E}$ has the form
    \begin{equation*}
      % https://q.uiver.app/#q=WzAsNCxbMCwwLCJhIl0sWzEsMCwiKFB0KV8hKGEpIl0sWzEsMSwiYiJdLFswLDEsImEiXSxbMCwxLCJcXHNpZ21hKGEsUHQpIiwwLHsic3R5bGUiOnsiYm9keSI6eyJuYW1lIjoiYmFycmVkIn19fV0sWzEsMiwiXFx0Z3RcXGVwc2lsb25fdCJdLFswLDMsIiIsMCx7ImxldmVsIjoyLCJzdHlsZSI6eyJoZWFkIjp7Im5hbWUiOiJub25lIn19fV0sWzMsMiwidCIsMix7InN0eWxlIjp7ImJvZHkiOnsibmFtZSI6ImJhcnJlZCJ9fX1dLFs0LDcsIlxcZXBzaWxvbl90IiwxLHsic2hvcnRlbiI6eyJzb3VyY2UiOjIwLCJ0YXJnZXQiOjIwfSwic3R5bGUiOnsiYm9keSI6eyJuYW1lIjoibm9uZSJ9LCJoZWFkIjp7Im5hbWUiOiJub25lIn19fV1d
      \begin{tikzcd}
        a & {(Pt)_!(a)} \\
        a & b
        \arrow[""{name=0, anchor=center, inner sep=0}, "{\sigma(a,Pt)}"{inner sep=.8ex}, "\shortmid"{marking}, from=1-1, to=1-2]
        \arrow[equals, from=1-1, to=2-1]
        \arrow["{\tgt\epsilon_t}", from=1-2, to=2-2]
        \arrow[""{name=1, anchor=center, inner sep=0}, "t"'{inner sep=.8ex}, "\shortmid"{marking}, from=2-1, to=2-2]
        \arrow["{\epsilon_t}"{description}, draw=none, from=0, to=1]
      \end{tikzcd}.
    \end{equation*}
  The assignment is then as follows.
    \begin{enumerate}
      \item A proarrow $t\colon a\proto b$ is sent to $(Pt,\tgt\epsilon_t)\colon (Pa,a)\proto (Pb,b)$;
      \item A cell $\alpha$ of $\dbl{E}$ as on the left is sent to the cell $P\alpha$ on the
        right:
        \begin{equation*}
          % https://q.uiver.app/#q=WzAsNCxbMCwwLCJhIl0sWzEsMCwiYiJdLFswLDEsImMiXSxbMSwxLCJkIl0sWzAsMSwidCIsMCx7InN0eWxlIjp7ImJvZHkiOnsibmFtZSI6ImJhcnJlZCJ9fX1dLFswLDIsInUiLDJdLFsyLDMsInMiLDIseyJzdHlsZSI6eyJib2R5Ijp7Im5hbWUiOiJiYXJyZWQifX19XSxbMSwzLCJ2Il0sWzQsNiwiXFxhbHBoYSIsMSx7InNob3J0ZW4iOnsic291cmNlIjoyMCwidGFyZ2V0IjoyMH0sInN0eWxlIjp7ImJvZHkiOnsibmFtZSI6Im5vbmUifSwiaGVhZCI6eyJuYW1lIjoibm9uZSJ9fX1dXQ==
          \begin{tikzcd}
            a & b \\
            c & d
            \arrow[""{name=0, anchor=center, inner sep=0}, "t"{inner sep=.8ex}, "\shortmid"{marking}, from=1-1, to=1-2]
            \arrow["u"', from=1-1, to=2-1]
            \arrow["v", from=1-2, to=2-2]
            \arrow[""{name=1, anchor=center, inner sep=0}, "s"'{inner sep=.8ex}, "\shortmid"{marking}, from=2-1, to=2-2]
            \arrow["\alpha"{description}, draw=none, from=0, to=1]
          \end{tikzcd}
          \qquad \mapsto \qquad
          % https://q.uiver.app/#q=WzAsNCxbMCwwLCIoUGEsYSkiXSxbMSwwLCIoUGIsYikiXSxbMCwxLCIoUGMsYykiXSxbMSwxLCIoUGQsZCkiXSxbMCwxLCIoUHQsXFx0Z3RcXGVwc2lsb25fdCkiLDAseyJzdHlsZSI6eyJib2R5Ijp7Im5hbWUiOiJiYXJyZWQifX19XSxbMCwyLCIoUHUsdSkiLDJdLFsyLDMsIihQcyxcXHRndFxcZXBzaWxvbl9zKSIsMix7InN0eWxlIjp7ImJvZHkiOnsibmFtZSI6ImJhcnJlZCJ9fX1dLFsxLDMsIihQdix2KSJdLFs0LDYsIlBcXGFscGhhIiwxLHsic2hvcnRlbiI6eyJzb3VyY2UiOjIwLCJ0YXJnZXQiOjIwfSwic3R5bGUiOnsiYm9keSI6eyJuYW1lIjoibm9uZSJ9LCJoZWFkIjp7Im5hbWUiOiJub25lIn19fV1d
          \begin{tikzcd}
            {(Pa,a)} & {(Pb,b)} \\
            {(Pc,c)} & {(Pd,d)}
            \arrow[""{name=0, anchor=center, inner sep=0}, "{(Pt,\tgt\epsilon_t)}"{inner sep=.8ex}, "\shortmid"{marking}, from=1-1, to=1-2]
            \arrow["{(Pu,u)}"', from=1-1, to=2-1]
            \arrow["{(Pv,v)}", from=1-2, to=2-2]
            \arrow[""{name=1, anchor=center, inner sep=0}, "{(Ps,\tgt\epsilon_s)}"'{inner sep=.8ex}, "\shortmid"{marking}, from=2-1, to=2-2]
            \arrow["{P\alpha}"{description}, draw=none, from=0, to=1]
          \end{tikzcd}.
        \end{equation*}
    \end{enumerate}
  The point now is to check well-definition in each case. For the proarrow 
  assignment, this is just to observe that the targets of the counit cells,
  $\tgt\varepsilon_t$,
  are \emph{vertical} in the sense that they are arrows of the appropriate 
  fiber categories $\dbl{E}_{Pb}$ (\cref{lemma:counit-equations}). For
  cell well-definition, the condition unpacked in 
  \cref{rmk:elements-cell-well-definition} needs to be verified. But this is a 
  consequence of the naturality of $\epsilon$.
\end{construction}

It remains to see now that these assignments can be boosted to the data of a
double functor. For this we need comparison cells for loose composition and
identities. The work will be in showing that these constructions are
well-defined; the upshot is that by construction functoriality and coherence
will for the most part follow automatically.

\begin{construction}[Loose composition comparisons]
  Adopt the same setup as in the previous construction. Let $s\colon a\proto b$ 
  and $t\colon b\proto c$ be composable proarrows of $\dbl{E}$. The proposal
  is that the comparison cell associated to this pair is 
    \begin{equation*}
      % https://q.uiver.app/#q=WzAsNSxbMCwwLCIoUGEsYSkiXSxbMSwwLCIoUGIsYikiXSxbMiwwLCIoUGMsYykiXSxbMiwxLCIoUGMsYykiXSxbMCwxLCIoUGEsYSkiXSxbMCwxLCIoUHMsXFx0Z3RcXGVwc2lsb25fcykiLDAseyJzdHlsZSI6eyJib2R5Ijp7Im5hbWUiOiJiYXJyZWQifX19XSxbMSwyLCIoUHQsXFx0Z3RcXGVwc2lsb25fdCkiLDAseyJzdHlsZSI6eyJib2R5Ijp7Im5hbWUiOiJiYXJyZWQifX19XSxbMiwzLCIiLDIseyJsZXZlbCI6Miwic3R5bGUiOnsiaGVhZCI6eyJuYW1lIjoibm9uZSJ9fX1dLFswLDQsIiIsMCx7ImxldmVsIjoyLCJzdHlsZSI6eyJoZWFkIjp7Im5hbWUiOiJub25lIn19fV0sWzQsMywiKFAoc1xcb2RvdCB0KSxcXHRndFxcZXBzaWxvbl97c1xcb2RvdCB0fSkiLDIseyJzdHlsZSI6eyJib2R5Ijp7Im5hbWUiOiJiYXJyZWQifX19XSxbMSw5LCJQX3tzLHR9IiwxLHsibGFiZWxfcG9zaXRpb24iOjQwLCJzaG9ydGVuIjp7InRhcmdldCI6MjB9LCJzdHlsZSI6eyJib2R5Ijp7Im5hbWUiOiJub25lIn0sImhlYWQiOnsibmFtZSI6Im5vbmUifX19XV0=
      \begin{tikzcd}
        {(Pa,a)} & {(Pb,b)} & {(Pc,c)} \\
        {(Pa,a)} && {(Pc,c)}
        \arrow["{(Ps,\tgt\epsilon_s)}"{inner sep=.8ex}, "\shortmid"{marking}, from=1-1, to=1-2]
        \arrow[equals, from=1-1, to=2-1]
        \arrow["{(Pt,\tgt\epsilon_t)}"{inner sep=.8ex}, "\shortmid"{marking}, from=1-2, to=1-3]
        \arrow[equals, from=1-3, to=2-3]
        \arrow[""{name=0, anchor=center, inner sep=0}, "{(P(s\odot t),\tgt\epsilon_{s\odot t})}"'{inner sep=.8ex}, "\shortmid"{marking}, from=2-1, to=2-3]
        \arrow["{P_{s,t}}"{description, pos=0.4}, draw=none, from=1-2, to=0]
      \end{tikzcd},
    \end{equation*}
  that is, essentially just the \emph{given} comparison $P_{s,t}$ coming with 
  the original double functor $P$.
  The question is whether this is well defined according to the
  condition of \cref{rmk:elements-cell-well-definition}. This condition
  substantially simplifies, however, owing to the fact that the proposed cell
  is tightly globular and thus the tight morphisms are identities. The 
  composition of the proarrows forming the domain of the proposed cell is 
    \begin{equation*}
      % https://q.uiver.app/#q=WzAsNCxbMCwwLCIoUHRcXG9kb3QgUHMpXyEoYSkiXSxbMSwwLCJQKHQpXyFQKHMpXyEoYSkiXSxbMywwLCJQKHQpXyEoYikiXSxbNCwwLCJjIl0sWzAsMSwiXFxjb25nIl0sWzEsMiwiKFB0KV8hKFxcdGd0XFxlcHNpbG9uX3MpIl0sWzIsMywiXFx0Z3RcXGVwc2lsb25fdCJdXQ==
      \begin{tikzcd}
        {(Pt\odot Ps)_!(a)} & {P(t)_!P(s)_!(a)} && {P(t)_!(b)} & c
        \arrow["\cong", from=1-1, to=1-2]
        \arrow["{(Pt)_!(\tgt\epsilon_s)}", from=1-2, to=1-4]
        \arrow["{\tgt\epsilon_t}", from=1-4, to=1-5]
      \end{tikzcd}.
    \end{equation*}
  The instances here of targets of components of $\epsilon$ at $s$ and at $t$
  are related to $\tgt\epsilon_{s\odot t}$ in the second component of the 
  codomain of the proposed cell by the algebra associativity condition
  in the definition of a cleavage (\cref{def:cleavage-for-disc-opfibration}).
  It is precisely this condition which thus provides well-definition of the 
  proposed cell. It was in fact the question of well-definition which led to 
  the phrasing of the algebra conditions for a cleavage in the first place.
\end{construction}

\begin{construction}[Loose identity comparisons]
  Likewise, take any object $a\in\dbl{E}$. The proposed unitor comparison cell
  is given by the unitor $P_a$ coming with $a$ as in 
    \begin{equation*}
      % https://q.uiver.app/#q=WzAsNCxbMCwwLCIoeCxhKSJdLFsyLDAsIih4LGEpIl0sWzIsMSwiKHgsYSkiXSxbMCwxLCIoeCxhKSJdLFswLDEsIihcXGlkX3gsIFxca2FwcGFfYV57LTF9KSIsMCx7InN0eWxlIjp7ImJvZHkiOnsibmFtZSI6ImJhcnJlZCJ9fX1dLFsxLDIsIiIsMCx7ImxldmVsIjoyLCJzdHlsZSI6eyJoZWFkIjp7Im5hbWUiOiJub25lIn19fV0sWzAsMywiIiwyLHsibGV2ZWwiOjIsInN0eWxlIjp7ImhlYWQiOnsibmFtZSI6Im5vbmUifX19XSxbMywyLCIoUFxcaWRfYSxcXHRndFxcZXBzaWxvbl97XFxpZF9hfSkiLDIseyJzdHlsZSI6eyJib2R5Ijp7Im5hbWUiOiJiYXJyZWQifX19XSxbNCw3LCJQYSIsMSx7InNob3J0ZW4iOnsic291cmNlIjoyMCwidGFyZ2V0IjoyMH0sInN0eWxlIjp7ImJvZHkiOnsibmFtZSI6Im5vbmUifSwiaGVhZCI6eyJuYW1lIjoibm9uZSJ9fX1dXQ==
      \begin{tikzcd}
        {(x,a)} && {(x,a)} \\
        {(x,a)} && {(x,a)}
        \arrow[""{name=0, anchor=center, inner sep=0}, "{(\id_x, \kappa_a^{-1})}"{inner sep=.8ex}, "\shortmid"{marking}, from=1-1, to=1-3]
        \arrow[equals, from=1-1, to=2-1]
        \arrow[equals, from=1-3, to=2-3]
        \arrow[""{name=1, anchor=center, inner sep=0}, "{(P\id_a,\tgt\epsilon_{\id_a})}"'{inner sep=.8ex}, "\shortmid"{marking}, from=2-1, to=2-3]
        \arrow["Pa"{description}, draw=none, from=0, to=1]
      \end{tikzcd}
    \end{equation*}
  It evidently needs to be seen that this is well-defined. We do this with some
  care, as unit conditions are all too easily glossed over. Unpacking the
  well-definition condition (\cref{rmk:elements-cell-well-definition}), we have
  a diagram of set functions
    \begin{equation*}
      % https://q.uiver.app/#q=WzAsNCxbMCwwLCJcXGRibHtFfV94KGEsYSkiXSxbMCwxLCJcXGRibHtFfV94KGEsYSkiXSxbMiwxLCJcXGRibHtFfV94KGEsYSkiXSxbMiwwLCJcXGRibHtFfV94KGEsYSkiXSxbMCwxLCJcXGRibHtFfV97UFxcaWRfYX0iLDJdLFsxLDIsIihcXHRndFxcZXBzaWxvbl97XFxpZF9hfSlfISIsMl0sWzMsMiwiKFxca2FwcGFfYV57LTF9KV4qIl1d
      \begin{tikzcd}
        {\dbl{E}_x(a,a)} && {\dbl{E}_x(a,a)} \\
        {\dbl{E}_x(a,a)} && {\dbl{E}_x(a,a)}
        \arrow["{\dbl{E}_{P\id_a}}"', from=1-1, to=2-1]
        \arrow["{(\kappa_a^{-1})^*}", from=1-3, to=2-3]
        \arrow["{(\tgt\epsilon_{\id_a})_!}"', from=2-1, to=2-3]
      \end{tikzcd}.
    \end{equation*}
  Chase the arrow $1_a$ around each side. From the definition of the action of 
  the pseudo-inverse construction on cells, we have in the counterclockwise 
  direction the composite $(\tgt\epsilon_{\id_a}) \tgt\sigma(1_a,P_a)$ and
  from the clockwise direction the component arrow $(\kappa_a)^{-1}$. That 
  these are equal is precisely the unitor algebra condition in the definition 
  of the cleavage $\sigma$ coming with $P$ 
  (\cref{def:cleavage-for-disc-opfibration}).
\end{construction}

With all the data assembled, we have now the second main result of the
subsection, namely, that $R$ is a \emph{weak} equivalence of double categories.

\begin{proposition}
  The assignments and comparison data from the previous constructions amount to 
  a double functor 
    \begin{equation*}
      R\colon\dbl{E}\to\Elt(\Tw(P))
    \end{equation*}
  that commutes with $P$ and the canonical projection from the elements
  construction to $\dbl{B}$. This is an equivalence in the sense that both $R_0$
  and $R_1$ are \emph{weak} equivalences of categories.
\end{proposition}
\begin{proof}
  The assignments clearly respect the source, target and loose identity functors.
  Moreover since the loose comparison cells are those of the original $P$, the 
  assignments are double-functorial since $P$ itself is a double functor. Again 
  the work on this front consisted in showing that the comparisons were 
  well-defined. As remarked above, $R_0$ is an equivalence. On the other hand, 
  evidently, $R_1$ is fully faithful by construction of $\Elt(-)$ applied to 
  any twisted functor. The claim is proved if we show that $R_1$ is essentially 
  surjective. To this end, let $(m,\overline m)\colon (x,a)\proto (y,b)$ denote 
  any proarrow of $\Elt(\Tw(P))$ where in particular $\overline m\colon m_!(a)
  \cong b$ is an isomorphism in the fiber $\dbl{E}_y$ and $m_!(a) = \tgt\sigma(a,m)$.
  The candidate isomorphism arises from taking the image of $m\colon a\proto b$ 
  under $R$ as constructed above:
    \begin{equation*}
      % https://q.uiver.app/#q=WzAsNCxbMCwwLCIoeCxhKSJdLFsyLDAsIih5LG1fIShhKSkiXSxbMCwxLCIoeCxhKSJdLFsyLDEsIih5LGIpIl0sWzAsMSwiKFBcXHNpZ21hKGEsbSksMSkiLDAseyJzdHlsZSI6eyJib2R5Ijp7Im5hbWUiOiJiYXJyZWQifX19XSxbMCwyLCIiLDIseyJsZXZlbCI6Miwic3R5bGUiOnsiaGVhZCI6eyJuYW1lIjoibm9uZSJ9fX1dLFsyLDMsIihtLFxcb3ZlcmxpbmUgbSkiLDIseyJzdHlsZSI6eyJib2R5Ijp7Im5hbWUiOiJiYXJyZWQifX19XSxbMSwzLCIoMSxcXG92ZXJsaW5lIG0pIl0sWzQsNiwiMV9tIiwxLHsic2hvcnRlbiI6eyJzb3VyY2UiOjIwLCJ0YXJnZXQiOjIwfSwic3R5bGUiOnsiYm9keSI6eyJuYW1lIjoibm9uZSJ9LCJoZWFkIjp7Im5hbWUiOiJub25lIn19fV1d
      \begin{tikzcd}
        {(x,a)} && {(y,m_!(a))} \\
        {(x,a)} && {(y,b)}
        \arrow[""{name=0, anchor=center, inner sep=0}, "{(P\sigma(a,m),1)}"{inner sep=.8ex}, "\shortmid"{marking}, from=1-1, to=1-3]
        \arrow[equals, from=1-1, to=2-1]
        \arrow["{(1,\overline m)}", from=1-3, to=2-3]
        \arrow[""{name=1, anchor=center, inner sep=0}, "{(m,\overline m)}"'{inner sep=.8ex}, "\shortmid"{marking}, from=2-1, to=2-3]
        \arrow["{1_m}"{description}, draw=none, from=0, to=1]
      \end{tikzcd}
    \end{equation*}
  which is well-defined since on the one hand $P\sigma(a,m) = m$ holds and also 
  $\tgt\epsilon_{\sigma(a,m)}= 1$ as noted in \cref{lemma:counit-equations}.
  Thus, the domain proarrow is well-defined and in the image of $R$. The cell 
  as well is clearly invertible since $\overline m$ is an iso.
\end{proof}

\begin{remark}[Right-wobbliness]
  An equivalence of double categories (a double functor with a pseudo-inverse) 
  is equivalently given by 
    \begin{enumerate}[noitemsep]
      \item an adjoint equivalence of double categories,
      \item equivalences at the $0$- and $1$- levels, and
      \item a fully faithful double functor that is essentially surjective on
        proarrows,
    \end{enumerate}
  \emph{provided} that the double categories involved are \emph{horizontally invariant}
  (\cite[\S\S{4.1.7, 4.4.4, 4.4.5}]{grandis2019}). What the previous proposition
  describes is clearly the third of these equivalent characterizations. But there 
  is at least one fatal issue that arises in trying to create a double functor 
  in the reverse direction at least in the general case. That is, given a 
  proarrow $(m,\overline m)\colon (x,a)\proto (y,b)$ of $\Elt(\Tw(P))$ as 
  described in the proof, we have a mere corner 
    \begin{equation*}
      % https://q.uiver.app/#q=WzAsMyxbMCwwLCJhIl0sWzEsMCwibV8hKGEpIl0sWzEsMSwiYiJdLFswLDEsIlxcc2lnbWEoYSxtKSIsMCx7InN0eWxlIjp7ImJvZHkiOnsibmFtZSI6ImJhcnJlZCJ9fX1dLFsxLDIsIlxcb3ZlcmxpbmUgbSJdXQ==
      \begin{tikzcd}
        a & {m_!(a)} \\
        & b
        \arrow["{\sigma(a,m)}"{inner sep=.8ex}, "\shortmid"{marking}, from=1-1, to=1-2]
        \arrow["{\overline m}", from=1-2, to=2-2]
      \end{tikzcd}
    \end{equation*}
  that does not obviously complete to a proarrow $a\proto b$ of $\dbl{E}$ with
  the appropriate target. That is, there is neither data nor structure available
  to derive an assignment in this direction that preserves targets on the nose
  as is required of a double functor. This is a persistent
  \emph{right-wobbliness} \cite{pare2015} endemic to the theory that is
  rectified potentially by certain further assumptions. For example, a fill or
  lift of the form
    \begin{equation*}
      % https://q.uiver.app/#q=WzAsNSxbMCwwLCJhIl0sWzEsMCwibV8hKGEpIl0sWzEsMSwiYiJdLFsyLDAsImIiXSxbMiwxLCJiIl0sWzAsMSwiXFxzaWdtYShhLG0pIiwwLHsic3R5bGUiOnsiYm9keSI6eyJuYW1lIjoiYmFycmVkIn19fV0sWzEsMiwiXFxvdmVybGluZSBtIl0sWzEsMywiPz8/IiwwLHsic3R5bGUiOnsiYm9keSI6eyJuYW1lIjoiYmFycmVkIn19fV0sWzIsNCwiXFxpZF9iIiwyLHsic3R5bGUiOnsiYm9keSI6eyJuYW1lIjoiYmFycmVkIn19fV0sWzMsNCwiMSJdLFs3LDgsIj8/PyIsMSx7InNob3J0ZW4iOnsic291cmNlIjoyMCwidGFyZ2V0IjoyMH0sInN0eWxlIjp7ImJvZHkiOnsibmFtZSI6Im5vbmUifSwiaGVhZCI6eyJuYW1lIjoibm9uZSJ9fX1dXQ==
      \begin{tikzcd}
        a & {m_!(a)} & b \\
        & b & b
        \arrow["{\sigma(a,m)}"{inner sep=.8ex}, "\shortmid"{marking}, from=1-1, to=1-2]
        \arrow[""{name=0, anchor=center, inner sep=0}, "{???}"{inner sep=.8ex}, "\shortmid"{marking}, from=1-2, to=1-3]
        \arrow["{\overline m}", from=1-2, to=2-2]
        \arrow["1", from=1-3, to=2-3]
        \arrow[""{name=1, anchor=center, inner sep=0}, "{\id_b}"'{inner sep=.8ex}, "\shortmid"{marking}, from=2-2, to=2-3]
        \arrow["{???}"{description}, draw=none, from=0, to=1]
      \end{tikzcd}
    \end{equation*}
  could be supplied by the assumptions of either horizontal invariance as above 
  or of \emph{isofibrancy} (which would certainly imply horizontal invariance).
  Almost all examples of interest would satisfy some such condition, but we see 
  no special reason to impose the assumption at this point, allowing for the 
  possibility that some assignments we would like to think of as equivalences 
  are just weak and in principle may not be able to be made strong without 
  further structure.
\end{remark}

\subsection{Pseudo-algebras into opfibrations}
\label{subsection:opfibrations-from-algebras}

Just as every cloven loosely discrete opfibration 
(\cref{def:cleavage-for-disc-opfibration}) comes with a pseudo-algebra structure 
describing the action of the cleavage, from any given such algebra, a loosely 
discrete opfibration can be recovered. It is the purpose of this subsection to 
show how this can be done.

In the following, for a fixed functor $P\colon \cat{A}\to\dbl{B}_0$ between ordinary
categories, let $\cat{A}^{\cong}$ denote the category whose objects are $P$-vertical
isomorphisms and whose morphisms are pairs of morphisms of $\cat{A}$ making
appropriate commutative squares. Such squares project back to $\cat{A}$ by
taking the left or right side of such a square yielding functors
  \begin{equation*}
    \cat{A} \xleftarrow{\src}\cat{A}^{\cong}\xrightarrow{\tgt}\cat{A}.
  \end{equation*}
With suitable loose composition and identities, these are the source and
target functors of a double category obtained as a full sub-double-category
of the commutative squares in $\cat{A}$.

\begin{construction}[Corners] \label{construction:corners}
  Given a double category $\dbl{B}$, let $P\colon \cat{A}\to\dbl{B}_0$ be the
  underlying functor for a pseudo-algebra of the pull-push monad on $\dbl{B}$
  (\cref{construction:pull-push-monad}). Denote the accompanying action and
  cells by $A$, $\mu$, and $\kappa$. For compactness, write $m_!(a)$ for $A(a,m)$ and
  similarly on arrows. In particular, note that $\mu$ and $\kappa$ are \emph{fibered
    over} $\dbl{B}_0$ via $P$. In this sense, the components of these cells are
  \emph{vertical} over objects of the base category $\dbl{B}_0$.

  Let $\Ctx(P)_1$ denote the apex of the pullback square given by
    \begin{equation*}
      % https://q.uiver.app/#q=WzAsNixbMCwxLCJcXGNhdHtBfSJdLFsxLDEsIlxcY2F0e0F9XFx0aW1lc197XFxkYmx7Qn1fMH1cXGRibHtCfV8xIl0sWzIsMSwiXFxjYXR7QX0iXSxbMywxLCJcXGNhdHtBfV57XFxjb25nfSJdLFs0LDEsIlxcY2F0e0F9Il0sWzIsMCwiXFxDdHgoUClfMSJdLFsxLDAsIlxccGlfMSJdLFsxLDIsIkEiLDJdLFszLDIsIlxcc3JjIl0sWzUsMV0sWzUsM10sWzUsMiwiIiwyLHsic3R5bGUiOnsibmFtZSI6ImNvcm5lciJ9fV0sWzMsNCwiXFx0Z3QiLDJdXQ==
      \begin{tikzcd}
        && {\Ctx(P)_1} \\
        {\cat{A}} & {\cat{A}\times_{\dbl{B}_0}\dbl{B}_1} & {\cat{A}} & {\cat{A}^{\cong}} & {\cat{A}}
        \arrow[from=1-3, to=2-2]
        \arrow["\lrcorner"{anchor=center, pos=0.125, rotate=-45}, draw=none, from=1-3, to=2-3]
        \arrow[from=1-3, to=2-4]
        \arrow["{\pi_1}", from=2-2, to=2-1]
        \arrow["A"', from=2-2, to=2-3]
        \arrow["\src", from=2-4, to=2-3]
        \arrow["\tgt"', from=2-4, to=2-5]
      \end{tikzcd}.
    \end{equation*}
  Thus, $\Ctx(P)_1$ is the category whose objects are triples 
  $(a,m,\overline{m})$, where $a\in \cat{A}$ is over $x$, $m\colon x\proto y$ is
  a proarrow in $\dbl{B}$, and $\overline m \colon m_!(a)\cong b$ is a
  $P$-vertical isomorphism. This data can be visualized as a kind of corner:
    \begin{equation*}
      % https://q.uiver.app/#q=WzAsMyxbMCwwLCJhIl0sWzEsMCwibV8hKGEpIl0sWzEsMSwiYiJdLFsxLDIsIlxcYmFyIG0iXV0=
      \begin{tikzcd}
        a & {m_!(a)} \\
        & b
        \arrow["{\bar m}", from=1-2, to=2-2]
      \end{tikzcd}
      \qquad\xmapsto{P}\qquad
      % https://q.uiver.app/#q=WzAsMyxbMCwwLCJ4Il0sWzEsMCwieSJdLFsxLDEsInkiXSxbMCwxLCJtIiwwLHsic3R5bGUiOnsiYm9keSI6eyJuYW1lIjoiYmFycmVkIn19fV0sWzEsMiwiIiwwLHsibGV2ZWwiOjIsInN0eWxlIjp7ImhlYWQiOnsibmFtZSI6Im5vbmUifX19XV0=
      \begin{tikzcd}
        x & y \\
        & y
        \arrow["m"{inner sep=.8ex}, "\shortmid"{marking}, from=1-1, to=1-2]
        \arrow[equals, from=1-2, to=2-2]
      \end{tikzcd}.
    \end{equation*}
  Via the left projection, we have a functor $P_1$ to $\dbl{B}_1$ as well as
  purported source and target functors:
    \begin{equation*}
      % https://q.uiver.app/#q=WzAsNSxbMCwwLCJcXEN0eChQKV8xIl0sWzIsMSwiXFxkYmx7Qn1fMSJdLFsyLDIsIlxcZGJse0J9XzAiXSxbMSwyLCJcXGNhdHtBfSJdLFsxLDEsIlxcY2F0e0F9XFx0aW1lc197XFxkYmx7Qn1fMH1cXGRibHtCfV8xIl0sWzAsMSwiUF8xIiwwLHsiY3VydmUiOi0zfV0sWzEsMiwiXFxzcmMiXSxbMCwzLCJcXHNyYyIsMix7ImN1cnZlIjo1fV0sWzMsMiwiUF8wIiwyXSxbNCwzLCJcXHBpXzEiLDJdLFs0LDEsIlxccGlfMiJdLFswLDRdLFs0LDIsIiIsMix7InN0eWxlIjp7Im5hbWUiOiJjb3JuZXIifX1dXQ==
      \begin{tikzcd}
        {\Ctx(P)_1} \\
        & {\cat{A}\times_{\dbl{B}_0}\dbl{B}_1} & {\dbl{B}_1} \\
        & {\cat{A}} & {\dbl{B}_0}
        \arrow[from=1-1, to=2-2]
        \arrow["{P_1}", curve={height=-18pt}, from=1-1, to=2-3]
        \arrow["\src"', curve={height=30pt}, from=1-1, to=3-2]
        \arrow["{\pi_2}", from=2-2, to=2-3]
        \arrow["{\pi_1}"', from=2-2, to=3-2]
        \arrow["\lrcorner"{anchor=center, pos=0.125}, draw=none, from=2-2, to=3-3]
        \arrow["\src", from=2-3, to=3-3]
        \arrow["{P_0}"', from=3-2, to=3-3]
      \end{tikzcd}
      \qquad\qquad\qquad
      % https://q.uiver.app/#q=WzAsMyxbMCwwLCJcXEN0eChQKV8xIl0sWzEsMSwiXFxjYXR7QX0iXSxbMCwxLCJcXGNhdHtBfVxcdGltZXNfe1xcZGJse0J9XzB9XFxkYmx7Qn1fMSJdLFswLDEsIlxcdGd0IiwwLHsiY3VydmUiOi0xfV0sWzIsMSwiQSIsMl0sWzAsMl1d
      \begin{tikzcd}
        {\Ctx(P)_1} \\
        {\cat{A}\times_{\dbl{B}_0}\dbl{B}_1} & {\cat{A}}
        \arrow[from=1-1, to=2-1]
        \arrow["\tgt", curve={height=-6pt}, from=1-1, to=2-2]
        \arrow["A"', from=2-1, to=2-2]
      \end{tikzcd}.
  \end{equation*}
  Together with loose composition and identities, $\Ctx(P)_1$ will then be the
  category of proarrows and cells of a double category with
  underlying category $\Ctx(P)_0 \coloneqq \cat{A}$ admitting a double functor
  to $\dbl{B}$. Anticipating this result, display the triples
  $(a,m,\overline m)$ giving the objects of $\Ctx(P)_1$ in a proarrow-style
  notation $(m,\overline{m})\colon a\proto b$. Loose identities are then those
  proarrows $(\id_{Pa},\kappa_a) \colon a\proto a$ using $\kappa_a$, a component of the unit
  isomorphism coming with the pseudo-algebra structure. Now, given proarrows
  $(m,\overline{m})\colon a\proto b$ and $(n,\overline{n})\colon b\proto c$ that are
  \emph{composable} in the sense that $\tgt(m) = \src (n)$ and
  $\overline m\colon m_!(a)\cong b$, define the \emph{loose composite} in a familiar
  category-of-elements style:
    \begin{equation*}
      (a,m,\overline{m})\odot (b,n,\overline{n}) \coloneqq (a,m\odot n,\overline{m\odot n}),
    \end{equation*}
  where
    \begin{equation*}
      \overline{m\odot n}:(m\odot n)_!(a)\cong n_!m_!(a) \xto{n_!(\bar m)} n_!(b) \xto{\bar n} c,
    \end{equation*}
  using the $a$-component of the algebra isomorphism $\mu$. Letting $\Ctx(P)_0
  \coloneqq\cat{A}$, these operations define a double category $\Ctx(P)$ and the
  associated functors $P_0$ and $P_1$ underlie a double functor 
  $P\colon \Ctx(P)\to\dbl{B}$ by construction.
\end{construction}

\begin{proposition}
  The double functor $P\colon \Ctx(P)\to\dbl{B}$ of \cref{construction:corners} is a cloven loosely discrete opfibration as in \cref{def:cleavage-for-disc-opfibration}.
\end{proposition}
\begin{proof}
  We shall first show that the projection morphism
    \begin{equation*}
      \langle\src,P_1\rangle\colon\Ctx(P)_1 \to \cat{A}\times_{\dbl{B}_0}\dbl{B}_1
    \end{equation*}
  is the right adjoint of a LARI split equivalence of categories. The left adjoint will be
    \begin{equation*}
      \sigma\colon\cat{A}\times_{\dbl{B}_0}\dbl{B}_1\to\Ctx(P)_1,
      \qquad\qquad
      \sigma(a,m) \;\coloneqq\; (m,1_{m_!(a)})\colon a\proto m_!(a),
    \end{equation*}
  and extended suitably to morphisms. By definition $\sigma$ is thus a right 
  inverse to the projection under consideration. This splitting provides the 
  unit of the adjunction underlying the purported equivalence. Required is thus 
  an invertible counit. Define its components in the following way. Given an 
  object of $\Ctx(P)_1$, say, $(m,\overline m)\colon a\proto b$, we have an 
  invertible morphism of $\Ctx(P)_1$, namely,
    \begin{equation*}
      \epsilon_{(a,m,\overline m)}
      \;\coloneqq\;
      (1_a, 1_m, (1_{m_!(a)}, \overline{m}))
      \;\leftrightsquigarrow\;
      % https://q.uiver.app/#q=WzAsNCxbMCwwLCJhIl0sWzEsMCwibV8hKGEpIl0sWzAsMSwiYSJdLFsxLDEsImIiXSxbMCwxLCIobSwxKSIsMCx7InN0eWxlIjp7ImJvZHkiOnsibmFtZSI6ImJhcnJlZCJ9fX1dLFswLDIsIiIsMix7ImxldmVsIjoyLCJzdHlsZSI6eyJoZWFkIjp7Im5hbWUiOiJub25lIn19fV0sWzIsMywiKG0sXFxvdmVybGluZSBtKSIsMix7InN0eWxlIjp7ImJvZHkiOnsibmFtZSI6ImJhcnJlZCJ9fX1dLFsxLDMsIlxcb3ZlcmxpbmUgbSJdLFs0LDYsIjFfbSIsMSx7InNob3J0ZW4iOnsic291cmNlIjoyMCwidGFyZ2V0IjoyMH0sInN0eWxlIjp7ImJvZHkiOnsibmFtZSI6Im5vbmUifSwiaGVhZCI6eyJuYW1lIjoibm9uZSJ9fX1dXQ==
      \begin{tikzcd}
        a & {m_!(a)} \\
        a & b
        \arrow[""{name=0, anchor=center, inner sep=0}, "{(m,1)}"{inner sep=.8ex}, "\shortmid"{marking}, from=1-1, to=1-2]
        \arrow[equals, from=1-1, to=2-1]
        \arrow["{\overline m}", from=1-2, to=2-2]
        \arrow[""{name=1, anchor=center, inner sep=0}, "{(m,\overline m)}"'{inner sep=.8ex}, "\shortmid"{marking}, from=2-1, to=2-2]
        \arrow["{1_m}"{description}, draw=none, from=0, to=1]
      \end{tikzcd}
    \end{equation*}
  which we take to be the counit component. These are suitably natural. Note 
  that $\overline m$ is of course invertible and $P$-vertical. This is thus an 
  appropriately fibered adjunction and indeed a LARI equivalence of categories. 
  It remains to check the coherence conditions of 
  \cref{def:cleavage-for-disc-opfibration}. First is that the
  composite $\tgt\sigma$ is the structure map for a pseudo algebra. But by
  definition this composite is just $A$, the structure map for the given 
  algebra. The further associativity and unit conditions follow by the 
  construction of $\epsilon$ and the definitions of the loose composition and 
  units in the double category $\Ctx(P)$.
\end{proof} 

The foregoing shows that from any pseudo-algebra for the canonical pull-push 
monad, an opfibration can be recovered. In fact, there is a more complete
correspondence between the two concepts.

\begin{corollary}
  Every cloven loosely discrete opfibration arises up to equivalence as the corner construction of a pseudo-algebra for the appropriate pull-push monad. In this sense, there is thus a correspondence 
    \begin{equation*}
      \lbrace\text{pseudo-algebras for the pull-push monad}\rbrace\quad\leftrightarrow\quad\lbrace\text{cloven loosely discrete opfibrations}\rbrace    
    \end{equation*}
  that is strict on the side of algebras and up to equivalence on opfibrations.
\end{corollary}
\begin{proof}
  Clearly in one direction, we assign the corner construction above to each pseudo-algebra to obtain an opfibration. In the reverse direction, we take each cloven opfibration to the associated pseudo-algebra. Now, on the one hand, starting with an algebra $P\colon\cat{A}\to\dbl{B}_0$, the algebra associated to $\Ctx(P)$ is exactly the original $P$. On the other hand, starting with an 
  opfibration $P\colon\dbl{E}\to\dbl{B}$, the associated algebra has structure morphism $P_0\colon\dbl{E}_0\to\dbl{B}_0$ and $\Ctx(P_0)_0 \coloneqq \dbl{E}_0$ and crucially $\Ctx(P_0)_1$ was shown above to be suitably equivalent to the same pullback that $\dbl{E}_1$ is equivalent to by the assumption of being an opfibration. These equivalences assemble into a suitably fibered equivalence
    \begin{equation*}
      % https://q.uiver.app/#q=WzAsMyxbMCwwLCJcXEN0eChQXzApIl0sWzIsMCwiXFxkYmx7RX0iXSxbMSwxLCJcXGRibHtCfSJdLFswLDEsIlxcc2ltZXEiXSxbMCwyLCJQXzAiLDJdLFsxLDIsIlAiXV0=
      \begin{tikzcd}
        {\Ctx(P_0)} && {\dbl{E}} \\
        & {\dbl{B}}
        \arrow["\simeq", from=1-1, to=1-3]
        \arrow["{P_0}"', from=1-1, to=2-2]
        \arrow["P", from=1-3, to=2-2]
      \end{tikzcd}
    \end{equation*}
  in the category of loosely discrete opfibrations.
\end{proof}

Combining this result with other work in the prior subsection, we deduce a final
result tying together the various correspondences.

\begin{corollary}
  Twisted copresheaves on a double category $\dbl{B}$
  correspond to pseudo-algebras for the pull-push monad on $\dbl{B}$.
\end{corollary}

This explains the elements-style definition of composition in the corners
construction: each occurs (essentially) as an elements construction of a twisted
copresheaf.

\begin{remark}[Outlook on double toposes]
  This corollary parallels the well-known results, summarized in the
  introduction, that discrete fibrations, ordinary fibrations, discrete
  2-fibrations and other notions of fibration occur as algebras (of some kind)
  for suitable pull-push monads. At least in the case of discrete fibrations,
  this leads to an elementary axiomatization of an internal (co)presheaf and a
  development of a completely internal elementary theory of \emph{filteredness}
  and \emph{flatness} (\cite[\S 2.5, \S 4.3]{johnstone1977}, \cite[\S B1.4, \S
  B3.2]{johnstone2002}), axiomatizing the well-known theory for Grothendieck
  toposes \cite[Ch. VII]{maclane1992} in a manner free of the assumption of a
  fixed set-theoretic universe.

  That is, the ordinary theory for Grothendieck toposes shows that geometric
  morphisms into a presheaf topos correspond to flat functors from the
  represented 1-category in the sense that there is an equivalence of categories
  \begin{equation*}
    \cat{Geom}(\cat{E},[\cat{C}^{\op},\Set])\simeq \cat{Flat}(\cat{C},\cat{E})
  \end{equation*}
  for any topos $\cat{E}$ \cite[Theorem VII.7.2]{maclane1992}. Such flat
  functors $\cat{C}\to\cat{E}$ are shown to be \emph{filtering} \cite[Theorem
  VII.9.1]{maclane1992}, which gives a tractable description of left-exactness
  in this context and leads to important applications in understanding morphisms of
  sites, classifying toposes, torsors, nonabelian cohomology, and geometric
  theories and their models. The point of the elementary theory is to replace
  $\Set$ by an ambient topos $\cat{F}$ thought of as the new base and as a
  self-contained universe of discourse. Likewise, $\cat{C}$ is replaced by an
  internal category. However, base-valued functors on $\cat{C}$ are not
  axiomatizable in the same manner as ordinary presheaves. Rather, the internal
  category acts on $\cat F$ via a pull-push monad, and $[\cat{C}^{\op},\Set]$ is
  replaced by the category of algebras for this monad. The notion of how such an
  algebra might be filtered is introduced and it is shown that such
  flat/filtering algebras correspond to geometric morphisms into the category of
  algebras for the pull-push monad. This is \emph{Diaconescu's theorem}
  \cite{diaconescu1975} and leads to many examples and applications facilitating
  our understanding of toposes and their properties \cite[\S
  4.3]{johnstone1977}.
  
  Thus, in summary, via our results on loosely discrete opfibrations, we are led
  directly to an elementary internal axiomatization of twisted representations
  and twisted copresheaves of a fixed double category viewed as a
  pseudo-category internal to the 2-category of categories. In this way the
  foundation is laid for a potential investigation of double-categorical
  \emph{filteredness} and \emph{flatness} and higher-dimensional \emph{geometric
    theories} which ought to be concepts central to both \emph{geometric} and
  \emph{elementary double toposes}.
\end{remark}

\paragraph{Acknowledgments}

Patterson was supported by the Air Force Office of Scientific Research (AFOSR)
Young Investigator Program (YIP) through Award FA9550-23-1-0133.

\printbibliography[heading=bibintoc]

\end{document}